\documentclass[12pt,letterpaper]{amsbook}
\usepackage{amsmath}
\usepackage{epic,eepic,latexsym, amssymb, amscd, amsfonts, xypic, euler}
\usepackage{color}
\usepackage{bbm}
\usepackage{enumerate}
\usepackage{setspace}

\usepackage{perpage} 
\MakePerPage{footnote}

\newcommand{\comments}[1]{}

\usepackage{hyperref}

 \newlength{\baseunit}               
 \newcount{\numlines}                
\newtheorem{defn}{Definition}
\newtheorem{remk}{Remark}
\newtheorem{thm}{Theorem}
\newtheorem{lem}{Lemma}

\newtheorem{cor}{Corollary}
\newtheorem{prop}{Proposition}

\newcommand{\environmentA}[1]{\noindent \underline{\textsc{#1}}}

\newcommand{\oleq}[1]{\overset{ {\scriptscriptstyle #1}}{\leq}}
\newcommand{\ogeq}[1]{\overset{ {\scriptscriptstyle #1}}{\geq}}

\def \N {\mathbb N}

\def \cA {\mathcal A}
\def \cQ {\mathcal Q}

\def \cB {\mathcal B}
\def \cC {\mathcal C}

\def \cI {\mathcal I}
\def \cL {\mathcal L}

\def \cJ {\mathcal J}

\def \cO {\mathcal O}
\def \cS {\mathcal S}
\def \cP {\mathcal P}

\def \cM {\mathcal M}

\def \oa {{\overline{\alpha}}}
\def \ob {{\overline{\beta}}}
\def \oA {{\overline{\mathcal A}}}

\def \oQ {{\overline{Q}}}
\def \oP {{\overline{P}}}

\def \hQ {{\widehat{Q}}}
\def \upQ {Q^{\up}}
\def \hmu {{\widehat{\mu}}}
\def \hP {{\widehat{P}}}

\def \path {\mbox{\textbf{path}}}

\def \lv {{\left\lvert}}
\def \rv {{\right\rvert}}

\def \ouch {{\text{\tiny ouch}}}
\def \ind {{\text{\tiny indicated}}}
\def \up {{\text{\tiny up}}}
\def \go {{\text{\tiny go}}}
\def \st {{\text{\tiny start}}}
\def \en {{\text{\tiny end}}}
\def \stay {{\text{\tiny stay}}}
\def \spl {{\text{\tiny split}}}
\def \cover {{\text{\tiny cov}}}
\def \near {{\text{\tiny near}}}
\def \loc {{\text{\tiny loc}}}
\def \new {{\text{\tiny new}}}
\def \old {{\text{\tiny old}}}
\def \spec {{\text{\tiny sp}}}

\def \bad {{\text{\tiny bad}}}

\def \extra {\text{extra}}

\def \cov {{\cI_\cover}}

\def \lv {\left\lvert}
\def \rv {\right\rvert}
\newcommand{\main}{\operatorname{main}}
\newcommand{\depth}{\operatorname{depth}}
\newcommand{\Vol}{\operatorname{Vol}}

\newcommand{\rt}{\operatorname{root}}

\newcommand{\Time}{\operatorname{Time}}
\newcommand{\Work}{\mathfrak{W}}
\newcommand{\Space}{\mathfrak{S}}
\newcommand{\NULL}{\operatorname{NULL}}

\newcommand{\sgn}{\operatorname{sgn}}
\newcommand{\dist}{\operatorname{dist}}
\newcommand{\gchild}{\operatorname{gochild}}
\newcommand{\schild}{\operatorname{staychild}}
\newcommand{\supp}{\operatorname{supp}}
\newcommand{\diam}{\operatorname{diam}}
\newcommand{\CT}{\operatorname{CT}}
\newcommand{\BD}{\operatorname{BD}}
\newcommand{\BT}{\operatorname{BT}}
\newcommand{\CZ}{\operatorname{CZ}}
\newcommand{\DC}{\operatorname{DC}}
\newcommand{\NR}{\operatorname{NR}}

\newcommand{\LLC}{\operatorname{LLC}}
\newcommand{\URC}{\operatorname{URC}}

\newcommand{\junk}{{\operatorname{junk}}}

\newcommand{\fin}{{\operatorname{  fin}}}

\newcommand{\ooline}[1]{\overline{\overline{#1}}}

\newcommand{\til}{\widetilde}

\newcommand{\Z}{\mathbb{Z}}

\newcommand{\X}{\mathbb{X}}

\newcommand{\R}{\mathbb{R}}
\newcommand{\T}{\mathbb{T}}

\newcommand{\desc}{\operatorname{Descendants}}
\newcommand{\ndesc}{\operatorname{Nondescendants}}
\newcommand{\DR}{\mathcal{DR}}

\numberwithin{remk}{section}
\numberwithin{equation}{section}
\numberwithin{thm}{section}
\numberwithin{defn}{section}
\numberwithin{lem}{section}
\numberwithin{cor}{section}
\numberwithin{prop}{section}
\numberwithin{alg}{section}

\begin{document}
\pagestyle{plain}
\title{Fitting a Sobolev function to data}
\author{{\large Charles Fefferman, Arie Israel, and Garving K. Luli\footnote{The first author is supported in part by NSF grant DMS-1265524 and AFOSR grant FA9550-12-1-0425\\
The second author is supported in part by an NSF postdoctoral fellowship, DMS-1103978. \\
The third author is supported in part by NSF grant DMS-1355968.}}}

 \maketitle

\begin{abstract} We exhibit an algorithm to solve the following extension problem: Given a finite set $E \subset \mathbb{R}^n$ and a function $f: E \rightarrow \mathbb{R}$, compute an extension $F$ in the Sobolev space $L^{m,p}(\mathbb{R}^n)$, $p>n$, with norm having the smallest possible order of magnitude, and secondly, compute the order of magnitude of the norm of $F$. Here, $L^{m,p}(\mathbb{R}^n)$ denotes the Sobolev space consisting of functions on $\mathbb{R}^n$ whose $m$th order partial derivatives belong to $L^p(\mathbb{R}^n)$. The running time of our algorithm is at most $C N \log N$, where $N$ denotes the cardinality of $E$, and $C$ is a constant depending only on $m$,$n$, and $p$.

 \end{abstract}

 \tableofcontents

\chapter{Introduction}

In this paper, we interpolate data by a function $F:\R^n \rightarrow \R$ whose Sobolev norm has the least possible order of magnitude. Our computations consist of efficient algorithms, to be implemented on an (idealized) von Neumann computer.

Our results are the analogues for Sobolev spaces of some of the main results of Fefferman-Klartag \cite{F6,FK1,FK2} on interpolation of data by functions in $C^m(\R^n)$.

Let us set up notation and definitions. Fix $m,n \geq 1$ and $1< p< \infty$. We work in the Sobolev space 
\begin{align}
\X=L^{m,p}(\R^n)\text{ with seminorm }||F||_{\X}=\left( \int_{\R^n} \sum_{|\alpha|=m} \left| \partial^\alpha F(x) \right|^p dx \right)^{1/p} \label{intro1}\end{align}
or 
\begin{align}
\X=W^{m,p}(\R^n)\text{ with norm }||F||_{\X}=\left( \int_{\R^n} \sum_{|\alpha|\leq m} \left| \partial^\alpha F(x) \right|^p dx \right)^{1/p}. \label{intro2}\end{align}

We make the assumption \begin{align} p>n, \label{intro3}\end{align} so that the Sobolev theorem tells us that  \begin{align} \X \subset C_{loc}^{m-1}(\R^n).\label{intro4}\end{align}

We write $c, C, C'$, etc. to denote ``universal constants," i.e., constants determined by $m,n,p$ in \eqref{intro1},\eqref{intro2}. These symbols may denote different universal constants in different occurrences. 

Now let 
 \begin{align}  E = \{z_1,\cdots,z_N\} \subset \R^n. \label{intro5} \end{align}

Then $\X(E)$ denotes the vector space of all real-valued functions on $E$, equipped with the norm (or seminorm) 
$$||f||_{\X(E)}=\inf \{||F||_{\X}: F\in \X, F=f \text{ on } E \}.$$

Let $A \geq 1$ be a constant. An ``$A$-optimal interpolant'' for a function $f \in \X(E)$ is a function $F \in \X$ that satisfies $F=f$ on $E$ and $||F||_{\X} \leq A \cdot ||f||_{\X(E)}$. 

Our goal here is to solve the following

\underline{Problems:}
\begin{itemize}
\item[(A)] Compute a $C$-optimal interpolant for a given function $f \in \X(E)$.
\item[(B)] Given $f \in \X(E)$, compute a number $|||f|||$ such that $$c|||f||| \leq ||f||_{\X(E)} \leq C|||f|||.$$
\end{itemize}

We owe the reader an explanation of what it means to ``compute a function" in Problem (A). First of all, our computations are performed on a computer with standard von Neumann architecture. We assume that each memory cell can store a single integer or real number. We study two distinct models of computation. In the first model (``infinite-precision") we assume that our computer deals with exact real numbers, without roundoff errors. Our second, more realistic model (``finite-precision") assumes that our machine handles only $S$-bit machine numbers for some fixed, large $S$. To work with the finite-precision model, we make a rigorous study of the roundoff errors arising in our algorithms. For simplicity, in this introduction, we restrict attention to the infinite-precision model.

To ``compute" a function $F \in C_{loc}^{m-1}(\R^n)$, the computer first performs ``one-time work", then answers ``queries." A query consists of a point $\underline{x} \in \R^n$, and the computer responds to a query $\underline{x}$ by computing $\partial^{\alpha}F(\underline{x})$ for all $|\alpha|\leq m-1$.

We want algorithms that make minimal use of the resources of our computer. For the computation of a function $F$ as in Problem (A), the relevant resources are 
\begin{itemize}
\item The number of computer operations used for the one-time work,
\item The number of memory cells used for all the work, and 
\item The number of computer operations used in responding to a query. 
\end{itemize}

We refer to these as the ``one-time work", the ``space" (or ``storage"), and the ``query work", respectively.

For the computation of the single number $|||f|||$ in Problem (B), the relevant computer resources are the number of operations used, and the number of memory cells required. We refer to these as, respectively, the ``work" and ``storage".

We are concerned with algorithms that work for arbitrary $f$ and $E$. If we allowed ourselves favorable geometric assumptions on $E$, our problems would be much easier. 

We can now state our results in their simplest form. Recall that $N$ denotes the number of points in our finite set $E$.

\begin{thm}\label{introThm1}
One can compute a $C$-optimal interpolant for a given function $f \in \X(E)$, with one-time work $\leq CN \log N$, storage $\leq CN$, and query work $\leq C \log N$.  
\end{thm}

\begin{thm}\label{introThm2}
Given $f\in \X(E)$, one can compute a number $|||f|||$ such that 
$$c|||f||| \leq ||f||_{\X(E)} \leq C|||f|||;$$
the computation uses work $\leq CN \log N$ and storage $\leq CN$.  
\end{thm}

Obviously, in Theorem \ref{introThm1}, the one-time work must be at least $N$, since we have to read the data; and the query work is at least $1$, since we must at least read the query. Similarly, in Theorem \ref{introThm2}, the work must be at least $N$. Hence, for trivial reasons, the work of our algorithms can be improved at most by a factor $\log N$. 

Very likely, the work and storage asserted above are best possible.

To prepare to state our results in their full strength, we recall the following results from our previous paper \cite{FIL1}. 

\begin{thm}[Extension Operators]\label{introThm3}
There exists a linear map $T:\X(E) \rightarrow \X$ such that $Tf$ is a $C$-optimal interpolant of $f$ for any $f \in \X(E)$. 
\end{thm}

\begin{thm}[Formula for the Norm]\label{introThm4}
There exist linear functionals $\xi_l:\X(E) \rightarrow \R$ ($l=1,\cdots,L$) such that 
\begin{itemize}
\item $L \leq CN$ and 
\item The quantity $$|||f|||=\left( \sum_{l=1}^L |\xi_l(f)|^p\right)^{1/p}$$ satisfies $$c|||f||| \leq ||f||_{\X(E)} \leq C|||f|||$$ for all $f \in \X(E)$. 
\end{itemize} 
\end{thm}

To prove Theorems \ref{introThm1} and \ref{introThm2}, we will compute the linear map $T$ and the functionals $\xi_l$ in Theorems \ref{introThm3} and \ref{introThm4}. To do so, we exploit a sparse structure for $T$ and $\xi_l$, established in \cite{FIL1}.

We recall the relevant definitions. 

Let $\Omega = \{ \omega_1, \cdots, \omega_{\nu_{\max}}\}$ be a finite list of linear functionals on $\X(E)$. Then we say that $\Omega$ is a ``set of assists" if each $\omega _{\nu }$ can be written as 
\begin{equation}
\omega _{\nu }\left( f\right) =\sum_{i=1}^{I_{\nu }}\mu _{\nu i}f\left(
z_{\nu i}\right) \quad\text{ (}f\in \X(E)\text{);}  \label{intro6}
\end{equation}
where $I_\nu \geq 1$, $\mu_{\nu i} \in \R$, $z_{\nu i} \in E$ are independent of $f$, and
\begin{equation}
\sum_{\nu =1}^{\nu _{\max }}I_{\nu }\leq CN.  \label{intro7}
\end{equation}

The point is that if (\ref{intro7}) holds, then for a given $f\in \X(E)$ we
can compute all the assists $\omega _{1}\left( f\right) ,\cdots ,\omega
_{\nu _{\max }}\left( f\right) $, using at most $CN$ computer operations. It
will be useful to precompute the $\omega _{\nu }\left( f\right) $, because
each of these quantities may be used many times in subsequent calculations.

Let $\Omega =\left\{ \omega _{1},\cdots ,\omega _{\nu _{\max }}\right\} $ be
a set of assists.

A linear functional 
\[
\xi :\X\left( E\right) \rightarrow \mathbb{R}
\]%
has \textquotedblleft $\Omega $-assisted bounded depth" if it can be written
in the form 
\begin{equation}
\xi \left( f\right) =\sum_{i=1}^{I}\lambda _{i}f\left( z_{i}\right)
+\sum_{j=1}^{J}\beta _{j}\omega _{\nu _{j}}\left( f\right) \text{ for all }%
f\in \X\left( E\right) ,  \label{intro8}
\end{equation}%
where $I,J,\lambda _{i},\beta _{j},\nu _{j}$ and $z_{i}\in E$ are
independent of $f$, and 
\begin{equation}
I+J\leq C.  \label{intro9}
\end{equation}

If (\ref{intro8}),(\ref{intro9}) hold, and if we have precomputed $\omega
_{1}\left( f\right) ,\cdots ,\omega _{\nu _{\max }}\left( f\right) $, then
we can calculate $\xi \left( f\right) $ using at most $C$ computer
operations.

We call (\ref{intro6}) and (\ref{intro8}) \textquotedblleft short forms"
of the assists $\omega _{\nu }$ and the functional $\xi $, respectively. Note that a functional may be written in short form in more than
one way.

A linear map $T:\X\left( E\right) \rightarrow C_{loc}^{m-1}\left( \mathbb{R}%
^{n}\right) $ will be said to have \textquotedblleft $\Omega $-assisted
bounded depth" if for each $x\in \mathbb{R}^{n}$ and each multiindex $\alpha 
$ of order $\left\vert \alpha \right\vert \leq m-1$, the linear functional 
\begin{equation}
\X(E)\ni f\longmapsto \partial ^{\alpha }Tf\left( x\right)   \label{intro10}
\end{equation}%
has $\Omega $-assisted bounded depth.

In \cite{FIL1}, we proved the following sharper version of Theorems \ref{introThm3} and \ref{introThm4}.

\begin{thm}
\label{introThm5}There exists a set of assists $\Omega =\left\{ \omega
_{1},\cdots ,\omega _{\nu _{\max }}\right\} $ such that the linear map $T$
in Theorem \ref{introThm3}, and the linear functionals $\xi _{1},\cdots
,\xi _{L}$ in Theorem  \ref{introThm4}, may be taken to have $\Omega $%
-assisted bounded depth.
\end{thm}

If we knew the assists $\omega _{1},\cdots ,\omega _{\nu _{\max }}$ and the
functionals $\xi _{1},\cdots ,\xi _{L}$ in their short form, then we could
easily compute $\left\vert \left\vert \left\vert f\right\vert \right\vert
\right\vert $ in Theorem \ref{introThm4} by first computing the $\omega
_{\nu }\left( f\right) $, then computing the $\xi _{l}\left( f\right)$.
The whole computation would require only $CN$ computer operations.

Similarly, suppose we knew the assists $\omega _{\nu }$ and the linear
functionals (\ref{intro10}) in their short form. 

Given $f\in \X\left( E\right) $, we could precompute $\omega _{1}\left(
f\right) ,\cdots ,\omega _{\nu _{\max }}\left( f\right) $ with at most $CN$
operations, after which we could answer queries:\ Given a query point $x\in 
\mathbb{R}^{n}$, we could compute $\partial ^{\alpha }Tf\left( x\right) $
(all $\left\vert \alpha \right\vert \leq m-1$) in at most $C$ operations.
Thus, we could give highly efficient solutions to Problems\ (A)\ and (B)
above.

Unfortunately, the proof of Theorem \ref{introThm5} in \cite{FIL1} is not constructive. It does not supply any formulas for the
assists $\omega _{\nu }$, the functionals $\xi _{l}$, or the operator $T$.
Our purpose here is to remedy this defect by proving the following result.

\begin{thm}[MAIN THEOREM]
\label{introMainTheorem}For suitable $\Omega =\left\{ \omega _{1},\cdots
,\omega _{\nu _{\max }}\right\} $, $T,\xi _{1},\cdots ,\xi _{L}$ as in
Theorems \ref{introThm3}, \ref{introThm4}, \ref{introThm5}, the assists $%
\omega _{\nu }$, and the functionals $\xi _{l}$ can all be computed in short
form, using work $\leq CN\log N$ and storage $\leq CN$. Moreover, after
one-time work $\leq CN\log N$ in space $\leq CN$, we can answer queries as
follows: 

A query consists of a point $x\in \mathbb{R}^{n}$. The response to
a query $x$ is a short-form description of the functional (\ref{intro10})
for each $\left\vert \alpha \right\vert \leq m-1.$ The work to answer a
query is $\leq C\log N.$
\end{thm}

To prove Theorem \ref{introMainTheorem}, we modify the proofs of Theorems \ref%
{introThm3}, \ref{introThm4}, \ref{introThm5} in \cite{FIL1}. Let us first review some of the ideas in \cite{FIL1}, and then explain some of the modifications needed for Theorem \ref%
{introMainTheorem}. Our discussion will be highly simplified, so that the basic
ideas are not obscured by technical details.

We introduce a bit more notation. If $F\in C_{loc}^{m-1}\left( \mathbb{R}%
^{n}\right) $ and $x\in \mathbb{R}^{n}$, then we write $J_{x}\left( F\right) 
$ (the \textquotedblleft jet" of $F$ at $x$) to denote the $\left(
m-1\right) ^{\text{st}}$ degree Taylor polynomial of $F$ at $x.$ Thus, $%
J_{x}\left( F\right) $ belongs to $\mathcal{P}$, the vector space of all
(real) polynomials of degree at most $\left( m-1\right) $ on $\mathbb{R}^{n}$%
.

We write $Q,Q^{\prime }$, etc. to denote cubes in $\mathbb{R}^{n}$ with
sides parallel to the coordinate axes. We write $\delta _{Q}$ to denote the
sidelength of a cube $Q.$

Our review of \cite{FIL1} starts with a local version of
our present Problem (A). Let $Q\subset \mathbb{R}^{n}$ be a cube, let $%
x_{0}\in Q$ be a point, and let $P_{0}\in \mathcal{P}$ be a jet. We pose the
following Local Interpolation Problem\ :\ 

\underline{$LIP\left( Q,E,f,x_{0},P_{0}\right) :$} Find a function $F\in
L^{m,p}\left( \mathbb{R}^{n}\right)$, depending linearly on $\left(
f,P_{0}\right) $, such that 
\begin{eqnarray*}
&&F =f\text{ on }E\cap Q\text{,}\\
&&J_{x_{0}}\left( F\right)  =P_{0}\text{, and } \\
&&\int_{Q}\sum_{\left\vert \alpha \right\vert = m}\left\vert \partial
^{\alpha }F\left( x\right) \right\vert ^{p}dx\text{ is as small as possible
up to a factor }C\text{.}
\end{eqnarray*}

If we can solve $LIP\left( Q,E,f,x_{0},P_{0}\right) $ whenever $Q$ is the
unit cube $Q^\circ$, then we can easily find a linear extension operator $T$
as in Theorem \ref{introThm3}. Moreover, careful inspection of our solution
to $LIP\left( Q^\circ,E,f,x_{0},P_{0}\right) $ in \cite{FIL1} yields also Theorems \ref{introThm4} and \ref{introThm5}. Thus, the
heart of the matter is to solve $LIP\left( Q^\circ,E,f,x_{0},P_{0}\right) $. 

To do so, we first associate to any point $x\in \mathbb{R}^{n}$ the crucial
convex set 
\[
\sigma \left( x,E\right) =\left\{ J_{x}\left( F\right) :F\in \X,\left\Vert
F\right\Vert _{\X}\leq 1,F=0\text{ on }E\right\}.
\]%
This set measures the ambiguity in $J_{x}\left( F\right) $ when we seek
functions $F\in \X$ satisfying $F=f$ on $E$, with control on $\left\Vert
F\right\Vert _{\X}$.

Using the geometry of the convex sets $\sigma \left( x,E\right)$, we will
attach \textquotedblleft labels" $\mathcal{A}$ to cubes $Q\subset \mathbb{R}%
^{n}$. A label is simply a set of multiindices of order $\leq m-1.$ Very
roughly speaking, we say that $Q$ is \textquotedblleft tagged" with $%
\mathcal{A}$ if either 

\begin{itemize}
\item $Q$ consists of at most one point in $E$, or 

\item The scaled monomial $y\mapsto \delta _{Q}^{\text{power}}\cdot \left(
y-x\right) ^{\alpha }$ belongs to $\sigma \left( x,E\right) $ for all $%
\alpha \in \mathcal{A}$ and $x\in E\cap Q$. 
\end{itemize}

If $Q$ is tagged with $\mathcal{A}$, then we are relatively free to modify $%
\partial ^{\alpha }F\left( x\right) $ for $\alpha \in \mathcal{A}$, $x\in E$
when we seek a solution $F$ to our local interpolation problem $LIP\left(
Q,E,f,x_{0},P_{0}\right)$. 

The notion of tagging gives rise to a Calder\'on-Zygmund decomposition $CZ\left( \mathcal{A}\right) $ of the unit cube $Q^\circ$ for each label $%
\mathcal{A}$. The cubes of $CZ\left( \mathcal{A}\right) $ are the maximal
dyadic subcubes of $Q^\circ$ that are tagged with $\mathcal{A}$.

There is a natural order relation $<$ on labels. If labels $\mathcal{A}$, $%
\mathcal{B}$ satisfy $\mathcal{A}<\mathcal{B}$, then the decomposition $%
CZ\left( \mathcal{A}\right) $ of $Q^\circ$ refines the decomposition $CZ\left( 
\mathcal{B}\right) $. The maximal label under $<$ is the empty set $\emptyset$,
and the Calder\'on-Zygmund decomposition $CZ\left( \emptyset \right) $
consists of a single cube $Q^\circ$. The minimal label under $<$ is the
set $\mathcal{M}$ of all multiindices of order $\leq m-1$. The
decomposition $CZ\left( \mathcal{M}\right) $ is so fine that each $Q\in
CZ\left( \mathcal{M}\right) $ contains at most one point of our finite set $%
E.$

We now use the decomposition $CZ\left( \mathcal{A}\right) $ to solve Local Interpolation Problems. By induction on the label $\mathcal{A}$ (with
respect to the order $<$), we solve the problem $LIP\left(
Q,E,f,x_{0},P_{0}\right) $ whenever $Q\in CZ\left( \mathcal{A}\right)$. 

\underline{In the base case}, $\mathcal{A}=\mathcal{M}$, the minimal label. 

Since any $Q\in CZ\left( \mathcal{M}\right) $ contains at most one point of $%
E$, our local interpolation problem $LIP\left( Q,E,f,x_{0},P_{0}\right) $ is
trivial.

\underline{For the induction step}, fix a label $\mathcal{A\not=M}$. Let $%
\mathcal{A}^{-}$ be the label immediately preceding $\mathcal{A}$ in the
order $<$. We make the inductive assumption that we can solve $LIP\left(
Q^{\prime },E,f,x^{\prime },P^{\prime }\right) $ whenever $Q^{\prime }\in
CZ\left( \mathcal{A}^{-}\right) $. Using this assumption, we solve $%
LIP\left( Q,E,f,x_{0},P_{0}\right) $ when $Q\in CZ\left( \mathcal{A}\right) $%
. To do so, we recall that $CZ\left( \mathcal{A}^{-}\right) $ refines $%
CZ\left( \mathcal{A}\right) $, hence our cube $Q$ is partitioned into
finitely many cubes $Q_{\nu }\in CZ\left( \mathcal{A}^{-}\right) $. For each 
$Q_{\nu }$, we carefully pick a point $x_{\nu }\in Q_{\nu }$ and a jet $%
P_{\nu }\in \mathcal{P}$. Our inductive assumption allows us to solve the
local problem 
\[
LIP\left( Q_{\nu },E,f,x_{\nu },P_{\nu }\right) 
\]%
for each $\nu $. Using a partition of unity, we patch together the solutions 
$F_{\nu }$ to the above local problems, and hope that the resulting function 
$F$ solves our problem $LIP\left( Q,E,f,x_{0},P_{0}\right) $. It works
provided we do a good job of picking the jets $P_{\nu }$. We refer the
reader to the introduction of \cite{FIL1} for some of
the ideas involved in picking the $P_{\nu }$. (See especially the discussion
in \cite{FIL1} of \textquotedblleft keystone cubes").

Thus, we can complete our induction on $\mathcal{A}$, and solve $LIP\left(
Q,E,f,x_{0},P_{0}\right) $ whenever $Q\in CZ\left( \mathcal{A}\right) $.

In particular, since $CZ\left( \mathcal{\emptyset }\right) $ consists of the
single cube $Q^\circ$, we have succeeded in solving any local interpolation
problem $LIP\left( Q,E,f,x_{0},P_{0}\right) $ with $Q=Q^\circ$. As explained
above, this allows us to deduce Theorems \ref{introThm3}, \ref{introThm4}, %
\ref{introThm5}. That's the good news.

The bad news is that we cannot tell whether a given cube $Q$ is tagged with
a given label $\mathcal{A}$, since that requires perfect knowledge of the
convex sets $\sigma \left( x,E\right) \subset \mathcal{P}$. Therefore, our
Calder\'on-Zygmund decomposition $CZ\left( \mathcal{A}\right) $ and our proofs
of Theorems \ref{introThm3}, \ref{introThm4}, \ref{introThm5} in \cite{FIL1} are non-constructive.

To overcome the obstacle, we introduce here a variant of our local
interpolation problem, a variant of the convex set $\sigma \left( x,E\right)
$, and a modified definition of tagging of a cube $Q$ with a label $\mathcal{%
A}$. We still cannot tell whether a given $Q$ is tagged with a given $%
\mathcal{A}$. However, using ideas from \cite{FIL1}, we
show how to \underline{test} $Q$ for tagging with $\mathcal{A}$. If $Q$ passes the test,
then it is tagged with $\mathcal{A}$. If $Q$ fails the test, then we do not
know whether $Q$ is tagged with $\mathcal{A}$, but we know that a somewhat
larger cube $Q^{\prime }\supset Q$ cannot be tagged with $\mathcal{A}$.

We show how to implement the above test by efficient algorithms. Moreover,
if we are given dyadic cubes $Q_{1}\subset Q_{2}\subset \cdots \subset Q_{\nu }$
such that $Q_{1}\cap E = \cdots = Q_{\nu }\cap E$, then we can test all the $Q_{i}$ simultaneously.
This idea is useful if our set $E$ involves vastly different lengthscales.
It provides a crucial speedup that allows us to bound the work by $N\log N$
as promised in Theorem \ref{introMainTheorem}.

Using the above tests, we produce a decomposition $CZ\left( \mathcal{A}%
\right) $ analogous to the decomposition defined in \cite{FIL1}. This allows a
constructive proof of Theorems \ref{introThm3}, \ref{introThm4}, \ref%
{introThm5}. To implement that proof by efficient algorithms and thus
establish Theorem \ref{introMainTheorem} requires additional ideas not
discussed in this introduction.

This concludes our sketch of the proof of Theorem \ref{introMainTheorem}. We
again warn the reader that it is highly oversimplified. The sections that
follow will give the correct version. In the next section, we start from
scratch.

We mention several open problems related to our work.

\begin{itemize}
\item In place of our standing assumption $p>n$, we may assume merely that $%
p>n/m$. The Sobolev theorem would then tell us that $L^{m,p}\left( \mathbb{R}%
^{n}\right) \subset C_{loc}^{0}\left( \mathbb{R}^{n}\right) $. Consequently,
any $F\in L^{m,p}\left( \mathbb{R}^{n}\right) $ may be restricted to a
finite set $E$, and our Problems\ (A) and (B) still make sense. It would be
very interesting to understand the problems of interpolation and extension
for $L^{m,p}\left( \mathbb{R}^{n}\right) $ and $W^{m,p}\left( \mathbb{R}%
^{n}\right) $ when $n/m<p\leq n$.

\item Is it possible to dispense with the assists $\Omega =\left\{ \omega
_{1},\cdots ,\omega _{\nu _{\max }}\right\} $ in Theorem \ref{introThm4},
and write each $\xi _{l}$ in the form 
\[
\xi _{l}\left( f\right) =\sum_{i=1}^{I}\beta _{l {i}}f\left(
z_{l {i}}\right) 
\]%
with $I$, $\beta _{l {i}}$ and $z_{l {i}}$ $\in E$ independent of $f$, and
with $\left\vert I\right\vert \leq C?$ Shvartsman \cite{S5} has proven this for $\X=L^{2,p}\left( \mathbb{R}^{2}\right) $ ($p>2$%
). Perhaps, it is true for general $L^{m,p}\left( \mathbb{R}^{n}\right) $.
The analogous result for interpolation by functions in $C^{m}\left( \mathbb{R%
}^{n}\right) $ is contained in Fefferman-Klartag \cite{FK2}. For the extension operators $T$ in Theorem \ref{introThm3}, one
cannot get away without assists; see \cite{FIL2}. 

These issues are connected with \textquotedblleft sparsification"; see 
\cite{BSS}.

\item We have constructed essentially optimal functions $F\in \X$ that agree
perfectly with a given function $f$ on $E.$ It would be natural to require
instead that $F$ agree with $f$ up to a given tolerance. More precisely,
given $f\in \X\left( E\right) $ and a positive function $\mu :E\rightarrow
(0,\infty ]$, we should compute a function $F\in \X$ that minimizes 
\begin{equation}
\left\Vert F\right\Vert _{\X}^{p}+\sum_{x\in E}\mu \left( x\right) \left\vert
F\left( x\right) -f\left( x\right) \right\vert ^{p},  \label{intro11}
\end{equation}%
up to a universal constant factor $C.$ (When $\mu \left( x\right) =+\infty $%
, we demand that $F\left( x\right) =f\left( x\right) $ and delete the
corresponding term from the above sum.) Compare with Fefferman-Klartag 
\cite{FK2}. 

It would be interesting to study the problem of optimizing (\ref{intro11}) for general $L^{m,p}\left( \mathbb{R}^{n}\right) $ and $W^{m,p}\left( \mathbb{R}^{n}\right) $ ($p>n$). The work of P. Shvartsman \cite{S6} on the Banach
space $L^{1,p}\left( \mathbb{R}^{n}\right) +L^{p}\left( \mathbb{R}^{n},d\mu
\right) $ is surely relevant here. 

\item It would be very interesting to produce practical algorithms that
(unlike our present algorithms) compute $C$-optimal interpolants for a
not-so-big universal constant $C.$ 
\end{itemize}

\bigskip 

This paper is a part of a literature on \textquotedblleft Whitney's
extension problem", going back over $3/4$ century and including
contributions by many authors. See e.g., H. Whitney \cite{W1,W2, W3}, G. Glaeser \cite{G}, Y. Brudnyi and P. Shvartsman \cite{BS3, BS1, BS4, BS2}, P. Shvartsman \cite{S3,S4,S7,S1}, J. Wells \cite{We1}, E. Le Gruyer \cite{LG1,LG2}, M. Hirn and E. Le Gruyer \cite{HL1}, C. Fefferman and B. Klartag \cite{FK1,FK2}, N. Zobin \cite{Z1,Z2}, as well as our own works \cite{F1,F5,F3,F4,F2,F6, I1, L1}.

We are grateful to Bernard Chazelle, Bo'az Klartag, Assaf Naor, Pavel Shvartsman, and Nahum Zobin for many enlightening conversations. We are grateful also to the American Institute of Mathematics, the College of William and Mary, the Fields Institute, and the Banff International Research Station for hosting workshops on Whitney's problems. The support of the AFOSR and the NSF is gratefully acknowledged.

Let us now begin the work of interpolating data.

\chapter{Preliminaries}

\section{Notation}
\label{sec_not}

Fix integers $m, n \geq 1$ and a real number $p > n $. We work in $\R^n$ with the  $\ell^\infty$ metric. Thus, given $x = (x_1,\cdots,x_n) \in \R^n $ we denote
\begin{equation*} \lvert x \rvert  := \max_{1 \leq i \leq n} \lvert x_i \rvert.
\end{equation*}
Given nonempty subsets $S,S' \subset \R^n$, we denote
\begin{align*}
& \dist (S,S') := \inf \{ \lvert x - y \rvert : x \in S, y \in S' \},\\
& \diam(S) := \sup \{ \lvert x - y \rvert : x,y \in S\}.
\end{align*}

A \emph{cube} takes the form
\[ Q = \bigl[ x_1 - \delta/2, x_1 + \delta/2 \bigr) \times \cdots \times \bigl[x_n - \delta/2,x_n + \delta/2 \bigr).\]
Let $x_Q := (x_1,\cdots,x_n)$ and $\delta_Q := \delta$ denote the center and sidelength of the cube $Q$, respectively. Let $A\cdot Q$ ($A > 0$) denote the $A$-dilate of $Q$ about its center. Hence, the cube $A \cdot Q$ has center $x_Q$ and sidelength $A \delta_Q$.

A \emph{dyadic cube} takes the form
\[Q = \bigl[ j_1 \cdot 2^{k}, (j_1 + 1) \cdot 2^{k} \bigr) \times \cdots \times \bigl[ j_n \cdot 2^{k}, (j_n + 1) \cdot 2^{k} \bigr)\]
for $j_1,\cdots,j_n, k \in \Z$.

We say that two dyadic cubes $Q$ and $Q'$ \underline{touch} either if $Q=Q'$, or if $Q$ is disjoint from $Q'$ but the boundaries $\partial Q$ and $\partial Q'$ have a nonempty intersection. We write $Q \leftrightarrow Q'$ to indicate that $Q$ touches $Q'$. 

We may bisect a dyadic cube $Q$ into $2^n$ dyadic subcubes of sidelength $\frac{1}{2} \delta_Q$ in the natural way. We call these subcubes the \underline{children} of $Q$. We write $Q^+$ to denote the \underline{parent} of $Q$, i.e., the unique dyadic cube for which $Q$ is a child of $Q^+$.

We let $\cP$ denote the vector space of real-valued $(m-1)$-st degree polynomials on $\R^n$, and we set $D := \dim \cP$. We identify $\cP$ with $\R^D$, by identifying $P \in \cP$ with $(\partial^\alpha P(0))_{|\alpha| \leq m -1}$.

Given $F \in C^{m-1}_{\loc}(\R^n)$ and a point $\underline{x} \in \R^n$, let $J_{\underline{x}} F \in \cP$ (the ``jet of $F$ at $\underline{x}$'') denote the $(m-1)$-st order Taylor polynomial 
\[ (J_{\underline{x}} F) (x) = \sum_{|\alpha| \leq m-1} \frac{1}{\alpha!} \partial^\alpha F(\underline{x}) \cdot (x-\underline{x})^\alpha.\]

Given $P, R \in \cP$, we define the product $ P \odot_x R = J_x(P \cdot R) \in \cP$. 

\noindent\textbf{Sobolev spaces.}

We work with the Sobolev space $\X = L^{m,p}(\R^n)$ with seminorm
\[
\| F \|_{\X} = \left( \int_{\R^n} \sum_{|\alpha| = m} \lvert \partial^\alpha F(x) \rvert^p dx \right)^{1/p}.
\]
We assume throughout that $p>n$. Given a connected domain $\Omega \subset \R^n$ with piecewise smooth boundary, let $\X(\Omega) = L^{m,p}(\Omega)$ be the Sobolev space consisting of functions $F : \Omega \rightarrow \R$ whose distributional derivatives $\partial^\alpha F$ (for $|\alpha| = m$) belong to $L^p(\Omega)$. On this space we define the seminorm
\[
\| F \|_{\X(\Omega)} = \left( \int_{\Omega} \sum_{|\alpha| = m} \lvert \partial^\alpha F(x) \rvert^p dx \right)^{1/p}.
\]
We may restrict attention to domains that are given as the union of two intersecting rectangular boxes. 
(A rectangular box is a  Cartesian product of left-closed, right-open intervals.)

\noindent\textbf{Lists}

A \emph{list} $\Xi$ is a collection of elements that can contain duplicates. Hence, for a list
\begin{equation}
\label{alist}
\Xi = \{ \xi_1,\cdots, \xi_L\}
\end{equation}
we may have $\xi_\ell = \xi_{\ell'}$ for distinct $\ell, \ell'$. We define $\# \left[ \Xi \right] =L$ for the list \eqref{alist}.  

Given a sequence $(a_\xi)$ of real numbers, indexed by elements $\xi$ in $\Xi$, we define
\[ \sum_{\xi \in \Xi} a_\xi = \sum_{\ell = 1}^L a_{\xi_{\ell}}.
\]

Given a sequence  $\Xi_1,\cdots,\Xi_M$ of lists, where $\Xi_m = \left\{ \xi^{[m]}_1,\cdots,\xi^{[m]}_{L_m} \right\}$ for each $1 \leq m \leq M$, we define the list $\Xi_1 \cup \cdots \cup \Xi_M$ by taking all the elements in the respective sublists together, namely
\[\Xi_1 \cup \cdots \cup \Xi_M := \left\{\xi_1^{[1]},\cdots,\xi_{L_1}^{[1]},\cdots, \xi_1^{[M]},\cdots, \xi^{[M]}_{L_M} \right\}.\]
We do not remove duplicate elements when forming the union of lists.

\noindent\textbf{Convention on constants.}

A \emph{universal constant} is a positive number determined by $m,n,$ and $p$. We use letters $C, c, C',$ etc, to denote universal constants. Let $t \in \R$. We use the symbol $C(t)$ to denote a positive number that depends only on $m,n, p, $ and $t$. A single letter or symbol may be used to denote different constants in separate occurrences.

We write $ A \lesssim B$ or $A = \mathcal{O}(B)$ to indicate the estimate $A \leq C B$, and we write $A \sim B$ to indicate the estimate $C^{-1} B \leq A \leq C B$. Here, the constant $C$ depends only on $m,n,$ and $p$.

Similarly, we write $ A \lesssim_t B$ or $A = \mathcal{O}_t(B)$ to indicate the estimate $A \leq C(t) \cdot B$, and we write $A \sim_t B$ to indicate the estimate $C(t)^{-1} \cdot B \leq A \leq C(t) \cdot B$. Here,  the constant $C(t)$ depends only on $m,n,p,$ and $t$.

\section{The Infinite-Precision Model of Computation}
\label{sec_moc1}

For infinite-precision, our model of computation consists of an idealized von Neumann computer \cite{V} able to work with exact real numbers. We assume that a single memory cell is capable of storing an arbitrary real number with perfect precision.

We assume that each of the following operations can be carried out using one unit of "work".

\begin{itemize}
\item We read the real number stored at a given address, or entered from an input device.

\item We write a real number from a register to a given memory address or to an output device.

\item Given real numbers $x$ and $y$, we return the numbers $x+y$, $x-y$, $xy$, $x/y$ (unless $y=0$), $\exp(x)$, and $\ln(y)$ (if $y>0$), and we decide whether $x<y$, $x>y$ or $x=y$.

\item Given a real number $x$, we return the greatest integer less than or equal to $x$.

\item Given dyadic intervals $I=[x,y)$ and $J=[a,b)$, both contained in $[0, \infty)$, we return the smallest dyadic interval containing both $I$ and $J$.
\end{itemize}

The above model of computation is subject to serious criticism, even without our assumption on the rapid processing of dyadic intervals. (See \cite{FK2,HS,Sch}.) Therefore, in a later section, we will give a model of computation in finite-precision. We believe that our finite-precision model faithfully reflects a subset of the capabilities of an actual computer (e.g. we don't assume any possibility of parallel processing).

Presumably, few readers will want to wade through the issues arising from implementing our algorithms in finite-precision. Hence, we first present our results assuming the above infinite-precision model of computation. We then explain how to modify our algorithms in order to succeed in the finite-precision model. These modifications are described in the Appendix.

\section{Basic Estimates on Sobolev Functions}

Let $\cP$ denote the vector space of $(m-1)^{\text{st}}$ degree polynomials.

Given a point $x \in \R^n$ and a number $\delta > 0$, we define
$$\lvert P \rvert_{x,\delta} := \max_{\lvert \beta \rvert \leq m-1} \lvert \partial^\beta P(x) \rvert \cdot \delta^{|\beta| + n/p - m} \qquad \mbox{for} \; P \in \cP.$$
Thus, $\lvert P \rvert_{x,\delta}$ is a norm on $\cP$. The unit ball of this norm is  $\cB(x,\delta) := \{ P \in \cP :  | P|_{x,\delta} \leq 1\}.$

We next present a few basic properties of the objects defined above.
\begin{lem}\label{pnorm} Let $Q \subset \R^n$ be a cube and let $K \geq 1$. For any polynomial $P \in \cP$, the following estimates hold:
\begin{align*}
&\delta_Q^{-m}\|P \|_{L^p(Q)} \sim \sum_{k=0}^{m-1} \delta_Q^{k-m} \| \nabla^k P \|_{L^p(Q)} \sim_K \lvert P \rvert_{x,\delta} \;\;\; \mbox{for} \; x \in K Q \; \mbox{and} \; \delta \in [ K^{-1} \delta_Q , K \delta_Q].\\ 
& \lvert P \rvert_{x',\delta'} \leq C(K) \cdot \lvert P \rvert_{x,\delta} \qquad\; \mbox{for} \; \lvert x - x' \rvert \leq K \delta' \; \mbox{and} \;\; \delta \leq K \delta'.\\
&\cB(x,\delta) \subseteq C(K) \cdot \cB(x',\delta') \quad  \mbox{for} \;\; \lvert x - x' \rvert \leq K \delta' \; \mbox{and} \; \delta \leq K \delta'.
\end{align*}
Here, $C(K)$ depends only on $m,n,p,$ and $K$.
\end{lem}

We present a few useful estimates on functions in $\X(Q) = L^{m,p}(Q)$, where $Q \subset \R^n$ is a cube. These estimates originate from the classical Sobolev inequality (Proposition \ref{si1}) and an interpolation inequality (Proposition \ref{int_ineq}).

\begin{prop}[Sobolev inequality]\label{si1}
Let $Q \subset \R^n$ be a cube, and let $F \in \X(Q)$. For any $x,y \in Q$, and any multiindex $\beta$ with $\lvert \beta \rvert \leq m -1$, we have
\begin{align*} & \left\lvert \partial^\beta ( J_x F - F)(y) \right\rvert \leq C \cdot \delta_Q^{m - |\beta| - n/p} \| F \|_{\X(Q)}.
\end{align*}
\end{prop}
\begin{proof}
The estimate is an easy consequence of the Sobolev embedding theorem and Taylor's theorem. 

We first review the Sobolev embedding theorem. Let $F \in \X(Q) = L^{m,p}(Q)$. Due to our standing assumption that $p>n$, we can use the Sobolev embedding theorem (see \cite{GT}), which implies that $F$ belongs to the H\"{o}lder space $C^{m-1,1- \frac{n}{p}}(Q)$\footnote{Here, $C^{m-1,\gamma}(Q)$ ($\gamma \in (0,1]$) is the H\"{o}lder space consisting of all functions $F : Q \rightarrow \R$ that satisfy the estimate $\lv \partial^\alpha F(x) - \partial^\alpha F(y) \rv \leq A \cdot \lv x - y \rv^{\gamma} $ for some $A < \infty$ and for all multiindices $\alpha$ with $\lv \alpha \rv = m-1$ and all $x,y \in Q$. The H\"{o}lder seminorm of $F$ in  $C^{m-1,\gamma}(Q)$ is defined to be the infimum of all such constants $A$.}, and moreover we have an estimate on the H\"{o}lder seminorm:
\begin{equation}
\label{sobembthm}
\| F \|_{C^{m-1,1-\frac{n}{p}}(Q)} \leq K \cdot \| F \|_{L^{m,p}(Q)}.
\end{equation}
Here, $K = K(m,n,p)$. In particular, $K$ is independent of  $Q$. We refer the reader to \cite{GT} for a proof of the estimate \eqref{sobembthm} when $Q = [0,1)^n$. We can prove \eqref{sobembthm} for the same choice of $K$ and a general $Q \subset \R^n$ using a standard rescaling argument. We provide details below.

Let $Q \subset \R^n$, and let $F \in L^{m,p}(Q)$. Let $\tau: \R^n \rightarrow \R^n$ be a transformation of the form $\tau(x) = R \cdot x + x_0$ ($R > 0$, $x_0 \in \R^n$) satisfying that $\tau$ maps $Q^\circ := [0,1)^n$ onto $Q$. We define a transformed function $\widetilde{F} = F \circ \tau : Q^\circ \mapsto \R$. The Sobolev and H\"{o}lder norms relevant to our discussion are transformed in a simple fashion. Indeed, by a change of variables we have
\begin{align*}
\| \widetilde{F} \|_{L^{m,p}(Q^\circ)} &= \left( \int_{Q^\circ} \lv \nabla^m( F(R \cdot x + x_0)) \rv^p dx \right)^{1/p} \\
& = R^{m-n/p} \left( \int_{Q} \lv \nabla^m F(x) \rv^p dx \right)^{1/p} \\
& = R^{m - n/p} \| F \|_{L^{m,p}(Q)}.
\end{align*}
Similarly, we have $\| \widetilde{F} \|_{C^{m-1,1-\frac{n}{p}}(Q^\circ)} = R^{m-n/p} \| F \|_{C^{m-1,1-\frac{n}{p}}(Q)}$ . 

We apply the known version of the estimate \eqref{sobembthm} to the function $\widetilde{F}$. Thus, we obtain $\| \widetilde{F} \|_{C^{m-1,1-\frac{n}{p}}(Q^\circ)} \leq K \cdot \| \widetilde{F} \|_{L^{m,p}(Q^\circ)}$. From the above equations we thus deduce that
\[
\| F \|_{C^{m-1,1-n/p}(Q)} \leq K \cdot \| F \|_{L^{m,p}(Q)}.
\]
This completes the proof of \eqref{sobembthm}.


Recall that $J_x F$ is the $(m-1)^{\text{st}}$ degree Taylor polyomial of $F$ at $x$. Hence, by definition,
\begin{equation}
\label{taylorpoly}
J_{x} F(y) = \displaystyle \sum_{|\alpha| \leq m-1} \frac{1}{\alpha!} \partial^\alpha F(x) \cdot (y-x)^\alpha \qquad (y \in \R^n).
\end{equation}
Taylor's theorem states that 
\[
\lv \partial^\beta (J_x F  - F)(y) \rv \leq C \cdot \lv x - y \rv^{m - 1 + \gamma - \lv \beta \rv} \cdot \| F \|_{C^{m-1,\gamma}(Q)}
\]
for any $F \in C^{m-1,\gamma}(Q)$, any $x,y \in Q$, and any multiindex $\beta$ with $\lv \beta \rv \leq m-1$. Here, $C = C(m,n,\gamma)$ depends only on $m$, $n$, and $\gamma$. We apply this estimate with $\gamma = 1-n/p$ and use \eqref{sobembthm} to see that
\[
\lv \partial^\beta (J_x F  - F)(y) \rv \leq C K \cdot \lv x - y \rv^{m-n/p - \lv \beta \rv} \cdot \| F \|_{L^{m,p}(Q)}
\]
for any $F \in L^{m,p}(Q)$. We bound $\lv x -y \rv \leq \delta_Q$ to prove the desired estimate and complete the proof of Proposition \ref{si1}.

\end{proof}

We can formulate the Sobolev inequality as an estimate involving the norms $\lv \cdot \rv_{x,\delta}$ introduced earlier. Indeed,
\[
\lvert J_x F - J_y F \rvert_{y,\delta_Q} \leq C \| F \|_{\X(Q)} \;\; \mbox{whenever} \; x,y \in Q.
\]

\begin{prop}\label{int_ineq}
Let $Q \subset \R^n$ be a cube, and let $F \in \X(Q)$. For any multiindex $\beta$ with $\lv \beta \rv \leq m$, we have
\[ 
\| \partial^\beta F \|_{L^p(Q)} \leq C \cdot \left[ \delta_Q^{-|\beta|} \| F\|_{L^p(Q)} + \delta_Q^{m - |\beta|} \|F\|_{\X(Q)} \right].
\]
\end{prop}
\begin{proof}
A standard scaling argument allows us to reduce to the case when $Q = [0,1)^n$. For a proof of this estimate, see \cite{GT}.
\end{proof}

\begin{lem}\label{si2}
Let $Q\subset \R^n$ be a cube, and let $F \in \X(Q)$. For any $x \in Q$, we have
\begin{align}
\label{s_ineq0} 
\delta_Q^{-m} \| F \|_{L^p(Q)} & \leq C \cdot \left[ \|F \|_{\X(Q)}  +   \sum_{|\beta| \leq m-1} \lvert \partial^\beta F(x) \rvert \delta_Q^{|\beta| + n/p - m} \right] \\
\notag{}
& \leq C' \cdot \left[ \| F \|_{\X(Q)} + \delta_Q^{-m} \| F\|_{L^p(Q)} \right].
\end{align}
For any cube $Q' \subset Q$, we have
\begin{equation}
 \delta_Q^{-m} \| F\|_{L^p(Q)} \leq C'' \cdot \left[ \delta_{Q'}^{-m} \| F \|_{L^p(Q')} + \| F\|_{\X(Q)} \right]. \label{s_ineq2}
\end{equation}
Here, $C,C',C''$ denote constants depending only on $m,n,$ and $p$.
\end{lem}
\begin{proof}

Let $Q \subset \R^n$, and let $x, y \in Q$. 

From the Sobolev inequality and the definition \eqref{taylorpoly} of the Taylor polyomial, we have
\begin{align*}
 \lvert F(y) \rvert &\leq \lvert F(y) - J_x F(y) \rvert + \lvert J_xF(y) \rvert \\
 & \lesssim \delta_Q^{m - \frac{n}{p}} \| F \|_{\X(Q)} + \sum_{|\beta| \leq m-1} \lvert \partial^\beta F(x) \rvert \cdot \lvert x - y \rvert^{|\beta|}.
 \end{align*}
 Hence,
 \begin{align*}
 \int_Q \lv F(y) \rv^p dy &\lesssim \delta_Q^{mp - n} \| F \|_{\X(Q)}^p \cdot \Vol(Q) + \sum_{\lv \beta \rv \leq m - 1} \lv \partial^\beta F(x) \rv^p \cdot \int_{Q} \lv x - y \rv^{\lv \beta \rv \cdot p } dy \\
 &\lesssim \delta_Q^{mp} \| F \|_{\X(Q)}^p + \sum_{\lv \beta \rv \leq m - 1} \lv \partial^\beta F(x) \rv^p \cdot \delta_Q^{\lv \beta \rv \cdot p + n}
 \end{align*}
We raise each side to the power $1/p$. This implies the first estimate in \eqref{s_ineq0}.

We now complete the proof of \eqref{s_ineq0}.
Let $\lvert \beta \rvert \leq m-1$, and let $x,y \in Q$. As before, from the Sobolev inequality and \eqref{taylorpoly} we have
\begin{align*}
\lv   \partial^\beta F(x) \rv &\leq \lv \partial^\beta F(x) - \partial^\beta (J_y F)(x) \rv + \lv  \partial^\beta (J_y F)(x) \rv \\
&\lesssim \delta_Q^{m - |\beta| - n/p} \| F\|_{\X(Q)} + \sum_{|\gamma| \leq m-1-|\beta|} \lv   \partial^{\beta + \gamma}F(y) \rv  \cdot \lvert x - y\rvert^{|\gamma|} \\
& \leq \delta_Q^{m - |\beta| - n/p} \| F\|_{\X(Q)} + \sum_{|\gamma| \leq m-1-|\beta|}  \lv    \partial^{\beta + \gamma}F(y)    \rv \cdot \delta_Q^{|\gamma|}.
\end{align*}
We raise this estimate to the power $p$ and average over $y \in Q$. Hence,
\begin{align*}
\lvert \partial^\beta F(x) \rvert \lesssim \delta_Q^{m - |\beta| - n/p} \| F \|_{\X(Q)} + \sum_{|\gamma| \leq m - 1 - |\beta| } \delta_Q^{|\gamma| - n/p} \| \partial^{\beta + \gamma} F \|_{L^p(Q)}.
\end{align*}
We apply Proposition \ref{int_ineq} to estimate the terms in the sum over $\gamma$. We see that $\|  \partial^{\beta + \gamma} F \|_{L^p(Q)}$ is bounded by $C \cdot \left[ \delta_Q^{m - \lv \beta \rv - \lv \gamma \rv} \cdot \| F \|_{\X(Q)} + \delta_Q^{- \lv \beta \rv - \lv \gamma \rv} \cdot \| F \|_{L^p(Q)} \right]$.  Thus, we conclude that 
\[
\lv \partial^\beta F(x) \rv  \lesssim \delta_Q^{m - |\beta| - n/p} \| F \|_{\X(Q)} + \delta_Q^{- | \beta| -n/p} \| F\|_{L^p(Q)}.
\]
We thus obtain the second estimate in \eqref{s_ineq0}.

We will finally prove the inequality \eqref{s_ineq2}. Let $Q' \subset Q$ be given cubes. Then \eqref{s_ineq0} implies that
\[ \delta_Q^{-m} \| F \|_{L^p(Q)} \lesssim \|F \|_{\X(Q)} +  \sum_{|\beta| \leq m-1} \lvert \partial^\beta F(x) \rvert \cdot \delta_{Q'}^{|\beta| + \frac{n}{p} - m} \quad \mbox{for any} \; x \in Q'.\]
(Here, we use that $\delta_{Q'} \leq \delta_Q$ and $|\beta| + \frac{n}{p} - m < 0$.) By averaging $p$-th powers in the above estimate, we see that
\[ \delta_Q^{-m} \| F \|_{L^p(Q)} \lesssim \|F \|_{\X(Q)} + \sum_{|\beta| \leq m-1} \delta_{Q'}^{|\beta| - m}  \| \partial^\beta F  \|_{L^p(Q')}.\]
Finally, we apply Proposition \ref{int_ineq} to estimate the terms in the sum over $\beta$.  Thus, we see that $\delta_{Q'}^{\lv \beta \rv - m} \| \partial^\beta F \|_{L^p(Q')}$ ($0 \leq \lv \beta \rv \leq m-1$) is bounded by $C \cdot \left[ \|F \|_{\X(Q)} + \delta_{Q'}^{-m} \| F \|_{L^p(Q')} \right]$ . This completes the proof of \eqref{s_ineq2}. This concludes the proof of Lemma \ref{si2}.
\end{proof}

\begin{lem}
\label{si4}
Let $Q', Q'' \subset \R^n$ be cubes with intersecting closures, and suppose that $\frac{1}{2} \delta_{Q''} \leq \delta_{Q'} \leq 2\delta_{Q''}$. Then for any $R',R'' \in \cP$ and any $H \in \X$, we have
\begin{align*} \delta_{Q'}^{-m} \|R' - R''\|_{L^p(Q')}  \lesssim \delta_{Q'}^{-m}  \| H - R' \|_{L^p(\frac{65}{64}Q')} + \delta_{Q''}^{-m} \| H - R'' \|_{L^p(\frac{65}{64}Q'')} &+ \| H\|_{\X(\frac{65}{64}Q')} \\
&+ \| H \|_{\X(\frac{65}{64}Q'')}.
\end{align*}
\end{lem}
\begin{proof}
First we write
\[ \delta_{Q'}^{-m} \|R' - R''\|_{L^p(Q')} \leq \delta_{Q'}^{-m} \|H- R'\|_{L^p(\frac{65}{64}Q')} + \delta_{Q'}^{-m} \|H - R'' \|_{L^p(\frac{65}{64}Q')}.\]
For any fixed $x \in \frac{65}{64}Q' \cap \frac{65}{64}Q''$, Lemma \ref{si2} gives that
\[\delta_{Q'}^{-m} \|H - R'' \|_{L^p(\frac{65}{64}Q')} \lesssim \| H \|_{\X(\frac{65}{64}Q')} +  \sum_{|\beta| \leq m-1} \lvert \partial^\beta (H-R'')(x) \rvert \delta_{Q'}^{|\beta| + \frac{n}{p} - m}\]
and also
\[\sum_{|\beta| \leq m-1} \lvert \partial^\beta (H-R'')(x) \rvert \delta_{Q'}^{|\beta| + \frac{n}{p} - m} \lesssim \| H \|_{\X(\frac{65}{64}Q'')} + \delta_{Q''}^{-m} \| H - R'' \|_{L^p(\frac{65}{64}Q'')}.\]
This implies the conclusion of the lemma.
\end{proof}

A \emph{rectangular box} $B \subset \R^n$ is a Cartesian product of coordinate intervals that are left-closed and right-open, where each interval has a nonempty interior. The length of each interval is a \emph{sidelength} of $B$. The \emph{aspect ratio} of $B$ is the ratio of the longest to shortest sidelength of $B$.

Let $K \geq 1$. Suppose that $B$ is a rectangular box with aspect ratio at most $K$. We can map a cube onto $B$ by applying a transformation of the form $\tau : (x_1,\cdots,x_n) \mapsto (\delta_1 x_1, \cdots, \delta_n x_n)$, with $1 \leq \lvert \delta_j \rvert \leq K$ for $j=1,\cdots,n$. Let $F \in \X(B)$. By precomposing $F$ with the transformation $\tau$, and using Proposition \ref{si1} (the Sobolev inequality), we see that
\begin{equation}
\label{si1a}
\lvert \partial^\beta(F - J_yF)(x) \rvert \leq C(K) \cdot \| F\|_{\X(B)} \lvert x - y \rvert^{m-n/p-|\beta|} \qquad \mbox{for all} \;  x,y \in B, \; |\beta| \leq m-1.
\end{equation}
We raise this estimate to the power $p$, and integrate over $x \in B$ to obtain
\begin{equation}
\label{si1b}
\| \partial^\beta (F - J_yF) \|_{L^p(B)} \leq C(K) \cdot \| F \|_{\X(B)} \diam(B)^{m-|\beta|}.
\end{equation}

\begin{lem}\label{si3}
Let $B_1,B_2$ be rectangular boxes with aspect ratio at most $K$ and with $B_1 \cap B_2 \neq \emptyset$. Then for any $F \in \X(B_1 \cup B_2)$, any $x,y \in B_1 \cup B_2$, and any $\beta$ with $\lvert \beta \rvert \leq m-1$, we have
\begin{equation}\label{stuff} \lvert \partial^\beta ( J_x F - F)(y) \rvert \leq C(K) \cdot \lvert x - y \rvert^{m - |\beta| - \frac{n}{p}} \cdot \left\{ \| F \|_{\X(B_1)} + \| F\|_{\X(B_2)} \right\}.
\end{equation}
\end{lem}
\begin{proof}
If either $x,y \in B_1$ or $x,y \in B_2$ then \eqref{stuff} follows from the estimate \eqref{si1a}. 

Otherwise, we may assume that $x \in B_1$ and $y \in B_2$. Pick $z \in B_1 \cap B_2$ with the property that $\lvert x - z \rvert \leq \lvert x - y\rvert$ and $\lvert y - z \lvert \leq \lvert x - y \rvert$, and note that
\begin{align*}
\lvert J_x F - J_y F \rvert_{y,|x-y|} &\leq \lvert J_x F - J_z F \rvert_{y,|x-y|} + \lvert J_z F - J_y F \rvert_{y,|x-y|} \\
& \lesssim \lvert J_x F - J_z F \rvert_{z,|x-z|} + \lvert J_z F - J_y F \rvert_{z,|y-z|} \qquad (\mbox{by Lemma \ref{pnorm}}).\\
\end{align*}
Now, \eqref{si1a} implies that
\begin{align*}
\lvert J_x F - J_y F \rvert_{y,|x-y|} 
& \leq C(K) \cdot \left\{ \| F \|_{\X(B_1)} + \| F\|_{\X(B_2)} \right\}.
\end{align*}
This completes the proof of the lemma.
\end{proof}

Our last result is a special case of the Jones extension theorem \cite{J1}.

\begin{prop}\label{jones}
Let $Q$ be a cube in $\R^n$. Then there exists a bounded linear map $T  : \X(Q) \rightarrow \X$, such that $T(F) = F$ on $Q$, and $\| T(F) \|_{\X} \leq C \| F\|_{\X(Q)}$ for each $F \in \X(Q)$. Here, the constant $C$ depends only on the parameters of the function space $\X$, i.e., on the numbers $m,n,p$.
\end{prop}

\section{Trace Norms}

Given a finite subset $E \subset \R^n$, let $\X(E)$ denote the space of all functions $f : E \rightarrow \R$, with the trace seminorm
$$\|f\|_{\X(E)} := \inf \{ \|F \|_{\X} : F \in \X , \;\; F = f \; \mbox{on} \; E\}.$$
Given a cube $Q \subseteq \R^n$, and given $(f,P) \in \X(E) \oplus \cP$, let
\begin{equation}\label{norm}
\|(f,P)\|_Q := \inf \Biggl\{ M \geq 0 \; : \;\; \exists \; F \in \X \;\;\; \mbox{s.t.} \;\;\;
\begin{aligned}
& \quad\quad F = f \; \mbox{on} \; E \cap Q  \\
& \|F\|_{\X(Q)} +  \delta_Q^{-m} \| F - P \|_{L^p(Q)} \leq M
\end{aligned} \Biggr\}.
\end{equation}
Note that $\|(f,P)\|_Q$ is a seminorm on the space $\X(E) \oplus \cP$. 

Let
\begin{equation}\label{sigma}
\sigma(Q) := \Biggl\{ P \in \cP \; : \;\; \exists \; \varphi \in \X \;\;\; \mbox{s.t.} \;\;\;
\begin{aligned}
& \quad \quad \varphi = 0 \; \mbox{on} \; E \cap Q  \\
&\|\varphi\|_{\X(Q)} +  \delta_Q^{-m} \| \varphi - P \|_{L^p(Q)} \leq 1
\end{aligned} \Biggr\}.
\end{equation}
Note that $\sigma(Q) \subset \cP$ is convex and symmetric ($P \in \sigma(Q) \implies - P \in \sigma(Q)$).

As an easy consequence of our definitions, we have the following result.

\begin{lem}\label{lem_normmon}
Given cubes $Q_1 \subset Q_2$ such that $\delta_{Q_2} \leq A \delta_{Q_1}$, we have $\| (f,P) \|_{Q_1} \leq C(A) \cdot \|(f,P) \|_{Q_2}$ and $\sigma(Q_2) \subset C(A) \cdot \sigma(Q_1)$. In fact, one can take $C(A) = A^m$.
\end{lem}

Our next result relates the convex sets $\sigma(Q_1)$ and $\sigma(Q_2)$, where $Q_1 \subset Q_2$ are cubes that may have vastly different sizes.

\begin{lem}\label{pre_lem0}
Given cubes $Q_1 \subset Q_2$, we have
\begin{equation} \label{t0}
\sigma(Q_2) \subset C \cdot \left[ \sigma(Q_1) + \cB(x_{Q_1}, \delta_{Q_2})\right].
\end{equation}
If additionally $Q_1 \cap E  = Q_2 \cap E$, then also
\begin{equation}
\label{t1}
c \cdot \left[ \sigma(Q_1) + \cB(x_{Q_1}, \delta_{Q_2}) \right] \subset \sigma(Q_2).
\end{equation}
\end{lem}

\begin{proof} Suppose that $R \in \sigma(Q_2)$. By definition, this means that there exists $\phi \in \X$ such that $\phi = 0$ on $Q_2 \cap E$, and $\| \phi \|_{\X(Q_2)} + \delta_{Q_2}^{-m} \| \phi - R \|_{L^p(Q_2)} \leq 1$. In particular, we have $\phi = 0$ on $Q_1 \cap E$. Also, the Sobolev inequality implies that
\[ \| \phi \|_{\X(Q_1)} + \delta_{Q_1}^{-m} \| \phi - J_{x_{Q_1}} \phi \|_{L^p(Q_1)} \leq C \| \phi \|_{\X(Q_1)} \leq C \| \phi \|_{\X(Q_2)} \leq C.\]
Hence, $ J_{x_{Q_1}} \phi \in C \sigma(Q_1)$.  (Recall that $x_{Q_1}$ is the center of the cube $Q_1$.)

Similarly, using the triangle inequality followed by the Sobolev inequality, we have
\begin{align*}
 \delta_{Q_2}^{-m} \| J_{x_{Q_1}} \phi - R \|_{L^p(Q_2)} & \leq  \delta_{Q_2}^{-m} \| \phi - R \|_{L^p(Q_2)}  +  \delta_{Q_2}^{-m} \| \phi - J_{x_{Q_1}} \phi \|_{L^p(Q_2)} \\
 & \leq  \delta_{Q_2}^{-m} \| \phi - R \|_{L^p(Q_2)} + C \| \phi \|_{\X(Q_2)} \leq 1 + C.
 \end{align*}
Thus, $J_{x_{Q_1}} \phi - R \in C \cdot \cB(x_{Q_1},\delta_{Q_2})$.

Hence, we have shown that $R = J_{x_{Q_1}} \phi + (R - J_{x_{Q_1}} \phi)  \in C \sigma(Q_1) + C \cB(x_{Q_1},\delta_{Q_2})$. Since $R \in \sigma(Q_2)$ was arbitrary, this proves the first inclusion \eqref{t0}.

We now assume that $Q_1 \cap E = Q_2 \cap E$ and prove the second inclusion \eqref{t1}. 

Let $R \in \sigma(Q_1) + \cB(x_{Q_1},\delta_{Q_2})$, i.e., suppose that $R = P + P^\#$ with $P \in \sigma(Q_1)$ and $P^\# \in \cB(x_{Q_1}, \delta_{Q_2} )$.

By definition, $P^\# \in \cB(x_{Q_1}, \delta_{Q_2} )$ implies that $\lvert P^\# \rvert_{x_{Q_1},\delta_{Q_2}} \leq 1$, hence
\begin{equation} \label{plum1}
\delta_{Q_2}^{-m} \| P^\# \|_{L^p(Q_2)} \leq C \qquad \mbox{(see Lemma \ref{pnorm})}.
\end{equation}

By definition, $P \in \sigma(Q_1)$ implies that there exists $\varphi \in \X$ such that $\varphi = 0$ on $Q_1 \cap E$ and $\| \varphi \|_{\X(Q_1)} + \delta_{Q_1}^{-m} \| \varphi - P \|_{L^p(Q_1)} \leq 1$. We now pick $\widetilde{\varphi} \in \X$ with $\widetilde{\varphi} = \varphi$ on $Q_1$ and $\| \widetilde{\varphi} \|_{\X} \leq C \| \varphi \|_{\X(Q_1)}$. (See Proposition \ref{jones}.) Note that $\widetilde{\varphi} = 0$ on $Q_1 \cap E = Q_2 \cap E$. Moreover,
\begin{align*}
 \delta_{Q_2}^{-m} \| \widetilde{\varphi} - R \|_{L^p(Q_2)} &\leq \delta_{Q_2}^{-m} \| \widetilde{\varphi} - P \|_{L^p(Q_2)} + \delta_{Q_2}^{-m} \| P^\# \|_{L^p(Q_2)} \\
&\leq \delta_{Q_2}^{-m} \| \widetilde{\varphi} - P \|_{L^p(Q_2)} + C   \;\qquad\qquad  \qquad \quad (\mbox{by} \; \eqref{plum1})\\
&\leq C\delta_{Q_1}^{-m} \| \widetilde{\varphi} - P \|_{L^p(Q_1)} + C \| \widetilde{\varphi} \|_{\X(Q_2)} + C  \quad (\mbox{by Lemma \ref{si2}})\\
& \leq C' \delta_{Q_1}^{-m} \| \varphi - P \|_{L^p(Q_1)} + C' \| \varphi \|_{\X(Q_1)} + C' \leq C''.
\end{align*}
Since we also have $\| \widetilde{\varphi} \|_{\X(Q_2)} \leq C$, it follows that $R \in C \sigma(Q_2)$. Since $R \in \sigma(Q_1) + \cB(x_{Q_1},\delta_{Q_2})$ was arbitrary, this proves \eqref{t1} and completes the proof of the lemma.
\end{proof}

\section{The Depth of Linear Maps}
\label{sec_depth}

Let $E = \{ z_1,\cdots,z_N\} \subset \R^n$. We fix this enumeration of $E$ for the rest of the paper.

A linear functional $\omega : X(E) \rightarrow \R$ may be written as
\begin{equation} \label{longform1}
\omega(f) = \sum_{j=1}^{N} \mu_j \cdot f(z_{j})  \qquad \mbox{for real coefficients} \; \mu_1,\cdots,\mu_{N}.
\end{equation}
That's the \underline{long form of $\omega$}. Let $\depth(\omega)$ (``the depth of $\omega$'') be the number of nonzero coefficients $\mu_j$ in \eqref{longform1}.

Suppose $\depth(\omega) = d$. Then let $1 \leq j_1 < j_2 < \cdots < j_d \leq N$ be the indices for which $\mu_j \neq 0$ above. Also, let $\widetilde{\mu}_k = \mu_{j_k}$ for $k=1,\cdots,d$. Then we can write $\omega$ in the form
\[\omega(f) = \sum_{k=1}^d \widetilde{\mu}_k \cdot f(z_{{j_k}}).\]
That's the \underline{short form of $\omega$}.

To store $\omega$ in its long form, we store $\mu_1,\cdots, \mu_{N}$.

To store $\omega$ in its short form, we store $d$, $\widetilde{\mu}_1,\cdots, \widetilde{\mu}_{d}$, and $j_1,\cdots,j_d$.

Let $\Omega = \{ \omega_1,\cdots,\omega_K\}$ be a list of linear functionals on $\X(E)$, each given in short form. Recall that a list may contain duplicates. Hence, we can have  $\omega_k = \omega_{k'}$ with $k \neq k'$. We store the list $\Omega$ by storing a list of pointers to the functionals in $\Omega$. If we have stored two lists of functionals $\Omega$ and $\Omega'$, then we can compute and store $\Omega \cup \Omega'$ using work at most $C \cdot \left[ \#(\Omega) + \#(\Omega') \right]$.

A functional $\xi : \X(E) \oplus \cP \rightarrow \R$ may be written in the form
\begin{equation}
\label{longform2}
\xi(f,P) = \lambda(P) + \sum_{j=1}^N \mu_j f(z_j)
\end{equation}
for coefficients $\mu_1,\cdots,\mu_N$ and a functional $\lambda : \cP \rightarrow \R$. That's the \underline{long form of $\xi$}.

Let $d \in \N$. A functional $\xi : \X(E) \oplus \cP \rightarrow \R$ has \underline{$\Omega$-assisted depth $d$} provided that
\begin{equation}
\label{shortform}
\xi(f,P) = \lambda(P) + \eta(f) + \sum_{\nu=1}^{\nu_{\max}} \gamma_\nu \omega_{k_\nu}(f)
\end{equation}
where $\eta : \X(E) \rightarrow \R$ is a linear functional, and $\depth(\eta) + \nu_{\max} \leq d$. That's a \underline{short form of $\xi$} in terms of the assists $\Omega$. Note that perhaps we can describe a given $\xi$ in many different ways in short form.

To store the long form of $\xi$, we store $\lambda$, $\mu_1,\cdots,\mu_N$. See \eqref{longform2}.

To store a short form of $\xi$ (in terms of the assists $\Omega$), we store $\lambda, \nu_{\max},$ $\gamma_1,\cdots,\gamma_{\nu_{\max}},$ $k_1,\cdots,k_{\nu_{\max}}$, and the short form of $\eta$. See \eqref{shortform}.

A linear map $S : \X(E) \oplus \cP \rightarrow \cP$ has \underline{$\Omega$-assisted depth $d$} provided that
\[
(f,P) \mapsto \partial^\alpha \left[ S(f,P) \right](0) \; \mbox{has} \; \Omega\mbox{-assisted depth} \; d, \; \mbox{for each} \; |\alpha| \leq m-1.
\]
To store a short form of the map $(f,P) \mapsto S(f,P)$, we store a short form of each of the linear functionals $(f,P) \mapsto \partial^\alpha \left[ S(f,P) \right](0)$ (for $\lv \alpha \rv \leq m-1$).

A linear map $T : \X(E) \oplus \cP \rightarrow \X$ has \underline{$\Omega$-assisted depth $d$} provided that
\[
(f,P) \mapsto \partial^\alpha \left[ T(f,P) \right](x) \; \mbox{has} \; \Omega\mbox{-assisted depth} \; d, \; \mbox{for each} \; x \in \R^n, \; |\alpha| \leq m-1.
\]
We can represent the map $T$ on a computer by giving an algorithm that accepts queries: A query consists of a point $x\in \R^n$. The response to a query is a short form of each of the linear functionals $(f,P) \mapsto \partial^\alpha \bigl[ T(f,P) \bigr](x)$ (for $\lv \alpha \rv \leq m-1$).

When we say that a linear functional $\omega$ has \underline{bounded depth}, we mean that its depth is bounded by a universal constant $C$.

When we say that a linear map $T$ (or linear functional $\xi$) has \underline{$\Omega$-assisted bounded depth}, we mean that $T$ (or $\xi$) has $\Omega$-assisted depth $d$, where $d$ is at most a universal constant $C$.

\comments{
Added the above lines
}

\section{Sets of Multi-indices} \label{sec_multi}

Let $\cM$ denote the collection of all multiindices $\alpha = (\alpha_1,\cdots,\alpha_n)$ of order $\lvert \alpha \rvert = \alpha_1 + \cdots + \alpha_n \leq m-1$. 

We define a total order relation $<$ on $\cM$ as follows: Given distinct $\alpha = (\alpha_1,\cdots,\alpha_n), \beta = (\beta_1,\cdots,\beta_n) \in \cM$, let $k \in \{1,\cdots,n\}$ be the maximal index such that $\alpha_1+\cdots + \alpha_k \neq \beta_1 + \cdots + \beta_k$. Then we write $\alpha < \beta$ if $\alpha_1+\cdots + \alpha_k <\beta_1 + \cdots + \beta_k$, and we write $\alpha > \beta$ otherwise.

We also define a total order relation $<$ on $2^{\cM}$. Given distinct subsets $\cA, \cB \subset \cM$, pick the minimal element $\alpha \in \cA \Delta \cB$ (with respect to the order relation defined above). Then we write $\cA < \cB$ if $\alpha \in \cA$, and we write $\cB < \cA$ otherwise. Here, $\cA \Delta \cB$ denotes the symmetric difference $(\cA \setminus \cB) \cup ( \cB \setminus \cA)$. Note that $\cM$ is minimal and that the empty set $\emptyset$ is maximal with respect to this order relation on $2^\cM$.

\begin{lem}\label{order_props}
The following properties hold.

\begin{itemize}
\item If $\alpha,\beta \in \cM$ and $\lvert \alpha \rvert < \lvert \beta \rvert$ then $\alpha < \beta$.
\item If $\alpha,\beta \in \cM$, $\alpha < \beta$ and $\lvert \gamma \rvert \leq m - 1 - \lvert \beta \rvert$, then $\alpha + \gamma < \beta + \gamma$.
\end{itemize}
\end{lem}

Given $\cA \subset \cM$, we say that $\cA$ is \underline{monotonic} if for every $\alpha \in \cA$ and $\gamma \in \cM$ with $\lvert \gamma \rvert \leq m - 1 - \lvert \alpha \rvert$, we have $\alpha + \gamma \in \cA$.

\begin{remk}\label{mon_rem}
Assume that $\cA\subset\cM$ is monotonic, $P \in \cP$, $x_0 \in \R^n$, and $\partial^\alpha P(x_0) = 0$ for all $\alpha \in \cA$. Then, for $x \in \R^n$ and $\alpha \in \cA$, we have
\[ \partial^\alpha P(x) = \sum_{| \gamma | \leq m-1-|\alpha|} \frac{1}{\gamma!} \partial^{\alpha + \gamma} P(x_0) \cdot (x-x_0)^\gamma = 0.\]
Hence, $\partial^\alpha P \equiv 0$ for any $\alpha \in \cA$.
\end{remk}

\section{Bases for the Space of Polynomials}
\label{bases_sec}

Let $\epsilon \in (0,1)$ be a given real number. We assume throughout this section that
$$\epsilon < \; \mbox{small enough constant determined by} \; m,n,p.$$

\subsection{Bases.}

Suppose we are given the following.
\begin{itemize}
\item A set of multiindices $\cA \subset \cM$.
\item A collection of polynomials $(P_\alpha)_{\alpha \in \cA}$ with each $P_\alpha \in \cP$.
\item A symmetric convex subset $\sigma \subset \cP$.
\item A point $x \in \R^n$, and real numbers $\Lambda \geq 1$, $\delta > 0$ (we call $\delta$ a ``lengthscale'').
\end{itemize}
We say that $(P_\alpha)_{\alpha \in \cA}$ forms an \underline{$(\cA,x,\epsilon,\delta)$-basis for $\sigma$} if the following conditions are met.
\begin{description}
\item[(B1)] $P_\alpha \in \epsilon \cdot \delta^{|\alpha| + n/p - m} \cdot \sigma$ \quad for all $\alpha \in \cA$.
\item[(B2)] $\partial^\beta P_\alpha(x) = \delta_{\alpha \beta} $ \qquad \quad for all $\alpha, \beta \in \cA$.
\item[(B3)] $|\partial^\beta P_\alpha(x)| \leq \epsilon \cdot \delta^{|\alpha| - |\beta|}$ \;\; for all $\alpha \in \cA$, $\beta \in \cM$, $\beta > \alpha$.
\end{description}
(Here, $\delta_{\alpha \beta}$ denotes the Kronecker delta: $\delta_{\alpha \beta} = 1$ if $\alpha = \beta$; $\delta_{\alpha \beta} = 0$ if $\alpha \neq \beta$.)
We say that $(P_\alpha)_{\alpha \in \cA}$ forms an \underline{$(\cA,x,\epsilon,\delta,\Lambda)$-basis for $\sigma$} if, in addition to \textbf{(B1)}-\textbf{(B3)}, the following condition is met.
\begin{description}
\item[(B4)] $|\partial^\beta P_\alpha(x)| \leq \Lambda \cdot \delta^{|\alpha| - |\beta|}$ \;\; for all $\alpha \in \cA$, $\beta \in \cM$.
\end{description}

\begin{remk}\label{basis_rem}
An $(\cA,x,\epsilon,\delta)$-basis is automatically an $(\cA,x,\epsilon',\delta')$-basis, for $\epsilon' \geq \epsilon$ and $\delta' \leq \delta$. An $(\cA,x,\epsilon,\delta,\Lambda)$-basis is automatically an $(\cA, x , \epsilon', \delta, \Lambda')$-basis, for $\epsilon' \geq \epsilon$ and $\Lambda' \geq \Lambda$. However, there is no simple relationship between an $(\cA,x,\epsilon,\delta,\Lambda)$-basis and an $(\cA,x,\epsilon,\delta',\Lambda)$-basis, due to the positive powers of $\delta$ appearing in condition \textbf{(B4)}.

Note that an $(\cA,x,\epsilon,\delta)$-basis is also an $(\cA,x, C^{m} \epsilon, C \delta)$-basis for any $C \geq 1$.
\end{remk}

\begin{remk}\label{rescale_base}
The notion of bases admits a natural rescaling, described below.

Given $P \in \cP$ define the polynomial $\tau_{x,\delta}(P) \in \cP$ by
\[\tau_{x,\delta}(P)(z) = P(\delta \cdot (z - x) + x).\]

Assume that $(P_\alpha)_{\alpha \in \cA}$ forms an $(\cA,x,\epsilon,\delta)$-basis for a symmetric convex set $\sigma \subset \cP$. Define the rescaled polynomials $\oP_\alpha = \delta^{-|\alpha|} \tau_{x,\delta}(P_\alpha)$ for $\alpha \in \cA$. Also define the convex set of polynomials
\[\overline{\sigma} = \{ \delta^{n/p - m} \tau_{x,\delta}(P) : P \in \sigma\}.\]
Then  \textbf{(B1-B3)} imply that $(\oP_{\alpha})_{\alpha \in \cA}$ forms an $(\cA,x,\epsilon,1)$-basis for $\overline{\sigma}$.

Similarly, under the assumption that $(P_\alpha)_{\alpha \in \cA}$ forms an $(\cA,x,\epsilon,\delta,\Lambda)$-basis for $\sigma$, we deduce that $(\oP_\alpha)_{\alpha \in \cA}$ forms an $(\cA,x,\epsilon,1,\Lambda)$-basis for $\overline{\sigma}$.

\end{remk}

\subsection{Tagged cubes} \label{tags_sec}

Assume that we are given a subset $E \subset \R^n$, a set of multiindices $\cA \subset \cM$, and a cube $Q \subset \R^n$.  

We say that $Q$ is \underline{tagged with $(\cA,\epsilon)$} provided that $\#(E \cap Q) \leq 1$ or there exists $\cA' \leq \cA$ such that
\[\sigma(Q) \; \mbox{has an} \; (\cA',x_Q,\epsilon,\delta_Q)\mbox{-basis} \qquad (\mbox{recall that } x_Q = \mbox{ center of } Q).\]

\begin{remk} \label{tag_rem} Note that every cube is tagged with $(\cA,\epsilon)$ with $\cA = \emptyset$. 
\end{remk}

Let real numbers $\eta \in (0,1)$,  $\Lambda \geq 1$, and $\delta>0$ be given. Let $\cA$ be a set of multiindices of order $\leq m-1$, and let $M = (M_{\alpha \beta})_{\alpha, \beta \in \cA}$ be a matrix (with real entries).

We say that $M$ is \underline{$(\eta,\Lambda, \delta)$-near triangular} provided that
\[
\lvert M_{\alpha \beta} - \delta_{\alpha \beta} \rvert \leq
\left\{
\begin{array}{ll}
\eta \cdot \delta^{|\alpha|-|\beta|} &: \mbox{if}\; \alpha, \beta \in \cA, \; \beta \geq \alpha \\
\Lambda \cdot \delta^{|\alpha|-|\beta|}  &: \mbox{if} \; \alpha, \beta \in \cA.\\
\end{array} \right.
\]

\begin{lem}\label{lem00}
If the matrices $M = (M_{\alpha\beta})_{\alpha,\beta \in \cA}, \widetilde{M} = (\widetilde{M}_{\alpha \beta})_{\alpha , \beta \in \cA}$ are $(\eta,\Lambda,\delta)$-near triangular and $(\widetilde{\eta},\widetilde{\Lambda},\delta)$-near triangular, respectively, then $M \widetilde{M}$ is $(\check{\eta}, \check{\Lambda}, \delta)$-near triangular if $ \check{\eta} < 1$, where $\check{\eta} = C \cdot ( \eta \widetilde{\Lambda} + \widetilde{\eta} \Lambda)$ and $\check{\Lambda} = C \Lambda \widetilde{\Lambda}$ for a universal constant $C$.
\end{lem}
\begin{proof}
Suppose that $\alpha,\beta \in \cA$ and $\beta > \alpha$. Then
\[(M \widetilde{M})_{\alpha \beta}= \sum_{\gamma \in \cA} M_{\alpha \gamma} \widetilde{M}_{\gamma \beta}.\]
If $\gamma \in \cA$ and $\gamma >\alpha$, then the corresponding term in the above sum is bounded in magnitude by $\eta  \widetilde{\Lambda} \delta^{|\alpha| - |\gamma|} \delta^{|\gamma| - |\beta|} = \eta  \widetilde{\Lambda}  \delta^{|\alpha| - |\beta|}$. Alternatively, if $\gamma \in \cA$ and $\gamma < \beta $ then the relevant term is bounded in magnitude by $\widetilde{\eta} \Lambda \delta^{|\alpha| - |\beta|}$. The total number of terms is at most $D = \dim(\cP)$, hence we see that $\lvert (M \widetilde{M})_{\alpha \beta} \rvert \leq \check{\eta} \delta^{|\alpha| - |\beta|}$ with $\check{\eta}$ as in the statement of the lemma.

Next, observe that
\[(M \widetilde{M})_{\alpha \alpha}= \sum_{\gamma \in \cA} M_{\alpha \gamma} \widetilde{M}_{\gamma \alpha}.\]
If $\gamma \in \cA$ and either $\gamma > \alpha$ or $\gamma < \alpha$, then the relevant term in the above sum is bounded in magnitude either by $\eta  \widetilde{\Lambda} \delta^{|\alpha| - |\gamma|} \delta^{|\gamma| - |\alpha|} = \eta \widetilde{\Lambda}$ or by $\widetilde{\eta} \Lambda \delta^{|\alpha| - |\gamma|} \delta^{|\gamma| - |\alpha|} = \widetilde{\eta} \Lambda$, respectively. If $\gamma = \alpha$, then the relevant term in the sum is equal to $M_{\alpha \alpha} \widetilde{M}_{\alpha \alpha} = (1 + \cO(\eta))(1+ \cO(\widetilde{\eta})) = 1 + \cO(\eta +  \widetilde{\eta})$. Hence, we find that $\lvert (M \widetilde{M})_{\alpha \alpha} - 1 \rvert \rvert \leq \check{\eta}$.

Finally, we assume that $\alpha,\beta \in \cA$ and $\beta < \alpha$. Then the estimates $\lvert M_{\alpha \gamma} \rvert \leq \Lambda  \delta^{|\alpha| - |\gamma|}$ and $\lvert \widetilde{M}_{\gamma \beta} \rvert \leq \widetilde{\Lambda} \delta^{|\gamma| - |\beta|}$ imply that $\lvert (M \widetilde{M})_{\alpha \beta} \rvert \leq \check{\Lambda} \delta^{|\alpha| - |\beta|}$ with $\check{\Lambda}$ as in the statement of the lemma.

This concludes the proof of Lemma \ref{lem00}. 

\end{proof}

\begin{lem}\label{lem0} Assume that $\eta \Lambda^D $ is less than a small enough constant depending on $m$ and $n$. Then the following holds.
\begin{itemize}
\item If the matrix $M = (M_{\alpha \beta})_{\alpha,\beta \in \cA}$ is $(\eta,\Lambda,\delta)$-near triangular, then $M$ is invertible and the inverse matrix $M^{-1}$ is $(C \eta \Lambda^D, C \Lambda^D,\delta)$-near triangular. 

\end{itemize}
Here, $D = \dim(\cP)$, and $C$ depends only on $m$ and $n$.
\end{lem}

\begin{proof}
Let $Y = (\delta_{ij} + X_{ij})_{i,j=1,\cdots,K}$ be a $K \times K$ matrix, where the $X_{ij}$ are variables. Let $(Y^{-1})_{ab}$ be the entries of $Y^{-1}$ ($a,b=1,\cdots,K$). Cramer's rule gives
\[\det Y = P(X) \;\; \mbox{and} \;\; (\det Y) \cdot \left[ (Y^{-1})_{ab} - \delta_{ab} \right] = P_{ab}(X),\]
where $P(X), P_{ab}(X)$ are $K$-th degree polynomials in $X=(X_{ij})_{i,j=1,\cdots,K}$. In $P$, $P_{ab}$, we separate the monomials containing only the variables $X_{ij}$ with $i<j$ from the monomials containing at least one variable $X_{ij}$ with $i \geq j$. We write $P = P_0 + P_1$ and $P_{ab} = P_{ab,0} + P_{ab,1}$, where the monomials in $P_0,P_{ab,0}$ contain only $X_{ij}$ with $i < j$, and the monomials in $P_1,P_{ab,1}$ contain at least one $X_{ij}$ with $i \geq j$.

Suppose that $X_{ij} = 0$ for $i \geq j$. Then $Y$ is upper triangular with $1$'s on the main diagonal, hence the same is true of $Y^{-1}$. It follows that $P_0 \equiv 1$,  and $P_{ab,0}  \equiv 0$ for $a \geq b$. Now we drop the assumption that $X_{ij} = 0$ for $i \geq j$, and assume instead that $\lvert X_{ij} \rvert \leq \eta$ for $i \geq j$ and $\lvert X_{ij} \rvert \leq \Lambda$ for all $i,j$. (Here, $0 < \eta < 1 \leq \Lambda$.)

We write $C,C',C''$, etc. to denote constants depending only on $K$. By examining each monomial separately, we see that
\[
\lvert P_1(X) \rvert, \lvert P_{ab,1}(X) \rvert \leq C \eta \Lambda^{K-1}, \; \mbox{and} \; \lvert P_{ab,0}(X) \rvert \leq C \Lambda^K.
\]
Combining these estimates with our knowledge of $P_0$ and $P_{ab,0}$ ($a \geq b$), we conclude that
\[
\lvert \det Y - 1 \rvert \leq C \eta \Lambda^{K-1}, \; \text{and} \; \lvert (\det Y) \cdot ((Y^{-1})_{ab} - \delta_{ab} ) \rvert \leq \left\{
\begin{array}{c}
C \eta \Lambda^{K-1} \; \mbox{if} \; a \geq b \\
C \Lambda^K, \; \mbox{all} \; a,b.
\end{array}
\right.
\]
This immediately implies the following result:
\begin{itemize}
\item \textbf{(*)} Let $Y = (Y_{ij})$ be a $K \times K$ matrix, satisfying
\[\lvert Y_{ij} - \delta_{ij} \rvert \leq \eta \; \mbox{for} \; i \geq j \;\; \mbox{and} \;\; \lvert Y_{ij} \rvert \leq \Lambda \; (\mbox{all} \; i,j),\]
where $\eta \Lambda^{K-1}$ is less than a small enough constant depending only on $K$.

Then the inverse matrix $Y^{-1} = ((Y^{-1})_{ab})$ satisfies
\[\lvert (Y^{-1})_{ab} - \delta_{ab} \rvert \leq C \eta \Lambda^{K-1} \; \mbox{for} \; a \geq b \;\; \mbox{and} \;\; \lvert (Y^{-1})_{ab} \rvert \leq C \Lambda^K \; (\mbox{all} \; a,b).\]
\end{itemize}

Now let $M=(M_{\alpha \beta})_{\alpha,\beta \in \cA}$ be an $(\eta,\Lambda,\delta)$-near triangular matrix. Let $\alpha_1 < \alpha_2 < \cdots < \alpha_K$ be the elements of $\cA$. Applying \textbf{(*)} to the matrix $Y_{ij} = \delta^{\lvert \alpha_i \rvert - \lvert \alpha_j \rvert} M_{\alpha_i \alpha_j}$, we find that $M^{-1}$ is $(C \eta \Lambda^{K-1},C \Lambda^K,\delta)$-near triangular, as claimed in Lemma \ref{lem0}.
\end{proof}

\begin{lem}\label{pre_lem1}
Let $x \in \R^n$. Suppose $\sigma_2 \subset C \cdot \left[ \sigma_1 + \cB(x,\delta) \right]$, and suppose $\sigma_2$ has an $(\cA,x,\epsilon,\delta,\Lambda)$-basis. Then $\sigma_1$ has an $(\cA,x,C \epsilon\Lambda, \delta, C \Lambda)$-basis. Here, $C$ depends only on $m$, $n$, and $p$.
\end{lem}
\begin{proof}
By rescaling, we may assume without loss of generality that $\delta=1$. (See Remark \ref{rescale_base}.)

Let $(\widetilde{P}_\alpha)_{\alpha \in \cA}$ be an $(\cA,x,\epsilon,1,\Lambda)$-basis for $\sigma_2$.

Then
\begin{itemize}
\item $\widetilde{P}_\alpha \in \epsilon \sigma_2 \subset C \epsilon \left[ \sigma_1 + \cB(x,1)\right] \qquad (\alpha \in \cA)$

\item $\partial^\beta \widetilde{P}_\alpha(x) = \delta_{\beta \alpha} \qquad (\beta, \alpha \in \cA)$

\item $\lvert \partial^\beta \widetilde{P}_\alpha(x) \rvert \leq \epsilon \qquad (\alpha \in \cA, \beta \in \cM, \beta > \alpha)$

\item $\lvert \partial^\beta \widetilde{P}_\alpha(x) \rvert \leq \Lambda \qquad (\alpha \in \cA, \beta \in \cM)$.
\end{itemize}

The first bullet point above gives $\widetilde{P}_\alpha = P_\alpha + (\widetilde{P}_\alpha - P_\alpha)$ with $P_\alpha \in C \epsilon \sigma_1$ (all $\alpha \in \cA$) and $\lvert \partial^\beta (\widetilde{P}_\alpha - P_\alpha)(x) \rvert \leq C \epsilon$ (all $\alpha \in \cA, \beta \in \cM$).

The four bullet point properties of the $\widetilde{P}_\alpha$ now yield the following properties of the $P_\alpha$.
\begin{itemize}
\item $P_\alpha \in C \epsilon \sigma_1 \qquad (\alpha \in \cA)$

\item $\lvert \partial^\beta P_\alpha(x) - \delta_{\beta \alpha} \rvert \leq C \epsilon \qquad (\beta, \alpha \in \cA)$

\item $\lvert \partial^\beta P_\alpha(x) \rvert \leq C\epsilon \qquad (\alpha \in \cA, \beta \in \cM, \beta > \alpha)$

\item $\lvert \partial^\beta P_\alpha(x) \rvert \leq C\Lambda \qquad (\alpha \in \cA, \beta \in \cM)$.
\end{itemize}
Inverting the matrix $(\partial^\beta P_\alpha(x))_{\beta,\alpha \in \cA}$, we obtain a matrix $(M_{\alpha \gamma})_{\alpha,\gamma \in \cA}$ such that
\begin{align*}
\sum_{\alpha \in \cA} \partial^\beta P_\alpha(x) \cdot M_{\alpha \gamma} &= \delta_{\beta \gamma} \quad (\beta, \gamma \in \cA)\\
 \lvert M_{\alpha \gamma} - \delta_{\alpha \gamma} \rvert &\leq C \epsilon \quad (\alpha, \gamma \in \cA).
\end{align*}

Set $P^\#_\gamma = \sum_{\alpha \in \cA} P_\alpha M_{\alpha \gamma}$ for $\gamma \in \cA$. Then
\begin{itemize}
\item $P^\#_\gamma \in C \epsilon \sigma_1$ \qquad ($\gamma \in \cA$)

\item $\partial^\beta P^\#_\gamma(x) = \delta_{\beta \gamma}$ \quad ($\beta,\gamma \in \cA$)

\item $\lvert \partial^\beta P^\#_\gamma(x) \rvert \leq C \Lambda$ \quad ($\gamma \in \cA, \beta \in \cM$).
\end{itemize}
For $\beta > \gamma$, we have
$$\partial^\beta P^\#_\gamma(x) = \sum_{\alpha \leq \gamma} \partial^\beta P_\alpha(x) M_{\alpha \gamma} + \sum_{\alpha > \gamma} \partial^\beta P_\alpha(x) M_{\alpha \gamma}.$$
For the sum over $\alpha \leq \gamma $, we note that $\beta > \gamma \geq \alpha$, hence
$$\lvert \partial^\beta P_\alpha(x) \rvert \leq C \epsilon, \;\;\; \mbox{whereas} \;\;\; \lvert M_{\alpha \gamma} \rvert \leq C. $$
For the sum over $\alpha > \gamma$, we note that
$$\lvert \partial^\beta P_\alpha(x) \rvert \leq \Lambda \;\;\; \mbox{and} \;\;\; \lvert M_{\alpha \gamma} \rvert \leq C \epsilon.$$
Therefore,
\begin{itemize}
\item $\lvert \partial^\beta P^\#_\gamma(x)\rvert \leq C \epsilon \Lambda$ \quad ($\gamma \in \cA, \beta \in \cM, \beta >  \gamma$).
\end{itemize}
Thus, the $(P^\#_\gamma)_{\gamma \in \cA}$ form an $(\cA,x,C \epsilon \Lambda, 1, C \Lambda)$-basis for $\sigma_1$, which is what we asserted, since $\delta = 1$.
\end{proof}

\begin{lem}
\label{pre_lem6}
Let $x \in \R^n$, $\epsilon > 0$, and $1 \leq Z \leq \epsilon^{-1/2}$ be given. Suppose that $Z$ exceeds a large enough universal constant. Let $(P_\alpha)_{\alpha \in \cA}$ be an $(\cA,x,\epsilon,\delta)$-basis for $\sigma$, with
\begin{equation}
\label{p6} \max \bigl\{ \lvert \partial^\beta P_\alpha(x) \rvert \delta^{|\beta|-|\alpha|}  : \;  \alpha \in \cA, \; \beta \in \cM \bigr\} \geq Z.
\end{equation}
Then $\sigma$ has an $(\cA',x,Z^{-\kappa},\delta)$-basis, with $\cA' < \cA$. Here, $\kappa > 0$ is a universal constant.
\end{lem}
\begin{proof}

By rescaling, we may assume without loss of generality that $\delta=1$. (See Remark \ref{rescale_base}.)

Our hypothesis tells us that $(P_\alpha)_{\alpha \in \cA}$ is an $(\cA,x,\epsilon,1)$-basis for $\sigma$, meaning that
\begin{align}
\label{p3} & P_\alpha \in \epsilon \cdot \sigma;\\
\label{p4} & \partial^\beta P_\alpha(x) = \delta_{\alpha \beta} && (\alpha,\beta \in \cA); \; \mbox{and}\\
\label{p5} & |\partial^\beta P_\alpha(x)| \leq \epsilon && (\alpha \in \cA, \;\beta \in \cM, \; \beta > \alpha).
\end{align}

Pick the minimal multiindex $\oa \in \cA$ with $\displaystyle \max_{\beta \in \cM} |\partial^\beta P_\oa(x)| \geq Z$. (See \eqref{p6}.) Thus,
\begin{equation}
\label{stuff1}  |\partial^\beta P_\alpha(x)| < Z, \;\; \mbox{for all} \; \beta \in \cM, \;\alpha \in \cA, \; \alpha < \oa,
\end{equation}
and there exists $\beta_0 \in \cM$ such that
\begin{equation}
\label{p6c} |\partial^{\beta_0} P_{\oa}(x)| = \max_{\beta \in \cM} |\partial^\beta P_\oa(x)| \geq Z.
\end{equation}
Note that $\beta_0 \neq \oa$ by \eqref{p4}, and $\beta_0 \leq \oa$ by \eqref{p5}. Thus, $\beta_0 <\oa$.

Let the elements of $\cM$ between $\beta_0$ and $\oa$ be ordered as follows:
$$\beta_0 < \beta_1 < \cdots < \beta_k = \oa.$$
Note that $k + 1 \leq \# \cM = D$. 

Pick $\overline{k} \in \{0,\ldots,k\}$ such that 
$$|\partial^{\beta_{\overline{k}}} P_\oa(x)| Z^{{\overline{k}}/(D+1)} \geq |\partial^{\beta_\ell} P_\oa(x)| Z^{\ell/(D+1)} \;\; \mbox{for all} \;  \ell \in\{0,\cdots,k\}.$$
In particular, setting $\ob = \beta_{\overline{k}}$, we have
\begin{align}\label{stuffs} &|\partial^{\ob} P_\oa(x)|   \geq Z^{-D/(D+1)} |\partial^{\beta_0} P_\oa(x)|  \ogeq{\eqref{p6c}} Z^{1/(D+1)},  \;\; \mbox{and} \\
& |\partial^{\ob} P_\oa(x)| \geq Z^{1/(D+1)} |\partial^{\beta_l} P_\oa(x)| \quad \mbox{for} \;\; \ell=\overline{k}+1,\ldots,k. \label{stuffs2}
\end{align}
If $\beta \in \cM$, $\beta > \oa$, then \eqref{p5} and \eqref{stuffs} give
$$|\partial^\beta P_{\oa}(x)| \leq \epsilon \leq 1 \leq Z^{-1/(D+1)} \lvert \partial^{\ob} P_\oa (x) \rvert .$$ 
Meanwhile, if $\beta \in \cM$, $\ob < \beta \leq \oa$, then \eqref{stuffs2} states that 
\[ |\partial^\beta P_{\oa}(x)| \leq Z^{-1/(D+1)} |\partial^{\ob} P_{\oa}(x)|.\]
Thus,
\begin{equation}
\label{p8} |\partial^\beta P_{\oa}(x)| \leq Z^{-1/(D+1)} |\partial^{\ob} P_{\oa}(x)| \qquad \mbox{for any} \; \beta \in \cM, \; \beta > \ob.
\end{equation}

Note that $|\partial^{\ob} P_{\oa}(x)|> 1$, thanks to \eqref{stuffs}. Hence, \eqref{p4} and \eqref{p5} show that
\begin{equation}
\label{p9} \ob < \oa \; \mbox{and} \; \ob \notin \cA.
\end{equation}
Set $P_{\ob} = P_{\oa} / \partial^{\ob} P_{\oa}(x)$. Then
\begin{align}
\label{p10} & P_{\ob} \in \epsilon \cdot \sigma, \qquad\qquad\qquad\qquad\qquad\qquad\qquad\quad\; \mbox{from} \; \eqref{p3} \; \mbox{and} \; |\partial^{\ob} P_{\oa}(x)|> 1;\\
\label{p13} & |\partial^\beta P_{\ob}(x)| \leq Z^{-1/(D+1)} \quad (\beta \in \cM, \; \beta > \ob), \qquad \mbox{from} \; \eqref{p8}; \\
\label{p14} & |\partial^\beta P_{\ob}(x)| \leq Z^{D/(1+D)} \quad\; (\beta \in \cM), \qquad\qquad\quad\; \mbox{from} \; \eqref{p6c} \; \mbox{and}\; \eqref{stuffs}; \; \mbox{and} \\
\label{p14a} & \partial^{\ob}P_{\ob}(x) = 1.
\end{align}

Now define
\begin{equation*}
P^\#_\ob := P_\ob - \sum_{\gamma \in \cA, \gamma < \ob} \partial^\gamma P_{\ob}(x) P_\gamma.
\end{equation*}
We derive some estimates on $P^\#_{\ob}$. From \eqref{p4} we see that
\begin{equation*}
\partial^\alpha P^\#_\ob(x) = \partial^\alpha P_\ob(x) - \sum_{\gamma \in \cA, \gamma < \ob} \partial^\gamma P_\ob(x) \delta_{\alpha \gamma} = 0 \qquad (\alpha \in \cA, \alpha < \ob).
\end{equation*}
Thanks to \eqref{p5}, \eqref{p14}, and \eqref{p14a}, we have
\begin{equation}
\label{eq21a}
|\partial^\ob P^\#_\ob(x) - 1| \leq \sum_{\gamma \in \cA, \gamma < \ob} \lvert \partial^\gamma P_{\ob}(x) \rvert \cdot \lvert \partial^\ob P_{\gamma}(x) \rvert \leq C Z^{D/(D+1)} \epsilon.
\end{equation}
Meanwhile, if $\beta > \ob$ and $\gamma < \ob$, then  $\beta > \gamma$. Hence, by \eqref{p5}, \eqref{p13}, and \eqref{p14}, we have
\begin{equation*}
|\partial^\beta P^\#_\ob(x)| \leq \lvert \partial^\beta P_\ob(x) \rvert + \sum_{\gamma \in \cA, \gamma < \ob} \lvert \partial^\gamma P_{\ob}(x) \rvert \cdot \lvert \partial^\beta P_{\gamma}(x) \rvert  \leq Z^{-\frac{1}{D+1}} + C Z^{\frac{D}{D+1}} \epsilon \quad (\beta \in \cM, \beta > \ob).
\end{equation*}
From \eqref{p3}, \eqref{p10}, and \eqref{p14}, we have
\begin{equation*}
P^\#_\ob \in  \epsilon \cdot \sigma + C \cdot Z^{D/(D+1)} \epsilon \cdot \sigma \subseteq \left( C Z^{D/(D+1)} \epsilon \right) \cdot \sigma.
\end{equation*}
For each $\gamma \in \cA$, $\gamma < \ob$, \eqref{p9} implies that $\gamma < \oa$. Hence, from \eqref{stuff1} and \eqref{p14} we have
\begin{equation*}
\lvert \partial^\beta P^\#_\ob(y) \rvert \leq C \cdot Z^{(2D+1)/(D+1)} \qquad (\beta \in\cM).
\end{equation*}

Since $\epsilon \leq Z^{-1}$, if $\epsilon$ is sufficiently small then \eqref{eq21a} implies that $\partial^\ob P^\#_\ob(x) \in [1/2,2]$. Hence, we may define $\hP_{\ob} = P^\#_\ob / \partial^\ob P^\#_\ob(x)$. The estimates written above show that
\begin{align}
\label{p22} & \hP_\ob \in \left(C \cdot Z^{D/(D+1)} \epsilon\right) \cdot \sigma; \\
\label{p23} & \partial^\beta \hP_\ob(x) = \delta_{\beta \ob} \qquad \; (\beta \in \cA, \beta < \ob \;\; \mbox{or} \;\; \beta = \ob);\\
\label{p24} & |\partial^\beta \hP_\ob(x)| \leq C \cdot Z^{-1/(D+1)} + C \cdot Z^{D/(D+1)} \epsilon \qquad ( \beta \in \cM, \; \beta > \ob); \; \mbox{and}\\
\label{p24a} & |\partial^\beta \hP_\ob(x)| \leq C \cdot Z^{(2D+1)/(D+1)} \qquad (\beta \in \cM).
\end{align}

For each $\alpha \in \cA$, $\alpha < \ob$, set $\hP_\alpha = P_\alpha - \partial^{\ob} P_\alpha(x) \hP_\ob$. Note that $|\partial^{\ob} P_\alpha(x) | \leq \epsilon \leq 1$, thanks to \eqref{p5}. From \eqref{p3} and \eqref{p22}, we have
\begin{equation}
\label{p26} \hP_\alpha \in \left( C Z^{D/(D+1)}\epsilon \right) \cdot \sigma.
\end{equation} 
From \eqref{p5} and \eqref{p24a}, we have
\begin{align}
 \label{p28} \vert \partial^\beta \hP_\alpha(x) \rvert &\leq \lvert \partial^\beta P_\alpha(x)\rvert + \lvert \partial^{\ob} P_\alpha(x)\rvert \cdot \lvert \partial^\beta \hP_\ob(x) \rvert \leq \epsilon  + \epsilon \cdot C Z^{(2D+1)/(D+1)} \\
& \leq C \epsilon \cdot Z^{(2D+1)/(D+1)} \qquad ( \beta \in \cM, \; \beta > \alpha).\notag{}
\end{align}
From \eqref{p4} and \eqref{p23}, we have
\begin{align}
\label{p27}  \partial^{\beta} \hP_\alpha(x) &= \partial^\beta P_\alpha(x) - \partial^{\ob} P_\alpha(x) \partial^\beta \hP_\ob(x) \\
&= \left\{
\begin{array}{ll}
 \delta_{\alpha \beta} - \partial^{\ob} P_\alpha(x) \delta_{\beta \ob} \quad\quad\;\;\;= \delta_{\alpha \beta}&: \mbox{if}\; \beta \in \cA, \; \beta < \ob \\
\partial^{\ob} P_\alpha(x) - \partial^{\ob} P_\alpha(x) \delta_{\ob \;\ob} = 0 &: \mbox{if} \; \beta = \ob\\
\end{array}
\right. \notag{} \\
& = \delta_{\alpha \beta} \qquad\qquad\qquad\qquad\qquad\qquad\qquad \mbox{if either} \; \beta < \ob \; \mbox{and} \; \beta \in \cA,\;  \mbox{or} \; \beta = \ob. \notag{}
\end{align}

Set $\oA = \{\alpha \in \cA : \alpha < \ob\} \cup \{\ob\}$. Then \eqref{p9} shows that the minimal element of $\cA \Delta \oA$ is $\ob$. Therefore, $\overline{\cA} < \cA$.

From \eqref{p22}-\eqref{p24} and \eqref{p26}-\eqref{p27} we deduce that 
\begin{equation*}
(\hP_\alpha)_{\alpha \in \oA} \; \mbox{is an} \; (\oA,x,C \cdot (Z^{-\frac{1}{D+1}} + Z^{\frac{2D+1}{D+1}} \cdot \epsilon ),1)\mbox{-basis for} \; \sigma.
\end{equation*}
Since $\epsilon \leq Z^{-2}$ and $\delta=1$, this implies the conclusion of Lemma \ref{pre_lem6}.
\end{proof}

\begin{lem}\label{pre_lem2}
There exist constants $\kappa_1, \kappa_2 \in (0,1]$ depending only on $m$, $n$, and $p$ such that the following holds.

Let $x \in \R^n$. Suppose that $\sigma$ has an $(\cA,x,\epsilon,\delta)$-basis. 

Then there exists a multiindex set $\cA ' \leq \cA$, and there exist numbers $\kappa' \in [\kappa_1,\kappa_2]$ and $\Lambda \geq 1$ with $\epsilon^{\kappa'} \Lambda^{100D} \leq \epsilon^{\kappa'/2}$, such that $\sigma$ has an $(\cA',x,\epsilon^{\kappa'},\delta,\Lambda)$-basis.

Here, $D = \dim \cP$.

\end{lem}

\begin{proof}

By rescaling, we may assume without loss of generality that $\delta=1$. (See Remark \ref{rescale_base}.)

Let $\cA_0 = \cA$ and $L = 2^D$, and let $\kappa \in (0,1)$ be as in Lemma \ref{pre_lem6}. Set
\[ \epsilon_0 = \epsilon, \; \epsilon_\ell = \epsilon^{\frac{\kappa^\ell}{(200D)^\ell}} \; \mbox{and} \; Z_\ell = \epsilon^{-\frac{\kappa^{\ell - 1}}{(200D)^\ell}}  \;\;\; \mbox{for} \; \ell = 1,\cdots, L.\]
Note that $(Z_\ell)^{- \kappa} = \epsilon_\ell$ and $Z_\ell \leq (\epsilon_{\ell-1})^{-1/2}$ for each $\ell\geq 1$.

Let $(P^{(0)}_\alpha)_{\alpha \in \cA}$ be an $(\cA_0,x,\epsilon_0,1)$-basis for $\sigma$. We carry out the following iterative procedure:
\begin{description}
\item[Stage 0] From Lemma \ref{pre_lem6}, we have either
\begin{description}
\item[Case A]
$\lvert \partial^\beta P^{(0)}_\alpha(x)  \rvert \leq Z_1 \quad \mbox{for all} \; \alpha \in \cA_0, \beta \in \cM$

or
\item[Case B] There exist polynomials $(P^{(1)}_\alpha)_{\alpha \in \cA_1}$, such that $(P^{(1)}_\alpha)_{\alpha \in \cA_1}$ is an $(\cA_1,x,\epsilon_1,1)$-basis for $\sigma$, for some $\cA_1 < \cA_0$.
\end{description}
In Case A we terminate. In  Case B, we pass to

\item[Stage 1] From Lemma \ref{pre_lem6}, we have either
\begin{description}
\item[Case A]
$\lvert \partial^\beta P^{(1)}_\alpha(x)  \rvert \leq Z_2 \quad \mbox{for all} \; \alpha \in \cA_1, \beta \in \cM$

or
\item[Case B] There exist polynomials $(P^{(2)}_\alpha)_{\alpha \in \cA_2}$, such that $(P^{(2)}_\alpha)_{\alpha \in \cA_2}$ is an  $(\cA_2,x,\epsilon_2,1)$-basis for $\sigma$, for some $\cA_2 < \cA_1$.
\end{description}
In Case A we terminate. In Case B, we pass to Stage 2, and so forth.
\end{description}

Since $\cA_0 > \cA_1 > \cA_2 > \cdots$ and $\# \{ \cA : \cA \subseteq \cM \} = L$, there exists $\ell \in \{0,\cdots,L-1\}$ such that Case A occurs in Stage $\ell$. Thus,
\[ (P^{(\ell)}_\alpha)_{\alpha \in \cA_\ell} \; \mbox{is an} \; (\cA_\ell,x,\epsilon_\ell,1)\mbox{-basis for} \; \sigma, \; \mbox{with}\]
\[ \lvert \partial^\beta P^{(\ell)}_\alpha(x) \rvert \leq Z_{\ell+1} \quad \mbox{for all} \; \alpha \in \cA_\ell, \beta \in \cM.\]
Note that
\begin{align*}
\epsilon_\ell \cdot Z_{\ell+1}^{100D} =  \epsilon^{\frac{\kappa^\ell}{(200D)^\ell} - \frac{\kappa^\ell}{(200D)^{\ell+1}} 100D}  = \epsilon^{\frac{\kappa^\ell}{2 (200D)^\ell}} = \sqrt{\epsilon_\ell}.
\end{align*}
Note that $\epsilon_\ell = \epsilon^{\kappa'}$ for $\kappa' = \kappa^\ell/(200D)^\ell$. We set $\Lambda = Z_{\ell+1}$. Then the above conditions imply the conclusion of Lemma \ref{pre_lem2}, since $\delta=1$.
\end{proof}

\begin{lem}\label{pre_lem3}
Let $x,y \in \R^n$; assume that $\lvert x - y \rvert \leq C \delta$. Suppose $\sigma$ has an $(\cA,x,\epsilon,\delta,\Lambda)$-basis. Assume that $\epsilon \Lambda^D$ is less than a small enough constant depending on $m$, $n$, and $p$.

Then $\sigma$ has an $(\cA,y, C \epsilon \Lambda^{2D+1},\delta, C \Lambda^{2D+1})$-basis.
\end{lem}
\begin{proof}

By rescaling, we may assume that $\delta=1$. (See Remark \ref{rescale_base}.)

Let $(P_\alpha)_{\alpha \in \cA}$ be an $(\cA,x,\epsilon,\delta, \Lambda)$-basis for $\sigma$. Thus,
\begin{itemize}
\item $P_\alpha \in \epsilon \sigma \qquad (\alpha \in \cA)$

\item $\partial^\beta P_\alpha (x) = \delta_{\beta \alpha} \quad (\beta,\alpha \in \cA)$

\item $\lvert \partial^\beta P_\alpha(x) \rvert \leq \epsilon \quad (\alpha \in \cA, \beta \in \cM, \beta > \alpha)$

\item $\lvert \partial^\beta P_\alpha(x) \rvert \leq \Lambda \quad (\alpha \in \cA, \beta \in \cM)$.
\end{itemize}
For $\beta \in \cM$ and $\alpha \in \cA$ with $\beta > \alpha$, we have
$$\lvert \partial^\beta P_\alpha(y) \rvert = \left\lvert \sum_{\gamma} \frac{1}{\gamma!} \partial^{\beta + \gamma} P_\alpha(x) \cdot (y-x)^\gamma \right\rvert \leq C \epsilon.$$
Also,
$$\partial^\alpha P_\alpha(y) = \partial^\alpha P_\alpha(x) + \sum_{\gamma \neq 0} \frac{1}{\gamma!} \partial^{\alpha + \gamma} P_\alpha(x) \cdot (y-x)^\gamma = 1 +  \mbox{Error}, \;\; \mbox{where} \; \lvert \mbox{Error} \rvert \leq C \epsilon.$$
On the other hand, for general $\beta \in \cM$ and $\alpha \in \cA$, we have
$$\lvert \partial^\beta P_\alpha(y) \rvert = \left\lvert \sum_{\gamma} \frac{1}{\gamma!} \partial^{\beta + \gamma} P_\alpha(x) \cdot (y-x)^\gamma \right\rvert \leq C \Lambda.$$
Thus, $(\partial^\beta P_\alpha(y))_{\beta,\alpha \in \cA}$ is a $(C \epsilon,C \Lambda,1)$-near triangular matrix. Therefore, 
\begin{equation} \label{inv1}
\mbox{the inverse} \; (M_{\alpha \gamma})_{\alpha,\gamma \in \cA} \; \mbox{of} \; (\partial^\beta P_\alpha(y))_{\beta,\alpha \in \cA} \; \mbox{is a} \;  (C \epsilon \Lambda^{2D},\Lambda^{2D}) \mbox{-near triangular matrix.}
\end{equation}

For each $\gamma \in \cA$, we define $P^\#_\gamma = \sum_{\alpha \in \cA} P_\alpha \cdot M_{\alpha \gamma}$.
From the properties of the $P_\alpha$, we read off the following.
\begin{itemize}
\item $\lvert \partial^\beta P_\gamma^\#(y) \rvert \leq C \Lambda^{2D+1} \quad (\beta \in \cM, \gamma \in \cA)$

\item $\partial^\beta P_\gamma^\#(y) = \delta_{\beta \gamma} \qquad (\beta,\gamma \in \cA)$

\item $P_\gamma^\# \in C \epsilon \Lambda^{2D+1} \sigma \qquad (\gamma \in \cA)$.
\end{itemize}
Finally, for each $\beta \in \cM$ and $\gamma \in \cA$ with $\beta > \gamma$, we have
\begin{align*}
\lvert \partial^\beta P^\#_\gamma(y) \rvert &\leq  \sum_{\alpha \leq \gamma} \lvert \partial^\beta P_\alpha(y) \rvert \cdot \lvert M_{\alpha \gamma} \rvert + \sum_{\alpha > \gamma} \lvert \partial^\beta P_\alpha(y)\rvert \cdot \lvert M_{\alpha \gamma} \rvert \\
& \leq \sum_{\alpha \leq \gamma} C \epsilon \cdot C \Lambda^{2D+1} + \sum_{\alpha > \gamma} C \Lambda \cdot C \epsilon \Lambda^{2D}  \quad (\mbox{see \eqref{inv1}}) \\
& \leq C \epsilon \Lambda^{2D+1}.
\end{align*}
Thus, $(P^\#_\gamma)_{\gamma \in \cA}$ is an $(\cA,y,C \epsilon \Lambda^{2D+1}, 1, C \Lambda^{2D+1})$-basis for $\sigma$.
\end{proof}

\begin{lem}\label{pre_lem4}
There exists $\kappa > 0$ depending only on $m$, $n$, and $p$, such that the following holds. Let $x,y \in \R^n$. Suppose that $\sigma$ has an $(\cA,x,\epsilon,\delta)$-basis and that $\lvert x - y \rvert \leq C \delta$. Then, there exists $\cA' \leq \cA$ such that  $\sigma$ has an $(\cA',y, \epsilon^{\kappa},\delta)$-basis.
\end{lem}
\begin{proof}
By Lemma \ref{pre_lem2}, there exist $\kappa' \in [ \kappa_1,\kappa_2]$, $\cA' \leq \cA$, and $\Lambda \geq 1$, such that
$$\sigma \; \mbox{has an} \; (\cA',x,\epsilon^{\kappa'},\delta,\Lambda)\mbox{-basis}, \; \mbox{and} \; \epsilon^{\kappa'} \Lambda^{100D} \leq \epsilon^{\kappa'/2}.$$
Here, $\kappa_1, \kappa_2 > 0$ are universal constants. Thus, $\sigma$ has an $(\cA', y, C \epsilon^{\kappa'} \Lambda^{2D+1}, \delta, C \Lambda^{2D+1})$-basis, due to Lemma \ref{pre_lem3}. 

Note that $C \epsilon^{\kappa'} \Lambda^{2D+1} \leq C \epsilon^{\kappa'/2} \leq \epsilon^{\kappa_1/4}$, if $\epsilon$ is less than a small enough universal constant. Hence, $\sigma$ has an $(\cA',y,\epsilon^{\kappa_1/4},\delta)$-basis. This completes the proof of Lemma \ref{pre_lem4}.
\end{proof}

\begin{lem}\label{pre_lem5}
Suppose that $Q' \subset Q$ and $Q$ is tagged with $(\cA,\epsilon)$. Then $Q'$ is tagged with $(\cA,\epsilon^\kappa)$, where $\kappa > 0$ depends only on $m$, $n$, and $p$.
\end{lem}
\begin{proof}
Let $Q' \subset Q$, and suppose $Q$ is tagged with $(\cA,\epsilon)$. Then either $\#(Q \cap E) \leq 1$ or $\sigma(Q)$ has an $(\cA',x_Q,\epsilon,\delta_Q)$-basis for some $\cA' \leq \cA$. (Recall that $x_Q$ is the center of $Q$.)

If $\#(Q \cap E) \leq 1$ then $\#(Q' \cap E) \leq 1$, hence $Q'$ is tagged with $(\cA,\epsilon^\kappa)$ for any $\kappa > 0$, which implies the conclusion of Lemma \ref{pre_lem5}.

Suppose instead that
\[\sigma(Q) \; \mbox{has an} \; (\cA',x_Q,\epsilon,\delta_Q)\mbox{-basis with} \; \cA' \leq \cA.\] 
Then Lemma \ref{pre_lem2} implies that there exist $\kappa' \in [ \kappa_1,\kappa_2]$, $\Lambda \geq 1$, and $\cA'' \leq \cA'$, such that 
\[ \sigma(Q) \; \mbox{has an} \; (\cA'',x_Q,\epsilon^{\kappa'},\delta_Q,\Lambda)\mbox{-basis}, \] 
with $\epsilon^{\kappa'} \cdot \Lambda^{100D} \leq \epsilon^{\kappa'/2}$. Here, $\kappa_1, \kappa_2 > 0$ are universal constants. 

We have $\lvert x_{Q'} - x_Q \rvert \leq \delta_Q$, since $x_{Q'} \in Q' \subset Q$ and $x_Q \in Q$. Hence, Lemma \ref{pre_lem3} implies that
$$\sigma(Q) \; \mbox{has an} \; (\cA'',x_{Q'}, C \epsilon^{\kappa'} \Lambda^{2D+1}, \delta_Q, C \Lambda^{2D+1})\mbox{-basis}.$$
Since Lemma \ref{pre_lem0} gives $\sigma(Q) \subset C \left[ \sigma(Q') + \cB(x_{Q'},\delta_Q) \right]$, Lemma \ref{pre_lem1} implies that
$$\sigma(Q') \; \mbox{has an} \; (\cA'',x_{Q'},C \epsilon^{\kappa'} \Lambda^{10D}, \delta_Q, C \Lambda^{10D})\mbox{-basis}.$$
Since $\delta_{Q'} \leq \delta_Q$ and $C \epsilon^{\kappa'} \Lambda^{100D} \leq C \epsilon^{\kappa'/2} \leq \epsilon^{\kappa'/4} \leq \epsilon^{\kappa_1/4}$, we see that
$$\sigma(Q') \; \mbox{has an} \; (\cA'',x_{Q'},\epsilon^{\kappa_1/4},\delta_{Q'})\mbox{-basis}.$$
Recall that $\cA'' \leq \cA' \leq \cA$. This completes the proof of Lemma \ref{pre_lem5}.
\end{proof}

\subsection{Computing a basis}\label{sec_compbase}

We fix $x \in \R^n$ and $\cA \subset \cM$.

Recall that $\cP$ is the vector space of polynomials on $\R^n$ of degree at most $m-1$, and $D = \dim \cP$. We identify $\cP$ with $\R^D$, by identifying $P \in \cP$ with $\left( \partial^\alpha P(x) \right)_{\beta \in \cM}$. We define 
\[
\lvert P \rvert_{x} = \sum_{\beta} \lvert \partial^\beta P(x) \rvert.
\] 

\comments{
p. 39: Added this line. Changed constants for ooline{sigma} below.
}

Suppose we are given $\Lambda \geq 1$. In this subsection, we write $c(\Lambda)$, $C(\Lambda)$, etc. to denote constants determined by $m$,$n$,$p$, and $\Lambda$. We write $c$, $C$, etc. to denote constants determined by $m$,$n$, and $p$.

Let $q$ be a nonnegative quadratic form on $\cP$; thus, $q(P) \geq 0$ for all $P \in \cP$. We are given a symmetric $D$ x $D$ matrix $(q_{\beta \gamma})_{\beta, \gamma \in \cM}$, with
\begin{equation}
\label{mtt0}
q(P) = \sum_{\beta, \gamma \in \cM} q_{\beta \gamma} \cdot \partial^\beta P(x)  \cdot \partial^\gamma P(x) \qquad \text{for} \; P \in \cP.
\end{equation}

Let $\ooline{\sigma} \subset \cP$ be a symmetric convex set with
\begin{equation}
\label{mt0}
\left\{ P \in \cP :  q(P) \leq \Lambda^{-1} \right\} \subset \ooline{\sigma} \subset \left\{ P \in \cP : q(P) \leq \Lambda \right\}.
\end{equation}

Given $x \in \R^n$, $(q_{\beta \gamma})_{\beta, \gamma \in \cM}$, $\delta \in (0,\infty)$, and $\cA \subset \cM$, we want to compute (approximately) the least $\eta$ for which there exists a collection $(P_\alpha)_{\alpha \in \cA}$ of $(m-1)$-st degree polynomials  such that
\begin{align}
\label{mt1}
& P_\alpha \in \eta^{1/2} \delta^{|\alpha| + n/p - m} \cdot \ooline{\sigma} \quad (\alpha \in \cA)\\
\label{mt2}
& \partial^\beta P_\alpha(x) = \delta_{\beta \alpha} \quad (\beta,\alpha \in \cA)\\
\label{mt3}
& \lvert \partial^\beta P_\alpha(x) \rvert \leq \eta^{1/2} \delta^{|\alpha| - |\beta|} \quad (\alpha \in \cA,\; \beta \in \cM, \; \beta > \alpha).
\end{align}


To compute such an $\eta$, we introduce the quadratic form
\begin{align}\label{we1}
M^\delta ((P_\alpha)_{\alpha \in \cA}) &:= \sum_{\alpha \in \cA} q(\delta^{m - n/p - |\alpha|} P_\alpha) + \sum_{\substack{\alpha \in \cA, \beta \in \cM \\ \beta > \alpha}} (\delta^{|\beta| - |\alpha|} \partial^\beta P_\alpha(x))^2 \\
&= \sum_{\alpha \in \cA} \sum_{\beta, \gamma \in \cM} \delta^{2(m - n/p - |\alpha|)} q_{\beta \gamma } \cdot \partial^\beta P_\alpha(x) \cdot \partial^\gamma P_\alpha(x) + \sum_{\substack{\alpha \in \cA, \beta \in \cM \\ \beta > \alpha}} (\delta^{|\beta| - |\alpha|} \partial^\beta P_\alpha(x) )^2, \notag{}
\end{align}
on the affine subspace
\begin{equation}
\label{affine_sub}
H := \left\{ \vec{P} = (P_\alpha)_{\alpha \in \cA} : \partial^\beta P_\alpha(x) = \delta_{\beta \alpha} \; \mbox{for} \; \alpha, \beta \in \cA \right\}.
\end{equation}
For fixed $q$, $\cA$, $x$, we denote 
\[\eta_{\min}(\delta) = \min_{\vec{P} \in H} M^\delta(\vec{P}), \] 
which we regard as a function of $\delta \in (0,\infty)$.

The definition of $\eta_{\min}(\delta)$ shows that 
\begin{equation}
\label{abc1}
\left\{
\begin{aligned}
&\mbox{we can satisfy \eqref{mt1}, \eqref{mt2}, \eqref{mt3} if } \eta > C(\Lambda) \cdot \eta_{\min}(\delta), \mbox{ but} \\
&\mbox{we cannot satisfy \eqref{mt1}, \eqref{mt2}, \eqref{mt3} if } \eta < c(\Lambda) \cdot \eta_{\min}(\delta).
\end{aligned}
\right.
\end{equation}
Hence,
\begin{equation}
\label{abc2}
\ooline{\sigma} \mbox{ has an } (\cA,x,\eta^{1/2},\delta)\mbox{-basis if } \eta > C(\Lambda) \cdot \eta_{\min}(\delta)\mbox{, but not if } \eta < c(\Lambda) \cdot \eta_{\min}(\delta).
\end{equation}
Moreover, 
\begin{equation}
\label{slowvariance}
\eta_{\min}(\delta_1) \leq \eta_{\min}(\delta_2) \leq \left(\frac{\delta_2}{\delta_1} \right)^{2m} \eta_{\min}(\delta_1) \;\; \mbox{for} \;\; \delta_1 \leq \delta_2,
\end{equation}
which follows at once from the definition of $\eta_{\min}$.

\comments{
Changed c,C to c(Lambda), C(Lambda)
}

We now compute an expression that approximates the function $\eta_{\min}(\delta)$. 

We identify the index set $\cI = \{ (\alpha,\beta) : \alpha \in \cA, \beta \in \cM \setminus \cA \}$ with $\{1,\cdots,J\}$ ($J = (\# \cA) \cdot (\# \cM - \# \cA)$) by fixing an enumeration of $\cI$. We introduce coordinates $w = (w_j)_{1 \leq j \leq J} = (w_{\alpha \beta})_{\alpha \in \cA, \beta \in \cM \setminus \cA} \in \R^J$ on the space $H$. We denote
\begin{equation}
\label{coord1} P^w_{\alpha}(z) := \frac{1}{\alpha!} (z-x)^\alpha + \sum_{\beta \in \cM \setminus \cA} \frac{1}{\beta!} w_{\alpha \beta} (z - x)^\beta \; \mbox{for} \; w \in \R^J.
\end{equation}
We identify
\begin{equation}
\label{coord2}
\vec{P}^w = (P_\alpha^w)_{\alpha \in \cA} \in H \;\; \mbox{with} \;\; w = (w_j)_{1 \leq j \leq J} = (w_{\alpha \beta})_{\alpha \in \cA, \beta \in \cM \setminus \cA} \in \R^J.
\end{equation}

We wish to minimize the quadratic function $\widetilde{M}^\delta(w) := M^\delta(\vec{P}^w)$ over $w \in \R^J$. We write
\begin{align}\label{we2}
\widetilde{M}^\delta(w) &= \sum_{i,j=1}^J A_{ij}^\delta w_i w_{j}  - 2 \sum_{j=1}^J b_j^\delta w_j + m^\delta \\
&= \langle A^\delta w, w \rangle - 2 \langle b^\delta, w \rangle + m^\delta. \notag{}
\end{align}
Here, we specify a symmetric matrix $A^\delta = (A_{i j}^\delta)$, vector $b^\delta = (b_j^\delta)$, and scalar $m^\delta$ -- all functions of $\delta > 0$. Here, we write $\langle \cdot, \cdot \rangle$ to denote the standard Euclidean inner product on $\R^J$. We express
\begin{align}
\label{we2.1}
A_{ij}^\delta &= \sum_{\mu,\nu} c^{ij}_{\mu \nu} \delta^{\mu + \nu/p}, \\
\label{we2.2}
b_j^\delta & = \sum_{\mu,\nu} c^{j}_{\mu \nu} \delta^{\mu + \nu/p}, \\
\label{we2.3}
m^\delta & = \sum_{\mu,\nu} c_{\mu \nu} \delta^{\mu + \nu/p}, 
\end{align}
for computable coefficients $c^{i j}_{\mu \nu}$, $c^{j}_{\mu \nu}$, and $c_{\mu \nu}$; here, the sums on $\mu$,$\nu$ are finite, and $\mu$,$\nu$ are integers. The coefficient matrix $(c^{ij}_{\mu \nu})$ is symmetric with respect to $(i,j)$. We compute these expressions by writing equation \eqref{we1} in $w$-coordinates. The quadratic function $\widetilde{M}^\delta$ is nonnegative, hence $A^\delta \geq 0$.

Let $\epsilon > 0$. We eventually send $\epsilon$ to zero. We define
\begin{equation}\label{we3}
\widetilde{M}^{\epsilon,\delta}(w) := \langle A^{\epsilon,\delta} w , w \rangle - 2 \langle b^\delta, w \rangle + m^\delta, \;\;\; \mbox{where} \;\; A^{\epsilon,\delta}_{i j} := A^\delta_{i j} + \epsilon \delta_{i j}.
\end{equation}
Note that $A^{\epsilon,\delta}$ is invertible because $A^\delta \geq 0$ and $\epsilon > 0$. Cramer's rule shows that
\begin{equation}
\label{pseudoinv}
(A^{\epsilon,\delta})^{-1}_{i j} = \frac{ [A^{\epsilon,\delta}]_{ ij} }{\det (A^{\epsilon,\delta})} =  \frac{\sum_{k, k'}  a^{ i j}_{k k'} \delta^{\lambda_k} \epsilon^{k'}} {\sum_{\ell , \ell'} b_{\ell \ell'} \delta^{\gamma_\ell} \epsilon^{\ell'}}
\end{equation}
for computable numbers $a^{ i j}_{k k'}$, $b_{\ell \ell'}$, $\lambda_k$, and $\gamma_\ell$; here, the sums on $k$,$k'$,$\ell$,$\ell'$ are finite, and $k$,$k'$,$\ell$,$\ell'$ are nonnegative integers. We write $[A^{\epsilon,\delta}]_{ij}$ to denote the $(i,j)$-cofactor of the matrix $A^{\epsilon,\delta}$.

The minimum of the quadratic function $\widetilde{M}^{\epsilon,\delta}(w)$ is achieved when $\nabla \widetilde{M}^{\epsilon,\delta}(w) = 0$, namely, for $w = w^{\epsilon,\delta} := (A^{\epsilon, \delta})^{-1} b^\delta$. From \eqref{we3} we see that the minimum value is 
\begin{align*}
\widetilde{M}^{\epsilon,\delta}(w^{\epsilon,\delta}) &= \left\langle A^{\epsilon,\delta} \left( A^{\epsilon,\delta}\right)^{-1} b^\delta , \left( A^{\epsilon,\delta}\right)^{-1} b^\delta \right\rangle - 2 \left\langle b^\delta, \left( A^{\epsilon,\delta}\right)^{-1} b^\delta \right\rangle + m^\delta  \\
&= - \left\langle b^\delta, \left( A^{\epsilon,\delta}\right)^{-1} b^\delta  \right\rangle + m^\delta.
\end{align*}
 Therefore, based on \eqref{pseudoinv} and based on the form of the vector $b^\delta$ and scalar $m^\delta$ written in \eqref{we2.2}, \eqref{we2.3}, we learn that
\begin{equation}
\label{rational0}
\min_{w \in \R^J} \widetilde{M}^{\epsilon,\delta}(w) =  \frac{\sum_{k,k'}  a_{kk'} \delta^{\lambda_k} \epsilon^{k'} }{  \sum_{\ell,\ell'} b_{\ell\ell'} \delta^{\gamma_\ell} \epsilon^{\ell'} }.
\end{equation}
for computable numbers $a_{kk'}$, $b_{k k'}$, $\lambda_k$, and $\gamma_\ell$. We abuse notation, since the exponents $\lambda_k$ in \eqref{rational0} might differ from the exponents $\lambda_k$ in \eqref{pseudoinv}. However, note that the denominator of \eqref{rational0} matches the expression in the denominator of \eqref{pseudoinv}. Also note that all exponents in \eqref{rational0} have the form $\mu + \nu/p$ for $\mu,\nu \in \Z$.

The minimum value of $\widetilde{M}^{\epsilon,\delta}(w)$ converges to the minimum value of $\widetilde{M}^\delta(w)$ as $\epsilon \rightarrow 0^+$. Hence,
\[
\eta_{\min}(\delta) = \min_{w \in \R^J} \widetilde{M}^\delta(w)  = \lim_{\epsilon \rightarrow 0^+} \frac{\sum_{k,k'}  a_{k k'} \delta^{\lambda_k} \epsilon^{k'} }{ \sum_{\ell, \ell'} b_{\ell\ell'} \delta^{\gamma_\ell} \epsilon^{\ell'} }.
\]
Canceling the smallest powers of $\epsilon$ from the numerator and denominator above, we obtain the formula
\begin{equation}
\label{rational1}
\eta_{\min}(\delta) = \frac{\sum_{k=1}^{K}  a_k \delta^{\lambda_k} }{\sum_{\ell=1}^{L} b_{\ell} \delta^{\gamma_\ell}}
\end{equation}
for nonzero coefficients $a_k$, $b_\ell$. All the coefficients and exponents in \eqref{rational1} can be computed using the numbers in \eqref{rational0}. Both $\lambda_k$ and $\gamma_\ell$ have the form $\mu + \nu/p$ with $\mu,\nu \in \Z$. Here, we abuse notation, since $\lambda_k$ and $\gamma_\ell$ may be different from before. By collecting terms, we may assume that
\begin{equation}
\label{distinct1}
\lvert \lambda_k - \lambda_{k'} \rvert \geq c  \quad \mbox{and} \quad \lvert \gamma_\ell - \gamma_{\ell'} \rvert \geq c \qquad \mbox{for} \; k \neq k', \; \ell \neq \ell'.
\end{equation}
Here, both $K$ and $L$ are bounded by a universal constant, and $c > 0$ is a universal constant.

We have thus obtained a computable expression for $\eta_{\min}(\delta)$.

We approximate $\eta_{\min}(\delta)$ with a piecewise-monomial function using the following procedure.

\label{pp01}

\environmentA{Procedure: Approximate Rational Function.}

Given nonzero numbers $a_k$, $b_\ell$ and numbers $\lambda_k$, $\gamma_\ell$ satisfying \eqref{distinct1}, let
\[
\eta_{\min}(\delta) = \frac{\sum_{k=1}^{K}  a_k \delta^{\lambda_k} }{\sum_{\ell=1}^{L} b_{\ell} \delta^{\gamma_\ell}}.
\]
We assume that $K$, $L$ are bounded by a universal constant, and $\eta_{\min}(\delta) \geq 0$ for $\delta \in (0,\infty)$. We further assume that $\eta_{\min}(\delta)$ satisfies \eqref{slowvariance}.

We compute a collection of pairwise disjoint intervals $I_\nu$ with $(0,\infty) = \cup_\nu I_\nu$, and we compute numbers $c_\nu$, $\lambda_\nu$ associated to each $I_\nu$, such that the function
\[\eta_*(\delta) :=  c_\nu \cdot \delta^{\lambda_\nu} \quad \; \mbox{for} \; \delta \in I_\nu\]
satisfies
\[c \cdot \eta_{\min}(\delta) \leq \eta_*(\delta) \leq C \cdot \eta_{\min}(\delta) \quad \mbox{for all} \; \delta \in (0,\infty).\] 
The algorithm requires work and storage at most $C$. In particular, the number of distinct intervals $I_\nu$ is at most $C$.

\begin{proof}[\underline{Explanation}]

We will analyze separately the numerator and denominator in the rational function $\eta_{\min}(\delta)$. We define
\begin{align*}
\cB &:= \bigcup_{k \neq k'} I_{k k'} \; , \quad \mbox{with} \\
& \hspace{2cm} I_{k k'} := \left\{ \delta \in (0,\infty) : 5^{-1} \cdot \lvert a_k \delta^{\lambda_k} \rvert \leq \lvert a_{k'} \delta^{\lambda_{k'}} \rvert \leq 5 \cdot \lvert a_k \delta^{\lambda_k} \rvert \right\},
\end{align*}
and similarly
\begin{align*}
\cC &:= \bigcup_{\ell \neq \ell'} J_{\ell \ell'} \; , \quad \mbox{with} \\
& \hspace{2cm} J_{\ell\ell'} := \left\{ \delta \in (0,\infty) : 5^{-1} \cdot \lvert b_\ell \delta^{\gamma_\ell} \rvert \leq \lvert b_{\ell'} \delta^{\gamma_{\ell'}} \rvert \leq 5 \cdot \lvert b_\ell \delta^{\gamma_\ell} \rvert \right\}.
\end{align*}


Let $I \subset (0,\infty) \setminus \cB$. For $\delta \in I$, all elements in the set $\left\{ \lvert a_k \delta^{\lambda_k}\rvert : 1 \leq k \leq K \right\}$ are nonzero and are separated by at least a factor of $5$. We choose $k$ such that $\lvert a_k \delta^{\lambda_k} \rvert$ is maximized. By continuity, the same $k$ must work for all $\delta \in I$. By summing a geometric series, we have 
\[ 
\sum_{k' \neq k} \lvert a_{k'} \delta^{\lambda_{k'}} \rvert < 2^{-1} \cdot  \lvert a_k \delta^{\lambda_k} \rvert \quad \mbox{for all} \; \delta \in I.
\] 
We obtain the analogous estimate for $\cC$ using a similar argument. Hence, for any interval $I \subset (0,\infty) \setminus (\cB \cup \cC)$,
\begin{equation}
\label{pickmon1}
\left\{ \;\;
\begin{aligned}
& \mbox{there exist unique }  k=k(I) \in \{1,\cdots,K\} \; \mbox{and} \; \ell = \ell(I) \in \{1,\cdots,L \} \\
& \; \mbox{such that}\; \lvert a_k \delta^{\lambda_k} \rvert >  2 \sum_{k' \neq k} \lvert a_{k'} \delta^{\lambda_{k'}} \rvert \; \mbox{and} \; \lvert b_{\ell} \delta^{\gamma_{\ell}} \rvert >  2 \sum_{\ell' \neq \ell} \lvert b_{\ell'} \delta^{\gamma_{\ell'}}\rvert \; \mbox{for all} \; \delta \in I.
\end{aligned}
\right.
\end{equation}
The fact that $k=k(I)$ and $\ell=\ell(I)$ are unique is obvious from the above statement.

The endpoints of each nonempty interval $I_{k k'} = [h_{kk'}^-,h_{kk'}^+]$ are solutions of the equations $\lvert a_{k'} \delta^{\lambda_{k'}} \rvert = 5 \cdot \lvert a_k \delta^{\lambda_k} \rvert$ and $\lvert a_{k'} \delta^{\lambda_{k'}} \rvert = 5^{-1} \cdot \lvert a_k \delta^{\lambda_k} \rvert$, namely $h_{kk'}^-$ and $h_{kk'}^+$ are among the numbers
\[
\delta_1 = \left( 5 \left\lvert \frac{a_k}{a_{k'}} \right\rvert \right)^{1/(\lambda_{k'} - \lambda_k)} \qquad \mbox{and} \qquad \delta_2 = \left( 5^{-1} \left\lvert \frac{a_k}{a_{k'}} \right\rvert \right)^{1/(\lambda_{k'} - \lambda_k)}.
\]
Thus, the endpoints $h_{kk'}^- ,h_{kk'}^+$ of the intervals $I_{k k'}$ are computable. Moreover, using \eqref{distinct1} we see that
\begin{align*}
\int_{\cB} \frac{dt}{t} \leq \sum_{k \neq k'} \int_{I_{k k'}} \frac{dt}{t} &= \sum_{k \neq k'} \log \left( \frac{h^+_{k k'}}{h_{k k'}^-} \right) \\
&= \sum_{k \neq k'} \frac{1}{\lvert \lambda_{k'} - \lambda_{k} \rvert} \log(25) \leq  A, \;\; \mbox{where} \; A = A(m,n,p),
\end{align*}
For a similar reason, the endpoints of the intervals $J_{\ell \ell'}$ are computable and
\[\int_{\cC} \frac{dt}{t} \leq  A.\]

We replace each pair of intersecting intervals among the $I_{k k'}$ and $J_{\ell \ell'}$ with their union. We continue until all the remaining intervals are pairwise disjoint. Thus, we can compute pairwise disjoint closed intervals $I^{\bad}_{\nu}$ and pairwise disjoint open intervals $I_{\mu}$ such that 
\[
\cB \cup \cC =  \bigcup_{\nu=1}^{\nu_{\max}} I^\bad_\nu, \qquad (0,\infty) \setminus (\cB \cup \cC) = \bigcup_{\mu=1}^{\mu_{\max}} I_\mu,
\]
and
\begin{equation}
\label{logbound}
 \int_{I^\bad_\nu} \frac{dt}{t} \leq \int_{\cB \cup \cC} \frac{dt}{t} \leq 2 A \qquad \mbox{for each} \; \nu.
\end{equation}
Note that $\nu_{\max} \leq \# \{ I_{kk'} \} + \# \{ J_{\ell \ell'} \} \leq K^2 + L^2$ and $\mu_{\max} = \nu_{\max} + 1$, hence $\nu_{\max}$ and $\mu_{\max}$ are bounded by a universal constant.

From \eqref{pickmon1}, there exist indices $k=k(\mu) \in \{1,\cdots,K\}$ and $\ell=\ell(\mu) \in \{1,\cdots,L\}$ for each $\mu$ such that
\begin{equation}
\label{qwe1}
\lvert a_k \delta^{\lambda_k} \rvert > 2 \sum_{k' \neq k} \lvert a_{k'} \delta^{\lambda_{k'}} \rvert \;\; \mbox{and} \;\; \lvert b_{\ell} \delta^{\gamma_{\ell}} \rvert > 2 \sum_{\ell' \neq \ell} \lvert b_{\ell'} \delta^{\gamma_{\ell'}} \rvert \qquad \mbox{for all } \delta \in I_\mu.
\end{equation}
We can compute $k(\mu)$ and $\ell(\mu)$ for $\mu=1,\cdots,\mu_{\max}$, using a brute-force search.

Let $\mu$ be given, and set $k=k(\mu)$ and $\ell = \ell(\mu)$. We have $\eta_{\min}(\delta) = \frac{N(\delta)}{D(\delta)}$, with 
\[
N(\delta) = a_k \delta^{\lambda_k}   + \sum_{k' \neq k} a_{k'} \delta^{\lambda_{k'}} \;\; \mbox{and} \;\; D(\delta) = b_\ell \delta^{\gamma_\ell}   + \sum_{\ell' \neq \ell} b_{\ell'} \delta^{\gamma_{\ell'}}.
\]
From \eqref{qwe1}, we have
\[
\frac{1}{2}  \leq \frac{ N(\delta) }{ a_k \delta^{\lambda_k} } \leq \frac{3}{2} \quad \mbox{and} \quad \frac{1}{2}  \leq \frac{ D(\delta) }{ b_\ell \delta^{\gamma_\ell} } \leq \frac{3}{2}\ \qquad \mbox{for all} \; \delta \in I_\mu.
\]
Hence,
\[ (1/4) \cdot \eta_{\min}(\delta) \leq \frac{a_{k} \delta^{\lambda_{k}}}{b_{\ell} \delta^{\gamma_{\ell}}} \leq (9/4) \cdot \eta_{\min}(\delta) \qquad \mbox{for all} \; \delta \in I_\mu.
\]

We fix $\delta_\nu \in I_\nu^\bad$ for each $\nu$. Note that $e^{- 2A} \leq \delta/\delta_\nu \leq e^{2A}$ for all $\delta \in I_\nu^\bad$, by \eqref{logbound}. Hence, \eqref{slowvariance} implies that
\[
c \cdot \eta_{\min}(\delta) \leq \eta_{\min}(\delta_\nu) \leq C \cdot \eta_{\min}(\delta) \qquad \mbox{for all} \; \delta \in I_\nu^\bad.
\]

We define
\begin{equation}
\label{etastar}
\eta_*(\delta) = 
\left\{
\begin{array}{c}
\frac{a_{k(\mu)} \delta^{\lambda_{k(\mu)}}}{b_{\ell(\mu)} \delta^{\gamma_{\ell(\mu)}}}
 \;\; \mbox{if} \; \delta \in I_\mu \\
 \eta_{\min}(\delta_\nu) \;\; \mbox{if} \; \delta \in I_\nu^\bad.
\end{array}
\right.
\end{equation}

From the previous two paragraphs, we see that $\eta_{\min}(\delta)$ and $\eta_*(\delta)$ differ by at most a universal constant factor for all $\delta \in (0,\infty)$. 

The above computations clearly require work at most $C$.

That completes our description of the procedure \textsc{Approximate Rational Function}.
\end{proof}

We have computed a piecewise-monomial function $\eta_*(\delta)$ that differs from $\eta_{\min}(\delta)$ by at most a constant factor. 
Thus, we see that the properties \eqref{abc2} and \eqref{slowvariance} of $\eta_{\min}(\delta)$ imply the first and second bullet points below.

\comments{
p. 45: Constants
}

\label{page_FBTCB}

\environmentA{Algorithm: Fit Basis to Convex Body.}

Given a nonnegative quadratic form $q$ on $\cP$, given a point $x \in \R^n$, and given a set $\cA \subset \cM$: We compute a partition of $(0,\infty)$ into at most $C$ intervals $I_\nu$, and for each $I_\nu$ we compute real numbers $\lambda_\nu, c_\nu$ with $c_\nu \geq 0$, such that the function
\[\eta_*(\delta) := c_\nu \cdot \delta^{\lambda_\nu} \quad \mbox{for} \;\; \delta \in I_\nu\]
has the following properties.
\begin{itemize}
\item Let $\ooline{\sigma} \subset \cP$ be a symmetric convex set that satisfies $\{ q \leq \Lambda^{-1} \} \subset \ooline{\sigma} \subset \{ q \leq \Lambda \}$ for a real number $\Lambda \geq 1$. Then, for any $\delta > 0$, $\ooline{\sigma}$ has an $(\cA,x,\eta^{1/2},\delta)$-basis if $\eta > C(\Lambda) \cdot \eta_*(\delta)$, but not if $\eta < c(\Lambda) \cdot \eta_*(\delta)$.
\item Moreover, $c  \cdot \eta_*(\delta_1) \leq \eta_*(\delta_2) \leq C \cdot \eta_*(\delta_1)$ whenever $\frac{1}{10} \delta_1 \leq \delta_2 \leq 10 \delta_1$.
\item The components of the piecewise-monomial function $\eta_*(\delta)$, i.e., the intervals $I_\nu$ and the numbers $\lambda_\nu, c_\nu$, can be computed using work and storage at most $C$.
\end{itemize}
Here, $c > 0$ and $C \geq 1$ are constants depending only on $m$,$n$, and $p$, while $c(\Lambda) > 0$ and $C(\Lambda) \geq 1$ are constants depending only on $m$,$n$,$p$, and $\Lambda$.

\comments{
Changes statement here
}

\section{Algorithms for Linear Functionals}\label{sec_lf}

\environmentA{Algorithm: Compress Norms.} 

Fix $1 < p < \infty$ and $D \geq 1$.
Given linear functionals $\mu_1,\ldots,\mu_L : \R^D \rightarrow \R$ ($L \geq 1$), we produce linear functionals $\mu_1^*,\ldots,\mu_D^* : \R^D \rightarrow \R$ such that
$$c \cdot \sum_{i=1}^D \lvert \mu_i^*(v) \rvert^p \leq \sum_{i=1}^L \lvert \mu_i(v) \rvert^p \leq C \cdot \sum_{i=1}^D \lvert \mu_i^*(v) \rvert^p \quad \mbox{for all} \; v \in \R^D.$$
The work and storage used to do so are at most $C' L$. Here, $c,C,C'$ depend only on $D$ and $p$.

\begin{proof}[\underline{Explanation}] 

In this explanation, $c,C,C'$, etc., depend only on $D$ and $p$.

We start with a few elementary estimates.
Let $(\Omega, d\mu)$ be a probability space. Then, for $f : \Omega \rightarrow \R$ measurable, the mean $\overline{f} = \int_\Omega f d \mu$ satisfies
\[\lvert \overline{f} \rvert \leq \left(\int_\Omega  \lvert f \rvert^p d \mu \right)^{1/p},\]
hence
\[\lvert \overline{f} \rvert^p + \int_\Omega \lvert f  - \overline{f} \rvert^p d \mu \leq C \int_\Omega \lvert f \rvert^p d \mu.\]
Also,
\[\int_\Omega \lvert f \rvert^p d \mu \leq C \lvert \overline{f} \rvert^p + C \int_\Omega \lvert f - \overline{f}\rvert^p d \mu.\]
Applying the above to the function $f - b$ for a constant $b$, we find that
\begin{equation}
\label{cn1}
c \left\{ \lvert \overline{f} - b \rvert^p + \int_{\Omega} \lvert f - \overline{f} \rvert^p d \mu \right\} \leq \int_\Omega \lvert f - b \rvert^p d \mu \leq  C \left\{ \lvert \overline{f} - b \rvert^p + \int_\Omega \lvert f - \overline{f} \rvert^p d \mu \right\}
\end{equation}
with $c,C > 0$ depending only on $p$.

We now return to the task of constructing $\mu^*_1, \cdots, \mu^*_D$.

We proceed by induction on $D$. In the base case $D=1$, the construction of $\mu^*_1$ is trivial. For the induction step, fix $D \geq 2$, and assume we can carry out \textsc{Compress Norms} in dimension $D-1$. We show how to carry out that algorithm in dimension $D$.

Let $\mu_1,\cdots,\mu_L  : \R^D \rightarrow \R$ be given linear functionals. We write 
\begin{equation}\label{qfin1}
\mu_i(v_1,\cdots,v_D) = \pm \left[ \beta_i v_D - \widetilde{\mu}_i(v_1,\cdots,v_{D-1}) \right]
\end{equation} with $\beta_i \geq 0$, and we let $I = \{i : \beta_i \neq 0 \}$. If $I$ is empty, then we succeed simply by setting $\mu_D^* = 0$ and invoking \textsc{Compress Norms} in dimension $D-1$. Suppose $I$ is non-empty. Let
\[
\mathbf{B} := \sum_{j \in I} \beta_j^p.
\]
We view $I$ as a probability space, with
\[\mbox{Prob}(i) := \beta_i^p/ \mathbf{B} \;\; \mbox{for} \; i \in I.\]
Then
\[\sum_{i \in I} \lvert \mu_i(v_1,\cdots,v_D) \rvert^p =  \mathbf{B} \cdot \sum_{i \in I} \mbox{Prob}(i) \cdot \left\lvert v_D - \beta_i^{-1} \widetilde{\mu}_i(v_1,\cdots,v_{D-1}) \right\rvert^p .\]
Invoking \eqref{cn1}, with $b = v_D$ and $f(i) = \beta_i^{-1} \widetilde{\mu}_i(v_1,\cdots,v_{D-1})$, we see that
\begin{align*}
c \sum_{i \in I} \lvert \mu_i(v_1,\cdots,v_D) \rvert^p & \leq   \mathbf{B} \cdot \biggl\{ \lvert v_D - \overline{\mu}(v_1,\cdots,v_{D-1}) \rvert^p \\
& \qquad + \sum_{i \in I} \mbox{Prob}(i) \cdot \left\lvert \overline{\mu}(v_1,\cdots,v_{D-1})  - \beta_i^{-1} \widetilde{\mu}_i(v_1,\cdots,v_{D-1}) \right\rvert^p  \biggr\}  \\
& \leq C \sum_{i \in I} \lvert \mu_i(v_1,\cdots,v_D)\rvert^p,
\end{align*}
where 
\begin{equation}
\label{mubar}
\overline{\mu}(v_1,\cdots,v_{D-1}) := \sum_{i \in I} \mbox{Prob}(i) \cdot \bigl\{ \beta_i^{-1} \widetilde{\mu}_i(v_1,\cdots,v_{D-1}) \bigr\}.
\end{equation}
Consequently, $\displaystyle \sum_{i=1}^L \lvert \mu_i(v_1,\cdots,v_D ) \rvert^p$ differs by at most a factor of $C$ from
\begin{align*}
&\mathbf{B} \cdot \lvert v_D - \overline{\mu}(v_1,\cdots,v_{D-1}) \rvert^p \\
& + \biggl\{ \mathbf{B} \cdot \sum_{i \in I} \mbox{Prob}(i) \cdot \lvert \overline{\mu}(v_1,\cdots,v_{D-1}) - \beta_i^{-1} \widetilde{\mu}_i(v_1,\cdots,v_{D-1}) \rvert^p \\
& \qquad + \sum_{i \notin I} \lvert \widetilde{\mu}_i(v_1,\cdots,v_{D-1}) \rvert^p \biggr\} \quad = \\
& \boxed{ \begin{aligned} \mathbf{B} \cdot \lvert v_D - \overline{\mu}(v_1,\cdots,v_{D-1}) \rvert^p + \biggl\{ \sum_{i =1}^L \lvert \beta_i \overline{\mu}(v_1,\cdots,v_{D-1}) - \widetilde{\mu}_i(v_1,\cdots,v_{D-1}) \rvert^p \biggr\}
\end{aligned}}
\end{align*}
where
\[\boxed{ \overline{\mu}(v_1,\cdots,v_{D-1}) := \mathbf{B}^{-1} \cdot \sum_{i=1}^L \beta_i^{p-1} \widetilde{\mu}_i(v_1,\cdots,v_{D-1})}.\]
Applying \textsc{Compress Norms} in dimension $D-1$ to the expression in curly brackets in the first box above, we obtain functionals $\mu^*_1, \cdots, \mu^*_{D-1} : \R^{D-1} \rightarrow \R$ such that $\displaystyle \sum_{i=1}^{D-1} \lvert \mu^*_i(v_1,\cdots,v_{D-1}) \rvert^p$ differs by at most a factor of $C$ from that expression in curly brackets.

Setting 
\begin{equation}
\label{muDstar} 
\mu_D^*(v_1,\cdots,v_D) := \mathbf{B}^{1/p} \cdot [ v_D - \overline{\mu}(v_1,\cdots,v_{D-1}) ],
\end{equation}
we see that $\displaystyle \sum_{i=1}^{D-1} \lvert \mu^*_i(v_1,\cdots,v_{D-1}) \rvert^p + \lvert \mu^*_D(v_1,\cdots,v_D) \rvert^p$ differs by at most a factor of $C$ from $\displaystyle \sum_{i=1}^{L} \lvert \mu_i(v_1,\cdots,v_{D}) \rvert^p$.

This completes our explanation of \textsc{Compress Norms}; note that the work and storage required are as promised.

\end{proof}

\environmentA{Algorithm: Optimize via Matrix.}

Given $1 < p < \infty$ and given a matrix $(a_{\ell j})_{\substack{1 \leq \ell \leq L \\ 1 \leq j \leq J}}$, we compute a matrix $(b_{j \ell})_{\substack{1 \leq j \leq J \\ 1 \leq \ell \leq L}}$ for which the following holds.

Let $y_1,\cdots, y_L$ be real numbers. Define
\[x_j^* = \sum_{\ell=1}^L b_{j \ell} y_\ell  \quad \mbox{for} \;\; j=1,\cdots,J.\]
Then, for any real numbers $x_1,\cdots,x_J$, we have
\[\sum_{\ell=1}^L \biggl\lvert y_\ell + \sum_{j=1}^J a_{\ell j} x_j^* \biggr\rvert^p \leq C_1 \cdot \sum_{\ell=1}^L \biggl\lvert y_\ell + \sum_{j=1}^J a_{\ell j} x_j \biggr\rvert^p \]
with $C_1$ depending only on $J$ and $p$.

The work and storage used to compute $(b_{j \ell})_{\substack{1 \leq j \leq J \\ 1 \leq \ell \leq L}}$ are at most $C L$, where  $C$ depends only on $J$.

\begin{proof}[\underline{Explanation}] 

We write $c$, $C$, $C'$, etc. to denote constants depending only on $J$, and $c(p)$, $C(p)$, etc. to denote constants depending only on $J$ and $p$.

\underline{For the case $J=1$} of our algorithm, we proceed as follows.

Let $(a_{\ell 1})_{1 \leq \ell \leq L}$ be a given matrix.

If $(a_{\ell1})_{1 \leq \ell \leq L}$ is identically zero, then the conclusion holds for any choice of $(b_{1 \ell})_{1 \leq \ell \leq L}$ if we take $C_1 = 1$. 

We suppose instead that $(a_{\ell1}) \neq (0)$. Let
\begin{equation}
\label{mn1}
b_{1\ell} :=  - \left(  \sum_{\ell'=1}^L \lvert a_{\ell' 1} \rvert^p \right)^{-1} \cdot \lvert a_{\ell 1} \rvert^{p-1}  \sgn(a_{\ell 1}) \quad \mbox{for}\; 1 \leq \ell \leq L,
\end{equation}
where $\sgn$ denotes the signum function: $\sgn(\eta) := 1$ for $\eta \geq 0$, $\sgn(\eta) := -1$ for $\eta < 0$. We compute the matrix $(b_{1\ell})_{1 \leq \ell \leq L}$ using work and storage at most $C L$.

For given real numbers $y_1,\cdots,y_L$ we set
\begin{equation*}
x^* = \sum_{\ell=1}^L b_{1\ell} y_\ell =  - \sum_{\ell : a_{\ell1} \neq 0} \frac{y_\ell \cdot \lvert a_{\ell 1} \rvert^{p}}{a_{\ell1} \cdot \sum_{\ell'} \lvert a_{\ell' 1} \rvert^p}.
\end{equation*}

Define a probability measure $d\mu$ and function $f$ on $\{1,\cdots,L\}$ by setting $d \mu(\ell) = \frac{ \lvert a_{\ell1} \rvert^p}{\sum_{k} \lvert a_{k1} \rvert^p}$, and setting $f(\ell) = y_\ell/a_{\ell1}$ if $a_{\ell1} \neq 0$ and $f(\ell) = 0$ otherwise. We then have $x^* = - \int f d \mu$. 

By applying \eqref{cn1} we see that $\int \lvert f + x^* \rvert^p d \mu \leq C(p) \int \lvert f + x \rvert^p d \mu$ for any $x \in \R$. Therefore,
\[
\sum_{\ell : a_{\ell1} \neq 0 } \lvert y_\ell/a_{\ell1} +  x^* \rvert^p \lvert a_{\ell1} \rvert^p \leq C(p) \sum_{\ell : a_{\ell1} \neq 0 } \lvert y_\ell/a_{\ell1} + x \rvert^p \lvert a_{\ell1} \rvert^p \quad \mbox{for any} \;  x \in \R.\]
This gives the desired conclusion in the case $J=1$.

\underline{For the general case}, we use induction on $J$.

Let $J \geq 2$, and let $1 < p < \infty$ and $(a_{\ell j})_{\substack{1 \leq \ell \leq L \\ 1 \leq j \leq J}}$ be given. Then
\begin{equation}\label{ovm1}
\sum_{\ell=1}^L \lvert y_\ell + \sum_{j=1}^J a_{\ell j} x_j \rvert^p = \sum_{\ell=1}^L \lvert \widehat{y}_\ell + \sum_{j=1}^{J-1} a_{\ell j} x_j \rvert^p
\end{equation}
with
\begin{equation}\label{ovm2}
\widehat{y}_\ell = y_\ell + a_{\ell J} \cdot x_J \qquad \mbox{for} \; \ell=1,\cdots,L.
\end{equation}

Applying our algorithm \textsc{Optimize via Matrix} recursively to $1 < p < \infty$ and $(a_{\ell j})_{\substack{1 \leq \ell \leq L \\ 1 \leq j \leq J-1}}$, we produce a matrix $(\widehat{b}_{j \ell})_{\substack{1 \leq j \leq J-1 \\ 1 \leq \ell \leq L}}$, for which the following holds.

\begin{itemize}
\item Let $\widehat{y}_1,\cdots,\widehat{y}_L$ be real numbers, and set
\begin{equation}\label{ovm3}
\widehat{x}_j = \sum_{\ell=1}^L \widehat{b}_{j \ell} \widehat{y}_\ell \quad \mbox{for} \; j=1,\cdots,J-1.
\end{equation}
Then, for any real numbers $x_1,\cdots,x_{J-1}$, we have
\begin{equation}\label{ovm4}
\sum_{\ell=1}^L \lvert \widehat{y}_\ell + \sum_{j=1}^{J-1} a_{\ell j} \widehat{x}_j \rvert^p \leq C(p) \sum_{\ell=1}^L \lvert \widehat{y}_\ell + \sum_{j=1}^{J-1} a_{\ell j} x_j \rvert^p.
\end{equation}
\end{itemize}

From \eqref{ovm1}-\eqref{ovm4}, we draw the following conclusion.

Let real numbers $y_1,\cdots,y_L$, and $x_1,\cdots,x_J$ be given. Define $\widehat{y}_1,\cdots, \widehat{y}_L$ by \eqref{ovm2}, next define $\widehat{x}_1,\cdots,\widehat{x}_{J-1}$ by \eqref{ovm3}, and finally set 
\begin{equation}\label{ovm5}
\widehat{x}_J = x_J.
\end{equation}
Then 
\begin{equation}\label{ovm6}
\sum_{\ell=1}^L \lvert y_\ell + \sum_{j=1}^J a_{\ell j} \widehat{x}_j \rvert^p \leq C(p) \sum_{\ell=1}^L \lvert y_\ell + \sum_{j=1}^J a_{\ell j} x_j \rvert^p
\end{equation}
and
\begin{equation}\label{ovm7}
\widehat{x}_j = \sum_{\ell=1}^L \widehat{b}_{j \ell} \cdot (y_\ell + a_{\ell J} \widehat{x}_J ) \quad \mbox{for} \; j=1,\cdots,J-1.
\end{equation}
Thus,
\begin{align}\label{ovm8}
&\widehat{x}_j = \sum_{\ell=1}^L \widehat{b}_{j \ell} y_\ell + g_j \widehat{x}_J \quad \mbox{for} \; j=1,\cdots, J-1, \; \mbox{where} \\
\label{ovm9}
&g_j = \sum_{\ell=1}^L \widehat{b}_{j \ell} a_{\ell J} \quad \mbox{for} \; j=1,\cdots, J-1.
\end{align}

Next, note that
\begin{align*}
y_\ell + \sum_{j=1}^{J} a_{\ell j} \widehat{x}_j & = y_\ell + \sum_{j=1}^{J-1} a_{\ell j} \left[ \sum_{\ell' = 1}^L \widehat{b}_{j \ell'} y_{\ell'} +  g_j \widehat{x}_J  \right]  + a_{\ell J} \widehat{x}_J \\
&= \left\{ y_\ell + \sum_{j=1}^{J-1} a_{\ell j} \sum_{\ell'=1}^L \widehat{b}_{j \ell'} y_{\ell'} \right\} + \left\{ a_{\ell J} + \sum_{j=1}^{J-1} a_{\ell j} g_j \right\} \widehat{x}_J.
\end{align*}
We set
\begin{equation}\label{ovm10}
y_\ell^{\ouch} = y_\ell + \sum_{j=1}^{J-1} a_{\ell j} \sum_{\ell'=1}^L \widehat{b}_{j \ell'} y_{\ell'} \quad \mbox{for} \; \ell=1,\cdots,L
\end{equation}
and
\begin{equation}\label{ovm11}
h_\ell = a_{\ell J} + \sum_{j=1}^{J-1} a_{\ell j}  g_j \quad \mbox{for} \; \ell=1,\cdots,L.
\end{equation}
Thus,
\begin{equation} \label{ovm12} 
\sum_{\ell=1}^L \lvert y_\ell + \sum_{j=1}^J a_{\ell j} \widehat{x}_j \rvert^p = \sum_{\ell=1}^L \lvert y^\ouch_\ell + h_\ell \widehat{x}_J \rvert^p.
\end{equation}
Here, \eqref{ovm12} holds whenever $\widehat{x}_1,\cdots,\widehat{x}_{J-1}$ are determined from $\widehat{x}_J$ via \eqref{ovm8}. 

Note that it is too expensive to compute $y^\ouch_\ell$ for all $\ell$ ($1 \leq \ell \leq L$); that computation would require $\sim L^2 J$ work. However, the $y^\ouch_\ell$ are determined by \eqref{ovm10}; they are independent of our choice of $\widehat{x}_J$.

Applying the known case $J=1$ of our algorithm \textsc{Optimize via Matrix}, we compute from the $h_\ell$ a vector of coefficients $\gamma_\ell$ ($1 \leq \ell \leq L$), for which the following holds.

Let $\check{y}_1,\cdots,\check{y}_L$ be real numbers. Set $\displaystyle \check{x} = \sum_{\ell=1}^L \gamma_\ell \check{y}_\ell$.
 Then
 \[\sum_{\ell=1}^L \lvert \check{y}_\ell + h_\ell \check{x} \rvert^p \leq C(p) \sum_{\ell=1}^L \lvert \check{y}_\ell + h_\ell \widehat{x} \rvert^p\]
 for any real number $\widehat{x}$.
 
Taking $\check{y}_\ell = y^\ouch_\ell$ for $\ell=1,\cdots,L$, we learn the following. Let
\begin{equation}\label{ovm13}
\check{x}_J = \sum_{\ell=1}^L \gamma_\ell y^\ouch_\ell
\end{equation}
and then define $\check{x}_1,\cdots,\check{x}_{J-1}$ from $\check{x}_J$ via \eqref{ovm8}, i.e.,
\begin{equation}\label{ovm14}
\check{x}_j = \sum_{\ell=1}^L \widehat{b}_{j \ell} y_\ell + g_j \check{x}_J \quad \mbox{for} \; j=1,\cdots,J-1.
\end{equation}
Then
\begin{equation}\label{ovm15}
\sum_{\ell=1}^L \lvert y_\ell + \sum_{j=1}^J a_{\ell j} \check{x}_j \rvert^p \leq C(p) \sum_{\ell=1}^L \lvert y_\ell + \sum_{j=1}^J a_{\ell j} \widehat{x}_j \rvert^p.
\end{equation}
(See \eqref{ovm12}.)

From \eqref{ovm6} and \eqref{ovm15}, we see that
\begin{equation}\label{ovm16}
\sum_{\ell=1}^L \lvert y_\ell + \sum_{j=1}^J a_{\ell j} \check{x}_j \rvert^p \leq C(p) \sum_{\ell=1}^L \lvert y_\ell + \sum_{j=1}^J a_{\ell j} x_j \rvert^p.
\end{equation}
Here, $\check{x}_1,\cdots,\check{x}_J$ are computed from \eqref{ovm13},\eqref{ovm14}; and $x_1,\cdots,x_J$ are arbitrary.

We produce efficient formulas for the $\check{x}_j$. Putting \eqref{ovm10} into \eqref{ovm13}, we find that
\begin{align*}
\check{x}_J & = \sum_{\ell=1}^L \gamma_\ell \cdot \left\{ y_\ell + \sum_{j=1}^{J-1} a_{\ell j} \sum_{\ell' = 1}^L \widehat{b}_{j \ell'} y_{\ell'} \right\} \\
& = \sum_{\ell=1}^L \gamma_\ell \cdot y_\ell + \sum_{\ell' = 1}^L \sum_{j=1}^{J-1}  \left[ \sum_{\ell=1}^L   \gamma_\ell a_{\ell j}  \right] \widehat{b}_{j \ell'} y_{\ell'} \\
& = \sum_{\ell=1}^L \left\{ \gamma_\ell + \sum_{j=1}^{J-1}  \left[ \sum_{\ell'=1}^L  \gamma_{\ell'} a_{\ell' j}  \right] \widehat{b}_{j \ell}  \right\} \cdot y_\ell.
\end{align*}
Therefore, setting
\begin{equation}\label{ovm17}
\Delta_j = \sum_{\ell=1}^L \gamma_\ell a_{\ell j} \quad \mbox{for} \; j=1,\cdots, J -1
\end{equation}
and
\begin{equation}\label{ovm18}
b^{\#\#}_{J\ell} = \gamma_\ell + \sum_{j=1}^{J-1} \Delta_j \widehat{b}_{j \ell} \quad \mbox{for} \; \ell=1,\cdots,L
\end{equation}
we find that
\begin{equation}\label{ovm19}
\check{x}_J = \sum_{\ell=1}^L b^{\# \#}_{J \ell} y_\ell.
\end{equation}
Substituting \eqref{ovm19} into \eqref{ovm14}, we find that
\[\check{x}_j = \sum_{\ell=1}^L \left\{ \widehat{b}_{j \ell} + g_j b_{J \ell}^{\# \#} \right\} y_\ell \quad \mbox{for} \; j=1,\cdots,J-1.\]
Thus, setting
\begin{equation}\label{ovm20}
b^{\# \#}_{j \ell} = \widehat{b}_{j \ell} + g_j b^{\# \#}_{J \ell} \quad \mbox{for} \; j=1,\cdots,J-1,
\end{equation}
we have
\begin{equation}\label{ovm21}
\check{x}_j = \sum_{\ell=1}^L b^{\# \#}_{j \ell} y_\ell \quad \mbox{for} \; j=1,\cdots,J-1.
\end{equation}
Recalling \eqref{ovm19}, we see that \eqref{ovm21} holds for $j=1,\cdots,J$. Thus, with $\check{x}_1,\cdots,\check{x}_{J}$ defined by \eqref{ovm21}, we have
\[\sum_{\ell=1}^L \lvert y_\ell + \sum_{j=1}^J a_{\ell j} \check{x}_j \rvert^p \leq C(p) \sum_{\ell=1}^L \lvert y_\ell + \sum_{j=1}^J a_{\ell j} x_j \rvert^p\]
for any real numbers $x_1,\cdots,x_J$. (See \eqref{ovm16}.)

So the matrix $(b^{\# \#}_{j \ell})_{\substack{1 \leq j \leq J \\ 1 \leq \ell \leq L}}$ is as promised in our algorithm.

Let us review the computation of $(b^{\# \#}_{j \ell})$.

\label{page60}

\begin{itemize}
\item Recursively, we apply \textsc{Optimize via Matrix} to the arguments $p$, $(a_{\ell j})_{\substack{1 \leq \ell \leq L \\ 1 \leq j \leq J-1}}$; this yields $(\widehat{b}_{j \ell})_{\substack{1 \leq j \leq J-1 \\ 1 \leq \ell \leq L}}$.
\item Next, we compute from \eqref{ovm9} the quantities
\[g_j = \sum_{\ell=1}^L \widehat{b}_{j \ell} a_{\ell J} \quad \mbox{for} \; j=1,\cdots,J-1.\]
\item We then compute from \eqref{ovm11} the numbers
\[h_\ell = a_{\ell J} + \sum_{j=1}^{J-1} a_{\ell j} g_j \quad \mbox{for} \; \ell=1,\cdots,L.\]
\item We apply the case $J=1$ of \textsc{Optimize via Matrix} to the $L \times 1$ matrix $(h_\ell)$, to produce the numbers $\gamma_\ell$ ($\ell=1,\cdots,L$).
\item From \eqref{ovm17} we then compute the numbers
\[\Delta_j  = \sum_{\ell=1}^L \gamma_\ell a_{\ell j} \quad \mbox{for} \; j=1,\cdots,J-1.\]
\item We set 
\[ b^{\# \#}_{J \ell} = \gamma_\ell + \sum_{j=1}^{J-1} \Delta_j \widehat{b}_{j \ell} \quad \mbox{for} \; \ell=1,\cdots,L.
\]
(See \eqref{ovm18}.)
\item Finally, we set
\[b^{\# \#}_{j \ell} = \widehat{b}_{j \ell} + g_j b^{\# \#}_{J \ell} \quad \mbox{for} \; j=1,\cdots,J-1 \; \mbox{and} \; \ell=1,\cdots,L.\]
(See \eqref{ovm20}.)
\end{itemize}
One can now check easily that our algorithm uses work and storage at most $C L$, as promised.

This concludes our explanation of the algorithm \textsc{Optimize via Matrix}.

\end{proof}

\chapter{Statement of the Main Technical Results}\label{sec_mainresults}

Suppose that $E \subset \frac{1}{32}Q^\circ$ is finite, where $Q^\circ = [0,1)^n$. We assume that $ N = \#(E) \geq 2$.

If $\cA \subsetneq \cM$, then let $\cA^- \subset \cM$ denote the maximal subset less than $\cA$. (See Section \ref{sec_multi} for the definition of the order relation $<$ on $2^\cM$.)

For each $\cA \subset \cM$, we will define the following.

\begin{itemize}
\item A decomposition $\CZ(\cA)$ of $Q^\circ$ into dyadic cubes.  We guarantee the following. \label{CZprops}
\begin{description}
\item[(CZ1)] If $Q,Q' \in \CZ(\cA)$ and $Q \leftrightarrow Q'$,  then  $\frac{1}{2} \delta_Q \leq \delta_{Q'} \leq 2 \delta_Q$ (``good geometry'').
\item[(CZ2)] If $Q \in \CZ(\cA)$ and $\delta_Q \leq c_*(\cA)$ then $S(\cA)Q$ is not tagged with $(\cA,\epsilon_1(\cA))$. \\
Moreover, $S(\cA) = 9$ for $\cA = \cM$.
\item[(CZ3)] In the case $\cA \neq \cM$:

If $Q \in \CZ(\cA)$, $Q' \in \CZ(\cA^-)$, $Q' \subset Q$ and $\delta_{Q'} \leq c_*(\cA) \delta_Q$ 

then the cube $3Q$ is tagged with $(\cA,\epsilon_2(\cA))$.
\item[(CZ4)] In the case $\cA = \cM$: 

If $Q \in \CZ(\cM)$, then $3Q$ is tagged with $(\cM,\epsilon_2(\cM))$.

\item[(CZ5)] In the case $\cA \neq \cM$: 

$\CZ(\cA^-)$ refines $\CZ(\cA)$. 

(We do not exclude the possibility that $\CZ(\cA^-) = \CZ(\cA)$.) 

\end{description}
Moreover, $c_*(\cA), \epsilon_1(\cA),\epsilon_2(\cA) \in (0,1)$ and $S(\cA) \geq 1$.

\item A collection $\CZ_{\main}(\cA)$ consisting of all cubes $Q \in \CZ(\cA)$ such that $\frac{65}{64}Q \cap E \neq \emptyset$.
\item For each $Q \in \CZ_{\main}(\cA)$, a list of functionals in short form $\Omega(Q,\cA) \subset \left[ \X(\frac{65}{64}Q \cap E) \right]^*$ (the ``assists'') such that
$$\sum_{Q \in \CZ_{\main}(\cA)} \sum_{\omega \in \Omega(Q,\cA)} \depth(\omega) \leq C N.$$

\item For each $Q \in \CZ_{\main}(\cA)$, a list of functionals $\Xi(Q,\cA) \subset \left[ \X(\frac{65}{64}Q \cap E) \oplus \cP \right]^*$, each having $\Omega(Q,\cA)$-assisted depth at most $C$. We guarantee that
$$ \sum_{Q \in \CZ_{\main}(\cA)} \# \bigl[ \Xi(Q,\cA) \bigr] \leq C N.$$
We set
$$M_{(Q,\cA)}(f,P) = \left( \sum_{\xi \in \Xi(Q,\cA)} \lvert \xi(f,P) \rvert^p \right)^{1/p}.$$
For each $(f,P) \in \X(\frac{65}{64}Q \cap E) \oplus \cP$, we guarantee that
$$c \cdot \|(f,P) \|_{(1+a(\cA))Q} \leq M_{(Q,\cA)}(f,P) \leq C \cdot \|(f,P) \|_{\frac{65}{64}Q}.$$
Here, $ 0 < a(\cA) \leq 1/64$. 

\item For each $Q \in \CZ_{\main}(\cA)$, a linear map $T_{(Q,\cA)} : \X(\frac{65}{64}Q \cap E) \oplus \cP \rightarrow \X$ with the following properties.
\begin{description}
\item[(E1)] $T_{(Q,\cA)}(f,P) = f$ on $(1+a(\cA))Q \cap E$ for each $(f,P)$.
\item[(E2)] $\|T_{(Q,\cA)}(f,P) \|_{\X((1+a(\cA))Q)}^p + \delta_Q^{-mp} \| T_{(Q,\cA)}(f,P) - P \|_{L^p((1+a(\cA))Q)}^p \leq C \left[ M_{(Q,\cA)}(f,P) \right]^p$ for each $(f,P)$.
\item[(E3)] $T_{(Q,\cA)}$ has $\Omega(Q,\cA)$-assisted depth at most $C$.
\end{description}
\item The constants $c_*(\cA), S(\cA),\epsilon_1(\cA),\epsilon_2(\cA),a(\cA),c,C$ depend only on $m,n,p$, and $\cA$. The constant $S(\cA)$ is a positive integer.
\end{itemize}

\begin{remk}\label{mr_rem0}
Note that both $\Xi(Q,\cA)$ and $\Omega(Q,\cA)$ are lists, hence they may contain more than one copy of the same linear functional. In the sums in the third and fourth bullet points, we include separate summands for each occurrence of a given functional $\xi \in \Xi(Q,\cA)$ or $\omega \in \Omega(Q,\cA)$. See Section \ref{sec_not} for more information on our notation concerning lists.
\end{remk}

To compute the objects defined above, we will produce the following algorithms.

\begin{itemize}
\item

\noindent \environmentA{Algorithm: CZ-Oracle.}
We perform one-time work at most $C N \log N$ in space $CN$, after which we can answer queries. A query consists of a point $\underline{x} \in Q^\circ$. The response to the query $\underline{x}$ is the list of all cubes $Q \in \CZ(\cA)$ such that $\underline{x} \in \frac{65}{64}Q$. To answer the query requires work and storage at most $C \log N$.

\item

\noindent  \environmentA{Algorithm: Compute Main-Cubes.}
With work at most $C N \log N$ in space $CN$, we compute the collection of cubes $\CZ_{\main}(\cA)$.

\item

\noindent  \environmentA{Algorithm: Compute Functionals.}
With work at most $C N \log N$ in space $CN$, we compute the following. 
\begin{itemize}
\item For each cube $Q \in \CZ_{\main}(\cA)$,
 the list of functionals $\Omega(Q,\cA)$, with each functional written in short form.
\item For each cube $Q \in \CZ_{\main}(\cA)$,
 the list of functionals $\Xi(Q,\cA)$, with each functional written in short form (in terms of the assists $\Omega(Q,\cA)$).
\end{itemize}

\item 
\noindent  \environmentA{Algorithm: Compute Extension Operator.} 
We perform one-time work at most $C N \log N$ in space $C N$, after which we can answer queries. A query consists of a cube $Q \in \CZ_{\main}(\cA)$ and a point $\underline{x} \in Q^\circ$. The response to the query $(Q,\underline{x})$ is a short form description of the $\Omega(Q,\cA)$-assisted bounded depth linear map
$$(f,P) \in \X\left(\frac{65}{64}Q \cap E \right) \oplus \cP \mapsto J_{\underline{x}}T_{(Q,\cA)}(f,P).$$
To answer the query requires work and storage at most $C \log N$.

\end{itemize}

\chapter{Data Structures}

\section{Algorithms for Dyadic Cubes}\label{algs1}

\subsection{Dyadic Cuboids}\label{sec_dc}

We define a \underline{dyadic interval} to be a subinterval $[a,b) \subset [0,\infty)$ where
$$a=\sum_{\nu=-\infty}^\infty \delta_\nu 2^\nu, \;\; b= \sum_{\nu = -\infty}^\infty \delta_\nu' 2^\nu, \;\;\; \mbox{each} \; \delta_\nu, \delta_\nu' = 0 \; \mbox{or} \; 1,$$
only finitely many $\delta_\nu, \delta_{\nu}'$ are nonzero, and for some $\mu$,
$$\delta_\nu = \delta'_\nu \; \mbox{for} \; \nu > \mu; \; \delta_\mu = 0, \; \delta_\mu' = 1; \;\; \delta_\nu = \delta_\nu' = 0 \; \mbox{for} \; \nu < \mu.$$

The \underline{dyadic cuboids} to be defined in a moment, will be Cartesian products of dyadic intervals. Thus, by definition, a dyadic cuboid $Q \subset \R^n$ will be a subset of $[0,\infty)^n$.

Fix a dimension $n$. A \underline{dyadic cuboid} $Q$ is a Cartesian product of the form
$$[a_1,b_1) \times \cdots \times [a_n,b_n) \subset \R^n,$$
where each $[a_i,b_i)$ is a dyadic interval, and one of the following holds:
\begin{description}
\item[(1)] All the $[a_i,b_i)$ have the same length, 

or for some $j$ ($1 \leq j < n$),
\item[(2)] All the $[a_i,b_i)$ ($1 \leq i \leq j$) have the same length, and each $[a_i,b_i)$ ($j < i \leq n$) has length $\frac{1}{2}(b_1 - a_1)$.
\end{description}

To \underline{bisect} the cuboid $Q$ means the following.

\noindent \textbf{Case 1}: Suppose $Q$ is as in \textbf{(1)}. Then we bisect $[a_n,b_n)$ into two dyadic intervals $I' = [a_n,\frac{a_n + b_n}{2})$ and $I'' = [\frac{a_n + b_n}{2},b_n)$.

To \underline{bisect} $Q$, we express $Q$ as the disjoint union of the two dyadic cuboids
$$Q' = [a_1,b_1) \times \cdots \times [a_{n-1},b_{n-1}) \times I'$$
and
$$Q'' = [a_1,b_1) \times \cdots \times [a_{n-1},b_{n-1}) \times I''$$
We call $Q'$ the \underline{lesser dyadic child} of $Q$,
and we call $Q''$ the \underline{greater dyadic child} of $Q$.

\noindent \textbf{Case 2}: For some $j$ ($1 \leq j < n$), suppose $Q$ is as in \textbf{(2)}.

Then we bisect $[a_j,b_j)$ into two dyadic intervals $I' = [a_j,\frac{a_j + b_j}{2})$ and $I'' = [\frac{a_j + b_j}{2},b_j)$.

To \underline{bisect} $Q$ is to express $Q$ as a disjoint union of the dyadic cuboids
$$Q' = [a_1,b_1) \times \cdots \times [a_{j-1},b_{j-1}) \times I' \times [a_{j+1},b_{j+1}) \times \cdots \times [a_n,b_n)$$
and
$$Q'' = [a_1,b_1) \times \cdots \times [a_{j-1},b_{j-1}) \times I'' \times [a_{j+1},b_{j+1}) \times \cdots \times [a_n,b_n)$$
We call $Q'$ the \underline{lesser dyadic child} of $Q$,
and we call $Q''$ the \underline{greater dyadic child} of $Q$.

To understand dyadic cuboids and their dyadic children, it is convenient to think of base $2$ expansions of real numbers. Let $\DR$ (``dyadic rationals'') be the set of all sums of the form $\displaystyle \sum_{\nu=-\infty}^\infty \delta(\nu) 2^\nu$, where finitely many $\delta(\nu)$ are equal to $1$, and all other $\delta(\nu)$ are equal to zero.

We define a map $\psi : \DR^n \rightarrow \DR $ as follows.

Let $x = (x_1,\cdots,x_n) \in \DR^n$, with each $\displaystyle x_i = \sum_{\nu= - \infty}^\infty \delta_i(\nu) 2^\nu$ as in the definition of $\DR$. Then we define $$\psi(x) = \sum_{\nu = -\infty}^\infty \sum_{i=1}^n \delta_i(\nu) 2^{\nu n + i} \in \DR.$$
Thus, $\psi$ is a $1$-$1$ correspondence between $\DR^n$ and $\DR$.

We define a $1$-$1$ correspondence between dyadic cuboids in $\R^n$ and dyadic intervals in $\R$, by the following rule:

The dyadic cuboid $Q \subset \R^n$ corresponds to the dyadic interval $I \subset \R$ if and only if $\DR \cap I = \psi((\DR)^n \cap Q)$. By thinking about base $2$ expansions of numbers, one sees easily that this is indeed a $1$-$1$ correspondence between dyadic cuboids in $\R^n$ and dyadic intervals in $\R$. Let us denote this $1$-$1$ correspondence by $I = \Psi(Q)$.

Suppose that $Q$ is a dyadic cuboid, with lesser dyadic child $Q'$ and greater dyadic child $Q''$. Then $\Psi(Q')$ and $\Psi(Q'')$ are the two dyadic children of the dyadic interval $\Psi(Q)$; and $\Psi(Q')$ lies to the left of $\Psi(Q'')$. Again, we leave to the reader the verification of this fact.

We now define a binary relation on dyadic cuboids, and another binary relation on dyadic intervals. We will see that these two relations are both order relations, and moreover, the two order relations are equivalent via the $1$-$1$ correspondence $\Psi$.

\noindent\underline{For cuboids:} Let $Q_1,Q_2$ be distinct dyadic cuboids.
\begin{itemize}
\item If $Q_1 \subset Q_2$, then we say that $Q_2 < Q_1$.
\item If $Q_2 \subset Q_1$, then we say that $Q_1 < Q_2$.
\item Suppose $Q_1$ and $Q_2$ are disjoint. Let $Q$ be the least common ancestor of $Q_1$ and $Q_2$ among dyadic cuboids. Let $Q'$, $Q''$ be the lesser and greater dyadic children of $Q$, respectively. Then one of $Q_1$, $Q_2$ is contained in $Q'$, and the other is contained in $Q''$.
\begin{itemize}
\item If $Q_1 \subset Q'$ and $Q_2 \subset Q''$, then we say that $Q_1 < Q_2$.
\item If $Q_2 \subset Q'$ and $Q_1 \subset Q''$, then we say that $Q_2 < Q_1$.
\end{itemize}
\end{itemize}

\noindent\underline{For dyadic intervals:} Let $I_1, I_2$ be distinct dyadic intervals.
\begin{itemize}
\item If $I_1 \subset I_2$, then we say that $I_2 < I_1$.
\item If $I_2 \subset I_1$, then we say that $I_1 < I_2$.
\item If $I_1$ and $I_2$ are disjoint, then let $I$ be the smallest dyadic interval containing $I_1$ and $I_2$. We bisect $I$ into $I'$ and $I''$, with $I'$ lying to the left of $I''$. Then one of $I_1$, $I_2$ is contained in $I'$, and the other is contained in $I''$.
\begin{itemize}
\item If $I_1 \subset I'$ and $I_2 \subset I''$, then we say that $I_1 < I_2$.
\item If $I_2 \subset I'$ and $I_1 \subset I''$, then we say that $I_2 < I_1$.
\end{itemize}
\end{itemize}

Thus, we have defined binary relations $<$ on dyadic cuboids, and on dyadic intervals. It is clear that these two relations correspond to each other via the $1$-$1$ correspondence $\Psi$ between dyadic cuboids and dyadic intervals.

Next, we check that $<$ is an order relation. To see this, it is most convenient to work with dyadic intervals. By examining each case mentioned above, we see that $[a_1,b_1) < [a_2,b_2)$ if and only if either $[a_1 < a_2]$ or $[a_1 = a_2 \;\; \mbox{and} \;\; b_2 < b_1]$.

This makes it obvious that $<$ is an order relation.

We will make use of the following

\begin{prop} \label{p31} \hfill
\begin{enumerate}
\item Let $I_1,I_2$ be dyadic intervals, and suppose $I_1 < I_2$. Then either $I_2 \subset I_1$, or $I_1 \cap I_2 = \emptyset$.
\item Let $I_1,I_2,I_3$ be dyadic intervals, and suppose $I_1 < I_2 < I_3$. If $I_1 \cap I_2 = \emptyset$ then $I_1 \cap I_3 = \emptyset$.
\end{enumerate}
\end{prop}
\begin{proof}
(1) holds simply because we cannot have $I_1 \subset I_2$ when $I_1 < I_2$.

To check (2), let $I_i = [a_i,b_i)$ for $i=1,2,3$. We have $a_1 \leq a_2$, and $[a_1,b_1) \cap [a_2,b_2) = \emptyset$. Hence $a_2 \geq b_1$. Since $I_2 < I_3$, we have also $a_3 \geq a_2$. Therefore, $a_3 \geq b_1$, and thus $[a_1,b_1) \cap  [a_3,b_3) = \emptyset$, as claimed.
\end{proof}

\begin{cor} \label{c31}
Let $I_1 < I_2 < \cdots < I_N$ for dyadic intervals $I_1,\cdots,I_N$ ($N \geq 2$). Then one of the following holds.
\begin{itemize}
\item All of $I_2, \cdots , I_N$ are contained in $I_1$.
\item All of $I_2,\cdots,I_N$ are disjoint from $I_1$.
\item For some $j$ ($2 \leq j < N$), we find that $I_2,\cdots,I_j \subset I_1$ and $I_{j+1},\cdots,I_N$ are disjoint from $I_1$.
\end{itemize}
\end{cor}

\begin{cor} \label{c32}
Let $Q_1,\cdots,Q_N$ be dyadic cuboids ($N \geq 2$), and suppose $Q_1 < Q_2 < \cdots < Q_N$. Then one of the following holds.
\begin{itemize}
\item All of $Q_2, \cdots , Q_N$ are contained in $Q_1$.
\item All of $Q_2,\cdots,Q_N$ are disjoint from $Q_1$.
\item For some $j$ ($2 \leq j <N$), we find that $Q_2,\cdots,Q_j \subset Q_1$ and $Q_{j+1},\cdots,Q_N$ are disjoint from $Q_1$.
\end{itemize}
\end{cor}

We briefly discuss the computer implementation of dyadic cuboids. In our infinite-precision model of computation (see Section \ref{sec_moc1}), a dyadic cuboid $Q$ can be stored using at most $C$ memory locations\footnote{In our finite-precision model of computation, we will deal only with dyadic cuboids whose sides have length between $2^{-CS_0}$ and $2^{+ C S_0}$; see Section \ref{sec_moc2}.}. 
We can compute the lesser and greater children of a given dyadic cuboid. We can determine whether two given dyadic cuboids are disjoint, and if not, then we can determine which one contains the other. We can also decide whether two given cuboids are equal. These operations require work and storage at most $C$. (In this paragraph $C$ denotes a constant depending only on the dimension $n$.) We note also that we can compute the least common ancestor of two given dyadic cuboids using work and storage at most $C$. (Recall that in our model of computation it takes a single operation to compute the smallest dyadic interval containing two given dyadic intervals.) It follows that we can compare two given cuboids under the order relation $<$ using work and storage at most $C$.

\subsection{Preliminary Definitions} \label{pdef_sec}

A \underline{B-Tree} is a rooted finite tree $T$ such that every node has zero, one or two children. We write $\rt(T)$ to denote the root of $T$. 

Let $x \in T$ be a node. Then $\desc(x,T)$ denotes the set of descendants of $x$ in the tree $T$, and $\ndesc(x,T)$ denotes the set $T \setminus \desc(x,T)$. 

If $x= \rt(T)$, then of course $\ndesc(x,T)$ is empty. 

If $x \neq \rt(T)$, then $\ndesc(x,T)$ is again a $B$-Tree, with the same root as $T$. 

In any case, $\desc(x,T)$ is a BTree with root $x$.

(Here, we adopt the convention that each node is a descendant of itself.)

For any BTree $T$, we write $\#(T)$ for the number of nodes in $T$.

When we implement a BTree $T$ in the computer, we store the nodes of $T$, a pointer to the root of $T$, and a pointer from each node of $T$ (except the root) to its parent. Also, we mark each node to indicate whether it is a leaf (recall that a leaf is a node with no children); and we mark each internal node (i.e. each non-leaf) with pointers to each of its children.

A \underline{binary tree} is a BTree such that each node has either zero or two children.

\noindent\underline{\textbf{DTrees and ADTrees}}

Fix $1 < p < \infty$, $n \geq 1$, $D \geq 1$.

A \underline{DTree} is a BTree $T$ each of whose nodes $x$ is identified with a dyadic cuboid $Q_x \subset \R^n$, such that the following hold.
\begin{itemize}
\item Let $y$ be a child of $x$ in $T$. Then $Q_y$ is a proper sub-cuboid of $Q_x$.
\item Let $y,z$ be distinct children of $x$ in $T$. Then $Q_y$ and $Q_z$ are disjoint.
\end{itemize}
An \underline{ADTree} is a DTree $T$ each of whose nodes $x$ is marked with $D$ linear functionals $\mu_1^x,\cdots,\mu_D^x$ on $\R^D$.

\label{pp5}

\noindent(``D'' in ``DTree'' stands for ``dyadic''; \\
``AD'' in ``ADTree'' stands for ``agumented dyadic''.)

\environmentA{Algorithm: BTree1.}

Given a BTree $T$ with $\#(T) \geq 2$, we produce a node $x_\spl \in T$, other than the root of $T$, such that
\[\#[\desc(x_\spl(T), T)] \leq \frac{9}{10} \#(T)\]
and
\[\#[ \ndesc(x_\spl(T),T))] \leq \frac{9}{10} \#(T). \]
The work and storage used to do so are at most $C \cdot \#(T)$ for a universal constant $C$.

\begin{proof}[\underline{Explanation}]

We first mark each node of $T$ with the number of its descendants. We then start with $\widetilde{x} = \rt(T)$. Initially, $\#[\desc(\widetilde{x},T)] = \#(T) > \frac{9}{10} \#(T)$.

\noindent While $\bigl( \#[\desc(\widetilde{x},T)] > \frac{9}{10} \#(T) \bigr)$

$\bigl\{$
/* Note that $\widetilde{x}$ cannot be a leaf of $T$, hence there are one or two children of $\widetilde{x}$. */

\noindent We let $\widetilde{y}$ be a child of $\widetilde{x}$ having as many descendants as possible (among the children of $\widetilde{x}$). We then set $\widetilde{x} := \widetilde{y}$.
$\bigr\}$

The above loop will terminate, since otherwise we would obtain an infinite descending sequence in the finite tree $T$.

When the loop terminates, we have
\begin{equation}
\label{bt0}
\#[\desc(\widetilde{x}, T)] \leq \frac{9}{10} \#(T).
\end{equation}
We will check also that $\widetilde{x}$ is not the root of $T$, and that 
\begin{equation}
\label{bt1}
\#[\ndesc(\widetilde{x}, T)] \leq \frac{9}{10} \#(T).
\end{equation}
Since the work and storage of the above procedures are at most $C \cdot \#(T)$, we can return $x_\spl(T) = \widetilde{x}$, and our algorithm will perform as promised.

Thus, it remains only to check that $\widetilde{x}$ isn't the root of $T$, and that \eqref{bt1} holds.

That $\widetilde{x}$ isn't the root of $T$ follows at once from \eqref{bt0}.

To check \eqref{bt1}, we note that $\widetilde{x}$ arose from its parent $\widetilde{x}^+$ by executing our loop for the last time. We have $\#[\desc(\widetilde{x}^+, T)] > \frac{9}{10} \#(T)$ since we executed the loop to produce $\widetilde{x}$ from $\widetilde{x}^+$.

Also,
\begin{align*}
\#[\desc(\widetilde{x}^+, T)] &= 1 + \sum_{y \; \mbox{\tiny children of} \; \widetilde{x}^+} \#[\desc(y,T)]\\
& \leq 1 + 2 \cdot \max \{ \#[\desc(y,T)] : \; y \; \mbox{children of} \; \widetilde{x}^+ \} \\
& = 1 + 2 \cdot \# [\desc(\widetilde{x},T)] \},
\end{align*}
since $\widetilde{x}^+$ has at most $2$ children, and since $\widetilde{x}$ has at least as many descendants as any child of $\widetilde{x}^+$.

Therefore, $1 + 2 \cdot \# [\desc(\widetilde{x},T)] > \frac{9}{10} \#(T)$, hence $\#[\desc(\widetilde{x},T)] > \frac{9}{20} \#(T) - \frac{1}{2}$. Since $\#(T) \geq 2$, we find that
\begin{align*}
\#[\desc(\widetilde{x},T)] > \frac{2}{20} \#(T) + \left(\frac{7}{20} \#(T) - \frac{1}{2}\right) > \frac{1}{10} \#(T).
\end{align*}
Thus, $\#[\desc(\widetilde{x},T)] > \frac{1}{10} \#(T)$, from which \eqref{bt1} follows at once.

This completes our explanation of the algorithm \textsc{BTree1}.

\end{proof}

\subsection{Control Trees}
\label{control_tree}

Let $T$ be a BTree. If $x,y \in T$, then we write $x \leq y$ if and only if $x$ is a descendant of $y$ in $T$.

Let $\widetilde{T}$ be a BTree. We call $\widetilde{T}$ a \underline{sub-tree} of $T$ if $\widetilde{T}$ consists of nodes in $T$ and if the following condition holds: for any nodes $x \leq y \leq z$ in $T$, if $x,z \in \widetilde{T}$ then $y \in \widetilde{T}$.

A \underline{control tree candidate} for $T$ is a finite binary tree $\T$ (i.e., each internal node of $\T$ has exactly two children), whose nodes are marked as follows.
\begin{itemize}
\item \underline{Let $\xi$ be any node of $\T$}. Then $\xi$ is marked by a pointer to a BTree called $\BT(\xi)$, which is a sub-tree of $T$. We mark $\xi$ with a pointer to a node $x_{\rt}(\xi) \in T$ which is the root of $\BT(\xi)$.
\item \underline{Let $\xi$ be any internal node of $\T$}. Then the two children of $\xi$ in $\T$ are marked separately as $\gchild(\xi)$ and $\schild(\xi)$; and the node $\xi$ is marked by a node $x_\spl(\xi) \in \BT(\xi)$.
\item \underline{Let $\xi$ be any leaf of $\T$}. Then $\xi$ is marked by a node $x_\ind (\xi) \in T$.
\end{itemize}

Let $T$ be a BTree. By induction on $\#(T)$, we define a particular control tree candidate for $T$, called the \underline{control tree} for $T$, to be denoted $\CT(T)$. The inductive definition of $\CT(T)$ proceeds as follows.

\noindent\underline{Base Case}: Suppose $\#(T) = 1$. Thus, $T$ consists of a single node $x_0$. We then take $\CT(T)$ to consist of a single node $\xi_0$, marked with the nodes $x_\ind(\xi_0) = x_0$ and $x_{\rt}(\xi_0) = x_0$, and also marked with a pointer to the BTree $\BT(\xi_0) = T$. Note that $\CT(T)$ is a control tree candidate for $T$. Thus we have defined $\CT(T)$ in the base case.

\noindent\underline{Induction Step}: Suppose $\#(T) \geq 2$, and suppose we have already defined $\CT(T')$ for any BTree $T'$ with fewer nodes than $T$. We then define $\CT(T)$ as follows.

We apply to the BTree $T$ the algorithm \textsc{BTree1}, to produce a node $x_\spl(T) \in T$. We know that $x_\spl(T)$ is not the root of $T$; and that 
\[  \#[\desc(x_\spl(T),T)] \leq \frac{9}{10} \#(T), \; \mbox{and}\]
\[\#[ \ndesc(x_\spl(T), T) ] \leq \frac{9}{10} \#(T).\]

Let $T_\go := \desc(x_\spl(T),T)$ and $T_\stay: = \ndesc(x_\spl(T),T)$. Then $T_\go$, $T_\stay$ are BTrees, with fewer nodes than $T$. By induction hypothesis, we have already defined the control trees $\CT(T_\go)$, $\CT(T_\stay)$. 

We define the tree $\T$ to consist of a root $\xi_0$, together with the two trees $\CT(T_\go)$, $\CT(T_\stay)$, where we take the two children of $\xi_0$ to be the roots of $\CT(T_\go)$ and of $\CT(T_\stay)$. Thus, $\T$ is a binary tree. We mark the nodes of $\T$ to form a control tree candidate for $T$, as follows:
\begin{itemize}
\item We keep the markings of the nodes of $\CT(T_\go)$ and $\CT(T_\stay)$, without change.
\item We mark the root of $\CT(T_\go)$ as $\gchild(\xi_0)$, and we mark the root of $\CT(T_\stay)$ as $\schild(\xi_0)$.
\item We mark the root $\xi_0$ with the node $x_\spl(\xi_0) = x_\spl(T)$.
\item We mark the root $\xi_0$ with a pointer to the tree $T$, i.e., we take $\BT(\xi_0) = T$ and $x_{\rt}(\xi_0) = $ the root of $T$.
\end{itemize}

We define $\CT(T)$ to be the marked tree $\T$. This concludes our inductive definition of $\CT(T)$. To show that $\CT(T)$ is a control-tree candidate for $T$, we just need to establish the following lemma. 

\begin{lem}\label{ggg_lem}
Let $T$ be a BTree. Then $\BT(\xi)$ is a sub-tree of $T$ for each $\xi \in \CT(T)$.
\end{lem}
\begin{proof}

The proof is by induction on $\#(T)$.

If $\#(T) = 1$ then $\BT(\xi) = T$ for the single node $\xi$ in $\CT(T)$. The result is immediate.

Suppose that $\#(T) \geq 2$, and suppose that lemma holds for all trees with fewer nodes than $T$. 

If $\xi$ is the root $\xi_0$ of $\CT(T)$, then $\BT(\xi) = T$ by definition, hence the  conclusion of the lemma is obvious.

If $\xi$ is not the root of $\CT(T)$ then $\xi$ is a node in either $\CT(T_{\go})$ or $\CT(T_{\stay})$.

Thus, by the inductive hypothesis, $\BT(\xi)$ is a sub-tree of either $T_{\go}$ or $T_{\stay}$. Since $T_{\go}$ and $T_{\stay}$ are sub-trees of $T$, we conclude that $\BT(\xi)$ is a sub-tree of $T$.

This concludes the proof of the lemma by induction.
\end{proof}

\begin{lem}\label{fff_lem}
The following properties of $\CT(T)$ hold.
\begin{enumerate}[(A)]
\item The number of nodes of $\CT(T)$ is $2 \#(T) - 1$.
\item Moreover, any descending sequence in $\CT(T)$ has length at most $1 + C \log(\#(T))$ for a universal constant $C$.
\item Finally, 
\[\sum_{\xi \in \CT(T)} \# \left( \BT(\xi) \right) \leq C \cdot \#(T) \cdot \left\{ \log_2(\#(T)) + 1 \right\} \]
for a universal constant $C$.
\end{enumerate}
\end{lem}
\begin{proof}
We prove  (A) by induction on $\#(T)$. If $\#(T) = 1$, then by definition $\#(\CT(T)) = 1$, so (A) holds in this case.

For the induction step, suppose we know (A) for all trees with fewer nodes than $T$ ($\#(T) \geq 2$). We establish (A) for $T$. Indeed, with $x_\spl(T)$ as in the definition of $\CT(T)$, we know that $\CT(T)$ consists of the root $\xi_0$, the control tree $\CT(\desc(x_\spl(T),T))$, and the control tree $\CT(\ndesc(x_\spl(T),T))$. Hence,
\[\# \CT(T) = 1 + \# \CT(\desc(x_\spl(T),T)) + \# \CT(\ndesc(x_\spl(T),T)),\]
whereas
\[ \# T =  \# \desc(x_\spl(T),T) + \# \ndesc(x_\spl(T),T).\]
Since $\# \desc(x_\spl(T),T)$, $\# \ndesc(x_\spl(T),T)$ are strictly less than $\#(T)$, the induction hypothesis gives
\[  \# \CT(\desc(x_\spl(T),T))  = 2 \cdot \# \desc(x_\spl(T),T) - 1\]
and
\[  \# \CT(\ndesc(x_\spl(T),T))  = 2 \cdot \# \ndesc(x_\spl(T),T) - 1.\]

Adding the above, we find that $\# \CT(T) - 1 = 2 \cdot \# T -  2$, proving (A) for the BTree $T$. This completes our induction and proves (A).

To prove (B), we check that $\# BT(\xi') \leq \frac{9}{10} \cdot \# \BT(\xi)$ whenever $\xi'$ is a child of $\xi$ in $\CT(T)$. Indeed, this follows from the definition of $\BT(\xi)$ and the defining property of $x_\spl(T)$ by an obvious induction on $\#(T)$. 

For a descending chain $\xi_0,\xi_1,\xi_2,\cdots,\xi_\ell$ in $\CT(T)$, we therefore have
\[1 \leq \# \BT(\xi_\ell) \leq (9/10)^\ell \# \BT(\xi_0) \leq (9/10)^\ell \cdot \# [T].\]
Thus, $\ell \leq \frac{\log (\# T)}{ \log(10/9)}$, proving (B).

To prove (C), we prove by induction on $\#(T)$ that
\[
(*) \qquad \sum_{\xi \in \CT(T)} \# \BT(\xi) \leq \# T \cdot \bigl\{ \log_{10/9}(\# T) + 1 \bigr\}.
\]
For $\#(T) =1$, this holds because $\BT(\xi_0) = T$ where $\xi_0$ is the one and only node of $\CT(T)$.

Assume $(*)$ holds for all BTrees with fewer nodes than $T$, where $T$ is a given BTree with $\#(T) \geq 2$. Then 
\begin{align*} 
(+) \quad & \sum_{\xi \in \CT(T)}  \# \BT(\xi) = \# \BT(\rt(\CT(T)))  + \sum_{\xi \in \CT(\desc(x_\spl(T), T))} \# \BT(\xi)  \\
& \qquad + \sum_{\xi \in \CT(\ndesc(x_\spl(T),T) )} \# \BT(\xi) \\
&\leq \#(T) + \# \desc(x_\spl(T),T) \cdot \left\{ 1 + \log_{10/9} \left[  \# \desc(x_\spl(T),T) \right] \right\} \\
& \qquad +  \# \ndesc(x_\spl(T),T) \cdot \left\{ 1 + \log_{10/9} \left[ \# \ndesc(x_\spl(T),T) \right] \right\} 
\end{align*}
by induction hypothesis.

We know that
\[1 + \log_{10/9} \left[ \# \desc(x_\spl(T),T) \right] \leq \log_{10/9} \# T\]
and that
\[1 + \log_{10/9} \left[ \# \ndesc(x_\spl(T),T)  \right] \leq \log_{10/9} \# T.\]
Hence, $(+)$ yields the estimate
\begin{align*} 
\sum_{\xi \in \CT(T)} \# \BT(\xi) &\leq \#T + \log_{10/9} (\#T) \cdot \# \desc(x_\spl(T),T) \\
& \qquad +  \log_{10/9} (\#T) \cdot \# \ndesc(x_\spl(T),T),
\end{align*}
thus proving $(*)$.

The proof of our lemma is complete.
\end{proof}

\environmentA{Algorithm: Make Control Tree (Deluxe edition).} 

Given a BTree $T$, we produce the control tree $\CT(T)$. The work and storage used to do so are at most $C \cdot \#(T) \cdot (1 + \log \#(T))$ for a universal constant $C$.

\begin{proof}[\underline{Explanation}]

We simply follow the definition in the obvious way. Where the definition proceeds by induction, the algorithm calls itself recursively.

The assertion about the work and storage follows from assertion (C) of Lemma \ref{fff_lem}, and also the bound on the running time of the algorithm \textsc{BTree1}, which is used as a sub-routine.

\end{proof}

We will not use the Deluxe edition explained above, because it uses too much storage.

\environmentA{Algorithm: Make Control Tree (Paperback edition).} 

Given a BTree $T$, we produce the tree $\CT(T)$ with all its markings \underline{except} for the BTrees $\BT(\xi)$ $(\xi \in \CT(T))$. For each $\xi \in \CT(T)$ we indicate whether $\BT(\xi)$ is a singleton.

The work used to do so is at most $C \cdot \#(T) \cdot (1 + \log \#(T))$, and the storage used is at most $C \cdot \#(T)$. Here, $C$ is a universal constant.

\begin{proof}[\underline{Explanation}]

We proceed as in the deluxe edition of the algorithm \textsc{Make Control Tree}, except that we delete $T$ when we are finished using it.

We spell out the details.

\noindent \underline{If $\#(T) = 1$}, then we take $\CT(T)$ to consist of a single node $\xi_0$, marked with $x_\ind(\xi_0) = x_{\rt}(\xi_0) = $ the one and only node of $T$. We indicate that the BTree $\BT(\xi_0)$ is a singleton.

\noindent \underline{If $\#(T) > 1$}, then we execute the algorithm \textsc{BTree1} to produce the node $x_\spl(T)$.

/* In a later variant of this algorithm, we insert code here */

We compute the trees $T' = \desc(x_\spl(T),T)$ with root $x_\spl(T)$, \\
and $T'' = \ndesc(x_\spl(T),T)$ with root $= \rt(T)$.

To produce the trees $T'$, $T''$ efficiently, we can simply erase the marking indicating $x_{\spl}(T)$ as a child of its parent in $T$, and then produce pointers to the roots of $T'$, $T''$. This destroys the tree $T$ after we no longer need it.

Recursively, we apply the paperback edition of \textsc{Make Control Tree} to $T'$ and $T''$. Thus, we obtain $\CT(T')$ and $\CT(T'')$ with all their markings, except for the markings $\BT(\xi')$ ($\xi' \in \CT(T')$) and $\BT(\xi'')$ ($\xi'' \in \CT(T'')$). These latter markings have not been computed (or rather, they were computed and then deleted). 

The tree $\CT(T)$ then consists of the two trees $\CT(T')$ and $\CT(T'')$, together with a root $\xi_0$. The children of $\xi_0$ are the roots of the two trees $\CT(T')$, $\CT(T'')$. We mark the root of $\CT(T')$ as $\gchild(\xi_0)$, and we mark the root of $\CT(T'')$ as $\schild(\xi_0)$. Also, we mark the root $\xi_0$ of $\CT(T)$ with the node $x_\spl(\xi_0) = x_\spl(T) \in T$.

We mark the root $\xi_0$ of $\CT(T)$ with the node $x_{\rt}(\xi_0) = $ the root of $T$. (Recall that $\BT(\xi_0) = T$.) We indicate that the BTree $\BT(\xi_0)$ is not a singleton. We do not mark the root $\xi_0$ with anything else to tell us what the tree $T$ was before we destroyed it.

This completes our description of the algorithm.

Let us check how much time and space are used.

Let $\Time(T)$ be the number of operations needed to execute the paperback algorithm for the tree $T$. Recalling that the algorithm \textsc{BTree1} uses work $C \#(T)$ to produce $x_\spl(T)$, we see that
\[\Time(T) \leq C \#(T) + \Time(T') + \Time(T''),\]
and we recall that $\#(T'), \#(T'') \leq \frac{9}{10} \#(T)$ and that $\#(T') + \#(T'') = \#(T)$. Hence, it follows by induction on $\#(T)$ that
\[\Time(T) \leq C \#(T) \cdot [1 + \log_{10/9} \#T].\]
Thus, the work required to execute our paperback algorithm is as promised.

Next, we study the storage used by our paperback algorithm, which we denote by $\Space(T)$.

Since we erase $T$, we see easily that
$$\Space(T) \leq \max \{ C \cdot \#(T), \Space(T') + \Space(T'') + C \},$$
i.e.,
\[[\Space(T) + C] \leq \max\{ C' \#(T), [\Space(T') + C] + [\Space(T'') + C].\]
Since $\#(T) = \#(T') + \#(T'')$, it follows by induction on $\#(T)$ that
\[[\Space(T) + C ] \leq C'' \#(T).\]
This proves that the storage used by our paperback algorithm is as promised.

\end{proof}

\label{pp6}

\environmentA{Algorithm: Make Control Tree (Hybrid version).} 

Given an ADTree $T$, with each node $x \in T$ marked by functionals $\mu_1^x,\cdots,\mu^x_D : \R^D \rightarrow \R$, we produce the control tree $\CT(T)$ with all its markings \underline{except} for the trees $\BT(\xi)$ ($\xi \in \CT(T)$). For each node $\xi \in \CT(T)$, we produce functionals $\mu^\xi_1,\cdots,\mu^\xi_D : \R^D \rightarrow \R$ such that
\[\sum_{x \in \BT(\xi)} \sum_{i=1}^D \lvert \mu^x_i(v) \rvert^p \;\;\; \mbox{and} \;\;\; \sum_{i=1}^D \lvert \mu^\xi_i(v) \rvert^p\]
differ by at most a factor $C(D,p)$ for any $v \in \R^D$. (This makes sense because $\BT(\xi)$ is a sub-tree of $T$ for each $\xi$.)

We mark each node $\xi \in \CT(T)$ with such functionals $\mu^\xi_1,\cdots,\mu^\xi_D$.

We mark each node $\xi \in \CT(T)$ to indicate whether $\BT(\xi)$ is a singleton.

The work and storage needed to execute this algorithm are at most $C \#(T) \cdot [ \log(\# T) + 1 ]$ and $C \#(T)$, respectively.

\begin{proof}[\underline{Explanation}]

We proceed as in the explanation of the paperback edition of \textsc{Make Control Tree}, with the following changes.

\begin{itemize}
\item If $\#(T) = 1$, then the tree $\CT(T)$ contains only the root node $\xi_0$. We set $\mu^{\xi_0}_i = \mu^x_i$ for $i=1,\cdots,D$, for the one and only one node $x \in T$.
\item  If $\#(T) \geq 2$, then we proceed as follows. Where we wrote \\
\noindent /* In a later variant of this algorithm, we insert code here */ \\
we now insert a call to  \textsc{Compress Norms} (Section \ref{sec_lf}). Thus, with work and storage at most $C \#(T)$, we produce functionals $\mu^*_1,\cdots,\mu^*_D : \R^D \rightarrow \R$ such that
\[
\sum_{x \in T} \sum_{i=1}^D \lvert \mu^x_i(v) \rvert^p \;\;\; \mbox{and} \;\;\; \sum_{i=1}^D \lvert \mu^*_i(v) \rvert^p
\]
differ by at most a factor $C(D,p)$ for any $v \in \R^D$.

Instead of recursively applying the paperback edition of \textsc{Make Control Tree}, we now recursively apply the hybrid version.

Just before the sentence ``This completes the description of the algorithm'' \\
we set $\mu^{\xi_0}_i = \mu_i^*$ for $i=1,\cdots,D$.

Since $\BT(\xi_0) = T$, our $\mu^{\xi_0}_i$ behave as we ask.

The work and storage needed to execute this algorithm are as promised.
\end{itemize}

\end{proof}

\subsection{Encapsulations}
\label{encaps}

Let $T$ be a DTree, and let $\CT(T)$ be its control tree. Recall that each node $\xi \in \CT(T)$ is marked with a BTree $\BT(\xi)$ consisting of nodes of $T$. Also, each node $x \in T$ is marked with a dyadic cuboid $Q_x$.

Let $Q$ be a dyadic cuboid. An \underline{encapsulation} of $Q$ is a set $S$ of nodes of $\CT(T)$, such that $\{x \in  T : Q_x \subset Q \}$ is the disjoint union of the sets $\BT(\xi)$ as $\xi$ varies over $S$.

\environmentA{Algorithm: Encapsulate.} 

Let $T$ be a DTree with $N$ nodes. After $C N ( 1 + \log N)$ one-time work in space $C N$, we can answer queries as follows:

A query consists of a dyadic cuboid $Q$.

The response to a query $Q$ is an encapsulation $S$ of $Q$, consisting of at most $C + C \log N $ nodes of $\CT(T)$.

The work and storage used to answer a query are at most $C+ C \log N$. Here, $C$ denotes a constant depending only on the dimension $n$.

\begin{proof}[\underline{Explanation}]

Suppose $T$ isn't a singleton. Let $x_\spl(T)$ be the node produced by  the algorithm \textsc{BTree1}. Write $T' = \desc(x_\spl(T),T)$ and $T'' = \ndesc(x_\spl(T),T)$.

We ask: For which nodes $x \in T$ do we have $Q_x \subset Q$ ? To answer this question, we compare the dyadic cuboids $Q$ and $Q_{x_\spl(T)}$. There are three cases:

\noindent\underline{Case 1}: $Q \subset Q_{x_\spl(T)}$.

In this case, we never have $Q_x \subset Q$ for an $x \in T''$.

Hence, in this case, $\{ x \in T : Q_x \subset Q\} = \{ x \in T' : Q_x \subset Q\}$. Thus, we have reduced matters from $T$ to $T'$.

\noindent\underline{Case 2}: $Q_{x_\spl(T)} \subsetneq Q$.

In this case, all $x \in T'$ satisfy $Q_x \subset Q$.

Therefore, in this case, 
\begin{align*}
\{x \in T : Q_x \subset Q\} &= \{x \in T'' : Q_x \subset Q\} \cup T' \\
& = \{ x \in T'' : Q_x \subset Q\} \cup \BT(\gchild(\mbox{root of } \CT(T))).
\end{align*}
Thus, we have reduced matters from $T$ to $T''$.

\noindent\underline{Case 3}: $Q_{x_\spl(T)} \cap Q = \emptyset$.

In this case, no $x \in T'$ satisfy $Q_x \subset Q$, hence $\{x \in T : Q_x \subset Q \} = \{ x \in T'' : Q_x \subset Q \}$. Again, we have reduced matters from $T$ to $T''$.

Cases 1,2,3 are the only possibilities, since $Q$ and $Q_{x_\spl(T)}$ are dyadic cuboids.

Thanks to the above remarks, the following procedure produces an encapsulation, when applied to $\widehat{\xi} = \rt(\CT(T))$.

\begin{itemize}
\item One-time work: Paperback edition of \textsc{Make Control Tree}.
\item Procedure Encap$(Q,\widehat{\xi})$: \\
\noindent /* Produces an encapsulation $S$ of $Q$ for the BTree $\BT(\widehat{\xi})$. */ 

\begin{itemize}
\item If $\widehat{\xi}$ is a leaf of $\CT(T)$, then Encap$(Q,\widehat{\xi})$ returns $\{\widehat{\xi}\}$ if $Q_{x_\ind(\widehat{\xi})} \subset Q$, and returns $\emptyset$ otherwise.

\item If $\widehat{\xi}$ is an internal node of $\CT(T)$, then let $\widehat{x} = x_{\spl}(\widehat{\xi})$, $\xi' = \gchild(\widehat{\xi})$, $\xi'' = \schild(\widehat{\xi})$.

\begin{itemize}
\item \underline{If $Q \subset Q_{\widehat{x}}$}, then return the set produced by (recursively) executing Encap$(Q,\xi')$.

\item \underline{If $Q_{\widehat{x}} \subsetneq Q$}, then return the union of $\{\xi'\}$ with the set produced by (recursively) executing Encap$(Q,\xi'')$.

\item \underline{If $Q_{\widehat{x}}\cap Q = \emptyset$}, then return the set produced by (recursively) executing Encap$(Q,\xi'')$.
\end{itemize}
\end{itemize}
\end{itemize}

Note that the one-time work here is simply that of the paperback edition of \textsc{Make Control Tree}; hence, we perform one-time work at most $C N (1+ \log N)$ in space $CN$.

Regarding the query work, note that the ``depth'' of the recursion (i.e., the number of recursive calls to Encap$(\cdot,\cdot)$) is at most $1+ C \log N$, since $\# ( \BT(\xi))$ decreases by at least a factor of $\frac{9}{10}$ each time we pass from a node $\xi$ to $\gchild(\xi)$ or $\schild(\xi)$.

Therefore, the work and storage used to answer a query are at most $C+ C \log N$. In particular, $\#(S) \leq C+ C \log N$, since it takes work at most $C+C \log N$ to write down $S$.

This completes our explanation of the algorithm \textsc{Encapsulate}.

\end{proof}

\label{pp7}

\environmentA{Algorithm: ADProcess.} 

Given an ADTree $T$, (recall that each node $x \in T$ is marked with a dyadic cuboid $Q_x$ and with linear functionals $\mu^x_1,\cdots,\mu^x_D : \R^D \rightarrow \R$) with $N$ nodes, we perform one-time work at most $C N ( 1 + \log N)$ in space $CN$, after which we can answer queries as follows:

A query consists of a dyadic cuboid $Q$.

The response to a query consists of linear functionals $\mu^Q_1,\cdots,\mu^Q_D  : \R^D \rightarrow \R$ such that
\[\sum_{i = 1}^D \lvert \mu^Q_i(v) \rvert^p \;\;\; \mbox{and} \;\; \sum_{x \in T, \; Q_x \subset Q} \sum_{i=1}^D \lvert \mu^x_i(v) \rvert^p\]
differ by at most a factor $C$ for any $v \in \R^D$.

The work and storage needed to respond to a query are at most $C \cdot (\log N + 1 )$.

Here, $C$ depends only on $p,n,D$ ($n=$ dimension of the cuboids).

\begin{proof}[\underline{Explanation}]

We perform the one-time work for the algorithm  \textsc{Encapsulate}, and we perform the hybrid version of the algorithm  \textsc{Make Control Tree}. Thus, we produce the control tree $\CT(T)$; each node $\xi$ is marked with linear functionals $\mu^\xi_1,\cdots,\mu^\xi_D : \R^D \rightarrow \R$ such that
\[\sum_{i = 1}^D \lvert \mu^\xi_i(v) \rvert^p \;\;\; \mbox{and} \;\; \sum_{x \in \BT(\xi)} \sum_{i=1}^D \lvert \mu^x_i(v) \rvert^p\]
differ by at most a factor of $C$ for any $v \in \R^D$.

Moreover, thanks to the query algorithm in  \textsc{Encapsulate}, we can answer queries as follows.

A query consists of a dyadic cuboid $Q$. The response to a query $Q$ consists of a set $S$ of at most $C + C \log N$ nodes in $\CT(T)$ such that $\{x \in T : Q_x \subset Q \}$ is the disjoint union over $\xi \in S$ of $\BT(\xi) \subset T$. \label{page61}

Given $v \in \R^D$, we have therefore
\[\sum_{x \in T : Q_x \subset Q} \sum_{i=1}^D \lvert \mu^x_i(v) \rvert^p  = \sum_{\xi \in S} \left\{ \sum_{x \in \BT(\xi)} \sum_{i=1}^D \lvert \mu^x_i(v) \rvert^p \right\}\]
which differs by at most a factor of $C$ from
\[(*) \qquad \sum_{\xi \in S} \sum_{i=1}^D\lvert \mu^\xi_i (v) \rvert^p.\]
Applying the algorithm \textsc{Compress Norms} (Section \ref{sec_lf}) to the linear functionals $\mu_i^\xi$ ($\xi \in S$; $i = 1 ,\cdots, D$), we obtain linear functionals $\mu_1^Q , \cdots, \mu_D^Q$, such that for any $v \in \R^D$, the quantity $(*)$ differs by at most a factor of $C$ from $\displaystyle \sum_{i=1}^D \lvert \mu_i^Q(v) \rvert^p$ . Thus, $\mu^Q_1,\cdots,\mu^Q_D$ satisfy the desired condition.

The one-time work of the above algorithm is at most $C N (1 + \log N)$, in space $CN$, thanks to our estimates of the one-time work of \textsc{Encapsulate}, and the work and storage of the hybrid version of \textsc{Make Control Tree}.

The query work of our algorithm is at most that of the algorithm \textsc{Encapsulate}, together with the work of applying \textsc{Compress Norms} to the $\mu^\xi_i$ ($\xi \in S$; $i=1,\cdots,D$).

Since $\#(S) \leq C  (\log N + 1 )$, the work of applying \textsc{Compress Norms} is at most $C  (\log N + 1)$. The query work of \textsc{Encapsulate} is also at most $C (\log N + 1)$.

Therefore, the total query work of our present algorithm is as promised.

This completes our explanation of the algorithm \textsc{ADProcess}.

\end{proof}

\subsection{Making a tree from a list of cuboids.}
\label{maketree_sec}

Fix a dimension $n$, and let $Q_1,\cdots,Q_N$ be a sequence of distinct dyadic cuboids in $\R^n$.

We assume that
\begin{equation} \label{list1}
Q_1 < Q_2 < \cdots < Q_N.
\end{equation}
(See Section \ref{sec_dc} for the definition of the order relation $<$.)

\environmentA{Algorithm: Make Forest.}

Given $ 1 \leq i_{\st} \leq i_{\en} \leq N$, we produce the following:
\begin{itemize}
\item A sorted list $i_1 < i_2 < \cdots < i_L$ consisting precisely of all the $i$ such that $i_\st \leq i \leq i_\en$, and such that there exists no $i'$ with $i_\st \leq i' \leq i_\en$ and $Q_i \subsetneq Q_{i'}$. The list is computed as a linked list: We do not compute an array $i_1,i_2,\cdots,i_L$. Rather, we compute the initial entry $i_1$, and we mark each $i_\nu$ with a pointer to its successor $i_{\nu+1}$ ($1 \leq \nu < L$). The last entry $i_L$ is marked with a $\NULL$ pointer.
 \item For each $i$ ($i_\st \leq i \leq i_\en$) we produce a pointer which is $\NULL$ if $i$ appears in the list $i_1,\cdots,i_L$, and otherwise indicates $i''$ such that among all $i'$ such that $i_\st \leq i' \leq i_\en$ and $Q_i \subsetneq Q_{i'}$, the cuboid $Q_{i''}$ is the smallest with respect to inclusion.
\end{itemize}

To do so requires at most $C \cdot (i_\en - i_\st + 1) \cdot (1 + \log N)$ work  and at most $C \cdot ( i_\en - i_\st + 1)$ storage, aside from that used to hold the list \eqref{list1}. Here, $C$ depends only on the dimension $n$.

\begin{proof}[\underline{Explanation}]

We proceed recursively, by induction on $i_\en - i_\st$.

\noindent\underline{In the base case}: $i_\en = i_\st$.

Then our task is trivial, and it takes work and storage at most $C$.

\noindent\underline{The induction step}: Suppose $i_\en > i_\st$.

Then, by Corollary \ref{c32}, we have one of the following cases.

\noindent\underline{Case 1}: $Q_{i_\st + 1}, \cdots,Q_{i_\en} \subset Q_{i_\st}$.

\noindent\underline{Case 2}: $Q_{i_\st + 1}, \cdots, Q_{i_\en}$ are all disjoint from $Q_{i_\st}$.

\noindent\underline{Case 3}: For some $j$ ($i_\st + 1 \leq j < i_\en$), we have $Q_{i_\st+1},\cdots,Q_j \subset Q_{i_\st}$, and $Q_{j+1},\cdots,Q_{i_\en}$ are disjoint from $Q_{i_\st}$.

We can determine which of these cases holds, simply by checking $Q_{i_\st+1} \cap Q_{i_\st}$ and $Q_{i_\en} \cap Q_{i_\st}$. Moreover, if Case 3 holds, then we can find $j$ by doing a binary search. This requires work at most $C \cdot (1 + \log N)$ and storage at most $C$, aside from the storage used to hold the given cuboids \eqref{list1}.

We describe how to proceed in Case 3. Later, we explain the modifications needed for Cases 1 and 2.

Suppose we are in Case 3, with $j$ known. Recursively, we apply the algorithm \textsc{Make Forest} to indices $i_\st+1$ and $j$ (in place of $i_\st$ and $i_\en$) and also to indices $j+1$ and $i_\en$  (in place of $i_\st$ and $i_\en$).

Thus we obtain the following:

\begin{align}
\label{maketree2}
&\widetilde{i}_1 < \widetilde{i}_2 < \cdots < \widetilde{i}_{\widetilde{L}}, \; \mbox{a linked list of all} \; \widetilde{i} \in \{i_\st+1,\cdots,j\}\\
&\qquad \mbox{such that} \; Q_{\widetilde{i}} \; \mbox{is maximal (under inclusion) among} \; Q_{i_\st+1},\cdots,Q_j. \notag{}
\end{align}

\begin{align}
\label{maketree3}
& \mbox{For each} \; \widetilde{i} \in \{i_\st+1,\cdots,j \}, \; \mbox{either a} \; \NULL \; \mbox{pointer indicating that}\\
& \qquad Q_{\widetilde{i}} \; \mbox{is maximal as in \eqref{maketree2}, or else a pointer to the} \; \widetilde{i}_+\in \{i_\st+1,\cdots,j\}\notag{}\\
& \qquad \mbox{such that} \; Q_{\widetilde{i}_+} \supsetneq Q_{\widetilde{i}} \; \mbox{with} \; Q_{\widetilde{i}_+} \; \mbox{as small as possible under inclusion}. \notag{}
\end{align}

\begin{align}
\label{maketree4}
& i_1^\# < i_2^\# < \cdots < i_{L^\#}^\#, \mbox{a linked list of all} \; i^\# \in \{j+1,\cdots,i_\en\}\\
&\qquad \mbox{such that} \; Q_{i^\#} \; \mbox{is maximal (under inclusion) among} \; Q_{j+1},\cdots,Q_{i_\en}. \notag{}
\end{align}

\begin{align}
\label{maketree5}
& \mbox{For each} \; i^\# \in \{j+1,\cdots,i_{\en}\}, \; \mbox{either a} \; \NULL \; \mbox{pointer indicating that}\\
& \qquad i^\# \; \mbox{appears in the list \eqref{maketree4}, or else a pointer to the} \; i^\#_+ \in \{j+1,\cdots,i_\en\}\notag{}\\
& \qquad \mbox{such that} \; Q_{i^\#_+} \supsetneq Q_{i^\#} \; \mbox{with} \; Q_{i^\#_+} \; \mbox{as small as possible under inclusion}. \notag{}
\end{align}

We now produce the desired output for $i_\st$, $i_\en$.

\begin{itemize}
\item Our linked list $i_1 < i_2 < \cdots < i_L$ consists of the list $i^\#_1 < i_2^\# < \cdots < i_{L^\#}^\#$, with $i_0^\# := i_\st$ added to the beginning of the list. (Since the lists in question are implemented here as linked lists, it takes work at most $C$ to add $i^\#_0$ to the list.)
\item Our pointers are as follows.

For $i \in \{j+1,\cdots,i_\en\}$, the pointers are precisely those produced in \eqref{maketree5}.

For $i \in \{i_\st+1,\cdots,j\}$, we take the pointers produced in \eqref{maketree3}. However, for each $\widetilde{i}_\ell$ in the linked list \eqref{maketree2}, we set the pointer associated to $\widetilde{i}_\ell$ (which was $\NULL$ in \eqref{maketree3}) so that it indicates $i_\st$.

For $i=i_\st$, we take a $\NULL$ pointer.
\end{itemize}

The work needed to implement the above bullet points is at most 
\[C + C \cdot \#\{ \mbox{Maximal cuboids under inclusion among} \; Q_{i_\st+1}, \cdots, Q_{j} \}.\]

This concludes our description of the algorithm in Case 3. It produces the desired information thanks to the inclusions and disjointness conditions that hold in Case 3.

Cases 1 and 2 are similar to Case 3, but easier.

In Case 1, there are no $\{j+1,\cdots,i_\en\}$ to deal with, so we omit all steps relevant to $\{j+1,\cdots,i_\en\}$.

Similarly, in Case 2, there are no $\{i_\st+1,\cdots,j\}$ to deal with, so we omit all steps relevant to $\{i_\st+1,\cdots,j\}$.

Thus, in all cases, our recursive algorithm works as promised, except that we have not yet estimated the work and storage needed to carry it out.

Regarding the work, which we call $\Work(i_\st,i_\en)$, we note that (in \underline{Case 3}) we have
\begin{align}
\label{maketree6}
\Work(i_\st,i_\en) \leq &  \; C \cdot (1+ \log N) + \Work(i_\st+1,j)  \\
 \notag{} & \qquad\qquad\qquad\;\; + \Work(j+1,i_\en)  \\
& + C  \cdot \# \bigl\{ \mbox{Maximal} \; Q_i \; \mbox{(under inclusion) among} \; Q_{i_\st+1},\cdots,Q_j \bigr\}. \notag{}
\end{align}
In \underline{Case 1} we have instead
\begin{align}
\label{maketree7}
\Work(i_\st,i_\en) \leq & \; C \cdot (1+ \log N) + \Work(i_\st+1,i_\en) \\
& + C  \cdot \# \bigl\{ \mbox{Maximal} \; Q_i \; \mbox{(under inclusion) among} \; Q_{i_\st+1},\cdots,Q_{i_\en} \bigr\}. \notag{}
\end{align}
In \underline{Case 2} we have
\begin{align}
\label{maketree8}
\Work(i_\st,i_\en) \leq & \; C \cdot (1+ \log N) + \Work(i_\st+1,i_\en)
\end{align}
since there are no pointers to modify in Case 2.

To analyze \eqref{maketree6},\eqref{maketree7},\eqref{maketree8}, we introduce $\NR(i_\st,i_\en) :=$ the number of cuboids $Q_i$ among $Q_{i_\st},\cdots,Q_{i_\en}$ that are not contained in any other $Q_{i'}$ among $Q_{i_\st},\cdots,Q_{i_\en}$. (``$\NR$'' stands for ``Number of Roots''.)

Then \eqref{maketree6},\eqref{maketree7},\eqref{maketree8} together with the known inclusions that follow from Corollary \ref{c32} tell us the following. 

In \underline{Case 3}, with $j$ as given in Case 3, we have
\begin{align}
\label{mt9}
\Work(i_\st,i_\en) \leq & \; C \cdot (1+ \log N) + \Work(i_\st+1,j) \\
&  + \Work(j+1,i_\en) + C \cdot \NR(i_\st+1, j) \notag{}
\end{align}
and
\begin{equation}
\label{mt10}
\NR(i_\st,i_\en) = 1 + \NR(j+1,i_\en).
\end{equation}

In \underline{Case 1}, we have instead
\begin{equation}
\label{mt11}
\Work(i_\st,i_\en) \leq C \cdot (1+ \log N) + \Work(i_\st+1,i_\en) + C \cdot \NR(i_\st+1, i_\en)
\end{equation}
and
\begin{equation}
\label{mt12}
\NR(i_\st,i_\en) = 1.
\end{equation}

In \underline{Case 2}, we have
\begin{equation}
\label{mt13}
\Work(i_\st,i_\en) \leq C \cdot (1+ \log N) + \Work(i_\st+1,i_\en)
\end{equation}
and
\begin{equation}
\label{mt14}
\NR(i_\st,i_\en) = 1 + \NR(i_\st+1,i_\en).
\end{equation}

Consequently, in \underline{Case 3} we have
\begin{align}
\label{mt15}
\bigl[ \Work(i_\st,i_\en)  + \underline{C}'  \NR(i_\st,i_\en) \bigr] \leq & \; \underline{C}'' \cdot (1 + \log N)  \\
& + \bigl[ \Work(i_\st + 1,j) + \underline{C}' \NR(i_\st + 1, j) \bigr] \notag{} \\
& + \bigl[ \Work(j+1,i_\en) + \underline{C}' \NR(j+1,i_\en)\bigr].\notag{}
\end{align}
(Here, we pick $\underline{C}'$ big, and then pick $\underline{C}''$ much bigger.)

In \underline{Cases 1 \& 2}, we have instead
\begin{align}
\label{mt16}
\bigl[ \Work(i_\st,i_\en) + \underline{C}' \NR(i_\st,i_\en) \bigr] \leq  & \; \underline{C}'' \cdot (1+ \log N)  \\
& + \bigl[ \Work(i_\st + 1, i_\en ) + \underline{C}' \NR(i_\st+1,i_\en) \bigr].\notag{}
\end{align}
Also, when $i_\st = i_\en$ we have
\begin{equation}
\label{mt18}
\Work(i_\st,i_\en) + \underline{C}' \NR(i_\st,i_\en) \leq \underline{C}''.
\end{equation}
Thanks to \eqref{mt15}$\cdots$\eqref{mt18}, induction on $i_\en - i_\st$ yields the estimate
\[\Work(i_\st,i_\en) + \underline{C}' \NR(i_\st,i_\en) \leq \underline{C}''' \cdot (1+ \log N ) \cdot (i_\en - i_\st + 1).\]
In particular,
\[\Work(i_\st,i_\en)  \leq C \cdot (1+ \log N ) \cdot (i_\en - i_\st + 1)\]
as promised in the statement of \textsc{Make Forest}.

Next, we analyze the storage needed to execute \textsc{Make Forest}.

We implement the pointers as global data (see the second bullet point in the statement of the algorithm); this requires space at most $C \cdot (i_\en - i_\st + 1)$. Let $\Space(i_\st,i_\en)$ be the space in which we can execute our algorithm \textsc{Make Forest}$(i_\st,i_\en)$, not counting the space needed for the pointers, but including the space needed to compute and display the linked list; see the first bullet point in the statement of the algorithm.

In \underline{Case 3}, we see that
\[\Space(i_\st,i_\en) \leq C + \Space(i_\st+1,j) + \Space(j,i_\en).\]

In \underline{Case 1} and in \underline{Case 2} we have instead
\[\Space(i_\st,i_\en) \leq C + \Space(i_\st+1,i_\en).\]

Since also $\Space(i_\st,i_\st) \leq C$, it therefore follows by induction that $\Space(i_\st,i_\en) \leq C (i_\en - i_\st + 1)$.

Therefore, the storage used in executing \textsc{Make Forest}$(i_\st,i_\en)$ is as promised.

This completes our explanation of that algorithm.

\end{proof}

\environmentA{Algorithm: Fill in gaps.} 

Suppose we are given a cuboid $\hQ$ and a list of pairwise disjoint cuboids $Q_{i_\st}, \cdots, Q_{i_\en}$, sorted so that $Q_{i_\st} < Q_{i_\st + 1} < \cdots < Q_{i_\en}$. Assume that each $Q_i \subset \hQ$.

We produce a DTree $T$ consisting of cuboids, with root $\hQ$, and with leaves $Q_{i_\st},\cdots,Q_{i_\en}$ (and with no other leaves). Each node of $T$ is either a leaf, the parent of a leaf, or has precisely two children.

The work and storage used to do so are at most $C \cdot (i_\en - i _\st + 1) \cdot \log (i_\en - i _\st + 2) $ and $C \cdot (i_\en - i _\st + 1)$, respectively. Here, $C$ depends only on the dimension $n$.

\begin{proof}[\underline{Explanation}]

If any $Q_i = \hQ$, then there is only one $Q_i$, so we take our DTree to consist only of $\hQ$. Suppose otherwise.

If $i_\en - i_\st + 1 \leq 2$, then we can take our DTree to have root $\hQ$, and take $\hQ$ to have children $Q_{i_\st},\cdots,Q_{i_\en}$. Thus, our present algorithm is trivial in this case.

Suppose instead $i_\en - i_\st + 1 \geq 3$. We take $Q^\#$ to be the least dyadic cuboid (under inclusion) that contains $Q_{i_\st}$ and $Q_{i_\en}$. Let $Q'$, $Q''$ be, respectively, the lesser and greater child of $Q^\#$ (as dyadic cuboids).

Since $Q_{i_\st} < \cdots < Q_{i_\en}$ and the $Q_j$ are pairwise disjoint, it follows that
\begin{itemize}
\item $Q^\# \subset\hQ$.
\item $Q_i \subset Q^\#$ for each $i = i_\st, \cdots, i_\en$.
\item There exists $j$ ($i_\st \leq j < i_\en$) such that
\[Q_{i_\st},\cdots,Q_j \subset Q' \quad \mbox{and} \quad Q_{j+1},\cdots,Q_{i_\en} \subset Q''.\]
\end{itemize}
(These properties are obvious for dyadic cuboids, since they are obvious for dyadic intervals; see our discussion of the $1$-$1$ correspondence $\Psi$ in Section \ref{sec_dc}.)

We can recursively apply the present algorithm to produce DTrees $T'$,$T''$ with roots $Q'$,$Q''$, respectively. The leaves of $T'$ are precisely $Q_{i_\st},\cdots,Q_j$; and the leaves of $T''$ are precisely $Q_{j+1},\cdots,Q_{i_\en}$.

Our DTree $T$ will consist of a root $\hQ$, together with $T'$ and $T''$. The two children of $\hQ$ will be $Q'$ and $Q''$. This $T$ is obviously as promised.

Aside from recursively calling on itself, the above algorithm uses work  at most $C \log (\#(T) + 1)$ and storage at most $C$ (not counting the storage used to hold $Q_{i_\st},\cdots, Q_{i_\en}$). The factor of $\log (\# (T) + 1)$ comes from a binary search used to find $j$.

Therefore, altogether, our algorithm uses work and storage at most $C \cdot \#(T) \log (\# (T) + 1)$ and $C \cdot \#(T)$, respectively, where $T$ is the DTree arising from the algorithm. Since the leaves of $T$ are precisely $Q_{i_\st},\cdots,Q_{i_\en}$, and since each node of $T$ (other than the leaves and parents of leaves) has $2$ children, we see that
\[ \#(T) \leq C \cdot \#(\mbox{leaves of} \; T) = C \cdot (i_\en - i_\st + 1).\]
Thus, the work and storage of the algorithm are as promised.

\end{proof}

\environmentA{Algorithm: Make DTree.} Suppose we are given a list of distinct cuboids $Q_1, \cdots, Q_N$. With work $\leq C N(1+ \log N)$ in space $C N$ we produce a DTree $T$ with the following properties.
\begin{itemize}
\item Each of the $Q_j$ is a node of $T$.
\item Each node of $T$ is marked as \underline{original} if and only if it is one of the $Q_j$. 
\item Furthermore, we mark each node $Q$ of $T$ to indicate either the smallest (under inclusion) original node containing $Q$, or else to indicate that no such original node exists.
\item The number of nodes of $T$ is at most $C N$.
\end{itemize}

\begin{proof}[\underline{Explanation}] 

We pick a cuboid $Q^{00}$ that strictly contains the cuboids $Q_1,\cdots, Q_N$.

Applying the algorithm \textsc{Make Forest}, we make a tree $T^{(1)}$ with nodes $Q^{00},Q_1,\cdots, Q_N$, such that $Q$ is a descendant of $Q'$ in $T^{(1)}$ if and only if $Q \subset Q'$. 

Note that $Q^{00}$ is the root of $T^{(1)}$. 

Recall that each non-root node in $T^{(1)}$ is marked with a pointer to its parent. We mark each internal node in $T^{(1)}$ with pointers to its children, using an obvious algorithm. We also mark each non-root node in $T^{(1)}$ to indicate that it is original. 

By repeatedly applying the algorithm \textsc{Fill in gaps} to a node $Q$ of $T^{(1)}$ together with a list of its children (sorted under $<$), we imbed the tree $T^{(1)}$ in a DTree $T^{(2)}$ that has at most $C N$ nodes. The nodes in $T^{(1)}$ retain their markings in $T^{(2)}$. At each stage, we indicate the smallest original node containing each of the newly generated nodes (if such an original node exists). Indeed, if $Q = Q^{00}$, then we mark each of the new nodes $Q'$ generated by \textsc{Fill in gaps} to indicate that $Q'$ is not contained in any original node. If $Q \neq Q^{00}$, then we mark each of the new nodes $Q'$ generated by \textsc{Fill in gaps} to indicate that the smallest original node containing $Q'$ is the node $Q$.

This completes the construction of the marked DTree $T^{(2)}$.

One can easily check that the algorithm satisfies the desired work and storage bounds.

\end{proof}

\environmentA{Algorithm: Compute Norms From Marked Cuboids.} 

\label{pp8}

Suppose we are given a list $Q_1,\cdots,Q_{N}$ of distinct dyadic cuboids in $\R^n$, with each cuboid $Q_i$ marked with linear functionals $\mu_1^{Q_i},\cdots,\mu_{L_i}^{Q_i} : \R^D \rightarrow \R$. Let $\widehat{N} = \sum_{i=1}^N (L_i+1)$. 

Given $1 < p < \infty$, we perform one-time work at most $C \widehat{N} ( 1 + \log \widehat{N})$ in space $C \widehat{N}$, after which we can answer queries as follows:

A query consists of a dyadic cuboid $Q$ in $\R^n$.

The response to the query $Q$ is a list of linear functionals $\hmu_1^Q, \cdots, \hmu_D^Q : \R^D \rightarrow \R$, for which we guarantee the estimate
$$c \cdot \sum_{j=1}^D \lvert \hmu^Q_j(v) \rvert^p \leq \sum_{Q_i \subset Q} \sum_{j=1}^{L_i} \lvert \mu_j^{Q_i}(v) \rvert^p \leq C \sum_{j=1}^D \lvert \hmu^Q_j(v) \rvert^p \qquad \mbox{for all} \; v \in \R^D.$$
The work and storage used to answer a query are at most $C \cdot ( 1 + \log N)$. Here, $c$ and $C$ depend only on $n$, $p$, and $D$.

\begin{proof}[\underline{Explanation}]  

For each $i=1,\cdots,N$, we apply \textsc{Compress Norms} (see Section \ref{sec_lf}) to the functionals $\mu^{Q_i}_j$ ($1 \leq j \leq L_i$). We obtain linear functionals $\overline{\mu}^{Q_i}_1, \cdots, \overline{\mu}^{Q_i}_D : \R^D \rightarrow \R$ such that
\begin{equation}\label{ccnorm}
c \cdot \sum_{j=1}^D \lvert \overline{\mu}_j^{Q_i}(v) \rvert^p \leq \sum_{j=1}^{L_i} \lvert \mu_j^{Q_i}(v) \rvert^p \leq C \cdot \sum_{j=1}^D \lvert \overline{\mu}_j^{Q_i}(v) \rvert^p \qquad \mbox{for all} \; v \in \R^D.
\end{equation}

We then construct a Dtree $T$ such that each cuboid $Q_j$ is a node of $T$. This requires an application of the algorithm \textsc{Make DTree}. We guarantee that $T$ has at most $C N$ nodes.

We mark each node $Q_i$ of $T$ ($1 \leq i \leq N$) with the linear functionals $\overline{\mu}_1^{Q_i},\cdots,\overline{\mu}_D^{Q_i} : \R^D \rightarrow \R$. We mark each node $\widetilde{Q}$ of $T$ that is not among $Q_1,\cdots,Q_N$ with linear functionals $\overline{\mu}_1^{\widetilde{Q}},\cdots,\overline{\mu}_D^{\widetilde{Q}} : \R^D \rightarrow \R$, all of which are simply zero. 

Equipped with these markings, $T$ becomes an ADTree. (See Section \ref{pdef_sec} for the definition of an ADTree.)

Applying the algorithm \textsc{ADProcess} to the ADTree $T$, we perform one-time work, after which we can answer queries. A query consists of a dyadic cuboid $Q$ in $\R^n$. The response to the query $Q$ is a set of linear functionals $\hmu_1^Q,\cdots,\hmu_D^Q :\R^D \rightarrow \R$ such that
$$c \sum_{j=1}^D \lvert \hmu^Q_j(v) \rvert^p \leq \sum_{Q_i \subset Q} \sum_{j=1}^D \lvert \overline{\mu}_j^{Q_i}(v) \rvert^p \leq C \sum_{j=1}^D \lvert \hmu^Q_j(v) \rvert^p \qquad \mbox{for all} \; v \in \R^D.$$
This estimate and \eqref{ccnorm} show that the functionals $\hmu_1^Q,\cdots,\hmu_D^Q$ satisfy the conclusion of the present algorithm. The work and storage used are easily seen to be as promised.

\end{proof}

\environmentA{Algorithm: Placing a Point Inside Target Cuboids.}

Given a list of dyadic cuboids $Q_1,\cdots,Q_N \subset Q^\circ$ (not necessarily pairwise disjoint), we perform one-time work $\leq C N (1 + \log N)$ in space $\leq CN$, after which we can answer queries as follows:

A query consists of a point $\underline{x} \in \R^n$. 

The response to a query $\underline{x}$ is either one of the $Q_i$ containing $\underline{x}$, or else a promise that no such $Q_i$ exists.

The work to answer a query is at most $C \cdot ( 1 + \log N)$.

Here, $C$ depends only on the dimension $n$.

\begin{proof}[\underline{Explanation}] 

We construct a DTree $T$ and its control tree $\T = \CT(T)$ with at most $C N$ nodes (using the paperback edition of \textsc{Make Control Tree}), such that each of the $Q_j$ is a node of $T$. We mark each node of $T$ as \underline{original} if and only if it is one of the $Q_j$. Furthermore, we mark each node $Q$ of $T$ to indicate either the smallest (under inclusion) original node containing $Q$, or else to indicate that no such original node exists. (See \textsc{Make DTree}.)

Recall that each internal node $\xi \in \T$ has two children, marked as $\gchild(\xi)$ and $\schild(\xi)$. Each internal node $\xi \in \T$ is also marked with a node $Q_{\spl}(\xi) \in T$. Each node $\xi \in \T$ is marked with a node $Q_{\rt}(\xi) = $ the root of the DTree $\BT(\xi)$.

Recall that we've marked each $\xi$ to say whether $\BT(\xi)$ is a singleton.

Recall that $\BT(\rt(\T)) = T$. Also recall that $\BT(\gchild(\xi)) = \desc(Q_{\spl}(\xi),\BT(\xi))$, and $\BT(\schild(\xi)) = \ndesc(Q_{\spl}(\xi),\BT(\xi))$ for each internal node $\xi \in \T$.

Let $\underline{x} \in \R^n$ and $\xi \in \T$ be given. We consider the following procedure.

\underline{\textsc{Procedure Find-Original-Node $(\underline{x},\xi)$:}} \emph{We determine whether $\underline{x}$ is contained in an original node in $\BT(\xi)$. If such an original node exists, then we exhibit one.}

The above procedure answers our query when applied to $\xi = \rt(\T)$.

We now assume that $\xi \in \T$ is arbitrary. We ask whether $\underline{x}$ is contained in an original node in $\BT(\xi)$.

To study our question, we compare $\underline{x}$ with the root $Q_{\rt}(\xi)$ of $\BT(\xi)$. If $\underline{x} \notin Q_{\rt}(\xi)$, the answer is obviously \underline{NO}.

(No original node $Q \in \BT(\xi)$ contains $\underline{x}$.)

Suppose $\underline{x} \in Q_{\rt}(\xi)$. If $Q_{\rt}(\xi)$ is \underline{original}, the answer is obviously \underline{YES}, and we exhibit $Q_{\rt}(\xi)$ as an original node containing $\underline{x}$.

Suppose $Q_{\rt}(\xi)$ is \underline{not original}. If $Q_{\rt}(\xi)$ is a leaf of $\BT(\xi)$ (i.e., $\BT(\xi)$ is a singleton), then the answer is obviously \underline{NO}.

Suppose $Q_{\rt}(\xi)$ isn't a leaf of $\BT(\xi)$ (i.e., $\BT(\xi)$ isn't a singleton). We then compare $\underline{x}$ with $Q_{\spl}(\xi) \in \BT(\xi)$.

If $\underline{x} \notin Q_{\spl}(\xi)$, then all the descendants of $Q_{\spl}(\xi)$ in $\BT(\xi)$ are irrelevant for our discussion, i.e., they can't possibly contain $\underline{x}$. Therefore, in this case it's enough to ask whether $\underline{x}$ is contained in an original cuboid of $\ndesc(Q_{\spl}(\xi),\BT(\xi))$.

Thus, in this case, we can pass from the root $\xi$ to $\xi^+ = \schild(\xi) \in \T$, and we can answer our question by recursion.

On the other hand, suppose $\underline{x} \in Q_{\spl}(\xi)$. We examine the following two situations.

\begin{itemize}
\item Suppose $Q_{\spl}(\xi)$ \underline{is contained} in an original node of $T$. Furthermore, suppose that the smallest original node $Q \in T$ containing $Q_{\spl}(\xi)$ is contained in $Q_{\rt}(\xi)$. Note that $Q_{\spl}(\xi)$ and $Q_{\rt}(\xi)$ are nodes in $\BT(\xi)$, and $Q_{\spl}(\xi) \leq Q \leq Q_{\rt}(\xi)$ (inclusion), hence  $Q$ is a node of $\BT(\xi)$  because $\BT(\xi)$ is a sub-tree of $T$ (see Lemma \ref{ggg_lem}). Thus, the answer to our question is \underline{YES}. We exhibit the original node $Q \in \BT(\xi)$ containing $\underline{x}$.

\item Suppose that either $Q_{\spl}(\xi)$ \underline{is not contained} in an original node of $T$, or that the smallest original node containing $Q_{\spl}(\xi)$ \underline{is not contained} in $Q_{\rt}(\xi)$. This means that none of the original nodes in $\BT(\xi)$ contain $Q_{\spl}(\xi)$. Then any original node in $\BT(\xi)$ that contains $\underline{x}$ must be a descendant of $Q_{\spl}(\xi)$.

Therefore, we may pass from $\xi$ to $\xi^- = \gchild(\xi) \in \T$, and we can answer our question by recursion.
\end{itemize}

So, in all cases, we can answer our question.

The one-time work is at most $C N (1 + \log N)$. The query work, apart from recursing, is at most $C$. We recurse at most $C \cdot (1+ \log N)$ times, since $\T$ has depth at most $C \cdot (1 + \log N)$. So, the query work is at most $C \cdot ( 1 + \log N)$, as desired.

This completes the description of the procedure \textsc{Find-Original-Node}. As mentioned before, this yields the algorithm \textsc{Placing a Point Inside Target Cuboids}.

\end{proof}

\section{The Callahan-Kosaraju Decomposition} \label{sec_CK}

Let $E \subset \mathbb{R}^n$ with $\#(E)=N \geq 2$. A \emph{well-separated pairs decomposition} (WSPD) of $E$ is a finite sequence of Cartesian products $E'_1\times E_1'',\ldots,E'_L\times E_L''$ contained in $E\times E$, with the following properties.

\begin{enumerate}
\item[(WSPD1)] Each pair $(x',x'') \in E\times E$ with $x'\not=x''$ belongs to exactly one of the sets $E'_\ell \times E_\ell ''$ (for $\ell=1,\ldots,L$). Moreover, $E_\ell'\cap E_\ell'' = \emptyset$ for $\ell=1,\ldots, L$.
\item[(WSPD2)] For each $\ell=1,\ldots, L$, we have $\diam(E_\ell') + \diam(E_\ell'') \leq 10^{-10} \cdot \dist(E_\ell',E_\ell'')$.
\item[(WSPD3)] We have $L \leq CN$, for some constant $C$ depending only on the dimension $n$.
\end{enumerate}

The next algorithm arises in the work of Callahan and Kosaraju in \cite{CK} (see also \cite{FK2}).

\environmentA{Algorithm: Make WSPD.} With work at most $CN \log N$ in space at most $C N$, we compute a WSPD for $E$ and we compute representative pairs $(x_\ell',x_\ell'')\in E_\ell'\times E_\ell''$ for $\ell = 1,\cdots,L$. 

We do not explain here what it means to ``compute'' a WSPD. This does not matter, however, since we will only need the representative pairs $(x_\ell',x_\ell'')$.

\section{The BBD Tree}\label{sec_bbd}

We recall a few of the results of Arya, Mount, Netanyahu, Silverman and
Wu in \cite{A1}.

Given $E \subset \R^n$ such that $\#(E) = N \geq 2$, and given $x \in \R^n$, we can enumerate the points of $E$ as $y_1,\cdots,y_N$ so that 
\[ \lvert x - y_1 \rvert \leq \lvert x - y_2 \rvert \leq \cdots \leq \lvert x - y_N \rvert.\] 
We define $d_k(x,E) = \lvert x - y_k \rvert$. The definition of $d_k(x,E)$ is clearly independent of the chosen enumeration.

The following result is contained in \cite{A1} (see also \cite{FK2}).

\begin{thm}\label{bbd_thm}
There exists an algorithm with the following properties:
\begin{itemize}
\item The algorithm receives as input a subset $E$ with $\#(E) = N \geq 2$. The algorithm performs one-time work at most $C N \log N$ using storage $C N $, after which the algorithm is prepared to answer queries.
\item A query consists of a point $x\in \mathbb{R}^n$.
\item The answer to a query $x$ consists of two distinct points $\tilde{x}_1,\tilde{x}_2 \in E$ with $|x-\tilde{x}_1|\leq 2 d_1(x,E)$ and  $|x-\tilde{x}_2|\leq 2 d_2(x,E)$.
\item The work required to answer a query is at most $C \log N $.
\item Here, $C$ depends only on the dimension $n$.
\end{itemize}
\end{thm}
The proof of Theorem \ref{bbd_thm} relies on a data structure called a BBD Tree, which is associated to $E$. As another application of the BBD Tree, we have the following algorithm (see \cite{FK2}).

\environmentA{Algorithm: RCZ.} Given real numbers $\lambda(x)$ ($x \in E$), we perform one-time work at most $ C N \log N$ in space $C N$, after which we can do the following.

Given a dyadic cuboid $Q$, we can compute the following numbers and points.
 \begin{itemize}
 \item $\#(E \cap Q)$.
 \item $\min \{ \lambda(x) : x \in E \cap Q\}$ (or a promise that $E \cap Q$ is empty).
 \item A representative point $x(Q) \in E \cap Q$ (if $E \cap Q \neq \emptyset$).
 \item $\diam(Q \cap E)$.
 \end{itemize}
This computation requires work at most $C \log N$

\begin{proof}[\underline{Explanation}] 
The computation follows directly from \textsc{Algorithm RCZ1} in Section 25 of \cite{FK2} and \textsc{Algorithm REP1} in Section 27 of \cite{FK2}. These algorithms explain how to compute the quantities $\min \left\{ \delta(x,\cA) : x \in E \cap Q \right\}$ and $\#(E \cap Q)$, and how to compute a representative point $x(Q) \in E \cap Q$ when $E \cap Q \neq \emptyset$. The numbers $\delta(x,\cA)$ ($x \in E$) in Section 25 of \cite{FK2} are treated as arbitrary given real numbers (except in Lemma 1 in Section 25 of \cite{FK2}, which is not used elsewhere in the relevant algorithms). See Section \ref{sec_czdecomp} of this paper for a related discussion.

For each $1 \leq i \leq n$ we define coordinate functions $\lambda_i(x) = x_i$ for $x = (x_1,\cdots,x_n) \in E$. Applying the computation in the first bullet point, we compute the quantities
\begin{align*}
&r_i := \min \{ x_i : x =(x_1,\cdots,x_n) \in E \cap Q \}, \; \mbox{and} \\
&s_i := \max \{ x_i : x =(x_1,\cdots,x_n) \in E \cap Q \} \quad (1 \leq i \leq n).
\end{align*}
Thus, we can compute $\diam(E \cap Q) = \max \left\{ \lvert s_i - r_i \rvert : i=1,\cdots,n \right\}$. (Recall that diameters are measured using the $\ell^\infty$ norm.)

\end{proof}

\begin{remk} \label{bbd_rem}
Let $Q$ be a dyadic cube. We can decide whether $3 Q \cap E$ is nonempty, and if so we compute $\min \left\{ \lambda(x) : x \in E \cap 3 Q \right\}$ and $\#(E \cap 3  Q)$. We use a divide and conquer strategy. We write $3Q$ as the disjoint union of $3^n$ dyadic cubes $Q_\nu$ of sidelength $\delta_Q$. We apply  \textsc{Algorithm: RCZ} to each $Q_\nu$. We can tell whether $E \cap 3Q$ is nonempty by checking whether $E \cap Q_\nu$ is nonempty for some $\nu$. We complete the computation using the formulas
\begin{align*}
& \min \left\{ \lambda(x) : x \in E \cap 3 Q \right\} = \min_\nu \min \left\{ \lambda(x) : x \in E \cap Q_\nu \right\} \\
& \#(E \cap 3Q) = \sum_\nu \# (E \cap Q_\nu).
\end{align*}
Similarly, we can compute $\min \left\{ \lambda(x) : x \in E \cap \frac{65}{64}Q \right\}$ and $\#(E \cap \frac{65}{64}Q)$ (if $E \cap \frac{65}{64} Q \neq \emptyset$). Here, we use the fact that $\frac{65}{64}Q$ is the disjoint union of $130^n$ dyadic cubes of sidelength $\frac{1}{128}\delta_Q$.

Hence, by replicating the argument in \textsc{Algorithm RCZ}, we can compute $\diam(E \cap 3Q)$ and $\diam(E \cap \frac{65}{64}Q)$.

All the above computations requires work at most $C \log N$ after the one-time work of \textsc{Algorithm: RCZ} has been carried out.
\end{remk}

\section{Clusters}\label{clusters}

Suppose we are given $E \subset \mathbb{R}^n$ with $\#(E)=N \geq 2$. 

Suppose we are given $A \geq 1$. Assume $A$ exceeds a large enough constant determined by $n$. Assume also that $A$ is an integer power of $2$.

In this section, let $C,c,C'$, etc. denote constants determined by $n$, and let $C(A),c(A),C'(A)$, etc. denote constants determined by $A$ and $n$.

These symbols may denote different constants in different occurrences.

Recall that we use the $l^\infty$-norm and $l^\infty$-metric on $\mathbb{R}^n$: For $x=(x_1,\ldots,x_n)\in \mathbb{R}^n$, we have $|x|=\max_{1\leq i \leq n}|x_i|$. 

A subset $S\subset E$ is called a \underline{cluster} if 
$$\#(S)\geq 2\text{ and } \dist(S,E\setminus S)\geq A^3 \cdot \diam (S),$$
where $\dist(S,E\setminus S)=\infty$ if $E\setminus S = \emptyset$.

A cluster $S$ is called a \underline{strong cluster} if
$$\dist(S,E\setminus S) \geq A^5\cdot \diam(S),$$
where $\dist(S,E\setminus S)=\infty$ if $E\setminus S = \emptyset$.

A cluster that is not a strong cluster is called a \underline{weak cluster}.

Let $S$ be any finite non-empty subset of $\mathbb{R}^n$. 

The lower left corner of $S$, denoted by $\LLC(S)$, is defined by
$$\LLC(S)= (x_1,\ldots,x_n)\in \mathbb{R}^n \; \mbox{where} \;  x_i = \min \left\{y_i: y=(y_1,\ldots,y_n)\in S\right\} \; \mbox{for} \; 1\leq i \leq n .$$

The upper right corner of $S$, denoted by $\URC(S)$, is defined by
$$\URC(S) = (x_1,\ldots,x_n)\in \mathbb{R}^n \; \mbox{where} \; x_i = \max \left\{y_i: y=(y_1,\ldots,y_n)\in S\right\} \; \mbox{for} \; 1\leq i \leq n .$$

Note that $\diam(S)= \lv \LLC(S)-\URC(S) \rv$, since we are using the $l^\infty$ metric on $\mathbb{R}^n$.

If $S\subset \mathbb{R}^n$ is finite and contains at least 2 points, then we define the \underline{descriptor cube} of $S$ to be the smallest dyadic cube $Q$ such that $Q$ contains $\LLC(S)$ and $3Q$ contains $\URC(S)$. We note that $S$ has one and only one descriptor cube.

We write $\DC(S)$ for the descriptor cube of $S$.

If $\LLC(S)=(x_1^{\downarrow},\ldots,x_n^{\downarrow})$ and $\URC(S)=(x_1^{\uparrow},\ldots,x_n^\uparrow)$, then any point $(x_1,\ldots,x_n)\in S$ satisfies $x_i^\downarrow\leq x_i\leq x_i^\uparrow$ for $1\leq i \leq n$. 

Hence, if $Q=I_1\times \cdots \times I_n$ is the descriptor cube of $S$, we have $x_i^\downarrow,x_i^\uparrow \in 3I_i$ for $1\leq i \leq n$, consequently, we have $x_i \in 3I_i$ for $1\leq i \leq n$. It follows that $x \in 3Q$.

Thus, if $Q=\DC(S)$, then $S \subset 3Q$. Moreover, $\diam(S)\geq c\delta_Q$, by the minimal property of $Q$. 

\environmentA{Algorithm: Find Descriptor Cube.} 

We perform one-time work $\leq CN\log N$ in space $CN$, after which we answer queries as follows:

A query consists of a dyadic cube $\underline{Q}$. The response to a query is as follows:

Either we guarantee that $\#(3 \underline{Q} \cap E) \leq 1$, or we guarantee that $\#(3 \underline{Q} \cap E) \geq 2$ and we compute the descriptor cube $\DC(3\underline{Q}\cap E)$ together with the points $\LLC(3\underline{Q}\cap E)$ and  $\URC(3\underline{Q}\cap E)$. The query work is at most $C\log N$. 

\begin{proof}[\underline{Explanation}] 

We perform the one-time work of the BBD tree, after which we can do the following:

We compute $\#(3 \underline{Q} \cap E)$ using \textsc{Algorithm: RCZ} (see Remark \ref{bbd_rem}). If $\#(3 \underline{Q} \cap E) \leq 1$ then we indicate as such and terminate the computation. Otherwise, we guarantee that $\#(3 \underline{Q} \cap E) \geq 2$ and proceed as follows.

Suppose we assign to each $x\in E$ a label $\lambda(x) \in \mathbb{R}$. After one-time work at most $CN \log(N) $ we can answer queries, as follows: 

A query is a dyadic cube $\underline{Q}$ and a response to a query is $\max_{x \in E \cap 3 \underline{Q}} \lambda(x)$ and $\min_{x \in E \cap 3\underline{Q}} \lambda(x)$. The query work is at most $C\log N$.  (See Remark \ref{bbd_rem}.)

Taking $\lambda(x)$ to be the $i^\text{th}$ coordinate of $x$ for each $x \in E$ and looping over all $i$, we see that we can perform one-time work at most $CN \log(N)$, after which, given any dyadic query cube $\underline{Q}$, we can compute $\LLC(3\underline{Q} \cap E)$ and $\URC(3 \underline{Q} \cap E)$ with work at most $C \log N$.

After we obtain $\LLC(3\underline{Q} \cap E)$ and $\URC(3 \underline{Q} \cap E)$, we can compute $\DC(3\underline{Q}\cap E)$ with work at most $C$.

This completes the explanation of the algorithm \textsc{Find Descriptor Cube}.

\end{proof}

\environmentA{Algorithm: Make Cluster Descriptors.}

We produce a list of dyadic cubes $Q_1^{\text{CD}},\ldots, Q_L^{\text{CD}}$, with the following properties:

\begin{itemize}
\item For each $l$, the set $S_l=3Q_l^{\text{CD}}\cap E$ is a cluster.
\item Every strong cluster is one of the the $S_l$ above.
\item For each $l$, the cube $Q_l^{\text{CD}}$ is the descriptor cube of $S_l$.
\item $L\leq CN$.
\item The cubes $Q_1^{\text{CD}},\ldots, Q_L^{\text{CD}}$ are all distinct. 
\end{itemize}

The algorithm uses work at most $CN \log N$ in space $CN$. 

\begin{proof}[\underline{Explanation}] 

We perform the one-time work to make representatives $( x_{\nu }^{\prime },x_{\nu }^{\prime \prime }) _{\nu =1,\cdots
,\nu _{\max }}$ of the Well-Separated Pairs Decomposition of $E$. Thus, $\nu
_{\max }\leq CN$ and for any $x^{\prime }$, $x^{\prime \prime }\in E$ with $%
x^{\prime }\not=x^{\prime \prime }$, there exists $\nu $ such that 
\[
\left\vert x^{\prime }-x_{\nu }^{\prime }\right\vert +\left\vert x^{\prime
\prime }-x_{\nu }^{\prime \prime }\right\vert \leq 10^{-10}\left\vert
x^{\prime }-x^{\prime \prime }\right\vert 
\]%
and 
\[
\left\vert x^{\prime }-x_{\nu }^{\prime }\right\vert +\left\vert x^{\prime
\prime }-x_{\nu }^{\prime \prime }\right\vert \leq 10^{-10}\left\vert x_{\nu
}^{\prime }-x_{\nu }^{\prime \prime }\right\vert \text{.}
\]
(See Section \ref{sec_CK}.)

We perform the one-time work of the algorithm \textsc{Find Descriptor Cube}.

We perform the one-time work of the BBD tree. After that, given a dyadic cube $\underline{Q}$, we can compute $\#(E\cap 3\underline{Q})$ in time $C\log N$. 
(See the algorithm \textsc{RCZ} in Section \ref{sec_bbd}.)

For each $\nu $, we let $S_{\nu }^{\text{cand}}=3Q_{\nu }^{\text{cand}}\cap E
$, where $Q_{\nu }^{\text{cand}}$ is a dyadic cube containing $x'_{\nu }$
with $2\left\vert x_{\nu }^{\prime }-x_{\nu }^{\prime \prime }\right\vert
\leq \delta _{Q_{\nu }^{\text{cand}}}\leq 8\left\vert x_{\nu }^{\prime
}-x_{\nu }^{\prime \prime }\right\vert $. 

Instead of computing $S_{\nu }^{\text{cand}}$ (which will take too much
work), we compute $Q_{\nu }^{\text{cand}}$.

To test whether we ``like'' $S_{\nu }^{\text{cand}}$, we test whether 
\begin{equation}
\#( E\cap 3Q_{\nu }^{\text{cand}}) =\#( E\cap 3\hat{Q}_{\nu
}^{\text{cand}}) \text{,}  \label{testnumber}
\end{equation}%
where  $\hat{Q}_{\nu }^{\text{cand}}\supset Q_{\nu }^{\text{cand}}$ is a dyadic cube with sidelength $A^4\delta _{Q_{\nu }^{\text{cand}}}$. (This
test takes work $C\log N$ after we perform the one-time work of the BBD
tree. Recall that $A$ is a power of $2$.)

If we like $S_{\nu }^{\text{cand}}$ (i.e., \eqref{testnumber} holds), then we apply the query algorithm within the
algorithm \textsc{Find Descriptor Cube} to find 
\[
\DC( S_{\nu }^{\text{cand}}) =\DC( 3Q_{\nu }^{\text{cand}}\cap
E) .
\]

We then add $\DC( S_{\nu }^{\text{cand}}) $ to the list of the
cubes $\left\{ Q_{l}^{\text{CD}}\right\} $.

If (\ref{testnumber}) does not hold, we do nothing further regarding $Q_{\nu
}^{\text{cand}}$.

Thus, we produce a list of cubes 
\[
Q_{1}^{\text{CD}},\cdots ,Q_{L}^{\text{CD}}\text{.}
\]%
Note that $L\leq \nu _{\max }\leq CN$, since each $Q_{\ell}^{\text{CD}}$ arises
from $Q_{\nu }^{\text{cand}}$ for some $1\leq \nu \leq \nu _{\max }$.

If we like $S_{\nu }^{\text{cand}}$, then $S_{\nu }^{\text{cand}}=3Q_{\nu }^{\text{cand}}\cap E
$ is a
cluster. Indeed, since 
\[
\#( E\cap 3Q_{\nu }^{\text{cand}}) =\#( E\cap 3\hat{Q}_{\nu
}^{\text{cand}}) \text{,}
\]%
we have 
\begin{eqnarray}
\dist( S_{\nu }^{\text{cand}}, E\setminus S_{\nu }^{\text{%
cand}})  &\geq &\dist( 3Q_{\nu }^{\text{cand}},\mathbb{R}%
^{n}\setminus 3\hat{Q}_{\nu }^{\text{cand}})  \label{ptk1}\\
&\geq &cA^{4}\delta _{Q_{\nu }^{\text{cand}}}\geq c^{\prime }A^{4}\cdot\diam%
( S_{\nu }^{\text{cand}}) \text{;}\notag
\end{eqnarray}%
in obtaining inequality \eqref{ptk1}, we used $S_{\nu }^{\text{cand}}=3Q_{\nu }^{\text{cand}}\cap E\subset 3Q_{\nu
}^{\text{cand}}$ and
\begin{equation}
\#\left\{ ( E\setminus S_{\nu }^{\text{cand}}) \cap 3\hat{Q}_{\nu
}^{\text{cand}}\right\}  = \#\left\{ ( E\cap 3\hat{Q}_{\nu }^{\text{%
cand}}) \setminus ( E\cap 3Q_{\nu }^{\text{cand}}) \right\} = 0.
\end{equation}%
This completes the proof that $S_{\nu }^{\text{cand}}$ is a cluster whenever (\ref{testnumber})\ holds.

Next, we show that whenever we like $S_{\nu }^{\text{cand}}=3Q_{\nu }^{\text{%
cand}}\cap E$ (i.e., (\ref{testnumber})\ holds), then $$\left\{
Q_{l}^{\text{CD}}:=\DC( S_{\nu }^{\text{cand}}) ,S_{l}:=S_{\nu }^{\text{cand%
}}\right\} $$ satisfies 
\begin{equation}
S_{l}=3Q_{l}^{\text{CD}}\cap E\text{.} \label{ptk0}
\end{equation}
Indeed, since $Q_{l}^{\text{CD}}=\DC( S_{\nu }^{\text{cand}}) $,
we have 
\begin{equation}
S_{l}=S_{\nu }^{\text{cand}}\subset 3Q_{l}^{\text{CD}}\cap E. \label{ptk2} \end{equation} On the
other hand, $\delta _{Q_{l}^{\text{CD}}}$ and $\delta _{Q_{\nu }^{\text{cand}%
}}$ are both comparable to $\left\vert x_{\nu }^{\prime }-x_{\nu }^{\prime
\prime }\right\vert $. Indeed, for $Q_\nu^\text{cand}$, this follows from the defining condition. For $Q_{l}^{\text{CD}}$, we have $x_{\nu }^{\prime }$,$x_{\nu }^{\prime \prime }\in
3Q_{\nu }^{\text{cand}}\cap E=S_{\nu }^{\text{cand}}\subset 3Q_{l}^{\text{CD}}$,
hence%
\[
\left\vert x_{\nu }^{\prime }-x_{\nu }^{\prime \prime }\right\vert \leq
C\delta _{Q_{l}^{\text{CD}}}.
\]%
Also, since $S_{\nu }^{\text{cand}}=3Q_{\nu }^{\,\text{cand}}\cap E$, we
have 
\[
\diam( S_{\nu }^{\text{cand}}) \leq 3\delta _{Q_{\nu }^{\text{cand}%
}}\leq C\left\vert x_{\nu }^{\prime }-x_{\nu }^{\prime \prime }\right\vert 
\text{,}
\]%
hence $Q_{l}^{\text{CD}}=\DC( S_{\nu }^{\text{cand}}) $ satisfies 
\[
\delta _{Q_{l}^{\text{CD}}}\leq C\cdot \diam( S_{\nu }^{\text{cand}}) \leq
C^{\prime }\left\vert x_{\nu }^{\prime }-x_{\nu }^{\prime \prime
}\right\vert \text{.}
\]%
Thus, as claimed, $\delta _{Q_{\nu }^{\text{cand}}}$ and $\delta _{Q_{l}^{\text{CD}}}$ are
comparable to $\left\vert x_{\nu }^{\prime }-x_{\nu }^{\prime \prime
}\right\vert $.

Furthermore, since $x_{\nu }^{\prime }\in 3Q_{\nu }^{\text{cand}}\cap
3Q_{l}^{\text{CD}}$, it follows that 
\[
3Q_{l}^{\text{CD}}\cap E\subset AQ_{\nu }^{\text{cand}}\cap E=3Q_{\nu }^{\text{cand}%
}\cap E\text{;}
\]%
the last equality holds since we like $S_{\nu }^{\text{cand}}$.

Therefore, 
\begin{equation}
3Q_{l}^{\text{CD}}\cap E\subset 3Q_{\nu }^{\text{cand}}\cap E=S_{\nu }^{\text{cand}%
}=S_{l}.  \label{ptk3}
\end{equation}

From (\ref{ptk2}) and (\ref{ptk3}), we obtain (\ref{ptk0}).

Since $Q_{l}^{\text{CD}}=\DC( S_{\nu }^{\text{cand}}) $, it now follows
that $Q_{l}^{\text{CD}}=\DC( S_{l}) $, and $S_{l}=3Q_{l}^{\text{CD}}\cap E.$

Since $S_{l}=S_{\nu }^{\text{cand}}$ and $S_{\nu }^{\text{cand}}$ is a
cluster, we have shown that $S_{l}$ is a cluster.

We have now proven the first, third, and fourth bullet points asserted in the specification of our algorithm \textsc{Make Cluster Descriptors}.

Next, we show the second bullet point: every strong cluster is one of the $S_l$. 

Indeed, let $S$ be a strong cluster. Thus, $\#( S) \geq 2$ and 
\[
\dist( S,E\setminus S) \geq A^{5}\cdot \diam( S) 
\text{.}
\]%
Fix $x^{\prime },x^{\prime \prime }\in S$ such that $\left\vert x^{\prime
}-x^{\prime \prime }\right\vert = \diam( S) $. Then we can find $%
\nu $ such that 
\[
\left\vert x_{\nu }^{\prime }-x^{\prime }\right\vert +\left\vert x_{\nu
}^{\prime \prime }-x^{\prime \prime }\right\vert \leq 10^{-10}\left\vert
x^{\prime }-x^{\prime \prime }\right\vert =10^{-10}\diam( S) 
\text{.}
\]%
We have 
\begin{equation}
\dist( x_{\nu }^{\prime },S) \leq \left\vert x_{\nu
}^{\prime }-x^{\prime }\right\vert \leq 10^{-10}\diam( S) 
\label{ptk4}
\end{equation}%
and 
\begin{equation}
\dist( x_{\nu }^{\prime \prime },S) \leq \left\vert x_{\nu
}^{\prime \prime }-x^{\prime \prime }\right\vert \leq 10^{-10}\diam%
( S) \text{.}  \label{ptk5}
\end{equation}%
Since $x_{\nu }^{\prime },x_{\nu }^{\prime \prime }\in E$ and dist$(
S,E\setminus S) \geq A^{5}\cdot \diam( S) $, from (\ref{ptk4})
and (\ref{ptk5}), we conclude that $x_{\nu }^{\prime },x_{\nu }^{\prime
\prime }\in S$.

Next, we show that $S=S_{\nu }^{\text{cand}}=3Q_{\nu }^{\text{cand}}\cap E$.

By definition of $Q_{\nu }^{\text{cand}}$, we have $x_{\nu }^{\prime }\in
Q_{\nu }^{\text{cand}}$ and 
\begin{equation}
2\left\vert x_{\nu }^{\prime }-x_{\nu }^{\prime \prime }\right\vert \leq \delta_{Q_{\nu }^{\text{cand}}}\leq 8\left\vert x_{\nu }^{\prime }-x_{\nu }^{\prime
\prime }\right\vert \text{.}  \label{ptk6}
\end{equation}%
Therefore, every point $z\in \mathbb{R}^{n}$ such that $\left\vert z-x_{\nu
}^{\prime }\right\vert \leq \frac{3}{2}\left\vert x_{\nu }^{\prime }-x_{\nu
}^{\prime \prime }\right\vert $ belongs to $3Q_{\nu }^{\text{cand}}$. On the other
hand, since $x_{\nu }^{\prime }\in S$, we know that every $x\in S$ satisfies 
\[
\left\vert x-x_{\nu }^{\prime }\right\vert \leq \diam( S)
=\left\vert x^{\prime }-x^{\prime \prime }\right\vert \leq \frac{3}{2}%
\left\vert x_{\nu }^{\prime }-x_{\nu }^{\prime \prime }\right\vert \text{.}
\]%
Therefore, $S\subset 3Q_{\nu }^{\text{cand}}.$ Since $S$ is a cluster, $%
S\subset E$. Thus, 
\[
S\subset 3Q_{\nu }^{\text{cand}}\cap E\text{.}
\]

If $S\not=3Q_{\nu }^{\text{cand}}\cap E$, then there would exist $\hat{x}\in
( 3Q_{\nu }^{\text{cand}}\cap E) \setminus S$. We would then have 
\begin{equation*}
\dist( E\setminus S,S)  \leq \left\vert \hat{x}-x_{\nu}^{\prime }\right\vert \leq C\delta _{Q_{\nu }^{\text{cand}}} \leq  C^{\prime }\left\vert x_{\nu }^{\prime }-x_{\nu }^{\prime \prime
}\right\vert  \leq C^{\prime} \diam( S) \text{,}
\end{equation*}
contradicting our assumption that $S$ is a strong cluster. 

This completes the proof that $S=S_{\nu }^{\text{cand}}=3Q_{\nu }^{\text{cand%
}}\cap E$, where the last equality follows by definition.

Next, we check that we like $S_{\nu }^{\text{cand}}$, i.e., that 
\[
\#( 3Q_{\nu }^{\text{cand}}\cap E) =\#( 3\hat{Q}_{\nu }^{%
\text{cand}}\cap E) \text{,}
\]%
where $\hat{Q}_{\nu }^{\text{cand}}\supset Q_{\nu }$ is a dyadic cube with $\delta _{\hat{Q}_{\nu }^{\text{cand}}}=A^{4}\delta
_{Q_{\nu }}$.

Indeed, suppose not. Then there would exist $\hat{x}\in ( 3\hat{Q}_{\nu
}^{\text{cand}}\cap E) \setminus ( 3Q_{\nu }^{\text{cand}}\cap
E) \subset E\setminus S$, where the last inclusion holds because $%
S=3Q_{\nu }^{\text{cand}}\cap E$. We have 
\begin{eqnarray*}
\dist( E\setminus S,S)  &\leq &\left\vert \hat{x}-x_{\nu
}^{\prime }\right\vert \leq C\delta_{\hat{Q}_{\nu }^{\text{cand}}} \\
&=&CA^4\delta _{Q_{\nu }^{\text{cand}}} \\
&\leq &C' A^{4}\left\vert x_{\nu }^{\prime }-x_{\nu }^{\prime \prime
}\right\vert \text{ \quad (see \eqref{ptk6})} \\
&\leq &C' A^{4}\cdot \diam( S)  \quad \text{ (since }x_{\nu
}^{\prime },x_{\nu }^{\prime \prime }\in S\text{),}
\end{eqnarray*}
contradicting our assumption that $S$ is a strong cluster. This completes
the proof that we like $S_{\nu }^{\text{cand}}$.

We now know that $\DC( S_{\nu }^{\text{cand}}) $ is one of the $%
Q_{l}^{\text{CD}}$, and that (for the same $l$), we have $S_{l}=S_{\nu }^{\text{cand%
}}=S$.

This proves the second bullet point asserted in our specification of the
algorithm \textsc{Make Cluster Descriptors}.

It remains to verify the last bullet point of the algorithm \textsc{Make
Cluster Descriptors}, i.e., the cubes $Q_{1}^{\text{CD}}$, $\cdots $, $Q_{L}^{\text{CD}}
$ are all distinct. Since $L\leq CN$, we can sort the $Q_{l}^{\text{CD}}$'s and
remove the duplicates with work $CN\log N$. 

Now all bullet points asserted in the specification of our algorithm hold.

The reader can easily check that the work and storage of our algorithm are
as promised.

\end{proof}

\begin{remk}
Note that the clusters $S_l$ produced (implicitly) by the above algorithm are all distinct, since their descriptor cubes $Q_l^{\text{CD}}$ are all distinct, and any cluster has one and only one descriptor cube. 
\end{remk}

\environmentA{Algorithm: Locate Relevant Cluster.} 

After performing the algorithm \textsc{Make Cluster Descriptors} and other one-time work, we can answer queries as follows:
A query consists of a point $\underline{x}\in \mathbb{R}^{n}$, for which
there exist a strong cluster $S$ and a point $x( S) \in S$ such
that 
\begin{equation}
A\diam( S) \leq \left\vert \underline{x}-x( S)
\right\vert \leq A^{-1}\dist( E\setminus S,S) \text{.}
\label{ptk9}
\end{equation}

We do not assume that $S$ or $x( S) $ is known.

The response to a query $\underline{x}$ is one of the descriptor cubes $%
Q_{l}^{\text{CD}}$ produced by the algorithm \textsc{Make Cluster Descriptors} such
that (\ref{ptk9}) holds for $S=S_{l}:=3Q_{l}^{\text{CD}}\cap E$ and for some $%
x( S) \in S$.

The one-time work is at most $CN\log N$ in space $CN$; the query work is at
most $C\log N$.

\begin{proof}[\underline{Explanation}]  

Suppose \eqref{ptk9} holds for some $S$ and for some $x( S) \in S$.

Assume that $\underline{x} \in S$. Then $A \cdot \diam(S) \leq \lv \underline{x} - x(S) \rv \leq \diam(S)$. This gives a contradiction if $A > 1$. Hence, we have shown that $\underline{x} \notin S$.

Assume next that $\underline{x} \in E \setminus S$. Then $\dist(E \setminus S, S) \leq \lv \underline{x} - x(S) \rv \leq A^{-1} \dist( E\setminus S,S) $. This gives a contradiction if $A > 1$. Hence, we have shown that $\underline{x} \notin E \setminus S$.

We have proven that $\underline{x} \notin E$.

Now,
\[
\dist( \underline{x},S) \leq \left\vert \underline{x}%
-x( S) \right\vert \leq A^{-1}\dist( E\setminus
S,S) \text{,}
\]%
and if $E \setminus S \neq \emptyset$ then
\begin{eqnarray*}
\dist( \underline{x},E\setminus S)  &\geq &\dist%
( E\setminus S,S) -\dist( \underline{x},S)  \\
&\geq &\dist( E\setminus S,S) -\left\vert \underline{x}-x( S) \right\vert \\
&\geq &\dist( E\setminus S,S) -A^{-1}\dist(E\setminus S,S)  \\
&\geq &(1/2)\dist( E\setminus S,S) \text{.}
\end{eqnarray*}
If $E \setminus S = \emptyset$, then by definition $\dist(\underline{x}, E\setminus S) = \infty$.

Therefore, 
\[
\dist( \underline{x},E\setminus S) \geq cA\cdot \dist( 
\underline{x},S) \text{,}
\]%
which yields 
\[
\dist( \underline{x},S) =\dist( \underline{x}%
,E) \text{.}
\]

Using the BBD tree, we compute a number $\Delta >0$ such that 
\[
8 \cdot \dist( \underline{x},E) \leq \Delta \leq 32 \cdot \dist( \underline{x},E),
\]
and such that $\Delta $ is a power of $2$.

We then produce the dyadic cube $Q^{\#}$ of sidelength $\Delta $ containing $%
\underline{x}$. 

We claim that $S=3Q^{\#}\cap E$. 

To see this, note that dist$( \underline{x},E\setminus S) \geq cA\cdot$%
dist$( \underline{x},S) =cA\cdot $dist$( \underline{x},E)
\geq c^{\prime }A\Delta $. On the other hand, $\underline{x}\in Q^{\#}$ and $\delta
_{Q^{\#}}=\Delta $.

Therefore, $3Q^{\#}\cap E\subset S$. If $3Q^{\#}\cap E \neq S$, then there
exists $\hat{x}\in S\setminus ( 3Q^{\#}\cap E) \subset S$. On the other
hand, since $\underline{x}\in Q^{\#}$ and $\delta _{Q^{\#}}=\Delta \geq 8$%
dist$( \underline{x},E) $, we know that $\frac{3}{2}Q^{\#}$
contains a point of $E$; say $\check{x}\in \frac{3}{2}Q^{\#}\cap E$. Note
that $\check{x}\in 3Q^{\#}\cap E\subset S$.

Thus, $\hat{x},\check{x}\in S$, with $\check{x}\in \frac{3}{2}Q^{\#}$ and $%
\hat{x}\not\in 3Q^{\#}$. Therefore, 
\begin{eqnarray*}
\diam( S)  &\geq &\left\vert \hat{x}-\check{x}\right\vert
\geq c\delta _{Q^{\#}}=c\Delta  \\
&\geq &c\cdot \dist( \underline{x},E) =c\cdot \dist( 
\underline{x},S)  \\
&\geq &c\left[ \left\vert \underline{x}-x( S) \right\vert - \diam( S) \right] \text{, since }x( S) \in S\text{,} \\
&\geq &c\left[ A\diam( S) -\diam( S) %
\right] \text{, by (\ref{ptk9}).}
\end{eqnarray*}%
Thus, $\diam(S) \geq c \cdot (A-1) \cdot \diam(
S) $. Since $S$ contains at least two points (because it is a
cluster), we have reached a contradiction. This completes the proof of our
claim that $S=3Q^{\#}\cap E$.

Since $S$ is a strong cluster, its descriptor cube $\DC( S) $ is
among the cubes $Q_{l}^{\text{CD}}$ produced by the algorithm \textsc{Make Cluster
Descriptors}.

We can compute $\DC( S) $ by applying the algorithm \textsc{Find
Descriptor Cube} to the query cube $Q^{\#}$; this produces $\DC(
S) $ because $S=3Q^{\#}\cap E$. 

Accordingly, our algorithm proceeds as follows:

\begin{itemize}
\item Compute $\Delta $, a power of $2$, such that $8$dist$( \underline{%
x},E) \leq \Delta \leq 32$dist$( \underline{x},E) $, using
the BBD tree.

\item Produce the dyadic cube $Q^{\#}$ of sidelength $\Delta $, containing $%
\underline{x}$.

\item Apply the query algorithm in \textsc{Find Descriptor Cube} to the
query cube $Q^{\#}$. This produces the cube $\DC( S) $.

\item By a binary search, locate $\DC( S) $ among the cubes $%
Q_{1}^{\text{CD}},\ldots ,Q_{L}^{\text{CD}}$ produced previously by the algorithm \textsc{%
Make Cluster Descriptors}. This produces one of the $Q_{l}^{\text{CD}}$, which is
the descriptor cube for the cluster $S$ in (\ref{ptk9}); in particular $%
S=S_{l}=3Q_{l}^{\text{CD}}\cap E$ satisfies (\ref{ptk9}), for some $x(
S) \in S$.
\end{itemize}

Thus, our algorithm does what we promised. The work and storage of the
algorithm are easily seen to be as promised also.

\end{proof}

\environmentA{Algorithm: Make Cluster Representatives.}

For each of the cubes $Q_{1}^{\text{CD}},\ldots ,Q_{L}^{\text{CD}}$ produced by the
algorithm \textsc{Make Cluster Descriptors}, we compute a point 
\[
x(S_{l}) \in S_{l}=3Q_{l}^{\text{CD}}\cap E\text{.}
\]%
The algorithm uses work $\leq CN\log N$ in space $CN$.

\begin{proof}[\underline{Explanation}] 

For each $Q_{l}^{\text{CD}}$ (a dyadic cube), we use the BBD tree to compute a
point $x(S_{l}) \in 3Q_{l}^{\text{CD}}\cap E$ (which we know to be non-empty,
since it is a cluster). We use the fact than $3Q_l^{\text{CD}}$ is the union of $3^n$ dyadic cubes; the required algorithm can be found in Section \ref{sec_bbd}.

The work and storage are as promised.

\end{proof}

Since every strong cluster is one of the $S_{l}$, we have computed a
representative point of every strong cluster and possibly also of some weak
clusters.

For each of the clusters $S= S_l :=  3Q_l^{\text{CD}} \cap E$ we define the \emph{halo} $H(S)$ by
\begin{equation}\label{halo}
H(S) =\left\{ y\in \mathbb{R}^{n}:A\cdot \diam(
S) <\left\vert y-x( S) \right\vert <A^{-1}\cdot \dist%
( E\setminus S,S) \right\}
\end{equation}
where $x(S)$ is the cluster representative produced by the algorithm \textsc{Make Cluster Representative}.

Finally, we recall Lemma 6.7 from \cite{FIL1}.
\begin{lem}
\label{path_lem}
Fix a cluster $S= 3Q_l^{\text{CD}} \cap E$. Suppose that $x \in H(S)$ and $x' \in H(S)$ satisfy $\lvert x - x(S) \rvert \geq \lvert x'- x(S) \rvert$. Assume furthermore that $x$ and $x'$ belong to the same connected component of $H(S)$. Then there exist a finite sequence of points $x_1,\cdots,x_{\underline{L}} \in H(S)$, and a positive integer $L_*$, with the following properties:

\begin{itemize}
\item $x_1 = x$ and $x_{\underline{L}} = x'$.
\item $\lvert x_{\ell+1} - x(S) \rvert \leq \lvert x_\ell - x(S) \rvert$ for $\ell = 1,\cdots, {\underline{L}} -1$.
\item $\lvert x_\ell - x_{\ell+1} \rvert \leq A^{-2} \lvert x_\ell - x(S) \rvert$ for $\ell = 1,\cdots , {\underline{L}} -1$.
\item $\lvert x_{\ell + L_*} - x(S) \rvert \leq (1-A^{-3}) \lvert x_\ell - x(S) \rvert$ for $1 \leq \ell \leq {\underline{L}} - L_*$.
\item $L_* \leq A^3$.
\end{itemize}
\end{lem}

\begin{remk} \label{rem_oned1}
Lemma 6.7 in \cite{FIL1} was stated incorrectly in dimension $n=1$. Here, we include the minor yet necessary modifications. The additional assumption that $x$ and $x'$ belong to the same connected component of $H(S)$ is required, since in dimension $n =1$ the halo $H(S)$ consists of two connected components. Clearly, if $x$ and $x'$ belong to distinct connected components there can be no finite sequence as in the statement of the lemma, and so this extra hypothesis is necessary. We shall not prove this lemma in the case $n=1$, as the argument is quite obvious. The proof of Lemma 6.7 in \cite{FIL1} remains valid when $n \geq 2$, since then the halos are connected.
\end{remk}

\section{Paths to Keystone Cubes}\label{sec_ptkc}

We assume we are given a finite subset $E\subset \mathbb{R}^{n}$, with $N=\#(
E) \geq 2$.

We are also given constants $K\geq 10$, $A\geq 10$. We assume that $A$ is greater than a large enough constant determined by $n$. We further assume that $A$ is a power of $2$ and that $K$ is an odd integer.

We write $c,C,C^{\prime }$, etc. to denote constants that depend only on the
dimension $n$; we write $c( K) $, $C( K) $, and $%
C^{\prime }( K) $, etc. to denote constants that depend only on $K$ and 
$n$; we write $c( A) $, $c( A) $, $C^{\prime }(
A) $, etc. to denote  constants that depend only on $A,K,n$. These symbols may denote different constants in different occurrences.

We suppose we are given a locally finite collection $\CZ$ consisting of dyadic
cubes that form a partition of $\mathbb{R}^{n}$. We do \textit{not} assume
that any list of $\CZ$ cubes is given; in fact, there are infinitely many $\CZ$
cubes. Rather, we assume that we have access to a $\CZ$-\textsc{Oracle}. Given a query
point $\underline{x}\in \mathbb{R}^{n}$, the $\CZ$-\textsc{Oracle} returns the one and
only $Q\in \CZ$ that contains $\underline{x}$. We do not count any calls to
the $\CZ$-\textsc{Oracle} in the one-time work or the query work of any of the algorithms presented here; we will instead keep track of the
\emph{number} of calls to the $\CZ$-\textsc{Oracle}.

We make the following assumptions on the decomposition $\CZ$.

\begin{itemize}
\item (Good geometry) If $Q,Q^{\prime }\in \CZ$ and $Q\leftrightarrow
Q^{\prime }$, then $\frac{1}{64}\delta _{Q}\leq \delta _{Q^{\prime }}\leq
64\delta _{Q}$.

\item ($E$ is nearby) For each $Q\in \CZ$, we have $\#( 9Q\cap E)
\geq 2$.
\end{itemize}

Due to good geometry, we see that there exists a constant $c_G > 0$, which is an integer power of $2$ depending only on the dimension $n$, such that, for any $Q,Q' \in \CZ$ we have
\begin{equation} \label{ggz}
(1+ 8c_G) Q \cap (1+ 8c_G)Q' \neq \emptyset \implies Q \leftrightarrow Q'.
\end{equation}

We next make a few definitions.

A finite sequence $\mathcal{S} = (Q_1,Q_2,\cdots,Q_{\underline{L}})$ consisting of $\CZ$ cubes is called a \emph{path} provided that any two consecutive cubes in the sequence touch. This property may be equivalently written as
\[ Q_1 \leftrightarrow Q_2 \leftrightarrow Q_3 \leftrightarrow \cdots \leftrightarrow Q_{\underline{L}}.\]
We sometimes say that the path $\mathcal{S}$ \emph{joins} $Q_1$ to $Q_{\underline{L}}$.

A path $\mathcal{S} = (Q_1,\cdots,Q_{\underline{L}})$ is called \emph{exponentially decreasing} if there exist constants $0 < c(A) < 1$ and $C(A) \geq 1$ such that
\[\delta_{Q_{\ell'}} \leq C(A) \cdot (1-c(A) )^{\ell' - \ell} \delta_{Q_\ell} \quad \mbox{for} \; 1 \leq \ell \leq \ell' \leq {\underline{L}}.\]
In particular, the constants $c(A)$ and $C(A)$ are assumed to be independent of the length ${\underline{L}}$ of the sequence $\mathcal{S}$.

A cube $Q\in \CZ$ is called \emph{keystone} provided $\delta
_{Q^{\prime }}\geq \delta _{Q}$ for each $Q^{\prime }\in \CZ$ that meets $KQ$%
. (Obviously, this definition depends on the choice of $K$, which will always be clear from the context.)

Recall that a subset $S\subset E$ is called a \underline{cluster} if $\#(S)\geq 2$ and $\dist(S,E\setminus S)\geq A^3 \cdot \diam (S)$.

There is a special collection of clusters $S = 3Q^{\text{CD}}_l \cap E$ ($1 \leq \ell \leq L$) arising in the algorithm \textsc{Make Cluster Descriptors}. We compute the representative point $x(S) \in S$ associated to each such cluster $S$ using the algorithm \textsc{Make Cluster Representatives}. For each such cluster, we let the halo $H(S)$ be defined as in \eqref{halo}.

From this point onward in the section, until further notice, we shall assume that $n \geq 2$. We make use of this assumption in several of the results that follow. Toward the end of this section we sketch the modifications required in dimension $n=1$.

Now suppose we are given $x\in ( 1+c_{G}) Q\cap H(S) $, where $Q\in \CZ$ and $S=3Q^{\text{CD}}_l \cap E$ for some fixed $1 \leq \ell \leq L$.

If $\delta _{Q}<\hat{c}\left\vert x-x( S) \right\vert $ for a
small enough constant $\hat{c}$, then since also $x\in ( 1+c_{G})
Q$, we have 
\begin{eqnarray*}
\dist( 10Q,x( S) )  &\geq &\left\vert x-x(
S) \right\vert -10\delta _{Q} \\
&\geq &(1/2) \cdot \left\vert x-x( S) \right\vert  \\
&\geq &(A/2) \cdot \diam( S) \text{,}
\end{eqnarray*}%
hence%
\[
\dist( 10Q,S) \geq \dist( 10Q,x( S)
) -\diam( S) \geq (A/4) \cdot \diam(S),\]
which implies $10Q\cap S=\emptyset $.

Additionally, note that 
\begin{equation}\label{sw1}
10Q\subset B( x( S) ,\left\vert x-x( S)
\right\vert +10\delta _{Q}) \subset B( x( S)
,2\left\vert x-x( S) \right\vert ),
\end{equation} since we assume $%
\delta _{Q}<\hat{c}\left\vert x-x( S) \right\vert$.

On the other hand, 
\[
\dist( x( S) ,E\setminus S) \geq \dist%
( S,E\setminus S) \geq A\left\vert x-x( S) \right\vert 
\]%
since we assume that $x\in H(S)$. 

Together with \eqref{sw1}, the above estimate tells us that $%
10Q\cap ( E\setminus S) =\emptyset $.

Since also $10Q\cap S=\emptyset $, we now know that $10Q\cap E=\emptyset $,
contradicting our assumption \textquotedblleft $E$ is nearby" for the cube $Q \in \CZ$. 

Thus, our assumption $\delta _{Q}<\hat{c}\left\vert x-x( S)
\right\vert $ must be false. We have proven the following.
\begin{equation}
\left[ 
\begin{array}{c}
\text{Suppose }Q\in \CZ\text{, }x\in ( 1+c_{G}) Q\cap H(
S) \text{, where }S = 3Q^{\text{CD}}_l \cap E. \\ 
\text{Then }\delta _{Q}\geq c\left\vert x-x( S) \right\vert \text{%
.}%
\end{array}%
\right]   \label{ptk12}
\end{equation}

\begin{remk}
Suppose $Q,Q^{\prime }\in \CZ$ and $x\in ( 1+c_{G}) Q,$ $x^{\prime
}\in ( 1+c_{G}) Q^{\prime }$. From \eqref{ggz} and the good geometry of $\CZ$ 
we see that 
\[
\delta _{Q^{\prime }}\leq C \cdot \left[ \delta _{Q}+\left\vert x-x^{\prime
}\right\vert \right].
\]%
This estimate and its analogue with $Q$ and $Q^{\prime }$ interchanged tell
us that 
\begin{equation}
c \cdot \bigl[ \delta _{Q}+\lvert x-x^{\prime }\rvert \bigr] \leq \bigl[
\delta _{Q^{\prime }}+\left\vert x-x^{\prime }\right\vert \bigr] \leq C \cdot
\bigl[ \delta _{Q^{\prime }}+\left\vert x-x^{\prime }\right\vert \bigr] 
\text{.}  \label{ptk11}
\end{equation}
\end{remk}

\begin{lem}\label{lem_cluster1}
Let $S = 3Q^{\text{CD}}_l \cap E$ be a cluster produced by the algorithm \textsc{Make Cluster Descriptors}, and let $x(S) \in S$ be its associated representative. Recall that
\begin{equation}
\label{halodefn}
H(S) =\left\{ y\in \mathbb{R}%
^{n}: A\cdot \diam( S) <\left\vert y-x( S)
\right\vert <A^{-1}\cdot \dist( E\setminus S,S) \right\}.
\end{equation}
Let $x\in ( 1+c_{G}) Q\cap H( S) $, with $Q\in \CZ$.
Finally let $Q^{CZ}( S) $ be the $\CZ$ cube containing $x(
S)$. Then $\delta _{Q}$ and $\delta _{Q^{CZ}( S) }+\left\vert x-x( S) \right\vert $ differ by at most a factor $C$.
\end{lem}
\begin{proof}
From \eqref{ptk11} we deduce that $\delta_{Q^{CZ}(S)} + \lvert x - x(S) \rvert$ and $\delta_Q + \lvert x - x(S) \rvert$ differ by at most a factor of $C$. From \eqref{ptk12}, we have $\delta_Q \geq c \lvert x - x(S) \rvert$. Combining these two estimates, we obtain the conclusion of the lemma.
\end{proof}

\begin{lem}\label{lem_cluster1a}
Let $S = 3Q^{\text{CD}}_l \cap E$ be a cluster produced by the algorithm \textsc{Make Cluster Descriptors}. Let $x,x' \in H(S)$, and let $Q,Q' \in \CZ$, with $x \in Q$ and $x' \in Q'$. If $\lvert x - x' \rvert \leq A^{-2}\lvert x - x(S) \rvert$, then $Q \leftrightarrow Q'$.
\end{lem}
\begin{proof}
By Lemma \ref{lem_cluster1}, we have 
\[ \lvert x' - x \rvert \leq A^{-2} \left[ \lvert x  - x(S) \rvert + \delta_{Q^{CZ}(S)} \right] \leq C A^{-2} \delta_Q.
\] 
Since $x \in Q$, we have $x' \in (1+c_G) Q$. Also $x' \in Q' \subset (1+c_G)Q'$. Thus $(1+c_G)Q \cap (1+c_G)Q' \neq \emptyset$, which implies that $Q \leftrightarrow Q'$. Here, we use \eqref{ggz}.
\end{proof}

We assume that we have done all the one-time work for the
algorithms in Section \ref{clusters}. Thus, the query algorithms from Section \ref{clusters} are at our disposal in the present section. Recall that the
one-time work just mentioned consists of work at most $CN\log N$ in space at
most $CN$.

\environmentA{Algorithm: Keystone-or-Not.} 

Given a cube $Q \in \CZ$, we produce one of the following outcomes:

\begin{enumerate}
\item[(KEY 1)] We guarantee that $Q$ is a keystone cube.

\item[(KEY 2)] We produce a cube $Q^{\prime }\in \CZ$ such that
\[
\delta _{Q^{\prime }}\leq \frac{1}{2}\delta _{Q} \; \mbox{and} \; Q^{\prime }\cap KQ\not=\emptyset \text{,}
\]%
and such that there exists an exponentially decreasing path of $\CZ$ cubes 
\[
Q=Q_1 \leftrightarrow Q_2 \leftrightarrow \ldots \leftrightarrow
Q_{L}=Q^{\prime },
\]%
with $L \leq C( K)$, and 
\begin{equation}
\label{exp_dec_ineq}
\delta _{Q_{l}} \leq C(K) \cdot ( 1-c( K) ) ^{l-l^{\prime }} \cdot \delta_{Q_{l^{\prime}}} \qquad \text{for} \;\; 1 \leq l^{\prime }\leq l\leq L.
\end{equation}
\end{enumerate}

The work, space and number of calls to the $\CZ$-\textsc{Oracle} required by the algorithm are bounded by a constant $C( K)$.

\begin{proof}[\underline{Explanation}] 

Since $K$ is an odd integer, we can partition $KQ$ into $K^{n}$ many dyadic cubes $\widetilde{Q%
}_{\nu }$, each of sidelength $\delta _{Q}$. For each $\widetilde{Q}_{\nu }$, we apply the $\CZ$-\textsc{Oracle} to determine $Q_{\nu }^{CZ}$, the $\CZ$ cube containing the center of $\widetilde{Q}_{\nu }$.

The cube $Q$ is keystone if and only if $Q_{\nu }^{CZ}\supseteq \widetilde{%
Q}_{\nu }$ for each $\nu $. Thus, we can test whether $Q$
is a keystone, using work at most $C( K)$ and using at most $C( K) $ calls to the $\CZ$-\textsc{Oracle}. If $Q$ is a keystone, then we are done.

Suppose $Q$ is not a keystone cube. We claim that there exists a path $Q=Q_1 \leftrightarrow Q_2 \leftrightarrow \ldots \leftrightarrow Q_{L}$ as in (KEY 2). This claim was essentially proven in Lemma 6.12 of \cite{FIL1}. The main difference is that the keystone cubes in \cite{FIL1} are defined with $K=100$, whereas here $K$ is an odd integer of size at least $10$ (to be fixed later). By making superficial modifications to the argument in Lemma 6.12 of \cite{FIL1} we prove our claim in the present setting.

To find an exponentially decreasing path as in (KEY 2), we first enumerate all the paths of $\CZ$ cubes $Q=Q_1 \leftrightarrow Q_2 \leftrightarrow \ldots
\leftrightarrow Q_{L}$ with $L\leq C( K) $, and we then ``test'' each path to see whether it satisfies the necessary
conditions: $\delta_{Q_L} \leq (1/2) \delta_Q$, $Q_L \cap K Q \neq \emptyset$, and \eqref{exp_dec_ineq}.

There are at most $C( K) $ such paths, and we can generate them
using work, space and calls to the $\CZ$-\textsc{Oracle} at most $C( K) $,
because, given a cube $\widetilde{Q}\in \CZ$, we can determine all the cubes $%
\widetilde{Q}^{\prime }\in \CZ$ such that $\widetilde{Q}\leftrightarrow \widetilde{Q}%
^{\prime }$, by using at most $C$ calls to the $\CZ$-\textsc{Oracle}. (Just query the $\CZ$-\textsc{Oracle} using as $\underline{x}$ the center of each dyadic cube $\widetilde{Q}%
^{\prime }$ such that $\widetilde{Q}\leftrightarrow \widetilde{Q}^{\prime }$ and $%
\frac{1}{64}\delta _{\widetilde{Q}}\leq \delta _{\widetilde{Q}^{\prime }}\leq 64 \delta
_{\widetilde{Q}}$.)

We can  ``test'' a given path using work and storage at most $C(K)$, and using no calls to the $\CZ$-\textsc{Oracle}.

This concludes the explanation of the algorithm \textsc{Keystone-or-Not}.

\end{proof}

\environmentA{Algorithm: List All Keystone Cubes.}

We produce a list $Q_{1}^{\#},\ldots ,Q_{L^{\#}}^{\#}$, consisting of all
the keystone cubes in $\CZ$. Each keystone cube appears once and only once in our
list. We have $L^{\#}\leq C( K) N$. The algorithm uses work at most $
C( K) N\log N$ in space $C( K) N$, and at most $%
C( K) N$ calls to the $\CZ$-\textsc{Oracle}.

\begin{proof}[\underline{Explanation}] 

Let $Q^{\#} \in \CZ$ be a keystone cube. Since $Q^\# \in \CZ$, there exists a point $x\in 9Q^{\#}\cap E$ due to our assumption that ``$E$ is nearby''. Let $Q_{x}$ be the $\CZ$ cube containing $x$. Then $\delta_{Q_{x}}\geq \delta _{Q^{\#}}$ because $Q^{\#}$ is keystone; moreover, $\delta _{Q_{x}}\leq C\delta _{Q^{\#}}$ because of good geometry. Hence, for each keystone cube $Q^{\#}$ there exists $x \in E$ such that
\begin{equation}
x \in 9Q^{\#} \text{ and }c\delta _{Q_{x}}\leq \delta _{Q^{\#}}\leq
\delta _{Q_{x}}.   \label{ptk10}
\end{equation}

To generate all the keystone cubes we may therefore proceed as follows.

We loop over all $x \in E$. For each $x \in E$, we produce the unique $\CZ$ cube $Q_{x}$ containing $x$, using the $\CZ$-\textsc{Oracle}. We then
list all the dyadic cubes $Q^{\#}$ satisfying \eqref{ptk10} (there are at most $C$ such $Q^\#$ for a fixed $x$). Finally, we apply the algorithm \textsc{Keystone-or-Not} to test each $Q^{\#}$ to see whether it is a keystone cube. We discard the cubes $Q^\#$ that are not keystone and retain the remaining cubes. Clearly, a single iteration of the loop requires at most $C(K)$ calls to the $\CZ$-\textsc{Oracle} and additional work at most $C(K)$.

Thus, with work and storage at most $C(K)N$, and with at most $%
C(K) N $ calls to the $\CZ$-\textsc{Oracle}, we produce a list $%
Q_{1}^{\#},\ldots ,Q_{L^{\#}}^{\#}$, with $
L^{\#}\leq C(K) N$, consisting of all the keystone cubes,
but possibly containing multiple copies of the same cube.

With work at most $C(K) N\log N$ in space $C(K) N$, we can
sort the list $Q_{1}^{\#},\ldots ,Q_{L^{\#}}^{\#}$ and remove duplicates.

This completes our explanation of the algorithm \textsc{List All Keystone
Cubes}.

\end{proof}

\bigskip 

\environmentA{Algorithm: Make Auxiliary Cubes.}

For each $Q_{l}^{\text{CD}}$, $S_{l}=3Q_{l}^{\text{CD}}\cap E$, produced by the algorithm 
\textsc{Make Cluster Descriptors}, we compute a point $x_{l}^{\extra}\in
H( S_{l}) $ such that 
\begin{equation} \label{extrapoint}
2A\cdot \diam( S_{l}) \leq \left\vert x_{l}^{\extra}-x(
S_{l}) \right\vert \leq 8A\cdot \diam( S_{l}) 
\end{equation}
and we compute $Q_{l}^{\extra}$, the $\CZ$ cube containing $x_{l}^{\extra}$.

The algorithm uses work $\leq C( A) N \log N$ in space $C(A) N$, and makes at most $C( A) N$ calls to the $\CZ$-\textsc{Oracle}.

\begin{proof}[\underline{Explanation}]  

For each $1 \leq l \leq L$ we compute $\diam(S_l)$; see Remark \ref{bbd_rem}. We choose $x \in \R^n$ satisfying
\[
2A\cdot \diam( S_{l}) \leq \left\vert x -x(
S_{l}) \right\vert \leq 8A\cdot \diam( S_{l}). 
\]
Note that we necessarily have $x \in H( S_l ) $. After picking such a point  $%
x= x_{l}^{\extra}$, we call the $\CZ$-\textsc{Oracle} to determine $Q_l^{\extra}$.

Recall that there are at most $CN$ distinct indices $l$; see the algorithm \textsc{Make Cluster Descriptors}. Thus, in the present algorithm, the work, storage and number of calls to the $\CZ$-\textsc{Oracle} are bounded as required.

\end{proof}

From Lemma \ref{lem_cluster1} and the definition of $Q_{l}^{\extra}$ ($1 \leq l \leq L$), we obtain
the following:

\begin{lem}\label{lem_cluster2} Assume that $n \geq 2$. Let $Q\in \CZ$ and $l \in \{1,\cdots,L\}$, and suppose that 
\begin{equation} \label{relevantcube} cA^{10}\cdot \diam( S_l ) <\left\vert x-x(S_l) \right\vert < CA^{-10}\cdot \dist( S_l ,E\setminus
S_l) \text{ for all }x\in ( 1+c_{G}) Q \text{.}
\end{equation}
Then there exists an exponentially decreasing path $\mathcal{S} = (Q_1,\cdots,Q_{\overline{\overline{J}}})$ joining $Q$ to $Q^{\extra}_l$.
\end{lem}

\begin{proof}

Denote the point $x^{\extra} = x^{\extra}_l$, the cube $Q^{\extra} = Q^{\extra}_l$, and the cluster $S= S_l$. By the conditions in the algorithm \textsc{Make Auxiliary Cubes}, we have $x^{\extra} \in Q^{\extra} \cap H(S)$.

We will construct an exponentially decreasing path $\mathcal{S} = (\widehat{Q}_1,\cdots,\widehat{Q}_{\overline{\overline{J}}})$ that joins $Q$ to $Q^\extra$.

Let $x \in Q$. From \eqref{relevantcube}, we see that $x \in Q \cap H(S)$.

From \eqref{extrapoint} and \eqref{relevantcube} we see that
\[ \lvert x^{\extra} - x(S) \rvert  \leq 8 A \diam(S) \leq c A^{10} \diam(S) \leq \lvert x - x(S) \rvert. \]
Lemma \ref{path_lem} implies that there exists a sequence of points $x_1,\cdots,x_J \in H(S)$ and an integer $J_* \geq 1$, satisfying the following bullet points.

\begin{itemize}
\item $x_1 = x$ and $x_J = x^{\extra}$.
\item $\lvert x_{j+1} - x(S) \rvert \leq \lvert x_j - x(S) \rvert$ for $j = 1,\cdots, J-1$.
\item $\lvert x_j - x_{j+1} \rvert \leq A^{-2} \lvert x_j - x(S) \rvert$ for $j = 1,\cdots , J-1$.
\item $\lvert x_{j + J_*} - x(S) \rvert \leq (1-A^{-3}) \lvert x_j - x(S) \rvert$ for $1 \leq j \leq J - J_*$.
\item $J_* \leq A^3$.
\end{itemize}
(Our assumption that $n \geq 2$ implies that the halo $H(S_l)$ has a single connected component. When $n=1$ we cannot use Lemma \ref{path_lem}, hence we will have to modify our approach.)

Since $\lvert x_j - x(S) \rvert$ is a non-increasing sequence, \eqref{extrapoint} and \eqref{relevantcube} imply that
\begin{align*}
2 A \diam(S) & \leq \lvert x^{\extra} - x(S) \rvert = \lvert x_J - x(S) \rvert \leq \lvert x_j - x(S) \rvert \leq \lvert x_1 - x(S) \rvert \\
&= \lvert x - x (S) \rvert \leq C A^{-10} \dist(S,E \setminus S) \leq A^{-1} \dist(S,E \setminus S) \quad \mbox{for} \; j=1,\cdots,J.
\end{align*}
Hence, $x_j \in H(S)$ for each $j=1,\cdots,J$.

Let $Q_j$ ($1 \leq j \leq J$) denote the $\CZ$ cube containing $x_j$. Thus $Q_1 = Q$ and $Q_J = Q^{\extra}$. (Recall that $x_1 = x \in Q$ and $x_J = x^{\extra} \in Q^{\extra}$.) 

By the third bullet point and by Lemma \ref{lem_cluster1a} we have
\begin{equation}\label{joins}
Q_j \leftrightarrow Q_{j+1} \; \mbox{for} \; j=1,\cdots,J-1.
\end{equation} 
Hence, $\delta_{Q_{j+1}}$ and $\delta_{Q_j}$ differ by at most a factor of $64$, thanks to good geometry.

Recall that $c_G > 0$ is a universal constant satisfying \eqref{ggz}. We now prove the following

\noindent\textbf{Claim:} If $1 \leq J_0 \leq J$ satisfies
\begin{equation}
\label{ss0}
\lvert x_{J_0} - x (S) \rvert \geq c_G \cdot \delta_{Q^{CZ}(S)},
\end{equation}
then
\begin{equation} \label{ss1}
\delta_{Q_{j'}} \leq C(A) \cdot (1-c(A))^{j' - j} \delta_{Q_j} \; \mbox{for all} \; 1 \leq j \leq j' \leq J_0.
\end{equation}

\noindent\textbf{Proof of Claim:} Since the sequence $\lvert x_j - x(S) \rvert$ ($1 \leq j \leq J$) is non-increasing, \eqref{ss0} implies that $\lvert x_{j} - x (S) \rvert \geq c_G \cdot \delta_{Q^{CZ}(S)}$ for all $1 \leq j \leq J_0$. Thus, Lemma \ref{lem_cluster1} implies that
\[c \cdot \lvert x_j - x(S) \rvert \leq \delta_{Q_j} \leq C \cdot \lvert x_j - x(S) \rvert \; \mbox{for} \; 1 \leq j \leq J_0.\]
This estimate and the fourth bullet point written above imply that
\[\delta_{Q_{J_* k}} \leq C \cdot (1-A^{-3})^{k-\ell} \cdot \delta_{Q_{J_* \ell}} \; \mbox{for} \; 1 \leq J_* \ell \leq J_* k \leq J_0.\]
Since $\delta_{Q_{j+1}}$ and $\delta_{Q_j}$ differ by at most a factor of $64$, and since $1 \leq J_* \leq A^3$, we obtain \eqref{ss1}. This completes the proof of the claim.

Since $\lvert x_j - x(S) \rvert$ ($1 \leq j \leq J$) is non-increasing, the following cases are exhaustive.

\begin{description}
\item[Case 1] $\lvert x_J - x(S) \rvert \geq c_G \cdot \delta_{Q^{CZ}(S)}$.
\item[Case 2] There exists $1 \leq \overline{J} \leq J-1$ such that
\begin{itemize}
\item $\lvert x_{\overline{J} + 1} - x(S) \rvert  < c_G \cdot \delta_{Q^{CZ}(S)}$, and
\item $\lvert x_{\overline{J}} - x(S) \rvert \geq c_G \cdot \delta_{Q^{CZ}(S)}$.
\end{itemize}
\item[Case 3] $\lvert x_1 - x(S) \rvert < c_G \cdot \delta_{Q^{CZ}(S)}$.
\end{description}

\underline{First, we consider \textbf{Case 1}.} In this case, \eqref{ss0} holds with $J_0 = J$. Hence, \eqref{ss1} implies that
\begin{equation*}
\delta_{Q_{j'}} \leq C(A) \cdot (1-c(A))^{j' - j} \delta_{Q_j} \; \mbox{for} \; 1 \leq j \leq j' \leq J.
\end{equation*}
We define the sequence $\mathcal{S} := (Q_1,\cdots,Q_J)$. The above estimate and \eqref{joins} show that $\mathcal{S}$ is an exponentially decreasing path joining $Q$ to $Q^{\extra}$.

\underline{Next, we consider \textbf{Case 2}.} In this case, \eqref{ss0} holds with $J_0 = \overline{J}$. Hence, \eqref{ss1} implies that
\begin{equation}\label{ss2}
\delta_{Q_{j'}} \leq C(A) \cdot (1-c(A))^{j' - j} \delta_{Q_j} \; \mbox{for} \; 1 \leq j \leq j' \leq \overline{J}.
\end{equation}
Recall that $x(S) \in Q^{CZ}(S)$, and $\lvert x_{\overline{J} + 1} - x(S) \rvert  < c_G \cdot \delta_{Q^{CZ}(S)}$. Hence, $x_{\overline{J}+1} \in (1 + 8c_G) Q^{CZ}(S)$. Since also $x_{\overline{J}+1} \in Q_{\overline{J} + 1}$ we see that $(1+8c_G) Q_{\overline{J}+1} \cap (1+8c_G) Q^{CZ}(S) \neq \emptyset$, hence $Q_{\overline{J} + 1} \leftrightarrow Q^{CZ}(S)$ thanks to \eqref{ggz}.

Moreover, since $\lvert x_j - x(S) \rvert$ is non-increasing, we have
\[ \lvert x^{\extra} - x(S) \rvert = \lvert x_J - x(S) \rvert \leq \lvert x_{\overline{J} + 1} - x(S) \rvert \leq c_G \cdot \delta_{Q^{CZ}(S)}.\]
Recall that $x(S) \in Q^{CZ}(S)$. Hence, the above estimate shows that $x^{\extra} \in (1+8c_G) Q^{CZ}(S)$. Since also $x^{\extra} \in Q^{\extra}$ we see that $(1+8c_G) Q^{CZ}(S) \cap (1+8c_G) Q^{\extra} \neq \emptyset$, hence $Q^{CZ}(S) \leftrightarrow Q^{\extra}$ thanks to \eqref{ggz}.

We define the sequence $\mathcal{S} := (\widehat{Q}_1,\widehat{Q}_2, \cdots,\widehat{Q}_{\overline{\overline{J}}}) := (Q_1, Q_2, \cdots, Q_{\overline{J}}, Q_{\overline{J} + 1}, Q^{CZ}(S), Q^{\extra})$. From \eqref{joins}, we see that $\widehat{Q}_j \leftrightarrow \widehat{Q}_{j+1}$ ($1 \leq j \leq \overline{\overline{J}} - 1$), hence $\delta_{\widehat{Q}_j}$ and $\delta_{\widehat{Q}_{j+1}}$ differ by at most a factor of $64$, thanks to good geometry. Combined with \eqref{ss2}, this shows that $\mathcal{S}$ is an exponentially decreasing path joining $Q$ to $Q^\extra$. 

\underline{Lastly, we consider \textbf{Case 3}.} In this case, $Q_1 \leftrightarrow Q^{CZ}(S)$ and $Q^{CZ}(S) \leftrightarrow Q^{\extra}$ as in the discussion of \textbf{Case 2}. Hence, the exponentially decreasing path $\mathcal{S} := (Q_1, Q^{CZ}(S), Q^{\extra})$ joins $Q$ to $Q^{\extra}$.

This completes the proof of Lemma \ref{lem_cluster2}.

\end{proof}

A cube $Q \in \CZ$ is called \underline{interstellar} provided that
\[
(1+3c_G)Q \cap E = \emptyset \;\; \mbox{and} \;\; \diam( A^{10}Q\cap E) \leq A^{-10}\delta_{Q}\text{.}
\]
Any cube $Q \in \CZ$ that is not interstellar will be called \underline{non-interstellar}.

\environmentA{Algorithm: Test an Interstellar Cube.}

We perform one-time work at most $C(A) N \log N$ in space $C(A) N$, after which we can answer queries.

A query consists of a cube $Q \in \CZ$.

The response to the query $Q$ is as follows. We first determine whether $Q$ is
interstellar. If it is, then we find an index $1 \leq l \leq L$ such that $S=S_{l}=3Q_{l}^{\text{CD}}\cap E
$ satisfies  
\begin{equation}
cA^{10}\cdot \diam( S) <\left\vert x-x( S) \right\vert  < C A^{-10}\cdot \dist( S, E \setminus S),
\label{ptk15}
\end{equation}
for all $x \in (1+c_G)Q$. The query work is at most $C(A) \log N$.

\begin{proof}[\underline{Explanation}]

We determine whether $Q$ is interstellar by computing $\diam( A^{10}Q\cap E)$ and by testing whether $(1+ 3c_{G}) Q\cap E=\emptyset$ using the BBD tree. We compute $\diam( A^{10}Q\cap E)$ using the algorithm \textsc{RCZ} in Section \ref{sec_bbd}. We use that $A^{10} Q$ can be expressed as the disjoint union of finitely many dyadic cubes of sidelength $\delta_Q/2$; see Remark \ref{bbd_rem}. Similarly, we  test whether $(
1+ 3c_{G}) Q\cap E=\emptyset $ using the BBD tree and the fact that $(1+3c_G)Q$ can be expressed as the disjoint union of finitely many dyadic cubes of sidelength $c_G \delta_Q/2$; again, see Remark \ref{bbd_rem}. This computation requires work and storage at most $C(A) \log N$.

The proof of Lemma 6.3 in \cite{FIL1} shows that if $Q$ is
interstellar, then for some cluster $S$ we have (\ref{ptk15}) for all $x\in
( 1+c_{G}) Q$. 

It follows that $S$ is a strong cluster and therefore $S$ is among the clusters $S_{l}
$ produced by the algorithm \textsc{Make Cluster Descriptors}. Hence, for any $x\in ( 1+c_{G}) Q$ we have $S=3\hat{Q}\cap E$%
, where $\hat{Q}$ is a dyadic cube such that $ \delta_{\hQ} \in \left[ 8\cdot \dist(x,E), 64\cdot \dist( x,E) \right]$ and $x\in \hQ$.

Therefore, given an interstellar $Q$, we can learn which $S_{l}$ satisfies (%
\ref{ptk15}) as follows.

\begin{itemize}
\item Let $x$ be the center of $Q$.

\item Compute $\dist(x,E) $ up to a factor of $2$, using the BBD tree.

\item Compute $\hat{Q}$, a dyadic cube containing $x$, with $8\cdot $dist$%
( x,E) \leq \delta _{\hat{Q}}\leq 64\cdot $dist$( x,E) 
$.

\item Using the algorithm \textsc{Find Descriptor Cube}, compute the
descriptor cube $Q^{\text{CD}}=\DC( 3\hat{Q}\cap E) $.

\item We know that $Q^{\text{CD}}$ will be one of our cubes $Q_{l}^{\text{CD}}$ produced
by the algorithm \textsc{Make Cluster Descriptors}. By a binary search, we
find $l$ such that $Q^{\text{CD}}=Q_{l}^{\text{CD}}$.

\item Thus, we find the $l$ for which $S=S_{l}=3Q_{l}^{\text{CD}}\cap E$ satisfies  (\ref%
{ptk15}).

\end{itemize}

This process takes work $\leq C( A) \log N$ and uses at most $C(A)$ calls to the $\CZ$-\textsc{Oracle}.

\end{proof}

\environmentA{Algorithm: List All Non-Interstellar Cubes.}

Using work at most $C(A) N \log N$ in space $C(A) N$ and making at most $C(A)N$ calls to the $\CZ$-\textsc{Oracle}, we produce all the non-interstellar cubes $Q \in \CZ$.

\begin{proof}[\underline{Explanation}] 

We compute representatives $(x_\nu',x_\nu'') \in E \times E$ ($\nu=1,\cdots,\nu_{\max}$) arising in the WSPD. These representatives have the property that, for each $(x',x'') \in E \times E \setminus \{ (x,x) : x \in E \}$, there exists $\nu$ such that 
\[\lvert x_\nu' - x' \rvert + \lvert x_\nu'' - x'' \rvert \leq \frac{1}{100} \lvert x' - x'' \rvert,\]
and $\nu_{\max} \leq C N$.

Let $Q \in \CZ$ be non-interstellar. Then either $(1+3 c_G) Q \cap E \neq \emptyset$ or $\diam(A^{10}Q \cap E) > A^{-10} \delta_Q$. In the second case, there exist two points $x',x'' \in A^{10}Q \cap E$  with $\lvert x' - x'' \rvert > A^{-10} \delta_Q$. Hence, in the second case there exists $\nu$ such that $x_\nu', x_\nu'' \in A^{11}Q \cap E$ and $\lvert x_\nu' - x_\nu'' \rvert > A^{-11} \delta_Q$.

We have shown the following: If $Q \in \CZ$ is non-interstellar, then
\[\exists x \in E \; \mbox{s.t.} \; x \in (1+3c_G)Q \;\; \mbox{\underline{or}} \; \; \exists \nu \; \mbox{s.t.} \; x_\nu',x_\nu'' \in A^{11}Q \cap E \; \mbox{and} \; \lvert x_\nu' - x_\nu'' \rvert > A^{-11} \delta_Q.\]

Thus, to list all the non-interstellar cubes we can proceed as follows:

For each $x \in E$, find all the cubes $Q \in \CZ$ such that $(1+3c_G)Q$ contains $x$. There are at most $C$ such cubes for each fixed $x$. We produce these cubes by making at most $C$ calls to the $\CZ$-\textsc{Oracle}.

For each $\nu =1,\cdots,\nu_{\max}$, find all the dyadic cubes $Q$ such that $x_\nu',x_\nu'' \in A^{11}Q \cap E$ and $\delta_Q \leq A^{11} \lvert x_\nu' - x_\nu'' \rvert$. There are at most $C(A)$ such cubes for each fixed $\nu$.

We have produced a list that contains all the non-interstellar cubes, and consists of at most $C(A) N$ dyadic cubes. We now pass through this list and remove any cubes that are interstellar. We then sort the remaining cubes and remove duplicates.

Thus we have computed the list of all non-interstellar cubes.

The reader may easily check that our algorithm performs as promised in terms of work, storage, and calls to the $\CZ$-\textsc{Oracle}.

\end{proof}

We now create a list \underline{USUAL-SUSPECTS} \label{pp1}, consisting of all keystone
cubes and all non-interstellar $\CZ$ cubes, and all the cubes $Q_{l}^{\extra}$
produced by the algorithm \textsc{Make Auxiliary Cubes}.

There are at most $C( A) N$ cubes in the list USUAL-SUSPECTS,
and we can produce the list using work $\leq C( A) N\log N$,
storage $\leq C( A) N$, and at most $C( A) N $ calls to
the $\CZ$-\textsc{Oracle}.

We assume that our list USUAL-SUSPECTS is sorted so that $Q,Q^{\prime }\in $
USUAL-SUSPECTS and $\delta _{Q}<\delta _{Q^{\prime }}\Rightarrow Q$ precedes 
$Q^{\prime }$ in the list USUAL-SUSPECTS. Thus, all the smallest cubes are located at the beginning of the list. We may also assume that the list USUAL-SUSPECTS contains no duplicates.

This can be achieved by doing extra work at most $C( A) N\log N$ in space 
$C( A) N$.

\bigskip 

In preparation for the next two algorithms, we prove a small lemma.

\begin{lem}\label{USUAL-SUSPECTS}
Fix constants $A^{\#}\geq 1$ and $0<a^{\#}<1$. Suppose we are given finite
sequences 
\begin{eqnarray*}
\mathcal{S}^{\left[ 1\right] } &:&\delta _{1}^{\left[ 1\right] },\delta
_{2}^{\left[ 1\right] },\cdots ,\delta _{L^{\left[ 1\right] }}^{\left[ 1%
\right] } \\
&&\vdots  \\
\mathcal{S}^{\left[ M\right] } &:&\delta_{1}^{\left[ M\right] },\delta _{2}^{\left[ M%
\right] },\cdots ,\delta _{L^{\left[ M\right] }}^{\left[ M\right] }
\end{eqnarray*}%
of positive real numbers. Assume 

\begin{itemize}
\item $\delta _{\ell}^{\left[ k\right] }\leq A^{\#}\cdot ( 1-a^{\#})
^{\ell - \ell^{\prime }}\delta _{\ell^{\prime }}^{\left[ k\right] }$ \quad for $1 \leq k \leq M$  \quad  and $1\leq \ell^{\prime }\leq \ell \leq L^{\left[ k\right] }$.

\item $\delta _{L^{\left[ k\right] }}^{\left[ k\right] }\leq (
1-a^{\#}) ^{L^{\left[ k\right]} - 1  }\delta _{1}^{\left[ k\right] }$  \quad  for $%
1 \leq k \leq M$.

\item $\delta_{1}^{\left[ k+1\right] }=\delta _{L^{\left[ k\right] }}^{%
\left[ k\right] }$  \quad  for $1 \leq k \leq M-1$.
\end{itemize}

Then the sequence 
\begin{eqnarray*}
&&\left( \delta _{1}^{\left[ 1\right] },\cdots ,\delta _{L^{\left[ 1\right]}}^{\left[ 1\right] },\delta_{2}^{\left[ 2\right] },\cdots ,\delta _{L^{\left[ 2\right] }}^{\left[ 2\right] },\delta _{2}^{\left[ 3\right] },\cdots,\delta _{L^{\left[ 3\right] }}^{\left[ 3\right] },\ldots, \delta_{L^{\left[ M-1 \right]}}^{\left[M-1\right]},  \delta_{2}^{\left[ M\right] },\cdots ,\delta _{L^{\left[ M\right] }}^{\left[ M\right]} \right)\\
&\equiv &\left( \delta _{1},\delta _{2},\ldots ,\delta_{J} \right) 
\end{eqnarray*}%
satisfies 
\[
\delta _{j}\leq (A^{\#})^2 \cdot ( 1-a^{\#}) ^{j-j^{\prime }}\delta
_{j^{\prime }}  \quad \text{ for } 1\leq j^{\prime }\leq j\leq J\text{.}
\]
\end{lem}

\begin{proof}
Let $1 \leq j' \leq j \leq J$. 

First, suppose that $\delta_j = \delta_\ell^{\left[k \right]}$ and $\delta_{j'} = \delta_{\ell'}^{\left[ k \right]}$ with $1 \leq k \leq M$ and $1 \leq \ell' \leq \ell \leq L^{\left[ k \right]}$. Then we have
\begin{align*}
\delta_j = \delta_{\ell}^{\left[ k\right] } & \leq A^{\#}\cdot ( 1-a^{\#})
^{\ell-\ell^{\prime }}\delta _{\ell^{\prime }}^{\left[ k\right] }\\ 
& = A^{\#} \cdot ( 1-a^{\#})^{j - j^\prime} \delta _{j^\prime}.
\end{align*}

We consider the remaining case. Namely, suppose that $\delta_j = \delta_\ell^{\left[k \right]}$ and $\delta_{j'} = \delta_{\ell'}^{\left[ k' \right]}$ with $1 \leq k' < k \leq M$, $1 \leq \ell \leq L^{\left[ k \right]}$, and $1 \leq \ell' \leq L^{\left[ k' \right]}$. Then we have
\begin{align*}
\delta_j = \delta_\ell^{\left[ k \right]} & \leq A^\# \cdot (1- a^\#)^{\ell - 1} \delta_1^{\left[ k \right]} \\
& = A^\# \cdot (1-a^\#)^{\ell -1} \delta_{L^{\left[ K-1 \right]}}^{\left[ k - 1 \right]} \\
& \leq  A^\# \cdot (1-a^\#)^{\ell -1}  (1-a^\#)^{L^{\left[k-1 \right] - 1}} \delta_1^{\left[ k -1 \right] }.
\end{align*}
Now, iterating the above reasoning we obtain
\begin{align*}
\delta_j = \delta_\ell^{\left[ k \right]} & \leq A^\# \cdot (1-a^\#)^{\ell -1}  (1-a^\#)^{L^{\left[k-1 \right] - 1}} \cdots (1-a^\#)^{L^{[k' +1]} - 1 } \delta_1^{\left[ k' + 1 \right] } \\
& = A^\# \cdot (1-a^\#)^{\ell -1}  (1-a^\#)^{L^{\left[k-1 \right] - 1}} \cdots (1-a^\#)^{L^{[k' +1]} - 1 } \delta_{L^{\left[k' \right]}}^{\left[ k' \right] } \\
& \leq (A^\#)^2 \cdot (1-a^\#)^{\ell -1}  (1-a^\#)^{L^{\left[k-1 \right] - 1}} \cdots (1-a^\#)^{L^{[k' +1]} - 1 } \cdot (1 - a^\#)^{L^{\left[ k' \right]} - \ell '} \delta_{\ell'}^{\left[ k' \right] } \\
& =  (A^\#)^2 \cdot (1-a^\#)^{j - j'}  \delta_{j'}.
\end{align*}
This completes the proof of the lemma.
\end{proof}

\environmentA{Algorithm: Mark Usual Suspects.}

We mark each cube $Q$ appearing in the list USUAL-SUSPECTS with a keystone cube 
$\mathcal{K}( Q) $ such that $Q$ is joined to $\mathcal{K}(Q)$ by an exponentially decreasing path. That is, there exists a finite sequence of $CZ
$ cubes 
\[
Q=Q_{1}\leftrightarrow Q_{2}\leftrightarrow \cdots \leftrightarrow
Q_{L( Q) }=\mathcal{K}( Q) 
\]%
such that 
\[
\delta _{Q_{\ell'}}\leq C( A) \cdot ( 1-c( A) )
^{\ell' - \ell}\delta _{Q_{\ell }}\text{ for } 1 \leq \ell \leq
\ell' \leq L( Q) .
\]%
We do not compute $Q_{1},\cdots ,Q_{L( Q) -1}$, but we guarantee
that they exist.

If $Q$ is keystone, then we guarantee that $\mathcal{K}( Q) =Q$.

The algorithm does work at most $C( A) N\log N$ in space $
C( A) N$, and it makes at most $C( A) N$ calls to the $\CZ$-\textsc{Oracle}.

\begin{proof}[\underline{Explanation}]

For each $Q$ in the list USUAL-SUSPECTS we will compute a keystone cube $\mathcal{K}(Q) \in \CZ$ such that we guarantee that there exists a list of sequences of $\CZ$ cubes
\begin{align*}
&\mathcal{S}_1 = \left(Q^{[1]}_1,\cdots,Q^{[1]}_{L^{[1]}} \right),  \\
&\mathcal{S}_2 = \left(Q^{[2]}_1,\cdots,Q^{[2]}_{L^{[2]}} \right), \\
&\qquad\qquad\qquad \vdots\\
&\mathcal{S}_M = \left(Q^{[M]}_1,\cdots,Q^{[M]}_{L^{[M]}} \right),
\end{align*}
with the following properties.
\begin{itemize}
\item The initial cube of $\mathcal{S}_1$ is $Q$ and the terminal cube of $\mathcal{S}_M$ is $\mathcal{K}(Q)$, i.e.,
\[
Q^{[1]}_1 = Q \;\; \mbox{and} \;\; Q^{[M]}_{L^{[M]}} = \mathcal{K}(Q).
\]
\item The terminal cube of a sequence matches the initial cube of its successor:
\[
Q^{[k]}_{L^{[k]}} = Q^{[k+1]}_1 \;\; \mbox{for} \;\; 1 \leq k \leq M-1.
\]
\item Each sequence is connected and exponentially decreasing:
\begin{equation*}
Q^{[k]}_{1}\leftrightarrow Q^{[k]}_{2}\leftrightarrow \cdots \leftrightarrow
Q^{[k]}_{L^{[k]}}
\end{equation*}
and
\begin{equation*}
\delta _{Q^{[k]}_{\ell}}\leq C_\# \cdot (1-c_\#)^{\ell-\ell^{\prime }}\delta _{Q^{[k]}_{\ell'}} \;\; \text{for}\; 1 \leq \ell^{\prime}  \leq \ell \leq L^{[k]}.
\end{equation*}
Moreover, we guarantee that
\begin{equation*}
\delta_{Q^{[k]}_{L^{[k]}}} \leq (1-c_\#)^{L^{[k]}-1} \cdot \delta_{Q^{[k]}_1}.
\end{equation*}
Here, $c_\# \in (0,1)$ and $C_\# \geq 1$ are controlled constants. In this discussion, a ``controlled constant'' is a constant that depends only on $A$, $K$, and $n$.
\end{itemize}

If these conditions hold, then we say that $Q$ and $\mathcal{K}(Q)$ are \emph{connected by a chain} $\mathcal{S}_1,\cdots,\mathcal{S}_M$ with constants $c_\#$ and $C_\#$.

By concatenating the sequences $\mathcal{S}_1,\cdots,\mathcal{S}_M$ we obtain a sequence of $\CZ$ cubes as in the algorithm \textsc{Mark Usual Suspects}, with $C(A) = (C_\#)^2$ and $c(A) = c_\#$. See Lemma \ref{USUAL-SUSPECTS}. Thus it suffices to compute a keystone cube $\mathcal{K}(Q)$ and verify the existence of a suitable chain for each $Q$ in the list USUAL-SUSPECTS.

Recall that the cubes in USUAL-SUSPECTS are sorted according to their size. We loop through all the $Q$ in USUAL-SUSPECTS, starting with the smallest cubes at the beginning of the list. We will compute $\mathcal{K}(Q)$ in the body of the loop, which is presented below.

We fix $Q$ in USUAL-SUSPECTS. By induction, we may assume that for each $Q'$ in USUAL-SUSPECTS with $\delta_{Q'} < \delta_Q$ we have computed a keystone cube $\mathcal{K}(Q')$ to which $Q'$ is connected by a chain with  constants $c_\#$ and $C_\#$.

We assume that $c_\#$ is less than a small enough controlled constant, and that $C_\#$ is greater than a large enough controlled constant. We will later pick $c_\#$ and $C_\#$ to be controlled constants, but not yet.

We perform the following procedure.

\noindent \underline{Main Procedure:}

\begin{itemize}
\item We initialize $Q^{[1]} = Q$.
\item Let $M_{\max}$ be a large enough integer determined by $A$, $K$, and $n$, to be picked later.
\item We perform the following loop: \\
For ($k=1,\cdots, M_{\max}  - 1$) \\
$\bigl\{$
\begin{itemize}
\item We execute the algorithm \textsc{Keystone-or-Not} to produce one of two outcomes. \\
\textbf{(Outcome A)} We guarantee that $Q^{[k]}$ is a keystone cube. We then return the cube $Q_{\text{out}} = Q^{[k]}$, indicating that it is a keystone cube, and terminate the loop.  \\
\textbf{(Outcome B)} We witness that $Q^{[k]}$ fails to be a keystone cube: We compute a cube $Q^{[k+1]} \in \CZ$ with $\delta_{Q^{[k+1]}} \leq \frac{1}{2} \delta_{Q^{[k]}}$ and $Q^{[k+1]} \cap K Q^{[k]} \neq \emptyset$, such that there exists a sequence of $\CZ$ cubes $\mathcal{S}_k = (Q^{[k]}_{1},  Q^{[k]}_{2},  \cdots, Q^{[k]}_{L^{[k]}} )$, with $L^{[k]} \leq C(K)$ and
\begin{equation*}
Q^{[k]}=Q^{[k]}_{1}\leftrightarrow Q^{[k]}_{2}\leftrightarrow \cdots \leftrightarrow
Q^{[k]}_{L^{[k]}}= Q^{[k+1]}
\end{equation*}
such that 
\begin{equation*}
\delta _{Q^{[k]}_{\ell}}\leq C( A) \cdot (1-c( A))^{\ell-\ell^{\prime }}\delta _{Q^{[k]}_{\ell'}} \;\; \text{for}\; 1 \leq \ell^{\prime}  \leq \ell \leq L^{[k]}.
\end{equation*}
Combining the estimate $\delta_{Q^{[k]}_{L^{[k]}}} \leq \frac{1}{2} \delta_{Q^{[k]}_1}$ with our bound on $L^{[k]}$, we see that
\begin{equation*}
\delta_{Q^{[k]}_{L^{[k]}}} \leq (1-c(A))^{L^{[k]}-1} \cdot \delta_{Q^{[k]}_1}.
\end{equation*}
Here, $c(A)$ and $C(A)$ are  controlled constants.

\item If $k=M_{\max} -1$ then we return the cube $Q_{\text{out}} = Q^{[M_{\max} ]}$. Using the algorithm \textsc{Keystone-or-Not}, we determine whether $Q_{\text{out}}$ is a keystone cube, and after indicating the result to the user we terminate the loop.
\end{itemize}
$\bigr\}$
\end{itemize}

We will now analyze the output of the Main Procedure.

Suppose that the Main Procedure returns $Q_{\text{out}} = Q^{[M_0]}$ with $1 \leq M_0 \leq M_{\max} $. Recall that the Main Procedure indicates whether $Q^{[M_0]}$ is keystone.

According to the construction in the Main Procedure, the following \underline{Main Condition} holds: there exists a chain connecting $Q=Q^{[1]}$ to $Q^{[M_0]}$ with constants $c(A)$ and $C(A)$. \\
(If $M_0=1$ then a trivial chain connects $Q$ to $Q^{[1]}=Q$.)

The construction proceeds in three cases below.

\underline{Case 1}: Suppose that $Q^{[M_0]}$ is keystone. According to the Main Condition, $Q$ is connected to the keystone cube $Q^{[M_0]}$ by a chain with  constants $c_\#$ and $C_\#$. Here, we assume that $c_\# \leq c(A)$ and $C_\# \geq C(A)$.
Therefore, we can define $\mathcal{K}(Q) := Q^{[M_0]}$ and the requisite properties listed in the bullet points  at the beginning of the explanation will be satisfied. This concludes the analysis in Case 1.

Note that, if $Q$ is keystone, then $M_0=1$ and $Q^{[1]} = Q$. To see this, just examine the Main Procedure. Hence, $\mathcal{K}(Q) = Q$ when $Q$ is keystone. This proves one of the conditions in the algorithm.

In the remaining cases, $Q^{[M_0]}$ is not a keystone cube. We then have $M_0=M_{\max} $ because the loop on $k$ cannot terminate early. Hence, by construction, $Q^{[M_{\max} ]}$ is \underline{not} a keystone cube, and 
\begin{equation}
\label{good1}
\delta_{Q^{[M_{\max} ]}} \leq 2^{-1} \cdot \delta_{Q^{[M_{\max} -1]}} \leq \cdots \leq 2^{-M_{\max} +1} \cdot \delta_{Q^{[1]}} = 2^{- M_{\max} + 1} \cdot \delta_Q.
\end{equation}

We can determine whether $Q^{[M_{\max} ]}$ appears in the list USUAL-SUSPECTS using a binary search. This takes work at most $C(A) \log N$.

\underline{Case 2}: Suppose that $Q^{[M_{\max} ]}$ is in the list USUAL-SUSPECTS. From \eqref{good1} and since $M_{\max} \geq 2$, we have $\delta_{Q^{[M_{\max}]}} \leq \frac{1}{2} \delta_Q$, hence $Q^{[M_{\max} ]}$ precedes $Q$ in the list USUAL-SUSPECTS. By induction hypothesis, we have computed a keystone cube $\mathcal{K}(Q^{[M_{\max} ]})$ to which $Q^{[M_{\max}]}$ is connected by a chain with constants $c_\#$ and $C_\#$. Moreover, another chain connects $Q$ to $Q^{[M_{\max} ]}$ (by the Main Condition).  By concatenating these chains, we see that $Q$ is connected to $\mathcal{K}(Q^{[M_{\max} ]})$ by a chain with constants $c_\#$ and $C_\#$. Here, we require that $c_\# \leq c(A)$ and $C_\# \geq C(A)$. We may thus define $\mathcal{K}(Q) :=\mathcal{K}( Q^{[M_{\max} ]})$ and the requisite properties are satisfied. This concludes the analysis in Case 2.

\underline{Case 3}: \label{pp2} Suppose that $Q^{[M_{\max} ]}$ is not in the list USUAL-SUSPECTS. Then $Q^{[M_{\max} ]}$ is interstellar, since all non-interstellar $\CZ$ cubes appear in the list USUAL-SUSPECTS. Using the algorithm \textsc{Test an
Interstellar Cube}, we determine a value of $l$ such that
\[
cA^{10}\cdot \diam( S_{l}) <\left\vert x-x(
S_{l}) \right\vert <CA^{-10}\cdot \dist( S_{l},E\setminus
S_{l}) \text{ for all }x\in ( 1+c_{G}) Q^{[M_{\max} ]} \text{.}
\]%
By definition, the cube $Q^{\text{fin}} :=Q_{l}^{\extra}$ appears in the
list USUAL-SUSPECTS.

By Lemma \ref{lem_cluster2} there exists a sequence of $\CZ$ cubes $Q_{1}\leftrightarrow\cdots \leftrightarrow Q_{L}$ such that $Q_1 = Q^{[M_{\max} ]}$, $Q_{L}=Q^{\text{fin}} $, and
\begin{equation*}
\delta _{Q_{\ell}}\leq C( A) \cdot ( 1-c( A) )
^{\ell-\ell^{\prime }}\delta _{Q_{\ell^{\prime }}}\text{ for } 1 \leq \ell^{\prime }\leq
\ell \leq L,
\end{equation*}
for controlled constants $c(A)$ and $C(A)$. 

Hence, $ \delta_{Q^{\text{fin}}}  \leq C(A) \delta_{Q^{[M_{\max} ]}} \leq C(A) 2^{-M_{\max} } \delta_Q$. We pick
\[
M_{\max} \geq \log_2(C(A)) + 1,
\]
and thus we obtain the estimate $ \delta_{Q^{\text{fin}}} \leq \frac{1}{2} \delta_Q$. 

Now, there exists a sequence of $\CZ$ cubes $\widetilde{Q}_1 \leftrightarrow \cdots \leftrightarrow \widetilde{Q}_{\widetilde{L}}$ such that $\widetilde{Q}_1 = Q$, $\widetilde{Q}_{\widetilde{L}} = Q^{[M_{\max}]}$, and $\widetilde{L} \leq C(A)$ for a controlled constant $C(A)$. This is a consequence of the construction in the Main Procedure. We concatenate the sequences $(\widetilde{Q}_\ell)_{1 \leq \ell \leq \widetilde{L}}$ and $(Q_\ell)_{1 \leq \ell \leq L}$. The resulting sequence $(\widehat{Q}_1, \cdots, \widehat{Q}_{\widehat{L}})$ satisfies
\begin{itemize}
\item $\widehat{Q}_1 = Q$, and $\widehat{Q}_{\widehat{L}} = Q^{\fin}$.
\item $\widehat{Q}_\ell \leftrightarrow  \widehat{Q}_{\ell+1}$ for $1 \leq \ell \leq \widehat{L}$.
\item $\delta_{\widehat{Q}_\ell} \leq C'(A) \cdot (1-c'(A))^{\ell - \ell'}  \cdot \delta_{Q_{\ell'}}$ for $1 \leq \ell' \leq \ell \leq \widehat{L}$.
\item $\delta_{\widehat{Q}_{\widehat{L}}} \leq \frac{1}{2} \delta_{\widehat{Q}_1}$.
\end{itemize}
Here, $c'(A)$ and $C'(A)$ are controlled constants. The last two bullet points imply that
\[\delta_{\widehat{Q}_{\widehat{L}}} \leq (1-c''(A))^{\widehat{L}-1} \delta_{Q_1}\]
for a controlled constant $c''(A) \leq c'(A)$. Hence, $Q$ is connected to $Q^\fin$ by a chain (in fact, the chain consists of a single sequence) with constants $c''(A)$ and $C'(A)$.

Moreover, since $ \delta_{Q^{\text{fin}}} \leq \frac{1}{2} \delta_Q$, we know that $Q^{\text{fin}}$ precedes $Q$ in the list USUAL-SUSPECTS.

By induction hypothesis, we have computed a keystone cube $\mathcal{K}(Q^{\fin})$ to which $Q^{\fin}$ is connected by a chain with constants $c_\#$ and $C_\#$. Moreover, as shown above,  $Q$ connects to $Q^{\fin}$ by a chain with  constants $c''(A)$ and $C'(A)$. Hence, $Q$ connects to $\mathcal{K}(Q^{\fin})$ by a chain with  constants $c_\#$ and $C_\#$. Here, we assume that $c_\# \leq c''(A)$ and $C_\# \geq C'(A)$.  We may thus define $\mathcal{K}( Q) :=\mathcal{K}(Q^{\fin})$ and the requisite properties are satisfied. This concludes the analysis in Case 3.

We review what we have achieved. By looping over all the cubes $Q$ in the list USUAL-SUSPECTS (sorted by size), we have computed for each $Q$ a keystone cube $\mathcal{K}(Q)$, and we have verified that $Q$ is connected to $\mathcal{K}(Q)$ by a chain with constants $c_\#$ and $C_\#$. We may choose $c_\#$ and $C_\#$ to be controlled constants. As mentioned before, by Lemma \ref{USUAL-SUSPECTS}, there thus exists an exponentially decreasing path connecting $Q$ to $\mathcal{K}(Q)$.

The reader may easily check that our algorithm performs as promised in terms
of the work, storage, and number of calls to the $\CZ$-\textsc{Oracle}

This concludes the explanation of the algorithm \textsc{Mark Usual Suspects}.

\end{proof}

\environmentA{Main Keystone Cube Algorithm.}

We perform one-time work, after which we can answer queries.

A query consists of a cube $\underline{Q}\in
\CZ$. The response to a query is a keystone cube $\mathcal{K}( 
\underline{Q}) $.

We guarantee the following:

\begin{itemize}
\item For each $\underline{Q}\in\CZ$ there is a finite sequence of $\CZ$
cubes 
\[
\underline{Q}=Q_{1}\leftrightarrow Q_{2}\leftrightarrow \cdots
\leftrightarrow Q_{\underline{L}}=\mathcal{K}( \underline{Q}) 
\text{ }
\]%
such that
\[
\delta _{Q_{\ell}}\leq C( A) \cdot ( 1-c( A) )
^{\ell-\ell^{\prime }}\delta _{Q_{\ell^{\prime }}}\text{ for } 1 \leq \ell^{\prime }\leq
\ell \leq \underline{L}.
\]

\item If $\underline{Q}$ is keystone, then $\mathcal{K}( \underline{Q}%
) =\underline{Q}$.

\item As part of the one-time work we compute a list called BORDER-DISPUTES, consisting of
pairs $( Q,Q^{\prime })$ with $Q,Q' \in \CZ$. A pair of $\CZ$ cubes $( Q,Q^{\prime }) $ belongs to BORDER-DISPUTES if and only if $\mathcal{K}( Q) \not=\mathcal{K}( Q^{\prime })$ and $Q\leftrightarrow Q'$. We guarantee that the list BORDER-DISPUTES consists of at most $C( A) \cdot N$ pairs of $\CZ$ cubes. 

\item The query work is at most $C( A) \log N$. The query work makes at most $C(A)$ additional calls to the $\CZ$-\textsc{Oracle}.
\item The one-time work is at most $C( A) N\log N$ in space $C(A) N$. The one-time work makes at most $C( A) N$ additional calls to the $\CZ$-\textsc{Oracle}.

\end{itemize}

\begin{proof}[\underline{Explanation}] 

As part of the one-time work, we execute the algorithm \textsc{Mark
Usual Suspects}. Hence, each cube $Q$ from the list USUAL-SUSPECTS is marked with a keystone cube $\mathcal{K}( Q) $, and we guarantee that there exists an exponentially decreasing path connecting $Q$ and $\mathcal{K}(Q)$. Furthermore, if $Q$ is keystone, then we guarantee that $\mathcal{K}( Q) =Q$.

We now explain the query algorithm.

Let $\underline{Q}$ be a $\CZ$ cube. By a binary search, we can check whether $
\underline{Q}$ belongs to the list USUAL-SUSPECTS. This requires work at most $C(A) \log N$.

If $\underline{Q} \in $ USUAL-SUSPECTS, then we have precomputed $\mathcal{K}( \underline{Q}) $ satisfying the first bullet point. 

Note that all the keystone cubes are among this list of USUAL-SUSPECTS. Hence, the second bullet point will always hold.

If $\underline{Q} \notin $ USUAL-SUSPECTS, then $\underline{Q}$ is interstellar, since all non-interstellar $\CZ$
cubes are among the USUAL-SUSPECTS. 

Applying the algorithm \textsc{Test an Interstellar Cube}, we compute an index $l$ for which
\[ c A^{10}\cdot \diam( S_l ) < \left\vert x-x( S_l)
\right\vert  < C A^{-10}\cdot \dist( S_l ,E\setminus S_l) \;\; \mbox{for all} \; x \in (1+c_G)\underline{Q},\]
where $S_{l}=3Q_{l}^{\text{CD}}\cap E$. \label{pp3} Hence, for this index $l$ we have
\begin{equation}
\label{cf1}
(1+c_G)\underline{Q} \subset H(S_l).
\end{equation}
(See the definition of the halo $H(S_l)$ in \eqref{halodefn}.) This computation requires work at most $C(A) \log N$ and uses at most $C(A)$ calls to the $\CZ$-\textsc{Oracle}.

By Lemma \ref{lem_cluster2}, there exists an exponentially decreasing
path of $\CZ$ cubes joining $\underline{Q}$ to $Q_{l}^{\extra}$; moreover, $Q_{l}^{\extra}$ is among the
USUAL-SUSPECTS. Therefore, we have precomputed a keystone cube $\mathcal{K}%
( Q_{l}^{\extra}) $, to which $Q_{l}^{\extra}$ may be joined by an
exponentially decreasing path.

We set $\mathcal{K}( \underline{Q}) :=\mathcal{K}(
Q_{l}^{\extra}) $. Note that $\underline{Q}$ is joined by an exponentially decreasing path to $Q_l^{\extra}$ and that $Q^{\extra}_l$ is joined by an exponentially decreasing path to $\mathcal{K}(Q_l^{\extra})$. Hence, there exists an exponentially decreasing path as in the first bullet point. The second bullet point holds vacuously. Indeed, all the
keystone cubes are among the USUAL-SUSPECTS, and $\underline{Q}$ is not
among the USUAL-SUSPECTS, hence $\underline{Q}$ is not keystone.

Thus, we have succeeded in responding to the query $\underline{Q}$.

We see that the work and the number of calls to the $\CZ$-\textsc{Oracle} in the query work are controlled as required.

This concludes our explanation of the query algorithm.

Next, we explain how to generate the list BORDER-DISPUTES. \label{pp4}

Suppose that $\underline{Q},\underline{\tilde{Q}}\in \CZ$, with $\underline{Q}%
\leftrightarrow \underline{\tilde{Q}}$. Assume that neither $\underline{Q}$
nor $\underline{\tilde{Q}}$ appears on the list of USUAL-SUSPECTS. Then our
query algorithm sets $\mathcal{K}( \underline{Q}) :=\mathcal{K}%
( Q_{l}^{\extra}) $ and $\mathcal{K}( \underline{\tilde{Q}}%
) :=\mathcal{K}( Q_{\tilde{l}}^{\extra}) $, where $(
1+c_{G}) \underline{Q}\subset H( S_{l}) $ and $(
1+c_{G}) \underline{\tilde{Q}}\subset H( S_{\tilde{l}}) $,
with the usual definitions $S_{l} = 3 Q_{l}^{\text{CD}} \cap E$ and $S_{\tilde{l}} = 3 Q_{\tilde{l}}^{\text{CD}} \cap E$. See \eqref{cf1}.

We recall Lemma 6.5 in Chapter 6 of \cite{FIL1}, which states that the halos $H(S)$ are pairwise disjoint as $S$ varies over all the clusters.

Since $( 1+c_{G}) \underline{Q}\cap ( 1+c_{G}) 
\underline{\tilde{Q}}\not=\emptyset $, it follows that $S_{l}=S_{\tilde{l}}$, hence $Q_{l}=\DC(
S_{l}) =\DC( S_{\tilde{l}}) =Q_{\tilde{l}}$, hence $l=\tilde{%
l}$.

Therefore, $\mathcal{K}( \underline{Q}) =\mathcal{K}(
Q_{l}^{\extra}) =\mathcal{K}( Q_{\tilde{l}}^{\extra}) =%
\mathcal{K}( \underline{\tilde{Q}}) $.

Consequently, whenever $Q,Q^{\prime }\in \CZ$ with $Q\leftrightarrow
Q^{\prime }$ and $\mathcal{K}( Q) \not=\mathcal{K}(
Q^{\prime }) $, either $Q$ or $Q^{\prime }$ is among the
USUAL-SUSPECTS.

Using our list USUAL-SUSPECTS and the $\CZ$-\textsc{Oracle}, we can easily generate a
list of all pairs of $\CZ$ cubes%
\begin{equation}
\left[ 
\begin{array}{l}
( Q,Q^{\prime }) \text{ such that }Q\leftrightarrow Q^{\prime }%
\text{ and} \\ 
Q\text{ or }Q^{\prime }\in \text{USUAL-SUSPECTS.}%
\end{array}%
\right]   \label{ptk16}
\end{equation}

There are at most $C( A) N$ such pairs, and we can generate them,
sort them and remove duplicates with work $\leq C( A) N\log N$
in space $C( A) N$, making at most $C( A) N$ calls to
the $\CZ$-\textsc{Oracle}.

Using our query algorithm, we can simply test each pair $( Q,Q^{\prime
}) $ satisfying (\ref{ptk16}) to determine whether $\mathcal{K}(
Q) =\mathcal{K}( Q^{\prime }) $.

This produces the list BORDER-DISPUTES, satisfying the third bullet point of
our algorithm. Since the are at most $C( A) N$ pairs satisfying (%
\ref{ptk16}), the fourth bullet point holds as well.

Note that we perform work $\leq C\left( A \right) N\log N$ in space $%
C( A) N$, and make at most $C( A) N$ calls to the $\CZ$-\textsc{Oracle}, in pruning the list (\ref{ptk16}) to make the list BORDER-DISPUTES.

This concludes our explanation of the \textsc{Main Keystone Cube Algorithm}.

\end{proof}

\subsubsection{The one-dimensional case.}
We now assume that $n=1$. As noted before, Lemma \ref{lem_cluster2} does not apply in this case. The cause of this failure is the fact that in one dimension the halo $H(S)$ has multiple connected components. Indeed,
\[H(S) = \{ y\in \R^{n}:  A\cdot \diam(S) <\lvert y-x( S) \rvert <A^{-1}\cdot \dist( E\setminus S,S) \}
\]
is the union of two disjoint intervals. 

We start by modifying the construction of $x_l^{\extra}$ and $Q_l^{\extra}$ from the algorithm \textsc{Make Auxiliary Cubes}. Instead of the points $x_l^{\extra}$, we define
\[x_l^{\extra, \pm} = x(S_l) \pm 4A \cdot \diam(S_l).\]
This definition yields the following result, just as before.
\begin{itemize}
\item \environmentA{Algorithm: Make Auxiliary Cubes (Version II).} For each $Q_{l}^{\text{CD}}$, $S_{l}=3Q_{l}^{\text{CD}}\cap E$, produced by the algorithm \textsc{Make Cluster Descriptors}, we compute two points $x_{l}^{\extra,-}, x_l^{\extra,+} \in
H( S_{l}) $ such that 
\[
2A\cdot \diam( S_{l}) \leq \left\vert x_{l}^{\extra,j}-x(
S_{l}) \right\vert \leq 8A\cdot \diam( S_{l})  \qquad \mbox{for} \;\; j \in \{ +,-\}.
\]
We guarantee that each of the connected components of $H(S_l)$ contains one of the points $x^{\extra,-}_l$, $x^{\extra,+}_l$. We also compute $Q_{l}^{\extra,-}$ and $Q_l^{\extra,+}$, the $\CZ$ cubes containing $x_{l}^{\extra,-}$ and $x_l^{\extra,+}$, respectively. The algorithm uses work at most $C( A) N \log N$ in space $C(A) N$, and makes at most $C( A) N$ calls to the $\CZ$-\textsc{Oracle}.
\end{itemize}

We require the following result, which is a modified version of Lemma \ref{lem_cluster2}.

\begin{lem} \label{lem_cluster2a}
Assume that $n = 1$. Let $Q\in \CZ$, and suppose that 
\begin{equation} \label{relevantcube_new} cA^{10}\cdot \diam( S_l ) <\left\vert x-x(S_l) \right\vert < CA^{-10}\cdot \dist( S_l ,E\setminus
S_l) \text{ for some} \; x \in ( 1+c_{G}) Q \text{.}
\end{equation}
Assume that $j \in \{+,-\}$ is chosen so that $x_l^{\extra,j}$ and $x$ belong to the same connected component of $H(S_l)$.

Then there exists an exponentially decreasing path $\mathcal{S} = (Q_1,\cdots,Q_{\overline{\overline{J}}})$ joining $Q$ to $Q^{\extra,j}_l$.
\end{lem}

To prove this lemma we mimic the proof of Lemma \ref{lem_cluster2}. Let $x$ and $x^{\extra,j}_l$ be as above. We apply Lemma \ref{path_lem} to the points $x$ and $x^{\extra,j}_l$, which belong to the same connected component of $H(S_l)$ according to hypothesis. Thus there exists a sequence $x_1,\cdots,x_J \in H(S_l)$ such that $x_1 = x$ and $x_1 = x^{\extra,j}_l$, which satisfies the remaining conditions described in the proof of Lemma \ref{lem_cluster2}. The remainder of the argument follows the proof of Lemma \ref{lem_cluster2} in an obvious way.

The remaining modifications necessary for the case $n=1$ are described below.

\begin{itemize}
\item (Following \textsc{Algorithm List All Non-Interstellar Cubes}, in the definition of USUAL-SUSPECTS.) \\
The list \underline{USUAL-SUSPECTS} consists of all keystone cubes and all non-interstellar cubes, and all the cubes $Q^{\extra,+}_l$, $Q_l^{\extra,-}$ produced by the algorithm \textsc{Make Auxiliary Cubes (Version II)}.
\item (The explanation of the algorithm \textsc{Mark Usual Suspects}, in the analysis of \underline{Case 2}.) \\
We know that $(1+c_G)Q' \subset H(S_l)$. We determine $j \in \{ +, -\}$ such that $x^{\extra,j}_l$ belongs to the connected component of $H(S_l)$ that contains $(1+c_G)Q'$. For that $l$ and that $j$, the cube $Q'' = Q_l^{\extra,j}$ appears in the list USUAL-SUSPECTS, and by Lemma \ref{lem_cluster2a} there exists a sequence of $\CZ$ cubes $Q' = Q_1 \leftrightarrow Q_2 \leftrightarrow \cdots \leftrightarrow Q_L = Q''$ such that ... 

We set $\mathcal{K}(Q) := \mathcal{K}(Q'')$.
\item (The explanation of the \textsc{Main Keystone Cube Algorithm}, in the analysis of the case in which $\underline{Q}$ is \underline{not} interstellar.) \\ We know that $(1+c_G) \underline{Q} \subset H(S_l)$. We determine $j \in \{ +, -\}$ such that $x^{\extra,j}_l$ belongs to the same connected component of $H(S_l)$ which contains $(1+c_G)\underline{Q}$. We know that $\underline{Q}$ can be joined by an exponentially decreasing path of $\CZ$ cubes to $Q_l^{\extra,j}$, and that $Q_l^{\extra,j}$ is among the USUAL-SUSPECTS. Therefore, we have precomputed a keystone cube $\mathcal{K}(Q_l^{\extra,j})$, to which $Q_l^{\extra,j}$ may be joined by an exponentially decreasing path.

We set $\mathcal{K}(\underline{Q}) : = \mathcal{K}(Q_l^{\extra,j})$.

\item (The explanation of the \textsc{Main Keystone Cube Algorithm}, in the definition of BORDER-DISPUTES.) \\
Then our query algorithm sets $\mathcal{K}(\underline{Q}) = \mathcal{K}(Q_l^{\extra,j})$ and $\mathcal{K}(\underline{\tilde{Q}}) = \mathcal{K}(Q_{\tilde{l}}^{\extra,\tilde{j}})$, where \\
$(1+c_G)\underline{Q}$ and $x^{\extra,j}_l$ are contained in the same connected component of $H(S_l)$, and where \\
$(1+c_G)\underline{\tilde{Q}}$ and $x^{\extra,\tilde{j}}_{\tilde{l}}$ are contained in the same connected component of $H(S_{\tilde{l}})$.

Since $( 1+c_{G}) \underline{Q}\cap ( 1+c_{G}) 
\underline{\tilde{Q}}\not=\emptyset $, while the halos $H( S) $
are pairwise disjoint as $S$ varies over all
clusters, it follows that $S_{l}=S_{\tilde{l}}$. Moreover, $(1+c_G)\underline{Q}$ and $( 1+c_{G}) 
\underline{\tilde{Q}}$ are contained in the same connected component of $H(S_l) = H(S_{\tilde{l}})$, hence $x^{\extra,j}_l = x^{\extra,\tilde{j}}_{\tilde{l}}$. Thus we have $l=\tilde{l}$ and $j = \tilde{j}$.

Therefore, $\mathcal{K}( \underline{Q}) =\mathcal{K}(
Q_{l}^{\extra,j}) =\mathcal{K}( Q_{\tilde{l}}^{\extra, \tilde{j}}) =%
\mathcal{K}( \underline{\tilde{Q}}) $.
\end{itemize}

This concludes the list of modifications required to treat the case $n=1$.

\section{CZ Decompositions}\label{sec_czdecomp}

\subsection{Preliminaries}

\begin{lem}\label{gg_lem_1}
Let $0 < \gamma < 1$ with $\gamma$ an integer power of two.

Let $\CZ$ be a collection of pairwise disjoint dyadic cubes. We assume either that $\CZ$ is a dyadic decomposition of a unit cube $Q^\circ$ or that $\CZ$ is a dyadic decomposition of $\R^n$.

Assume that for all $Q,Q' \in \CZ$ with $Q \leftrightarrow Q'$, we have $\gamma \delta_{Q'} \leq \delta_Q \leq \gamma^{-1} \delta_{Q'}$.

Then, for any $Q,Q' \in \CZ$ with $(1+\gamma/2)Q \cap (1+\gamma/2)Q' \neq \emptyset$, we have $Q \leftrightarrow Q'$.

\end{lem}

\begin{proof}
We assume that $\CZ$ is a dyadic decomposition of $Q^\circ$, where $Q^\circ$ is a unit cube. The case in which $\CZ$ is a dyadic decomposition of $\R^n$ is treated similarly.

Let $Q,Q' \in \CZ$ satisfy $(1+\gamma/2)Q \cap (1+\gamma/2)Q' \neq \emptyset$ and $\delta_Q \geq \delta_{Q'}$. For the sake of contradiction suppose that $Q$ and $Q'$ do not meet. That is, we assume that the closure of $Q$ is disjoint from the closure of $Q'$.

Fix a point $z \in (1+\gamma/2)Q \cap (1+\gamma/2)Q'$.

Now,
\begin{align*}
d(Q,Q') &= \inf_{\substack{x \in Q \\ x' \in Q'}} \lvert x -x' \rvert \\
&\leq \inf_{\substack{x \in Q \\ x' \in Q'}} \lvert x -z\rvert + \lvert x' - z \rvert \\
&\leq (\gamma/4) \delta_Q + (\gamma/4) \delta_{Q'},
\end{align*}
where in the last inequality we use that $z \in (1+\gamma/2)Q$ and $z \in (1+\gamma/2)Q'$. (Recall that we use the  $\ell^\infty$ norm on $\R^n$.) Since $\delta_{Q'} \leq \delta_Q$, we conclude that $d(Q,Q') \leq (\gamma/2) \delta_Q$.

Consider the subset
\[\mathcal{D}_Q = \bigcup \left\{ \overline{Q} : \overline{Q} \in \CZ, \; \overline{Q} \leftrightarrow Q \right\} \subset Q^\circ.\]
According to good geometry, each of the above $\overline{Q}$ satisfies $ \delta_{\overline{Q}} \geq \gamma \delta_Q$. Thus, because the cubes in $\CZ$ are a partition of $Q^\circ$, we have
\[ \bigl\{ y \in Q^\circ : d(y,Q) \leq (3\gamma/4) \cdot \delta_Q \bigr\} \subset \mathcal{D}_Q.
\]
Hence, since $d(Q,Q') \leq (\gamma/2) \delta_Q$, we know that $\mathcal{D}_Q$ intersects $Q'$. Therefore, there exists $\overline{Q} \in \CZ$ with  $\overline{Q} \leftrightarrow Q$ and $\overline{Q} \cap Q' \neq \emptyset$. Hence, because the cubes in $\CZ$ are pairwise disjoint, we must have $\overline{Q} = Q'$. Thus, $Q' \leftrightarrow Q$, which contradicts our assumption that $Q$ and $Q'$ do not meet. This completes the proof of the lemma by contradiction.

\end{proof}

\subsection{Review of known results}
\label{sec_rkr}

We review several results from Sections 20-26 in \cite{FK2}.

In those sections, we are given the following data (see Section 20 in \cite{FK2}):
\begin{itemize}
\item A finite subset $E \subset \R^n$, with $\#(E) = N$, $N \geq 2$.
\item A real number $A_2 \geq 1$, assumed to be an integer power of $2$.
\item For each $x \in E$ and $\cA \subset \cM$, a positive real number $\delta(x,\cA)$. \footnote{We recall from \cite{FK2} that Lemma 5 in Section 20 there makes use of a particular choice of the $\delta(x,\cA)$, but that lemma has no effect on anything else in Sections 20-26 of \cite{FK2}. Again, see the remarks in the first few paragraphs of Section 20 in \cite{FK2}.} 
\end{itemize}

These define a family of Calder\'on-Zygmund decompositions of $\R^n$, called $\CZ(\cA)$, indexed by subsets $\cA \subset \cM$.

Here, $\CZ(\cA)$ consists of the maximal dyadic cubes $Q \subset \R^n$ of sidelength $\delta_Q \leq A_2^{-1}$ such that either
\begin{enumerate}[(a)]
\item $\#(5Q \cap E) \leq 1$

or
\item for some $\cA' \leq \cA$ we have $\delta(x,\cA') \geq A_2 \delta_Q$ for all $x \in E \cap 5Q$.
\end{enumerate}

The following algorithm is presented in Section 26 of \cite{FK2}:

Given $A_2$, $E$, $(\delta(x,\cA))_{\substack{ \cA \subset \cM \\ x \in E}}$, we perform one-time work at most $C N \log N$ in space $CN$, after which we can answer queries as follows.

A query consists of a subset $\cA \subset \cM$ and a point $\underline{x} \in \R^n$. The response to the query $(\cA,\underline{x})$ is the one and only one cube $Q \in \CZ(\cA)$ containing $\underline{x}$. The work to answer a query is at most $C \log N$. Here, $C$ depends only on $m$ and $n$.

We will make a slight change here by replacing $5Q$ by $3Q$ in the definition of $\CZ(\cA)$ (see (a) and (b) above). This change affects nothing significant in the relevant discussion in \cite{FK2}.

The only point worth mentioning is the proof of good geometry. Lemma 2 in Section 21 of \cite{FK2} asserts that if $(1+2c_G)Q \cap (1+2c_G) Q' \neq \emptyset$ with $Q,Q' \in \CZ(\cA)$, then $\frac{1}{2} \delta_Q \leq \delta_{Q'} \leq 2 \delta_Q$. Here, $c_G > 0$ is a small constant depending only on the dimension $n$.

The proof of that lemma requires slight changes; the argument given in \cite{FK2} shows that $\frac{1}{2} \delta_Q \leq \delta_{Q'} \leq 2 \delta_Q$ for any $Q,Q' \in \CZ(\cA)$ with $Q \leftrightarrow Q'$. Thus, applying Lemma \ref{gg_lem_1} (with $\gamma = 1/2$), we see that
\[\bigl[ Q,Q' \in \CZ(\cA), (1+2c_G)Q \cap (1+2c_G)Q' \neq \emptyset \bigr] \implies \frac{1}{2} \delta_Q \leq \delta_{Q'} \leq 2 \delta_Q.\]
This proves the ``good geometry'' of the cubes in $\CZ(\cA)$.

\subsection{A Calder\'on-Zygmund Oracle} \label{sec_czoracle}

We assume we are given the following data.
\begin{itemize}
\item We are given a finite set $E \subset Q^\circ$, with $Q^\circ \subset \R^n$ a dyadic cube of unit sidelength; we assume that $\#(E) = N$, $N \geq 2$.
\item We are given a number $\Delta(x) \in (0,1] $ for each $x \in E$. We denote $\vec{\Delta} = (\Delta(x))_{x \in E}$. 
\end{itemize}

Given the data above, we define a Calder\'on-Zygmund decomposition $\CZ(\vec{\Delta})$ of $Q^\circ$ as follows: $\CZ(\vec{\Delta})$ consists of the maximal dyadic cubes $Q \subset Q^\circ$ such that either $\#(E \cap 3Q) \leq 1$ or $\Delta(x) \geq \delta_Q$ for all $x \in E \cap 3Q$.

\environmentA{Algorithm: Plain Vanilla CZ-Oracle.} 

Given $E$, $\vec{\Delta}$ as above, we perform one-time work at most $ C N \log N$ in space $CN$, after which we can answer queries. 

A query consists of a point $\underline{x} \in Q^\circ$. The response to the query $\underline{x}$ is the one and only one cube $Q \in \CZ(\vec{\Delta})$ containing $\underline{x}$. 

The work to answer a query is at most $C \log N$. Here, $C$ depends only on the dimension $n$.

\begin{proof}[\underline{Explanation}]

We take $A_2 =1$ and $\delta(x,\cA) = \Delta(x)$ for each $x \in E$, $\cA \subset \cM$, and we apply the query algorithm given in the previous section.

\end{proof}

\begin{remk}\label{rem_cw}
As a special case, we can apply the \textsc{Plain Vanilla CZ-Oracle} to the ``classic Whitney decomposition'' of $Q^\circ$, which consists of the maximal dyadic subcubes $Q \subset Q^\circ$ such that $\#(E \cap 3Q) \leq 1$. In fact, we need only pick $\delta_{\text{\tiny small}}$ with $0 < \delta_{\text{\tiny small}} < \frac{1}{100} \min \{ \lvert x -y\rvert : x,y \in E, x \neq y \}$, and then take $\Delta(x) = \delta_{\text{\tiny small}}$ for all $x \in E$.

Such a number $\delta_{\text{\tiny small}}$ may be computed with one-time work $\leq C N$ once we have the Well-Separated Pairs Decomposition available. The classic Whitney decomposition coincides with $\CZ(\vec{\Delta})$. We will use the Oracle for this decomposition in a later section.
\end{remk}

We close this section with an easy generalization of the \textsc{Plain Vanilla CZ-Oracle}.

\label{page62}
We assume we have already defined a decomposition $\CZ_{\text{\tiny old}}$ of $Q^\circ$ consisting of pairwise disjoint dyadic subcubes. We make the following assumptions:
\begin{itemize}
\item If $Q \subset Q^\circ$ is a dyadic subcube and $\#(E \cap 3Q) \leq 1$, then $Q$ is contained in a cube of $\CZ_{\text{\tiny old}}$.
\item Good geometry: If $Q,Q' \in \CZ_\old$ and $Q \leftrightarrow Q'$ then $\frac{1}{2} \delta_Q \leq \delta_{Q'} \leq 2 \delta_Q$.
\item We have available a $\CZ_{\text{\tiny old}}$-\textsc{Oracle}: Given a query point $\underline{x} \in Q^\circ$, the $\CZ_{\text{\tiny old}}$-\textsc{Oracle} returns the one and only one cube $Q \in \CZ_{\text{\tiny old}}$ containing $\underline{x}$. 
\end{itemize}


\label{cz_new}
We define a decomposition $\CZ_{\new}$ of $Q^\circ$ to consist of the maximal dyadic cubes $Q \subset Q^\circ$ such that either $Q \in \CZ_{\text{\tiny old}}$ or $\Delta(x) \geq \delta_Q$ for all $x \in E \cap 3Q$.

We clearly see that the decomposition $\CZ_\new$ has good geometry, namely
\[
\mbox{If} \; Q,Q' \in \CZ_\new \; \mbox{and} \; Q \leftrightarrow Q' \;\; \mbox{then} \;\; \frac{1}{2} \delta_Q \leq \delta_{Q'} \leq 2 \delta_Q.
\]
Applying Lemma \ref{gg_lem_1} with $\gamma = 1/2$, we obtain
\begin{equation} \label{gg_1}
\mbox{If} \;  Q,Q' \in \CZ_\new \; \mbox{and} \; \frac{65}{64}Q \cap \frac{65}{64}Q' \neq \emptyset, \; \mbox{then} \; Q \leftrightarrow Q' \; \mbox{and} \; \frac{1}{2} \delta_Q \leq \delta_{Q'} \leq 2 \delta_Q.
\end{equation}
We finish this section with the following algorithm.

\environmentA{Algorithm: Glorified CZ-Oracle.}

Given $E$ and $\vec{\Delta}$ as above, we perform one-time work at most $C N \log N$ in space $C N$, after which we can answer queries. 

A query consists of a point $\underline{x} \in Q^\circ$. The response to the query $\underline{x}$ is a list containing all the cubes $Q \in \CZ_{\text{\tiny \new}}$ such that $\underline{x} \in \frac{65}{64} Q$.

We answer the query using at most $C \log N$ computer operations as well as at most $C$ calls to the $\CZ_{\mbox{\tiny old}}$-\textsc{Oracle}.

\begin{proof}[\underline{Explanation}]  

First, given $\underline{x} \in Q^\circ$, we show how to compute the unique cube $Q_{\underline{x}} \in \CZ_{\text{\tiny new}}$ containing $\underline{x}$. In fact, $Q_{\underline{x}}$ is simply the larger of the following two cubes:
\begin{itemize}
\item The cube returned by the $\CZ_{\mbox{\tiny old}}$-\textsc{Oracle} in response to the query $\underline{x}$.
\item The cube returned by the \textsc{Plain Vanilla CZ-Oracle} applied to $\vec{\Delta} = (\Delta(x))_{x \in E}$ in response to the query $\underline{x}$.
\end{itemize}
The above computation requires work at most $C \log N$ and one call to the $\CZ_{\old}$-\textsc{Oracle}.

Next, given $\underline{x} \in Q^\circ$, we show how to compute a list of all $Q \in \CZ_{\text{\tiny new}}$ such that $\underline{x} \in \frac{65}{64} Q$. To do so, we first compute the cube $Q_{\underline{x}} \in \CZ_{\text{\tiny new}}$ containing $\underline{x}$. By \eqref{gg_1}, our desired list of cubes consists only of dyadic cubes $Q \subset Q^\circ$ such that 
\[ \left[ \underline{x} \in \frac{65}{64} Q \;\; \mbox{and} \;\; \frac{1}{2} \delta_{Q_{\underline{x}}} \leq \delta_Q \leq 2 \delta_{Q_{\underline{x}}} \right].\]

There are at most $C$ such cubes, and we can easily list them all.

Now, we test each such $Q$ to see whether $Q \in \CZ_{\text{\tiny new}}$. To do that, we just compute the one and only one cube $\hQ \in \CZ_{\text{\tiny new}}$ containing the center of $Q$, and we check whether $\hQ = Q$.

This completes our description of the \textsc{Glorified CZ-Oracle}. It's trivial to check that the algorithm works, and that the one-time work, query work, the storage, and the number of calls to the $\CZ_{\text{\tiny old}}$-\textsc{Oracle} are as promised.

\end{proof}

\subsection{Basic algorithms} \label{czalg_sec}

In the present section and in the next section (Section \ref{sec_pou}), we assume that we are given the following.
\begin{itemize}
\item A finite set $E \subset \frac{1}{32} Q^\circ$, with $Q^\circ$ a dyadic cube of unit sidelength in $\R^n$, such that $N := \#(E) \geq 2$.
\item A collection $\CZ$ consisting of dyadic cubes $Q \subset \R^n$. We assume that $\CZ$ is locally finite, i.e., any given compact set $S \subset \R^n$ intersects a finite number of cubes $Q \in \CZ$. Furthermore, we assume 
\[\text{Good geometry: \; If \;} Q \leftrightarrow Q' \; \mbox{and} \; Q,Q' \in \CZ, \; \mbox{then} \; \frac{1}{8} \delta_Q \leq \delta_{Q'} \leq 8 \delta_Q.\]
\item We assume we are either in\\
\textbf{Setting 1:} The cubes in $\CZ$ partition $Q^\circ$, or \\
\textbf{Setting 2:} The cubes in $\CZ$ partition $\R^n$.

\item We assume that a $\CZ$-\textsc{Oracle} is available. The $\CZ$-\textsc{Oracle} accepts queries.  In \textbf{Setting 1}, a query consists of a point $\underline{x} \in Q^\circ$; in \textbf{Setting 2}, a query consists of a point $\underline{x} \in \R^n$. Given a query $\underline{x}$, the $\CZ$-\textsc{Oracle} produces a list of the cubes $Q \in \CZ$ such that $\underline{x} \in \frac{65}{64} Q$. This requires work at most $C \cdot \log N$.
\end{itemize}

We see that $\CZ$ satisfies the hypotheses of Lemma \ref{gg_lem_1} with $\gamma = 1/8$. Thus, 
\begin{equation} \label{gg_2}
\mbox{If} \; \frac{65}{64}Q \cap \frac{65}{64}Q' \neq \emptyset \;\mbox{and} \; Q,Q' \in \CZ, \; \mbox{then} \; Q \leftrightarrow Q' \; \mbox{and} \; \frac{1}{8} \delta_Q \leq \delta_{Q'} \leq 8 \delta_Q.
\end{equation}

Let $Q \in \CZ$ and $\overline{Q} \in \CZ \setminus \{ Q \}$ be given. Let $x_Q$ be the center of $Q$. Suppose that $ \frac{65}{64} \overline{Q} \cap B(x_Q,\delta) \neq \emptyset$ for some $0 < \delta < \frac{1}{64} \min \{ \delta_{\oQ}, \delta_Q\}$. Then $ \frac{65}{64} \overline{Q} \cap \frac{65}{64}Q  \neq \emptyset$. From \eqref{gg_2}, we see that $\delta_Q$ and $\delta_{\overline{Q}}$ differ by at most a factor of $16$. Thus, because $Q$ and $\overline{Q}$ are disjoint dyadic cubes, we have
\[d(x_Q,\overline{Q}) \geq \frac{1}{2} \delta_Q \geq \frac{1}{32} \delta_{\overline{Q}}.\]
(Recall, distances are measured using the $\ell^\infty$ metric.) However, our assumption that $ \frac{65}{64} \overline{Q} \cap B(x_Q,\delta) \neq \emptyset$ implies that $d(x_Q, \overline{Q}) \leq \delta + \frac{1}{64}\delta_{\overline{Q}} < \frac{1}{32}\delta_{\oQ}$. This contradiction establishes
\begin{equation}
\label{gg_3}
\mbox{If} \; Q, \overline{Q} \in \CZ  \; \mbox{and} \; Q \neq \overline{Q}, \; \mbox{then} \; B(x_Q,\delta) \cap \frac{65}{64} \overline{Q} = \emptyset \quad \mbox{for} \; \delta < \frac{1}{64} \min \{ \delta_{\oQ}, \delta_Q\}.
\end{equation}

Under the above assumptions, we give the following algorithms.

\environmentA{Algorithm: Find Neighbors.} We can answer queries as follows. A query consists of a cube $Q \in \CZ$. The response to the query $Q$ is the list of all cubes $Q' \in \CZ$ such that $Q' \leftrightarrow Q$. To answer the query requires work at most $C \log N$.

\begin{proof}[\underline{Explanation}]

We first explain how to test whether a given dyadic cube $Q' \subset \R^n$ belongs to the collection $\CZ$. In \textbf{Setting 1}, it is necessary that $Q' \subset Q^\circ$. Assuming that this is the case, we examine the center $x_{Q'}$ of $Q'$. We query the $\CZ$-\textsc{Oracle} on $x_{Q'}$ to produce the list of all cubes $Q \in \CZ$ with $x_{Q'} \in \frac{65}{64}Q$. This list contains at most $C$ cubes, thanks to Good Geometry. Note that $Q'$ belongs to this list if and only if $Q'$ belongs to $\CZ$. We can check the former condition using work at most $C$. 

Let $Q \in \CZ$ be given. We wish to list all the $Q' \in \CZ$ such that $Q' \leftrightarrow Q$. According to Good Geometry, each such $Q'$ also satisfies $\frac{1}{8} \delta_Q \leq \delta_{Q'} \leq 8 \delta_Q$. We can list all the dyadic cubes $Q'$ with $Q' \leftrightarrow Q$ and $\frac{1}{8}\delta_Q \leq \delta_{Q'} \leq 8 \delta_Q$. We remove from this list those cubes that do not belong to $\CZ$. We return a list of the remaining cubes.

This completes our description of the algorithm \textsc{Find Neighbors}. It's easy to check that the algorithm operates as promised, and that the amount of work is as promised.

\end{proof}

\environmentA{Algorithm: Find Main-Cubes.} After one-time work at most $C N \log N$ in space $CN$, we produce the collection of cubes $\CZ_{\main} := \{ Q \in \CZ : \frac{65}{64}Q \cap E \neq \emptyset\}$. We mark each cube $Q \in \CZ_{\main}$ with a point $x(Q) \in \frac{65}{64}Q \cap E$.

\begin{proof}[\underline{Explanation}]

We loop over $x \in E$. For each point $x \in E$, we list all the cubes $Q \in \CZ$ such that $x \in \frac{65}{64}Q$. This requires $N$ calls to the $\CZ$-\textsc{Oracle}. For each $Q$ obtained above, we set $x(Q) := x$ for the relevant $x$. Thus, we produce a list of all the cubes in $CZ_{\main}$, possibly containing duplicates. After sorting this list, we can find and remove duplicates, and obtain our desired list of the cubes $Q \in CZ_{\main}$ marked by points $x(Q)$.

\end{proof}

\subsection{Partitions of unity}\label{sec_pou}

Aside from a decomposition $\CZ$ satisfying the conditions laid out in Section \ref{czalg_sec}, we assume that we are given a cube $\hQ \subset \R^n$, and real numbers $0 < \overline{r} \leq 1/64$ and $A \geq 1$. We are also given a finite subcollection $\cQ \subset \CZ$ with the following properties:
\begin{align}
\label{covers}
& \mbox{For each} \;x \in \hQ \; \mbox{we have} \; x \in (1+\overline{r}/2)Q \; \mbox{for some} \; Q \in \cQ.\\
\label{sizebd}
&\delta_Q \leq A \delta_{\hQ} \; \mbox{for each} \; Q \in \cQ.
\end{align}
(We do not assume here that $\hQ$ is dyadic.)

By \eqref{gg_2}, we see that the collection $\{\frac{65}{64}Q : Q \in \cQ\}$ has \emph{bounded overlap}, meaning that for each $Q \in \cQ$ there are at most $C$ cubes $Q' \in \cQ$ such that $\frac{65}{64}Q \cap \frac{65}{64}Q' \neq \emptyset$. Here, $C$ depends only on the dimension $n$.

For each $Q \in \CZ$ we choose a cutoff function $\widetilde{\theta}_Q \in C^m(\R^n)$ such that
\begin{itemize}
\item $\supp ( \widetilde{\theta}_Q) \subset (1+\frac{3\overline{r}}{4})Q$ .
\item $\widetilde{\theta}_Q \geq 0$ on $\R^n$.
\item $\widetilde{\theta}_Q \geq 1/2$ on $\left(1+\frac{\overline{r}}{2}\right)Q$.
\item $\lvert \partial^\alpha \widetilde{\theta}_Q(x) \rvert \leq C(\overline{r}) \cdot \delta_Q^{-|\alpha|}$ for $x \in \R^n$, $|\alpha| \leq m.$
\end{itemize}
We choose $\widetilde{\theta}_Q$ to depend only on $Q$ and $\overline{r}$.

We assume the existence of a query algorithm for $\widetilde{\theta}_Q$. For instance, we can take $\widetilde{\theta}_Q$ to be a tensor product of univariate splines, in which case the next algorithm is trivial.

\label{pp9}

\environmentA{Algorithm: Compute Cutoff Function.} Given a cube $Q \in \CZ$, a point $\underline{x} \in Q^\circ$, and $0 < \overline{r} \leq 1/64$, we compute the jet $J_{\underline{x}} (\widetilde{\theta}_Q)$ using work and storage at most $C$.

In the next lemma we use the cutoff functions $\widetilde{\theta}_Q$ to construct a \emph{partition of unity}.

Recall that $x_Q$ denotes the center of a cube $Q$.

\begin{lem} \label{pou_lem}
There exists $\theta^\hQ_Q \in C^m(\R^n)$ for each $Q \in \cQ$, such that \begin{align} 
\label{suppprop}
& \supp \theta^\hQ_Q \subset \left(1+ \frac{3\overline{r}}{4}\right)Q,\\
\label{derbds}
&\lvert \partial^\alpha \theta^\hQ_Q (x) \rvert \leq C(\overline{r}) \cdot \delta_Q^{-|\alpha|} \; \mbox{for} \;x \in \R^n, \; |\alpha| \leq m,\\
\label{sumtoone} 
&1 = \sum_{Q \in \cQ} \theta^\hQ_Q \;\;\; \mbox{on} \; \hQ.
\end{align}
Moreover, $\theta^\hQ_Q = 1$ near $x_Q$ and $\theta^\hQ_Q = 0$ near $x_{Q'}$ for all $Q' \in \cQ \setminus \{Q\}$. 

Here, the constant $C(\overline{r})$ depends only on $\overline{r}$, $m$ and $n$.

\end{lem}
\begin{proof}

We set
$$\Psi(x) = \sum_{\overline{Q} \in \cQ} \widetilde{\theta}_{\overline{Q}}(x) \;\; \mbox{for} \; x \in \R^n.$$
Because $\widetilde{\theta}_{\overline{Q}} \geq 1/2$ on $(1+\overline{r}/2){\overline{Q}}$, the condition \eqref{covers} implies that $\Psi \geq 1/2$ on $\hQ$. We can easily see that
\begin{align} \label{cut3}
\lvert \partial^\alpha \Psi(x) \rvert\leq C(\overline{r}) \delta_Q^{-|\alpha|} \qquad \mbox{for} \; x \in \frac{65}{64}Q, Q \in \cQ, \;  |\alpha| \leq m.
\end{align}
More precisely, since $\supp(\widetilde{\theta}_{\overline{Q}}) \subset \frac{65}{64}\overline{Q}$, any cube $\overline{Q} \in \cQ$ that contributes to the sum defining $\Psi(x)$ must satisfy $\frac{65}{64}Q \cap \frac{65}{64} \overline{Q}  \neq \emptyset$. (Recall that $x \in \frac{65}{64}Q$.) Moreover, $\delta_{\overline{Q}}$ and $\delta_Q$ differ by at most a factor of $16$ for any such $\overline{Q}$; see \eqref{gg_2}. Hence, \eqref{cut3} follows from \eqref{derbds} and from the fact that the sum defining $\Psi(x)$ has at most $C$ nonzero terms for each fixed $x$, a consequence of the bounded overlap of the cubes in $\cQ$.

Let $\eta \in C^m([0,\infty))$ be a function with $\eta(t) \geq 1/4$ for $t \in [0, 1/2)$, and $ \eta(t) = t$ for $t \geq 1/2$. 

Let $Q \in \cQ$. We define
$$\theta^\hQ_Q(x) :=  \frac{\widetilde{\theta}_Q(x)}{\eta \circ \Psi(x)}, \quad \mbox{a function in} \; C^m(\R^n).$$
Clearly, $ \supp \theta^\hQ_Q \subset \supp \widetilde{\theta}_Q \subset (1+3\overline{r}/4)Q$. Moreover, \eqref{cut3} implies that
\[ \lvert \partial^\alpha \left[ \eta \circ \Psi \right] (x) \rvert \leq C(\overline{r}) \delta_Q^{-|\alpha|} \;\; \mbox{for} \; x \in \frac{65}{64}Q, \; |\alpha| \leq m.
\] 
Using that $\eta \circ \Psi(x) \geq 1/4$, we obtain
\[\lvert \partial^\alpha \theta^\hQ_Q(x) \rvert = \left\lvert \partial^\alpha \left[ \frac{\widetilde{\theta}_Q}{\eta \circ \Psi} \right] (x) \right\rvert \leq C(\overline{r}) \delta_Q^{-|\alpha|}  \;\; \mbox{for} \; x \in \frac{65}{64}Q, \; |\alpha| \leq m.\]

Finally, note that
\begin{equation} \label{pou2} \sum_{Q \in \cQ}  \theta^\hQ_Q(x) = \sum_{Q \in \cQ} \frac{\widetilde{\theta}_Q(x)}{\eta \circ \Psi(x)} = \sum_{Q \in \cQ} \frac{\widetilde{\theta}_Q(x)}{ \Psi(x)} = 1 \quad \mbox{for} \; x \in \hQ.
\end{equation}
(Here, we use the fact that $\eta \circ \Psi(x) = \Psi(x)$, since $\Psi \geq 1/2$ on $\hQ$.)

Recall that $\theta_{\overline{Q}}^\hQ \geq 0$ and that $\supp(\theta^\hQ_{\overline{Q}}) \subset \frac{65}{64} \overline{Q}$ for each $\overline{Q} \in \cQ$. Thus, from \eqref{gg_3} and \eqref{pou2} we deduce that there exists $\delta > 0$ such that $\theta^\hQ_Q = 1$ on $B(x_Q,\delta)$ and $\theta^\hQ_Q=0$ on $B(x_{Q'},\delta)$ for every $Q' \in \cQ \setminus\{Q\}$.

This completes the proof of the lemma.
\end{proof}

\begin{lem}\label{patch_lem}
Given a function $F_Q \in \X((1+\overline{r})Q \cap \hQ)$ for each $Q \in \cQ$, we define
$$F := \sum_{Q \in \cQ} F_Q \theta^\hQ_Q \;\;\; \mbox{on} \; \hQ, \;\; \mbox{with} \; \theta_Q^{\hQ} \; \mbox{as in Lemma \ref{pou_lem}}.$$
Then, given a polynomial $P_Q \in \cP$ and a point $y_Q \in Q$ for each $Q \in \cQ$, we have
\begin{align}
\|F\|_{\X(\hQ)}^p \leq & \; C(A, \overline{r}) \cdot \biggl[ \;\sum_{Q \in \cQ} \bigl[ \| F_Q \|_{\X((1+\overline{r})Q \cap \hQ)}^p  + \delta_Q^{-mp} \| F_Q - P_Q \|^p_{L^p((1+\overline{r})Q \cap \hQ)} \bigr]  \label{patch_norm} \\
& \qquad\qquad + \sum_{\substack{ Q,Q' \in \cQ \\ (1+\overline{r}) Q \cap (1+\overline{r})Q' \neq \emptyset}} \sum_{|\beta| \leq m-1} \delta_Q^{(|\beta| -m)p + n} \lvert \partial^\beta (P_Q - P_{Q'})(y_Q) \rvert^p \biggr]. \notag{}
\end{align}
Here, the constant $C(A, \overline{r})$ depends only on $\overline{r}$, $A$, $m$, $n$, and $p$.
\end{lem}

\begin{proof}

Let $Q' \in \cQ$ be given, with $\hQ \cap (1 + \frac{\overline{r}}{2}) Q' \neq \emptyset$.

Let $x \in \hQ \cap (1+ \frac{\overline{r}}{2})Q'$. Recall that $\displaystyle \sum_{Q \in \cQ} \theta^\hQ_Q = 1$ on $\hQ$, hence
\[F = F_{Q'} + \sum_{Q \in \cQ} \theta^\hQ_Q \cdot (F_Q - F_{Q'}) \;\;\; \mbox{on} \; \hQ. \]
Differentiating the above equation, for $|\alpha| = m$ we have
\[\partial^\alpha F(x) = \partial^\alpha F_{Q'}(x) +  \sum_{\substack{Q \in \cQ \\ (1+\frac{3\overline{r}}{4})Q \ni x }} \sum_{\beta + \gamma = \alpha} \mbox{coeff}(\beta, \gamma) \cdot \partial^\beta ( F_Q - F_{Q'})(x) \cdot \partial^\gamma \theta^\hQ_Q(x).\]
There are at most $C$ nonzero terms in the above sum, thanks to bounded overlap of $\{ (1+\overline{r})Q : Q \in \cQ\}$. 

Let $Q \in \cQ$ be such that $(1+3\overline{r}/4)Q \ni x$. Note that $x \in \hQ \cap (1+\overline{r})Q$ and $x \in \hQ \cap (1+\overline{r})Q'$. 

In the above sum, if $|\beta| = m$ then $\beta = \alpha$ and $\gamma=0$. These terms are bounded in magnitude by $| \partial^{\alpha} F_Q(x) | + | \partial^\alpha F_{Q'}(x) |$.

In the above sum, if $| \beta | \leq m-1$ then we have
\[| \partial^\beta (F_Q - F_{Q'})(x)| \leq | \partial^\beta (F_Q - P_{Q})(x)| + |\partial^\beta(P_Q - P_{Q'})(x)| + | \partial^\beta(F_{Q'} - P_{Q'})(x)|.\]
Since $\hQ \cap (1+3\overline{r}/4)Q \neq \emptyset$ and $\delta_Q \leq A \delta_{\hQ}$ (see \eqref{sizebd}), the sidelengths of the rectangular box $\hQ \cap (1+\overline{r})Q$ are comparable to $\delta_Q$ (up to a constant factor depending on $\overline{r}$, $A$ and $n$). Similarly, the sidelengths of the rectangular box $\hQ \cap (1+\overline{r})Q'$ are comparable to $\delta_{Q'}$. Thus, by an easy rescaling argument, Lemma \ref{si2} shows that
\begin{align*}
\lvert \partial^\beta (F_Q - P_{Q})(x) \rvert &\leq C(A, \overline{r}) \cdot \left( \delta_Q^{-|\beta| - \frac{n}{p}} \| F_Q - P_Q \|_{L^p((1+\overline{r})Q \cap \hQ) } + \delta_Q^{m-|\beta| - \frac{n}{p}}\|F_Q\|_{\X((1+\overline{r})Q  \cap \hQ)} \right).\\
\lvert \partial^\beta (F_{Q'} - P_{Q'})(x) \rvert &\leq C(A, \overline{r}) \cdot \left( \delta_{Q'}^{-|\beta| - \frac{n}{p}} \| F_{Q'} - P_{Q'} \|_{L^p((1+\overline{r})Q' \cap \hQ ) } + \delta_{Q'}^{m-|\beta| - \frac{n}{p}}\|F_{Q'}\|_{\X((1+\overline{r})Q'  \cap \hQ )} \right).
\end{align*}

If $\beta + \gamma = \alpha$ then $|\gamma| = m - |\beta|$, hence $| \partial^\gamma \theta^\hQ_Q| \leq C(A,\overline{r}) \cdot \delta_Q^{-|\gamma|} = C(A, \overline{r}) \cdot \delta_Q^{|\beta| - m}$. Hence,
\begin{align*}
|\partial^\alpha F(x)| & \leq C(A, \overline{r}) \sum_{\substack{Q \in \cQ \\ (1+\overline{r})Q \ni x}} \biggl[ |\partial^\alpha F_Q(x)| +  \delta_Q^{-m- \frac{n}{p}} \| F_Q - P_Q \|_{L^p((1+\overline{r})Q \cap \hQ)} + \delta_Q^{-\frac{n}{p}} \| F_Q \|_{\X((1+\overline{r})Q \cap \hQ)} \\
& \qquad  + \sum_{|\beta| \leq m-1} | \partial^\beta (P_Q - P_{Q'})(x) | \cdot \delta_Q^{|\beta| - m} \biggr]\\
& \qquad\qquad (\mbox{note that the cube} \; Q' \; \mbox{enters into the above sum})\\
& \leq C(A, \overline{r})  \sum_{\substack{Q \in \cQ \\ (1+\overline{r})Q \ni x}} \biggl[ |\partial^\alpha F_Q(x)| +  \delta_Q^{-m- \frac{n}{p}} \| F_Q - P_Q \|_{L^p((1+\overline{r})Q \cap \hQ)} + \delta_Q^{-\frac{n}{p}} \| F_Q \|_{\X((1+\overline{r})Q \cap \hQ)} \\
& \qquad  + \sum_{|\beta| \leq m-1} | \partial^\beta (P_Q - P_{Q'})(y_Q) | \cdot \delta_Q^{|\beta| - m} \biggr] \\
& \qquad (\mbox{note that } y_Q \in Q \mbox{ and } x \in (1+\overline{r})Q, \; \mbox{hence} \; |y_Q - x| \leq C \delta_Q; \\
& \qquad\qquad \; \mbox{thus, the above inequality follows from Lemma \ref{pnorm}} ).
\end{align*}
We now raise each side to the power $p$, integrate over $\hQ \cap (1+\frac{\overline{r}}{2})Q'$, and sum over $|\alpha| = m$. Thus we obtain
\begin{align*}
\| F\|_{\X(\hQ \cap (1+\frac{\overline{r}}{2})Q')}^p & \leq C(A, \overline{r})  \sum_{\substack{Q \in \cQ \\ (1+\overline{r})Q' \cap (1+\overline{r})Q \neq \emptyset}} \biggl[ \|F_Q\|_{\X((1+\overline{r})Q \cap \hQ)}^p +  \delta_Q^{-mp} \| F_Q - P_Q \|^p_{L^p((1+\overline{r})Q \cap \hQ)} \\
& \qquad\qquad \qquad +  \sum_{|\beta| \leq m-1} \lvert \partial^\beta (P_Q - P_{Q'})(y_Q) \rvert^p \cdot \delta_Q^{(|\beta| - m)p + n} \biggr].
\end{align*}
Finally, summing over $Q' \in \cQ$, we obtain the conclusion of the lemma, thanks to \eqref{covers} and the bounded overlap and good geometry of $\cQ$.
\end{proof}

\chapter{Proof of the Main Technical Results}

We will prove the Main Technical Results by induction on $\cA$  (see Chapter \ref{sec_mainresults}). Recall the order relation $<$ on multiindex sets $\cA \subset \cM$ defined in Section \ref{sec_multi}. In particular, recall that $\cA = \cM$ is minimal under $<$.

Fix a finite subset $E \subset \frac{1}{32}Q^\circ$, where $Q^\circ$ denotes the unit cube $[0,1)^n$. We assume that $N = \#(E) \geq 2$.

\section{Starting the Induction} \label{sec_start}

We first establish the base case of the induction. This corresponds to proving the Main Technical Results for $\cA = \cM$.  (See Chapter \ref{sec_mainresults}.) 

Let $\CZ(\cM)$ be the collection of maximal dyadic cubes $Q \subset Q^\circ$ such that $\#(E \cap 3Q) \leq 1$. 

Using one time-work at most $C N \log N$ in space $CN$, we produce a $\CZ(\cM)$-\textsc{Oracle} that answers queries as follows.
\begin{itemize}
\item A query consists of a point $\underline{x} \in Q^\circ$.
\item The response to the query $\underline{x}$ is a list of all the cubes $Q \in \CZ(\cM)$ such that $\underline{x} \in \frac{65}{64} Q$.
\item The work and storage required to answer a query are at most $C \log N$.
\end{itemize}
We simply apply the \textsc{Plain Vanilla CZ-Oracle} from Section \ref{sec_czoracle}; see Remark \ref{rem_cw}.

Since $\#(E) \geq 2$ and $E \subset Q^\circ$, the collection $\CZ(\cM)$ does not contain the cube $Q^\circ$. Therefore, each $Q \in \CZ(\cM)$ is a strict subcube of $Q^\circ$, hence $Q$ has a dyadic parent $Q^+ \subset Q^\circ$ such that $\#(3Q^+ \cap E) \geq 2$ (because $Q$ is maximal), and so in particular
\begin{equation}
\label{Enearby} \#(9Q \cap E) \geq 2 \quad \mbox{for all} \; Q \in \CZ(\cM).
\end{equation}

\begin{lem}\label{baselem1} If $Q,Q' \in \CZ(\cM)$ and $Q \leftrightarrow Q'$ then $\frac{1}{2} \delta_Q \leq \delta_{Q'} \leq 2 \delta_Q$.
\end{lem}
\begin{proof}
We proceed by contradiction. Suppose that $Q \leftrightarrow Q'$ and $\delta_Q \leq \frac{1}{4} \delta_{Q'}$ for some $Q,Q' \in \CZ(\cM)$. Then $3Q^+ \subset 3Q'$, and hence $\#(E \cap 3Q^+) \leq \#(E \cap 3 Q') \leq 1$. However, this contradicts that $\#(3Q^+ \cap E) \geq 2$, completing the proof of the lemma.
\end{proof}

\begin{lem}\label{baselem2} There exists $\epsilon_1 > 0$, depending only on $m,n,$ and $p$, such that $9 Q$ is not tagged with $(\cM, \epsilon_1)$ for any $Q \in \CZ(\cM)$.
\end{lem}
\begin{proof} Assume that $\epsilon_1 \in (0,1)$ is less than a small enough universal constant.

Let $Q \in \CZ(\cM)$.  It suffices to show that $\sigma(9Q)$ does not have an $(\cM,x_Q,\epsilon_1,\delta_{9Q})$-basis, thanks to \eqref{Enearby}. 

We argue by contradiction. Suppose that $(P_\alpha)_{\alpha \in \cM}$ is an $(\cM,x_Q,\epsilon_1,\delta_{9Q})$-basis for $\sigma(9Q)$. Therefore, $P_0(x_Q) = 1$, and $ \partial^\alpha P_0(x_Q) = 0$ for $\alpha \in \cM$, $\alpha \neq 0$. In other words, $P_0 \equiv 1$.

Moreover, there exists $\varphi_0 \in \X$ such that $\varphi_0 = 0$ on $E \cap 9 Q$ and
$$\| \varphi_0 \|_{\X(9Q)} + \delta_{9Q}^{-m} \| \varphi_0 - P_0 \|_{L^p(9Q)} \leq \epsilon_1 \delta_{9Q}^{n/p - m}.$$
We know that $\#(E \cap 9Q ) \geq 2$. Fix $x \in E \cap 9 Q$. By Lemma \ref{si2} we have
\[ \delta_Q^{n/p-m} \cdot \lvert \varphi_0(x) - P_0(x) \rvert \leq C \cdot \left\{ \| \varphi_0 \|_{\X(9Q)} + \delta_Q^{-m} \| \varphi_0 - P_0 \|_{L^p(9Q)} \right\} \leq C'  \epsilon_1 \delta_{9Q}^{n/p - m}.\]
But $\varphi_0(x) = 0$, and thus $\lvert P_0(x) \rvert \leq C'' \epsilon_1$. However, if we take $\epsilon_1 < 1/C''$, then this inequality contradicts the fact that $P_0 \equiv 1$.
\end{proof}

Recall that $\# (3 Q \cap E) \leq 1$ for each $Q \in \CZ(\cM)$. This implies the next result.

\begin{lem} \label{baselem3} If $Q \in \CZ(\cM)$ then $3Q$ is tagged with $(\cM,1/2)$.
\end{lem}

We have thus established properties \textbf{(CZ1-CZ5)} for the decomposition $\CZ(\cM)$. Indeed, \textbf{(CZ1)}, \textbf{(CZ2)}, and \textbf{(CZ4)} are consequences of Lemmas \ref{baselem1},  \ref{baselem2},  and \ref{baselem3}, respectively. Note that \textbf{(CZ3)} and \textbf{(CZ5)} are vacuously true because we are treating the base case $\cA = \cM$.

We next associate an extension operator and a linear functional to each of the ``non-trivial'' cubes in $\CZ(\cM)$.

More precisely, we define $\CZ_{\main}(\cM) := \left\{ Q \in \CZ(\cM) : (65/64)Q \cap E \neq \emptyset \right\}$. For each $Q \in \CZ_{\main}(\cM)$ there is a unique point $x(Q) \in E \cap \frac{65}{64}Q$. (Recall that $\#(E \cap 3Q) \leq 1$ for each $Q \in \CZ(\cM)$.)

For each $Q \in \CZ_{\main}(\cM)$, we define the following objects:
\begin{itemize}
\item A linear map $T_{(Q,\cM)} : \X(\frac{65}{64}Q \cap E) \oplus \cP \rightarrow \X$ given by
\begin{equation} \label{base1}
T_{(Q,\cM)}(f,P) = P + f(x(Q)) - P(x(Q)).
\end{equation}
\item A list $\Xi(Q,\cM) = \{ \xi_Q \}$, where
\begin{equation} \label{base2}
\xi_Q(f,P) = \bigl( f(x(Q)) - P(x(Q)) \bigr) \cdot \delta_Q^{n/p-m}.
\end{equation}
\item A list of assist functionals $\Omega(Q,\cM)$, which we take to be empty.
\end{itemize}

Clearly, the functional $\xi_Q$ and map $T_{(Q,\cM)}$ both have $\Omega(Q,\cM)$-assisted bounded depth (bounded depth).

\environmentA{Algorithm: Find Main-Cubes and Compute Extension Operators (Base Case).}

We compute a list of the cubes in $\CZ_{\main}(\cM)$. For each $Q  \in \CZ_{\main}(\cM)$, we compute a short form description of the bounded depth functional 
\[ \xi_Q : \X\left(\frac{65}{64}Q \cap E\right) \oplus \cP \rightarrow \R.\]

We give a query algorithm, which requires work at most $C \log N$ to answer queries. A query consists of a cube $\underline{Q} \in \CZ_{\main}(\cM)$ and point $\underline{x} \in Q^\circ$. The response to the query $(\underline{Q},\underline{x})$ is a short form description of the linear map
$$(f,P) \mapsto J_{\underline{x}} T_{(\underline{Q},\cM)}(f,P).$$

These computations require one-time work at most $C N \log N$ in space $C N$.

\begin{proof}[\underline{Explanation}] 

We compute a list of cubes $Q \in \CZ_{\main}(\cM)$ and associated points $x(Q) \in E \cap \frac{65}{64}Q$. This computation requires work at most $C N \log N$ in space $CN$; see the algorithm \textsc{Find Main-Cubes} in Section \ref{czalg_sec}.

For each $Q$ in the list $\CZ_{\main}(\cM)$, we compute the linear functional
\begin{align*}
&\xi_Q(f,P) = \bigl\{ f(x(Q)) - P(x(Q)) \bigr\} \cdot \delta_Q^{n/p-m}.
\end{align*}
There are at most $C N$ such functionals, and we compute each one using work and storage at most $C$.

Given $(\underline{Q}, \underline{x}) \in \CZ_{\main}(\cM) \times Q^\circ$, we use a binary search to determine the position of $\underline{Q}$ in the list $\CZ_{\main}$. We then compute the linear map
\[ (f,P) \mapsto J_{\underline{x}}T_{(\underline{Q},\cM)}(f,P) = P + f(x(\underline{Q})) - P(x(\underline{Q})).
\]
This requires work at most $C \log N$ per query.
\end{proof}

\begin{lem} \label{lem_eo} There exists $C \geq 1$, depending only on $m,n,$ and $p$, such that for each $Q \in \CZ_{\main}(\cM)$, the following properties hold.
\begin{itemize}
\item $T_{(Q,\cM)}(f,P) = f$ on $\frac{65}{64}Q \cap E$.
\vspace{0.2cm}
\item $\| T_{(Q,\cM)}(f,P) \|_{\X(\frac{65}{64}Q)} + \delta_Q^{-m} \| T_{(Q,\cM)}(f,P) - P \|_{L^p(\frac{65}{64}Q)} \leq C \cdot \lv \xi_Q(f,P) \rv$.
\vspace{0.2cm}
\item $C^{-1} \cdot \| (f,P) \|_{\frac{65}{64}Q} \leq \lvert \xi_Q(f,P) \rvert \leq C  \cdot \| (f,P) \|_{\frac{65}{64}Q}.$
\end{itemize}
\end{lem}
\begin{proof}
Note that $E \cap \frac{65}{64}Q = \{x(Q) \}$ and $T_{(Q,\cM)}(f,P)(x(Q)) = f(x(Q))$ for each $Q \in \CZ_{\main}(\cM)$. This implies the first bullet point.

Recall that $T_{(Q,\cM)}(f,P) \in \cP$, hence $\| T_{(Q,\cM)}(f,P) \|_{\X(\frac{65}{64}Q)} = 0$. Moreover,
\begin{align*}
\delta_Q^{-m} \| T_{(Q,\cM)}(f,P) - P \|_{L^p(\frac{65}{64}Q)} &= \delta_Q^{-m} \| f(x(Q)) - P(x(Q)) \|_{L^p(\frac{65}{64}Q)} \\
&\leq C \delta_Q^{-m + n/p} \lvert f(x(Q)) - P(x(Q)) \rvert = C \lvert \xi_Q(f,P) \rvert.
\end{align*}
This implies the second bullet point.

From the first and second bullet points we have
\[ \| (f,P) \|_{\frac{65}{64}Q} \leq \| T_{(Q,\cM)}(f,P) \|_{\X(\frac{65}{64}Q)} + (\delta_{\frac{65}{64}Q})^{-m} \| T_{(Q,\cM)}(f,P) - P \|_{L^p(\frac{65}{64}Q)} \leq C\lvert \xi_Q(f,P) \rvert.
\]

Let $F \in \X$ satisfy $F = f$ on $\frac{65}{64} Q \cap E$. Then the Sobolev inequality implies that
\[
\delta_Q^{n/p - m} \lvert (f-P)(x(Q)) \rvert = \delta_Q^{n/p - m}  \lvert (F - P)(x(Q)) \rvert \leq C \cdot \left( \| F\|_{\X(\frac{65}{64}Q)} + \delta_Q^{-m} \| F - P \|_{L^p(\frac{65}{64}Q)} \right)
\]
Taking the infimum over such $F$, we obtain the estimate $\lvert \xi_Q(f,P) \rvert \leq C \| (f,P) \|_{\frac{65}{64}Q}$. Thus we obtain the third bullet point, and this completes the proof of the lemma.
\end{proof}

This completes the proof of the base case of the induction. In the next section we start to prove the induction step.

\section{The Induction Step} 
\label{sec_ind_step}
Fix a set of multiindices $\cA \subset \cM$ with 
\begin{equation} \label{not_base}
\cA \neq \cM.
\end{equation} 
We assume by induction that we have already carried out the Main Technical Results for each $\cA' < \cA$. Our goal is to find suitable constants $a(\cA)$, $\epsilon_1(\cA)$, $\epsilon_2(\cA)$, $c_*(\cA)$, $S(\cA)$ and to carry out the Main Technical Results for $\cA$.

Let $\cA^- < \cA$ be the maximal mutiindex set with respect to the order relation $<$ on $2^\cM$. (See Section \ref{sec_multi} for the definition of the order relation.) By induction hypothesis, we have already carried out the Main Technical Results for $\cA^-$. (See Chapter \ref{sec_mainresults}.) We have thus produced the following:
\begin{itemize}
\item A decomposition $\CZ(\cA^-)$ of $Q^\circ$ into dyadic cubes, with the following properties.
\begin{itemize}
\item If $\frac{65}{64}Q \cap \frac{65}{64}Q' \neq \emptyset$ with $Q,Q' \in \CZ(\cA^-)$, then $Q \leftrightarrow Q'$ and $\frac{1}{2} \delta_{Q'} \leq \delta_Q \leq 2 \delta_{Q'}$.
\item The collection of cubes $\{ \frac{65}{64}Q : Q \in \CZ(\cA^-) \}$ has \emph{bounded overlap}, meaning that there exists a constant $C=C(n)$ such that, for each $Q \in \CZ(\cA^-)$ there are at most $C$ cubes $Q' \in \CZ(\cA^-)$ with $(65/64)Q \cap (65/64)Q' \neq \emptyset$.
\item From \eqref{Enearby} and since $\CZ(\cM)$ refines $\CZ(\cA^-)$ (see the induction hypothesis) we know that
\begin{equation}
\label{Enearby2}
(\text{``}E\text{ is nearby''}) \;\; \#(E \cap 9Q) \geq 2 \; \mbox{for each} \; Q \in \CZ(\cA^-).
\end{equation}
\end{itemize}
\item An oracle that accepts queries $\underline{x} \in Q^\circ$ and returns the list of all cubes $Q \in \CZ(\cA^-)$ such that $\underline{x} \in \frac{65}{64}Q$.
\item A list $\CZ_{\main}(\cA^-)$ consisting of all the $Q \in \CZ(\cA^-)$ such that $\frac{65}{64}Q \cap E \neq \emptyset$.
\item For each $Q \in \CZ_{\main}(\cA^-)$, a list of assists $\Omega(Q,\cA^-) \subset \left[ \X( E) \right]^*$.
\item For each $Q \in \CZ_{\main}(\cA^-)$, a list of $\Omega(Q,\cA^-)$-assisted bounded depth linear functionals $\Xi(Q,\cA^-) \subset \left[ \X(\frac{65}{64}Q \cap E) \oplus \cP \right]^*$ written in short form, as well as a linear extension operator $$T_{(Q,\cA^-)} : \X\left(\frac{65}{64}Q \cap E\right) \oplus \cP \rightarrow \X,$$
which we ``compute'' in the sense that (after one-time work) we can answer queries: In response to a query $\underline{x} \in Q^\circ$ we return a short form description of the $\Omega(Q,\cA^-)$-assisted bounded depth linear map
$$(f,P) \mapsto J_{\underline{x}} T_{(Q,\cA^-)}(f,P).$$
\end{itemize}
These objects and algorithms have good properties as part of the induction assumption on $\cA^-$. We listed some of these properties just above. The remaining properties are mentioned later, as required.

We denote
\begin{equation}
\label{a_defn}
a = a(\cA^-), \; \mbox{the geometric constant used in the Main Technical Results for} \; \cA^-.
\end{equation}

\environmentA{Algorithm: Approximate Old Trace Norm.} 

For each $Q \in \CZ_{\main}(\cA^-)$, we compute linear functionals $\xi_1^Q,\ldots,\xi_D^Q$ on $\cP$, such that
\begin{equation}\label{ineq2}  \sum_{\xi \in \Xi(Q,\cA^-)} \lvert \xi(0,P) \rvert^p \quad \mbox{and} \quad \sum_{i=1}^D \lvert \xi^Q_i(P) \rvert^p \qquad (P \in \cP)
\end{equation}
differ by at most a factor of $C$. We carry this out using work and storage $\leq C N$.

\begin{proof}[\underline{Explanation}]

For each $\xi$ in the list $\Xi(Q,\cA^-)$, we compute the map $P \mapsto \xi(0,P)$ using work and storage at most $C$, by examining the short form description of $(f,P) \mapsto \xi(f,P)$ that has been computed. Applying \textsc{Compress Norms} (see Section \ref{sec_lf}), we compute linear functionals $\xi^Q_1, \cdots, \xi^Q_D$ such that \eqref{ineq2} holds, using work and storage at most $C \cdot \# \left[ \Xi(Q,\cA^-) \right]$. By the inductive hypothesis, we know that the sum of $\# \left[ \Xi(Q,\cA^-) \right]$ over all $Q \in \CZ_{\main}(\cA^-)$ is bounded by $C N$, hence the work and storage guarantees are met.

\end{proof}

\subsection{The Non-monotonic Case}\label{nonmon_sec}

Here, we assume that $\cA \subset \cM$ is not monotonic and prove the Main Technical Results for $\cA$. See Section \ref{sec_multi} for the definition of monotonic sets.

We define $\CZ(\cA) = \CZ(\cA^-)$ and 
\[
\epsilon_2(\cA) = \epsilon_2(\cA^-), \;\; c_*(\cA) = c_*(\cA^-), \;\; a(\cA) = a(\cA^-), \;\; \mbox{and} \; S(\cA) = S(\cA^-).
\]
The constant $\epsilon_1(\cA)$ is chosen later in this section.

We define 
\begin{align*}
&\Omega(Q,\cA) := \Omega(Q,\cA^-), \;\; \Xi(Q,\cA) := \Xi(Q,\cA^-) \;  \mbox{and} \\ 
&T_{(Q,\cA)} := T_{(Q,\cA^-)} \; \;  \mbox{for each} \; Q \in \CZ_{\main}(\cA) = \CZ_{\main}(\cA^-).
\end{align*}
The properties of $\Omega(Q,\cA)$, $\Xi(Q,\cA)$ and $T_{(Q,\cA)}$ asserted in the Main Technical Results for $\cA$ are immediate from the corresponding properties of $\Omega(Q,\cA^-)$, $\Xi(Q,\cA^-)$ and $T_{(Q,\cA^-)}$ asserted in the Main Technical Results for $\cA^-$.

Next, we prove properties \textbf{(CZ1-CZ5)} for the label $\cA$.

Note that \textbf{(CZ1)} for $\cA$ follows from \textbf{(CZ1)} for $\cA^-$. Also, note that \textbf{(CZ5)} for $\cA$ holds because $\CZ(\cA) = \CZ(\cA^-)$.

Note that \textbf{(CZ3)} for $\cA$ holds vacuously: There do not exist cubes $Q \in \CZ(\cA)$, $Q' \in \CZ(\cA^-)$ which satisfy the hypotheses of \textbf{(CZ3)}. This follows because $\CZ(\cA) = \CZ(\cA^-)$.

We need not check \textbf{(CZ4)}, since $\cA \neq \cM$; see \eqref{not_base}.

It remains to prove \textbf{(CZ2)} for $\cA$, which we accomplish in the next lemma. We determine $\epsilon_1(\cA) = \epsilon_1$ in the lemma below.

\begin{lem} \label{secondprop_lem}
There exists a universal constant $\epsilon_1>0$ such that the following holds. Suppose that $Q \in \CZ(\cA)$ and $\delta_Q \leq c_*(\cA)$. Then $S(\cA) Q$ is not tagged with $(\cA,\epsilon_1)$.
\end{lem}

\begin{proof}

We assume that $\epsilon_1 > 0$ is less than a small enough universal constant.

Let $Q \in \CZ(\cA)$ satisfy $\delta_Q \leq c_*(\cA)$. Assume for the sake of contradiction that $S(\cA)Q$ is tagged with $(\cA,\epsilon_1)$.

If $\#(S(\cA)Q \cap E) \leq 1$ then $S(\cA^-)Q = S(\cA)Q$ is tagged with $(\cA^-,\epsilon_1(\cA^-))$. However, this contradicts the induction hypothesis. Hence, we may assume from now on that $\#(S(\cA) Q  \cap E) \geq 2$. Thus, 
\[ \sigma(S(\cA)Q) \; \mbox{has an} \; (\widetilde{\cA},x_Q,\epsilon_1,\delta_{S(\cA)Q})\mbox{-basis for some} \; \widetilde{\cA} \leq \cA.
\]
Hence, Lemma \ref{pre_lem2} implies that there exists $\kappa \in [ \kappa_1,\kappa_2] $ such that
\[ \sigma(S(\cA)Q) \; \mbox{has an} \; (\cA',x_Q,\epsilon_1^{\kappa},\delta_{S(\cA)Q}, \Lambda)\mbox{-basis}, \; \mbox{with} \;\cA' \leq \cA \; \mbox{and} \; \epsilon_1^{\kappa} \Lambda^{100D} \leq \epsilon_1^{\kappa/2}.\]
Here, $\kappa_1, \kappa_2 > 0$ are universal constants. 

Suppose for the moment that $\cA' < \cA$. Then $S(\cA)Q$ is tagged with $(\cA^-,\epsilon_1^{\kappa})$. Note that $\epsilon_1^{\kappa} \leq \epsilon_1^{\kappa_1} \leq \epsilon_1(\cA^-)$, for small enough $\epsilon_1$. Thus, $S(\cA^-)Q = S(\cA)Q$ is tagged with $(\cA^-,\epsilon_1(\cA^-))$. However, this contradicts the induction hypothesis. Hence,
\[ \sigma(S(\cA)Q) \; \mbox{has an} \; (\cA,x_Q,\epsilon_1^{\kappa},\delta_{S(\cA)Q}, \Lambda)\mbox{-basis}.\]
Thus, there exists $(P_\alpha)_{\alpha \in \cA}$ with
\begin{equation} \label{insigma}
P_\alpha \in \epsilon_1^{\kappa} \cdot  (\delta_{S(\cA) Q})^{\frac{n}{p} + |\alpha| - m }  \sigma(S(\cA)Q) \qquad (\alpha \in \cA)
\end{equation}
\begin{itemize}
\item $\partial^\beta P_\alpha (x_Q) = \delta_{\beta \alpha} \qquad\qquad\qquad\qquad\qquad (\beta,\alpha \in \cA)$

\item $\lvert \partial^\beta P_\alpha(x_Q) \rvert \leq \epsilon_1^\kappa \cdot (\delta_{S(\cA)Q})^{|\alpha| - | \beta| } \qquad\quad (\alpha \in \cA, \beta \in \cM, \beta > \alpha)$

\item $\lvert \partial^\beta P_\alpha(x_Q) \rvert \leq \Lambda \cdot (\delta_{S(\cA)Q})^{|\alpha| - | \beta| } \qquad\quad (\alpha \in \cA, \beta \in \cM)$.
\end{itemize}

We are assuming that $\cA$ is not monotonic. Thus we can pick $\alpha_0 \in \cA$ and $\gamma \in \cM$ such that $\alpha_0 + \gamma \in \cM \setminus \cA$. We define
\[\oa = \alpha_0 + \gamma \quad \mbox{and} \quad \oA = \cA \cup \{\oa\}.\]
Note that $\oA \Delta \cA = \{\oa\}$ with $\oa \in \oA$. Consequently, $\oA < \cA$.

We define $P_\oa = P_{\alpha_0} \odot_{x_Q} q$, where $q(y) = \frac{\alpha_0 !}{\oa !} (y - x_Q)^\gamma $. That is,
\[P_\oa(y) = \frac{\alpha_0!}{\oa!} \sum_{|\omega| \leq m - 1 - |\gamma| } \frac{1}{\omega!} \partial^\omega P_{\alpha_0}(x_Q) (y - x_Q)^{\omega+ \gamma}.\]
Note that $P_\oa = q \cdot P_{\alpha_0}^{\main}$, where
\[P_{\alpha_0}^{\main} = \sum_{|\omega| \leq m - 1 - |\gamma|} \frac{1}{\omega!} \partial^\omega P_{\alpha_0}(x_Q) (y-x_Q)^\omega, \]
and that
\[R_{\alpha_0}  := P_{\alpha_0} -  P_{\alpha_0}^{\main}  = \sum_{|\omega| > m - 1 - |\gamma|} \frac{1}{\omega!} \partial^\omega P_{\alpha_0}(x_Q) (y-x_Q)^\omega.\]
In the above sum for $R_{\alpha_0}$, since $|\omega| > m-1 - |\gamma| \geq |\alpha_0|$ we have $\omega > \alpha_0$, and so $\lvert \partial^\omega P_{\alpha_0}(x_Q) \rvert \leq C \epsilon_1^\kappa \delta_Q^{|\alpha_0| - |\omega|}$. Consequently, $\| R_{\alpha_0} \|_{L^p(S(\cA)Q)} \leq C' \epsilon_1^\kappa \delta_Q^{n/p + |\alpha_0|}$.

The bullet point properties of $P_{\alpha_0}$ now yield the following properties of $P_\oa$. 
\begin{itemize}
\item $\partial^\oa P_\oa (x_Q) = 1$
\item $\lvert \partial^\beta P_\oa(x_Q) \rvert \leq C \epsilon_1^\kappa \delta_Q^{|\oa| - |\beta|} \qquad\quad (\beta \in \cM, \beta > \oa)$
\item $\lvert \partial^\beta P_\oa(x_Q) \rvert \leq  C \Lambda \delta_Q^{|\oa| - |\beta|}  \qquad\quad (\beta \in \cM)$.
\end{itemize}
We now show that 
\begin{itemize}
\item $P_{\oa} \in C \epsilon_1^\kappa \cdot (\delta_Q)^{\frac{n}{p} + |\oa | - m} \cdot \sigma(S(\cA)Q)$.
\end{itemize}
To start, \eqref{insigma} implies that there exists $\varphi \in \X$ with $\varphi = 0$ on $S(\cA) Q \cap E$ and
\[\| \varphi \|_{\X(S(\cA)Q)} + \delta_Q^{-m} \| \varphi - P_{\alpha_0} \|_{L^p(S(\cA)Q)}  \leq C \epsilon_1^\kappa \cdot (\delta_Q)^{\frac{n}{p} + |\alpha_0| - m }.\]
Applying $\| R_{\alpha_0} \|_{L^p(S(\cA)Q)} \leq C \epsilon_1^\kappa \delta_Q^{\frac{n}{p} + |\alpha_0|}$, we see that
\[ \| \varphi - P_{\alpha_0}^{\main} \|_{\X(S(\cA)Q)} + \delta_Q^{-m} \| \varphi - P^{\main}_{\alpha_0} \|_{L^p(S(\cA)Q)}  \leq C \epsilon_1^\kappa \cdot (\delta_Q)^{\frac{n}{p} + |\alpha_0| - m }.\]
Moreover, the Leibniz Rule shows that
\begin{align*}
\| q\cdot ( \varphi - P_{\alpha_0}^{\main}) & \|_{\X(S(\cA)Q)}  + \delta_Q^{-m} \| q \cdot (\varphi  - P^{\main}_{\alpha_0}) \|_{L^p(S(\cA)Q)} \\
& \leq C\sum_{k=0}^{m}  (\delta_Q)^{k + |\gamma| - m }  \| \nabla^k (\varphi - P_{\alpha_0}^{\main}) \|_{L^p(S(\cA)Q)} \\
& \leq C \cdot (\delta_Q)^{|\gamma|} \left( \|\varphi -  P_{\alpha_0}^{\main} \|_{\X(S(\cA)Q)} + \delta_Q^{-m} \| \varphi - P_{\alpha_0}^{\main} \|_{L^p(S(\cA)Q)} \right) \\
& \qquad\qquad\qquad\qquad\qquad\qquad (\mbox{by Lemma} \; \ref{si2}) \\
& \leq C \epsilon_1^\kappa \cdot (\delta_Q)^{|\gamma|} (\delta_Q)^{\frac{n}{p} + |\alpha_0| - m } = C \epsilon_1^\kappa \cdot (\delta_Q)^{\frac{n}{p} + |\oa| - m }.
\end{align*}
Note that $q \cdot P_{\alpha_0}^{\main} \in \cP$, hence
\[\| q\cdot \varphi  \|_{\X(S(\cA)Q)} + \delta_Q^{-m} \| q \cdot \varphi  - P_{\oa} \|_{L^p(S(\cA)Q)} \leq C \epsilon_1^\kappa (\delta_Q)^{\frac{n}{p} + |\oa| - m }. \]
Since $q \cdot \varphi = 0$ on $S(\cA)Q \cap E$, we have shown that $P_{\oa} \in C \epsilon_1^\kappa \cdot (\delta_Q)^{\frac{n}{p} + |\oa | - m} \cdot \sigma(S(\cA)Q)$. This proves all the bullet point properties of $P_\oa$.

The bullet point properties of the $P_\alpha$ ($\alpha \in \oA$) imply that $(\partial^\beta P_\alpha(x_Q))_{\alpha,\beta \in \oA}$ is $(C \epsilon_1^\kappa, C \Lambda,\delta_Q)$-near triangular. Inverting the matrix $(\partial^\beta P_\alpha(x_Q))_{\alpha,\beta \in \oA}$, we obtain a matrix $(M_{\alpha \omega})_{\alpha,\omega \in \oA}$ such that
\[\sum_{\alpha \in \oA} \partial^\beta P_\alpha(x_Q) M_{\alpha \omega} = \delta_{\beta \omega} \;\;\; (\beta,\omega \in \oA)\]
and
\[
\lvert M_{\alpha \omega} - \delta_{\alpha \omega} \rvert \leq
\left\{
\begin{array}{ll}
C \epsilon_1^\kappa \Lambda^D \cdot \delta_Q^{|\omega| - |\alpha|} &: \mbox{if}\; \alpha, \omega \in \oA, \; \alpha \geq \omega \\
C \Lambda^D \cdot \delta_Q^{|\omega| - |\alpha|}  &: \mbox{if} \; \alpha, \omega \in \oA.\\
\end{array} \right. \]

Set $P^\#_\omega = \sum_{\alpha \in \oA} P_\alpha M_{\alpha \omega}$. The bullet point properties of $(P_\alpha)_{\alpha \in \oA}$ imply that
\begin{itemize}
\item $ P^\#_\omega \in C \epsilon_1^\kappa \cdot \Lambda^{2D} \cdot \delta_Q^{n/p + |\omega| - m} \sigma(S(\cA) Q)$ \qquad $(\omega \in \oA)$
\item $ \partial^\beta P^\#_\omega(x_Q) = \delta_{\beta \omega}$ \qquad\qquad\qquad\qquad $(\beta,\omega \in \oA)$
\end{itemize}
For $\omega \in \oA$ and $\beta \in \cM$ with $\beta > \omega$, we write
\begin{equation} \label{broken} \partial^\beta P^\#_\omega(x_Q) = \sum_{\alpha < \beta} \partial^\beta P_\alpha(x_Q) M_{\alpha \omega} + \sum_{\alpha \geq \beta} \partial^\beta P_\alpha(x_Q) M_{\alpha \omega}.\end{equation}
An arbitrary term in the first sum in \eqref{broken} is bounded by $\left[ C \epsilon_1^\kappa \delta_Q^{|\alpha| - |\beta|} \right] \cdot \left[ C \Lambda^D \delta_Q^{|\omega| - |\alpha|} \right]$. Hence, this sum is at most $ C' \epsilon_1^\kappa \Lambda^D \delta_Q^{|\omega| - |\beta|}$.

If $\alpha \geq \beta$, then $\alpha > \omega$, since $\beta > \omega$. Thus, an arbitrary term in the second sum in \eqref{broken} is bounded by $\left[ C \Lambda \delta_Q^{|\alpha| - |\beta|} \right] \cdot \left[ C \epsilon_1^\kappa \Lambda^{D} \cdot \delta_Q^{|\omega| - |\alpha|} \right]$. Hence, this sum is at most $C' \epsilon_1^\kappa \Lambda^{D+1} \delta_Q^{|\omega| - |\beta|}$.

Thus,
\begin{itemize}
\item $ \lvert \partial^\beta P^\#_\omega(x_Q) \rvert \leq C \epsilon_1^\kappa \Lambda^{2D} \delta_Q^{|\omega| - |\beta|}$ \qquad\quad $(\beta \in \cM, \omega \in \oA, \beta > \omega)$.
\end{itemize}

According to the bullet point properties of $(P^\#_\omega)_{\omega \in \oA}$, we see that \\
$\sigma(S(\cA)Q)$ has an $(\oA, x_Q, C \epsilon^\kappa_1 \Lambda^{2D}, \delta_Q)$-basis, hence \\
$\sigma(S(\cA)Q)$ has an $(\oA,x_Q,C'\epsilon^\kappa_1 \Lambda^{2D}, \delta_{S(\cA)Q})$-basis. (See Remark \ref{basis_rem}.)

For small enough $\epsilon_1$ we have $C' \epsilon_1^\kappa \Lambda^{2D} \leq C' \epsilon_1^{\kappa/2} \leq \epsilon_1^{\kappa_1/4} \leq \epsilon_1(\cA^-)$, hence
\[\sigma(S(\cA)Q) \; \mbox{has an} \; (\oA, x_Q, \epsilon_1(\cA^-), \delta_{S(\cA)Q})\mbox{-basis}.\]
Hence, $S(\cA)Q$ is tagged with $(\cA^-, \epsilon_1(\cA^-))$. However, since $\delta_Q \leq c_*(\cA^-)$ and $S(\cA) = S(\cA^-)$, this contradicts the induction hypothesis.

This completes the contradiction, and with it, the proof of the lemma.

\end{proof}

We have thus proven the Main Technical Results for $\cA$ in the non-monotonic case.

\subsection{The Monotonic Case} \label{sec_ocz}

From this point onward, we assume that $\cA$ is monotonic. (See Section \ref{sec_multi} for the definition of monotonic multiindex sets.) We drop this assumption when we prove our main theorem in Chapter \ref{final_chap}. We will now begin the task of carrying out the induction step by proving the Main Technical Results for $\cA$. (See Chapter \ref{sec_mainresults}.)

We begin by treating a preliminary case.

\begin{lem} \label{mini_lem} Suppose that $\delta_Q \geq \frac{1}{4}$ for all $Q \in \CZ(\cA^-)$. Then the Main Technical Results for $\cA^-$ imply the Main Technical Results for $\cA$.
\end{lem} 
\begin{proof}
We take $\CZ(\cA)$ to equal $\CZ(\cA^-)$. The other objects and algorithms in the Main Technical Results for $\cA$ are copies of the corresponding objects and algorithms in the Main Technical Results for $\cA^-$. 
\end{proof}

By making at most $C$ calls to the $\CZ(\cA^-)$-\textsc{Oracle}, we can check whether the hypothesis of Lemma \ref{mini_lem} holds. This takes one-time work at most $C \log N$. In the sequel, we assume that we are in the case that
\begin{equation} \label{smallish_cube}
\delta_Q \leq 1/8 \;\; \mbox{for some} \; Q \in \CZ(\cA^-).
\end{equation} 
Recall that the decomposition $\CZ(\cA^-)$ has the following properties:
\begin{itemize}
\item $\CZ(\cA^-)$ is a finite partition of $Q^\circ = [0,1)^n$ into pairwise disjoint dyadic subcubes.
\item If $Q,Q' \in \CZ(\cA^-)$ and $Q \leftrightarrow Q'$ then $\delta_Q/ \delta_{Q'} \in \{1/2,1,2\}$.
\item If $Q \in \CZ(\cA^-)$ then $\#(9Q \cap E) \geq 2$. (See \eqref{Enearby2}.)
\end{itemize}

\begin{lem}\label{bdry_lem}
If $Q \in \CZ(\cA^-)$ and $\dist(Q, \R^n \setminus Q^\circ) = 0$ then $\delta_Q \in \{ \frac{1}{2}, \frac{1}{4}, \frac{1}{8} \}$.
\end{lem}
\begin{proof}
Let $Q \in \CZ(\cA^-)$ with $\dist(Q, \R^n \setminus Q^\circ) = 0$. 

Recall that $\delta_{Q} \neq 1$, because $\CZ(\cA^-) \neq \{ Q^\circ \}$ (see \eqref{smallish_cube}).

We need to show that $\delta_Q \geq \frac{1}{8}$. For the sake of contradiction assume that $\delta_Q \leq \frac{1}{16}$. Then since $\dist(Q, \R^n \setminus Q^\circ) = 0$, we have $9 Q \subset \R^n \setminus \frac{1}{10}Q^\circ$, hence $9Q \subset \R^n \setminus E$. But $\#( E \cap 9Q) \geq 2$, according to the above bullet points. This contradiction completes the proof of Lemma \ref{bdry_lem}.
\end{proof}

We now pass from the decomposition $\CZ(\cA^-)$ of $Q^\circ$ to a decomposition $\overline{\CZ}(\cA^-)$ of $\R^n$.

\begin{prop} \label{newcz_prop}
There exists a decomposition $\overline{\CZ}(\cA^-)$ of $\R^n$ into pairwise disjoint dyadic cubes, with the following properties:
\begin{enumerate}[(a)]
\item $\CZ(\cA^-) \subset \overline{\CZ}(\cA^-)$.
\item If $Q,Q' \in \overline{\CZ}(\cA^-)$ and $Q \leftrightarrow Q'$ then $\frac{1}{8} \delta_{Q'} \leq \delta_Q \leq 8 \delta_{Q'}$ (``good geometry''). Moreover, the collection of cubes $\{ \frac{65}{64} Q : Q \in \overline{\CZ}(\cA^-)\}$ has bounded overlap (each cube intersects a bounded number of other cubes).
\item If $Q \in \overline{\CZ}(\cA^-) \setminus \CZ(\cA^-)$, then $\frac{65}{64}Q \cap E = \emptyset$.
\item If $Q \in  \overline{\CZ}(\cA^-) \setminus \CZ(\cA^-)$, then $100 Q$ intersects cubes in $\CZ(\cA^-)$ with sidelength less than $\delta_Q$.
\item If $Q \in \overline{\CZ}(\cA^-) \setminus \CZ(\cA^-)$, then $\delta_Q \geq 1$.
\item If $Q \in \overline{\CZ}(\cA^-)$ then $\#(9Q \cap E) \geq 2$.
\end{enumerate}
We produce a $\overline{\CZ}(\cA^-)$-\textsc{Oracle}. The $\overline{\CZ}(\cA^-)$-\textsc{Oracle} accepts a query consisting of a point $\underline{x} \in \R^n$. The response to a query $\underline{x}$ is the list of cubes $Q \in \overline{\CZ}(\cA^-)$ such that $\underline{x} \in \frac{65}{64}Q$. The work and storage required to answer a query are at most $C \log N$.

\end{prop}
\begin{proof}
Let $\cQ$ consist of the maximal dyadic cubes $Q \subset \R^n$ satisfying the condition [$\delta_Q \leq 1$ or $0 \notin 2Q$]. A dyadic cube $Q \subset \R^n$ belongs to $\cQ$ if and only if
\begin{equation}\label{qq6}
\delta_Q = 1 \; \mbox{or} \; 0 \notin 2Q,
\end{equation} 
and
\begin{equation}\label{qq7}
\delta_{Q^+} \geq 2 \; \mbox{and} \; 0 \in 2Q^+.
\end{equation}
Here, as usual, $Q^+$ denotes the parent of a dyadic cube $Q$.

For any $x \in \R^n$, there exists a dyadic cube $\overline{Q}$ containing $x$ such that $\delta_{\overline{Q}} \geq 2$ and $0 \in 2\overline{Q}$. Hence, each $x \in \R^n$ is contained in some cube $Q \in \cQ$. Hence, $\cQ$ partitions $\R^n$ into pairwise disjoint dyadic cubes.

Note that the cube $Q^\circ = [0,1)^n$ belongs to $\cQ$.

We now establish good geometry of $\cQ$ (with constant $1/4$). We prove that if $Q,Q' \in \cQ$ and $Q \leftrightarrow Q'$ then $\frac{1}{4} \delta_{Q'} \leq \delta_Q \leq 4 \delta_{Q'}$. 

Assume for the sake of contradiction that there exist cubes $Q,Q' \in \cQ$ with $\delta_{Q} \leq \frac{1}{8} \delta_{Q'}$ and $Q \leftrightarrow Q'$. By \eqref{qq7}, we have $\delta_{Q^+} \geq 2$ and $0 \in 2 Q^+$. Moreover, note that $2Q^+ \subset 2 Q'$ (since $Q \leftrightarrow Q'$ and $\delta_{Q} \leq \frac{1}{8} \delta_{Q'}$, hence $Q^+ \leftrightarrow Q'$ and $\delta_{Q^+} \leq \frac{1}{4} \delta_{Q'}$). Hence, $0 \in 2Q'$. Moreover, $\delta_{Q'} \geq 4 \delta_{Q^+} \geq 8$. However, since $Q' \in \cQ$, the analogue of \eqref{qq6} with $Q$ replaced by $Q'$ must hold. This yields a contradiction. This completes the proof that the cubes in $\cQ$ have good geometry.

We define the collection $\overline{\CZ}(\cA^-)$ to consist of all the cubes $Q \in \cQ$ except for $Q = Q^\circ$, together with all the cubes $Q \in \CZ(\cA^-)$. Since $\cQ$ partitions $\R^n$ and $\CZ(\cA^-)$ partitions $Q^\circ$, we see that $\overline{\CZ}(\cA^-)$ partitions $\R^n$ into pairwise disjoint dyadic cubes. Moreover, property (a) clearly holds.

If $Q \in \cQ$, $Q' \in \CZ(\cA^-)$, and $Q \leftrightarrow Q'$, then both $Q$ and $Q'$ touch the boundary of $Q^\circ$.

We prove the claim that $\cQ$ contains all $4^{n}$ of the dyadic cubes $Q \subset [-2,2)^n$ with $\delta_Q = 1$. Indeed, we have $Q^+ \subset [-2,2)^n$, $\delta_{Q^+} = 2$ and $0 \in 2Q^+$ for any such $Q$. Hence, each $Q$ satisfies \eqref{qq6} and \eqref{qq7}, which implies that $Q$ belongs to $\cQ$. This proves our claim. Hence, in particular, any $Q \in \cQ$ that intersects the boundary of $Q^\circ = [0,1)^n$ must satisfy $\delta_Q = 1$.

Moreover, by Lemma \ref{bdry_lem}, any $Q' \in \CZ(\cA^-)$ that intersects the boundary of $Q^\circ$ must satisfy $\delta_{Q'} \in \{1/2,1/4,1/8\}$.

Hence, the previous two statements imply that for any $Q \in \cQ$ and $Q' \in \CZ(\cA^-)$ with $Q \leftrightarrow Q'$ we have $\frac{1}{8} \delta_Q \leq \delta_{Q'} \leq \delta_Q$.

Finally, for $Q,Q' \in \CZ(\cA^-)$ with $Q \leftrightarrow Q'$, we have $\frac{1}{2} \delta_Q \leq \delta_{Q'} \leq 2 \delta_Q$, by good geometry of the cubes in $\CZ(\cA^-)$.

Recall that the cubes in $\cQ$ satisfy good geometry (with constant $1/4$).

Thus, combining the previous three statements, for any $Q,Q' \in \overline{\CZ}(\cA^-)$ with $Q \leftrightarrow Q'$, we have $\frac{1}{8} \delta_{Q'} \leq \delta_Q \leq 8 \delta_{Q'}$.

The above property shows that $\overline{\CZ}(\cA^-)$ satisfies the hypothesis of Lemma \ref{gg_lem_1} with $\gamma=1/8$. Hence, for $Q,Q' \in \overline{\CZ}(\cA^-)$ with $\frac{65}{64}Q \cap \frac{65}{64}Q' \neq \emptyset$, we have $Q \leftrightarrow Q'$. It follows that the collection $\{ \frac{65}{64} Q : Q \in \overline{\CZ}(\cA^-)\}$ has bounded overlap. This completes the proof of property (b).

From \eqref{qq7}, each $Q \in \overline{\CZ}(\cA^-) \setminus \CZ(\cA^-)$ satisfies $\delta_Q \geq 1$. This proves property (e).

We now prove property (c). Let $Q \in \overline{\CZ}(\cA^-) \setminus \CZ(\cA^-)$. Then $Q \in \cQ$ and $Q \subset \R^n \setminus Q^\circ$. According to property (e), there are only two cases to consider
\begin{itemize}
\item If $\delta_Q = 1$, then $\frac{65}{64}Q$ cannot intersect $\frac{1}{32}Q^\circ$ (because $Q \cap Q^\circ = \emptyset$ and $\delta_Q = \delta_{Q^\circ}=1$). Since $E \subset \frac{1}{32}Q^\circ$, we conclude that $\frac{65}{64}Q \cap E = \emptyset$. 
\item If $\delta_Q \geq 2$, then $ 0 \notin 2Q$ thanks to \eqref{qq6}. Assume for the sake of contradiction that $\frac{65}{64}Q \cap \frac{1}{32} Q^\circ \neq \emptyset$. Since $Q$ and $Q^\circ$ are disjoint, we conclude that $\frac{1}{64} \delta_Q \geq \frac{1}{4} \implies \delta_Q \geq 16$. Hence, $0 \in 2Q$ (since $\frac{65}{64}Q \cap \frac{1}{32} Q^\circ \neq \emptyset$ and $\delta_Q \geq 16$, and $Q^\circ = [0,1)^n$). Hence, we derive a contradiction. Thus, $\frac{65}{64}Q$ cannot intersect $\frac{1}{32}Q^\circ$. Since $E \subset \frac{1}{32}Q^\circ$, we conclude that $\frac{65}{64}Q \cap E = \emptyset$.
\end{itemize}
This completes the proof of property (c).

Property (d) is easy to prove. Let $Q \in \cQ$. Then $0 \in 2Q^+$ thanks to \eqref{qq7}. Hence, $0 \in 6Q$ (since $2Q^+ \subset 6Q$). Since $\delta_Q \geq 1$, this implies $Q^\circ \subset 9Q$ (recall that $Q^\circ = [0,1)^n$). Together with \eqref{smallish_cube}, this implies property (d).

We now prove property (f). Let $Q \in \overline{\CZ}(\cA^-)$ be given. 

If $Q \in \CZ(\cA^-)$ then $\#( 9Q \cap E ) \geq 2$, thanks to \eqref{Enearby2}.

If $Q \in \cQ$, then $9Q \supset Q^\circ$, hence $\#(9Q \cap E) = \#(E) \geq 2$.

This concludes the proof of property (f).

We prepare to describe the construction of the $\overline{\CZ}(\cA^-)$-\textsc{Oracle}.

We can determine whether a dyadic cube $Q \subset \R^n$ belongs to $\overline{\CZ}(\cA^-)$ using work and storage at most $C \log N$. We explain the procedure below.

Let $Q \subset \R^n$ be given.

First, suppose that $Q \subset Q^\circ$. Then $Q \in \overline{\CZ}(\cA^-)$ if and only if $Q \in \CZ(\cA^-)$. We can determine whether $Q \in \CZ(\cA^-)$ by using the $\CZ(\cA^-)$-\textsc{Oracle} to produce a list of all the cubes $Q' \in \CZ(\cA^-)$ satisfying $x_Q \in \frac{65}{64}Q'$. (Recall, $x_Q$ denotes the center of $Q$.) Then $Q \in \CZ(\cA^-)$ if and only if $Q$ belongs to the aforementioned list. Thus, in this case, we can determine whether $Q \in \overline{\CZ}(\cA^-)$ using work at most $C \log N$.

Next, suppose that $Q^\circ \subsetneq Q$. Then $Q$ can never belong to $\overline{\CZ}(\cA^-)$.

Lastly, suppose that $Q \subset \R^n \setminus Q^\circ$. Then $Q \in \overline{\CZ}(\cA^-)$ if and only if $Q \in \cQ$. Recall from \eqref{qq6} and \eqref{qq7} that $Q \in \cQ$ if and only if $[\delta_Q = 1 \; \mbox{or} \; 0 \notin 2Q]$ and $[\delta_{Q^+} \geq 2 \; \mbox{and} \; 0 \in 2 Q^+]$. We can check each of these conditions using at most $C$ computer operations. Thus, in this case we can determine whether $Q \in \overline{\CZ}(\cA^-)$ using work at most $C$.

Hence, we can determine whether a given cube belongs to $\overline{\CZ}(\cA^-)$ using work at most $C \log N$.

We next explain how to compute the unique cube $Q_{\underline{x}} \in \overline{\CZ}(\cA^-)$ containing $\underline{x}$. It will then not be difficult to produce a list of the cubes $Q \in \overline{\CZ}(\cA^-)$ satisfying $\underline{x} \in \frac{65}{64} Q$. We describe this step at the very end.

We check whether or not $\underline{x} \in Q^\circ$. We split into cases depending on the result.

First, suppose that $\underline{x} \in Q^\circ$. We then compute the cube $\overline{Q} \in \CZ(\cA^-)$ containing $\underline{x}$ using the $\CZ(\cA^-)$-\textsc{Oracle}. We set $Q_{\underline{x}} = \overline{Q}$.

Now suppose that $\underline{x} \in \R^n \setminus Q^\circ$.

Let $Q$ be the unique cube in $\cQ \setminus \{Q^\circ\}$ containing $\underline{x}$. We will explain how to compute $Q$.

We compute the dyadic cube $\widetilde{Q} \subset \R^n$ such that $\delta_{\widetilde{Q}}=1$ and $\underline{x} \in \widetilde{Q}$. 

We test to see whether $0 \in 2 \widetilde{Q}$. We can do that using at most $C$ computer operations.

If $0 \in 2 \widetilde{Q}$ then $\widetilde{Q}$ is a maximal dyadic cube satisfying the condition [$\delta_{\widetilde{Q}} \leq 1$ or $0 \notin 2 \widetilde{Q}$]. Hence, in that case, $\widetilde{Q}$ is the unique cube in $\cQ$ containing $\underline{x}$. We set $Q_{\underline{x}} = \widetilde{Q}$.

Now suppose that $0 \notin 2 \widetilde{Q}$. Thus, $\widetilde{Q}$ satisfies \eqref{qq6}. Since $Q$ and $\widetilde{Q}$ are intersecting dyadic cubes (they both contain $\underline{x}$), and since $Q$ is maximal with respect to the property \eqref{qq6}, we conclude that $\widetilde{Q} \subset Q$.

Assume that $\widetilde{Q} = Q$. Then $0 \notin 2 Q$, by assumption. On the other hand, suppose that $\widetilde{Q} \subsetneq Q$. Then $\delta_Q > 1$ (since $\delta_{\widetilde{Q}} = 1$). Since $Q$ satisfies \eqref{qq6}, we conclude that $0 \notin 2Q$. 

Thus, in the case where $0 \notin 2 \widetilde{Q}$, we know that $0 \notin 2Q$. Since $\underline{x} \in Q$ this shows that $\lvert \underline{x} \rvert \geq \frac{1}{4} \delta_Q$. Moreover, since $Q$ satisfies \eqref{qq7} we know that $0 \in 9Q$. Hence,
\begin{equation} \label{qq8} \frac{1}{4} \delta_Q \leq \lvert \underline{x}  \rvert   \leq 9 \delta_Q
\end{equation}
for the unique cube $Q \in \cQ$ containing $\underline{x}$. 

There are no more than $C$ dyadic cubes $Q \subset \R^n$ satisfying \eqref{qq8} with $\underline{x} \in Q$; moreover, it takes work at most $C$ to list all these cubes. We examine each cube and test to see whether it belongs to $\overline{\CZ}(\cA^-)$. We set aside the unique cube $Q$ that passes the test. We set $Q_{\underline{x}} = Q$.

We have just explained how to compute the cube $Q_{\underline{x}} \in \overline{\CZ}(\cA^-)$ containing a given point $\underline{x} \in \R^n$. The work requires is at most $C \log N$. We now explain how to construct the $\overline{\CZ}(\cA^-)$-\textsc{Oracle}. 

Suppose that $\overline{Q} \in \overline{\CZ}(\cA^-)$ satisfies $\underline{x} \in \frac{65}{64}\overline{Q}$. Then
\begin{equation}\label{gg8}
\overline{Q} \leftrightarrow Q_{\underline{x}} \; \mbox{and} \; \frac{1}{8} \delta_{Q_{\underline{x}}} \leq \delta_{\overline{Q}} \leq 8 \delta_{Q_{\underline{x}}}.
\end{equation}
This is a consequence of condition (b) in Proposition \ref{newcz_prop} and an application of Lemma \ref{gg_lem_1} (with $\gamma = 1/8$).

We produce a list of all the dyadic cubes $\overline{Q}$ that satisfy both \eqref{gg8} and $\underline{x} \in \frac{65}{64}\overline{Q}$. There are at most $C$ such cubes and it takes work at most $C$ to list them all. We examine each cube $Q$ to see whether it belongs to $\overline{\CZ}(\cA^-)$. We return the list of all those cubes that belong to $\overline{\CZ}(\cA^-)$.

This completes the description of the $\overline{\CZ}(\cA^-)$-\textsc{Oracle}. This completes the proof of the proposition.

\end{proof}

\subsection{Keystone Cubes}
\label{key_sec}

We define integer constants
\begin{equation} \label{consts}
\left\{
\begin{aligned}
&S_0 := S(\cA^-), \\
&S_1 := \mbox{the smallest integer greater than} \; 100, 10^5 \cdot S_0, \; \mbox{and} \;  2 \cdot \left[ c_*(\cA^-)\right]^{-1}, \\
&S_2 :=  \mbox{the smallest odd integer greater than} \; 10^{5}  S_1.
\end{aligned}
\right.
\end{equation}
We let $\epsilon > 0$ be a small parameter. We assume in what follows that
\begin{equation}
\label{smallepsassump}
\epsilon>0 \; \mbox{is less than a small enough universal constant}.
\end{equation}
We eventually fix $\epsilon$ to be a universal constant, but only much later in the proof. We will take $\epsilon_2(\cA) = \epsilon^\kappa$ and $\epsilon_1(\cA) = \epsilon^{1/\kappa}$ for a small universal constant $\kappa$. The discussion of the final choice of the numerical constants relevant to the Main Results for $\cA$ occurs in Section \ref{closing_remarks}. See also \eqref{fix_con}.

We next define the keystone cubes associated to the decomposition $\overline{\CZ}(\cA^-)$. We will prove a few basic properties of the keystone cubes and introduce the relevant algorithms.

\begin{defn}
A cube $Q^{\#} \in \overline{\CZ}(\cA^-)$ is keystone if and only if $\delta_{Q} \geq \delta_{Q^\#}$ for every $Q \in \overline{\CZ}(\cA^-)$ that meets $S_2 Q^{\#}$.

\end{defn}

\begin{lem}\label{key_geom}
The collection $\{ S_1 Q^\# : Q^\# \in \overline{\CZ}(\cA^-) \; \mbox{keystone} \}$ has bounded overlap. Moreover, each keystone cube $Q^\# \in \overline{\CZ}(\cA^-)$ belongs to $\CZ(\cA^-)$.
\end{lem}
\begin{proof}
Suppose that $Q_1^\#,Q^\#_2$ are keystone cubes such that $S_1 Q^\#_1 \cap S_1  Q_2^\# \neq \emptyset$ and $\delta_{Q^\#_1} \leq \delta_{Q^\#_2}$. Then $Q^\#_1 \cap S_2 Q^\#_2 \neq \emptyset$, since $S_2 \geq 10^5 S_1$. Therefore, $\delta_{Q^\#_1} \geq \delta_{Q^\#_2}$, by definition of the keystone cubes.

Consequently, $\delta_{Q_1^\#} = \delta_{Q^\#_2}$ whenever $S_1 Q^\#_1 \cap S_1 Q_2^\# \neq \emptyset$. Thus, no more than $C$ cubes $S_1Q_2^\#$ can intersect any given cube $S_1 Q^\#_1$. This implies the first conclusion of Lemma \ref{key_geom}.

Finally, observe that no cube in $\overline{\CZ}(\cA^-) \setminus \CZ(\cA^-)$ can be keystone, thanks to condition (d) in Proposition \ref{newcz_prop} and the fact that $S_2 \geq 100$. This completes the proof of the lemma.
\end{proof}

The definition of keystone cubes written above agrees with the definition in Section \ref{sec_ptkc}, where we set $K = S_2$ and let $A$ be a large universal constant in Section \ref{sec_ptkc}. The \textsc{Main Keystone Cube Algorithm} in Section \ref{sec_ptkc} says the following. Given $Q \in \CZ(\cA^-)$, we can compute a keystone cube $\mathcal{K}(Q) \in \CZ(\cA^-)$ such that the following condition holds.

There exists a sequence $\mathcal{S} = (Q_1,Q_2\cdots,Q_{\underline{L}})$ of $\overline{\CZ}(\cA^-)$ cubes such that
\[Q = Q_1 \leftrightarrow Q_2 \leftrightarrow \cdots \leftrightarrow Q_{\underline{L}} = \mathcal{K}(Q),\]
and such that
\[\delta_{Q_\ell} \leq C \cdot (1 - c)^{\ell - \ell'} \delta_{Q_{\ell'}} \quad \mbox{for} \; 1 \leq \ell' \leq \ell \leq \underline{L}.\]
We do not compute the sequence $\mathcal{S}$, we just claim its existence.

We now modify the sequence $\mathcal{S}$ to consist only of cubes from $\CZ(\cA^-)$ while maintaining the important properties of $\mathcal{S}$.

We first discuss the case in dimension $n = 1$. We let $\mathcal{S}'$ denote the sequence formed by omitting from $\mathcal{S}$ all the cubes that belong to $\overline{\CZ}(\cA^-) \setminus \CZ(\cA^-)$. Recall that all the cubes in $\CZ(\cA^-)$ are contained in $Q^\circ = [0,1)$ and all the cubes in $\overline{\CZ}(\cA^-)$ are contained in $\R \setminus [0,1)$. Consider a maximal subsequence $Q_{k_1},\cdots,Q_{k_2}$ of cubes in $\mathcal{S}$ that belong to $\overline{\CZ}(\cA^-)$. Then, by connectedness, each $Q_k$ ($k_1 \leq k \leq k_2$)  is contained in either $[1,\infty)$ or $(-\infty,0)$. Assume for sake of definiteness that each $Q_k$ is contained in $[1,\infty)$. Then $Q_{k_1-1}$ and $Q_{k_2+1}$ are the same cube in $\CZ(\cA^-)$, namely the unique cube in $\CZ(\cA^-)$ that meets the endpoint $x=1$. (This is because the sequence must exit and reenter $[0,1)$ using the same cube that borders the endpoint $x=1$.) Thus we can remove the aforementioned subsequence from $\mathcal{S}$ and obtain a connected path of cubes. The resulting sequence is exponentially decreasing with the same constants $C$ and $c$ above. The same argument shows that we can remove every maximal subsequence of $\mathcal{S}$ consisting of cubes in $\overline{\CZ}(\cA^-) \setminus \CZ(\cA^-)$.

We now handle the case when the dimension $n$ is at least $2$.

Suppose that some of the cubes in $\mathcal{S}$ belong to $\overline{\CZ}(\cA^-)$. Let $Q_{k_1}$ and $Q_{k_2}$ denote the first and last cubes in the sequence $\mathcal{S}$ belonging to $\overline{\CZ}(\cA^-) \setminus \CZ(\cA^-)$. Let $ \mathcal{S}_{\text{sub}} = (Q_{k_1},\cdots,Q_{k_2} )$ denote the corresponding subsequence of $\mathcal{S}$.

We know that $Q_1 = Q$ and $Q_{\underline{L}} = \mathcal{K}(Q)$ both belong to $\CZ(\cA^-)$. Hence, $1 < k_1 \leq k_2 < \underline{L}$.

Note that both $Q_{k_1-1}$ and $Q_{k_2-1}$ intersect the boundary of $Q^\circ$ and belong to $\CZ(\cA^-)$.


We join $Q_{k_1 - 1}$ and $Q_{k_2 +1}$ with a sequence $\mathcal{S}'_{\text{sub}} = (\widetilde{Q}_{k_1},\cdots,\widetilde{Q}_{k_3})$ with the following properties. 
\begin{itemize}
\item The cubes $\widetilde{Q}_k \in \CZ(\cA^-)$ intersect the boundary of $Q^\circ$. 
\item $\widetilde{Q}_{k_1} \leftrightarrow Q_{k_1 - 1}$, $\widetilde{Q}_{k_3} \leftrightarrow Q_{k_2 + 1}$, and
\[
\widetilde{Q}_k \leftrightarrow \widetilde{Q}_{k+1} \;\; \mbox{for} \; k_1 \leq k \leq k_3 - 1,
\]
\item $k_3 - k_1$ is bounded by a universal constant.
\item Each $\widetilde{Q}_k$ has sidelength between $1/2$ and $1/8$.
\end{itemize}
These properties can be arranged due to Lemma \ref{bdry_lem}.

We replace the subsequence $\mathcal{S}_{\text{sub}}$ with the sequence $\mathcal{S}_{\text{sub}}'$ in $\mathcal{S}$. We obtain a sequence $\mathcal{S}' = (\widetilde{Q}_1, \widetilde{Q}_2, \cdots, \widetilde{Q}_{L})$ of cubes in $\CZ(\cA^-)$ such that $\widetilde{Q}_1 = Q$ and $\widetilde{Q}_{L} = \mathcal{K}(Q)$; moreover, 
\[
\widetilde{Q}_{\ell} \leftrightarrow \widetilde{Q}_{\ell +1} \;\; (1 \leq \ell \leq L - 1) \; \mbox{ and } \; \delta_{\widetilde{Q}_{\ell}} \leq C' \cdot (1-c')^{\ell - \ell'} \delta_{\widetilde{Q}_{\ell'}} \;\; (1 \leq \ell' \leq \ell \leq L).
\]
Indeed, the fact that $\mathcal{S}'$ satisfies the exponentially decreasing property follows directly from the construction: We removed a subsequence of connected cubes in $\mathcal{S}$ and replaced it with a subsequence of bounded length consisting of cubes of size $\in \{1/2,1/4,1/8\}$. This has no effect on the fact that the sidelengths are exponentially decreasing in the sense of the above estimate. 

Hence, the sequence $\mathcal{S}'$ joining $Q$ and $\mathcal{K}(Q)$ is exponentially decreasing.

We never actually compute the sequences $\mathcal{S}$ or $\mathcal{S}'$, we just claim their existence.

Using the above analysis, the \textsc{Main Keystone Cube Algorithm} and the algorithm \textsc{List All Keystone Cubes} in Section \ref{sec_ptkc}, we obtain the following result.

\environmentA{Algorithm: Keystone-Oracle.} 

After one-time work at most $C N \log N$ in space $CN$ we produce the following outcomes:

\begin{itemize}
\item We list all the keystone cubes $Q^\#$ in $\CZ(\cA^-)$.
\item We can answer queries: A query consists of a cube $Q \in \CZ(\cA^-)$, and the response to a query $Q$ is a keystone cube $\mathcal{K}(Q)$ to which $Q$ is connected by an exponentially decreasing path
$$Q = \widetilde{Q}_1 \leftrightarrow \widetilde{Q}_2 \leftrightarrow \cdots \leftrightarrow \widetilde{Q}_L = \mathcal{K}(Q)$$
with
$$\delta_{\widetilde{Q}_\ell} \leq C \cdot (1-c)^{\ell-\ell'} \delta_{\widetilde{Q}_{\ell'}} \;\; \mbox{for} \; 1 \leq \ell' \leq \ell \leq L.$$
We guarantee that $\widetilde{Q}_\ell \in \CZ(\cA^-)$ and $\frac{65}{64} \widetilde{Q}_\ell \subset C Q$ for each $\ell$. We guarantee that $S_1 \mathcal{K}(Q) \subset CQ$; also that $\mathcal{K}(Q) = Q$ if $Q$ is keystone. The work required to answer a query is at most $C \log N$.
\item We list all $(Q',Q'') \in \CZ(\cA^-) \times \CZ(\cA^-)$ such that $Q' \leftrightarrow Q''$ and $\mathcal{K}(Q') \neq \mathcal{K}(Q'')$. Let $\BD(\cA^-)$ (the ``border disputes'') denote the set of all such pairs $(Q',Q'')$. We guarantee that the cardinality of $\BD(\cA^-)$ is at most $CN$.

\end{itemize}

\begin{remk} \label{key_rem0}
Let $\widetilde{Q}_1 \leftrightarrow \cdots \leftrightarrow \widetilde{Q}_L$ be as above. For fixed $Q'$, we can have $Q' = \widetilde{Q}_\ell$ for at most $C$ distinct $\ell$. This is because the path $\widetilde{Q}_1 \leftrightarrow \widetilde{Q}_2 \leftrightarrow \cdots \leftrightarrow \widetilde{Q}_L$ is exponentially decreasing.
\end{remk}

\begin{remk} \label{key_remk1}
We do not attempt to compute the sequence of cubes $\widetilde{Q}_1,\cdots,\widetilde{Q}_{L-1}$ - we only guarantee that this sequence exists, and we guarantee that we can compute the keystone cube $\mathcal{K}(Q) = \widetilde{Q}_L$ located at the end of the sequence.
\end{remk}

\section{An Approximation to the Sigma} \label{sec_test}

We begin the proof of the Main Technical Results for $\cA$. We recall that $\cA \subsetneq \cM$ is a monotonic set. 

In Sections \ref{sec_ocz} and \ref{key_sec}, we have defined a dyadic decomposition $\overline{\CZ}(\cA^-)$ of $\R^n$ and a notion of keystone cubes in $\overline{\CZ}(\cA^-)$. We cannot compute all the cubes in $\overline{\CZ}(\cA^-)$ since there are infinitely many. Instead, we have access to a $\overline{\CZ}(\cA^-)$-\textsc{Oracle} and the \textsc{Keystone-Oracle}.

The integer constants $S_0,S_1.S_2$ relating to the keystone cubes are defined in \eqref{consts}.

According to the Main Technical Results for $\cA^-$ (see Chapter \ref{sec_mainresults}), for each $Q \in \CZ_{\main}(\cA^-)$ the functional
\begin{equation}
\label{mdefn} M_{(Q,\cA^-)}(f,R) := \left(\sum_{\xi \in \Xi(Q,\cA^-)} \lvert \xi(f,R) \rvert^p \right)^{1/p}
\end{equation}
satisfies
\begin{equation}
\label{n_appx}
c \| (f,R) \|_{(1+ a) Q} \leq M_{(Q,\cA^-)}(f,R) \leq C \| (f,R)\|_{\frac{65}{64}Q}.
\end{equation}
Recall that $a$ is the constant $a(\cA^-)$ from the Main Technical Results; see \eqref{a_defn}.

For each $Q \in \overline{\CZ}(\cA^-)\setminus \CZ_{\main}(\cA^-)$, we define $\Xi(Q,\cA^-) := \emptyset$ and $M_{(Q,\cA^-)}(f,R) := 0$. By definition of the collection $\CZ_{\main}(\cA^-)$ and by property (c) in Proposition \ref{newcz_prop}, we have $\frac{65}{64}Q \cap E = \emptyset$. Thus, $\| (f,R)\|_{(1+a)Q} = 0$ for any $Q \in \overline{\CZ}(\cA^-)\setminus \CZ_{\main}(\cA^-)$. Thus, we see that \eqref{n_appx} holds for all $Q \in \overline{\CZ}(\cA^-)$.

\subsection{Assigning Jets to Keystone Cubes} \label{sec_assign}

Let $Q^\# \in \CZ(\cA^-)$ be a keystone cube. We define its associated $\CZ$ cubes to be the collection
\begin{equation}
\label{Idefn}
\cI(Q^\#) := \bigl\{ Q \in \overline{\CZ}(\cA^-): Q \cap S_0 Q^\# \neq \emptyset \bigr\}.
\end{equation}
We note that the cubes in $\cI(Q^\#)$ belong to $\overline{\CZ}(\cA^-)$ rather than $\CZ(\cA^-)$. Hence, the cubes in $\cI(Q^\#)$ are contained in $\R^n$, and may not be contained in $Q^\circ$.

\begin{lem}
\label{touch_lem}
Let $A \geq 1$ be given. Assume that $Q, \oQ \in \overline{\CZ}(\cA^-)$ and $Q \cap A  \oQ \neq \emptyset$. Then
\begin{align} \label{touch1} &\delta_Q \leq 10^3 A \delta_{\oQ}, \; \mbox{and} \\
\label{touch2} & \frac{65}{64}Q \subset 10^5 A  \oQ.
\end{align}
\end{lem}
\begin{proof}
We first prove \eqref{touch1}. Assume for the sake of contradiction that $\delta_Q > 10^3 A  \delta_{\oQ}$ for some $Q,\oQ \in \overline{\CZ}(\cA^-)$ with $Q \cap A  \oQ \neq \emptyset$. Then $\frac{65}{64}Q \cap \oQ \neq \emptyset$. However, this contradicts the good geometry of the cubes in $\overline{\CZ}(\cA^-)$ (see Proposition \ref{newcz_prop}). This completes the proof of \eqref{touch1} by contradiction. Lastly, \eqref{touch2} follows from \eqref{touch1} and our assumption that $Q \cap A  \oQ \neq \emptyset$.
\end{proof}

By Lemma \ref{touch_lem}, the $\CZ$ cubes associated to a given $Q^\#$ satisfy the following property: for each $Q \in \cI(Q^\#)$ we have
\begin{align} 
\label{touch1_a} &\delta_Q \leq 10^3 S_0 \cdot \delta_{Q^\#}, \; \mbox{and} \\
\label{touch2_a} & \frac{65}{64}Q \subset S_1 Q^\#.
\end{align}
(Recall \eqref{consts} which states that $S_1 \geq 10^5 S_0$.)

\begin{remk} \label{key_rem1}
The definition of keystone cubes shows that $\delta_Q \geq \delta_{Q^\#}$ whenever $Q \in \cI(Q^\#)$. Hence, \eqref{touch1_a} implies that the cardinality of $\cI(Q^\#)$ is bounded by $C$ for each keystone cube $Q^\#$.

If $Q \in \cI(Q_1^\#)$ and $Q \in \cI(Q_2^\#)$ then \eqref{touch2_a} implies that $S_1 Q^\#_1 \cap S_1 Q^\#_2 \supset \frac{65}{64}Q$. Recall that the cubes $S_1 Q^\#$ ($Q^\#$ keystone) have bounded overlap. (See Lemma \ref{key_geom}.) Thus, each $Q \in \overline{\CZ}(\cA^-)$ belongs to $\cI(Q^\#)$ for at most $C$ distinct keystone cubes $Q^\#$.
\end{remk}

\environmentA{Algorithm: Make New Assists and Assign Keystone Jets.}

For each keystone cube $Q^\#$, we compute a list of new assists $\Omega^{\new}(Q^\#) \subset \left[ \X(S_1 Q^\# \cap E) \right]^*$, written in short form, and we produce an $\Omega^\new(Q^\#)$-assisted bounded depth linear map $R^\#_{Q^\#} : \X(S_1 Q^\# \cap E) \oplus \cP \rightarrow \cP$, written in short form. Furthermore, we guarantee that the following conditions are met.
\begin{itemize}
\item The sum of $\depth(\omega)$ over all $\omega$ in $\Omega^\new(Q^\#)$, and over all keystone cubes $Q^\#$, is bounded by $C N$.
\end{itemize}
Given $(f,P) \in \X(S_1 Q^\# \cap E) \oplus \cP$ , denote $R^\# = R^\#_{Q^\#}(f,P)$. 

\begin{itemize}
\item Then $\partial^\alpha \left( R^\# - P \right) \equiv 0$ for all $\alpha \in \cA$ (recall, $\cA$ is monotonic; see Remark \ref{mon_rem}).
\item Let $R \in \cP$, with $ \partial^\beta \left(R - P \right) \equiv 0$ for all $\beta\in \cA$. Then 
\begin{equation}
\label{eqstuff}
\sum_{Q \in \cI(Q^\#)} \sum_{\xi \in \Xi(Q,\cA^-)} \lvert \xi(f,R^\#) \rvert^p \leq C \sum_{Q \in \cI(Q^\#)} \sum_{\xi \in \Xi(Q,\cA^-)} \lvert \xi(f,R) \rvert^p.
\end{equation}
\end{itemize}
To compute the assists $\Omega^\new(Q^\#)$ and the short form of the maps $R^\#_{Q^\#}$ (for all the keystone cubes $Q^\#$) requires work at most $C N \log N$, and storage at most $C N$.

\begin{proof}[\underline{Explanation}]

Given $P \in \cP$, we define $V_P$ to be the affine subspace consisting of all polynomials $R \in \cP$ satisfying $\partial^\alpha ( R - P) \equiv 0$ for all $\alpha \in \cA$.  We note that $R \in V_P$ $\iff$ $\partial^\alpha (R-P)(0) = 0$ for all $\alpha \in \cA$, since $\cA$ is monotonic.

We introduce coordinates on $V_P$, defined by
$$w = (w_1,\cdots,w_k) \in \R^k \implies R_{w}(x) = \sum_{\alpha \in \cA} \frac{\partial^\alpha P(0)}{\alpha!} x^\alpha + \sum_{j=1}^k w_j \cdot \frac{x^{\alpha_j}}{\alpha_j!},$$
where we write $\cM \setminus \cA = \{ \alpha_1,\cdots,\alpha_k \}$. 

We consider the sum of the $p$-th powers of the functionals $\xi(f,R_w)$ over all $\xi \in \Xi(Q,\cA^-)$ and $Q \in \cI(Q^\#)$. We want to minimize this expression with respect to $w \in \R^k$. We can approximately solve this minimization problem using the algorithm \textsc{Optimize via Matrix} from Section \ref{sec_lf}. We describe the process below.

For each $Q \in \CZ_{\main}(\cA^-)$ with $Q \cap S_0 Q^\# \neq \emptyset$, we have $\delta_{Q^\#} \leq \delta_Q$ by definition of keystone cubes. Hence, from \eqref{touch1_a} we have
\begin{align}
\label{influence}
Q \cap S_0 Q^\# \neq \emptyset \;  \mbox{and} &\; \delta_{Q^\#} \leq \delta_Q \leq C \cdot \delta_{Q^\#} \\
\notag{}
& \mbox{for a universal constant} \; C.
\end{align}
We list all the dyadic cubes $Q$ that satisfy \eqref{influence}. There are at most $C$ cubes in this list. We test each $Q$ to see whether it belongs to $\CZ_{\main}(\cA^-)$. Thus, we can compute the list
\[
\mathfrak{L} = \bigl\{ Q \in \CZ_{\main}(\cA^-) : \; Q \cap S_0 Q^\# \neq \emptyset \bigr\}.
\] 
The list $\mathfrak{L}$ contains all the cubes $Q$ that participate in \eqref{eqstuff} for which $\Xi(Q,\cA^-) \neq \emptyset$. \\
(Recall that $\Xi(Q,\cA^-) = \emptyset$ for $Q \in  \overline{\CZ}(\cA^-) \setminus \CZ_{\main}(\cA^-)$.)

We list all the functionals $\xi$ appearing in $\Xi(Q,\cA^-)$ for some $Q \in \mathfrak{L}$. From the Main Technical Results for $\cA^-$ (see Chapter \ref{sec_mainresults}), each such $\xi$ is given in the form
\[
\xi(f,R_w) = \lambda(f) + \sum_{a=1}^I \mu_a \cdot \omega_a(f)  + \check{\lambda}((\partial^\alpha P(0))_{\alpha \in \cA}) + \sum_{j=1}^k \check{\mu}_j \cdot w_j,
\]
where $\lambda$ and $\check{\lambda}$ are linear functionals; $\omega_a \in \Omega(Q,\cA^-)$ for some $Q \in \mathfrak{L}$;  $\mu_a$ and $\check{\mu}_j$ are real coefficients; and $\depth(\lambda) = \mathcal{O}(1)$, $I = \mathcal{O}(1)$. In this discussion, we write $X = \mathcal{O}(Y)$ to indicate that $X \leq C Y$ for a universal constant $C$.

Processing each functional $\xi$ this way takes work $\mathcal{O}(1)$ per functional. Thus, with work $\mathcal{O} (L)$ (see below),
we obtain a list of all the above $\xi$'s, written as
\begin{align}
\label{xi_short}
\xi_\ell(f,R_w) &= \lambda_\ell(f) + \sum_{a=1}^{I_\ell} \mu_{\ell a} \omega_{\ell a}(f) + \check{\lambda}_\ell((\partial^\alpha P(0))_{\alpha \in \cA}) + \sum_{j=1}^k \check{\mu}_{\ell j} w_j \\
& \qquad\qquad\qquad \mbox{for} \; \ell =1,\cdots,L; \; \mbox{here,} \; L = \sum_{Q \in \cI(Q^\#)} \# \bigl[ \Xi(Q,\cA^-) \bigr]. \notag{}
\end{align}
Here, each $I_\ell = \mathcal{O}(1)$, each $\lambda_\ell$ has bounded depth, and each $\omega_{\ell a}$ belongs to $\Omega(Q_{\ell a},\cA^- )$ for some $Q_{\ell a} \in \mathfrak{L}$. Of course the $Q_{\ell a}$ need not be distinct, and $k \leq \dim(\cP) = D$.

Processing the functionals $w \mapsto \xi_\ell(f,R_w)$ in \eqref{xi_short} with the algorithm \textsc{Optimize via Matrix} (see Section \ref{sec_lf}), we compute a matrix $(b_{j\ell})$ with the following property. The sum of the absolute values of the $p$-th powers of the functionals $\xi_\ell(f,R_w)$ ($1 \leq \ell \leq L$) is essentially minimized for fixed $f$, $(\partial^\alpha P(0))_{\alpha \in \cA}$ by setting $w = w^* = (w^*_1,\cdots,w^*_k)$, where
\begin{align}
\label{wc11}
w_j^* &= \sum_{\ell=1}^L b_{j\ell} \left[ \lambda_\ell(f) + \sum_{a=1}^{I_\ell} \mu_{\ell a} \omega_{\ell a}(f) + \check{\lambda}_\ell((\partial^\alpha P(0))_{\alpha \in \cA}) \right] \\
& = \left\{  \sum_{\ell=1}^L b_{j\ell} \left[ \lambda_\ell(f) + \sum_{a=1}^{I_\ell} \mu_{\ell a} \omega_{\ell a}(f) \right]\right\} + \sum_{\ell=1}^L b_{j\ell} \cdot \check{\lambda}_\ell((\partial^\alpha P(0))_{\alpha \in \cA}) \notag{} \\
& \equiv  \qquad\qquad \omega^{\new}_j(f) \qquad\qquad\qquad\quad\; + \quad \lambda^{\new}_j((\partial^\alpha P(0))_{\alpha \in \cA}). \notag{}
\end{align}
We have thus defined new assists $\omega_j^\new$ and new functionals $\lambda_j^\new$.

We may therefore take $R_{Q^\#}^\#(f,P) := R_{w^*}$ with $w^*_j = \omega^{\new}_j(f) + \lambda^{\new}_j((\partial^\alpha P(0))_{\alpha \in \cA})$ ($1 \leq j \leq k$) and we obtain the estimate \eqref{eqstuff}. Note that $R_{Q^\#}^\#$ has assisted bounded depth, with assists $\omega_j^{\new}$ ($j=1,\cdots,k$). Indeed, 
\begin{equation}
\label{derivativesform}
\partial^\alpha \left[ R^\#_{Q^\#}(f,P) \right](0) = \left\{
\begin{array}{ccc}
&\omega^{\new}_j(f) + \lambda^{\new}_j((\partial^\alpha P(0))_{\alpha \in \cA}) &\mbox{if} \; \alpha = \alpha_j, \; j \in \{1,\cdots,k\} \\
& \partial^\alpha P(0) & \mbox{if} \; \alpha \in \cA.
\end{array}
\right.
\end{equation}

We can compute the new functionals $\lambda^\new_j$ ($1 \leq j \leq k$) using the obvious method. This requires work 
\[ 
\mathcal{O}(L) = \mathcal{O}\biggl( \sum_{Q \in \cI(Q^\#)} \# \bigl[ \Xi(Q,\cA^-) \bigr]\biggr).
\]

We will now express the new assists $\omega^{\new}_j$ in short form. 

We write $\omega_j^\new = \omega_j^{\new,1} + \omega_j^{\new,2}$, where $\omega_j^{\new,1}$ and $\omega_j^{\new,2}$ are defined below (see \eqref{first_func} and \eqref{last_fnc}).

Each $\lambda_\ell(f)$ has bounded depth, so the functional
\begin{equation}
\label{first_func}
\omega_j^{\new,1} : f \mapsto \sum_{\ell=1}^L b_{j\ell}   \cdot   \lambda_\ell(f)
\end{equation}
can be computed in short form using \\
work $\displaystyle \mathcal{O}(L \log L) = \mathcal{O}\biggl( \log N \cdot \sum_{Q \in \cI(Q^\#)} \# \left[ \Xi(Q,\cA^-) \right]\biggr)$ and storage $\displaystyle \mathcal{O}(L) = \mathcal{O}\biggl( \sum_{Q \in \cI(Q^\#)} \# \left[ \Xi(Q,\cA^-) \right] \biggr)$. 

This computation follows by a simple sorting procedure. We provide details below.
\begin{itemize}
\item Recall that $\lambda_\ell$ has bounded depth and is given in short form (without assists):
\begin{equation}
\label{cl}
\lambda_\ell(f) = \sum_{x \in S_\ell} c_\ell(x) \cdot f(x), \;\; \mbox{where} \; \# (S_\ell) \leq C.
\end{equation}
Thus, we can express the functional \eqref{first_func} as
\begin{align}
\label{new_form1}
\omega_j^{\new,1} : f \mapsto \sum_{x \in S} d_j(x) \cdot f(x), \;\; \mbox{where} \; &S = \bigcup_{\ell=1}^L S_\ell \; \mbox{and} \\
\notag{}
& d_j(x) = \sum_{\ell=1}^L b_{j \ell} \cdot c_\ell(x) \quad \mbox{for} \; x \in S.
\end{align}
We compute the weights $d_j(x)$ by sorting. More precisely, we sort the points of $S$. We make an array indexed by $S$. We initialize the array to have all zero entries. We loop over $\ell = 1,\cdots,L$, and we loop over all the points $y \in S_\ell$. We determine the position of each $y$ in the list $S$, and we add the number $b_{j \ell} \cdot c_\ell(y)$ at the relevant position in the array. This requires work at most $C \log (S) \leq C\log L$ for a fixed pair $(\ell, y)$. Hence, the total work required is at most $C L \log L$, since
the number of relevant pairs $(\ell,y)$ is $\sum_{\ell=1}^L \#(S_\ell) \leq \sum_{\ell=1}^L C \leq C L$. Similarly, the total storage is at most $ C \sum_{\ell=1}^L \#(S_\ell) \leq CL$. 
\end{itemize}
Thus, we can compute the functional \eqref{first_func} using work $\cO(L \log L)$ and storage $\cO(L)$.

It remains to compute the functional
\begin{equation}
\label{last_fnc}
\omega_j^{\new,2} : f \mapsto \sum_{\ell=1}^L b_{j\ell} \sum_{a=1}^{I_\ell} \mu_{\ell a} \omega_{\ell a}(f) \quad \mbox{in short form.} \; (\mbox{Recall, each} \; I_\ell = \mathcal{O}(1).)
\end{equation}
We recall that each $\omega_{\ell a}$ belongs to $\displaystyle  \bigcup_{Q \in \cI(Q^\#)} \Omega(Q,\cA^-)$.

We can express the functional \eqref{last_fnc} in the form
\begin{equation} \label{new_form2}
\omega_j^{\new,2} : f \mapsto \sum_{\omega \in \bigcup_{Q \in \cI(Q^\#)} \Omega(Q,\cA^-)} q_{j\omega} \cdot \omega(f),
\end{equation}
with\\
work $\displaystyle \mathcal{O}(L \log L) = \mathcal{O}\biggl( \log N \cdot \sum_{Q \in \cI(Q^\#)} \# \left[ \Xi(Q,\cA^-) \right]\biggr)$ and storage $\displaystyle \mathcal{O}(L) = \mathcal{O}\biggl( \sum_{Q \in \cI(Q^\#)} \# \left[ \Xi(Q,\cA^-) \right] \biggr)$. We can compute the relevant numbers $q_{j\omega}$ by sorting, since 
\begin{equation}
\label{q_defn}
q_{j \omega} = \sum_{(\ell,a) : \omega_{\ell a} = \omega} b_{j \ell} \cdot \mu_{\ell a} .
\end{equation}

Finally, once our functional is in the form \eqref{new_form2}, we can easily write it in short form
\begin{equation}
\label{new_form3}
\omega_j^{\new,2} : f \mapsto \sum_{x \in S} k_j(x) \cdot f(x)
\end{equation}
with 
\[\mbox{work} \; \mathcal{O}\left( \log N \cdot \sum_{Q \in \cI(Q^\#)} \sum_{\omega \in \Omega(Q,\cA^-)} \depth(\omega) \right) \; \mbox{and storage} \; \mathcal{O}\left( \sum_{Q \in \cI(Q^\#)} \sum_{\omega \in \Omega(Q,\cA^-) )} \depth(\omega)\right).\]
Again, we perform a sort to carry this out.

We compute $\omega_j^\new = \omega_j^{\new,1} + \omega_j^{\new,2}$ in short form by adding the short form expressions \eqref{new_form1} and \eqref{new_form3}.

Altogether, we obtain the new assists $\omega_j^\new$ and the new functionals $\lambda_j^\new$ for a given $Q^\#$ using work at most
\[C \cdot (\log N) \cdot \left[ \sum_{Q \in \cI(Q^\#)} \left\{ 1 + \# \bigl[ \Xi(Q,\cA^-) \bigr] + \sum_{\omega \in \Omega(Q,\cA^-)} \depth(\omega) \right\} \right]\]
and storage at most
\[C \cdot \left[ \sum_{Q \in \cI(Q^\#)} \left\{ 1 + \# \bigl[ \Xi(Q,\cA^-) \bigr] + \sum_{\omega \in \Omega(Q,\cA^-)} \depth(\omega) \right\} \right].\]
(Again, recall that $\Xi(Q,\cA^-)  = \Omega(Q,\cA^-) = \emptyset$ for any $Q \in  \overline{\CZ}(\cA^-) \setminus \CZ_{\main}(\cA^-)$.)

Each $Q \in \CZ_{\main}(\cA^-)$ belongs to $\cI(Q^\#)$ for at most $C$ distinct $Q^\#$ (see Remark \ref{key_rem1}). Therefore, we can compute the new assists and the new functionals for all the keystone cubes $Q^\#$ using work at most
\begin{align*}
C \cdot (\log N) \cdot \biggl[ & \# \{ \; \text{Keystone Cubes} \; Q^\# \} + \sum_{Q \in \CZ_{\main}(\cA^-)}  \# \bigl[ \Xi(Q,\cA^-) \bigr] \\
& + \sum_{Q \in \CZ_{\main}(\cA^-)}  \sum_{\omega \in \Omega(Q,\cA^-)} \depth(\omega) \biggr],
\end{align*}
which is at most $C N \log N$ by the induction hypothesis and the statement of the \textsc{Keystone-Oracle} (which guarantees that the number of keystone cubes is bounded by $C N$). Similarly, we see that all the new assists can be computed using storage at most $C N$.

Finally, note that there are at most $D$ new assists for each given keystone cube $Q^\# \in \CZ(\cA^-)$, and each such assist has depth at most $\#(S_1 Q^\# \cap E)$. By the bounded overlap of the cubes $S_1 Q^\#$ (see Lemma \ref{key_geom}), we see that the sum of the depths of all the new assists is at most $C \cdot \#(E) = C N$.

This completes the explanation of the algorithm.

\end{proof}

Let $Q^\# \in \CZ(\cA^-)$ be a keystone cube. For each $(f,R) \in \X(S_1 Q^\# \cap E) \oplus \cP$, we define
\begin{align} \label{sharp_norm} \left[ M^\#_{Q^\#}(f,R) \right]^p &:= \sum_{\substack{Q \in \cI(Q^\#)}} \sum_{\xi \in \Xi(Q,\cA^-)} \lvert \xi(f,R) \rvert^p  = \sum_{\substack{ Q \in \cI(Q^\#)}}  \left[ M_{(Q,\cA^-)}(f,R)\right]^p.
\end{align}
The terms $ \xi(f,R)$ appearing above are well-defined, since \eqref{touch2_a} states that $\frac{65}{64}Q \subset S_1 Q^\#$ for each $Q \in \cI(Q^\#)$. (Recall that the domain of each functional $\xi$ in $\Xi(Q,\cA^-)$ is $\X((65/64)Q \cap E) \oplus \cP$.)

We now show that the ``keystone functional'' defined in \eqref{sharp_norm} is comparable to the trace semi-norm near the given keystone cube.

\begin{lem}\label{key_lem1} Let $Q^\#$ be a keystone cube. Then 
\[ c \cdot \|(f, R)\|_{S_0 Q^\#} \leq M^\#_{Q^\#}(f,R) \leq C \cdot  \|(f,R)\|_{S_1 Q^\#}\]
for all $(f,R) \in \X(S_1 Q^\# \cap E) \oplus \cP$.
\end{lem}

\begin{proof}
From \eqref{touch2_a} we learn that $(1+a)Q \subset (65/64)Q \subset S_1 Q^\#$ for any $Q \in \cI(Q^\#)$. (Recall that $a = a(\cA^-) \leq 1/64$; see \eqref{a_defn}.) 

Let $(f,R) \in \X(S_1 Q^\# \cap E) \oplus \cP$ be given.

For each $Q \in \cI(Q^\#)$, by definition of the seminorm $\| (\cdot,\cdot)\|_{(1+a)Q}$, we can choose $F_Q \in \X$ with $F_Q = f$ on $E \cap (1+a)Q$ and
\[ \|F_Q\|_{\X((1+a)Q)} + \delta_{(1+a)Q}^{-m} \|F_Q - R \|_{L^p((1+a)Q)} \leq 2 \cdot \| (f,R) \|_{(1+a)Q}.\]
From the left-hand estimate in \eqref{n_appx} we have $ \| (f,R) \|_{(1+a)Q} \leq C \cdot M_{(Q,\cA^-)}(f,R)$. Thus, by definition \eqref{sharp_norm} we have
\begin{equation} \label{yuk}\|F_Q\|_{\X((1+a)Q)} + \delta_Q^{-m} \|F_Q - R \|_{L^p((1+a)Q)} \leq C \cdot M_{Q^\#}^\#(f,R) \quad \mbox{for} \; Q \in \cI(Q^\#). 
\end{equation}

The assumptions in Sections \ref{czalg_sec} and \ref{sec_pou} are valid, where
\begin{itemize}
\item $\CZ = \overline{\CZ}(\cA^-)$ and $\cQ = \cI(Q^\#)$.
\item $\hQ = S_0 Q^\#$;
\item $\overline{r} = a$, and $A = C$ for a large enough universal constant $C$.
\end{itemize}
We exhibited a $\overline{\CZ}(\cA^-)$-\textsc{Oracle} and proved good geometry for $\overline{\CZ}(\cA^-)$, which is a decomposition of $\R^n$, in Proposition \ref{newcz_prop}. We proved the properties of the collection $\cQ$ stated in Section \ref{sec_pou}: by definition of $\cI(Q^\#)$ in \eqref{Idefn} and since $\overline{\CZ}(\cA^-)$ is a partition of $\R^n$, we obtain \eqref{covers}; also, from \eqref{touch1_a} we obtain \eqref{sizebd}. 

We may thus apply the results stated in Section \ref{sec_pou}.

By Lemma \ref{pou_lem}, there exists a partition of unity $\theta_Q^{Q^\#} \in C^m(\R^n)$ ($Q \in \cI(Q^\#)$) such that
\[ \sum_{Q \in \cI(Q^\#)} \theta^{Q^\#}_Q = 1 \;\;  \mbox{on} \;\; S_0 Q^\#,\]
where $\supp \theta^{Q^\#}_Q \subset (1+a)Q$ and $ \lvert \partial^\alpha \theta^{Q^\#}_Q(x)  \rvert \leq C \cdot \delta_Q^{-|\alpha|}$ for $x \in (1+a)Q$, $|\alpha| \leq m$. 

Define
\[F := \sum_{Q \in \cI(Q^\#)} F_Q \cdot \theta^{Q^\#}_Q.\]
Since the cardinality of $\cI(Q^\#)$ is at most $C$ (see Remark \ref{key_rem1}), Lemma \ref{patch_lem} and \eqref{yuk} show that
\begin{equation*}
\|F\|_{\X(S_0Q^\#)} \leq C \cdot M_{Q^\#}^\#(f,R).
\end{equation*}
Moreover, since $\delta_{Q^\#} \leq \delta_Q \leq C \delta_{Q^\#}$ for all $Q \in \cI(Q^\#)$, we have
\begin{align*} \delta_{Q^\#}^{-m} \| F - R \|_{L^p(S_0 Q^\#)} &\leq C  \sum_{Q \in \cI(Q^\#)}  \delta_Q^{-m} \| F_Q - R \|_{L^{p}((1+a)Q)} \| \theta_Q^{Q^\#} \|_{L^\infty((1+a)Q)}  \\ 
& \leq C \cdot M_{Q^\#}^\#(f,R) \qquad \qquad \qquad \qquad (\mbox{see \eqref{yuk}}).
\end{align*}

Because $\supp( \theta_Q^{Q^\#}) \subset (1+a)Q$ and $F_Q = f$ on $E \cap (1+a)Q$ we see that $F = f $ on $E \cap S_0 Q^\#$. Hence, the above estimates imply that
$$\|(f,R)\|_{S_0 Q^\#} \leq \|F\|_{\X(S_0Q^\#)} + \delta_{S_0Q^\#}^{-m} \| F - R \|_{L^p(S_0 Q^\#)}  \lesssim M_{Q^\#}^\#(f,R).$$
This proves one inequality in the statement of the lemma.

Next, using the right-hand estimate in \eqref{n_appx}, we see that
\begin{equation*}
\left[ M_{Q^\#}^\#(f,R) \right]^p = \sum_{\substack { Q \in \cI(Q^\#)  \\ Q \in \CZ_{\main}(\cA^-) } } \left[ M_{(Q,\cA^-)}(f,R) \right]^p \leq C \sum_{Q \in \cI(Q^\#)} \|(f,R)\|_{\frac{65}{64}Q}^p.
\end{equation*}
Recall that $\cI(Q^\#)$ contains at most $C$ cubes. Thus, by Lemma \ref{lem_normmon}, where we use the estimate $\delta_{Q^\#} \leq \delta_Q \leq C \delta_{Q^\#}$ and that $\frac{65}{64}Q \subset S_1 Q^\#$ for each $Q \in \cI(Q^\#)$, the right-hand side in the above estimate is bounded by $C \cdot \|(f,R)\|_{S_1 Q^\#}^p$. This completes the proof of the lemma.
\end{proof}

The parameter $\epsilon>0$ now makes its first appearance. Recall that $\epsilon$ is assumed to be less than a small enough universal constant. See \eqref{smallepsassump}.

\begin{prop}\label{key_prop1}
Let $\hQ$ be a dyadic subcube of $Q^\circ$, such that $3\hQ$ is tagged with $(\cA,\epsilon)$.

Assume also that $Q^\# \in \CZ(\cA^-)$ is a keystone cube, and that $S_1 Q^\# \subseteq \frac{65}{64}\hQ$.

Suppose that $H \in \X$ satisfies $H = f$ on $E \cap S_1Q^\#$ and $\partial^\alpha H(x_{Q^\#}) = \partial^\alpha P(x_{Q^\#})$ for all $\alpha \in \cA$. Then
\begin{equation} \label{imp} 
\delta_{Q^\#}^{-m} \| H - R_{Q^\#}^\#  \|_{L^p(S_1 Q^\#)} \leq C \cdot \| H \|_{\X(S_1 Q^\#)}.
\end{equation}
Here, $C \geq 1$ is a universal constant; and $R^\#_{Q^\#} = R^\#_{Q^\#}(f,P)$. \\
(See the algorithm \textsc{Make New Assists and Assign Keystone Jets}.)
\end{prop}
\begin{proof}

Recall that $S_0 Q^\# \subset \frac{65}{64}\hQ$ and $3 \hQ$ is tagged with $(\cA,\epsilon)$. Thus, Lemma \ref{pre_lem5} implies that $S_0 Q^\#$ is tagged with $(\cA, \epsilon^\kappa)$ for some universal constant $\kappa > 0$.

Recall that $S_1 \geq 2 [c_*(\cA^-)]^{-1}$; see \eqref{consts}. Thus, since $S_1 Q^\# \subset \frac{65}{64}\hQ$ we have 
\begin{equation}\label{small_cube}
\delta_{Q^\#} \leq S_1^{-1} \delta_{\frac{65}{64}  \hQ}  \leq c_*(\cA^-).
\end{equation}
Hence, the induction hypothesis implies that 
\begin{equation} \label{nottagged}
S_0 Q^\# \; \mbox{is not tagged with} \; (\cA',\epsilon_1(\cA^-) ) \; \mbox{for any} \; \cA' < \cA.
\end{equation}
In particular, $S_0 Q^\#$ is not tagged with $(\cA',\epsilon^\kappa)$ for any $\cA' < \cA$.

Since $S_0 Q^\#$ is tagged with $(\cA,\epsilon^\kappa)$ but not with $(\cA',\epsilon^\kappa)$ for any $\cA' < \cA$, we see that $\sigma(S_0 Q^\#)$ has an $(\cA,x_{Q^\#},\epsilon^\kappa,\delta_{S_0Q^\#})$-basis. Thus there exist polynomials $(P_\alpha)_{\alpha \in \cA}$ such that
\begin{align} 
\label{dun1} &P_\alpha \in \epsilon^\kappa \left[ \delta_{S_0Q^\#} \right]^{|\alpha| + n/p - m} \cdot \sigma(S_0 Q^\#) \qquad \mbox{for} \; \alpha \in \cA, \\
\label{dun2}  \partial^\beta &P_\alpha(x_{Q^\#}) = \delta_{\alpha \beta} \qquad\qquad\qquad\qquad\qquad\;\;          \mbox{for} \; \beta,\alpha \in \cA,\\
\label{dun3}  \lvert\partial^\beta &P_\alpha(x_{Q^\#}) \rvert  \leq \epsilon^\kappa \left[\delta_{S_0Q^\#} \right]^{|\alpha| - |\beta|}  \qquad\qquad\quad \mbox{for} \; \beta \in \cM,\; \alpha \in \cA, \; \beta > \alpha.
\end{align}

To start, we prove the  following statement.
\begin{itemize}
\item Suppose that 
\begin{align} 
\label{insigma0}
&R \in \sigma(S_0 Q^\#), \; \mbox{and} \\
\label{zeroder}
&\partial^\alpha R(x_{Q^\#}) = 0 \; \mbox{for all} \; \alpha \in \cA.
\end{align}
Then, for some $W = W(m,n,p) \geq 0 $ we have
\begin{equation} \label{smallsigma} \lvert \partial^\beta R(x_{Q^\#}) \rvert \leq W \cdot \delta_{Q^\#}^{m - n/p - |\beta|} \quad \mbox{for} \; \beta \in \cM.\end{equation}
\end{itemize}

For the sake of contradiction, suppose that \eqref{insigma0} and \eqref{zeroder} do not imply \eqref{smallsigma}. Then, for some large constant $\widehat{W} \geq 0$, which will be determined later, there exists $R \in \sigma(S_0 Q^\#)$ satisfying \eqref{zeroder} and
\begin{equation} \label{eq_22}  \max_{\beta \in \cM} \; \lvert \partial^\beta R(x_{Q^\#}) \rvert \cdot (\delta_{Q^\#})^{ \frac{n}{p} + |\beta| - m } = \widehat{W}.
\end{equation}

For each integer $\ell \geq 0$, define the multi-index set
\[\Delta_\ell = \left\{\beta \in \cM : |\partial^\beta R(x_{Q^\#})| \cdot (\delta_{Q^\#})^{\frac{n}{p} + |\beta| - m} \geq \widehat{W}^{(2^{-\ell})}  \right\}.\] 
Note that $\Delta_0 \neq \emptyset$ thanks to \eqref{eq_22}, and also $\Delta_{\ell} \subset \Delta_{\ell + 1}$ for $\ell\geq 0$.

Since $\# \cM = D$ and $\Delta_\ell \subset \cM$ is an increasing sequence, there is an index $\ell_* \in \{ 0 , \cdots, D\}$ such that $\Delta_{\ell_*} = \Delta_{\ell_* + 1}$. Pick the maximal element $\oa \in \Delta_{\ell_*}$ (under the standard order on multi-indices defined in Section \ref{sec_multi}). Since $\oa \in \Delta_{\ell_*}$ we have
\begin{equation} \label{eq_23} \lvert \partial^\oa R(x_{Q^\#}) \rvert \cdot (\delta_{Q^\#})^{ \frac{n}{p} + |\oa| - m } \geq \widehat{W}^{(2^{-\ell_*})}.
\end{equation}
Now, if $\beta \in \cM$ and $\beta > \oa$, then $\beta \notin \Delta_{\ell_*} = \Delta_{\ell_*+1}$ by the maximality of $\oa$. Hence,
\begin{equation}\label{eq_24}
\lvert \partial^\beta R(x_{Q^\#}) \rvert \cdot (\delta_{Q^\#})^{ \frac{n}{p} + |\beta| - m }  \leq \widehat{W}^{(2^{-\ell_* - 1})} \quad \mbox{for each} \; \beta \in \cM \; \mbox{with} \; \beta > \oa.
\end{equation}
We define $Z$ and $A$ by setting $\widehat{W} = Z^A$, with $A = 2^{\ell_*}$ and $0 \leq \ell_* \leq D$. Then \eqref{eq_22} - \eqref{eq_24} state that
\begin{align*}
|\partial^{\oa} R(x_{Q^\#})| &\geq Z \cdot [\delta_{Q^\#}]^{m-n/p-|\oa|}, \\
|\partial^{\beta} R(x_{Q^\#})| &\leq Z^{1/2} \cdot [\delta_{Q^\#}]^{m-n/p-|\beta|} \quad \mbox{for} \;  \beta \in \cM, \;\beta > \oa, \;\; \mbox{and} \\
|\partial^{\beta} R(x_{Q^\#})| &\leq Z^{A} \cdot [\delta_{Q^\#}]^{m-n/p-|\beta|} \quad\;\; \mbox{for} \; \beta \in \cM.
\end{align*}
Define $\oP_{\oa} := ( \partial^\oa R(x_{Q^\#}) )^{-1} \cdot R$. Then \eqref{zeroder} implies that
\begin{equation} \label{e725} \partial^\beta \oP_{\oa}(x_{Q^\#}) = \delta_{\alpha \beta} \;\; \mbox{for} \;  \beta \in \cA \cup \{\oa\}.
\end{equation}
Since $R \in \sigma(S_0 Q^\#)$,
\begin{equation} 
\label{e726} \oP_{\oa} \in Z^{-1} \left[ \delta_{Q^\#} \right]^{|\oa| + n/p - m}  \cdot \sigma(S_0 Q^\#),
\end{equation}
and also
\begin{align}
\label{e727} \lvert \partial^\beta \oP_{\oa}(x_{Q^\#}) \rvert &\leq Z^{-1/2} \cdot [\delta_{Q^\#}]^{|\oa| - |\beta|} \qquad \mbox{for} \; \beta \in \cM, \; \beta > \oa, \; \mbox{and} \\
\label{e728} \lvert \partial^\beta \oP_{\oa}(x_{Q^\#}) \rvert &\leq Z^A \cdot [\delta_{Q^\#}]^{|\oa| - |\beta|} \qquad\; \mbox{for} \; \beta \in \cM.
\end{align}
Set $\oA := \{ \alpha \in \cA : \alpha < \oa\} \cup \{\oa\}$.  Since $\partial^{\alpha} R(x_{Q^\#}) = 0$ for all $\alpha \in \cA$, we see that $\oa \notin \cA$. Thus, 
\begin{equation} \label{sss}
\oA < \cA.
\end{equation}
For $\alpha \in \cA$ with $ \alpha < \oa$, we set 
$$\oP_{\alpha} := P_\alpha - \partial^{\oa} P_\alpha(x_{Q^\#})  \oP_{\oa}.$$ 
Note that
\begin{align}
\label{e729} & \oP_\alpha \in C \epsilon^\kappa Z^A \left[\delta_{Q^\#}\right]^{|\alpha| + n/p - m}  \cdot \sigma(S_0 Q^\#) \qquad \qquad  (\mbox{by \eqref{dun1}, \eqref{dun3}, \eqref{e726}}), \\
\label{e730} & \partial^\beta \oP_\alpha(x_{Q^\#}) = \delta_{\alpha \beta} \qquad\qquad\qquad\;\;\; \mbox{for} \;\; \beta \in \oA \qquad\qquad\;\; (\mbox{by \eqref{dun2}, \eqref{e725}}), \; \mbox{and} \\
\label{e731} & \lvert \partial^\beta \oP_\alpha(x_{Q^\#}) \rvert \leq C \epsilon^\kappa Z^A \left[\delta_{Q^\#}\right]^{|\alpha| - |\beta|} \;\; \mbox{for} \;\; \beta \in \cM, \; \beta > \alpha \quad (\mbox{by \eqref{dun3}, \eqref{e728}}).
\end{align}

Examining \eqref{e725}-\eqref{e731}, we see that 
\begin{equation}
\label{e732}
(\oP_\alpha)_{\alpha \in \oA} \; \mbox{forms an} \; (\oA,x_{Q^\#}, C \cdot (S_0)^{m}  \max \{\epsilon^\kappa Z^A, Z^{-1/2}\}, \delta_{S_0 Q^\#}) \mbox{-basis for} \; \sigma(S_0 Q^\#).
\end{equation}
Recall that $Z = \widehat{W}^{1/A}$ with $A = 2^{\ell_*}$ and $0 \leq \ell_* \leq D$. We pick $\widehat{W}$ to be a large enough universal constant, and assume that $\epsilon$ is less than a small enough universal constant. Then \eqref{sss} and \eqref{e732} imply that $S_0 Q^\#$ is tagged with $(\cA^-, \epsilon_1(\cA^-) )$. But this contradicts \eqref{nottagged}. This completes our proof that \eqref{smallsigma} holds whenever the polynomial $R$ satisfies \eqref{insigma0} and \eqref{zeroder}.

We now prove the main assertion in Proposition \ref{key_prop1}. \label{pp10} Suppose that $H \in \X$ satisfies $H = f$ on $E \cap S_1 Q^\#$ and $\partial^\alpha H (x_{Q^\#}) = \partial^\alpha P(x_{Q^\#})$ for all $\alpha \in \cA$. Then $\partial^\alpha (J_{x_{Q^\#}}H -  P) \equiv 0$ for all $\alpha \in \cA$. (Recall, $\cA$ is monotonic; see Remark \ref{mon_rem}.) We apply the estimate \eqref{eqstuff} followed by Lemma \ref{key_lem1}, and hence, we see that
\begin{align*}
M^\#_{Q^\#}(f,R^\#_{Q^\#}) &\leq C \cdot M^\#_{Q^\#}(f,J_{x_{Q^\#}}H) \\ 
&\leq C \|(f,J_{x_{Q^\#}}H) \|_{S_1 Q^\#} \leq C \|H\|_{\X(S_1Q^\#)},
\end{align*}
which implies that
\[
M^\#_{Q^\#}(0, R^\#_{Q^\#} - J_{x_{Q^\#}}H) \leq C \|H\|_{X(S_1Q^\#)}.
\]
Thus, Lemma \ref{key_lem1} implies that $\|(0 , R^\#_{Q^\#} - J_{x_{Q^\#}}H ) \|_{S_0Q^\#} \leq C \|H\|_{\X(S_1Q^\#)}$, hence 
\[
R^\#_{Q^\#} - J_{x_{Q^\#}}H \in C \|H\|_{\X(S_1 Q^\#)}  \cdot \sigma(S_0 Q^\#).
\]
By the defining properties of $R^\#_{Q^\#}$ (see the algorithm \textsc{Make New Assists and Assign Keystone Jets}), and by our assumption on $\partial^\alpha H(x_{Q^\#})$, we have
\[\partial^\alpha ( R^\#_{Q^\#} - J_{x_{Q^\#}}H)(x_{Q^\#}) = \partial^\alpha(P - P)(x_{Q^\#}) = 0 \quad \mbox{for all} \; \alpha \in \cA.
\]
Thus, \eqref{smallsigma} shows that
\[
| \partial^\beta(J_{x_{Q^\#}}H - R^\#_{Q^\#})(x_{Q^\#})| \leq C \cdot (\delta_{Q^\#})^{m-n/p-|\beta|} \|H\|_{\X(S_1 Q^\#)} \quad \mbox{for all} \; \beta \in \cM.
\]
Hence, by the Sobolev inequality we have
\begin{equation*}
\delta_{Q^\#}^{-m} \| H - R_{Q^\#}^\#  \|_{L^p(S_1 Q^\#)} \leq C \cdot \| H \|_{\X(S_1 Q^\#)}.
\end{equation*}
That proves \eqref{imp} and completes the proof of Proposition \ref{key_prop1}.
\end{proof}

\subsection{Marked Cubes} \label{mk_cubes}

We summarize various objects that we have computed in previous sections of the paper. This is meant to serve as a reference for the reader.

\begin{itemize}
\item \underline{The main cubes}: We compute the collection of cubes $Q \in \CZ_{\main}(\cA^-)$, each marked with pointers to the following objects.
\begin{itemize}
\item The list $\Omega(Q, \cA^-)$ of assist functionals on $\X(\frac{65}{64} Q \cap E)$, expressed in short form.
\item The list $\Xi(Q,\cA^-)$ of functionals on $\X(\frac{65}{64}Q \cap E) \oplus \cP$, which have $\Omega(Q,\cA^-)$-assisted bounded depth, expressed in short form in terms of assists $\Omega(Q,\cA^-)$.
\item The list of functionals $\xi^Q_1,\cdots,\xi^Q_D$ on  $\cP$.
\end{itemize}
(See the Main Technical Results for $\cA^-$ and the algorithm \textsc{Approximate Old Trace Norm} in Section \ref{sec_ind_step}.)
\item \underline{The keystone cubes}: We list all the keystone cubes $Q^\#$ for $\CZ(\cA^-)$, each marked with pointers to the following objects.
\begin{itemize}
\item The list $\Omega^\new(Q^\#)$ of new assist functionals on $\X(S_1 Q^\# \cap E)$, expressed in short form.
\item The linear map $R^\#_{Q^\#} : \X(S_1 Q^\# \cap E) \oplus \cP \rightarrow \cP$, which has $\Omega^\new(Q^\#)$-assisted bounded depth, and is expressed in short form in terms of assists $\Omega^\new(Q^\#)$.
\end{itemize}
(See \textsc{Make New Assists and Assign Keystone Jets} in Section \ref{sec_assign}.)
\item \underline{The border-dispute pairs}: We list all the border-dispute pairs $(Q',Q'') \in \BD(\cA^-)$. \\
(See the \textsc{Keystone-Oracle}  in Section \ref{key_sec}.)
\end{itemize}
We store these cubes in memory along with their markings.

\subsection{Testing Cubes} \label{sec_testcube}

Let $\hQ$ be a dyadic subcube of $Q^\circ$. Since $\CZ(\cA^-)$ is a dyadic decomposition of $Q^\circ$, one and only one of the following alternatives holds.
\begin{enumerate}[(A)]
\item $\hQ$  is a disjoint union of cubes from $\CZ(\cA^-)$.
\item $\hQ$  is strictly contained in one of the cubes of $\CZ(\cA^-)$.
\end{enumerate}

\begin{defn}
\label{testing_defn}
Let $\hQ \subset Q^\circ$ be a dyadic cube. If alternative (A) holds, we call $\hQ$ a \underline{testing cube}. 

Let $0 < \lambda < 1$. We say that a testing cube $\hQ$ is \underline{$\lambda$-simple} if $\delta_Q \geq \lambda \cdot \delta_{\hQ}$ for any $Q \in \CZ(\cA^-)$ with $Q \subset (65/64)\hQ$.

\end{defn}

We introduce a geometric parameter
\begin{equation} 
\label{tdefn}
t_G \in \R, \; \mbox{which is an integer power of two}.
\end{equation}
We assume that
$$ 0 < t_G < c, \; \mbox{where} \; c \; \mbox{is a small enough constant determined by} \; m,n,p.$$
We will later determine $a(\cA)$ to be an appropriate constant depending on $t_G$. For the main conditions satisfied by $a(\cA)$, see the fourth and fifth bullet points in Chapter \ref{sec_mainresults}. Near the end of Section \ref{sec_test} we determine $t_G$ to be a constant depending only on $m$, $n$, and $p$ - but not yet. 

We recall that $a = a(\cA^-)$ is a fixed universal constant.

\begin{lem}\label{lem_cover} Let $\hQ$ be a testing cube. Assume that $t_G > 0$ is less than a small enough universal constant. The following properties hold.
\begin{itemize}
\item There exists a constant $a_\new > 0$, depending only on $t_G,m,n,p$, such that the cube $(1+a_\new)\hQ$ is contained in the union of the cubes $(1+\frac{a}{2})Q$ over all $Q \in \CZ(\cA^-)$ with $Q \subset (1+t_G)\hQ$.
\item If $Q \in \CZ(\cA^-)$ and $Q \subset (1+100t_G)\hQ$, then $\frac{65}{64}Q \subset \frac{65}{64}\hQ$.
\end{itemize}
\end{lem}

\begin{proof}
We assume $a_\new$ is less than a small enough constant determined by $t_G$, $m$, $n$, and $p$. We will later fix $a_\new$ to be a constant depending only on $t_G$, $m$, $n$, and $p$, but not yet.

Let $x \in (1+a_\new)\hQ$ be given. We will produce a cube $Q \in \CZ(\cA^-)$ with $Q \subset (1+t_G)\hQ$ such that $x \in (1+ \frac{a}{2})Q$, thus proving the first bullet point.

Pick a point $x_\near \in \hQ$ with $\lvert x_\near - x \rvert \leq a_\new \delta_\hQ$. (Recall that we use the $\ell^\infty$ metric on $\R^n$.) 

Since the cubes in $\CZ(\cA^-)$ partition $Q^\circ$, one of the following cases must occur

\noindent \textbf{Case 1:} There exists $Q_1 \in \CZ(\cA^-)$ with $\delta_{Q_1} \leq (t_G/40) \delta_\hQ$ such that $x \in Q_1$.

Because $x \in (1+a_\new)\hQ$, we have $Q_1 \subset (1+a_\new+\frac{t_G}{10})\hQ$ in Case 1. Therefore, $Q_1 \subset (1+t_G)\hQ$. (Here, we assume that $a_\new \leq \frac{9t_G}{10}$.)

\noindent \textbf{Case 2:} There exists $Q_2 \in \CZ(\cA^-)$ with $\delta_{Q_2} > (t_G/40)\delta_\hQ$ and $x \in Q_2$.

Because $\hQ$ is a testing cube, there exists $Q \in \CZ(\cA^-)$ such that $Q \subset \hQ$ and $x_\near \in Q$. 

Moreover, note that
\[\lvert x - x_\near \rvert \leq a_\new \delta_\hQ \leq \frac{40 a_\new}{t_G} \delta_{Q_2} \leq \frac{a}{8} \delta_{Q_2}.\]
(Here, we assume that $a_\new \leq \frac{a t_G}{320}$.)  The above estimate and the fact that $x \in Q_2$ imply that $x_\near \in (1+a)Q_2$. Since $x_\near \in Q$, we have $\delta_{Q_2} \leq 2 \delta_Q$ by good geometry. Therefore, $\lvert x - x_\near \rvert \leq \frac{a}{4} \delta_{Q}$. Consequently, since $x_\near \in Q$ we have $x \in (1+\frac{a}{2})Q$.

\noindent \textbf{Case 3:} $x \in \R^n \setminus Q^\circ$. 

Because $\hQ$ is a testing cube, there exists $Q \in \CZ(\cA^-)$ with $Q \subset \hQ$ and $x_\near \in Q$. Note that
\[ \dist(Q, \R^n \setminus Q^\circ) \leq | x - x_\near | \leq a_\new \delta_\hQ \leq a_\new. \] 
If $a_\new < 10^{-3}$, then Lemma \ref{bdry_lem} implies that $\delta_Q \in \{1/2,1/4,1/8\}$, hence $\lvert x - x_\near \rvert \leq 8 a_\new \delta_Q$. Since $x_\near \in Q$, we see that $x \in (1+ 100a_\new) Q \subset (1+a/2)Q$.

Thus, in all cases we have produced some cube $Q' \in \CZ(\cA^-)$ such that $Q' \subset (1+t_G) \hQ$ and $x \in (1+a/2)Q'$. Here, $x \in (1+a_\new)\hQ$ is arbitrary. We now fix $a_\new$ to be a small enough constant depending on $t_G$, $m$, $n$, and $p$. This completes the proof of the first bullet point.

We now prove the second bullet point. 

We assume we are given a cube $Q \in \CZ(\cA^-)$ with $Q \subset (1+100t_G)\hQ$.

Since $\hQ$ is a testing cube, either $Q \subset \hQ$ or $Q \subset (1+100t_G)\hQ \setminus \hQ$. 

In the former case, clearly $\frac{65}{64} Q \subset \frac{65}{64} \hQ$. 

In the latter case, we have $\delta_Q \leq 50t_G \delta_\hQ$ and so $\frac{65}{64}Q \subset (1+ 1000t_G) \hQ \subset \frac{65}{64} \hQ$. 

This proves the second bullet point and completes the proof of the lemma.
\end{proof}

\subsection{Testing Functionals}
\label{sec_as}

We recall that we have computed linear maps $R^\#_{Q^\#}$ associated to the keystone cubes $Q^\#$ in $\CZ(\cA^-)$. See the algorithm \textsc{Make New Assists and Assign Keystone Jets}. 

We assume we are given a parameter $t_G$ as in \eqref{tdefn}.

We assume we are given a testing cube $\hQ \subset Q^\circ$. (See Definition \ref{testing_defn}.)

For each $Q \in \CZ(\cA^-)$ with $Q \subseteq (1+100 t_G) \hQ$, we define
\begin{equation}\label{jet1}
   R^{\hQ}_Q(f,P) := \left\{
     \begin{array}{lr}
       P & : \delta_Q \geq t_G \delta_{\hQ}\\
       R^\#_{\mathcal{K}(Q)}(f,P) & : \delta_Q < t_G \delta_{\hQ}
     \end{array}   \qquad (\mbox{for any} \; (f,P) \in \X((65/64)\hQ \cap E) \oplus \cP).
   \right.
\end{equation}
We guarantee that $S_1 \mathcal{K}(Q) \subset C Q$ as in the \textsc{Keystone-Oracle} in Section \ref{key_sec}. If $\delta_Q < t_G \delta_{\hQ}$, then $CQ \subset (1+Ct_G) \hQ$. For small enough $t_G$, we conclude that 
\begin{equation} \label{subcc}
S_1 \mathcal{K}(Q) \subset ( 65/64) \hQ.
\end{equation}
This shows that the map $R^\hQ_Q$ is well-defined.

We define the ``testing functional'' $[ M_{\hQ}(f,P)]^p$ to be the sum of the following terms.

\begin{align} \label{i} \text{\bf (I)} = &\text{ the sum of } \bigl[ M_{(Q,\cA^-)}(f,R^\hQ_Q(f,P)) \bigr]^p = \sum_{\xi \in \Xi(Q,\cA^-)} \bigl\lvert \xi \bigl(f,R^\hQ_Q(f,P) \bigr) \bigr\rvert^p \\
&\qquad \text{over all } Q \in \CZ_{\main}(\cA^-) \text{ such that } Q \subset (1+t_G)\hQ. \notag{} \\
\label{ii} \text{\bf (II)} = &\text{ the sum of } \sum_{ \beta \in \cM} \delta_{Q'}^{n-(m-|\beta|)p} \left\lvert \partial^\beta \left[ R^\hQ_{Q'}(f,P) - R^\hQ_{Q''}(f,P) \right](x_{Q'}) \right\rvert^p \\
& \qquad \text{over all} \; (Q', Q'') \in \BD(\cA^-) \; \text{such that} \; Q' \subset (1+t_G)\hQ, \; \delta_{Q'} < t_G \delta_{\hQ}. \notag{} \\
\label{iii} \text{\bf (III)} = &\text{ the sum of } \sum_{\beta \in \cM} \delta_Q^{n-(m-|\beta|)p} \left\lvert \partial^\beta \left[ R^\hQ_Q(f,P) - P \right](x_Q) \right\rvert^p \\
& \qquad \text{over all} \; Q \in \CZ(\cA^-) \; \text{such that} \; Q \subset (1+t_G)\hQ,\; \delta_Q \geq t_G^2 \delta_{\hQ}. \notag{} \\
\label{iv} \text{\bf (IV)} = &\text{ the sum of } \sum_{\beta \in \cM} \delta_{\hQ}^{n-(m-|\beta|)p} \left\lvert \partial^\beta \left[ R^\hQ_{Q_{\spec}}(f,P) - P \right](x_{\hQ}) \right\rvert^p \\
& \qquad \text{for a single (arbitrarily chosen)} \; Q_{\spec} \in \CZ(\cA^-) \; \text{contained in} \; \hQ. \notag{}
\end{align}
(Note that $Q'' \subset (1+100t_G)\hQ$ in \eqref{ii}, thanks to the good geometry of cubes in $\CZ(\cA^-)$; hence the sum \textbf{(II)} is well-defined.)

Thus we have defined a functional $M_{\hQ}(f,P)$. Although $M_{\hQ}(f,P)$ depends on the parameter $t_G$, we leave this dependence implicit in our notation for the sake of brevity.

For each testing cube $\hQ$, we define
\begin{equation}\label{s1}\ooline{\sigma}(\hQ) = \left\{ P \in \cP : M_{\hQ}(0,P) \leq 1 \right\}.\end{equation}

\environmentA{Algorithm: Approximate New Trace Norm.}\label{alg_ANII}

Given a number $t_G > 0$ as in \eqref{tdefn}, we perform one-time work at most $C(t_G) N \log N$ in space $C(t_G)N$, after which we can answer queries.

A query consists of a testing cube $\hQ$.

The response to the query $\hQ$ is a list $\mu_1^{\hQ},\ldots,\mu_D^{\hQ}$ of linear functionals on $\cP$ such that 
\begin{equation}
\label{appp}
c \left[ M_{\hQ}(0,P) \right]^p \leq  \sum_{i=1}^D \lvert \mu_i^{\hQ}(P) \rvert^p \leq C \left[ M_{\hQ}(0,P)\right]^p.
\end{equation}
Define a quadratic form on $\cP$ by 
\begin{equation}
\label{quadform}
q_{\hQ}(P) := \sum_{i=1}^D \lvert \mu_i^{\hQ}(P) \rvert^2.
\end{equation}
This quadratic form satisfies
\begin{equation}
\label{qform_bd1}
c \left[ M_{\hQ}(0,P) \right]^2 \leq  q_{\hQ}(P) \leq C \left[ M_{\hQ}(0,P) \right]^2.
\end{equation}
In particular,
\begin{equation}
\label{qform_bd2}
\{ q_\hQ \leq c \} \subset \ooline{\sigma}(\hQ) \subset \{ q_\hQ \leq C\}.
\end{equation}
The work required to answer a query is at most $C(t_G) \log N$.

\begin{proof}[\underline{Explanation}] 

For each keystone cube $Q^\# \in \CZ(\cA^-)$ and each $\beta \in \cM$, we have stored a short form description of the $\Omega^\new(Q^\#)$-assisted bounded depth linear functional $(f,P) \mapsto \partial^\beta \left[R^\#_{Q^\#}(f,P)\right](0)$. This corresponds to an expansion
\[\partial^\beta \left[R^\#_{Q^\#}(f,P)\right](0) = \lambda_{(Q^\#,\beta)}(f) + \overline{\lambda}_{(Q^\#,\beta)}(P); \]
here, $\lambda_{(Q^\#,\beta)}(f)$ and $\overline{\lambda}_{(Q^\#,\beta)}(P)$ are linear functionals (with $\lambda_{(Q^\#,\beta)}$ given in short form in terms of some set of assists). We mark each keystone cube $Q^\#$ with the linear map 
\begin{equation}
\label{polymap}
P \mapsto R^\#_{Q^\#}(0, P) = \sum_{\beta \in \cM} \overline{\lambda}_{(Q^\#,\beta)}(P) \cdot \frac{1}{\beta!} x^\beta.
\end{equation}
This requires work and storage at most $C$ for each $Q^\#$. (We simply produce the functionals $\overline{\lambda}_{(Q^\#,\beta)} : \cP \rightarrow \R$  for all $\beta \in \cM$.) The number of keystone cubes is at most $CN$, hence this computation requires total work at most $C N$.

We now perform the marking procedure described below.
\begin{itemize}
\item For each cube $Q \in \CZ_{\main}(\cA^-)$, we mark $Q$ with the linear functionals
$$\xi_{(Q,i)}(P) := \xi_i^Q\left(R^\#_{\mathcal{K}(Q)}(0,P)\right) \qquad (i=1,\cdots,D).$$
To compute these functionals we simply compose linear maps that were already computed. The functionals $\xi_i^Q$ on $\cP$ satisfy \eqref{ineq2}, and are computed using the algorithm \textsc{Approximate Old Trace Norm}. We produce the keystone cube $\mathcal{K}(Q)$ using the \textsc{Keystone-Oracle}. We locate the map $P \mapsto R^{\#}_{\mathcal{K}(Q)}(0,P)$ using a binary search.

This requires work at most $C \log N$ for each given $Q \in \CZ_{\main}(\cA^-)$. (The binary search requires work at most $C \log N$.)

\item For each border-dispute pair $(Q', Q'') \in \BD(\cA^-)$, we mark $Q'$ with linear functionals
$$\xi_{(Q',Q'',\beta)}(P) :=  \delta_{Q'}^{n/p - m+|\beta|} \partial^\beta \left\{ R^\#_{\mathcal{K}(Q')}(0,P) - R^\#_{\mathcal{K}(Q'')}(0,P) \right\}(x_{Q'}) \qquad (\beta \in \cM).$$
The linear maps $P \mapsto R^{\#}_{\mathcal{K}(Q')}(0,P)$ and $P \mapsto R^{\#}_{\mathcal{K}(Q'')}(0,P)$ are computed using the \textsc{Keystone-Oracle} and a binary search, as in the previous bullet point.

This requires work at most $C \log N$ for each given $(Q',Q'') \in \BD(\cA^-)$. 
\end{itemize}

Each relevant cube is marked with at most $\mathcal{O}(1)$ functionals by the above bullet points. Since the number of cubes $Q$ and $Q'$ arising above is at most $C N$, the marking procedure requires work at most $C N \log N$ in space $C N$.

We perform the one-time work for the algorithm \textsc{Compute Norms From Marked Cuboids} on the marked cubes $Q$, $Q'$ arising above, which is at most $C N \log N$ work in space $C N$. Again, we use the fact that each cube is marked by $\mathcal{O}(1)$ functionals. This concludes the one-time work for the present algorithm.

We now explain the query work. Suppose that $\hQ$ is a given testing cube (a query).

We partition $(1+t_G)\hQ$ into dyadic cubes $ Q_1,\ldots,Q_L \subset \R^n$ such that $\delta_{Q_\ell} = (t_G/4)\delta_{\hQ}$. Note that $L = L(t_G)$ is a constant determined by $n$ and $t_G$. (Recall that $0 < t_G < 1$ is  an integer power of 2; see \eqref{tdefn}.)

Note that
\begin{align}\label{fact}
& Q \in \CZ(\cA^-), \; Q \subset (1+t_G)\hQ, \; \mbox{and} \;  \delta_Q \leq (t_G/4) \delta_\hQ \iff \\
& Q \in \CZ(\cA^-) \; \mbox{and} \; Q \subset Q_\ell \;\; \mbox{for some} \; \ell \in \{1,\ldots,L\}.
\notag{}
\end{align}

Next, we apply the query algorithm from \textsc{Compute Norms From Marked Cuboids} with each cube $Q_\ell$ used as a query ($\ell=1,\cdots, L$). We obtain linear functionals $\mu^{Q_\ell}_1,\cdots, \mu^{Q_\ell}_D$ on $\cP$ such that
\begin{equation}
\label{wc1}
c  \sum_{k=1}^D \lvert \mu^{Q_\ell}_k(P) \rvert^p  \leq  \sum_{\substack{ Q \in \CZ(\cA^-) \\ \text{linear functional} \; \xi } } \bigl\{ \lvert \xi(P) \rvert^p : Q \subset Q_\ell,  \; Q \; \mbox{marked with} \; \xi \bigr\} \leq C  \sum_{k=1}^D \lvert \mu^{Q_\ell}_k(P) \rvert^p.
\end{equation}
This requires work and storage at most $C \log N$ for each fixed $\ell$, and total work and storage at most $C(t_G) \log N$. Summing the above estimate from $\ell=1,\ldots,L$ and using \eqref{fact}, we learn that
\begin{align}
\sum_{\ell = 1}^L \sum_{k=1}^D\lvert \mu^{Q_\ell}_k(P) \rvert^p \sim & \sum   \biggl\{ \sum_{i=1}^D \lvert \xi^Q_i(R^\#_{\mathcal{K}(Q)}(0,P)) \rvert^p : Q \in \CZ_{\main}(\cA^-), \; Q \subset (1+t_G)\hQ, \; \delta_Q \leq \frac{t_G}{4} \delta_{\hQ} \biggr\}.\notag{} \\
\label{e763} 
& + \sum \biggl\{   \lvert\xi_{(Q',Q'',\beta)}(P)\rvert^p : \;  (Q',Q'') \in \BD(\cA^-), \; Q' \subset (1+t_G)\hQ,  \\
& \qquad\qquad\qquad\qquad\qquad\qquad\qquad\qquad\qquad\qquad\; \delta_{Q'} \leq \frac{t_G}{4}\delta_{\hQ}, \; \beta \in \cM \biggr\} \notag{}\\
& =: \mathfrak{S}_1 + \mathfrak{S}_2. \notag{}
\end{align}

We now compute the functionals described below.
\begin{align*}
&\mathbf{(F_1)} \;\; \boxed{\mu_k^{Q_\ell}(P)} \hspace{2cm} \text{ for }k=1,\cdots,D, \; \ell=1,\cdots,L.  \\
&\mathbf{(F_2)} \;\; \boxed{\xi^Q_i(R^\#_{\mathcal{K}(Q)}(0,P))}     \quad \text{ for } i=1,\cdots,D, \; Q \in \CZ_{\main}(\cA^-), \; Q \subset (1+t_G)\hQ, \; \delta_Q = \frac{t_G}{2} \delta_{\hQ}. \\
&\mathbf{(F_3)} \;\; \boxed{\xi^Q_i(P)}     \hspace{2cm}\; \text{ for }  i=1,\cdots,D, \; Q \in \CZ_{\main}(\cA^-), \; Q \subset (1+t_G)\hQ, \; \delta_Q \geq t_G \delta_{\hQ}. \\
&\mathbf{(F_4)} \;\; \boxed{\xi_{(Q',Q'',\beta)}(P)}  \qquad\; \text{ for }  \beta \in \cM, \; (Q',Q'') \in \BD(\cA^-), \\
& \hspace{10cm}  \delta_{Q'} = \frac{t_G}{2}\delta_{\hQ}, \; \delta_{Q''} \leq \frac{t_G}{2} \delta_\hQ. \\
&\mathbf{(F_5)} \;\; \boxed{\delta_Q^{n/p- m + |\beta|} \left\{ \partial^\beta (R^\hQ_Q(0,P) - P)(x_Q)\right\}}  \qquad \text{ for } \beta \in \cM, \; Q \in \CZ(\cA^-), \\
&  \hspace{10cm}  Q \subset (1+t_G)\hQ, \; \delta_Q \geq t_G^2 \delta_{\hQ}. \\
&\mathbf{(F_6)} \;\; \boxed{\delta_{\hQ}^{n/p-m+|\beta|} \left\{ \partial^\beta (R^\hQ_{Q_{\spec}}(0,P) - P)(x_{\hQ}) \right\}} \quad\text{ for } \beta \in \cM. \\
\end{align*}

The number of functionals listed here is at most $C(t_G)$. To compute these functionals, we proceed as follows.

We have already produced the functionals in $\mathbf{(F_1)}$ that satisfy \eqref{e763}.

We can compute the functionals arising in $\mathbf{(F_6)}$. If $\delta_{Q_{\spec}} \geq t_G \delta_\hQ$ then the functionals in $\mathbf{(F_6)}$ vanish identically. If instead $\delta_{Q_{\spec}} < t_G \delta_\hQ$ then the map $R^\hQ_{Q_{\spec}}(0,P) = R^\#_{\mathcal{K}(Q_{\spec})}(0,P)$ has been computed, and we easily produce the expression in $\mathbf{(F_6)}$. 

Next, we loop over all dyadic cubes $Q \subset (1+t_G)\hQ$ with $\delta_Q \geq t_G^2 \delta_\hQ$. For each such $Q$, we do the following. 

If $\delta_Q = \frac{t_G}{2} \delta_\hQ$ and $Q \in \CZ_{\main}(\cA^-)$ then we compute the functional in $\mathbf{(F_2)}$. 

If $\delta_Q \geq t_G \delta_\hQ$ and  $Q \in \CZ_{\main}(\cA^-)$ then we compute the functional in $\mathbf{(F_3)}$. 

If $Q \in \CZ(\cA^-)$ then we can compute the functionals in $\mathbf{(F_5)}$. These functionals are identically zero whenever $\delta_Q \geq t_G \delta_\hQ$. Otherwise, since we have already computed the map $R^\hQ_{Q}(0,P) = R^\#_{\mathcal{K}(Q)}(0,P)$, we can easily compute the expression in $\mathbf{(F_5)}$. That concludes the loop over $Q$.

Finally, we loop over the dyadic cubes $Q' \subset (1+t_G)\hQ$ with  $ \delta_{Q'} = \frac{t_G}{2} \delta_\hQ$. If $Q' \in \CZ(\cA^-)$, then we loop over $Q'' \in \CZ(\cA^-)$ such that $Q'' \leftrightarrow Q'$. If $\delta_{Q''} \leq (t_G/2) \delta_\hQ$ and $\mathcal{K}(Q'') \neq \mathcal{K}(Q')$ then we compute the functionals arising in $\mathbf{(F_4)}$. That concludes the loop over $Q'$.

Thus we have computed all the functionals arising in $\mathbf{(F_1)}$-$\mathbf{(F_6)}$. We define the functional $\bigl[X(P)\bigr]^p$ to be the sum of the $p$-th powers of all these functionals.

We will now show that $\bigl[X(P)\bigr]^p$ well approximates $\left[ M_\hQ(0,P) \right]^p$.

The sum of the $p$-th powers of the functionals arising in $\mathbf{(F_1)}$ is estimated in \eqref{e763}. We obtain from this the estimate
\begin{align} 
\label{Xeq}
\bigl[X(P)\bigr]^p \sim &  \sum   \biggl\{ \sum_{i=1}^D \lvert \xi^Q_i(R^\hQ_Q(0,P)) \rvert^p : Q \in \CZ_{\main}(\cA^-), \; Q \subset (1+t_G)\hQ\biggr\} \\ 
\notag{}
 & \quad +   \sum \biggl\{   \lvert\xi_{(Q',Q'',\beta)}(P)\rvert^p : \;  (Q',Q'') \in \BD(\cA^-), \; Q' \subset (1+t_G)\hQ, \; \\
 \notag{}
& \hspace{6cm} \delta_{Q'} \leq \frac{t_G}{2} \delta_{\hQ}, \; \delta_{Q''} \leq \frac{t_G}{2} \delta_\hQ, \; \beta \in \cM \biggr\}\\
\notag{}
 & \quad +  \mathfrak{S}_3 + \mathfrak{S}_4.
\end{align}
Here, $\mathfrak{S}_3$ and $\mathfrak{S}_4$ are the terms \textbf{(III)} and \textbf{(IV)}, respectively, with $f$ set to $0$ (see \eqref{iii} and \eqref{iv}). 
Let us explain how we obtained this formula. The sum of the term $\mathfrak{S}_1$ in \eqref{e763} and the sum of the $p$-th powers of all the functionals in $\mathbf{(F_2)}$ and $\mathbf{(F_3)}$ is equal to the first line in \eqref{Xeq}. (Recall the definition of $R_Q^\hQ$ in \eqref{jet1}.) The sum of the term $\mathfrak{S}_2$ in \eqref{e763} and the sum of the $p$-th powers of all the functionals in $\mathbf{(F_4)}$ is equal to the second line in \eqref{Xeq}. The sum of the $p$-th powers of all the functionals in $\mathbf{(F_5)}$ and $\mathbf{(F_6)}$ is equal to the third line in \eqref{Xeq}, i.e., the quantity $\mathfrak{S}_3 + \mathfrak{S}_4$.

The sum in the first line in \eqref{Xeq} is comparable to the term \textbf{(I)} with $f \equiv 0$ (see \eqref{i}), thanks to the estimate \eqref{ineq2}. Note that the sum in the second line in \eqref{Xeq} is equal to the term \textbf{(II)} with $f \equiv 0$ (see \eqref{ii}) minus all the summands in \textbf{(II)} with $\delta_{Q'} = (t_G/2) \delta_\hQ$ and $\delta_{Q''} = t_G \delta_\hQ$. (Recall that by good geometry the sidelengths of $Q'$ and $Q''$ can differ by at most a factor of $2$.) However, these discarded summands appear also in the term \textbf{(III)}. Thus, $\left[X(P)\right]^p$ is comparable to the sum of the terms \textbf{(I)},\textbf{(II)},\textbf{(III)},\textbf{(IV)} (with $f \equiv 0$). Thus, in summary, we have
$$c \cdot \left[ M_\hQ(0,P)\right]^p \leq \bigl[ X(P) \bigr]^p \leq C \cdot \left[ M_\hQ(0,P) \right]^p$$
for universal constants $c > 0$ and $C \geq 1$.

Processing the functionals in $\mathbf{(F_1)}$-$\mathbf{(F_6)}$ using the algorithm \textsc{Compress Norms} (Section \ref{sec_lf}), we compute functionals $\mu_1^{\hQ},\ldots,\mu_D^{\hQ}$ on $\cP$ such that
$$c \cdot \sum_{i=1}^D \lvert \mu_i^\hQ(P) \rvert^p \leq \left[X(P)\right]^p \leq C \cdot\sum_{i=1}^D \lvert \mu_i^\hQ(P) \rvert^p.$$
The previous two estimates imply the desired estimate \eqref{appp}.

The estimate \eqref{qform_bd1}, concerning the quadratic form $q_{\hQ}(P)$ defined in \eqref{quadform}, follows because the $\ell_p$ and $\ell_2$ norms on the space $\R^D$ are comparable up to a constant factor depending on $D$, which is, in turn, a universal constant. (Recall that $D = \dim(\cP)$ depends only on $m$ and $n$.) The pair of inclusions in \eqref{qform_bd2} follows directly from \eqref{qform_bd1} and the definition of $\ooline{\sigma}(\hQ)$ in \eqref{s1}.

This completes the description of the query work, which consists of at most $C(t_G) \log N$ computer operations.

This completes the explanation of the algorithm \textsc{Approximate New Trace Norm}.

\end{proof}

\subsection{Computing Data Associated to a Testing Cube} \label{sec_computingstuff}

Let $\hQ$ be a testing cube (see Definition \ref{testing_defn}), and let $t_G > 0$ be as in \eqref{tdefn}.

The \underline{supporting data for $\hQ$} consists of the following:
\begin{itemize}
\item[(SD1)] Pointers to the cubes $Q \in \CZ_{\main}(\cA^-)$ with $Q \subset (1+t_G)\hQ$. \\
(These are the cubes appearing in the sum \textbf{(I)}; see \eqref{i}.)
\item[(SD2)] Pointers to the pairs $(Q', Q'') \in \BD(\cA^-)$ with $Q' \subset (1+t_G)\hQ$ and $\delta_{Q'} < t_G \delta_{\hQ}$. \\
(These are the pairs of cubes appearing in the sum \textbf{(II)}; see \eqref{ii}.)
\item[(SD3)] Pointers to the cubes $Q \in \CZ(\cA^-)$ with $Q \subset (1+t_G)\hQ$ and $\delta_Q \geq t_G^2 \delta_{\hQ}$. \\
(These are the cubes appearing in the sum \textbf{(III)}; see \eqref{iii}.)
\item[(SD4)] A pointer to a cube $Q_{\spec} \in \CZ(\cA^-)$ with $Q_{\spec} \subset \hQ$. \\
(This cube appears in \textbf{(IV)}; see \eqref{iv}.)
\item[(SD5)] Pointers to the keystone cubes $Q^\#$ of $\CZ(\cA^-)$ with $S_1 Q^\# \subset (65/64)\hQ$. \\
(See \eqref{consts} for the definition of $S_1$.)
\end{itemize}

We are given markings as in Section \ref{mk_cubes}. Each cube $Q \in \CZ_{\main}(\cA^-)$ is marked with pointers to the lists $\Omega(Q,\cA^-)$ and $\Xi(Q,\cA^-)$, and each keystone cube $Q^\# \in \CZ(\cA^-)$ is marked with a pointer to the list $\Omega^{\new}(Q^\#)$.  We define
\begin{align} \label{test_assists}
\Omega(\hQ) := &\biggl[ \bigcup \left\{ \Omega(Q,\cA^-) : Q \in \CZ_{\main}(\cA^-), Q \subset (1+t_G)\hQ \right\} \biggr] \cup \\
& \biggl[  \bigcup \left\{ \Omega^{\new}(Q^\#) : Q^\# \in \CZ(\cA^-) \; \mbox{keystone},  \; S_1 Q^\# \subset (65/64) \hQ\right\} \biggr]. \notag{}
\end{align}
Using the supporting data for $\hQ$ and the above markings, we produce a list of all the functionals in $\Omega(\hQ)$. To form the list \eqref{test_assists}, we examine all the relevant $Q$ and $Q^\#$, and we copy each assist functional $\omega$ from $\Omega(Q,\cA^-)$ or $\Omega^{\new}(Q^\#)$ into a location in memory. The work and space required are bounded by the sum of the depths of all the $\omega$ that are copied. We make no attempt to remove duplicates in the list \eqref{test_assists}. For more details about our notation concerning unions of lists, see Section \ref{sec_not}. We summarize the procedure in the following algorithm.

\environmentA{Algorithm: Compute New Assists.}

Given a testing cube $\hQ$, and given the supporting data for $\hQ$, we compute a list of all the functionals in $\Omega(\hQ)$. We mark all the functionals that appear in the lists $\Omega(Q,\cA^-)$ (for $Q \in \CZ_{\main}(\cA^-)$, $Q \subset (1+t_G)\hQ$ in the supporting data) and $\Omega^{\new}(Q^\#)$ (for $Q^\#$ keystone, $S_1 Q^\# \subset (65/64)\hQ$ in the supporting data) with pointers to their position in the list $\Omega(\hQ)$. 
This requires work at most
\begin{align} \label{work1}
\Work_1(\hQ) = C  \log N  \cdot \biggl[ 1 &+ \sum_{\substack{ Q \in \CZ_{\main}(\cA^-) \\ Q \subset (1+t_G)\hQ}}  \sum_{\omega \in \Omega(Q,\cA^-)}\depth(\omega)  \\
& \hspace{1cm} + \sum_{\substack{ \mbox{\tiny keystone} \; Q^\# \in \CZ(\cA^-)  \\ S_1 Q^\# \subset \frac{65}{64} \hQ  }} \sum_{\omega \in \Omega^\new(Q^\#)} \depth(\omega) \biggr] \notag{}
\end{align}
and storage at most
\begin{align} \label{space1} \Space_1(\hQ) = C \cdot \biggl[ 1 &+ \sum_{\substack{ Q \in \CZ_{\main}(\cA^-) \\ Q \subset (1+t_G)\hQ}}  \sum_{\omega \in \Omega(Q,\cA^-)}\depth(\omega)  \\
& \hspace{1cm}  + \sum_{\substack{ \mbox{\tiny keystone} \; Q^\# \in \CZ(\cA^-)  \\ S_1 Q^\# \subset \frac{65}{64} \hQ  }} \sum_{\omega \in \Omega^\new(Q^\#)} \depth(\omega) \biggr]. \notag{}
\end{align}

\begin{remk}
\label{remk_translate} 
Let $\hQ$ be a testing cube, and let $\xi$ be a linear functional that has $\Omega(Q,\cA^-)$-assisted bounded depth for $Q \in \CZ_{\main}(\cA^-)$, $Q \subset (1+t_G)\hQ$ relevant to the supporting data for $\hQ$. Then $\xi$ has $\Omega(\hQ)$-assisted bounded depth, since $\Omega(Q, \cA^-)$ is a sublist of $\Omega(\hQ)$. If $\xi$ is given in short form in terms of the assists $\Omega(Q, \cA^-)$, then we can convert $\xi$ into a short form in terms of the assists $\Omega(\hQ)$. That is because we have marked each  functional in $\Omega(Q, \cA^-)$ with a pointer to its position in the list $\Omega(\hQ)$. The conversion requires a constant amount of work once we have carried out the algorithm \textsc{Compute New Assists} for the given $\hQ$. 

Similarly, let $\xi$ be a linear functional that has $\Omega^\new(Q^\#)$-assisted bounded depth, for some $Q^\#$ relevant to the supporting data for $\hQ$. Given a short form description of $\xi$ in terms of the assists $\Omega^\new(Q^\#)$, we  can express $\xi$ in short form in terms of the assists $\Omega(\hQ)$ using a constant amount of work.
\end{remk}

\environmentA{Algorithm: Compute Supporting Map.}

We perform one-time work at most $C N \log N$ in space $CN$, after which we can answer queries as follows.

A query consists of a testing cube $\hQ$, its supporting data, and a cube $Q \in \CZ(\cA^-)$ with $Q \subset (1+100t_G)\hQ$.

The response to a query $(\hQ,Q)$ is a short form description of the linear map $R^\hQ_Q : \X((65/64)\hQ \cap E) \oplus \cP \rightarrow \cP$ in terms of the assists $\Omega(\hQ)$ (see \eqref{jet1}).

The work and storage required to answer a query are at most $C \log N$.

(Here, we do not count the storage used to hold the supporting data for $\hQ$.)

\begin{proof}[\underline{Explanation}]  We simply use the definition in \eqref{jet1}.

We first test to see whether $\delta_Q < t_G \delta_\hQ$ or $\delta_Q \geq t_G \delta_\hQ$.

In the first case when $\delta_Q \geq t_G \delta_{\hQ}$, we have $R^\hQ_Q(f,P) = P$, and we produce a short-form description of this map.

In the second case when $\delta_Q < t_G \delta_{\hQ}$, we compute the map $R^\hQ_Q = R^\#_{\mathcal{K}(Q)}$ as follows.

First, we compute the keystone cube $Q^\# = \mathcal{K}(Q)$ using the \textsc{Keystone-Oracle}. Recall that $S_1 Q^\# \subset \frac{65}{64}\hQ$ (see \eqref{subcc}). We locate $Q^\#$ in the list of pointers appearing in (SD5) using a binary search. We have already computed the $\Omega^\new(Q^\#)$-assisted bounded depth linear map 
\[ R^\#_{Q^\#} : \X((65/64)\hQ \cap E) \oplus \cP \rightarrow \cP\]
in short form in terms of the assists $\Omega^\new(Q^\#)$, as described in Section \ref{mk_cubes}. Thanks to Remark \ref{remk_translate} we can express $R^\#_{Q^\#}$ in short form in terms of the assists $\Omega(\hQ)$.

Thus we have computed the desired expression for $R_Q^\hQ = R^\#_{Q^\#}$ in the second case.

That concludes the explanation of the algorithm.

\end{proof}

\environmentA{Algorithm: Compute New Assisted Functionals.}

Given a testing cube $\hQ$ and its supporting data, we produce a list $\Xi(\hQ)$ consisting of $\Omega(\hQ)$-assisted bounded depth functionals on $\X(\frac{65}{64}\hQ \cap E)  \oplus \cP$, with each functional written in short form, such that the following hold.
\begin{itemize}
\item $ \displaystyle \left[ M_{\hQ}(f,P) \right]^p =   \sum_{\xi \in \Xi(\hQ)} \lvert \xi(f,P) \rvert^p$ for each $(f,P) \in \X(\frac{65}{64}\hQ \cap E)  \oplus \cP$.
\item Denote
\begin{equation} \label{bd1} \mathfrak{N}(\hQ) := \# \bigl\{ (Q',Q'') \in \BD(\cA^-) : \; Q' \subset (1+t_G)\hQ, \;  \delta_{Q'} < t_G \delta_{\hQ} \bigr\}.
\end{equation}
We carry out the preceding computation using work at most 
\begin{equation} \label{work2}
\Work_2(\hQ) := C(t_G) \cdot \log N \cdot  \biggl[ 1 + \mathfrak{N}(\hQ) + \sum_{\substack{Q \in \CZ_{\main}(\cA^-) \\ Q \subset (1+t_G)\hQ }} \# \bigl[ \Xi(Q,\cA^-) \bigr]  \biggr]
\end{equation}
in space
\begin{equation} \label{space2} \Space_2(\hQ) := C(t_G) \cdot  \biggl[ 1 + \mathfrak{N}(\hQ) + \sum_{\substack{Q \in \CZ_{\main}(\cA^-) \\ Q \subset (1+t_G)\hQ }} \# \bigl[ \Xi(Q,\cA^-) \bigr] \biggr].
\end{equation}
In particular, $\# \Xi(\hQ) \leq \Space_2(\hQ)$.
\end{itemize}
(Again, we don't count the space used to hold the supporting data for $\hQ$.)

\begin{proof}[\underline{Explanation}]  

We compute the list $\Xi(\hQ)$ of all the functionals appearing in the sums \textbf{(I)}-\textbf{(IV)} in \eqref{i}-\eqref{iv}.

We loop over all the cubes $Q \in \CZ_{\main}(\cA^-)$ with $Q \subset (1+t_G)\hQ$ (as in (SD1)).

We form the functionals
\begin{equation}\label{i1}
(f,P) \mapsto \xi(f,R^\hQ_Q(f,P))   \quad (\text{for} \; \xi \in \Xi(Q,\cA^-)).
\end{equation}
The linear maps $R^\hQ_{Q}$ are written in short form in terms of the assists $\Omega(\hQ)$ (see the algorithm \textsc{Compute Supporting Map}). The functionals $\xi \in \Xi(Q,\cA^-)$ are written in short form in terms of assists $\Omega(Q,\cA^-)$. We can write the functionals $\xi \in \Xi(Q,\cA^-)$ in short form in terms of the assists $\Omega(\hQ)$ (see Remark \ref{remk_translate}). Hence, we can express each functional in \eqref{i1} in short form in terms of assists $\Omega(\hQ)$. This requires work at most $C$ for each $\xi$.

That concludes the loop on $Q$.

We now loop over all pairs $(Q',Q'') \in \BD(\cA^-)$ with $Q' \subset (1+t_G)\hQ$ and $\delta_{Q'} < t_G \delta_\hQ$ (as in (SD2)). For each such pair, we compute the $\Omega(\hQ)$-assisted bounded depth linear maps $R^\hQ_{Q'}$ and $R^\hQ_{Q''}$ in short form. We form the functionals
\begin{equation}
\label{ii1} \delta_{Q'}^{n/p- m + |\beta|} \left\{ \partial^\beta (R^\hQ_{Q'}(f,P) - R^\hQ_{Q''}(f,P))(x_{Q'}) \right\} \qquad (\mbox{for} \; \beta \in \cM).
\end{equation}
That concludes the loop on $(Q',Q'')$.

We loop over all the cubes $Q \in \CZ(\cA^-)$ such that $Q \subset (1+t_G)\hQ$ and $\delta_Q \geq t_G^2 \delta_\hQ$ (as in (SD3)). We form the functionals
\begin{equation} \label{iii1}
\delta_Q^{n/p- m + |\beta|} \left\{ \partial^\beta (R^\hQ_Q(f,P) - P)(x_Q)\right\}\qquad (\mbox{for} \;\beta \in \cM).
\end{equation}
That concludes the loop on $Q$.

We form the functionals
\begin{equation}\label{iv1}
\delta_\hQ^{n/p- m + |\beta|} \left\{ \partial^\beta (R^\hQ_{Q_\spec}(f,P) - P)(x_\hQ)\right\} \qquad (\mbox{for} \;\beta \in \cM).
\end{equation}
Here, we use the cube $Q_\spec$ in (SD4). 

Let $\Xi(\hQ)$ denote the list of functionals arising in \eqref{i1}-\eqref{iv1}. All these functionals have $\Omega(\hQ)$-assisted bounded depth and are expressed in short form in terms of assists $\Omega(\hQ)$. Comparing with \eqref{i}-\eqref{iv}, we see that $[M_\hQ(f,P)]^p$ is equal to the sum of $\lvert \xi(f,P) \rvert^p$ over all $\xi \in \Xi(\hQ)$. Clearly, the number of functionals in $\Xi(\hQ)$ is bounded by
\begin{equation*}
C(t_G) \cdot  \biggl[ 1 + \mathfrak{N}(\hQ) + \sum_{\substack{Q \in \CZ_{\main}(\cA^-) \\ Q \subset (1+t_G)\hQ }} \# \bigl[ \Xi(Q,\cA^-) \bigr] \biggr].
\end{equation*}
Since we perform work at most $C \log N$ (using storage at most $C$) to compute each functional, the total work and storage used by our algorithm are at most $\Work_2(\hQ)$ and $\Space_2(\hQ)$, respectively.

\end{proof}

\noindent\textbf{The extension operator.}

Given a testing cube $\hQ$, the \underline{covering cubes} for $\hQ$ are
\begin{equation}\label{scrI} \cov(\hQ) := \bigl\{ Q \in \CZ(\cA^-) : Q \subset (1+t_G)\hQ\bigr\}.\end{equation}

We assume that
\begin{equation}
\label{small_t}
t_G \; \mbox{satisfies the hypothesis of Lemma \ref{lem_cover}}.
\end{equation}
We do not fix $t_G$ just yet. Let $a_\new = a_\new(t_G)$ be as in Lemma \ref{lem_cover}. Thus, 
\begin{equation}
\label{inert}
(1+a_\new)\hQ \; \mbox{is contained in the union of the cubes} \; (1+a/2)Q \; \mbox{as} \; Q \; \mbox{ranges over} \; \cov(\hQ).
\end{equation}

The assumptions in Sections \ref{czalg_sec} and \ref{sec_pou} are valid, where
\begin{itemize}
\item $\CZ = \CZ(\cA^-)$ and $\cQ = \cov(\hQ)$.
\item The cube called $\hQ$ in Sections \ref{czalg_sec} and \ref{sec_pou} given by the cube $(1+a_\new)\hQ$ as in the present section.
\item $\overline{r} = a$, and $A = C$ for a large enough universal constant $C$.
\end{itemize}
The good geometry of $\CZ(\cA^-)$ and the existence of a $\CZ(\cA^-)$-\textsc{Oracle} follow from the Main Technical Results for $\cA^-$ (see Chapter \ref{sec_mainresults}). Regarding the conditions in Section \ref{sec_pou}: condition \eqref{covers} is stated in \eqref{inert}, while condition \eqref{sizebd} is a direct consequence of the definition of $\cov(\hQ)$. 

We may thus apply the results stated in Section \ref{sec_pou}.

By Lemma \ref{pou_lem}, there exist cutoff functions $\theta^\hQ_Q \in C^m(\R^n)$ such that
\begin{align} 
\label{pou701} 
& \sum_{Q \in \cov(\hQ)}  \theta^\hQ_Q  = 1 \; \mbox{on} \; (1+a_\new)\hQ,\\
\label{pou702} 
&\supp ( \theta_Q^\hQ) \subset (1+a)Q \;\; \mbox{and} \;\;  \lvert \partial^\alpha \theta^\hQ_Q \rvert \leq C \cdot \delta_Q^{-|\alpha|} \;\; \mbox{for} \; |\alpha| \leq m, \; \mbox{and} \\
\label{pou703}
& \theta_Q^\hQ = 1 \; \mbox{near} \; x_Q, \; \mbox{and} \; \theta_Q^\hQ = 0 \; \mbox{near} \; x_{Q'} \; \mbox{for each} \; Q' \in \cov(\hQ) \setminus \{Q\}.
\end{align}

\environmentA{Algorithm: Compute POU.} \label{alg_pou}

After one-time work at most $C N \log N$ in space $CN$, we can answer queries as follows.

A query consists of a testing cube $\hQ$ and a point $\underline{x} \in Q^\circ$. 

The response to the query $(\hQ, \underline{x})$ is a list of all the cubes $Q_1,\cdots,Q_L \in \cov(\hQ)$ (with $Q_1,\cdots,Q_L$ all distinct) such that $\underline{x} \in \frac{65}{64} Q_\ell$, and the list of polynomials $J_{\underline{x}} \theta^\hQ_{Q_1},\cdots, J_{\underline{x}} \theta^\hQ_{Q_L}$.

To answer a query requires work and storage at most $C \log N$.

\begin{proof}[\underline{Explanation}] 

We list all the cubes $Q \in \CZ(\cA^-)$ for which $\underline{x} \in \frac{65}{64}Q$ using the $\CZ(\cA^-)$-\textsc{Oracle}. We then discard any cubes that are not contained in $(1+t_G)\hQ$. The remaining cubes give the desired list $Q_1,\cdots,Q_L$.

We now compute the jet $J_{\underline{x}} \theta^\hQ_{Q_\ell}$ for each $\ell$. From the proof of Lemma \ref{pou_lem}, we have
\[ \theta^\hQ_{Q_\ell} = \widetilde{\theta}_{Q_\ell} \cdot \bigl[ \eta \circ \Psi \bigr]^{-1}, \qquad \mbox{where} \;\; \Psi = \sum_{Q \in \cov(\hQ)} \widetilde{\theta}_{Q}.\]
Applying the algorithm \textsc{Compute Cutoff Function} (Section \ref{sec_pou}), we compute the jet $J_{\underline{x}} \widetilde{\theta}_{Q_\ell}$ for each $\ell = 1 ,\cdots,L$. 
We can compute a formula for $\partial^\alpha J_{\underline{x}} ( \theta_{Q_\ell}^\hQ)(\underline{x})$ given a formula for the jet $J_{\underline{x}} ( \eta \circ \Psi )$. Indeed, by the Leibniz rule, $\partial^\alpha J_{\underline{x}} ( \theta_{Q_\ell}^\hQ)(\underline{x})$ ($\lv \alpha \rv \leq m-1$) is given by a rational function of the derivatives $\partial^\beta J_{\underline{x}}( \widetilde{\theta}_{Q_\ell})(\underline{x})$ and $\partial^\beta J_{\underline{x}} (\eta \circ \Psi)(\underline{x})$ ($\lv \beta \rv \leq m-1$).

Since each $\widetilde{\theta}_{Q_\ell}$ is supported on $\frac{65}{64}Q_\ell$, we have
\[J_{\underline{x}} \Psi = \sum_{\ell=1}^L J_{\underline{x}} \widetilde{\theta}_{Q_\ell}. \]
Recall from the proof of Lemma \ref{pou_lem} that the function $\eta : [0,\infty) \rightarrow \R$ satisfies $\eta(t) \geq 1/4$ for $t \in [0,1/2)$, and $\eta(t) = t$ for $t \in [1/2,\infty)$. Given $t_* \geq 0$ and $k \leq m$, we assume that the number $\frac{d^k \eta}{d t^k}(t_*)$ can be computed using work and storage at most $C$. This can be achieved by taking $\eta$ to be a suitable spline function. Thus, the jet $J_{\underline{x}} (\eta \circ \Psi)$ can be computed using the chain rule. 

Thus, we can compute the jets $J_{\underline{x}} ( \theta^\hQ_{Q_\ell})$ using work and storage at most $C$ once we know the list $Q_1,\cdots,Q_L$.

\end{proof}

Let $(f,P) \in \X(\frac{65}{64}\hQ \cap E) \oplus \cP$ be given.

For ease of notation we write $R^\hQ_Q = R^\hQ_Q(f,P)$ for the polyomial defined in \eqref{jet1} (the dependence on $(f,P)$ should be understood).

For each $Q \in \cov(\hQ)$ we define
\begin{equation}
\label{test_eo_aux}
   F^{\hQ}_Q := 
   \left\{
     \begin{array}{lr}
       T_{(Q,\cA^-)}(f,R^\hQ_Q) & : \mbox{if} \; \frac{65}{64}Q \cap E \neq \emptyset  \\
       R^\hQ_Q & :  \mbox{if} \; \frac{65}{64}Q \cap E = \emptyset.
     \end{array}
   \right.
\end{equation}
Note that the function $F_Q^\hQ \in \X$ is well-defined. According to the Main Technical Results for $\cA^-$ (see Chapter \ref{sec_mainresults}) we have
\begin{equation}\label{local_ext}
F_Q^\hQ = f \; \mbox{on} \; (1+a)Q \cap E.
\end{equation}
\begin{equation} \label{local_est} 
\| F_Q^\hQ \|_{\X((1+a)Q)} + \delta_Q^{-m} \| F_Q^\hQ - R_Q^\hQ \|_{L^p((1+a)Q)} \leq  \left\{
\begin{array}{lr}
       C M_{(Q,\cA^-)}(f,R_Q^\hQ) & : \mbox{if} \; \frac{65}{64}Q \cap E \neq \emptyset   \\
       0 & :  \mbox{if} \; \frac{65}{64}Q \cap E = \emptyset
     \end{array}
 \right.
\end{equation}
(Recall that $a = a(\cA^-) \leq 1/64$; see \eqref{a_defn}.)

Finally, we define
\begin{equation}\label{extopdefn}
T_{\hQ}(f,P) := \sum_{Q \in \cov(\hQ)} F^\hQ_Q \cdot \theta^\hQ_Q \in \X, \;\; \mbox{with} \; \theta_Q^\hQ \; \mbox{as in \eqref{pou701}-\eqref{pou703}}.
\end{equation}

\environmentA{Algorithm: Compute New Extension Operator.}

We perform one-time work at most $C N \log N$ in space $CN$, after which we can answer queries.

A query consists of a testing cube $\hQ$, the supporting data for $\hQ$, and a point $\underline{x} \in Q^\circ$.

The response to the query $\underline{x}$ is a short form description of the $\Omega(\hQ)$-assisted bounded depth linear map
$$(f,P)  \mapsto J_{\underline{x}} T_{\hQ}(f,P).$$
To answer a query requires work at most $C \log N$.

\begin{proof}[\underline{Explanation}]  

We compute a list of the cubes $Q_1,\ldots,Q_L \in \cov(\hQ)$ (with $Q_1,\cdots,Q_L$ all distinct) such that $\underline{x} \in \frac{65}{64}Q_\ell$, and a list of the jets $J_{\underline{x}} \theta^\hQ_{Q_1}, \cdots , J_{\underline{x}} \theta^\hQ_{Q_L}$. See the algorithm \textsc{Compute POU}. Recall that $L \leq C$.

Recall that $\supp(\theta_Q^\hQ) \subset \frac{65}{64}Q$. Therefore,
\begin{equation}\label{expansion} J_{\underline{x}} T_{\hQ}(f,P) = \sum_{\ell=1}^{L} J_{\underline{x}} \theta_{Q_\ell}^\hQ \odot_{\underline{x}} J_{\underline{x}} T_{(Q_\ell,\cA^-)}(f,R^\hQ_{Q_\ell}(f,P)).
\end{equation}

For each $\ell = 1,\cdots, L$, we compute (see below) the map
\begin{equation} \label{map1}
(f,R)  \mapsto J_{\underline{x}} T_{(Q_\ell,\cA^-)}(f,R) \qquad ((f,R) \in \X((65/64)\hQ \cap E) \oplus \cP).
\end{equation}
We recall the definition \eqref{scrI} of $\cov(\hQ)$. Since $Q_\ell \in \cov(\hQ)$, we have $Q_\ell \subset (1+t_G)\hQ$. Thus, Lemma \ref{lem_cover} implies that $\frac{65}{64}Q_\ell \subset \frac{65}{64}\hQ$, hence the map \eqref{map1} is well-defined.

If $\frac{65}{64} Q_\ell \cap E = \emptyset$ then $J_{\underline{x}} T_{(Q_\ell,\cA^-)}(f,R) = R$. Otherwise, if $\frac{65}{64} Q_\ell \cap E \neq \emptyset$, then we can compute the map \eqref{map1} in short form in terms of the assists $\Omega(Q, \cA^-)$, thanks to the Main Technical Results for $\cA^-$. We check whether $\frac{65}{64}Q_\ell \cap E \neq \emptyset$, by checking whether $Q_\ell$ appears in the list $\CZ_{\main}(\cA^-)$ using a binary search. We write each of the maps \eqref{map1} in short form in terms of the assists $\Omega(\hQ)$ (see Remark \ref{remk_translate}).

We compute a short form description of the $\Omega(\hQ)$-assisted bounded depth map
\[ R_{Q_\ell}^\hQ : \X\left(\frac{65}{64}\hQ \cap E \right) \oplus \cP \rightarrow \cP \quad \text{for} \; \ell=1,\cdots,L.
\]
We use the algorithm \textsc{Compute Supporting Map} (see Section \ref{sec_computingstuff}).

Substituting $R= R_{Q_\ell}^\hQ(f,P)$ in the formula for each of the maps \eqref{map1}, we can express the map $(f,P) \mapsto J_{\underline{x}} T_{\hQ}(f,P)$ in short form in terms of the assists $\Omega(\hQ)$ using \eqref{expansion}.

The query work is clearly bounded by $C \log N$.

\end{proof}

\subsection{The main estimates} \label{themainresult}

We first prove a few properties of the extension operator $T_{\hQ}$ defined in \eqref{extopdefn}.

Let $a_\new = a_\new(t_G)$ be as in Lemma \ref{lem_cover}.

\begin{prop}\label{prop_bddextop}
Let $\hQ$ be a testing cube, and let $(f,P) \in \X(\frac{65}{64}\hQ \cap E) \oplus \cP$. Then the following properties hold.
\begin{itemize}
\item $T_{\hQ}(f,P) = f$ on $(1+a_\new)\hQ \cap E$.
\item $\| T_{\hQ}(f,P) \|_{\X( (1+a_\new)\hQ)} + \delta_\hQ^{-m}  \| T_{\hQ}(f,P) - P \|_{L^p((1+a_\new)\hQ)} \leq C \cdot M_{\hQ}(f,P) .$
\end{itemize}
Here, the constant $C \geq 1$ depends only on $m$, $n$, and $p$.
\end{prop}
\begin{proof}
The first bullet point follows from \eqref{pou701}-\eqref{pou703}, \eqref{local_ext} and \eqref{extopdefn}. We now prove the second bullet point.

Recall that we defined the collection of cubes $\cov(\hQ)$ in \eqref{scrI}.

For ease of notation, we set $\overline{a} = a_\new$ throughout the proof.

We apply Lemma \ref{patch_lem} to the cube $(1+\overline{a})\hQ$, the collection $\cov(\hQ)$, the functions $F_Q^\hQ$, the polynomials $R_Q^\hQ$, and the partition of unity $\theta^\hQ_Q$ (defined for all $Q \in \cov(\hQ)$). Thus, for $G := T_{\hQ}(f,P)$ we have
\begin{align*}
\| G \|_{\X( (1+\overline{a})\hQ)}^p \lesssim & \sum_{Q \in \cov(\hQ)} \left[  \|F_Q^\hQ\|_{\X((1+a)Q)}^p + \delta_Q^{-mp} \| F_Q^\hQ - R_Q^\hQ\|^p_{L^p((1+a)Q)} \right] \\ 
& \qquad\qquad + \sum_{\substack{ Q',Q'' \in \cov(\hQ) \\  Q' \leftrightarrow Q''}} \sum_{|\beta| \leq m-1} \delta_{Q'}^{(|\beta|-m)p + n} \lvert \partial^\beta(R^\hQ_{Q'} - R^\hQ_{Q''})(x_{Q'}) \rvert^p.
\end{align*}
From \eqref{local_est}, this implies that
\begin{align}\label{midd}
\| G\|_{\X( (1+\overline{a})\hQ)}^p \lesssim & \sum_{\substack{Q \in \cov(\hQ) \\ \frac{65}{64}Q \cap E \neq \emptyset} } \left[ M_{(Q,\cA^-)}(f,R_Q^\hQ) \right]^p \\ 
& +  \sum_{\substack{ Q',Q'' \in \cov(\hQ) \\  Q' \leftrightarrow Q''}}  \sum_{|\beta| \leq m-1} \delta_{Q'}^{(|\beta|-m)p + n} \lvert \partial^\beta (R^\hQ_{Q'} - R^\hQ_{Q''})(x_{Q'}) \rvert^p \notag{} \\ 
& = A_1(f,P) + A_2(f,P).\notag{}
\end{align}

The expression $A_1(f,P)$ is equal to the sum of the terms in \eqref{i}.

We now analyze the expression $A_2(f,P)$. Suppose that $Q',Q'' \in \cov(\hQ)$ and $Q' \leftrightarrow Q''$. Then one of the following cases must occur.
\begin{enumerate}[(A)]
\item Both $\delta_{Q'}$ and $\delta_{Q''}$ are less than  $t_G \cdot \delta_\hQ$, and $\mathcal{K}(Q') = \mathcal{K}(Q'')$;
\item Both $\delta_{Q'}$ and $\delta_{Q''}$ are less than  $t_G \cdot \delta_\hQ$, and $\mathcal{K}(Q') \neq \mathcal{K}(Q'')$;
\item Both $\delta_{Q'}$ and $\delta_{Q''}$ are at least $t_G \cdot \delta_\hQ$;
\item Exactly one of $\delta_{Q'}$ and $\delta_{Q''}$ is at least $t_G \cdot \delta_\hQ$.
\end{enumerate}

If (A) or (C) occurs, then $R_{Q'}^\hQ = R_{Q''}^\hQ$, hence the summand in $A_2(f,P)$ vanishes identically; see \eqref{jet1}.

On the other hand, suppose that (D) occurs. Since $Q'$ and $Q''$ play symmetric roles in the summand from the second sum in \eqref{midd} (switching $Q'$ and $Q''$ does not change the order of magnitude of this term) we may assume that $\delta_{Q''} \geq t_G \cdot \delta_\hQ$. Hence, $R^\hQ_{Q''} = P$; see \eqref{jet1}. Since $Q' \leftrightarrow Q''$, we also have $\delta_{Q'} \geq \frac{1}{2} \delta_{Q''} \geq \frac{t_G}{2} \delta_\hQ \geq t_G^2 \delta_\hQ$ by good geometry.

The previous three paragraphs imply the following estimate:
\begin{align*}
 A_2(f,P) \lesssim  &\sum \bigl\{ \mbox{terms in \eqref{ii}} : Q', Q'' \subset (1+t_G)\hQ, \;\; \delta_{Q'}, \delta_{Q''} < t_G \cdot \delta_\hQ, \; \\
 & \qquad\qquad\qquad\qquad\qquad Q' \leftrightarrow Q'', \;\; \mathcal{K}(Q') \neq \mathcal{K}(Q'') \bigr\} \\
 + & \sum \bigl\{ \mbox{terms in \eqref{iii}} : Q \subset (1+t_G)\hQ,\; \delta_Q \geq t_G^2 \cdot \delta_\hQ \bigr\}.
\end{align*}
Thus, we have shown that
\begin{equation} \label{mainbound1} 
\|G\|_{\X((1+\overline{a})\hQ)} \leq C \cdot M_\hQ(f,P).
\end{equation}

We now estimate $\| G - P \|_{L^p((1+\overline{a})\hQ)}$. Let $Q_{\spec}$ be as in \eqref{iv}. Observe that
\begin{align*}
\delta_\hQ^{-m} \| G - P \|_{L^p((1+\overline{a})\hQ)} \lesssim   \; & \delta_\hQ^{-m} \| G - J_{x_{Q_\spec}} G \|_{L^p((1+\overline{a})\hQ)} \\ 
 + & \delta_\hQ^{-m}  \|J_{x_{Q_\spec}} G - R_{Q_\spec}^\hQ\|_{L^p((1+\overline{a})\hQ)}\\
 + & \delta_\hQ^{-m} \|R_{Q_\spec}^\hQ -P \|_{L^p((1+\overline{a})\hQ)}.
\end{align*}
We estimate each term on the right-hand side above. First, by the Sobolev inequality, 
$$\delta_\hQ^{-m} \| G - J_{x_{Q_\spec}} G \|_{L^p((1+\overline{a})\hQ)} \leq C \| G \|_{\X((1+\overline{a})\hQ)}.$$
Next, note that $J_{x_{Q_\spec}} G = J_{x_{Q_\spec}} F^\hQ_{Q_\spec}$ (see \eqref{pou703}). Thus,
\begin{align*}
\delta_\hQ^{-m} \|J_{x_{Q_\spec}} G - R_{Q_\spec}^\hQ\|_{L^p((1+\overline{a})\hQ)} & = \delta_\hQ^{-m} \|J_{x_{Q_\spec}} F_{Q_\spec}^\hQ - R_{Q_\spec}^\hQ\|_{L^p((1+\overline{a})\hQ)} \\
& \sim \lvert J_{x_{Q_\spec}} F_{Q_\spec}^\hQ - R_{Q_\spec}^\hQ \rvert_{x_{\hQ}, \delta_\hQ} \qquad \qquad (\mbox{by Lemma \ref{pnorm}})\\
& \lesssim \lvert J_{x_{Q_\spec}} F_{Q_\spec}^\hQ - R_{Q_\spec}^\hQ \rvert_{x_{Q_\spec}, \delta_{Q_\spec}}  \quad\quad \;\; (\mbox{by Lemma \ref{pnorm}}) \\
& \sim \delta_{Q_\spec}^{-m} \|J_{x_{Q_\spec}} F_{Q_\spec}^\hQ - R_{Q_\spec}^\hQ\|_{L^p(Q_\spec)}  \quad (\mbox{by Lemma \ref{pnorm}}) \\
& \lesssim \| F_{Q_\spec}^\hQ\|_{\X(Q_\spec)} + \delta_{Q_\spec}^{-m} \| F_{Q_\spec}^\hQ - R_{Q_\spec}^\hQ \|_{L^p(Q_\spec)}.
\end{align*}
According to \eqref{local_est}, this implies that
$$\delta_\hQ^{-m} \|J_{x_{Q_\spec}} G - R_{Q_\spec}^\hQ\|_{L^p((1+\overline{a})\hQ)}  \leq C \cdot M_\hQ(f,P).$$ 
Finally,
\begin{align*}
\delta_\hQ^{-m} \| R_{Q_\spec}^\hQ -P \|_{L^p((1+\overline{a})\hQ)} & \sim \sum_{|\beta| \leq m-1} \delta_\hQ^{(|\beta| - m) + n/p} \lvert \partial^\beta (R_{Q_\spec}^\hQ - P)(x_{\hQ}) \rvert \;\;\; (\mbox{by Lemma \ref{pnorm}}) \\
& \leq C \cdot M_{\hQ}(f,P) \qquad \qquad\qquad\qquad\qquad\quad\;\;  (\mbox{by \eqref{iv}}).
\end{align*}
We combine \eqref{mainbound1} and the above estimates to obtain
\begin{equation}\label{mainbound2}
\|G\|_{\X((1+\overline{a})\hQ)}  + \delta_\hQ^{-m} \| G - P \|_{L^p((1+\overline{a})\hQ)} \leq C \cdot M_\hQ(f,P).
\end{equation}
This completes the proof of Proposition \ref{prop_bddextop}.
\end{proof}

We next prove the following result.

\begin{prop}\label{prop_normappx}
Let $\hQ$ be a testing cube, and let $(f,P) \in \X(\frac{65}{64}\hQ \cap E) \oplus \cP$. Then the following inequalities hold.
\begin{description}
\item[Unconditional inequality] $\| (f,P) \|_{(1+a_\new)\hQ} \leq C \cdot M_{\hQ}(f,P).$
\item[Conditional inequality] If $3\hQ$ is tagged with $(\cA,\epsilon)$, then 
\[ M_{\hQ}(f,P) \leq C(t_G) \cdot (1/\epsilon) \cdot \|(f,P) \|_{\frac{65}{64}\hQ}.\]
\end{description}
\end{prop}

The rest of Section \ref{sec_appx} is devoted to the proof of Proposition \ref{prop_normappx}. We set $\overline{a} = a_\new$ for the remainder of the section for ease of notation.

The unconditional inequality in Proposition \ref{prop_normappx} follows easily from Proposition \ref{prop_bddextop}. Indeed, Proposition \ref{prop_bddextop} states that $T_{\hQ}(f,P) = f$ on $(1+\overline{a})\hQ \cap E$. Hence, by definition of the trace norm,
$$\| (f,P) \|_{(1+\overline{a})\hQ} \leq \| T_{\hQ}(f,P) \|_{\X( (1+\overline{a})\hQ)} + \delta_\hQ^{-m} \| T_{\hQ}(f,P) - P \|_{L^p((1+\overline{a})\hQ)}.$$
Again thanks to Proposition \ref{prop_bddextop}, the right-hand side is bounded by $C \cdot M_{\hQ}(f,P)$, which proves the unconditional inequality.

We now begin the proof of the conditional inequality in Proposition \ref{prop_normappx}. We assume that 
\begin{equation}
\label{tagged}
3 \hQ \; \mbox{is tagged with} \; (\cA,\epsilon).
\end{equation} 
and
\begin{equation} 
\label{tassump}
t_G \leq \eta, \;\; \mbox{where} \; \eta = \min \left\{ c_*(\cA^-), \bigl[ 100 \cdot S(\cA^-) \bigr]^{-1} \right\}
\end{equation}

Now, we consider two separate cases: either $\hQ$ is $\eta$-simple or $\hQ$ is not $\eta$-simple. For the definition of simple testing cubes, see Definition \ref{testing_defn}.

The conditional inequality is easy to prove in the former case.

\begin{lem}\label{lem_simple}
Suppose that a testing cube $\hQ$ is $\eta$-simple with $\eta \geq t_G$. \\
Then $M_{\hQ}(f,P) \leq C(t_G) \cdot \|(f,P) \|_{\frac{65}{64}\hQ}$, where $C(t_G)$ depends only on $m$, $n$, $p$, and $t_G$.
\end{lem}
\begin{proof}
We examine the definition of $\bigl[M_\hQ(f,P) \bigr]^p$ as a sum of terms \textbf{(I)}-\textbf{(IV)} (see \eqref{i}-\eqref{iv}).

Suppose that $Q \in \CZ(\cA^-)$ with $Q \subset (1+100t_G)\hQ$. Then $Q \subset \frac{65}{64}\hQ$ for small enough $t_G$. Our assumption that $\hQ$ is $\eta$-simple  with $\eta \geq t_G$ implies that $\delta_Q \geq  t_G \delta_\hQ$. Hence, from \eqref{jet1} we see that
\[  Q \in \CZ(\cA^-), \; Q \subset (1+100t_G)\hQ \implies R^\hQ_Q(f,P) = P.\] 

\label{van_fn}

For every $Q', Q''$ as in \eqref{ii}, by good geometry of $\CZ(\cA^-)$ we have $Q',Q'' \subset (1+100t_G) \hQ$, hence $R^\hQ_{Q'} = R^\hQ_{Q''} = P$ in \textbf{(II)}. Similarly, for each $Q$ in \eqref{iii}, we have $R^\hQ_Q = P$ in \textbf{(III)}. Similarly, $R^\hQ_{Q_{\spec}} = P$ in \textbf{(IV)}. Hence, the terms \textbf{(II)},\textbf{(III)}, and \textbf{(IV)}, all vanish, and thus $\bigl[ M_{\hQ}(f,P) \bigr]^p = \text{\textbf{(I)}}$.

We estimate the remaining term \textbf{(I)} (see \eqref{i}).

Let $Q \in \CZ_{\main}(\cA^-)$ satisfy $Q \subset (1+t_G)\hQ$. We will bound each of the summands $\left[M_{(Q,\cA)}(f,R_Q^\hQ) \right]^p$, which are relevant to the term \textbf{(I)}. As above, we have $R_Q^\hQ(f,P) = P$. Note that $\frac{65}{64}Q \subset \frac{65}{64}\hQ$ by Lemma \ref{lem_cover}. From the right-hand estimate in \eqref{n_appx}, Lemma \ref{lem_normmon}, and the estimate $\delta_Q \geq t_G \delta_\hQ$ (which follows because $\hQ$ is $\eta$-simple with $\eta \geq t_G$), we have
\[M_{(Q,\cA^-)}(f,P) \leq C \cdot \|(f,P)\|_{\frac{65}{64}Q} \leq C(t_G) \cdot \|(f,P)\|_{\frac{65}{64}\hQ}.\]
Therefore, each summand $\left[M_{(Q,\cA)}(f,R_Q^\hQ) \right]^p$ relevant to \textbf{(I)} is bounded by $C(t_G)^p \cdot \|(f,P)\|^p_{\frac{65}{64}\hQ}$.

Since $\hQ$ is $\eta$-simple, we can have $Q \subset (1+t_G)\hQ$ for no more than $C(t_G)$ of the cubes $Q \in \CZ(\cA^-)$. Hence, no more than $C(t_G)$ many cubes $Q$ arise in  \eqref{i}. Hence, by summing the estimates just obtained, we learn that $\bigl[ M_{\hQ}(f,P)  \bigr]^p \leq C(t_G) \cdot \|(f,P)\|^p_{\frac{65}{64}\hQ}$. This completes the proof of Lemma \ref{lem_simple}.
\end{proof}

If $\hQ$ is $\eta$-simple with $\eta = \min \left\{ c_*(\cA^-), \left[ 100 S(\cA^-)\right]^{-1} \right\}$, then the conditional inequality follows from Lemma \ref{lem_simple}. Here, note that the assumption \eqref{tassump} implies the hypotheses of Lemma \ref{lem_simple}.

Thus, in proving the conditional inequality, we may assume that
\begin{equation}\label{notsimple} 
\hQ \; \mbox{is \underline{not}} \; \eta\mbox{-simple,} \; \mbox{with} \; \eta = \min \{ c_*(\cA^-), \left[ 100 S(\cA^-)\right]^{-1} \}.
\end{equation} 
This is the latter, more difficult case in the dichotomy mentioned before. By definition, in this case, there exists a cube $\overline{Q} \in \CZ(\cA^-)$ with $\overline{Q} \subset \frac{65}{64}\hQ$ and  $\delta_{\overline{Q}} \leq \eta \cdot \delta_\hQ$. Hence, we note that $S(\cA^-) \overline{Q} \subset 3 \hQ$. Moreover, we have $\delta_\oQ \leq c_*(\cA^-) \delta_\hQ \leq c_*(\cA^-)$.

Thus, \textbf{(CZ2)} in the Main Technical Results for $\cA^-$ (see Chapter \ref{sec_mainresults}) implies that
\begin{equation} \label{nott} S(\cA^-) \overline{Q} \; \mbox{is not tagged with} \; (\cA^-,\epsilon_1(\cA^-)).
\end{equation}
Hence, in particular, we have $ \#(E \cap 3 \hQ) \geq \#(E \cap S(\cA^-) \overline{Q}) \geq 2$. 

Now, from \eqref{tagged} we know that $3 \hQ$ is tagged with $(\cA,\epsilon)$. Hence, since $\#(E \cap 3\hQ) \geq 2$, we know that $\sigma(3 \hQ)$ has an $(\cA',x_{\hQ},\epsilon,\delta_{3\hQ})\mbox{-basis}$, for some $\cA' \leq \cA$.

We next apply Lemma \ref{pre_lem2} to the convex set $\sigma = \sigma(3 \hQ)$. Thus, we can guarantee that there exist numbers $\Lambda \geq 1$, and  $\kappa_1 \leq {\overline{\kappa}} \leq \kappa_2$, and a multiindex set $\cA'' \leq \cA'$, such that
\begin{equation*} \sigma(3 \hQ) \; \mbox{has an} \; (\cA'',x_{\hQ},\epsilon^{\overline{\kappa}},\delta_{3\hQ},\Lambda)\mbox{-basis}, \; \mbox{where} \; \epsilon^{\overline{\kappa}} \cdot \Lambda^{100D} \leq \epsilon^{{\overline{\kappa}}/2}.\end{equation*}
Here, $\kappa_1, \kappa_2 \in (0,1]$ are universal constants. Hence, $3\hQ$ is tagged with $(\cA'', \epsilon^{\kappa_1/2})$, which implies that $S(\cA^-)\overline{Q}$ is tagged with $(\cA'', \epsilon^{\kappa'})$ for a universal constant $\kappa'>0$. (Here, we use that $S(\cA^-) \oQ \subset 3 \hQ$; see Lemma \ref{pre_lem5}). Comparing this statement and \eqref{nott}, we deduce that $\cA'' = \cA$. In summary,
\begin{equation} \label{eq_b}  
\sigma(3 \hQ) \;\; \mbox{has an} \;  (\cA,x_{\hQ},\epsilon^{\overline{\kappa}},\delta_{3\hQ},\Lambda)\mbox{-basis}, \;  \mbox{where} \; \epsilon^{\overline{\kappa}} \cdot \Lambda^{100D} \leq \epsilon^{{\overline{\kappa}}/2}.
\end{equation}
The assumptions \eqref{tagged}-\eqref{eq_b} will be used in the remainder of this section. We finish the section by completing the proof of the conditional inequality in Proposition \ref{prop_normappx} and by deriving a useful corollary.

The next result represents a main step in the proof of the conditional inequality.

\begin{prop}\label{prop_A} Assume that \eqref{tagged}-\eqref{eq_b} hold. Then there exists an $H \in \X$ such that
\begin{itemize}
\item $H = f$ on $E \cap \frac{65}{64}\hQ$.
\item $\partial^\alpha H(x_Q) = \partial^\alpha P(x_Q)$ for each $\alpha \in \cA$ and $Q \in \CZ(\cA^-)$ such that $Q \subset \frac{65}{64}\hQ$.
\item $\|H \|_{\X(\frac{65}{64}\hQ)} + \delta_{\hQ}^{-m} \|H - P \|_{L^p(\frac{65}{64}\hQ)} \leq C \Lambda^{2D+1} \cdot \|(f,P)\|_{\frac{65}{64}\hQ}.$
\end{itemize}
Here, $C \geq 1$ depends only on $m$, $n$, and $p$.
\end{prop}
\begin{proof}

We set
$$\cJ(\hQ) := \left\{ Q \in \overline{\CZ}(\cA^-) : Q \cap \frac{65}{64} \hQ \neq \emptyset \right\}.$$
Recall that the cubes in $\{ (65/64)Q : Q \in \cJ(\hQ) \}$ have bounded overlap, and that the cubes in $\cJ(\hQ)$ have good geometry, i.e., 
\[ 
\text{\bf (GG)} \;\; \mbox{If} \; Q,Q' \in \cJ(\hQ) \; \mbox{and} \; Q \leftrightarrow Q' \; \mbox{then} \; \frac{1}{16} \cdot \delta_Q \leq \delta_{Q'} \leq 16 \cdot \delta_Q.
\] 
This follows from Proposition \ref{newcz_prop}, since $ \cJ(\hQ) \subset \overline{\CZ}(\cA^-)$. We now prove that
\begin{equation}
\label{uws1}
\delta_Q \leq C \cdot \delta_{\hQ} \;\; \mbox{for each} \; Q \in \cJ(\hQ).
\end{equation}
For the sake of contradiction, assume that $\delta_Q \geq 10^5 \delta_\hQ$ for some $Q \in \cJ(\hQ)$. By definition of $\cJ(\hQ)$, we have $Q \cap \frac{65}{64}\hQ \neq \emptyset$.  Hence, since  $\delta_Q \geq 10^5 \delta_\hQ$, we see that there exists $x \in \frac{65}{64}Q \cap \hQ$. Now, $\hQ$ is partitioned into cubes in $\CZ(\cA^-)$, since $\hQ \subset Q^\circ$ is a testing cube. Thus, we can pick $Q_* \in \CZ(\cA^-)$ with $x \in Q_*$ and $Q_* \subset \hQ$. Note that $x \in \frac{65}{64}Q \cap Q_*$. By good geometry of the cubes in $\overline{CZ}(\cA^-)$, we conclude that $\delta_Q \leq 16 \cdot \delta_{Q_*}$.  Hence, $\delta_Q \leq 16 \delta_{Q_*}  \leq 16 \delta_\hQ < 10^5 \delta_\hQ$. This gives a contradiction and completes the proof of \eqref{uws1}.

For each $Q \in \cJ(\hQ)$ we select $y_Q \in Q \cap \frac{65}{64}\hQ$ such that
\begin{equation}\label{centered}
\mbox{if} \; Q \subset \frac{65}{64}\hQ \; \mbox{then} \; y_Q = x_Q \; (\mbox{the center of} \; Q).
\end{equation}

By definition of the trace seminorm $\|(\cdot, \cdot)\|_{\frac{65}{64}\hQ}$, there exists a function $F \in \X$ with
\begin{align}
&\|F \|_{\X(\frac{65}{64}\hQ)} + \delta_{\frac{65}{64}\hQ}^{-m} \|F - P \|_{L^p(\frac{65}{64}\hQ)} \leq 2 \|(f,P)\|_{\frac{65}{64}\hQ}, \; \mbox{and}\label{rrr2}\\
& F = f \;\; \mbox{on} \;\; \frac{65}{64} \hQ \cap E. \label{rrr1}
\end{align}

\noindent\textbf{Part I: Defining local basis functions}.

By \eqref{eq_b}, there exist $P_\alpha \in \cP$ and $\varphi_\alpha  \in \X$ such that
\begin{align}
& \|\varphi_\alpha\|_{\X(3 \hQ)} + \delta_{3\hQ}^{-m} \| \varphi_\alpha - P_\alpha \|_{L^p(3 \hQ)} \leq  \epsilon^{\overline{\kappa}} \delta_{3\hQ}^{|\alpha| + n/p - m} \quad (\alpha \in \cA), \label{qw2} \\
& \varphi_\alpha = 0 \;\mbox{ on} \; E \cap 3 \hQ \qquad\quad\qquad\qquad\qquad\qquad\qquad (\alpha \in \cA),\label{qw1} \\
&\lvert \partial^\beta P_\alpha(x_\hQ) - \delta_{\alpha \beta}  \rvert  \leq \left\{ \begin{array}{lr}
       \epsilon^{\overline{\kappa}} \delta_{3\hQ}^{|\alpha| - |\beta|}  & : \mbox{if} \;\beta \geq \alpha \label{qw4} \\
       \Lambda \delta_{3\hQ}^{|\alpha| - |\beta|} & : \mbox{if} \;\beta < \alpha
     \end{array}   \right. \qquad (\alpha \in \cA, \beta \in \cM), \\
& \partial^\beta P_\alpha(x_\hQ) = \delta_{\alpha \beta} \qquad\qquad\qquad\qquad\qquad\qquad\qquad (\alpha,\beta \in \cA). \label{qw3}
\end{align}
Moreover, by \eqref{qw2} and \eqref{qw4},
\begin{equation} \label{qw5}
\| \varphi_\alpha\|_{L^p(3 \hQ)} \leq \| \varphi_\alpha - P_\alpha \|_{L^p(3 \hQ)} + \| P_\alpha \|_{L^p(3 \hQ)} \leq C \epsilon^{\overline{\kappa}} \delta_\hQ^{|\alpha| + \frac{n}{p}} + C \Lambda \delta_\hQ^{|\alpha| + \frac{n}{p}}.
\end{equation}

For each $\beta \in \cM$, and each $Q \in \cJ(\hQ)$, we have
\begin{equation*}
  \lvert \partial^\beta \varphi_\alpha(y_Q) - \delta_{\alpha \beta} \rvert  \leq  \lvert \partial^\beta (\varphi_\alpha - P_\alpha)(y_Q) \rvert + \lvert \partial^\beta P_\alpha(y_Q) - \partial^\beta P_\alpha(x_\hQ) \rvert + \lvert \partial^\beta P_\alpha(x_\hQ) - \delta_{\alpha\beta}\rvert.
\end{equation*}

Moreover, Lemma \ref{si2} implies that
$$\delta_\hQ^{|\beta| + n/p - m} \lvert \partial^\beta (\varphi_\alpha - P_\alpha)(y_Q) \rvert \leq C \left( \delta_\hQ^{-m} \|\varphi_\alpha - P_\alpha \|_{L^p(3 \hQ)} + \|\varphi_\alpha \|_{\X(3 \hQ)} \right) \leq C \epsilon^{\overline{\kappa}} \delta_\hQ^{|\alpha| +n/p - m}.$$
(Recall that $y_Q \in \frac{65}{64}\hQ \subset 3 \hQ$.) Next, by a Taylor expansion we see that
\begin{align*}
 \lvert \partial^\beta P_\alpha(y_Q) - \partial^\beta P_\alpha(x_\hQ) \rvert &= \left\lvert \sum_{0 < |\gamma| \leq m-|\beta|-1} \frac{\partial^{\beta+\gamma}P_\alpha(x_\hQ)}{\gamma!} \cdot (y_Q - x_\hQ)^\gamma \right\rvert \\
&\leq \left\{ \begin{array}{lr}
       C \epsilon^{\overline{\kappa}} \delta_\hQ^{|\alpha| - |\beta|}  & : \mbox{if} \; \beta \geq \alpha \\
       C \Lambda \delta_\hQ^{|\alpha| - |\beta|} & : \mbox{if} \; \beta < \alpha
     \end{array}
   \right. \quad (\mbox{see \eqref{qw4}}).
\end{align*}
The previous three estimates and \eqref{qw4} show that
\begin{equation} \label{qw4a}
  \lvert \partial^\beta \varphi_\alpha(y_Q) - \delta_{\alpha \beta} \rvert \leq \left\{
     \begin{array}{lr}
       C \epsilon^{\overline{\kappa}} \delta_\hQ^{|\alpha| - |\beta|}  & : \beta \geq \alpha \\
       C \Lambda \delta_\hQ^{|\alpha| - |\beta|} & : \beta < \alpha
     \end{array}
   \right. \quad (\mbox{for} \; \alpha \in \cA, \beta \in \cM).
\end{equation}
In particular, the matrix $(\partial^\beta \varphi_\alpha(y_Q))_{\alpha,\beta \in \cA}$ is $(C \epsilon^{\overline{\kappa}} , C \Lambda, \delta_\hQ)$-near triangular (with $\epsilon^{\overline{\kappa}} \Lambda^{100D} \leq \epsilon^{{\overline{\kappa}}/2}$); hence, by Lemma \ref{lem0} it has an inverse matrix $( A_{\gamma \alpha}^Q )_{\gamma,\alpha \in \cA}$ such that
\begin{equation} \label{inverse} \sum_{\alpha \in \cA} A_{\gamma \alpha}^Q \cdot \partial^\beta \varphi_\alpha(y_Q)  = \delta_{\beta\gamma} \qquad  (\mbox{for all} \; \beta,\gamma \in \cA);
\end{equation}
 \begin{equation}\label{ineq11}
  \lvert A^Q_{\gamma \alpha} - \delta_{\gamma \alpha} \rvert  \leq \left\{
     \begin{array}{lr}
       C \epsilon^{\overline{\kappa}} \Lambda^{D} \delta_\hQ^{|\gamma| - |\alpha|}  & : \mbox{if} \; \alpha \geq \gamma \\
       C \Lambda^{D} \delta_\hQ^{|\gamma| - |\alpha|} & : \mbox{if} \; \alpha < \gamma
     \end{array}
   \right.  \qquad (\mbox{for all} \; \alpha,\gamma \in \cA).
\end{equation} 
We define
\begin{equation} \label{um} \varphi_{\alpha}^Q := \sum_{\beta \in \cA} A^Q_{\alpha \beta} \varphi_\beta \quad \mbox{on} \; \R^n \qquad (\mbox{for each} \; \alpha \in \cA). \\
\end{equation}
For any $Q' \in \cJ(\hQ)$, we can write
\begin{equation}  \label{bass} \begin{aligned} 
\varphi^{Q'}_\alpha = \sum_{\gamma \in \cA} \omega^{QQ'}_{\alpha \gamma} \varphi^{Q}_\gamma, \;\; \mbox{where} \; \omega^{QQ'} := A^{Q'} \cdot \left[ A^Q \right]^{-1}.
\end{aligned}
\end{equation}

For each $Q \in \cJ(\hQ)$, we have
\begin{align}
&\|\varphi^Q_\alpha \|_{\X(3 \hQ) } \leq C \epsilon^{\overline{\kappa}} \Lambda^{D} \delta_\hQ^{|\alpha| + n/p - m}.  \label{qe1} \\
&\varphi_\alpha^Q = 0 \; \mbox{on} \; E \cap 3\hQ. \label{qe2} \\
&\partial^\gamma \varphi_\alpha^Q (y_Q) = \delta_{\gamma \alpha} \qquad\qquad\qquad\qquad\;\; \mbox{for} \;\gamma \in \cA. \label{qe3} \\
&|\partial^\gamma \varphi_\alpha^Q(y_Q) | \leq C \epsilon^{\overline{\kappa}} \Lambda^{D+1} \delta_\hQ^{|\alpha| - |\gamma|} \qquad\;\;\; \mbox{for} \; \gamma \in \cM, \gamma > \alpha. \label{qe4} \\
&|\partial^\gamma \varphi_\alpha^Q(y_Q)| \leq C \Lambda^{D+1} \delta_\hQ^{|\alpha| - |\gamma|}  \qquad\quad\;\; \mbox{for} \; \gamma \in \cM. \label{qe5}
\end{align}
Here, \eqref{qe1}, \eqref{qe2}, \eqref{qe3} are immediate consequences of \eqref{qw2} and \eqref{ineq11}, \eqref{qw1}, and \eqref{inverse}, respectively. Moreover, \eqref{qe4} and \eqref{qe5} are both consequences of \eqref{qw4a} and \eqref{ineq11} (see Lemma \ref{lem00}).

We now show that there exists $Z > 0$, depending only on $m$, $n$, and $p$, such that
\begin{equation} \label{locineq} \lvert \partial^\gamma \varphi_\alpha^Q (y_Q) \rvert \leq \; Z \Lambda^{D+1} \delta_{S(\cA^-)Q}^{|\alpha| - |\gamma|} \qquad \mbox{for} \; \alpha \in \cA, \; \gamma \in \cM.
\end{equation}
For the sake of contradiction, assume that \eqref{locineq} fails to hold for some number $Z \geq 1$. We assume that $Z$ exceeds a large enough constant determined by $m$, $n$, and $p$. We later take $Z = Z(m,n,p)$, but not yet. We assume that $\epsilon$ is less than a small enough constant determined by $Z$, $m$, $n$, and $p$.

If $\delta_Q \geq \min \{c_*(\cA^-), 1/16 , [100 S(\cA^-)]^{-1} \}  \cdot \delta_\hQ$ then since also $\delta_Q \leq C \delta_{\hQ}$ (see \eqref{uws1}), the estimate \eqref{locineq} follows from \eqref{qe5}.

Alternatively, assume that 
\[ \delta_Q <  \min\{c_*(\cA^-), 1/16 ,  [100 S(\cA^-)]^{-1}\} \cdot \delta_\hQ.
\]
Thus, we have $S(\cA^-) Q \subseteq 3 \hQ$, since $Q \cap \frac{65}{64}\hQ \neq \emptyset$. Therefore, \eqref{qe2} implies that $ \varphi_\alpha^Q = 0 \; \mbox{on} \; E \cap S(\cA^-) Q$. Moreover, the Sobolev inequality and \eqref{qe1} imply that
\begin{align*}
\| \varphi_\alpha^Q \|_{\X(S(\cA^-)Q)} + \delta_{S(\cA^-)Q}^{-m} \| \varphi_\alpha^Q - J_{y_Q}\varphi_\alpha^Q \|_{L^p(S(\cA^-) Q)} &\leq C \| \varphi_\alpha^Q \|_{\X(S(\cA^-)Q)} \\ 
& \leq C \| \varphi_\alpha^Q \|_{\X(3\hQ)} \\
&\leq C \epsilon^{\overline{\kappa}} \Lambda^{D} \delta_{3 \hQ}^{|\alpha| + n/p - m} \\
&\leq C \epsilon^{\overline{\kappa}} \Lambda^{D} \delta_{S(\cA^-)Q}^{|\alpha| + n/p - m}.
\end{align*}
(In the last inequality, we use that $\delta_{S(\cA^-)Q} \leq \delta_{3\hQ}$ and $\lvert \alpha \rvert + n/p - m < 0$.)

Thus, from the previous paragraph and \eqref{qe3},\eqref{qe4}, we see that $(J_{y_Q} \varphi_\alpha^Q)_{\alpha \in \cA}$ is an $(\cA,y_Q,C \epsilon^{\overline{\kappa}} \Lambda^{D+1},\delta_{S(\cA^-)Q})$-basis for $\sigma(S(\cA^-)Q)$. 

Note that $C\epsilon^{\overline{\kappa}} \Lambda^{D+1} \leq C \epsilon^{{\overline{\kappa}}/2} \leq \epsilon^{\kappa_1/4}$ as long as $\epsilon$ is less than a small enough universal constant. Hence, $(J_{y_Q} \varphi_\alpha^Q)_{\alpha \in \cA}$ is an $(\cA,y_Q, \epsilon^{\kappa_1/4},\delta_{S(\cA^-)Q})$-basis for $\sigma(S(\cA^-)Q)$.

We are assuming that \eqref{locineq} does not hold, hence
\[ \max \left\{ \lvert \partial^\gamma \varphi_\alpha^Q (y_Q) \rvert \delta_{S(\cA^-)Q}^{|\gamma| - |\alpha|}  :  \alpha \in \cA, \; \gamma \in \cM  \right\} \geq Z.\] 
If $Z$ exceeds a large enough universal constant, and if $\epsilon^{\kappa_1/4} < Z^{-2}$, then from Lemma \ref{pre_lem6} we deduce that
$$\sigma(S(\cA^-) Q) \; \mbox{has an} \; (\cA',y_Q, Z^{-\kappa}, \delta_{S(\cA^-)Q})\mbox{-basis, with} \; \cA' < \cA.$$
Hence, $\sigma(S(\cA^-) Q)$ has an $(\cA'',x_Q, Z^{-\kappa'}, \delta_{S(\cA^-)Q})$-basis for some $\cA'' \leq \cA'$, thanks to Lemma \ref{pre_lem4}. (Here we use that $y_Q \in Q$ and $x_Q \in Q$, so $\lvert x_Q - y_Q \rvert \leq 2 \delta_Q$.) Here, $\kappa$ and $\kappa'$ are universal constants.

If $Z$ is chosen to be a large enough universal constant, we conclude that
$$\sigma(S(\cA^-)Q) \; \mbox{has an} \; (\cA'',x_Q,\epsilon_1(\cA^-), \delta_{S(\cA^-)Q})\mbox{-basis}.$$
Hence, $S(\cA^-) Q$ is tagged with $(\cA^-,\epsilon_1(\cA^-))$.
 
Recall that $Q \in \overline{\CZ}(\cA^-)$. In fact, since $\delta_Q \leq (1/16) \delta_\hQ \leq (1/16)$, condition (e) in Proposition \ref{newcz_prop} shows that $Q \in \CZ(\cA^-)$.

Since $\delta_Q \leq c_*(\cA^-)$, the previous two paragraphs contradict property \textbf{(CZ2)} of $\CZ(\cA^-)$ in Chapter \ref{sec_mainresults}. This completes the proof of \eqref{locineq} by contradiction. This concludes our analysis of the basis functions $(\varphi^Q_\alpha)_{\alpha \in \cA}$.

\noindent\textbf{Part II: Modifying the extension}.

Since $\overline{\CZ}(\cA^-)$ forms a partition of $\R^n$, we have
\begin{equation}
\label{qq1}
\frac{65}{64}\hQ \subset \bigcup_{Q \in \cJ(\hQ)} Q.
\end{equation}

The assumptions in Sections \ref{czalg_sec} and \ref{sec_pou} are valid, where
\begin{itemize}
\item $\CZ = \overline{\CZ}(\cA^-)$ and $\cQ = \cJ(\hQ)$.
\item The cube called $\hQ$ in Section \ref{czalg_sec} and Section \ref{sec_pou} given by the cube $\frac{65}{64}\hQ$ from the present setting.
\item $\overline{r} = a$, and $A = C$ for a large enough universal constant $C$.
\end{itemize}
Indeed, $\overline{\CZ}(\cA^-)$ is a decomposition of $\R^n$ into dyadic cubes that satisfies good geometry (see Proposition \ref{newcz_prop}). Regarding the conditions in Section \ref{sec_pou}: conditions \eqref{covers} and \eqref{sizebd} follow from \eqref{qq1} and \eqref{uws1}, respectively.

Thus, we may apply the results in Section \ref{sec_pou}.

By Lemma \ref{pou_lem}, there exists $\theta_Q \in C^m(\R^n)$ for $Q \in \cJ(\hQ)$ with
\begin{align*}
&(a) \;  \sum_{Q \in \cJ(\hQ)} \theta_Q = 1 \; \mbox{on} \; \frac{65}{64}\hQ, \\
&(b) \;\; \theta_Q = 1 \; \mbox{near} \; x_Q, \; \theta_Q = 0 \; \mbox{near} \; x_{Q'} \; \mbox{for} \; Q' \in \cJ(\hQ) \setminus\{Q\}, \\
& (c) \;\; \| \partial^\alpha \theta_Q \|_{L^\infty((1+a)Q)} \leq C \cdot \delta_{Q}^{-|\alpha|} \;\; \mbox{for} \; |\alpha| \leq m, \; \mbox{and} \;\; (d) \;\; \supp \theta_Q \subset (1+a)Q.
\end{align*}

We set $H := F + \widetilde{F}$ on $\R^n$, where
\begin{equation} \label{tt}
\begin{aligned}
\widetilde{F}(x) :=& \sum_{Q \in \cJ(\hQ)} \sum_{\alpha,\beta \in \cA}   \theta_{Q}(x) \cdot  \varphi_\beta(x) \cdot A^{Q} _{\alpha \beta} \cdot [\partial^\alpha(P-F)(y_{Q})] \\
=& \sum_{Q \in \cJ(\hQ)} \sum_{\alpha \in \cA} \theta_{Q}(x) \cdot \varphi^{Q}_\alpha(x) \cdot [\partial^\alpha(P-F)(y_{Q})].
\end{aligned}
\end{equation}
Note that $H$ belongs to $\X$.

Since $\varphi_\alpha = 0$ on $E \cap 3 \hQ$, we see that $\widetilde{F} = 0$ on $E \cap 3 \hQ$, hence $H = f$ on $E \cap \frac{65}{64}\hQ$; see \eqref{rrr1}. This proves the first bullet point in Proposition \ref{prop_A}.

Suppose that $Q \in \CZ(\cA^-)$ and $Q \subset \frac{65}{64}\hQ$. Then $y_Q = x_Q$, thanks to \eqref{centered}. Thus, property (b) of $\{\theta_Q\}$ states that $\theta_Q \equiv 1$ near $y_Q$, and $\theta_{Q'} \equiv 0$ near $y_Q$ for any $Q' \in \cJ(\hQ) \setminus \{Q\}$. Therefore, \eqref{qe3} and \eqref{tt} give
\begin{align*}
\partial^\gamma \widetilde{F}(y_Q) &= \sum_{\alpha \in \cA} \partial^\gamma \varphi_\alpha^Q(y_Q)  \cdot \partial^\alpha(P-F)(y_Q) \\
&= \sum_{\alpha \in \cA} \delta_{\alpha \gamma } \cdot \partial^\alpha(P-F)(y_Q)  = \partial^\gamma (P - F)(y_Q) \;\;\; \mbox{for each}  \; \gamma \in \cA.\end{align*}
Hence, $\partial^\gamma H(y_Q) = \partial^\gamma F(y_Q) + \partial^\gamma \widetilde{F}(y_Q) = \partial^\gamma P(y_Q)$ (with $y_Q = x_Q$). This proves the second bullet point in Proposition \ref{prop_A}.

\noindent \textbf{Part III: Estimating the norm}.

From property (a) of the partition of unity $(\theta_Q)$, we may write
\begin{align} 
& H = \sum_{Q \in \cJ(\hQ)} F_Q \cdot \theta_Q \;\; \mbox{on} \;\; \frac{65}{64}\hQ, \label{dd} \\ 
& \qquad \qquad \mbox{where} \;\; F_{Q} = F + \sum_{\alpha \in \cA} \varphi^{Q}_\alpha \cdot \partial^\alpha(P-F)(y_Q) \;\; \mbox{on} \; \frac{65}{64}\hQ. \notag{}
\end{align}

Before we estimate the semi-norm $\| H\|_{\X(\frac{65}{64}\hQ)}$, we present several estimates.

First, by the right-hand inequality in \eqref{s_ineq0} and by \eqref{rrr2},
\begin{align} 
\delta_\hQ^{|\alpha| + \frac{n}{p} - m} \lvert \partial^\alpha (F - P)(y_{Q}) \rvert &\leq C \cdot \left( \delta_\hQ^{-m} \|F - P \|_{L^p(\frac{65}{64} \hQ)} + \|F - P \|_{\X(\frac{65}{64} \hQ)} \right) \notag{} \\
&  \leq C' \cdot \|(f,P)\|_{\frac{65}{64}\hQ} \quad \mbox{for} \; \alpha \in \cM.\label{a1}
\end{align}

Given $Q,Q' \in \cJ(\hQ)$ such that $Q \leftrightarrow Q'$, define the rectangular boxes
\[B_1 = (1+a)Q \cap (65/64)\hQ \;\; \mbox{and} \;\; B_2 = (1+a)Q' \cap (65/64)\hQ.
\]

Since $Q \in \cJ(\hQ)$, we know that $Q \cap \frac{65}{64}\hQ \neq \emptyset$ and $\delta_Q \leq C \delta_{\hQ}$ (see \eqref{uws1}). Hence, $B_1$ is a product of $n$ intervals whose lengths are between $c \delta_Q$ and $C \delta_Q$, for universal constants $c$ and $C$. Thus, the sidelengths of $B_1$ are between $c \delta_Q$ and $C \delta_Q$.

Similarly, the sidelengths of $B_2$ are between $c \delta_{Q'}$ and $C \delta_{Q'}$.

Note that $\delta_Q$ and $\delta_{Q'}$ differ by at most a factor of $64$ thanks to good geometry.

We know that $(1+a)Q \cap (1+a)Q' \neq \emptyset$ because $Q \leftrightarrow Q'$. Since $B_1$ and $B_2$ are nonempty, the collection of cubes $\left\{ (1+a)Q, (1+a)Q', \frac{65}{64}\hQ \right\}$ have nonempty pairwise intersections, hence we conclude that there is a common point in these three cubes. \footnote{This follows from the fact that if three intervals have nonempty pairwise intersections then the three intervals share a point in common.}
Thus, $B_1 \cap B_2 \neq \emptyset$.

We have proven the following claim.

\noindent \textbf{Claim.} For any $Q, Q' \in \cJ(\hQ)$ with $Q \leftrightarrow Q'$, all the sides of the boxes
\[ B_1 = (1+a)Q \cap (65/64) \hQ \; \mbox{and} \; B_2 = (1+a)Q' \cap (65/64)\hQ
\] 
are between $c \delta_Q$ and $C \delta_Q$ for universal constants $c$ and $C$. Hence, in particular, $B_1$ and $B_2$ have aspect ratio at most a universal constant. Moreover, $B_1 \cap B_2 \neq \emptyset$. Hence, the hypotheses of Lemma \ref{si3} hold with $K$ a universal constant.

For each $\beta \in \cA$, we have
$$\partial^\beta(J_{y_{Q'}} F_{Q'} - P)(y_{Q'}) = \partial^\beta (F_{Q'} -  P)(y_{Q'}) = 0 \qquad (\mbox{see \eqref{qe3}, \eqref{dd}}).$$
Thus, $\partial^\beta(J_{y_{Q'}} F_{Q'} - P)(y_Q) = 0$ (recall, $\cA$ is monotonic; see Remark \ref{mon_rem}). Hence, $\lvert \partial^\beta (F_{Q'} - P)(y_Q) \rvert  = \lvert \partial^\beta (F_{Q'} - J_{y_{Q'}} F_{Q'})(y_Q) \rvert $, and so Lemma \ref{si3} implies that
\begin{equation} \label{eq201}\lvert \partial^\beta (F_{Q'} - P)(y_Q) \rvert   \lesssim \delta_Q^{m-|\beta| - \frac{n}{p}}  \left( \|F_{Q'}\|_{\X((1+a)Q \cap \frac{65}{64} \hQ)} + \|F_{Q'}\|_{\X((1+a)Q'  \cap \frac{65}{64} \hQ)} \right), \;\; \beta \in \cA.
\end{equation}
(Here, we use that $\lvert y_Q - y_{Q'} \rvert \leq C \delta_Q$, which is a consequence of $Q \leftrightarrow Q'$ and the good geometry of $\cJ(\hQ)$.)

From \eqref{bass} and \eqref{dd}, for $\ob \in \cM$ we have
\begin{align}
\lvert \partial^{\ob} (F_Q - F_{Q'})(y_{Q}) \rvert &= \Bigl\lvert \sum_{\beta \in \cA}  \partial^\beta (F-P)(y_{Q}) \partial^{\ob} \varphi^Q_{\beta}(y_{Q}) - \sum_{\alpha, \beta \in \cA}  \partial^\alpha (F-P)(y_{Q'})\omega^{QQ'}_{\alpha \beta}  \partial^{\ob} \varphi^Q_{\beta}(y_{Q}) \Bigr\rvert  \notag{}\\
&\leq \sum_{\beta \in \cA} \left| \partial^{\ob} \varphi^Q_{\beta}(y_{Q}) \right| \left[ \Bigl\lvert \partial^\beta (F-P)(y_{Q}) - \sum_{\alpha \in \cA} \partial^\alpha (F-P)(y_{Q'})  \omega^{QQ'}_{\alpha \beta}  \Bigr\rvert \right] \notag{}\\
& =  \sum_{\beta \in \cA} \left| \partial^{\ob} \varphi^Q_{\beta}(y_{Q}) \right| \left[ \Bigl\lvert \partial^\beta (F-P)(y_{Q}) - \sum_{\alpha \in \cA} \partial^\alpha (F-P)(y_{Q'})  \partial^\beta \varphi_\alpha^{Q'}(y_Q)  \Bigr\rvert \right] \notag{}\\
&\qquad \qquad (\mbox{note that} \; \omega_{\alpha\beta}^{QQ'} = \partial^\beta \varphi_\alpha^{Q'}(y_Q); \mbox{see \eqref{bass} and \eqref{qe3}}) \notag{} \\
&= \sum_{\beta \in \cA} \left| \partial^{\ob} \varphi^Q_{\beta}(y_{Q}) \right| \left[ \Bigl\lvert \partial^\beta (F_{Q'} - P)(y_{Q}) \Bigr\rvert \right] \notag{} \\
&\qquad \qquad (\mbox{see \eqref{dd}}) \notag{}\\
&\leq C \Lambda^{D+1} \sum_{\beta \in \cA} \delta_{Q}^{|\beta| - |\ob|} \delta_{Q}^{m-n/p-|\beta|} \bigl[ \|F_{Q'}\|_{\X((1+a) Q  \cap \frac{65}{64} \hQ)} + \|F_{Q'}\|_{\X((1+a) Q'  \cap \frac{65}{64} \hQ)} \bigr] \notag{} \\
&\qquad \qquad (\mbox{see \eqref{locineq} and \eqref{eq201}}) \notag{}\\
&\leq C \Lambda^{D+1} \cdot \delta_{Q}^{m-|\ob| - \frac{n}{p}} \bigl[ \|F_{Q'}\|_{\X((1+a) Q  \cap \frac{65}{64} \hQ)} + \|F_{Q'}\|_{\X((1+a) Q'  \cap \frac{65}{64} \hQ)} \bigr]. \label{a2}
\end{align}

We are now prepared to estimate $\| H\|_{\X(\frac{65}{64}\hQ)}$.

Applying Lemma \ref{patch_lem}, we see that
\begin{align}
\|H\|_{\X(\frac{65}{64}\hQ)}^p \lesssim & \sum_{Q \in \cJ(\hQ)}  \| F_Q \|_{\X((1+a)Q \cap \frac{65}{64}\hQ)}^p \notag{} \\
& + \sum_{Q \in \cJ(\hQ)} \delta_Q^{-mp} \| F_Q - J_{y_Q} F_Q \|^p_{L^p((1+a)Q \cap \frac{65}{64}\hQ)} \label{aa0} \\
& + \sum_{\substack{ Q,Q' \in \cJ(\hQ) \\ Q \leftrightarrow Q'}} \sum_{|\beta| \leq m-1} \delta_Q^{(|\beta| -m)p + n} \lvert \partial^\beta (J_{y_Q} F_Q - J_{y_{Q'}} F_{Q'})(y_Q) \rvert^p. \label{aa1}
\end{align}
(Here, we take $P_Q = J_{y_Q} F_Q$ in our application of Lemma \ref{patch_lem}.)

First we estimate the terms in \eqref{aa0}. From \eqref{si1b}, we obtain
\begin{equation*}
\delta_Q^{-m} \| F_Q - J_{y_Q} F_Q \|_{L^p((1+a)Q \cap \frac{65}{64}\hQ)} \leq C \| F_Q \|_{\X((1+a)Q \cap \frac{65}{64}\hQ)}.
\end{equation*}
(Here, we use that all the sides of the box $(1+a)Q \cap \frac{65}{64}\hQ$ are comparable to $\delta_Q$; see the previous \textbf{Claim}.)

Next we estimate the terms in \eqref{aa1}. For any $Q,Q' \in \cJ(\hQ)$ with $Q \leftrightarrow Q'$, we have
\begin{align*}
\lvert \partial^\beta \bigl[ J_{y_Q} F_Q  & - J_{y_{Q'}} F_{Q'} \bigr](y_Q) \rvert = \lvert \partial^\beta \bigl[ F_Q - J_{y_{Q'}} F_{Q'} \bigr](y_Q)  \rvert \notag{} \\
& \leq \lvert \partial^\beta \bigl[ F_Q - F_{Q'} \bigr](y_Q) \rvert + \lvert \partial^\beta \bigl[ F_{Q'} - J_{y_{Q'}} F_{Q'} \bigr](y_Q) \rvert  \notag{} \\
& \lesssim  \lvert \partial^\beta \bigl[ F_Q - F_{Q'} \bigr](y_Q) \rvert  + \delta_Q^{m - |\beta| - \frac{n}{p}} \bigl[\|F_{Q'} \|_{\X((1+a)Q \cap \frac{65}{64}\hQ)} + \| F_{Q'} \|_{\X((1+a)Q' \cap \frac{65}{64}\hQ)} \bigr].
\end{align*}
(Here, in the last inequality we use Lemma \ref{si3}.)

Using our previous estimates on \eqref{aa0} and \eqref{aa1}, we obtain
\begin{align*} \|H\|_{\X(\frac{65}{64}\hQ)}^p & \lesssim  \sum_{\substack{ Q,Q' \in \cJ(\hQ) \\ Q \leftrightarrow Q'} } \biggl[ \| F_{Q} \|^p_{\X((1+a)Q \cap \frac{65}{64}\hQ)}  +  \| F_{Q'} \|^p_{\X((1+a)Q \cap \frac{65}{64}\hQ)} \\
& \qquad\qquad\qquad\qquad\qquad\qquad\qquad\qquad\qquad +   \sum_{\beta \in \cM} \lvert \partial^\beta (F_Q - F_{Q'})(y_Q) \rvert^p \delta_Q^{(|\beta| - m)p + n}  \biggr] \\
&\leq C \Lambda^{(D+1)p} \sum_{ \substack{ Q, Q' \in \cJ(\hQ) \\ Q \leftrightarrow Q' }} \Bigl[ \| F_Q\|_{\X((1+a)Q  \cap \frac{65}{64} \hQ)}^p + \| F_Q\|_{\X((1+a)Q'  \cap \frac{65}{64} \hQ)}^p \Bigr]. \\
& \qquad\qquad\qquad (\mbox{by \eqref{a2}})\\
&\leq C \Lambda^{(D+1)p} \sum_{ \substack{ Q, Q' \in \cJ(\hQ) \\ Q \leftrightarrow Q' }} \biggl[ \|F\|^p_{\X((1+a)Q \cap \frac{65}{64} \hQ)}  +  \|F\|^p_{\X((1+a)Q' \cap \frac{65}{64} \hQ)} \\
& \qquad\qquad\qquad + \sum_{\alpha, \beta \in \cA} \lvert A^Q_{\alpha \beta} \rvert^p \cdot  \lvert \partial^\alpha (F - P)(y_{Q})\rvert^p \cdot \bigl( \| \varphi_\beta \|^p_{\X((1+a)Q \cap \frac{65}{64} \hQ)} + \| \varphi_\beta \|^p_{\X((1+a)Q' \cap \frac{65}{64} \hQ)} \bigr) \biggr] \\
& \qquad\qquad\qquad (\mbox{by \eqref{um} and \eqref{dd}}) \\
&\leq C \Lambda^{(2D+1)p}  \sum_{ \substack{ Q, Q' \in \cJ(\hQ) \\ Q \leftrightarrow Q' }}  \Bigl[ \| F \|_{\X((1+a)Q  \cap \frac{65}{64} \hQ)}^p + \sum_{\beta \in \cA} \delta_\hQ^{(m-|\beta|)p-n}  \|(f,P)\|^p_{\frac{65}{64}\hQ} \| \varphi_\beta \|^p_{\X((1+a)Q  \cap \frac{65}{64} \hQ)}  \Bigr] \\
& \qquad\qquad\qquad (\mbox{by \eqref{ineq11} and \eqref{a1}})\\
& \leq C \Lambda^{(2D+1)p}  \cdot \|(f,P)\|^p_{\frac{65}{64}\hQ} \left[ 1  +  \sum_{\beta \in \cA}  \| \varphi_\beta \|^p_{\X(3\hQ)} \delta_\hQ^{(m-|\beta|)p-n} \right]  \leq C  \Lambda^{(2D+1)p} \cdot  \|(f,P)\|^p_{\frac{65}{64}\hQ} \\
& \qquad\qquad\qquad (\mbox{by bounded overlap of the collection } \{ (1+a)Q : Q \in \cJ(\hQ)\}, \\
& \qquad\qquad\qquad\qquad\qquad\qquad\qquad\qquad\qquad\qquad\qquad\qquad \mbox{and by \eqref{rrr2}, \eqref{qw2}}).
\end{align*}
This concludes our estimation of $\|H\|_{\X(\frac{65}{64}\hQ)}$.

Writing $H - P = (F-P) + \widetilde{F}$, we also obtain
\begin{align*}
\|H - P \|_{L^p(\frac{65}{64}\hQ)}^p &\lesssim \| F - P \|^p_{L^p(\frac{65}{64}\hQ)} + \sum_{Q \in \cJ(\hQ)} \sum_{\alpha,\beta \in \cA} \| \varphi_\beta \|_{L^p((1+a)Q  \cap \frac{65}{64} \hQ)}^p \cdot \lvert A^Q_{\alpha \beta} \rvert^p \cdot | \partial^\alpha (P - F)(y_Q)|^p \\
& \qquad\qquad (\mbox{see \eqref{tt}} )\\
& \lesssim \| F - P \|^p_{L^p(\frac{65}{64}\hQ)} + \Lambda^{Dp} \cdot \sum_{Q \in \cJ(\hQ)}\sum_{\beta \in \cA} \|\varphi_\beta\|^p_{L^p((1+a)Q  \cap \frac{65}{64} \hQ)} \delta_\hQ^{(m-|\beta|)p - n} \|(f,P)\|^p_{\frac{65}{64}\hQ} \\
& \qquad\qquad  (\mbox{by \eqref{ineq11} and \eqref{a1}})\\
& \lesssim \| F - P \|^p_{L^p(\frac{65}{64}\hQ)} + \Lambda^{Dp} \cdot \sum_{\beta \in \cA} \|\varphi_\beta\|_{L^p(3\hQ)}^p \delta_\hQ^{(m-|\beta|)p - n} \|(f,P)\|^p_{\frac{65}{64}\hQ} \\
&\qquad\qquad   (\mbox{by bounded overlap of the collection} \; \{ (1+a)Q : Q \in \cJ(\hQ) \} ) \\
& \lesssim \Lambda^{(D+1)p} \cdot [\delta_\hQ]^{mp} \|(f,P)\|^p_{\frac{65}{64}\hQ} \qquad\qquad  (\mbox{by \eqref{rrr2} and \eqref{qw5}}).
\end{align*}
Adding together the previous two estimates, we have
$$\|H\|_{\X(\frac{65}{64}\hQ)} + \delta_\hQ^{-m}\|H-P\|_{L^p(\frac{65}{64}\hQ)} \leq C\Lambda^{2D+1} \|(f,P)\|_{\frac{65}{64}\hQ}.$$ This completes the proof of Proposition \ref{prop_A}
\end{proof}

We recall several facts, and set some notation  for the rest of this section.
\begin{itemize}
\item Suppose $Q \in \CZ(\cA^-)$ and $Q \subset (1+100t_G)\hQ$. Then $\frac{65}{64} Q \subset \frac{65}{64}\hQ$ (see Lemma \ref{lem_cover}).
\item Suppose that $Q', Q'' \in \CZ(\cA^-)$, $Q' \subset (1+t_G)\hQ$, $Q' \leftrightarrow Q''$, and $\delta_{Q'} < t_G \cdot \delta_\hQ$. Then $Q'' \subset (1+100t_G)\hQ$ (by good geometry).
\item The sums below are indexed over cubes $Q \in \CZ(\cA^-)$, and over pairs of cubes $(Q',Q'') \in \CZ(\cA^-) \times \CZ(\cA^-)$; and often for ease of notation we choose to not make this indexing explicit.
\item For ease of notation, we write $R_Q^\hQ = R_Q^\hQ(f,P)$ (the dependence on $(f,P)$ should be understood).
\item We are assuming that $t_G$ is less than a small constant determined by $m$, $n$, and $p$, to be picked later, whereas the universal constant $S_1$ has already been picked. (See \eqref{consts}.)
\end{itemize}

\begin{prop} \label{sob_prop} Given $H \in \X$, and given $\{ R_{Q^\#} :  Q^\# \; \mbox{keystone} \} \subset \cP$, the following inequality holds:
\begin{align}
\sum_{\substack{ Q \subset (1+100t_G)\hQ \\  \delta_Q < t_G \cdot \delta_\hQ}} \delta_Q^{-mp} & \| H - R_{\mathcal{K}(Q)} \|^p_{L^p(\frac{65}{64}Q)} \notag{} \\
& \lesssim \sum_{\substack{ Q^\# \; \tiny{\mbox{keystone}}  \\ S_1 Q^\# \subset \frac{65}{64}\hQ}} [\delta_{Q^\#}]^{-mp} \| H - R_{Q^\#} \|^p_{L^{p}(\frac{65}{64} Q^\#)} +  \|H\|_{\X(\frac{65}{64}\hQ)}^p. \label{important}
\end{align}
\end{prop}
\begin{proof}

Let $Q \in \CZ(\cA^-)$ satisfy 
\begin{equation}
\label{cond1}
Q \subset (1+100t_G)\hQ \;\; \mbox{and} \;\; \delta_Q < t_G \cdot \delta_\hQ.
\end{equation}
Then there exists an exponentially decreasing path connecting $Q$ and $\mathcal{K}(Q)$, as promised by the \textsc{Keystone-Oracle}. We denote this path by 
\[ Q = Q(1) \leftrightarrow Q(2) \leftrightarrow \cdots \leftrightarrow Q(L_Q) = \mathcal{K}(Q).\]
Recall that 
\begin{equation}\label{edp1} 
\delta_{Q(\ell')} \leq C \cdot (1-c)^{\ell' - \ell} \cdot \delta_{Q(\ell)} \; \mbox{for} \; \ell' \geq \ell;
\end{equation}
also $Q(\ell) \subset CQ$, and $S_1 \mathcal{K}(Q) \subset C Q$, for a universal constant $C$. From \eqref{cond1} we conclude that $\frac{65}{64}CQ \subset \frac{65}{64}\hQ$, as long as $t_G$ is sufficiently small. Therefore,
\begin{equation}
\label{edp2}
\frac{65}{64}Q(\ell) \subset \frac{65}{64} \hQ \;\; \mbox{for all} \; \ell = 1,\cdots,L_Q, \; \mbox{and} \; S_1\mathcal{K}(Q) \subset \frac{65}{64}\hQ.
\end{equation}
In particular, note that $\frac{65}{64}Q \subset \frac{65}{64} \hQ$.

Fix an arbitrary number $\eta \in (0,1-n/p)$ depending only on $n$ and $p$. By the triangle inequality,
\begin{align*}
 \delta_Q^{-mp}  \|H - R_{\mathcal{K}(Q)}\|_{L^p(\frac{65}{64}Q)}^p \lesssim \;  & \delta_{Q}^{-mp} \|H - J_{x_{Q}} H \|^p_{L^p(\frac{65}{64}Q)}  + \delta_{Q}^{-mp}\|J_{x_{\mathcal{K}(Q)}} H - R_{\mathcal{K}(Q)}\|_{L^p(\frac{65}{64}Q)}^p  \\
&  + \delta_Q ^{-mp} \left\| \sum_{\ell=1}^{L_Q - 1} \left( J_{x_{Q(\ell)} } H -  J_{x_{Q(\ell+1)}} H\right) \delta_{Q(\ell)}^{- \eta}  \delta_{Q(\ell)}^{+\eta} \right\|^p_{L^p(\frac{65}{64}Q)};
\end{align*}
here, H\"older's inequality shows that
\begin{align*}
\delta_Q ^{-mp} & \left\|   \sum_{\ell=1}^{L_Q - 1} \left( J_{x_{Q(\ell)} } H - J_{x_{Q(\ell + 1)} } H\right) \delta_{Q(\ell)}^{- \eta}  \delta_{Q(\ell)}^{+\eta} \right\|^p_{L^p(\frac{65}{64}Q)} \\
& \leq \delta_Q^{- mp}  \left(\sum_{\ell=1}^{L_Q - 1} \delta_{Q(\ell)}^{\eta p'} \right)^{p/p'}  \sum_{\ell=1}^{L_Q - 1} \delta_{Q(\ell)}^{-\eta p} \| J_{x_{Q(\ell)} } H - J_{x_{Q(\ell + 1)} } H \|_{L^p(\frac{65}{64}Q)}^p \\
& \lesssim \delta_Q^{\eta p - m p}  \sum_{\ell=1}^{L_Q - 1} \delta_{Q(\ell)}^{-\eta p} \| J_{x_{Q(\ell)} } H - J_{x_{Q(\ell + 1)} } H \|_{L^p(\frac{65}{64}Q)}^p \quad (\mbox{by \eqref{edp1}});
\end{align*}
also, the Sobolev inequality shows that $\delta_{Q}^{-mp} \|H - J_{x_Q}H\|^p_{L^p(\frac{65}{64}Q)} \leq  C \|H\|^p_{\X(\frac{65}{64}Q)}$. Combining these estimates, we have:
\begin{align*}
 \delta_Q^{-mp} \| H - R_{\mathcal{K}(Q)} \|^p_{L^p(\frac{65}{64}Q)} \lesssim  \;  \|H\|_{\X(\frac{65}{64}Q)}^p & + \delta_Q^{-mp} \| J_{x_{\mathcal{K}(Q)}} H - R_{\mathcal{K}(Q)}\|_{L^p(\frac{65}{64}Q)}^p \\
 & + \delta_Q^{\eta p - mp} \sum_{\ell=1}^{L_Q-1} \delta_{Q(\ell)}^{-\eta p} \| J_{x_{Q(\ell)}} H - J_{x_{Q(\ell+1)}} H \|_{L^p(\frac{65}{64}Q)}^p.
 \end{align*}
Using Lemma \ref{pnorm}, we find that
\begin{align*}
\delta_Q^{-mp} \| H - R_{\mathcal{K}(Q)} \|^p_{L^p(\frac{65}{64}Q)} \lesssim \; & \|H\|_{\X(\frac{65}{64}Q)}^p  + \delta_Q^{-mp } \sum_{|\beta| \leq m-1} \lvert \partial^\beta(J_{x_{\mathcal{K}(Q)}} H - R_{\mathcal{K}(Q)})(x_{\mathcal{K}(Q)}) \rvert^p  \delta_Q^{|\beta|p + n} \\
& + \delta_Q^{\eta p - mp} \sum_{\ell=1}^{L_Q-1} \delta_{Q(\ell)}^{-\eta p} \sum_{|\beta| \leq m-1}  \lvert \partial^\beta(J_{x_{Q(\ell)}} H - J_{x_{Q(\ell+1)}} H)(x_{Q(\ell)}) \rvert^p \delta_Q^{|\beta|p + n}.
\end{align*}

Let $X$ denote the sum of $\delta_Q^{-mp} \| H -  R_{\mathcal{K}(Q)} \|_{L^p(\frac{65}{64}Q)}^p$ over all $Q \in \CZ(\cA^-)$ with $Q \subset (1+100t_G)\hQ$ and $\delta_Q < t_G \delta_\hQ$. 

We now sum the previous estimate over $Q$. We denote $Q^\# = \mathcal{K}(Q)$, $Q' = Q(\ell)$, and $Q'' = Q(\ell+1)$, and we switch the order of summation in our sum. Using \eqref{edp2}, we see that
\begin{align*}
X \lesssim & \sum_{\frac{65}{64}Q \subset \frac{65}{64}\hQ} \| H \|_{\X(\frac{65}{64}Q)}^p \\
& + \sum_{\substack{Q^\# \; \tiny{\mbox{keystone}} \\  S_1 Q^\# \subset \frac{65}{64}\hQ }} \sum_{|\beta| \leq m-1} \lvert \partial^\beta( J_{x_{Q^\#}} H - R_{Q^\#})(x_{Q^\#})\rvert^p \sum_{\substack{ \frac{65}{64}Q \subset \frac{65}{64} \hQ \\ \mathcal{K}(Q) = Q^\#}} \delta_Q^{(|\beta| - m)p + n} \\
& + \sum_{\substack{Q' \leftrightarrow Q'' \\  \frac{65}{64}Q', \frac{65}{64}Q'' \subset \frac{65}{64} \hQ }}\delta_{Q'}^{- \eta p} \sum_{|\beta| \leq m-1}  \lvert \partial^\beta( J_{x_{Q'}} H - J_{x_{Q''}}H)(x_{Q'})\rvert^p \sum_{ \frac{65}{64}Q \subset \frac{65}{64} \hQ} \; \sum_{\substack{\ell : Q(\ell) = Q'} } \delta_Q^{(\eta + |\beta| - m)p + n}.
\end{align*}
Now, for fixed $Q^\#$ we have
\begin{equation*}
\sum_{\substack{ \frac{65}{64}Q \subset \frac{65}{64} \hQ \\ \mathcal{K}(Q) = Q^\#}} \delta_Q^{(|\beta| - m)p + n} \leq \sum_{\substack{ Q \; \tiny \mbox{dyadic} \\  Q^\# \subset CQ}} \delta_Q^{(|\beta| - m)p + n} \leq C \cdot \left[\delta_{Q^\#}\right]^{(|\beta| - m)p + n}.
\end{equation*}
Also, from Remark \ref{key_rem0}, which can be found after the \textsc{Keystone-Oracle}, for fixed $Q'$ we have
\begin{align*}
\sum_{ \frac{65}{64}Q \subset \frac{65}{64} \hQ} \; \sum_{\substack{\ell : Q(\ell) = Q'} } \delta_Q^{(\eta + |\beta| - m)p + n}  & \leq C \sum_{\substack{ \frac{65}{64}Q \subset \frac{65}{64} \hQ \\ \exists \ell, \; Q(\ell) = Q' }} \delta_Q^{(\eta + |\beta| - m)p + n}  \\
&\leq C \sum_{\substack{ Q \; \tiny \mbox{dyadic} \\  Q' \subset CQ}} \delta_Q^{(\eta + |\beta| - m)p + n} \leq C \cdot \left[\delta_{Q'}\right]^{( \eta + |\beta| - m)p + n}.
\end{align*}

Therefore,
\begin{align*}
X \lesssim  \sum_{\frac{65}{64}Q \subset \frac{65}{64}\hQ} \| H \|_{\X(\frac{65}{64}Q)}^p  & + \sum_{\substack{Q^\# \;\tiny{\mbox{keystone}} \\ S_1 Q^\# \subset \frac{65}{64}\hQ }} \sum_{|\beta| \leq m-1} \lvert \partial^\beta( J_{x_{Q^\#}} H - R_{Q^\#})(x_{Q^\#})\rvert^p  \cdot \left[ \delta_{Q^\#}\right]^{(|\beta| - m)p + n} \\
& + \sum_{\substack{ Q' \leftrightarrow Q'' \\ \frac{65}{64}Q', \frac{65}{64}Q'' \subset \frac{65}{64} \hQ}}  \sum_{|\beta| \leq m-1}  \lvert \partial^\beta( J_{x_{Q'}} H - J_{x_{Q''}}H)(x_{Q'})\rvert^p \cdot \left[ \delta_{Q'}\right]^{(|\beta| - m)p + n}.
\end{align*}
Next, for fixed $Q^\#$, Lemma \ref{si2} implies that
\[  \sum_{|\beta| \leq m-1} \lvert \partial^\beta( J_{x_{Q^\#}} H - R_{Q^\#})(x_{Q^\#})\rvert^p  \cdot \left[ \delta_{Q^\#}\right]^{(|\beta| - m)p + n}  \lesssim  \| H  \|_{\X(\frac{65}{64} Q^\#)}^p + \delta_{Q^\#}^{-mp} \| H - R_{Q^\#} \|_{L^p(\frac{65}{64} Q^\#)}^p.\]
Moreover, applying Lemma \ref{si3} with $B_1 = \frac{65}{64}Q'$ and $B_2 = \frac{65}{64}Q''$, where $Q',Q'' \in \CZ(\cA^-)$ and $Q' \leftrightarrow Q''$, we obtain the estimate
\[\sum_{|\beta| \leq m-1}  \lvert \partial^\beta( J_{x_{Q'}} H - J_{x_{Q''}}H)(x_{Q'})\rvert^p \cdot \left[ \delta_{Q'}\right]^{(|\beta| - m)p + n} \lesssim \| H \|_{\X(\frac{65}{64}Q')}^p + \|H\|_{\X(\frac{65}{64}Q'')}^p.\]
(Here, we use the fact that $\lvert x_{Q'} - x_{Q''} \rvert \leq C \delta_{Q'} $ and that $\frac{65}{64}Q' \cap \frac{65}{64}Q'' \neq \emptyset$.)

Thus
\begin{align*}
X \lesssim  \; & \sum_{\frac{65}{64}Q \subset \frac{65}{64}\hQ} \| H \|_{\X(\frac{65}{64}Q)}^p   + \sum_{\substack{ Q' \leftrightarrow Q'' \\ \frac{65}{64}Q', \frac{65}{64}Q'' \subset \frac{65}{64} \hQ}} \biggl[  \| H \|_{\X(\frac{65}{64}Q')}^p +  \| H \|_{\X(\frac{65}{64}Q'')}^p \biggr] \\
& + \sum_{\substack{Q^\# \;\tiny{\mbox{keystone}} \\ S_1 Q^\# \subset \frac{65}{64}\hQ }}   \| H \|_{\X(\frac{65}{64}Q^\#)}^p  + \sum_{\substack{Q^\# \;\tiny{\mbox{keystone}} \\ S_1 Q^\# \subset \frac{65}{64}\hQ }}   \delta_{Q^\#}^{-mp} \| H - R_{Q^\#}\|_{L^p(\frac{65}{64}Q^\#)}^p.
\end{align*}

From the bounded overlap of the cubes $\frac{65}{64}Q$, $Q \in \CZ(\cA^-)$, the previous estimate implies that
\begin{align*}
X \lesssim & \; \| H \|_{\X(\frac{65}{64} \hQ)}^p + \sum_{\substack{Q^\# \;\tiny{\mbox{keystone}} \\ S_1 Q^\# \subset \frac{65}{64}\hQ }}   \delta_{Q^\#}^{-mp} \| H - R_{Q^\#}\|_{L^p(\frac{65}{64}Q^\#)}^p.
\end{align*}
This completes the proof of Proposition \ref{sob_prop}.
\end{proof}

We are now prepared to prove the conditional inequality.

\label{page1}

We seek an estimate on $\bigl[ M_\hQ(f,P) \bigr]^p$, which is the sum of the terms \textbf{(I)}-\textbf{(IV)} (see \eqref{i}-\eqref{iv}). We first apply Lemma \ref{pnorm} to estimate the summands appearing in \eqref{ii}, \eqref{iii}, \eqref{iv}. We also replace \eqref{ii} by a sum over a larger collection of pairs $(Q',Q'')$ (see below). Thus, we obtain
\begin{align}
  \label{e772}
\bigl[ M_\hQ(f,P) \bigr]^p & \leq C(t_G) \cdot \biggl[  \sum_{\substack{ Q \in \CZ_{\main}(\cA^-) \\ Q \subset (1+100t_G)\hQ}} \bigl[ M_{(Q,\cA^-)}(f,R^\hQ_Q) \bigr]^p + \sum_{\substack{Q',Q'' \in \CZ(\cA^-) \\ Q',Q'' \subset (1+100t_G)\hQ \\ Q' \leftrightarrow Q''}} \delta_{Q'}^{-m p}  \|R^\hQ_{Q'} - R^\hQ_{Q''} \|_{L^p(Q')}^p \\ 
& + \sum_{\substack{ Q \in \CZ(\cA^-) \\ Q \subset (1+100t_G)\hQ \\ \delta_Q \geq t_G^2 \cdot \delta_\hQ}} \delta_{\hQ}^{-m p} \|R_Q^\hQ - P \|_{L^p(\hQ)}^p + \delta_{\hQ}^{-m p} \| R^\hQ_{Q_\spec} - P \|_{L^p(\hQ)}^p \biggr]. \notag{} 
\end{align}

We pick a function $H$ as in Proposition \ref{prop_A}. Our estimates proceed in three stages below.

\noindent{\textbf{Stage I:}} We bound the relevant summands in  \eqref{e772}.

We consider $Q,Q',Q'' \in \CZ(\cA^-)$ that satisfy $ Q,Q',Q''  \subset (1+100t_G)\hQ$ and $Q' \leftrightarrow Q''$. We impose either the assumption $Q \in \CZ_{\main}(\cA^-)$ or the assumption $\delta_Q \geq t_G^2 \cdot \delta_\hQ$, depending on which expression we seek to bound.

Assume first that $Q \in \CZ_{\main}(\cA^-)$. Then the right-hand estimate in \eqref{n_appx} implies that
\begin{equation}
\label{e772a}
M_{(Q,\cA^-)}(f,R_Q^\hQ) \leq C \cdot \| (f , R_Q^\hQ) \|_{\frac{65}{64}Q} \leq C \cdot \left[ \|H \|_{\X(\frac{65}{64}Q)} + \delta_Q^{-m} \|H - R_Q^\hQ \|_{L^p(\frac{65}{64}Q)} \right].
\end{equation}
Here, in the last inequality, we use the definition of the trace seminorm and recall that $H = f$ on $E \cap \frac{65}{64}\hQ$.

On the other hand, assume that $\delta_Q \geq t_G^2 \delta_\hQ$. We first apply the triangle inequality and next apply estimate \eqref{s_ineq2} from Lemma \ref{si2} (note that $\frac{65}{64}Q \subset \frac{65}{64}\hQ$, as shown in Lemma \ref{lem_cover}). Thus, we have
\begin{align*}
\delta_{\hQ}^{-m} \|R_{Q}^\hQ - P \|_{L^p(\hQ)} & \leq \delta_\hQ^{-m} \| H - P \|_{L^p(\frac{65}{64}\hQ)} + \delta_\hQ^{-m} \|R^\hQ_Q - H \|_{L^p(\frac{65}{64}\hQ)} \\
& \leq \delta_{\hQ}^{-m} \| H - P \|_{L^p(\frac{65}{64}\hQ)} + C \cdot \left[ \delta_{Q}^{-m}  \|R_{Q}^\hQ - H \|_{L^p(\frac{65}{64}Q)} + \| H \|_{\X(\frac{65}{64}\hQ)} \right].
\end{align*}

We now consider the summands indexed by pairs $(Q',Q'')$. Lemma \ref{si4} implies that
\begin{align*} \delta_{Q'}^{-m} \|R_{Q'}^\hQ - R_{Q''}^\hQ\|_{L^p(Q')}  \leq C \cdot \biggl[  \delta_{Q'}^{-m}  \|R_{Q'}^\hQ - H \|_{L^p(\frac{65}{64}Q')} & + \delta_{Q''}^{-m} \| H - R_{Q''}^\hQ \|_{L^p(\frac{65}{64}Q'')} \\
& + \| H\|_{\X(\frac{65}{64}Q')} + \| H \|_{\X(\frac{65}{64}Q'')} \biggr].
\end{align*}

Furthermore, since $Q_{\spec} \subset \hQ$, we have  $\frac{65}{64}Q_{\spec} \subset \frac{65}{64}\hQ$. Hence, Lemma \ref{si2} implies that
\begin{align*}
\delta_\hQ^{-m} \|P - R_{Q_\spec}^\hQ \|_{L^p(\frac{65}{64}\hQ)} & \leq \delta_\hQ^{-m} \| H - P \|_{L^p(\frac{65}{64}\hQ)} + \delta_\hQ^{-m} \|R^\hQ_{Q_\spec} - H \|_{L^p(\frac{65}{64}\hQ)} \\
&\leq \delta_\hQ^{-m}  \| H - P \|_{L^p(\frac{65}{64}\hQ)} + C \cdot \left[ \delta_{Q_\spec}^{-m} \| R_{Q_\spec}^\hQ - H \|_{L^p(\frac{65}{64}Q_\spec)}  + \|H \|_{\X(\frac{65}{64}\hQ)} \right].
\end{align*}

We combine \eqref{e772} and the previous four estimates to obtain
\begin{align} 
\bigl[ M_{\hQ}(f,P)\bigr]^p \leq C(t_G) \cdot \biggl( \|H \|_{\X(\frac{65}{64}\hQ)}^p & +  \delta_\hQ^{-mp}  \| H - P \|_{L^p(\frac{65}{64}\hQ)}^p  \notag{} \\
& + \sum_{Q \subset (1+100t_G)\hQ } \bigl[ \| H\|_{\X(\frac{65}{64}Q)}^p + \delta_Q^{-mp} \| H - R^\hQ_Q \|_{L^p(\frac{65}{64}Q)}^p \bigr] \biggr).\label{e773}
\end{align}

\noindent{\textbf{Stage II:}} Observe that
$$\sum_{Q \subset (1+100t_G)\hQ }  \|H \|_{\X(\frac{65}{64}Q)}^p \leq C \cdot \|H \|_{\X(\frac{65}{64}\hQ)}^p.$$
Indeed, we have $\frac{65}{64}Q \subset \frac{65}{64}\hQ$ for any cube $Q \in \CZ(\cA^-)$ arising above (see Lemma \ref{lem_cover}); hence, the desired estimate is a consequence of the fact that the cubes $\frac{65}{64}Q$, with $Q \in \CZ(\cA^-)$, have bounded overlap.

The number of cubes $Q \in \CZ(\cA^-)$ such that $Q \subset (1+100t_G)\hQ$ and $\delta_Q \geq t_G \delta_\hQ$ is bounded by a constant $C(t_G)$. Hence,
$$ \sum_{\substack{ Q \subset (1+100t_G)\hQ \\ \delta_Q \geq t_G \delta_\hQ}} \delta_Q^{-mp} \| H - R_Q^\hQ\|_{L^p(\frac{65}{64}Q)}^p  \leq C(t_G) \cdot \delta_\hQ^{-mp} \| H - P \|_{L^p(\frac{65}{64}\hQ)}^p \quad (\mbox{see \eqref{jet1}}).$$
On the other hand,
$$\sum_{\substack{Q \subset (1+100t_G)\hQ \\ \delta_Q < t_G  \delta_\hQ}} \delta_Q^{-mp} \| H - R_Q^\hQ\|_{L^p(\frac{65}{64}Q)}^p  =  \sum_{\substack{ Q \subset (1+100t_G)\hQ \\ \delta_Q < t_G  \delta_\hQ}} \delta_Q^{-mp} \| H - R^\#_{\mathcal{K}(Q)} \|_{L^p(\frac{65}{64}Q)}^p \quad (\mbox{see \eqref{jet1}}).$$

We combine \eqref{e773} and the previous three estimates to obtain
\begin{align} 
\left[M_{\hQ}(f,P)\right]^p & \leq C(t_G) \cdot \biggl(  \|H \|_{\X(\frac{65}{64}\hQ)}^p + \delta_\hQ^{-mp}  \| H - P \|_{L^p(\frac{65}{64}\hQ)}^p  + \sum_{\substack{ Q \subset (1+100t_G) \hQ \\ \delta_Q < t_G \delta_\hQ}} \delta_Q^{-mp} \| H - R^\#_{\mathcal{K}(Q)} \|_{L^p(\frac{65}{64}Q)}^p \biggr) \notag{} \\
& \leq C(t_G) \cdot \biggl(  \|H \|_{\X(\frac{65}{64}\hQ)}^p +  \delta_\hQ^{-mp} \| H - P \|_{L^p(\frac{65}{64}\hQ)}^p  + \sum_{\substack{ Q^\# \; \tiny{\mbox{keystone}} \\ \; S_1 Q^\# \subset \frac{65}{64}\hQ}} (\delta_{Q^\#})^{-mp} \| H - R^\#_{Q^\#} \|^p_{L^{p}(S_1Q^\#)} \biggr) \label{e774} \\
& \quad\qquad\qquad\qquad\qquad\qquad\qquad\qquad\qquad\qquad\quad\; (\mbox{see Proposition \ref{sob_prop}}).\notag{}
\end{align}

\label{stageiii}

\noindent{\textbf{Stage III:}} Let $Q^\# \in \CZ(\cA^-)$ be a keystone cube with $S_1 Q^\# \subset \frac{65}{64}\hQ$. Then, as stated in Proposition \ref{prop_A}, we have $\partial^\alpha H(x_{Q^\#}) = \partial^\alpha P(x_{Q^\#})$ for all $\alpha \in \cA$. Thus, by Proposition \ref{key_prop1}, we have
\begin{equation}
\label{wc3}
(\delta_{Q^\#})^{-mp}  \| H - R^\#_{Q^\#} \|^p_{L^{p}(S_1 Q^\#)} \lesssim \| H \|^p_{\X(S_1 Q^\#)}.
\end{equation}
From Lemma \ref{key_geom}, we recall that the cubes $S_1 Q^\#$ ($Q^\#$ keystone) have bounded overlap. Thus, \eqref{e774} and \eqref{wc3} imply that
\begin{align}
\label{wc4}
M_{\hQ}(f,P)^p & \leq C(t_G) \cdot \bigl( \|H \|_{\X(\frac{65}{64}\hQ)}^p +  \| H - P \|_{L^p(\frac{65}{64}\hQ)}^p \delta_\hQ^{-mp} \bigr) \\
& \leq C(t_G) \cdot \Lambda^{(2D+1)p} \cdot \|(f,P)\|_{\frac{65}{64}\hQ}^p \qquad\qquad\qquad\;\;\; (\mbox{see Proposition \ref{prop_A}}). \notag{}
\end{align}
Recall that $\epsilon^{\overline{\kappa}} \Lambda^{100D} \leq \epsilon^{{\overline{\kappa}}/2}$ and $\overline{\kappa} \leq \kappa_2 \leq 1$ (see \eqref{eq_b}). Hence, $\Lambda^{2D+1} \leq \epsilon^{-{\overline{\kappa}}/2} \leq \epsilon^{-1}$. This shows that
\[M_{\hQ}(f,P) \leq C(t_G) \cdot (1/\epsilon ) \cdot \|(f,P)\|_{\frac{65}{64}\hQ} \] 
This completes the proof of the conditional inequality. This completes the proof of Proposition \ref{prop_normappx}.

\hfill $\qed$

We fix $t_G > 0$, depending only on $m$, $n$, and $p$, small enough so that the above results hold. Since we have fixed the constant $t_G$, all the previous constants of the form $C(t_G)$ or $c(t_G)$ become universal constants $C$ or $c$. In particular, the constant $a_\new=a_\new(t_G)$ from Lemma \ref{lem_cover} depends only on $m$, $n$, and $p$. We set
\begin{equation} \label{newa_defn}
a(\cA) = a_\new.
\end{equation}

Recall the definition of the convex set $\ooline{\sigma}(\hQ)$ in \eqref{s1}.

Just for the moment, let $\epsilon=\epsilon_0$ be a small enough constant depending only on $m$, $n$, and $p$. From Proposition \ref{prop_normappx} we obtain the following result.

\begin{prop} \label{inc_prop}

There exist universal constants $\epsilon_0 > 0$ and $C \geq 1$ such that the following holds.

Let $\hQ$ be a testing cube. Then the following conclusions hold.

\begin{description}
\item[Unconditional Inequality] $\| (f,P) \|_{(1+a(\cA))\hQ} \leq C M_{\hQ}(f,P).$
\item[Conditional Inequality] If $3\hQ$ is tagged with $(\cA,\epsilon_0)$, then $M_{\hQ}(f,P) \leq C \|(f,P) \|_{\frac{65}{64}\hQ}$.
\end{description}

\begin{description}
\item[Unconditional inclusion] $\ooline{\sigma}(\hQ) \subset C \sigma((1+a(\cA))\hQ).$
\item[Conditional inclusion]  If $3\hQ$ is tagged with $(\cA,\epsilon_0)$, then
$\sigma(\frac{65}{64}\hQ) \subset C \ooline{\sigma}(\hQ).$
\end{description}

\end{prop}

Once again, let $\epsilon$ be a small parameter. As usual, we assume that $\epsilon$ is less than a small enough constant depending only on $m$, $n$, and $p$.

\section{Tools to Fill the Gap Between Geometrically Interesting Cubes}
\label{tools_sec}

For the results in this section, the reader may wish to review the definition of testing cubes (see Definition \ref{testing_defn}).

\begin{prop}\label{tool1}
Let $\hQ$ be a testing cube.

If 
$$\left[ \#\left( \frac{65}{64} \hQ \cap E \right) \leq 1 \;\; \mbox{or} \;\; \ooline{\sigma}(\hQ) \;\; \mbox{has an} \;\; (\cA',x_{\hQ},\epsilon,\delta_{\hQ})\mbox{-basis for some} \; \cA' \leq \cA \right]$$
then $(1+ a(\cA)) \hQ$ is tagged with $(\cA,\epsilon^\kappa)$. Otherwise, no cube containing $3 \hQ$ is tagged with $(\cA,\epsilon^{1/\kappa})$. Here, $\kappa$ is a universal constant.
\end{prop}
\begin{proof}
If $\#(\frac{65}{64}\hQ \cap E) \leq 1$, then $(1+a(\cA))\hQ$ is tagged with $(\cA,\epsilon)$.

Suppose $\ooline{\sigma}(\hQ)$ has an $(\cA',x_{\hQ},\epsilon,\delta_{\hQ})$-basis for some $\cA' \leq \cA$. Proposition \ref{inc_prop} implies that $\ooline{\sigma}(\hQ) \subset C \sigma((1+a(\cA))\hQ)$. Thus, $\sigma((1+a(\cA))\hQ)$ has an $(\cA',x_{\hQ},C \epsilon,\delta_{\hQ})\mbox{-basis}$.

Therefore, $(1+a(\cA))\hQ$ is tagged with $(\cA,\epsilon^\kappa)$. Here, we can arrange that $C \epsilon \leq \epsilon^\kappa$ by taking $\epsilon$ sufficiently small.

This proves the first part of Proposition \ref{tool1}.

On the other hand, suppose $Q \supset 3 \hQ$ and suppose $Q$ is tagged with $(\cA, \epsilon^{1/\kappa'})$, for some $\kappa' > 0$ to be picked below. Then $3 \hQ$ is tagged with $(\cA, \epsilon^{\kappa/\kappa'})$, thanks to Lemma \ref{pre_lem5}. Hence, from Proposition \ref{inc_prop} we see that
\begin{equation} \label{dd1} \sigma\left(\frac{65}{64} \hQ\right) \subset C \cdot \ooline{\sigma}(\hQ).
\end{equation}

Recall that $\frac{65}{64} \hQ \subset Q$ and that $Q$ is tagged with $(\cA, \epsilon^{1/\kappa'})$. Thus, Lemma \ref{pre_lem5} shows that $\frac{65}{64}\hQ$ is tagged with $(\cA, \epsilon^{\kappa/\kappa'})$. This means that either $\#(\frac{65}{64}\hQ \cap E) \leq 1$ or $\sigma(\frac{65}{64}\hQ)$ has an $(\cA',x_{\hQ},\epsilon^{\kappa/\kappa'},\delta_\hQ)$-basis, with $\cA' \leq \cA$. Thus, \eqref{dd1} implies that
\[\#\left(\frac{65}{64}\hQ \cap E\right) \leq 1 \;\; \mbox{or} \;\; \ooline{\sigma}(\hQ)  \; \mbox{has an} \; (\cA',x_{\hQ},C \epsilon^{\kappa/\kappa'},\delta_\hQ)\mbox{-basis}.\]
Hence, either $\#(\frac{65}{64}\hQ \cap E) \leq 1$ or $\ooline{\sigma}(\hQ)$ has an $(\cA',x_{\hQ}, \epsilon, \delta_\hQ)$-basis for some $\cA' \leq \cA$. Here, we have determined $\kappa ' = \kappa/2$, with $\kappa$ as in Lemma \ref{pre_lem5}; note that $C \epsilon^{\kappa/\kappa'} = C \epsilon^2 \leq \epsilon$.

This completes the proof of Proposition \ref{tool1}.
\end{proof}

Suppose that $\hQ_1 \subset \hQ_2$ are testing cubes.  We want to understand the tagging of $3 \hQ_2$ in terms of the convex symmetric set $\ooline{\sigma}(\hQ_1)$.

\begin{prop} \label{tag_prop1}
Suppose that $\hQ_1 \subset \hQ_2$ are testing cubes. We assume that 
\begin{equation}
\label{l1}\#(3\hQ_2 \cap E) \geq 2
\end{equation}
and
\begin{equation}
\label{l2}
(1+a(\cA))\hQ_1 \cap E = 3 \hQ_2 \cap E.
\end{equation}

If $\ooline{\sigma}(\hQ_1)$ has an $(\cA',x_{\hQ_1},\epsilon,\delta_{\hQ_2})$-basis, then $3 \hQ_2$ is tagged with $(\cA',\epsilon^\kappa)$ for a universal constant $\kappa$.
\end{prop}

\begin{proof}

Let $(P_\alpha)_{\alpha \in \cA'}$ be an $(\cA',x_{\hQ_1},\epsilon,\delta_{\hQ_2})$-basis for $\ooline{\sigma}(\hQ_1)$. Thus,
\begin{align}
\label{tag1} & P_\alpha \in \epsilon \delta_{\hQ_2}^{-(m - \frac{n}{p} - |\alpha|)} \ooline{\sigma}(\hQ_1) \qquad (\alpha \in \cA') \\
\label{tag2} & \partial^\beta P_\alpha(x_{\hQ_1}) = \delta_{\beta \alpha} \qquad\qquad\quad (\beta,\alpha \in \cA')\\
\label{tag3} & \lvert \partial^\beta P_\alpha(x_{\hQ_1}) \rvert \leq \epsilon \delta_{\hQ_2}^{|\alpha| - |\beta|} \qquad\quad (\alpha \in \cA', \; \beta \in \cM, \; \beta > \alpha).
\end{align}
The unconditional inclusion and \eqref{t1} show that
\begin{equation*}
\ooline{\sigma}(\hQ_1) \subset C \sigma((1+a(\cA)) \hQ_1) \subset C \left[ \sigma((1+a(\cA)) \hQ_1) + \cB(x_{\hQ_1}, 3 \delta_{\hQ_2}) \right] \subset C' \sigma(3 \hQ_2).
\end{equation*}
(We can apply \eqref{t1} from Lemma  \ref{pre_lem0}, since we assume here that $(1+a(\cA)) \hQ_1 \cap E = 3 \hQ_2 \cap E$.)

Thus, \eqref{tag1} implies that
\[
P_\alpha \in C \epsilon \delta_{3 \hQ_2}^{-(m - \frac{n}{p} - |\alpha|)} \cdot \sigma(3 \hQ_2) \qquad (\alpha \in \cA').
\]

Together with \eqref{tag2} and \eqref{tag3}, this shows that $(P_\alpha)_{\alpha \in \cA'}$ is an $(\cA',x_{\hQ_1}, C \epsilon, \delta_{3 \hQ_2})$-basis for $\sigma(3 \hQ_2)$.

It follows from Lemma \ref{pre_lem4} that $\sigma(3 \hQ_2)$ has an $(\cA'', x_{\hQ_2}, \epsilon^\kappa, \delta_{3 \hQ_2})$-basis, for some $\cA'' \leq \cA'$. By definition, this means that the cube $3 \hQ_2$ is tagged with $(\cA',\epsilon^\kappa)$, completing the proof of Proposition \ref{tag_prop1}.
\end{proof}

\begin{cor}
\label{cor_tag1}
Suppose that $\hQ_1 \subset \hQ_2$ are testing cubes and that \eqref{l1}, \eqref{l2} hold.

Suppose $\ooline{\sigma}(\hQ_1)$ has an $(\cA',x_{\hQ_1},\epsilon,\delta_{\hQ_2})$-basis for some $\cA' \leq \cA$. Then $3 \hQ_2$ is tagged with $(\cA,\epsilon^\kappa)$ for a universal constant $\kappa$.
\end{cor}

\begin{proof}
Proposition \ref{tag_prop1} tells us that $3 \hQ_2$ is tagged with $(\cA',\epsilon^\kappa)$. This trivially implies that $3 \hQ_2$ is tagged with $(\cA,\epsilon^\kappa)$.
\end{proof}

\begin{prop}\label{tag_prop2}
Suppose that $\hQ_1 \subset \hQ_2$ are testing cubes and that \eqref{l1}, \eqref{l2} hold.

Suppose  $3 \hQ_2$ is tagged with $(\cA,\epsilon)$. Then $\ooline{\sigma}(\hQ_1)$ has an $(\cA',x_{\hQ_1},\epsilon^{\kappa'},\delta_{\hQ_2})$-basis, for some $\cA' \leq \cA$. Here, $\kappa'$ is a universal constant.
\end{prop}

\begin{proof}
We have $3 \hQ_1 \subset 3 \hQ_2$, so Lemma \ref{pre_lem5} tells us that $3 \hQ_1$ is tagged with $(\cA,\epsilon^\kappa)$. Hence, by the conditional inclusion, we have
\begin{equation}\label{tag4}
c\cdot \sigma\left(\frac{65}{64} \hQ_1\right) \subset \ooline{\sigma}(\hQ_1).
\end{equation}
Next note that $\frac{65}{64}\hQ_1 \cap E = 3 \hQ_2 \cap E$, and that $\frac{65}{64}\hQ_1 \subset 3 \hQ_2$. Therefore, Lemma \ref{pre_lem0} gives the inclusion
\begin{equation} \label{tag5}
\sigma(3 \hQ_2) \subset C \cdot \left[ \sigma\left(\frac{65}{64}\hQ_1\right) + \cB(x_{\hQ_2},\delta_{3 \hQ_2})  \right] .
\end{equation}
(Since $\lvert x_{\hQ_1} - x_{\hQ_2} \rvert \leq \delta_{\hQ_2}$, it follows that $\cB(x_{\hQ_1}, \delta_{3 \hQ_2}) \subset C \cB(x_{\hQ_2},\delta_{3 \hQ_2})$. This shows that \eqref{tag5} follows from the conclusion of Lemma \ref{pre_lem0}.)

Now, $3 \hQ_2$ is assumed to be tagged with $(\cA,\epsilon)$, and $\#(3 \hQ_2 \cap E)$ is assumed to be at least $2$. Hence, by definition, $\sigma(3 \hQ_2)$ has an $(\cA',x_{\hQ_2},\epsilon,\delta_{3 \hQ_2})$-basis for some $\cA' \leq \cA$.

By Lemma \ref{pre_lem2},
\begin{equation}\label{tag6}
\sigma(3 \hQ_2) \; \mbox{has an} \; (\cA'',x_{\hQ_2}, \epsilon^\kappa, \delta_{3 \hQ_2}, \Lambda)\mbox{-basis, for some } \cA'' \leq \cA' \leq \cA,
\end{equation}
such that $\epsilon^\kappa \Lambda^{100D} \leq \epsilon^{\kappa/2}$ and $\kappa \in [\kappa_1, \kappa_2]$. Here, $\kappa_1, \kappa_2 > 0$ are universal constants.

Inclusions \eqref{tag4} and \eqref{tag5} show that
\begin{equation} \label{tag7}
\sigma( 3 \hQ_2) \subset C'' \cdot \left[ \ooline{\sigma}(\hQ_1) + \cB(x_{\hQ_2}, \delta_{3 \hQ_2}) \right].
\end{equation}
From \eqref{tag6}, \eqref{tag7}, and Lemma \ref{pre_lem1}, we see that
\[
\ooline{\sigma}(\hQ_1) \; \mbox{has an} \; (\cA'', x_{\hQ_2}, C\epsilon^\kappa \Lambda, \delta_{3 \hQ_2}, C \Lambda)\mbox{-basis}.
\]
We now apply Lemma \ref{pre_lem3}. Thus,
\[
\ooline{\sigma}(\hQ_1) \; \mbox{has an} \; (\cA'', x_{\hQ_1}, C\epsilon^\kappa \Lambda^{2D+2}, \delta_{3 \hQ_2}, C \Lambda^{2D+1})\mbox{-basis}.
\]
Since $C \epsilon^{\kappa} \Lambda^{2D+2} \leq \epsilon^{\kappa/3} \leq \epsilon^{\kappa_1/3}$, it follows that
\begin{equation} \label{tag8}
\ooline{\sigma}(\hQ_1) \; \mbox{has an} \; (\cA'',x_{\hQ_1},\epsilon^{\kappa_1/3},\delta_{\hQ_2})\mbox{-basis}.
\end{equation}
(Here, the passage from $\delta_{3 \hQ_2}$ to $\delta_{\hQ_2}$ is harmless; it just increases the constant ``$C$'' in $C \epsilon^{\kappa} \Lambda^{2D+2}$.)

Since $\cA'' \leq \cA$, \eqref{tag8} is the conclusion of Proposition \ref{tag_prop2}.
\end{proof}

Combining the results of Propositions \ref{tag_prop1}, \ref{tag_prop2}, we now prove the following.

\begin{prop} \label{tag_prop3} Suppose that $\hQ_1 \subset \hQ_2$ are testing cubes and that \eqref{l1}, \eqref{l2} hold. Then
\begin{enumerate}[(A)]
\item If $\ooline{\sigma}(\hQ_1)$ has an $(\cA',x_{\hQ_1},\epsilon,\delta_{\hQ_2})$-basis for some $\cA' \leq \cA$, then $(1+a(\cA)) \hQ_2$ is tagged with $(\cA,\epsilon^\kappa)$.
\item If some cube containing $3 \hQ_2$ is tagged with $(\cA,\epsilon)$, then $\ooline{\sigma}(\hQ_1)$ has an $(\cA',x_{\hQ_1},\epsilon^\kappa, \delta_{\hQ_2})$-basis for some $\cA' \leq \cA$.
\end{enumerate}
Here, $\kappa$ is a universal constant.
\end{prop}

\begin{proof}
First we check $(A)$. If $\ooline{\sigma}(\hQ_1)$ has an $(\cA',x_{\hQ_1},\epsilon,\delta_{\hQ_2})$-basis with $\cA' \leq \cA$, then according to Corollary \ref{cor_tag1}, the  cube $3 \hQ_2$ is tagged with $(\cA,\epsilon^\kappa)$. Hence, by Lemma \ref{pre_lem5}, $(1+a(\cA))\hQ_2$ is tagged with $(\cA,\epsilon^{\kappa'})$, completing the proof of (A).

To check (B), let $Q' \supset 3 \hQ_2$ be tagged with $(\cA,\epsilon)$. By Lemma \ref{pre_lem5}, $3 \hQ_2$ is tagged with $(\cA,\epsilon^\kappa)$. Hence, by Proposition \ref{tag_prop2}, $\ooline{\sigma}(\hQ_1)$ has an $(\cA',x_{\hQ_1},\epsilon^{\kappa'}, \delta_{\hQ_2})$-basis, for some $\cA' \leq \cA$. 

This completes the proof of (B). 
\end{proof}

We apply (A) with $\epsilon$ unchanged, and (B) with $\epsilon$ replaced by $\epsilon^{1/\kappa}$. Thus we obtain the following result.

\begin{prop}\label{tool2}
Let $\hQ \subset Q$ be testing cubes. 

Assume that $\#(3 Q \cap E) \geq 2$ and that $(1+a(\cA))\hQ \cap E = 3 Q \cap E$.

Then the following hold, for a universal constant $\kappa$.
\begin{enumerate}[(A)]
\item If $\ooline{\sigma}(\hQ)$ has an $(\cA', x_{\hQ},\epsilon,\delta_Q)$-basis for some $\cA' \leq \cA$, then $(1+a(\cA))Q$ is tagged with $(\cA,\epsilon^\kappa)$.
\item If $\ooline{\sigma}(\hQ)$ \underline{does not} have an $(\cA', x_{\hQ},\epsilon,\delta_Q)$-basis for any $\cA' \leq \cA$, then no cube containing $3Q$ is tagged with $(\cA,\epsilon^{1/\kappa})$.
\end{enumerate}
\end{prop}

The final result in this section is the following algorithm.

\environmentA{Algorithm: Optimize Basis.}

We perform one time work at most $C N \log N$ in space $C N$, after which we can answer queries as follows.

A query consists of a testing cube $\hQ$ and a set $\cA \subset \cM$

The response to the query $(\hQ,\cA)$ consists of a collection of pairwise disjoint intervals $I_\ell$ and numbers $a_\ell, \lambda_\ell$ ($\ell=1,\cdots,\ell_{\max}$), such that the following conditions hold.

\begin{itemize}
\item $\bigcup_\ell I_\ell = (0,\infty)$ and $\ell_{\max} \leq C$.
\item Let $\eta^{(\hQ,\cA)}(\delta) := a_\ell \delta^{\lambda_\ell}$ for $\delta \in I_\ell$. Then we have:
\begin{description}
\item[(A1)] For each $\delta \in (0,\infty)$ there exists $\cA' \leq \cA$ such that $\ooline{\sigma}(\hQ)$ has an $(\cA',x_{\hQ},\eta^{1/2},\delta)$-basis for all $\eta > C \cdot \eta^{(\hQ,\cA)}(\delta)$.

\item[(A2)] For each $\delta \in (0,\infty)$ and any $\cA' \leq \cA$, $\ooline{\sigma}(\hQ)$ \underline{does not} have an $(\cA',x_{\hQ}, \eta^{1/2} ,\delta)$-basis with $ \eta < c \cdot \eta^{(\hQ,\cA)}(\delta)$.
\item[(A3)] $c \cdot \eta^{(\hQ,\cA)}(\delta_1) \leq \eta^{(\hQ,\cA)}(\delta_2) \leq C \cdot \eta^{(\hQ,\cA)}(\delta_1)$ whenever $\frac{1}{10} \delta_1 \leq \delta_2 \leq 10 \delta_1$.

\end{description}
\item To answer a query requires work at most $C \log N$.
\end{itemize}

\begin{proof}[\underline{Explanation}] 

We compute a quadratic form $q_{\hQ}$ on $\cP$ such that there exist universal constants $c > 0$ and $C \geq 1$ so that $\{ q_{\hQ} \leq c \} \subset \ooline{\sigma}(\hQ) \subset \{ q_{\hQ} \leq C \}$. (See the algorithm \textsc{Approximate New Trace Norm} in Section \ref{themainresult}.)

Processing the quadratic form $q^\hQ$ using  \textsc{Fit Basis to Convex Body} (see Section \ref{sec_compbase}), we compute  a piecewise-monomial function $\eta_{*}^{(\hQ,\cA')}(\cdot)$ for each $\cA' \leq \cA$. We guarantee that $\ooline{\sigma}(\hQ)$ has an $(\cA',x_\hQ,\eta^{1/2},\delta)$-basis for all $\eta > C \cdot \eta_{*}^{(\hQ,\cA')}(\delta)$, but that $\ooline{\sigma}(\hQ)$ does not have an $(\cA',x_\hQ,\eta^{1/2},\delta)$-basis for any $\eta < c\cdot\eta_{*}^{(\hQ,\cA')}(\delta)$.

We define
\[\eta^{(\hQ,\cA)}(\delta) =  \eta(\delta) = \min_{\cA' \leq \cA} \eta_*^{(\hQ,\cA')}(\delta) \quad \mbox{for} \; \delta \in (0,\infty).\]
It follows that $\ooline{\sigma}(\hQ)$ has an $(\cA',x_\hQ,\eta^{1/2},\delta)$-basis for some $\cA' \leq \cA$ whenever $\eta > C \cdot \eta(\delta)$, but that $\ooline{\sigma}(\hQ)$ does not have an $(\cA',x_\hQ,\eta^{1/2},\delta)$-basis for any $\cA' \leq \cA$ whenever $\eta < c \cdot \eta(\delta)$. Thus we have proven \textbf{(A1)} and \textbf{(A2)}.

Recall that $c \cdot \eta^{(\hQ,\cA')}(\delta_1) \leq  \eta^{(\hQ,\cA')}(\delta_2) \leq C \cdot \eta^{(\hQ,\cA')}(\delta_1)$ for $\frac{1}{10} \delta_1 \leq \delta_2 \leq 10 \delta_1$. Taking the minimum with respect to $\cA' \leq \cA$ in this inequality, we prove \textbf{(A3)}.

Recall that $\eta_*^{(\hQ,\cA')}(\delta) = a_{\ell,\cA'} \delta^{\lambda_{\ell,\cA'}}$ for $\delta \in I_{\ell,\cA'}$, where the intervals $I_{\ell, \cA' }$ ($\ell = 1,\cdots, \ell_{\max}(\cA')$) form a partition of $(0,\infty)$, for each $\cA' \leq \cA$. Here, $\ell_{\max}(\cA')$ is bounded by a universal constant.

Thus we can partition $(0,\infty)$ into intervals $I_\ell$ ($\ell = 1,\cdots,\underline{\ell_{\max}}$), for which there exist real numbers $a_\ell, \lambda_\ell$ such that $\eta(\delta) = a_\ell \delta^{\lambda_\ell}$ for $\delta \in I_\ell$. Moreover, $\underline{\ell_{\max}}$ is at most some universal constant. This follows because, for fixed real numbers $a,b,\lambda,\gamma$, the equation $a \delta^\lambda = b \delta^{\gamma}$ is satisfied either for at most one $\delta$ or for all $\delta \in (0,\infty)$. To compute the intervals $I_\ell$ and the numbers $a_\ell,\lambda_\ell$ we solve at most $C$ equations of the above type, and we make at most $C$ comparisons between the functions $\eta^{(\hQ,\cA')}(\delta)$ ($\cA' \leq \cA$) to compute the minimum value on each of the relevant intervals. This completes the explanation of our algorithm.

\end{proof}

\section{Computing Lengthscales}
\label{sec_cl}

We say that a dyadic cube $Q \subset \R^n$ is \underline{geometrically interesting} provided that $\diam(3 Q \cap E) \geq \lambda \delta_Q$, where we set $\lambda := 1/40$.

\environmentA{Algorithm: Compute Interesting Cubes.}
We produce a tree $T$ consisting of all the cubes $Q \in \CZ(\cA^-)$ that contain points of $E$, together with all testing cubes $\hQ$ for which $\diam(3 \hQ \cap E) \geq \lambda \delta_{\hQ}$; as well as the unit cube $Q^\circ$.

Here, $T$ is a tree with respect to inclusion. We mark each internal node $Q \in T $ with pointers to its children, and we mark each node $Q \in T$ (except for the root) with a pointer to its parent.

The number of nodes in ${T}$ is at most $C N$, and ${T}$ can be computed with work at most $C N \log N$ in space $C N$.

We note that all the nodes of $T$ are testing cubes. (This is immediate from the definition of testing cubes - see Definition \ref{testing_defn}.)
\begin{proof}[\underline{Explanation}]  

We perform the one-time work of the BBD Tree (see Theorem \ref{bbd_thm}). Also, we compute representatives arising in the well-separated pairs decomposition using the algorithm \textsc{Make WSPD} from Section \ref{sec_CK}. Thus, we compute a sequence of tuples $(x_\nu',x_\nu'') \in E \times E$ ($\nu = 1,\ldots,\nu_{\max}$) such that, for each $(x',x'') \in E \times E \setminus \{ (x,x) : x \in E\}$ there exists $\nu$ such that
\[\lvert x_\nu' - x' \rvert + \lvert x_{\nu}'' - x'' \rvert \leq 10^{-10} \lvert x' - x'' \rvert, \]
and $\nu_{\max} \leq C N$.

We execute the following loop:
\begin{itemize}
\item For each $\nu = 1,\cdots,\nu_{\max}$, we compute the sequence of all dyadic cubes $\widetilde{Q}$ such that $x_\nu', x_\nu'' \in 5\widetilde{Q}$ and $\lvert x_\nu' - x_\nu'' \rvert \geq  \frac{\lambda}{2} \delta_{\widetilde{Q}}$. (There are at most $C$ such cubes for each $\nu$.)
\end{itemize}
We denote the sequence of all cubes produced above, for all $\nu$, by $Q_1,\cdots,Q_K$. We remove duplicates by sorting, which requires work at most $C N \log N$. Note that we have $K \leq C N$. 

Let $Q$ be a geometrically interesting cube. By definition, there exist $x',x'' \in 3Q \cap E$ with $\lvert x' - x'' \rvert \geq \lambda \delta_Q$. Hence, there is some $\nu$ such that
\[ \lvert x_\nu' - x_\nu'' \rvert \geq \frac{9}{10} \lvert x' - x'' \rvert \geq (\lambda/2) \delta_Q
\]
and
\[
\lvert x_\nu' - x' \rvert + \lvert x_\nu'' - x'' \rvert \leq \frac{1}{10} \lvert x' - x'' \rvert \leq  \frac{\delta_{3Q}}{10}.
\] 
Therefore, $x_\nu',x_\nu'' \in 5Q$, and hence $Q$ belongs to the list $Q_1,\cdots,Q_K$. 

We have proven that all geometrically interesting cubes belong to the list $Q_1,\cdots, Q_K$.

For each $k=1,\cdots,K$, we compute $\diam(3Q_k \cap E)$ using the BBD Tree. (See Remark \ref{bbd_rem}.) If $\diam(3Q_k \cap E) < \lambda \delta_{Q_k}$, then we remove $Q_k$ from our list. We also compute the cube in $\CZ(\cA^-)$ that contains the center of $Q_k$. If this cube strictly contains $Q_k$ then we remove $Q_k$ from our list. (This means that $Q_k$ is not a testing cube.)

We denote the sequence of surviving cubes by $\widetilde{Q}_1,\cdots,\widetilde{Q}_{\widetilde{K}}$. As shown above, these are all the testing cubes that are geometrically interesting.

We form a list of all the cubes $Q \in \CZ(\cA^-)$ that contain points of $E$, the cubes $\widetilde{Q}_1,\cdots,\widetilde{Q}_{\widetilde{K}}$, and the unit cube $Q^\circ$. There are at most $C N$ such cubes. By sorting, we can remove duplicates. We organize this list into a tree $T$ using the algorithm \textsc{Make Forest} (see Section \ref{maketree_sec}). We obtain a tree (rather than a forest) because all the cubes have a common ancestor, namely $Q^\circ$. This algorithm marks $Q^\circ$ as the root of $T$, and marks each non-root node with a pointer to its parent. In addition, we mark each internal node of $T$ with pointers to its children.

One can easily check that the work and storage of our algorithm are as promised.
\end{proof}

\begin{lem} \label{ff_lem}
Let $Q \subset Q^\circ$ be dyadic, with $\delta_Q \leq \frac{1}{4}$. Suppose that $3Q \cap E \neq \emptyset$ and $\diam(3Q^{++} \cap E) < \lambda \delta_{Q^{++}}$. 

Then $3Q^{++} \cap E = 3Q^+ \cap E$.  Here, $Q^{++}$ denotes the dyadic parent of the dyadic parent of $Q$.
\end{lem}
\begin{proof}
For the sake of contradiction, suppose that there exists $x \in E$ with $x \in 3Q^{++}$ and $x \notin  3Q^+$. Thus, for each $y \in E \cap 3Q$ we have
\[ \diam(3Q^{++} \cap E) \geq \lvert x - y \rvert \geq \dist(\R^n \setminus 3Q^+, 3Q) \geq \frac{\delta_{Q}}{10} = \frac{\delta_{Q^{++}}}{40} = \lambda \delta_{Q^{++}}.\]
This yields a contradiction, completing the proof of the lemma.
\end{proof}

\begin{lem}\label{gg_lem}
Let $Q_1 \subset Q_2$ be dyadic cubes such that $Q_2$ is the parent of $Q_1$ in the tree $T$. Let $a > 0$ be given. Let $Q_1^\up$ and $Q_2^{\text{down}}$ be dyadic cubes.

Assume that $Q_1 \subsetneq Q_1^\up \subsetneq Q_2^{\text{down}} \subsetneq Q_2$ with $\delta_{Q_2} \geq \Lambda \delta_{ Q_2^{\text{down}} } \geq \Lambda^2 \delta_{Q_{1}^\up} \geq \Lambda^3 \delta_{Q_1}$ for some $\Lambda \geq 2$.

If $\Lambda$ exceeds a large enough constant determined by $a$ and $n$, then $(1+a) Q_1^\up \cap E = 3 Q_2^{\text{down}}  \cap E$.
\end{lem}
\begin{proof}

If $\Lambda \geq 4$, then since $Q_1 \subset Q_2^{\text{down}}$ and $\delta_{Q_1} \leq \frac{1}{\Lambda^2} \delta_{Q_2^{\text{down}}}$, we have $Q_1^{+++} \subset Q_2^{\text{down}}$.

Fix a sequence of dyadic cubes $Q_{1,1} \subset Q_{1,2} \subset \cdots \subset Q_{1,K}$ with
\[Q_{1,1} = (Q_1)^{+++}, \; Q_{1,K} = Q_2^{\text{down}}, \;\; \mbox{and} \;\; Q_{1,k} = ( Q_{1,k-1})^+ \; \mbox{for} \; 2 \leq k \leq K.
\]

Since $Q_1$ is a testing cube (recall that all the nodes of $T$ are testing cubes), it follows by definition that $Q_1$ contains a cube in $\CZ(\cA^-)$. Thus, thanks to \eqref{Enearby2}, the set $9 Q_1 \cap E$ is nonempty. We have $3 Q_{1,k} \supset 3Q_{1,1} = 3 Q_1^{+++} \supset 9 Q_1$ for any $1 \leq k \leq K$. Hence, $3 Q_{1,k} \cap E \neq \emptyset$ for any  $1 \leq k \leq K$. Moreover, note that $Q_1 \subsetneq Q_{1,k} \subsetneq Q_2$, because $Q_1^{+++}$ and $Q_2^{\text{down}}$ are strictly contained between $Q_1$ and $Q_2$. Since $Q_2$ is the parent of $Q_1$ in the tree $T$, which contains all the geometrically interesting testing cubes, we learn that $Q_{1,k}$ is \underline{not} geometrically interesting, for each $1 \leq k \leq K$. In particular, we see that $Q_{1,k}^{++} = Q_{1,k+2}$ is not geometrically interesting, hence $\diam(Q_{1,k}^{++} \cap E) < \lambda \delta_{Q_{1,k}^{++}}$, for all $1 \leq k \leq K-2$. Thus, the hypotheses of Lemma \ref{ff_lem} are satisfied by $Q = Q_{1,k}$ for each $1 \leq k \leq K - 2$. We conclude that
\[
3Q_{1,2} \cap E = 3 Q_{1,3} \cap E = \cdots = 3 Q_{1,K} \cap E.
\]
That is, $3 Q_1^{++++} \cap E = 3 Q_2^{\text{down}} \cap E$.

Recall that $Q_1 \subset Q_1^{\up}$ are dyadic cubes with $\delta_{Q_1^\up} \geq \Lambda \delta_{Q_1}$. It follows that $3Q_1^{++++} \subset (1+a)Q_1^\up$ if $\Lambda$ is much larger than $a^{-1}$. Therefore, $3Q_2^{\text{down}} \cap E \subset (1+a) Q_1^\up \cap E$.  Moreover, the reverse inclusion follows because $Q_1^\up \subset Q_2^{\text{down}}$. Therefore, $3 Q_2^{\text{down}} \cap E = (1+a)\upQ_1 \cap E$.
\end{proof}

\subsection{Finding Enough Tagged Cubes}
\label{sec_fetc}

We produce the following algorithm.

\environmentA{Algorithm: Compute Critical Testing Cubes.}

Given $\epsilon > 0$ less than a small enough universal constant, we produce a collection $\widehat{\mathcal{Q}}_\epsilon$ of testing cubes with the following properties.
\begin{enumerate}[(a)]
\item Each point $x \in E$ belongs to some cube $\hQ_x \in \widehat{\mathcal{Q}}_\epsilon$.
\item The number of cubes belonging to $\widehat{\mathcal{Q}}_\epsilon$ is bounded by $C \cdot N$.
\item If $\widehat{Q} \in \widehat{\mathcal{Q}}_\epsilon$ strictly contains a cube in $\CZ(\cA^-)$, then $(1+a(\cA))\widehat{Q}$ is tagged with $(\cA, \epsilon^\kappa)$.
\item If $\widehat{Q} \in \widehat{\mathcal{Q}}_\epsilon$ and $\delta_{\widehat{Q}} \leq c^*$, then no cube containing $S \widehat{Q}$ is tagged with $(\cA, \epsilon^{1/\kappa})$.
\end{enumerate}
Here, $c^* > 0$ and $S \geq 1$ are integer powers of $2$, depending only on $m$, $n$, $p$; also, $\kappa \in (0,1)$ is a universal constant. The algorithm requires work at most $C N \log N$ in space $CN$.

\begin{proof}[\underline{Explanation}]  

We introduce a large parameter $\Lambda = 2^{\scriptsize\mbox{integer}} \geq 1$. We later pick $\Lambda$ to be a constant determined by $m$,$n$,$p$, but not yet. We assume that  $\Lambda$ exceeds a large enough constant determined by $m$,$n$,$p$, and that $\epsilon$ is less than a small enough constant determined by $\Lambda$,$m$,$n$,$p$.

We let $\kappa_0, \cdots, \kappa_{20} \in (0,1)$ be constants to be determined later. We assume that $\kappa_0$ is less than a small enough constant determined by $m$,$n$,$p$, and that $\kappa_{j+1} \leq \kappa_j^{100}$ for $j=0,\cdots, 19$.

We first describe the construction of $\widehat{\mathcal{Q}}_\epsilon$.

Let $T$ be the tree constructed in the algorithm \textsc{Compute Interesting Cubes}.

We initialize $\widehat{\mathcal{Q}}_\epsilon$ to be the empty collection. Next, for each cube $Q_1 \in {T}$ other than the root, we perform Steps 0-3 below.

\begin{itemize}

\item \underline{\textsf{Step 0}}: We find the parent $Q_2$ of $Q_1$ in the tree ${T}$.

\item \underline{\textsf{Step 1}}:  If $\delta_{Q_1} \leq \Lambda^{-20} \delta_{Q_2}$, then we do the following. 

Let $\upQ_1$ be the dyadic cube with $Q_1 \subset \upQ_1$ and $\delta_{\upQ_1} = \Lambda \cdot \delta_{Q_1}$.

\label{pp20}
We compute the function $\eta^{(\hQ^\up_1,\cA)}(\delta)$ using the algorithm \textsc{Optimize Basis} (see Section \ref{tools_sec}). We determine whether or not there exists a number $\delta \in [ \Lambda^{10} \delta_{Q_1}, \Lambda^{-10} \delta_{Q_2}]$ with the property that 
\[ \epsilon^{1/\kappa_5} \leq \eta^{(\upQ_1,\cA)}(\delta) \leq \epsilon^{\kappa_5}.\]
If such a $\delta$ exists, we can easily find one. Moreover, we can then find a dyadic cube $Q$ such that 
\[ Q_1 \subset Q \subset Q_2, \; \delta/2 \leq  \delta_Q \leq 2 \delta, \; \mbox{and} \; \Lambda^{10} \delta_{Q_1} \leq \delta_Q \leq \Lambda^{-10} \delta_{Q_2}.
\]
We add $Q$ to the collection $\widehat{\mathcal{Q}}_\epsilon$.
 Note that $c \eta^{(\upQ_1,\cA)}(\delta) \leq \eta^{(\upQ_1,\cA)}( \delta_Q) \leq C \eta^{(\upQ_1,\cA)}(\delta)$, thanks to condition \textbf{(A3)} in the algorithm \textsc{Optimize Basis}. Thus, we can guarantee that
\begin{equation} \label{step1}\left[ \epsilon^{1/\kappa_6} \leq \eta^{(\upQ_1,\cA)}( \delta_Q)  \right]\;\; \mbox{and} \;\; \left[ \eta^{(\upQ_1,\cA)}(\delta_Q) \leq \epsilon^{\kappa_6} \right].
\end{equation}

\item \underline{\textsf{Step 2}}:  We examine each dyadic cube $Q$ with $Q_1 \subset Q \subset Q_2$, $\delta_Q \leq \Lambda^{-10}$, and \\
$\left[ \delta_Q \leq \Lambda^{10} \delta_{Q_1} \; \mbox{or} \; \delta_Q \geq \Lambda^{-10} \delta_{Q_2} \right]$.

We can compute $\#(E \cap \frac{65}{64}Q)$ using work at most $C \log N$; see Remark \ref{bbd_rem}.

Let $Q^\up$ be the dyadic cube with $Q \subset Q^\up$ and $\delta_{Q^\up} = \Lambda \delta_Q$. We determine whether or not
\begin{equation} \label{step2}
\left[ \epsilon^{1/\kappa_5} \leq  \eta^{(Q^\up,\cA)}(\delta_{Q^\up}) \right] \; \mbox{and} \; \left[ \#\left( \frac{65}{64} Q \cap E\right) \leq 1 \;\; \mbox{or} \;\; \eta^{(Q,\cA)}(\delta_Q) \leq \epsilon^{\kappa_5} \right].
\end{equation}

We add $Q$ to the collection $\widehat{\mathcal{Q}}_\epsilon$ if and only if \eqref{step2} holds.

\item \underline{\textsf{Step 3}}: We examine each dyadic cube $Q$ with $Q_1 \subset Q \subset Q_2$ and $\delta_Q \geq \Lambda^{-10}$.

We can compute $\#(E \cap \frac{65}{64}Q)$ using work at most $C \log N$; see Remark \ref{bbd_rem}.

For each such $Q$, we determine whether or not
\begin{equation} \label{step3} \left[ \#\left(\frac{65}{64}Q \cap E\right) \leq 1 \;\; \mbox{or} \;\; \eta^{(Q,\cA)}(\delta_Q) \leq \epsilon^{\kappa_5} \right].
\end{equation}
We add $Q$ to the collection $\widehat{\mathcal{Q}}_\epsilon$ if and only if \eqref{step3} holds.

\end{itemize}

Finally, we perform Steps 4-6 below.

\begin{itemize}
\item \underline{\textsf{Step 4}}:  We check whether or not
\begin{equation} \label{step4} \left[ \eta^{(Q^\circ,\cA)}(\delta_{Q^\circ}) \leq \epsilon^{\kappa_5} \right].
\end{equation}
We add $Q^\circ$ to the collection $\widehat{\mathcal{Q}}_\epsilon$ if and only if  \eqref{step4} holds.

\item \underline{\textsf{Step 5}}:  We examine all dyadic cubes $Q \subset Q^\circ$ such that $\delta_Q \geq \Lambda^{-10}$. 

We can test whether $Q \in \CZ(\cA^-)$ by querying the $\CZ(\cA^-)$-\textsc{Oracle} on the center of $Q$. We add $Q$ to the collection $\widehat{\mathcal{Q}}_\epsilon$ if and only if $Q \in \CZ(\cA^-)$.

\item \underline{\textsf{Step 6}}: We examine all cubes $Q \in \CZ(\cA^-)$ such that $\delta_Q \leq \Lambda^{-10}$ and $Q \cap E \neq \emptyset$.

Let $Q^\up$ be the dyadic cube with $Q \subset Q^\up$ and $\delta_{Q^\up} = \Lambda \delta_Q$. We determine whether or not
\begin{equation} \label{step6} \left[ \epsilon^{1/\kappa_5} \leq \eta^{(Q^\up,\cA)}(\delta_{Q^\up}) \right].
\end{equation}
We add $Q$ to the collection $\widehat{\mathcal{Q}}_\epsilon$ if and only if \eqref{step6} holds.
\end{itemize}

This completes the construction of $\widehat{\cQ}_\epsilon$. We examined at most $C(\Lambda) N$ cubes, and performed work at most $C \log N$ on each cube. Hence, the computation required work at most $C(\Lambda) N \log N$ in space $C(\Lambda) N$. We later choose $\Lambda$ to be a constant depending only on $m$, $n$, and $p$. We have thus not exceeded the work and storage guarantees of \textsc{Compute Critical Testing Cubes}. Moreover, we have $\# ( \mathcal{\hQ}_\epsilon) \leq C(\Lambda) \cdot N$, which implies condition (b).

If $Q$ belongs to $\mathcal{\cQ}_\epsilon$, then $Q$ was chosen in one of the six steps above (not including Step 0). We will examine the six cases separately and prove conditions (c) and (d) for the cube $Q$.

\noindent\textbf{Analysis of Step 1.}
Suppose that $Q$ was chosen in Step 1. Then $Q$ satisfies \eqref{step1}.

We use properties \textbf{(A1)} and \textbf{(A2)} of the function $\eta^{(\upQ_1,\cA)}$ from the algorithm \textsc{Optimize Basis}.

From \textbf{(A2)} and \eqref{step1}, we find that $\ooline{\sigma}(\upQ_1)$ does not have an $(\cA',x_{\upQ_1},\epsilon^{1/\kappa_7},\delta_Q)$-basis for any $\cA' \leq \cA$.

Since $Q_1 \subset Q$ are testing cubes, and  $\delta_Q \geq \Lambda^{10} \delta_{Q_1}$, we have $\#(3Q \cap E) \geq \#(9Q_1 \cap E) \geq 2$.

Also note that $(1+a(\cA))\upQ_1 \cap E = 3 Q \cap E$ if $\Lambda$ is greater than some constant determined by $m$, $n$, $p$; see Lemma \ref{gg_lem}.

Hence, Proposition \ref{tool2} implies that no cube containing $3 Q$ is tagged with $(\cA, \epsilon^{1/\kappa_8})$. This proves property (d).

To prove property (c), note that \textbf{(A1)} and \eqref{step1} imply that
\[\ooline{\sigma}(\upQ_1) \mbox{ has an } (\cA',x_{\upQ_1},\epsilon^{\kappa_7}, \delta_Q)\mbox{-basis for some} \; \cA' \leq \cA.
\]
Thus, Proposition \ref{tool2} shows that $(1+a(\cA))Q$ is tagged with $(\cA, \epsilon^{\kappa_8})$.

\noindent\textbf{Analysis of Step 2.}
Suppose that $Q$ was chosen in Step 2, and let $Q^\up$ be as in Step 2. Then $Q$ and $Q^\up$ satisfy \eqref{step2}.

We use properties \textbf{(A1)} and \textbf{(A2)} of the functions $\eta^{(Q,\cA)}$ and $\eta^{(Q^\up,\cA)}$ from the algorithm \textsc{Optimize Basis}.

Since $Q \subset Q^\up$ are testing cubes, and $\delta_{Q^\up} = \Lambda \delta_{Q}$, we have $\#(E \cap \frac{65}{64}Q^\up) \geq \#(E \cap 9 Q) \geq 2$ for sufficiently large $\Lambda$.

From \textbf{(A2)} and \eqref{step2}, we find that $\ooline{\sigma}(Q^\up)$ does not have an $(\cA',x_{Q^\up},\epsilon^{1/\kappa_6},\delta_{Q^\up})$-basis for any $\cA' \leq \cA$. Thus, Proposition \ref{tool1} implies that no cube containing $3Q^\up$ is tagged with $(\cA,\epsilon^{1/\kappa_7})$. In particular, since $3Q^\up \subset 100\Lambda Q$, we find that no cube containing $100 \Lambda Q$ is tagged with $(\cA,\epsilon^{1/\kappa_7})$. This proves property (d).

From \textbf{(A1)} and \eqref{step2}, we find that either $\#(\frac{65}{64}Q \cap E ) \leq 1$ or $\ooline{\sigma}(Q)$ has an $(\cA',x_{Q},\epsilon^{\kappa_6}, \delta_Q)$-basis for some $\cA' \leq \cA$. Thus, Proposition \ref{tool1} implies that $(1+a(\cA))Q$ is tagged with $(\cA, \epsilon^{\kappa_7})$. This proves property (c).

\noindent\textbf{Analysis of Step 3.} Note that (d) holds vacuously for all the cubes $Q \in \mathcal{\hQ}_\epsilon$ chosen in Step 3, assuming that $c^* \leq \Lambda^{-10}$.

As in the analysis of Step 2, \eqref{step3} implies that $(1+a(\cA))Q$ is tagged with $(\cA, \epsilon^{\kappa_7})$. This implies property (c) for any $Q$ picked in Step 3.

\noindent\textbf{Analysis of Step 4.} Suppose that $Q^\circ$ was chosen in Step 4. Note that (d) holds vacuously for $Q^\circ$.

As in the analysis of Step 2, \eqref{step4} shows that $(1+a(\cA))Q^\circ$ is tagged with $(\cA, \epsilon^{\kappa_7})$. This implies property (c) for $Q^\circ$.

\noindent\textbf{Analysis of Step 5.} We may assume that $c^* \leq \Lambda^{-10}$. Therefore, (c) and (d) are vacuously true for all the cubes $Q \in \mathcal{\hQ}_\epsilon$ chosen in Step 5.

\noindent\textbf{Analysis of Step 6.} Suppose that $Q$ was chosen in Step 6. Note that (c) holds vacuously for $Q$, since $Q \in \CZ(\cA^-)$.

Since $\delta_{Q^\up} = \Lambda \delta_{Q}$, and since $Q \subset Q^\up$ are testing cubes, we have $\#(E \cap \frac{65}{64}Q^\up) \geq \#(E \cap 9 Q) \geq 2$.

By \eqref{step6}  and property \textbf{(A2)} of the function $\eta^{(\upQ,\cA)}$ stated in \textsc{Optimize Basis}, we find that $\ooline{\sigma}(Q^\up)$ does not have an $(\cA',x_{Q^\up},\epsilon^{1/\kappa_6},\delta_{Q^\up})$-basis for any $\cA' \leq \cA$. Then Proposition \ref{tool1} guarantees that no cube containing $3Q^\up$ is tagged with $(\cA,\epsilon^{1/\kappa_7})$. Therefore, since $3Q^\up \subset 100 \Lambda Q$, we find that no cube containing $100 \Lambda Q$ is tagged with $(\cA,\epsilon^{1/\kappa_7})$. This implies property (d) for $Q$, and concludes the analysis of Step 6.

This completes the proof of (c) and (d) in all cases. An inspection of our argument shows that we may take $c^* = \Lambda^{-10}$ and $S = 128 \Lambda$.

Next we prove property (a).

Let $x \in E$ be given. Consider the finite sequence of cubes $Q_\nu \in {T}$ such that
\begin{equation} \label{chain} x \in Q_0 \subsetneq Q_1 \subsetneq \cdots \subsetneq Q_{\nu_{\max}} = Q^\circ,
\end{equation}
where  $Q_0 \in \CZ(\cA^-)$ and $Q_{\nu+1}$ is the parent of $Q_\nu$ in ${T}$. (We do not attempt to compute this sequence.)

We will show that there exists $Q' \in \widehat{\cQ}_\epsilon$ with $Q_0 \subset Q' \subset Q_{\nu_{\max}}$. This will complete the proof of (a).

Note that one of the following cases must occur.
\begin{description}
\item[The First Extreme Case] For all dyadic cubes $Q$ such that $Q_0 \subset Q \subset Q_{\nu_{\max}}$, the cube $3Q$ is tagged with $(\cA,\epsilon)$.
\item[The Second Extreme Case] For all dyadic cubes $Q$ such that $Q_0 \subset Q \subset Q_{\nu_{\max}}$, the cube $3Q$ is not tagged with $(\cA,\epsilon)$.
\item[The Main Case] For some dyadic cube $Q$ such that $Q_0 \subset Q \subsetneq Q_{\nu_{\max}}$ we find that exactly one of $3Q$, $3Q^+$ is tagged with $(\cA,\epsilon)$.
\end{description}

\noindent\underline{In the First Extreme Case}: 
\begin{equation}
\label{extreme1}
3 Q^\circ \; \mbox{is tagged with} \; (\cA,\epsilon).
\end{equation}

Notice that $\#(\frac{65}{64}Q^\circ \cap E) = \#(E) \geq 2$. From \eqref{extreme1} and Proposition \ref{tool1}, we see that $\ooline{\sigma}(Q^\circ)$ has an $(\cA',x_{Q^\circ},\epsilon^{\kappa_1},\delta_{Q^\circ})$-basis for some $\cA' \leq \cA$. Then property \textbf{(A2)} from \textsc{Optimize Basis} shows that $\eta^{(Q^\circ,\cA)}(\delta_{Q^\circ}) \leq \epsilon^{\kappa_5}$. Therefore, we decided to include $Q^\circ$ in $\widehat{\cQ}_\epsilon$ in Step 4.

This completes the analysis in the First Extreme Case.

\noindent\underline{In the Second Extreme Case}: 
\begin{equation}
\label{extreme2}
3 Q_0  \; \mbox{is not tagged with} \; (\cA,\epsilon).
\end{equation}

If $\delta_{Q_0} \geq \Lambda^{-10}$, then we decided to include $Q_0$ in $\widehat{\cQ}_\epsilon$ in Step 5.

Otherwise, suppose that $\delta_{Q_0} < \Lambda^{-10}$.

Let $Q^\up_0$ be a dyadic cube with $Q_0 \subset Q^\up_0 \subset Q^\circ$ and $\delta_{Q^\up_0} = \Lambda \delta_{Q_0}$.

Note that $3Q_0 \subset (1+a(\cA))Q^\up_0$, if $\Lambda$ is sufficiently large. Then \eqref{extreme2} and Lemma \ref{pre_lem5} imply that $(1+a(\cA))Q_0^\up$ is not tagged with $(\cA,\epsilon^{1/\kappa_1})$. Hence, Proposition \ref{tool1} shows that
\[\ooline{\sigma}(Q^\up_0) \; \mbox{does not have an} \; (\cA',x_{Q^\up_0}, \epsilon^{1/\kappa_2},\delta_{Q^\up_0})\mbox{-basis} \; \mbox{for any} \; \cA' \leq \cA.\]
Thus, property \textbf{(A1)} from \textsc{Optimize Basis} shows that $\eta^{(Q^\up_0,\cA)}(\delta_{Q^\up_0}) \geq \epsilon^{1/\kappa_5}$. Therefore, we decided to include $Q_0$ in $\widehat{\cQ}_\epsilon$ in Step 6. (Recall that $x \in Q_0$, hence $E \cap Q_0 \neq \emptyset$.)

This completes the analysis in the Second Extreme Case.

\noindent\underline{In the Main Case}: Exactly one of $3Q$, $3Q^+$ is tagged with $(\cA,\epsilon)$, thus
\begin{equation}
\label{main1} 3Q \; \mbox{is tagged with} \; (\cA,\epsilon^{\kappa_0}) \;\; (\mbox{see Lemma} \; \ref{pre_lem5}),
\end{equation}
and
\begin{equation}
\label{main2} 3Q^+ \; \mbox{is not tagged with} \; (\cA,\epsilon^{1/\kappa_0}) \;\; (\mbox{again, see Lemma} \; \ref{pre_lem5}).
\end{equation}
We now consider three subcases of the Main Case.
\begin{itemize}
\item The Geometrically Interesting (``GI'') subcase: For some $\nu$,
\begin{equation} \label{gi0} Q_\nu \subset Q \subset Q_{\nu+1}, \;\; \left[\delta_Q \leq \Lambda^{10} \delta_{Q_\nu}  \; \mbox{or} \;   \delta_Q \geq \Lambda^{-10} \delta_{Q_{\nu+1}}  \right], \; \mbox{and} \;  \delta_{Q} \leq \Lambda^{-10}.
\end{equation}

\item The Geometrically Uninteresting (``GUI'') subcase: For some $\nu$,
\begin{equation} \label{gui3} Q_\nu \subset Q \subset Q_{\nu+1} \; \mbox{and} \; \Lambda^{10} \delta_{Q_\nu} \leq \delta_Q \leq \Lambda^{-10} \delta_{Q_{\nu+1}}.\end{equation}
\item The Near-Maximal (``NM'') subcase: 
\begin{equation} \label{nm0} \delta_Q \geq \Lambda^{-10}.
\end{equation}

\end{itemize}

\noindent\underline{First consider the GI subcase}.

From Proposition \ref{tool1} and \eqref{main1} we see that
$$\#\left(\frac{65}{64}Q \cap E\right) \leq 1 \;\; \mbox{or} \;\; \ooline{\sigma}(Q) \; \mbox{has an} \; (\cA',x_Q,\epsilon^{\kappa_1},\delta_Q)\mbox{-basis} \; \mbox{for some} \; \cA' \leq \cA.$$
Thus, by property \textbf{(A2)} from \textsc{Optimize Basis},
\begin{equation}
\label{eq01} \#\left(\frac{65}{64}Q \cap E\right) \leq 1 \;\; \mbox{or} \;\; \eta^{(Q,\cA)}(\delta_Q) \leq \epsilon^{\kappa_5}.
\end{equation}

Pick $Q^\up$ (dyadic) such that $Q \subset Q^\up \subset Q^\circ$ and $\delta_{Q^\up} = \Lambda \delta_Q$. (Recall that $\delta_Q \leq \Lambda^{-10}$.) Then $3Q^+ \subset (1+a(\cA))Q^\up$, assuming that $\Lambda$ is sufficiently large. Thus, \eqref{main2} shows that
\begin{equation*} (1+a(\cA))Q^\up \; \mbox{is not tagged with} \; (\cA,\epsilon^{1/\kappa_1}) \qquad (\mbox{see Lemma \ref{pre_lem5}}).
\end{equation*}
Therefore, Proposition \ref{tool1} gives that
\[\ooline{\sigma}(Q^\up) \; \mbox{doesn't have an} \; (\cA',x_{Q^\up}, \epsilon^{1/\kappa_2},\delta_{Q^\up})\mbox{-basis} \; \mbox{for any} \; \cA' \leq \cA.\]
Hence, using property \textbf{(A1)} from \textsc{Optimize Basis},
\begin{equation}
\label{eq02} \eta^{(Q^\up,\cA)}(\delta_{Q^\up}) \geq \epsilon^{1/\kappa_5}
\end{equation}

From \eqref{eq01} and \eqref{eq02}, we see that $Q$ was included in $\widehat{\cQ}_\epsilon$  in Step 2. This completes the analysis in the GI subcase.

\noindent\underline{Next consider the GUI subcase}.

Since $Q_\nu \subset Q$ are testing cubes, and $\delta_Q \geq \Lambda^{10} \delta_{Q_\nu}$, we have  $\#(E \cap 3Q) \geq \#(E \cap 9 Q_\nu) \geq 2$.

Let $\upQ_\nu$ denote the dyadic cube with $Q_\nu \subset \upQ_\nu$ and $\delta_{\upQ_\nu} = \Lambda \cdot \delta_{Q_\nu}$

Note that $(1+a(\cA))\upQ_\nu \cap E = 3 Q \cap E$, as long as $\Lambda \geq C$ for a large enough universal constant $C$ (see Lemma \ref{gg_lem}).

From Proposition \ref{tool2} and assumption \eqref{main1} (from the Main Case), we see that
\[\ooline{\sigma}(\upQ_\nu) \; \mbox{has an} \; (\cA',x_{\upQ_\nu},\epsilon^{\kappa_1},\delta_Q)\mbox{-basis}, \; \mbox{for some} \; \cA' \leq \cA.\]
Hence, condition \textbf{(A2)} in the algorithm \textsc{Optimize Basis} implies that
\begin{equation}
\label{eq4} 
\eta^{(\upQ_\nu,\cA)}(\delta_Q) \leq \epsilon^{\kappa_5}.
\end{equation}

Let $Q^\up$ be a dyadic cube with $Q \subset Q^\up \subsetneq Q_{\nu+1}$ and $\delta_{Q^\up} = \Lambda \cdot \delta_Q$. Such a dyadic cube exists because we are assuming that $\delta_{Q} \leq \Lambda^{-10} \delta_{Q_{\nu+1}}$. Then $ (1+a(\cA)) \upQ_\nu \cap E = 3Q^\up \cap E$ thanks to Lemma \ref{gg_lem}.

For large enough $\Lambda$, we have $3Q^+ \subset (1+a(\cA))Q^\up$. Thus, Lemma \ref{pre_lem5} and \eqref{main2} imply that $(1+a(\cA))Q^\up$ is not tagged with  $(\cA,\epsilon^{1/\kappa_1})$. We apply conclusion (A) in Proposition \ref{tool2} to the testing cubes $\upQ_\nu \subset \upQ$ 
in order to deduce that
\[ \ooline{\sigma}(\upQ_\nu) \; \mbox{does not have an} \; (\cA',x_{\upQ_\nu}, \epsilon^{1/\kappa_2},\delta_{Q^\up})\mbox{-basis}, \; \mbox{for any} \; \cA' \leq \cA.\]
Thus, property \textbf{(A1)} from \textsc{Optimize Basis} shows that $\eta^{(\upQ_\nu,\cA)}(\delta_{Q^\up}) \geq \epsilon^{1/\kappa_3}$. Moreover, property \textbf{(A3)} from \textsc{Optimize Basis} implies that $\eta^{(\upQ_\nu,\cA)}(\delta_{Q^\up}) \leq C(\Lambda) \eta^{(\upQ_\nu,\cA)}(\delta_Q)$, hence we have
\begin{equation} \label{eq5}
\eta^{(\upQ_\nu,\cA)}(\delta_Q) \geq \epsilon^{1/\kappa_5}.
\end{equation}

We are assuming that $\delta_Q \in [ \Lambda^{10} \delta_{Q_\nu} , \Lambda^{-10} \delta_{Q_{\nu+1}}]$ (from the GUI subcase). Hence, from \eqref{eq4} and \eqref{eq5}, we see that in Step 1 we included in $\widehat{\cQ}_\epsilon$ a dyadic cube $Q'$ such that $Q_\nu \subset Q' \subset Q_{\nu+1}$. This completes the analysis in the GUI subcase.

\noindent\underline{Finally, consider the NM subcase}.

From Proposition \ref{tool1} and \eqref{main1} we have
$$\#\left(\frac{65}{64}Q \cap E\right) \leq 1 \;\; \mbox{or} \;\; \ooline{\sigma}(Q) \; \mbox{has an} \; (\cA',x_Q,\epsilon^{\kappa_1},\delta_Q)\mbox{-basis} \; \mbox{for some} \; \cA' \leq \cA.$$
Then property \textbf{(A2)} from \textsc{Optimize Basis} implies that
\begin{equation}
\label{eq6} 
\#\left(\frac{65}{64}Q \cap E\right) \leq 1 \;\; \mbox{or} \;\; \eta^{(Q,\cA)}(\delta_Q) \leq \epsilon^{\kappa_5}.
\end{equation}
Thus, we included $Q$ in $\widehat{\cQ}_\epsilon$ in Step 3. This completes the analysis in the NM subcase. 

Thus, in all the cases, we see that there exists $Q' \in \widehat{\cQ}_\epsilon$ with $Q_0 \subset Q' \subset Q_{\nu_{\max}}$, where $Q_0$ is the unique cube in $\CZ(\cA^-)$ containing the point $x \in E$ . As mentioned before, this completes the proof of (a).

We fix a large enough constant $\Lambda = 2^{J} \geq 1$, depending only on $m$, $n$, and $p$.

This completes the explanation of the algorithm \textsc{Compute Critical Testing Cubes}.

\end{proof}

\subsection{Lengthscales}

\label{sec_cl3}

Using the algorithm \textsc{Compute Critical Testing Cubes} from the previous section, we compute a collection $\widehat{\mathcal{Q}}_\epsilon$ consisting of dyadic subcubes of $Q^\circ$. We proved that each point of $E$ belongs to a cube in $\widehat{\mathcal{Q}}_\epsilon$. Applying the algorithm \textsc{Placing a Point Inside Target Cuboids} (see Section \ref{maketree_sec}), we obtain the following algorithm.

\environmentA{Algorithm: Compute Lengthscales.}

For each $x \in E$ we compute a cube $Q_x \in \widehat{\mathcal{Q}}_\epsilon$ containing $x$. This requires work at most $C N \log N$ in space $C N$.

We write $c^* > 0$ and $S \geq 1$ for the universal constants from the algorithm \textsc{Compute Critical Testing Cubes}. The conclusion of this algorithm implies the next result.

\begin{prop} \label{lengthscales_prop} For each $x \in E$, the following properties hold.
\begin{description}
\item[(LS1)] Suppose that $Q_x$ strictly contains a cube of $\CZ(\cA^-)$. Then $(1+a(\cA))Q_x$ is tagged with $(\cA,\epsilon^\kappa)$.
\item[(LS2)] Suppose that $\delta_{Q_x} \leq c^*$. Then no cube containing $S Q_x$ is tagged with $(\cA,\epsilon^{1/\kappa})$.
\end{description}
Here, $\kappa > 0$ is a small universal constant.
\end{prop}

\section{Passing from Lengthscales to CZ Decompositions}

\label{sec_cz}

For each $x \in E$ we compute the sidelength
\begin{equation}\label{defn_Delta}
\Delta_{\cA}(x) := \delta_{Q_x}.
\end{equation}
Here, we compute the cube $Q_x$ using the algorithm \textsc{Compute Lengthscales}. Recall that $x \in Q_x$ for each $x \in E$. Since $Q_x \subset Q^\circ$, we know that 
\begin{equation} \label{ls_prop1}
\Delta_\cA(x) \in (0,1] \; \text{for all} \; x \in E.
\end{equation}

Let $Q \subset Q^\circ$ be a testing cube. We say that $Q$ is $OK(\cA)$ provided that either $Q \in \CZ(\cA^-)$ or $\Delta_{\cA}(x) \geq K \delta_Q$ for all $x \in E \cap 3Q$, where we set $K := \frac{10^9}{a(\cA)}$. Recall that we have defined the constant $a(\cA)$ in equation \eqref{newa_defn}. In particular, since $a(\cA) \leq 1$, we see that $K \geq 1$.

We define a Calder\'on-Zygmund decomposition $\CZ(\cA)$ of the unit cube $Q^\circ$ to consist of the maximal dyadic subcubes $Q \subset Q^\circ$ that are $OK(\cA)$.

We will prove properties \textbf{(CZ1-CZ5)} in the Main Technical Results for $\cA$.

First, however, we produce a $\CZ(\cA)$-\textsc{Oracle} as described in Chapter \ref{sec_mainresults}. The decomposition $\CZ(\cA)$ coincides with the decomposition $\CZ_\new$ from Section \ref{sec_czoracle}, where we use $\CZ_{\old} = \CZ(\cA^-)$ and $\Delta(x) = \Delta_\cA(x)/K$ in the notation therein. Note that $\Delta(x) \in (0,1]$ for each $x \in E$, hence the assumptions in Section \ref{sec_czoracle} hold. The \textsc{Glorified CZ-Oracle} coincides with the $\CZ(\cA)$-\textsc{Oracle} described in Chapter \ref{sec_mainresults}.

\begin{prop}\label{prop561}
The collection $\CZ(\cA)$ is a partition of $Q^\circ$ into pairwise disjoint dyadic subcubes.
\end{prop}
\begin{proof}
Each point $x \in Q^\circ$ belongs to some cube $Q_0 \in \CZ(\cA^-)$. Note that $Q_0$ is $OK(\cA)$, and hence $Q_0$ is contained in a maximal dyadic subcube $Q \subset Q^\circ$ that is also $OK(\cA)$. Thus, each point $x \in Q^\circ$ is contained in some cube $Q \in \CZ(\cA)$. 

Any two distinct cubes $Q, Q' \in \CZ(\cA)$ are dyadic, hence either $Q, Q'$ are disjoint or one of $Q, Q'$ contains the other. The latter case cannot occur, by definition of $\CZ(\cA)$. It follows that the cubes in $\CZ(\cA)$ are pairwise disjoint.
\end{proof}

Our previous decomposition $\CZ(\cA^-)$ clearly refines $\CZ(\cA)$. This establishes property \textbf{(CZ5)} for $\cA$. We now prove the remaining properties \textbf{(CZ1-CZ4)}.

We prove property \textbf{(CZ1)} in the next result.

\begin{prop}\label{prop_gg}
The cubes in $\CZ(\cA)$ have good geometry.
\end{prop}
\begin{proof}
For the sake of contradiction suppose that there are cubes $Q,Q' \in \CZ(\cA)$ such that $Q \leftrightarrow Q'$ and $\delta_Q \leq \frac{1}{4} \delta_{Q'}$. It follows that $3Q^+ \subset 3Q'$. 

First, suppose $Q' \in \CZ(\cA^-)$. Since $\CZ(\cA^-)$ refines $\CZ(\cA)$, there exists a cube $Q'' \in \CZ(\cA^-)$ with $Q'' \subset Q$ and $Q'' \leftrightarrow Q'$. Note that $\delta_{Q''} \leq \delta_Q \leq \frac{1}{4} \delta_{Q'}$. But this contradicts our assumption that the cubes in $\CZ(\cA^-)$ satisfy good geometry.

Next, suppose $Q' \notin \CZ(\cA^-)$. By definition of $\CZ(\cA)$ we know that $Q^+$ is not $OK(\cA)$, hence there exists $x \in E \cap 3Q^+$ with $\Delta_{\cA}(x) < K \delta_{Q^+}$. Thus, $x \in E \cap 3Q'$ and $\Delta_{\cA}(x) < K \delta_{Q'}$. Since also $Q' \notin \CZ(\cA^-)$ we see that $Q'$ is not $OK(\cA)$. But this contradicts our assumption that $Q' \in \CZ(\cA)$.
\end{proof}

\begin{prop}\label{mainprops}
There exists a universal constant $c_* > 0$ such that, for any $Q \in \CZ(\cA)$, the following conditions hold.
\begin{enumerate}[(a)]
\item If $Q$ is not $c_*$-simple then $3Q$ is tagged with $(\cA, \epsilon^\kappa)$.
\item If $\delta_Q \leq c_*$ then $WQ$ is not tagged with $(\cA,\epsilon^{1/\kappa})$.
\end{enumerate}
Here, $\kappa > 0$ and $W \in \N$ are universal constants.
\end{prop}
\begin{proof}

We choose $c_*$ much smaller than the constant $c^*$ from Proposition \ref{lengthscales_prop}. 

We now prove (a). Assume that $Q \in \CZ(\cA)$ is not $c_*$-simple. Then there exists $\overline{Q} \in \CZ(\cA^-)$ with $\overline{Q} \subset \frac{65}{64}Q$ and $\delta_{\overline{Q}} \leq c_* \delta_Q$. For small enough $c_*$ this implies that $9\overline{Q} \subset 3 Q$. Recall \eqref{Enearby2}, which implies that $9 \overline{Q} \cap E \neq \emptyset$, hence $3 Q \cap E \neq \emptyset$. We fix $x \in E \cap 3Q$.

We have $\delta_{Q_x} = \Delta_{\cA}(x) \geq K \delta_Q$ with $K = 10^9/a(\cA)$, because $Q$ is $OK(\cA)$ and $Q \notin \CZ(\cA^-)$; see also \eqref{defn_Delta}. 

For any $y \in 3Q$ we have $\lvert y - x \rvert \leq 3 \delta_Q \leq \frac{3}{K} \delta_{Q_x}$, because $x \in 3Q$. Moreover, $\lvert x - x_{Q_x} \rvert \leq \frac{1}{2} \delta_{Q_x}$, because $x \in Q_x$. (Recall that $x_{Q_x}$ is the center of $Q_x$.) Thus,
\begin{align*}
\lvert y - x_{Q_x} \rvert &\leq \left[ \frac{1}{2} + \frac{3}{K} \right] \delta_{Q_x} \\
& \leq \frac{1}{2} \left[1 + a(\cA) \right] \delta_{Q_x} \quad \mbox{for any} \; y \in 3Q.
\end{align*}
(Here, we use that $K = 10^9/a(\cA)$.) Hence,
\begin{equation}\label{qq5}
3Q \subset (1+a(\cA))Q_x.
\end{equation}
(Recall that we  are working with the $\ell^\infty$ metric.)

We now prove that $Q_x$ strictly contains a cube of $\CZ(\cA^-)$. Assume for the sake of contradiction that $Q_x$ is contained in a cube in $\CZ(\cA^-)$. (For a dyadic cube this is the only alternative.) Since $Q_x$ is a testing cube, it follows that $Q_x$ belongs to $\CZ(\cA^-)$; see Section \ref{sec_testcube} where the notion of a testing cube is defined. From \textbf{(CZ5)} we see that $\CZ(\cA^-)$ refines $\CZ(\cA)$, hence there exists $\widetilde{Q} \in \CZ(\cA^-)$ with $\widetilde{Q} \subset Q$. From \eqref{qq5} we have
\[\widetilde{Q} \subset 3 Q \subset (1+a(\cA)) Q_x \subset \frac{65}{64}Q_x.\]
Since $\widetilde{Q} \in \CZ(\cA^-)$ and $Q_x \in \CZ(\cA^-)$, from good geometry of the cubes in $\CZ(\cA^-)$ and from Lemma \ref{gg_lem_1}, we deduce that $\frac{1}{2} \delta_{\widetilde{Q}} \leq \delta_{Q_x} \leq 2 \delta_{\widetilde{Q}}$. Hence, because the cubes in $\CZ(\cA^-)$ are pairwise disjoint and dyadic, we must have $Q_x = \widetilde{Q}$. Thus, we have
\[(1+a(\cA))Q_x \subsetneq 3\widetilde{Q} \subset 3Q.\]
However, this contradicts \eqref{qq5}. This completes our proof that $Q_x$ strictly contains a cube of $\CZ(\cA^-)$.

Hence, from \textbf{(LS1)} in Proposition \ref{lengthscales_prop} we deduce that $(1+a(\cA))Q_x$ is tagged with $(\cA,\epsilon^\kappa)$; hence, $3Q$ is tagged with $(\cA,\epsilon^{\kappa'})$ for some universal constant $\kappa'$, thanks to Lemma \ref{pre_lem5} and \eqref{qq5}. This completes the proof of (a).

We now prove (b). Let $S$ be the universal constant in Proposition \ref{lengthscales_prop}. Assume that $Q \in \CZ(\cA)$ satisfies $\delta_Q \leq c_*$. 

Since $Q^+$ is not $OK(\cA)$, there exists $x \in E \cap 3Q^+$ such that 
\[ \delta_{Q_x} = \Delta_{\cA}(x) \leq K \delta_{Q^+}.
\]
Hence, because $x \in Q_x$ and $x \in 3Q^+$, we have $SQ_x \subset WQ$ for a large enough integer constant $W \geq 1$ depending only on $K$ and $S$. Recall that $K = 10^9/a(\cA)$ is a universal constant. Hence, we can choose $W$ to be a universal constant. Therefore, $\delta_{Q_x} \leq \frac{W}{S} \delta_Q \leq \frac{W}{S} c_* \leq c^*$. Here, we assume that $c_* \leq \frac{S}{W} c^*$. 

From \textbf{(LS2)} in Proposition \ref{lengthscales_prop} it follows that $WQ$ is not tagged with $(\cA,\epsilon^{1/\kappa})$.

This completes the proof of the proposition.
\end{proof}

We have proven \textbf{(CZ2)} and \textbf{(CZ3)} in the Main Technical Results for $\cA$, where we set
\begin{equation}
\label{fix_con}
\left\{
\begin{aligned}
&c_*(\cA) = c_*/2, \\
&S(\cA) = W, \\
&\epsilon_1(\cA) = \epsilon^{1/\kappa}, \\
&\epsilon_2(\cA) = \epsilon^\kappa.
\end{aligned}
\right.
\end{equation}
Here, $\kappa$ and $W$ are as in Proposition \ref{mainprops}. Note that \textbf{(CZ4)} holds vacuously, since we are assuming that $\cA \neq \cM$. We have thus proven \textbf{(CZ1-CZ5)} for the label $\cA$. We will later pick $\epsilon$ to be a small enough universal constant, at which point $\epsilon_1(\cA)$ and $\epsilon_2(\cA)$ will be determined once and for all.

We let $\CZ_{\main}(\cA)$ denote the collection of all cubes $Q \in \CZ(\cA)$ that satisfy $\frac{65}{64}Q \cap E \neq \emptyset$. We note that the collection $\{ \frac{65}{64}Q : Q \in \CZ(\cA)\}$ has bounded overlap, thanks to the good geometry of the cubes in $\CZ(\cA)$ (see Lemma \ref{prop_gg}). Hence,
\begin{equation}
\label{maincubesbd}
\# \left( \CZ_{\main}(\cA) \right) \leq C \cdot N.
\end{equation}

\section{Completing the Induction}
\label{compind_sec}

In the previous section, we defined a decomposition $\CZ(\cA)$ and gave a $\CZ(\cA)$-\textsc{Oracle}. Here, we construct the remaining objects in the Main Results for $\cA$.

We can compute a list of all the cubes $Q$ in $ \CZ_{\main}(\cA)$. We list  all the cubes $Q \in \CZ(\cA)$ that satisfy $E \cap \frac{65}{64}Q \neq \emptyset$ using the algorithm \textsc{Find Main-Cubes} in Section \ref{czalg_sec}.

We now show that for each $\hQ \in  \CZ_{\main}(\cA)$ we can efficiently collect all the ingredients we need to compute the assists, functionals, and local extension operator relevant to $\hQ, \cA$. 

Recall the notion of \underline{supporting data} associated to a testing cube; see Section \ref{sec_computingstuff}.

\environmentA{Algorithm: Produce All Supporting Data}

We produce the supporting data for each cube $\hQ$ in $\CZ_{\main}(\cA)$, using work at most $C N \log N$ in space $C N$.

\begin{proof}[\underline{Explanation}] 

We produce the cubes $Q, Q_{\spec}, Q^\#$ and the pairs of cubes $(Q',Q'')$ that arise in (SD1)-(SD5) in Section \ref{sec_computingstuff} for some testing cube $\hQ \in \CZ_{\main}(\cA)$.

For each $Q \in \CZ_{\main}(\cA^-)$, we apply the $\CZ(\cA)$-\textsc{Oracle} to find the cube $\hQ \in \CZ(\cA)$ that contains $Q$, as well as all the cubes $\hQ' \in \CZ(\cA)$ such that $\hQ' \leftrightarrow \hQ$. For each such $\hQ$ (or $\hQ'$), we check whether $\hQ$ (or $\hQ'$) appears in the list $\CZ_{\main}(\cA)$; if it does, then we check whether $Q \subset (1+t_G) \hQ$ (or $ (1+t_G) \hQ'$). If so, then we add the cube $Q$ to the list of cubes in (SD1) relevant to the testing cube $\hQ$ (or $\hQ'$).

Similarly, for each pair $(Q',Q'') \in \CZ(\cA^-) \times \CZ(\cA^-)$ such that $Q' \leftrightarrow Q''$ but $\mathcal{K}(Q') \neq \mathcal{K}(Q'')$ (the ``border disputes''), we look for all possible $\hQ \in \CZ_{\main}(\cA)$ such that $(Q',Q'')$ arises in (SD2) for the testing cube $\hQ$. That is, we look for all the $\hQ \in \CZ_{\main}(\cA)$ such that $Q' \subset (1+t_G)\hQ$ and $\delta_{Q'} < t_G \delta_{\hQ}$.

To find all the $\hQ$ as above, we need only search among the cubes $\hQ' =$ the cube of $\CZ(\cA)$ containing $Q'$, and the cubes of $\CZ(\cA)$ that touch $\hQ'$. We obtain all those cubes by making at most $C$ calls to the $\CZ(\cA)$-\textsc{Oracle} and doing additional work at most $C$.

We check each $\hQ$ obtained as above to see whether $\hQ \in \CZ_{\main}(\cA)$, and if so whether also $\hQ$ has the desired relationship with $Q'$. For each surviving $\hQ$, we add $(Q',Q'')$ to the list of cubes in (SD2) relevant to that $\hQ$.

To find all the $Q \in \CZ(\cA^-)$ that arise in (SD3), we loop over all the $\hQ \in \CZ_{\main}(\cA)$. For each fixed $\hQ$, we examine all the dyadic cubes $Q \subset (1+t_G)\hQ$ such that $\delta_Q \geq t_G^2 \delta_{\hQ}$. (There are only $C$ such $Q$.) We test $Q$ to see whether it belongs to $\CZ(\cA^-)$; if so, then we add $Q$ to the list of cubes in (SD3) relevant to $\hQ$.

For the supporting data in (SD4), we can loop over all $\hQ \in \CZ_{\main}(\cA)$. For each such $\hQ$, we can just take $Q_{\spec}$ to be the $\CZ(\cA^-)$-cube containing the center of $\hQ$.

Finally, we loop over all keystone cubes $Q^\#$ of $\CZ(\cA^-)$. For each such $Q^\#$, we look for all the $\hQ \in \CZ_{\main}(\cA)$ such that $S_1 Q^\# \subset (65/64)\hQ$.

To find all the $\hQ$ as above, we need only search among the cubes $\hQ' =$ the cube of $\CZ(\cA)$ containing $Q^\#$, and the cubes of $\CZ(\cA)$ that touch $\hQ'$. We obtain all those cubes by making at most $C$ calls to the $\CZ(\cA)$-\textsc{Oracle} and doing additional work at most $C$.

We check each $\hQ$ obtained as above to see whether $\hQ \in \CZ_{\main}(\cA)$, and if so whether also $\hQ$ has the desired relationship with $Q^\#$. If those conditions are satisfied, then we add $Q^\#$ to the list of cubes in (SD5) relevant to $\hQ$.

Once we have carried out the above, then for each $\hQ \in \CZ_{\main}(\cA)$, we have a list of all the cubes $Q, Q_{\spec},Q^\#$ and of all the pairs of cubes $(Q',Q'')$ relevant to the supporting data (SD1)-(SD5) for the given $\hQ$. Again, see Section \ref{sec_computingstuff}.

This uses work $\mathcal{O}(N \log N)$ in space $\mathcal{O}(N)$. This completes our explanation  of the algorithm \textsc{Produce All Supporting Data}.

\end{proof}

Next, we will define lists $\Omega(\hQ,\cA) \subset \left[ \X(\frac{65}{64}\hQ \cap E) \right]^*$ and $\Xi(\hQ,\cA) \subset \left[ \X(\frac{65}{64}\hQ \cap E) \oplus \cP \right]^*$ and  also a linear extension operator $T_{(\hQ,\cA)} : \X(\frac{65}{64}\hQ \cap E) \oplus \cP \rightarrow \X$ for each $\hQ \in \CZ_{\main}(\cA)$. We will prove that these objects satisfy the properties laid out in the third, fourth and fifth bullet points in the Main Technical Results for $\cA$ (see Chapter \ref{sec_mainresults}).

For each $\hQ \in \CZ_{\main}(\cA)$, we can define
\[
M_{\hQ}(f,P)  = \left( \sum_{\xi \in \Xi(\hQ,\cA)} \lvert \xi(f,P) \rvert^p \right)^{1/p}.
\]
We need to prove the estimates in the fourth bullet point in the Main Technical Results for $\cA$. These estimates are
\begin{equation}
\label{q1} c \|(f,P) \|_{(1+a(\cA))\hQ} \leq M_{\hQ}(f,P)
\end{equation}
and
\begin{equation}
\label{q2} 
M_{\hQ}(f,P) \leq C \| (f,P) \|_{\frac{65}{64}\hQ}.
\end{equation}

Recall that a testing cube $\hQ$ is called $\lambda$-simple if for every $Q \in \CZ(\cA^-)$ with $Q \subset \frac{65}{64} \hQ$ we have $\delta_Q \geq \lambda \cdot \delta_{\hQ}$. We can determine whether a given cube $\hQ$ is $\lambda$-simple using work at most $C(\lambda)$, and at most $C(\lambda)$ calls to the $\CZ(\cA^-)$-\textsc{Oracle}. Here, $C(\lambda)$ is a constant depending only on $\lambda$ and $n$.

Let $c_*$ be the universal constant in Proposition \ref{mainprops}.

We loop over all the cubes $\hQ \in \CZ_{\main}(\cA)$. We can determine in time $\cO(\log N)$ whether $\hQ$ is $c_*$-simple. (Recall that a call to the $\CZ(\cA^-)$-\textsc{Oracle} requires work $\cO(\log N)$.) The body of our loop separates into two cases depending on the result of the test.

\subsection{Case I: Non-simple cubes} \label{notsimple_sec} We suppose that $\hQ \in \CZ_{\main}(\cA)$ is \underline{not} $c_*$-simple (the non-simple case). We will explain how to construct the objects in the Main Technical Results for $\cA$ relevant to $\hQ$.

We have already computed the supporting data for all the cubes in $\CZ_{\main}(\cA)$. By executing the algorithms \textsc{Compute New Assists} and \textsc{Compute New Assisted Functionals} (see Section \ref{sec_computingstuff}), we can compute 
\begin{enumerate}
\item[(a)] A list of \emph{assist functionals}: $\Omega(\hQ) \subset \left[ \X\left(E \cap (65/64) \hQ\right) \right]^*$ (see \eqref{test_assists}), and 
\item[(b)] A list of \emph{assisted functionals}: $\Xi(\hQ) \subset \left[ \X\left(E \cap (65/64)\hQ\right) \oplus \cP \right]^*$. 
\end{enumerate}
Each functional $\xi \in \Xi(\hQ)$ has $\Omega(\hQ)$-assisted bounded depth, and is written in short form in terms of the assists $\Omega(\hQ)$.

We define $\Omega(\hQ,\cA) := \Omega(\hQ)$, $\Xi(\hQ,\cA) := \Xi(\hQ)$, and
\[M_{\hQ}(f,P) = \left(\sum_{\xi \in \Xi(\hQ)} \lvert \xi(f,P) \rvert^p \right)^{1/p}.\]

We now prove the estimates \eqref{q1} and \eqref{q2}.

The estimate \eqref{q1} is a direct consequence of the unconditional inequality in Proposition \ref{inc_prop}.

Since $\hQ \in \CZ_{\main}(\cA)$ and $\hQ$ is not $c_*$-simple, we know that $3\hQ$ is tagged with $(\cA,\epsilon^\kappa)$ (see Proposition \ref{mainprops}). We may assume that $\epsilon^\kappa \leq \epsilon_0$, with $\epsilon_0$ as in Proposition \ref{inc_prop}. Thus, $3\hQ$ is tagged with $(\cA,\epsilon_0)$. Hence, the conditional inequality in Proposition \ref{inc_prop} implies the estimate \eqref{q2}.

Next, we estimate how much work and storage are used to compute the lists $\Omega(\hQ,\cA)$ and $\Xi(\hQ,\cA)$ for \underline{all} the non-simple cubes $\hQ \in \CZ_{\main}(\cA)$. We will prove that the \underline{total work} is at most $C N \log N$ and that the storage used is at most $C N$.

We examine the algorithms \textsc{Compute New Assists} and \textsc{Compute New Assisted Functionals} (see Section \ref{sec_computingstuff}). We see that we can compute all the lists $\Omega(\hQ,\cA)$ and $\Xi(\hQ,\cA)$ for all the non-simple cubes $\hQ \in \CZ_{\main}(\cA)$, using total work at most
\begin{align*}
&\sum_{\hQ \in \CZ_{\main}(\cA)} \bigl\{ \Work_1(\hQ) + \Work_2(\hQ) \bigr\} \\
& \hspace{.5cm} \leq C \log N \cdot  \sum_{\hQ \in \CZ_{\main}(\cA)} \biggl\{ 1 + \sum_{\substack{Q \in \CZ_{\main}(\cA^-) \\ Q \subset (1+t_G)\hQ }}  \left[ \sum_{\omega \in \Omega(Q,\cA^-)} \depth(\omega) + \# \bigl( \Xi(Q,\cA^-) \bigr)  \right] \\
& \hspace{4cm} + \sum_{\substack{ \text{keystone } Q^\# \in \CZ(\cA^-) \\ S_1 Q^\# \subset \frac{65}{64}\hQ}} \sum_{\omega \in \Omega^{\new}(Q^\#)} \depth(\omega) \\
&\hspace{4cm} + \# \bigl\{ (Q',Q'') \in \BD(\cA^-) : \; Q' \subset (1+t_G)\hQ, \;  \delta_{Q'} < t_G \delta_{\hQ} \bigr\}  \biggr\}.
\end{align*}
See \eqref{work1}, \eqref{bd1}, and \eqref{work2}, for the definitions of the quantities $\Work_1(\hQ)$ and $\Work_2(\hQ)$. Recall that $t_G$ is now a fixed universal constant, and so $C(t_G)$ in \eqref{work2} is a universal constant $C$.

Each cube $Q$ in $\CZ_{\main}(\cA^-)$, each keystone cube $Q^\# \in \CZ(\cA^-)$, and each pair $(Q',Q'') \in \BD(\cA^-)$ participates above for at most $C$ distinct $\hQ$ in $\CZ(\cA)$. This follows because the collection $\left\{ (65/64)\hQ : \hQ \in \CZ(\cA) \right\}$ has bounded overlap, which follows from the good geometry of $\CZ(\cA)$. Thus, by reversing the order of summation in the above expression, we see that the total work is bounded by
\begin{align*}
 C \cdot \log N \cdot \biggl\{  \# \bigl( \CZ_{\main}(\cA) \bigr) & +  \sum_{Q \in \CZ_{\main}(\cA^-) } \left[ \sum_{\omega \in \Omega(Q,\cA^-)} \depth(\omega)   +  \# \bigl( \Xi(Q,\cA^-) \bigr)  \right] \\
 & + \sum_{\text{keystone } Q^\# \in \CZ(\cA^-)}  \sum_{\omega \in \Omega^{\new}(Q^\#)} \depth(\omega)\\
& + \# \bigl( \BD(\cA^-) \bigr) \biggr\}.
\end{align*}
According to the Main Technical Results for $\cA^-$ and \eqref{maincubesbd}, the sum of terms inside the curly brackets in the first line above is bounded by $CN$. According to the algorithm \textsc{Make New Assists and Assign Keystone Jets}, the term on the second line above is bounded by $CN$. According to the \textsc{Keystone-Oracle}, the term on the last line above is bounded by $CN$. Hence, with work at most $C N \log N$, we can compute the lists $\Omega(\hQ)$ and $\Xi(\hQ)$ for all the non-simple cubes $\hQ \in \CZ_{\main}(\cA)$.

Similarly, we see that the computation of the lists $\Omega(\hQ)$ and $\Xi(\hQ)$  for \underline{all} the non-simple cubes $\hQ \in \CZ_{\main}(\cA)$ requires space at most $C N$.

Next, we explain how to define a linear extension operator associated to a non-simple $\hQ \in \CZ_{\main}(\cA)$ as in the Main Technical Results for $\cA$.

We define the map $T_{\hQ} : \X(E \cap \frac{65}{64}\hQ) \oplus \cP \rightarrow \X$ as in Proposition \ref{prop_bddextop}, and set $T_{(\hQ,\cA)} := T_{\hQ}$. 

We perform the one-time work of the algorithm \textsc{Compute New Extension Operator} (see Section \ref{sec_computingstuff}). We thus obtain a query algorithm for $T_\hQ$. Given $\underline{x} \in Q^\circ$, we can compute a short form description of the the $\Omega(\hQ)$-assisted bounded depth linear functional
\[
(f,P) \mapsto \partial^\beta \left[ J_{\underline{x}} T_{\hQ}(f,P) \right](\underline{x}) \;\;\;\; \mbox{for every} \; \beta \in \cM.
\]
This computation requires work at most $C \log N$ per query point.

\label{pp21}Proposition \ref{prop_bddextop} states that $T_\hQ(f,P) = f$ on $(1+a(\cA)) \hQ \cap E$, and
\[\| T_\hQ(f,P) \|_{\X((1+a(\cA))\hQ)} + \| T_\hQ(f,P) - P \|_{L^p((1+a(\cA))\hQ)} \leq C \cdot M_\hQ(f,P)
\]
for any $(f,P) \in \X(\frac{65}{64} \hQ \cap E) \oplus \cP$, where
\[M_{\hQ}(f,P) = \left(\sum_{\xi \in \Xi(\hQ)} \lvert \xi(f,P) \rvert^p \right)^{1/p}.\]
This proves \textbf{(E1)} and \textbf{(E2)} in the Main Technical Results for $\cA$.

We have thus treated all the non-simple cubes in $\CZ_{\main}(\cA)$.

\subsection{Case II: Simple cubes} \label{simple_sec} We suppose that $\hQ \in \CZ_{\main}(\cA)$ \underline{is} $c_*$-simple. We will explain how to construct the objects in the Main Technical Results for $\cA$ relevant to $\hQ$.

We have computed lists $\Omega(Q,\cA^-)$ and $\Xi(Q,\cA^-)$ of linear functionals on $\X(E \cap (65/64)Q)$ and $\X(E \cap (65/64)Q) \oplus \cP$, respectively, for each $Q \in \CZ_{\main}(\cA^-)$. See the Main Technical Results for $\cA^-$. Each functional in $\Xi(Q,\cA^-)$ has $\Omega(Q,\cA^-)$-assisted bounded depth and is given in short form.

From \eqref{mdefn} and \eqref{n_appx}, we know that
\begin{equation}\label{simple_1}
M_{(Q,\cA^-)}(f,R) := \left(\sum_{\xi \in \Xi(Q,\cA^-)} \lvert \xi(f,R) \rvert^p \right)^{1/p}
\end{equation}
satisfies
\begin{equation}\label{simple_2}
c \cdot \| (f,R) \|_{(1+ a) Q} \leq M_{(Q,\cA^-)}(f,R) \leq C \cdot \| (f,R)\|_{\frac{65}{64}Q}.
\end{equation}
Here, $a : = a(\cA^-) \in (0,1/64]$ is a universal constant in the Main Technical Results for $\cA^-$.

Recall that we have fixed a universal constant $t_G \in (0, 1/64]$ satisfying \eqref{small_t}.

We define
\begin{align*}
&\Omega(\hQ,\cA) := \bigcup_{\substack{Q \in \CZ_{\main}(\cA^-) \\ Q \subset (1+t_G)\hQ}} \Omega(Q,\cA^-).\\
&\Xi(\hQ,\cA) := \bigcup_{\substack{Q \in \CZ_{\main}(\cA^-) \\ Q \subset (1+t_G)\hQ}} \Xi(Q,\cA^-).\\
\end{align*}

Each $Q \in \CZ_{\main}(\cA^-)$ participates above for at most $C$ distinct $\hQ \in \CZ_{\main}(\cA)$. This is a consequence of the bounded overlap of $\{ \frac{65}{64}\hQ : \hQ \in \CZ_{\main}(\cA)\}$, since $t_G \leq \frac{1}{64}$. We can thus compute the lists $\Omega(\hQ,\cA)$ for all $c_*$-simple cubes $\hQ \in \CZ_{\main}(\cA)$, using work at most
\[
C \cdot \sum_{Q \in \CZ_{\main}(\cA^-)} \sum_{\omega \in \Omega(Q,\cA^-)} \depth(\omega) \leq CN,
\]
and we can compute the lists $\Xi(\hQ,\cA)$ for all $c_*$-simple cubes $\hQ \in \CZ_{\main}(\cA)$, using work at most
\[
C \cdot \sum_{Q \in \CZ_{\main}(\cA^-)} \biggl\{ 1 + \# \bigl[ \Xi(Q,\cA^-) \bigr] \biggr\} \leq CN.
\]
(The upper bound by $CN$ on these sums is stated in the Main Technical Results for $\cA^-$.) We do not attempt to remove duplicates from the lists $\Omega(\hQ,\cA)$ and $\Xi(\hQ,\cA)$, which are computed simply by copying.

When we copy the functionals in the list $\Omega(Q,\cA^-)$, for $Q \in \CZ_{\main}(\cA^-)$, $Q \subset (1+t_G)\hQ$, into the list $\Omega(\hQ,\cA)$, we mark each functional in $\Omega(Q,\cA^-)$ ($Q \in \CZ_{\main}(\cA^-)$, $Q \subset (1+t_G)\hQ$) with a pointer to its position in the list $\Omega(\hQ,\cA)$. This requires total extra work at most $CN$.

Each functional $\xi \in \Xi(\hQ,\cA)$ has $\Omega(Q,\cA^-)$-assisted bounded depth for some $Q \in \CZ_{\main}(\cA^-)$ with $Q \subset (1+t_G)\hQ$, hence $\xi$ has $\Omega(\hQ,\cA)$-assisted bounded depth, because $\Omega(Q,\cA^-)$ is a sublist of $\Omega(\hQ,\cA)$. We can compute a short form of $\xi$ in terms of the assists $\Omega(\hQ,\cA)$ by using the pointers from $\Omega(Q,\cA^-)$ into $\Omega(\hQ,\cA)$ (see Remark \ref{remk_translate}). This requires a constant amount of work per functional $\xi$. We assume that this work was carried out when we formed the lists $\Xi(\hQ,\cA)$.

We fix $\hQ \in \CZ_{\main}(\cA)$ such that $\hQ$ is $c_*$-simple.

As in the Main Technical Results for $\cA$, we define
\begin{equation}\label{drum1} \left[ M_{(\hQ,\cA)}(f,P) \right]^p  = \sum_{\xi \in \Xi(\hQ,\cA)} \lvert \xi(f,P) \rvert^p =  \sum_{\substack{ Q \in \CZ_{\main}(\cA^-) \\  Q \subset (1+t_G)\hQ} } \left[ M_{(Q,\cA^-)}(f,P) \right]^p.
\end{equation}

We next define an extension operator $T_{(\hQ,\cA)} : \X(E \cap (65/64)\hQ) \oplus \cP \rightarrow \X$. We follow an argument in Section \ref{sec_computingstuff}. 

We define the \emph{covering cubes}
\begin{equation*} \cov(\hQ) := \bigl\{ Q \in \CZ(\cA^-) : Q \subset (1+t_G)\hQ\bigr\}.\end{equation*}

Thanks to our assumption \eqref{small_t}, we can choose a universal constant $a_\new = a_\new(t_G)$ satisfying the conclusion of Lemma \ref{lem_cover}. Hence, since $\hQ$ is a testing cube, we obtain the following \\
\noindent\textbf{Covering Property}: The cube $(1+a_\new)\hQ$ is contained in the union of the cubes $(1+a/2) Q$ over all $Q \in \CZ(\cA^-)$ such that $Q \subset (1+t_G)\hQ$.

Recall that we have defined $a(\cA) = a_\new$ in \eqref{newa_defn}.

We pick cutoff functions $\theta_Q^\hQ \in C^m(\R^n)$, for each $Q \in \cov(\hQ)$, with 
\begin{equation}
\label{pou_new}
\left\{
\begin{aligned} 
& \sum_{Q \in \cov(\hQ)}  \theta^\hQ_Q  = 1 \; \mbox{on} \; (1+a_\new)\hQ,\\
&\supp ( \theta_Q^\hQ) \subset (1+a)Q \;\; \mbox{and} \;\;  \lvert \partial^\alpha \theta^\hQ_Q \rvert \leq C \cdot \delta_Q^{-|\alpha|} \;\; \mbox{for} \; |\alpha| \leq m, \; \mbox{and} \\
& \theta_Q^\hQ = 1 \; \mbox{near} \; x_Q, \; \mbox{and} \; \theta_Q^\hQ = 0 \; \mbox{near} \; x_{Q'} \; \mbox{for each} \; Q' \in \cov(\hQ) \setminus \{Q\}.
\end{aligned}
\right.
\end{equation}

For each $Q \in \cov(\hQ)$ we define
\begin{equation}\label{simple_lf}
   F^{\hQ}_Q := 
   \left\{
     \begin{array}{lr}
       T_{(Q,\cA^-)}(f,P) & : \mbox{if} \; \frac{65}{64}Q \cap E \neq \emptyset  \\
       P & :  \mbox{if} \; \frac{65}{64}Q \cap E = \emptyset.
     \end{array}
   \right.
\end{equation}
We define a linear map $T_{(\hQ,\cA)} : \X(E \cap \frac{65}{64}\hQ) \oplus \cP \rightarrow \X$ by the formula
\begin{equation}\label{simple_extop}
T_{(\hQ,\cA)}(f,P) := \sum_{Q \in \cov(\hQ)} F^\hQ_Q \cdot \theta^\hQ_Q.
\end{equation}
(Compare to \eqref{extopdefn}.)

Here, the maps $T_{(Q,\cA^-)}$ are as in the Main Technical Results for $\cA^-$; see Chapter \ref{sec_mainresults}. Each $T_{(Q,\cA^-)}$ has $\Omega(Q,\cA^-)$ assisted bounded depth, hence  $T_{(Q,\cA^-)}$ has $\Omega(\hQ,\cA)$-assisted bounded depth, since by definition $\Omega(Q,\cA^-)$ is a subslist of $\Omega(\hQ,\cA)$ for each $Q \in \cov(\hQ)$. 

Therefore, each $T_{(\hQ,\cA)}$ has $\Omega(\hQ,\cA)$-assisted bounded depth. We also give a query algorithm for $T_{(\hQ,\cA)}$: Given $\underline{x} \in Q^\circ$, we compute the map $(f,P) \mapsto J_{\underline{x}} T_{(\hQ,\cA)}(f,P)$ in short form in terms of the assists $\Omega(\hQ,\cA)$. We leave details to the reader.

\begin{prop}\label{simple_extprops}
Let $(f,P) \in \X(\frac{65}{64}\hQ \cap E) \oplus \cP$. Then the following properties hold.
\begin{itemize}
\item $T_{(\hQ,\cA)}(f,P) = f$ on $(1+a_\new)\hQ \cap E$.
\item $\|T_{(\hQ,\cA)}(f,P) \|_{\X( (1+a_\new)\hQ)} + \delta_\hQ^{-m}  \| T_{(\hQ,\cA)}(f,P) - P \|_{L^p((1+a_\new)\hQ)} \leq C \cdot M_{(\hQ, \cA)}(f,P)$.
\end{itemize}
\end{prop}
\begin{proof}
The proof is analogous to the proof of Proposition \ref{prop_bddextop}, except much easier. We spell out the details.

For ease of notation, we set $\overline{a}  = a_\new$.

The definition of the linear map in \eqref{simple_extop} is the same as that in \eqref{extopdefn}, except that the polynomials $R_{Q}^\hQ$ used in the functions $F_Q^\hQ$ in \eqref{test_eo_aux} are replaced by $P$ (compare \eqref{test_eo_aux} and \eqref{simple_lf}). Thus, to prove our proposition, we may follow parts of the reasoning in the proof of Proposition \ref{prop_bddextop}, as long as we substitute $R_Q^\hQ$ everywhere with $P$. 

The functions $F_Q^\hQ$ in \eqref{simple_lf} satisfy 
\begin{equation}
\label{new_local}
\left\{
\begin{aligned}
&F_Q^\hQ = f \; \mbox{on} \; (1+a)Q \cap E \\
& \| F_Q^\hQ \|_{\X((1+a)Q)} + \delta_Q^{-m} \| F_Q^\hQ - P \|_{L^p((1+a)Q)} \leq  \left\{
\begin{array}{lr}
       C M_{(Q,\cA^-)}(f,P) & : \mbox{if} \; \frac{65}{64}Q \cap E \neq \emptyset   \\
       0 & :  \mbox{if} \; \frac{65}{64}Q \cap E = \emptyset
     \end{array}
 \right.
 \end{aligned}
 \right.
\end{equation}
This follows from the Main Technical Results for $\cA^-$.

Thus, the function $T_{(\hQ,\cA)}(f,P)$ defined in \eqref{simple_extop} satisfies the first bullet point of Proposition \ref{simple_extprops}. This is a consequence of the first and second conditions in \eqref{pou_new}, and the first condition in \eqref{new_local}.

We now prove the second bullet point of Proposition \ref{simple_extprops}.

Let $G = T_{(\hQ,\cA)}(f,P)$.

The equation \eqref{midd} holds in the present setting if we replace $R_Q^\hQ$ with $P$, for the same reason as before. (Here, we use the \textbf{Covering Property}.) Moreover, when we replace $R_Q^\hQ$ with $P$, the term $A_2(f,P)$ vanishes. Thus, we have
\[\| G \|^p_{\X((1+\overline{a})\hQ)} \lesssim \sum_{\substack{ Q \in \cov(\hQ) \\ \frac{65}{64}Q \cap E \neq \emptyset} }   \left[ M_{(Q,\cA^-)}(f,P) \right]^p. 
\]
By definition, the right-hand side is equal to $\left[ M_{(\hQ,\cA)}(f,P) \right]^p$ (see \eqref{drum1}). Thus we have proven
\[
\| G \|_{\X((1+\overline{a})\hQ)} \leq C \cdot M_{(\hQ,\cA)}(f,P).
\]

It remains to show that $\| G - P \|_{L^p((1+\overline{a})\hQ)} \leq C \cdot M_{(\hQ,\cA)}(f,P)$. We proceed directly without referring to the previous arguments. Using  \eqref{simple_extop} and the first condition in \eqref{pou_new}, we have
\[
G - P = \sum_{Q \in \cov(\hQ)} \theta_Q^\hQ \cdot ( F_Q^\hQ - P ) \qquad \mbox{on} \; (1+\overline{a})\hQ.
\]
Recall that $\theta_Q^\hQ$ is supported on $(1+a) Q$ and $\lvert \theta_Q^\hQ \rvert \leq C$ (see \eqref{pou_new}). Since $\hQ$ is $c_*$-simple, at most $C$ cubes $Q$ contribute to the above sum, and $\delta_Q \geq c_* \delta_{\hQ}$ for each $Q$. Hence,
\[ (\delta_{\hQ})^{-m} \| G - P \|_{L^p((1+\overline{a})\hQ)}^p \leq C \sum_{Q \in \cov(\hQ)} (\delta_Q)^{-m} \| F_Q^\hQ - P \|^p_{L^p((1+a)Q)}.\]
Hence, using \eqref{new_local}, we have
\[ (\delta_{\hQ})^{-m} \| G - P \|_{L^p((1+\overline{a})\hQ)}^p \leq C \sum_{Q \in \cov(\hQ)}  \left[ M_{(Q,\cA^-)}(f,P) \right]^p = C \cdot \left[M_{(\hQ,\cA)}(f,P) \right]^p.\]
This completes the proof of the second bullet point in Proposition \ref{simple_extprops}.
\end{proof}

\begin{lem}\label{simple_lem1}
We have
\[
c \|(f,P) \|_{(1+a_\new)\hQ} \leq M_{(\hQ,\cA)}(f,P) \leq C \| (f,P) \|_{\frac{65}{64}\hQ}.
\]
\end{lem}
\begin{proof}
The inequality $ \|(f,P) \|_{(1+a_\new)\hQ} \leq C M_{(\hQ,\cA)}(f,P)$ is an easy consequence of Proposition \ref{simple_extprops} and the definition of the trace seminorm. Thus, the only task is to prove the second inequality, $M_{(\hQ,\cA)}(f,P) \leq C \| (f,P) \|_{\frac{65}{64}\hQ}$.

First, the upper bound in \eqref{simple_2} implies that
\[
\left[ M_{(\hQ,\cA)}(f,P) \right]^p \leq C \sum_{\substack{ Q \in \CZ_{\main}(\cA^-) \\  Q \subset (1+t_G)\hQ} } \| (f,P) \|_{\frac{65}{64} Q}^p.\]
Since $\hQ$ is $c_*$-simple, each cube $Q$ relevant to the above sum satisfies $\delta_Q \geq c_* \delta_\hQ$. Moreover, Lemma \ref{lem_cover} implies that $\frac{65}{64}Q \subset \frac{65}{64}\hQ$ for each relevant $Q$; recall \eqref{small_t}. Hence, each term $\| (f,P) \|_{\frac{65}{64} Q}$ is bounded by $C \| (f,P) \|_{\frac{65}{64} \hQ}$ thanks to Lemma \ref{lem_normmon}. Moreover, the number of terms is at most a universal constant, hence
\[M_{(\hQ,\cA)}(f,P) \leq C \| (f,P) \|_{\frac{65}{64} \hQ}.\]
This completes the proof of the lemma.
\end{proof}

We have produced lists $\Omega(\hQ,\cA)$ and $\Xi(\hQ,\cA)$, and we have defined a linear map $T_{(\hQ,\cA)}$ that satisfy the conditions in the Main Technical Results for $\cA$ (see Chapter \ref{sec_mainresults}), for every $\hQ \in \CZ_{\main}(\cA)$ that is $c_*$-simple. We have remarked that one can easily produce a query algorithm for $T_{(\hQ,\cA)}$. We have performed these computations using work at most $C N \log N$ in space $CN$.

We have thus  treated all the simple cubes in $\CZ_{\main}(\cA)$

\subsection{Closing remarks}
\label{closing_remarks}

All the previously defined objects satisfy the conditions set down in Chapter \ref{sec_mainresults} with many of the constants depending on $\epsilon$, and with $\epsilon_2(\cA) = \epsilon^\kappa$, $\epsilon_1(\cA) = \epsilon^{1/\kappa}$. We have computed a list of assists $\Omega(\hQ, \cA)$, and a list of assisted functionals $\Xi(\hQ,\cA)$, and we have given a query algorithm for a linear map $T_{(\hQ,\cA)}$ for each $\hQ \in \CZ_{\main}(\cA)$, using one-time work at most $C N \log N$ in space $CN$. In particular, the bound on the required space implies that
\[
\left\{
\begin{aligned}
 &\sum_{\hQ \in \CZ_{\main}(\cA)}  \sum_{\omega \in \Omega(\hQ,\cA)} \depth(\omega)  \leq CN, \;\; \mbox{and}  \\
 &\sum_{\hQ \in \CZ_{\main}(\cA)}  \#\bigl[ \Xi(\hQ,\cA)\bigr] \leq CN.
\end{aligned}
\right.
\]

We now fix $\epsilon$ to be a universal constant, small enough so that the previous results hold. That completes the Induction Step, and thus we have achieved the Main Technical Results for $\cA$.

\chapter{Proofs of the Main Theorems}
\label{final_chap}

\section{Extension in Homogeneous Sobolev Spaces}\label{hom_sec}

In this section we prove our main theorem concerning homogeneous Sobolev spaces $\X = L^{m,p}(\R^n)$ ($p > n$), which reads as follows.

\begin{thm} \label{main_thm_hom}
Let $E \subset \R^n$ satisfy $N = \#(E) \geq 2$.
\begin{itemize}
\item We produce lists $\Omega$ and $\Xi$, consisting of functionals on $\X(E) = \{ f : E \rightarrow \R\}$, with the following properties.
\begin{itemize}
\item The sum of $\depth(\omega)$ over all $\omega \in \Omega$ is bounded by $C N$. The number of functionals in $\Xi$ is at most $C N$.
\item Each functional $\xi$ in $\Xi$ has $\Omega$-assisted depth at most $C$. The functionals in $\Omega$ and $\Xi$ are represented in their short form.
\item For all $f \in \X(E)$ we have
\[ c \| f\|_{\X(E)} \leq \left[ \sum_{\xi \in \Xi} \lvert \xi(f) \rvert^p  \right]^{1/p}  \leq C \| f\|_{\X(E)}.\]
\end{itemize}
\end{itemize}

Moreover, there exists a linear map $T : \X(E) \rightarrow \X$ with the following properties.
\begin{itemize}
\item $T$ has $\Omega$-assisted depth at most $C$.
\item $T f = f$ on $E$ and $\| Tf \|_{\X} \leq C \| f \|_{\X(E)}$ for all $f \in \X(E)$.
\item We produce a query algorithm that operates as follows.

Given a point $\underline{x} \in \R^n$, we compute a short form description of the $\Omega$-assisted bounded depth linear map $ \X(E) \ni f \mapsto J_{\underline{x}}\left( Tf  \right) \in \cP$ using work and storage at most $C \log N$.
\end{itemize}

The  computations above require one-time work at most $C N \log N$ in space $C N$.
\end{thm}

By translating and rescaling, we may assume without loss of generality that $E \subset \frac{1}{32}Q^\circ$, with $Q^\circ = [0,1)^n$.

\label{pp22}
We deduce Theorem \ref{main_thm_hom} from the Main Technical Results for $\cA = \emptyset$. Recall that we have achieved the following (see Chapter \ref{sec_mainresults}).

\begin{itemize}
\item There is a decomposition $\CZ$ of $Q^\circ$ into dyadic cubes. Every point $x \in Q^\circ$ belongs to a unique cube $Q_x \in \CZ$.
\item We produce a $\CZ$-\textsc{Oracle}.

The $\CZ$-\textsc{Oracle} accepts a query point $\underline{x} \in Q^\circ$. The response to a query $\underline{x}$ is the list of all $Q \in \CZ$ such that $\frac{65}{64} Q$ contains $\underline{x}$. The work and storage required to answer a query are at most $C \log N$.
\item If $Q,Q' \in \CZ$ and $Q \leftrightarrow Q'$ then $\frac{1}{2} \delta_Q \leq \delta_{Q'} \leq 2 \delta_Q$.
\item Each point $x \in \R^n$ is contained in at most $C$ of the cubes $\frac{65}{64}Q$, $Q \in \CZ$.\\
(This is an easy consequence of the previous bullet point.)
\item If $Q \in \CZ$ and $\delta_Q \leq c_*$ then $S Q$ is not tagged with $(\emptyset, \epsilon_1)$.

Recall that any cube is tagged with $(\emptyset,\epsilon_1)$; see Remark \ref{tag_rem}. Thus, we learn that $\delta_Q > c_* $ for each $Q \in \CZ$. In particular, the cardinality of $\CZ$  is bounded by some universal constant $C$.

\item Next, we recall the various assists, functionals and local extension operators described in the Main Technical Results.

For each $Q \in \CZ$ with $\frac{65}{64}Q \cap E \neq \emptyset$, we compute a list of assists $\Omega(Q)$  and assisted functionals $\Xi(Q) \subset \left[ \X(E \cap \frac{65}{64}Q) \oplus \cP\right]^*$. Each $\xi$ in $\Xi(Q)$ has $\Omega(Q)$-assisted bounded depth. We have
\begin{equation}
\label{mm0}
\sum_{\xi \in \Xi(Q)} \lvert \xi(f,P) \rvert^p \leq C \cdot \|(f,P)\|_{\frac{65}{64}Q}^p.
\end{equation}
We compute these lists of functionals using one-time work at most $C N \log N$ in space $C N$.

We also define an $\Omega(Q)$-assisted bounded depth linear map $T_Q : \X(\frac{65}{64}Q \cap E) \oplus \cP \rightarrow \X$ such that
\begin{equation}
\label{m1} T_Q(f,P) = f \quad \mbox{on} \; E \cap (1+a)Q
\end{equation}
and
\begin{equation}
\label{m2} \| T_Q(f,P) \|_{\X((1+a)Q)}^p + \delta_Q^{-mp} \| T_Q(f,P) - P \|_{L^p((1+a)Q)}^p \leq C \cdot \sum_{\xi \in \Xi(Q)} \lvert \xi(f,P) \rvert^p.
\end{equation}
Given a query $\underline{x} \in Q^\circ$, we can compute the linear map $(f,P) \mapsto J_{\underline{x}}T_Q(f,P)$ in short form in terms of the assists $\Omega(Q)$, using work at most $C \log N$.

\end{itemize}

This completes the description of the objects from Chapter \ref{sec_mainresults}.

We list the cubes in $\CZ$ with the following procedure. For each $x \in (c_*/10) \Z^n \cap Q^\circ$, we use the $\CZ$-\textsc{Oracle} to list all the cubes $Q$ in $\CZ$ such that $x \in \frac{65}{64}Q$. Each $Q\in\CZ$ contains at least one point in $(c_*/10) \Z^n\cap Q^\circ$ (because $Q \subset Q^\circ$ and $\delta_Q > c_*$), hence each $Q \in \CZ$ arises in an aforementioned list for some $x$. We concatenate these lists and then sort the resulting list to remove duplicate cubes. 

We now construct a suitable partition of unity adapted to the decomposition $\CZ$.

Let $a :=a(\cA)$ with $\cA = \emptyset$, as in Chapter \ref{sec_mainresults}. Recall that
\begin{equation}
\label{star1}
0 < a \leq 1/64.
\end{equation}

For each $Q \in \CZ$, let $\widetilde{\theta}_Q \in C^m(\R^n)$ be a function such that
\begin{enumerate}
\item $0 \leq \widetilde{\theta}_Q \leq 1$ on $\R^n$, 
\item $\widetilde{\theta}_Q \geq 1/2$ on $Q$,
\item $\widetilde{\theta}_Q = 0$ outside $ (1+a)Q$,
\item $\lvert \partial^\beta \widetilde{\theta}_Q(x) \rvert \leq C$ for $x \in \R^n$, $\lvert \beta \rvert \leq m$.
\end{enumerate}
Also, let $\eta : [0,\infty) \rightarrow \R$ be a $C^m$ function such that
\begin{enumerate}
\item $\eta(t) \geq 1/4$ for $t \geq 0$,
\item $\eta(t) = t$ for $t \geq 1/2$,
\item $\lvert (d/dt)^k \eta(t) \rvert \leq C$ for $t \geq 0$, $k \leq m$.
\end{enumerate} 
We can satisfy these conditions by choosing $\widetilde{\theta}_Q$ and $\eta$ to be appropriate spline functions. We assume that the following queries can be answered using work at most $C$.

\environmentA{Algorithm: Compute Auxiliary Functions.} Given $Q \in \CZ$ and $\underline{x} \in \R^n$, we can compute the jet $J_{\underline{x}}(\widetilde{\theta}_Q)$. Given $t_* \geq 0$ and an integer $0 \leq k \leq m$, we can compute $\frac{d^k \eta}{d t^k}(t_*)$.

For each $Q \in \CZ$ we define
\[\theta_Q(x) = \frac{\widetilde{\theta}_Q(x)}{\eta \circ \Psi(x)}, \;\;\; \mbox{where} \; \Psi(x) = \sum_{Q \in \CZ} \widetilde{\theta}_Q(x).\]
Clearly, we can answer the following query using work at most $C$.

\environmentA{Algorithm: Compute POU2.} Given  $Q \in \CZ$ and $\underline{x} \in \R^n$, we compute the jet $J_{\underline{x}} (\theta_Q)$.

We can easily prove the following properties. (See the proof of Lemma \ref{pou_lem}.)
\begin{enumerate}
\item $\theta_Q \in C^m(\R^n)$ is well-defined, by property (1) of $\widetilde{\theta}_Q$ and property (1) of $\eta$.
\item $ \theta_Q = 0$ outside $(1+a)Q$, by property (3) of $\widetilde{\theta}_Q$.
\item $\lvert \partial^\alpha \theta_Q(x) \rvert \leq C$ for $x \in \R^n$, $| \alpha | \leq m$, by property (4) of $\widetilde{\theta}_Q$ and properties (1),(3) of $\eta$.
\item $\displaystyle \sum_{Q \in \CZ} \theta_Q = 1$ on $Q^\circ$.
\end{enumerate}
To prove property (4), recall that the cubes in $\CZ$ cover $Q^\circ$. Hence, properties (1) and (2) above imply that $ \Psi(x) = \sum_{Q \in \CZ} \widetilde{\theta}_Q(x) \geq 1/2$ for $x \in Q^\circ$. Hence, property (2) of $\eta$ implies that $\eta \circ \Psi(x) = \Psi(x)$ for $x \in Q^\circ$, which implies that
\[\sum_{Q \in \CZ} \theta_Q(x) = \frac{\sum_{Q \in \CZ} \widetilde{\theta}_Q(x)}{\Psi(x)} = 1 \quad \mbox{for} \; x \in  Q^\circ.\]
This completes the proof of property (4) of $\{\theta_Q\}_{Q \in \CZ}$.

We let $\Xi^\circ \subset (\X(E) \oplus \cP)^*$ be the union of the  lists $\Xi(Q)$ for all $Q \in \CZ$ such that $E \cap \frac{65}{64} Q\neq \emptyset$.  Similarly, we let $\Omega^\circ$ be the union of the lists $\Omega(Q)$ for all $Q \in \CZ$ such that $E \cap \frac{65}{64} Q\neq \emptyset$. Hence,
\[\sum_{\xi^\circ \in \Xi^\circ} \lvert \xi^\circ (f,P) \rvert^p = \sum_{Q \in \CZ} \sum_{\xi \in \Xi(Q)}\lvert \xi(f,P) \rvert^p.\]
The functionals in $\Xi^\circ$ have $\Omega^\circ$-assisted bounded depth. We can express each functional $\xi$ in $\Xi^\circ$ in short form in terms of the assists $\Omega^\circ$ by sorting. This requires work at most $C \log N$ for each $\xi \in \Xi^\circ$. Because there are at most $C N$ functionals in $\Xi^\circ$, this requires work at most $C N \log N$ in total.

Given $(f,P) \in \X(E) \oplus \cP$ we define
\begin{equation}
\label{star2}
T^\circ(f,P)  = \sum_{Q \in \CZ} \theta_Q \cdot F_Q,
\end{equation}
where $F_Q = T_Q(f,P)$ whenever $\frac{65}{64}Q \cap E \neq \emptyset$, and $F_Q = P$ whenever $\frac{65}{64}Q \cap E = \emptyset$.

\begin{prop}\label{mprop}
The following hold.
\begin{itemize}
\item The sum of $\depth(\omega)$ over all $\omega \in \Omega^\circ$ is bounded by $C N$. \\
The cardinality of $\Xi^\circ$ is bounded by $C N$.
\item Given $\underline{x} \in Q^\circ$, we can compute a short form description of the $\Omega^\circ$-assisted bounded depth linear map 
\[  \X(E) \oplus \cP  \ni (f,P)  \mapsto J_{\underline{x}} T^\circ(f,P) \in \cP\]
using work and storage at most $C \log N$.
\item  Given $(f,P) \in \X(E) \oplus \cP$ we have $T^\circ(f,P) = f$ on $E$, and
\[ \| T^\circ(f,P) \|_{\X(Q^\circ)}^p + \| T^\circ(f,P) - P \|^p_{L^p(Q^\circ)} \leq C \sum_{\xi \in \Xi^\circ} \lvert \xi(f,P) \rvert^p.\]
\item Given $(f,P) \in \X(E) \oplus \cP$, we have
\[ \sum_{\xi \in \Xi^\circ} \lvert \xi(f,P) \rvert^p \leq C \cdot \| (f,P) \|_{\frac{65}{64}Q^\circ}^p.\]
\end{itemize}
\end{prop}

\begin{proof}

From Chapter  \ref{sec_mainresults}, recall that 
\[ \sum_{\substack{Q \in \CZ \\ \frac{65}{64}Q \cap E \neq \emptyset}} \sum_{\omega \in \Omega(Q)} \depth(\omega) \leq C N \quad \mbox{and} \quad \sum_{\substack{ Q \in \CZ \\ \frac{65}{64}Q \cap E \neq \emptyset } }\# \bigl[ \Xi(Q) \bigr] \leq C N. \]
This implies the conclusion of the first bullet point.

We fix a query point  $\underline{x} \in Q^\circ$. Then we have
\begin{equation} \label{jet_form}
J_{\underline{x}} T^\circ(f,P) =  \sum_{\substack{Q \in \CZ \\ \frac{65}{64}Q \cap E \neq \emptyset}} J_{\underline{x}} \theta_Q \odot_{\underline{x}} J_{\underline{x}}  T_Q(f,P)  +  \sum_{\substack{Q \in \CZ \\ \frac{65}{64}Q \cap E = \emptyset}}  J_{\underline{x}} \theta_Q \odot_{\underline{x}} P.
\end{equation}
Recall that we have computed a list of all the $Q \in \CZ$, and that there are at most $C$ such cubes.  We loop over all $Q \in \CZ$, and perform the steps below.
\begin{itemize}
\item \textbf{Step 1}: For each $\alpha \in \cM$ we compute $\partial^\alpha( J_{\underline{x}} \theta_Q)(\underline{x})$ using \textsc{Compute POU2}.
\item \textbf{Step 2}: If $\frac{65}{64}Q \cap E \neq \emptyset$, then we compute the linear map $(f,P) \mapsto J_{\underline{x}} T_Q(f,P)$ in short form in terms of the assists $\Omega(Q)$ (see the Main Technical Results for $\cA = \emptyset$). This means that for each $\alpha \in \cM$ we compute linear functionals $\lambda^{Q,\alpha} : \cP \rightarrow \R$ and $\eta^{Q, \alpha} : \X(E) \rightarrow \R$, assists $\omega^{Q,\alpha}_1,\cdots,\omega^{Q,\alpha}_d \in \Omega(Q)$, and numbers $\gamma^{Q,\alpha}_1,\cdots, \gamma^{Q,\alpha}_d \in \R$, such that 
\[
\partial^\alpha( J_{\underline{x}} T_Q(f,P))(0) = \lambda^{Q,\alpha}(P) + \eta^{Q,\alpha}(f) + \sum_{k=1}^d \gamma_k^{Q,\alpha} \cdot \omega^{Q,\alpha}_k(f).
\]
We guarantee that $\depth(\eta^{Q,\alpha}) + d$ is at most a universal constant $C$. 

From the first bullet point in Proposition \ref{mprop}, we know that $\Omega^\circ = \cup_{Q \in \CZ} \Omega(Q)$ contains at most $CN$ functionals. When we formed the list $\Omega^\circ$ by concatenation, we assume that we marked each functional in $\Omega(Q)$ with a pointer to its position in the list $\Omega^\circ$.  This requires additional one-time work work at most $CN$. Thus, the previous formula gives a short form representation of the functional $(f,P) \mapsto \partial^\alpha( J_{\underline{x}} T_Q(f,P))(0)$ in terms of the assists $\Omega^\circ$. Therefore, by Taylor's theorem we can compute a short form representation of the functional $(f,P) \mapsto \partial^\alpha( J_{\underline{x}} T_Q(f,P))(\underline{x})$ in terms of the assists $\Omega^\circ$.

From the definition of the product $\odot_{\underline{x}}$ and the computation in \textbf{Step 1}, for each $\alpha \in \cM$ we can compute a short form of the functional 
\[
(f,P) \mapsto \partial^\alpha( J_{\underline{x}} \theta_Q \odot_{\underline{x}} J_{\underline{x}} T_Q(f,P) )(\underline{x})
\]
in terms of the assists $\Omega^\circ$.

\item \textbf{Step 3}: If $\frac{65}{64} Q \cap E = \emptyset$,  then for each $\alpha \in \cM$ we compute a short form of the functional $(f,P) \mapsto \partial^\alpha ( J_{\underline{x}} \theta_Q \odot_{\underline{x}} P )(\underline{x})$. Here, we use Taylor's theorem to compute the change-of-coordinate map $(\partial^\alpha P(0))_{\alpha \in \cM} \mapsto (\partial^\alpha P(\underline{x}))_{\alpha \in \cM}$. Thus, the desired computation is a consequence of the definition of the product $\odot_{\underline{x}}$ and the result of \textbf{Step 1}. This concludes the loop over $Q$.
\end{itemize}
For each $\alpha \in \cM$, we compute a short form of the functional $(f,P) \mapsto \partial^\alpha(J_{\underline{x}} T^\circ(f,P) )(\underline{x})$ in terms of the assists $\Omega^\circ$ by adding together the short form representations of the functionals determined at the end of \textbf{Step 2} and \textbf{Step 3} (see the formula \eqref{jet_form}). Therefore, we can compute a short form of the functional $(f,P) \mapsto \partial^\alpha(J_{\underline{x}} T^\circ(f,P) )(0)$ in terms of the  assists $\Omega^\circ$. This is a consequence of Taylor's theorem and the previous computation. The reader may easily check that the above computation requires work at most $C \log N$ per query $\underline{x} \in Q^\circ$. This completes the proof of the second bullet point in  Proposition \ref{mprop}.

Fix $x \in E$. Then
\begin{equation}
\label{m3}
T^\circ(f,P)(x) = \sum_{\substack{Q \in \CZ \\ \frac{65}{64}Q \cap E \neq \emptyset}} \theta_Q(x) \cdot T_Q(f,P)(x)  +  \sum_{\substack{Q \in \CZ \\ \frac{65}{64}Q \cap E = \emptyset}}  \theta_Q(x) \cdot P(x).
\end{equation}
Recall that $\theta_Q$ is supported on the cube $(1+a)Q$, which is contained in $\frac{65}{64}Q$. (See \eqref{star1}.)

For the $Q$ arising in the second sum in \eqref{m3} we learn that $\theta_Q(x) = 0$, since the support of $\theta_Q$ does not intersect $E$. 

For the $Q$ arising in the first sum in \eqref{m3}, if $x \in (1+a)Q$ then $T_Q(f,P)(x) = f(x)$. Otherwise, if $x \notin (1+a)Q$ then $\theta_Q(x) = 0$, due to the support properties of $\theta_Q$.

Hence, $T^\circ(f,P)(x) = \sum_{Q \in \CZ} \theta_Q(x) f(x) = f(x)$, due to the fact that $\sum_{Q \in \CZ} \theta_Q = 1$ on $E$ (recall that $E \subset Q^\circ$). Hence, $T^\circ(f,P) = f$ on $E$, as desired.

We apply Lemma \ref{patch_lem}, where $P_Q = P$ and $F_Q$ ($Q \in \CZ$) is determined just below \eqref{star2}. The conditions on $\CZ$ in Section \ref{czalg_sec} are given in the Main Technical Results in Chapter \ref{sec_mainresults}. Here, for $\hQ = Q^\circ$, the conditions \eqref{covers} and \eqref{sizebd} in Section \ref{sec_pou} are obvious (since $Q^\circ$ equals the union of all $Q \in \CZ$; also, $\delta_{Q^\circ} = 1$  and $\delta_Q \leq 1$ for all $Q \in \CZ$). Thus, by Lemma \ref{patch_lem}, we have
\begin{align}
\label{m4}
\|T^\circ(f,P) \|_{\X(Q^\circ)}^p & \lesssim  \sum_{Q \in \CZ} \left[ \| F_Q \|_{\X((1+a)Q)}^p + \delta_Q^{-mp}\| F_Q - P \|_{L^p((1+a)Q)}^p \right] \\
\notag{}
& =  \sum_{\substack{ Q \in \CZ \\ \frac{65}{64}Q \cap E \neq \emptyset}} \left[ \| T_Q(f,P) \|_{\X((1+a)Q)}^p + \delta_Q^{-mp} \|T_Q(f,P) - P \|^p_{L^p((1+a)Q)} \right] \\
\notag{} 
&\lesssim \sum_{\substack{ Q \in \CZ \\ \frac{65}{64}Q \cap E \neq \emptyset}}  \sum_{\xi \in \Xi(Q)} \lvert \xi(f,P) \rvert^p \quad (\mbox{see \eqref{m2}}).
\end{align}
(All the terms with $\frac{65}{64}Q \cap E = \emptyset$ vanish, since $F_Q = P$ is an $(m-1)$-st degree polynomial.)

Since $ \sum_{Q \in \CZ} \theta_Q = 1$ on $Q^\circ$, we have
\[T^\circ(f,P) - P  = T^\circ(f,P) -  \sum_{Q \in \CZ} \theta_Q \cdot P = \sum_{\substack{ Q \in \CZ \\ \frac{65}{64}Q \cap E \neq \emptyset}}  \left( T_Q(f,P) - P \right) \cdot \theta_Q \;\; \mbox{on} \; Q^\circ.\]
There are at most $C$ terms in the above sum. Thus, since each $\theta_Q$ is supported on $(1+a)Q$ and $\| \theta_Q \|_{L^\infty} \leq C$, we have
\begin{equation}\label{m5} \| T^\circ(f,P) - P \|^p_{L^p(Q^\circ)} \lesssim \sum_{\substack{ Q \in \CZ \\ \frac{65}{64}Q \cap E \neq \emptyset}} \| T_Q(f,P) - P \|_{L^p((1+a)Q)}^p \lesssim \sum_{\substack{ Q \in \CZ \\ \frac{65}{64}Q \cap E \neq \emptyset}}  \sum_{\xi \in \Xi(Q)} \lvert \xi(f,P) \rvert^p.
\end{equation}
Here, in the last inequality we used \eqref{m2}. (Recall that $\delta_Q \leq 1$ whenever $Q \in \CZ$.)

Summing \eqref{m4} and \eqref{m5} shows that
\[ \| T^\circ(f,P) \|_{\X(Q^\circ)}^p + \| T^\circ(f,P) - P \|_{L^p(Q^\circ)}^p \lesssim \sum_{\substack{ Q \in \CZ \\ \frac{65}{64}Q \cap E \neq \emptyset}}  \sum_{\xi \in \Xi(Q)} \lvert \xi(f,P) \rvert^p.\]
The right-hand expression is equal to $\sum_{\xi \in \Xi^\circ } \lvert \xi(f,P) \rvert^p$. This completes the proof of the third bullet point in Proposition \ref{mprop}.

\label{pp23}
From \eqref{mm0}, recall that
\begin{equation}\label{snuff22}
 \sum_{\substack{ Q \in \CZ \\ \frac{65}{64}Q \cap E \neq \emptyset}}  \sum_{\xi \in \Xi(Q)} \lvert \xi(f,P) \rvert^p \leq C \sum_{\substack{ Q \in \CZ \\ \frac{65}{64}Q \cap E \neq \emptyset}} \| (f,P) \|_{\frac{65}{64}Q}^p.
 \end{equation}

We have $\| (f,P )\|_{\frac{65}{64}Q} \lesssim \|(f,P)\|_{\frac{65}{64}Q^\circ}$ for each $Q \in \CZ$. Here, we apply Lemma \ref{lem_normmon} and use the fact that $\delta_Q \geq c_*$ for all $Q \in \CZ$. Since the cardinality of $\CZ$ is bounded by a universal constant, we conclude that
\[ \sum_{\substack{ Q \in \CZ \\ \frac{65}{64}Q \cap E \neq \emptyset}}  \sum_{\xi \in \Xi(Q)} \lvert \xi(f,P) \rvert^p \leq C \cdot \| (f,P)\|^p_{\frac{65}{64}Q^\circ}.\]
This implies the fourth bullet point in Proposition \ref{mprop}. This completes the proof of Proposition \ref{mprop}.

\end{proof}

We will now construct the various assists, functionals, and the extension operator from Theorem \ref{main_thm_hom}.

\noindent\textbf{\underline{Computing a near-optimal jet:}}

Each functional $\xi_\ell \in \Xi^\circ$ is given in the form
\begin{align}
\label{xi_short2}
\xi_\ell(f,R) &= \lambda_\ell(f) + \sum_{a=1}^{I_\ell} \mu_{\ell a} \omega_{\ell a}(f) + \sum_{\alpha \in \cM} \check{\mu}_{\ell \alpha} \cdot  \partial^\alpha R (0) \\
& \qquad\qquad\qquad \mbox{for} \; \ell =1,\cdots,L; \; \mbox{here,} \; L = \#(\Xi^\circ) \leq CN. \notag{}
\end{align}
Here, $\omega_{\ell a} \in \Omega^\circ$; $\lambda_\ell$ is a linear functional; $\mu_{\ell a}$ and $\check{\mu}_{\ell \alpha}$ are real coefficients; and $\depth(\lambda_\ell) = \mathcal{O}(1)$, $I_\ell = \mathcal{O}(1)$. In this discussion, we write $X = \mathcal{O}(Y)$ to indicate that $X \leq C Y$ for a universal constant $C$.

Applying the algorithm \textsc{Optimize via Matrix}, we find a matrix $(b_{\alpha \ell})_{\substack{ \alpha \in \cM \\ \ell = 1,\cdots, L }}$ such that the sum of the $p$-th powers of the $\lvert \xi_\ell(f,R) \rvert$ ($\ell=1,\cdots,L$) in \eqref{xi_short2} is essentially minimized for fixed $f$ by setting
\begin{align}
\label{min_func}
\partial^\alpha R(0) &= \sum_{\ell=1}^L b_{\alpha \ell} \left[ \lambda_\ell(f) + \sum_{a=1}^{I_\ell} \mu_{\ell a} \omega_{\ell a}(f) \right] \\
& \equiv  \qquad\qquad \omega^{\new}_\alpha(f).
\notag{}
\end{align}

We express the functionals $\omega^\new_\alpha$ in short form. We first compute real coefficients $(\mu_{\alpha x})_{x \in E}$ ($\alpha \in \cM$) so that
\begin{equation}
\label{min_func2}
\omega^\new_\alpha(f) = \sum_{x \in E} \mu_{\alpha x} \cdot f(x).
\end{equation}
We achieve this by summing all the coefficients $b_{\alpha \ell} \cdot \mu_{\ell a}$ in \eqref{min_func} that correspond to the same functional $\omega = \omega_{\ell a}$. (We accomplish this by sorting over $\Omega^\circ$.) We can convert the resulting expression into the form \eqref{min_func2}, by sorting over $E$. Hence, we can express each $\omega^\new_\alpha$ in short form using work $ \mathcal{O} (N \log N)$, since $\displaystyle \sum_{\omega \in \Omega^\circ} \depth(\omega) \leq CN$.

Define the map $ \mathfrak{R} : \X(E) \rightarrow \cP$ by the formula $\mathfrak{R}(f)(x) = \sum_{\alpha \in \cM} \frac{1}{\alpha !} \omega^{\new}_\alpha(f) \cdot x^\alpha$. Hence, we have the key condition
\begin{equation} \label{mm1}
\sum_{\xi \in \Xi^\circ} \lvert \xi(f, \mathfrak{R}(f)) \rvert^p \leq C \cdot \sum_{\xi \in \Xi^\circ} \lvert \xi(f,R) \rvert^p \quad \mbox{for any} \; R \in \cP.
\end{equation}

\noindent\textbf{\underline{Picking a cutoff function:}}

We choose a function $\theta^\circ \in C^m(\R^n)$ such that
\begin{enumerate}
\item $\theta^\circ = 0$ on $\R^n \setminus Q^\circ$,
\item $\theta^\circ = 1$ on $E$,
\item $\lvert \partial^\alpha \theta^\circ(x) \rvert \leq C $ for $x \in \R^n$, $\lvert \alpha \rvert \leq m$,
\item Given $\underline{x} \in \R^n$, we can compute $J_{\underline{x}} \theta^\circ$ using work and storage at most $C$.
\end{enumerate}
We can arrange these conditions by taking $\theta^\circ$ to be a spline function that equals $1$ on $\frac{1}{32}Q^\circ$ and equals $0$ on $\R^n \setminus Q^\circ$; recall that $E \subset \frac{1}{32}Q^\circ$.

\noindent\textbf{\underline{Main definitions:}}
\begin{itemize}
\item Let the list $\Omega \subset (\X(E))^*$ consist of all the functionals $\omega$ in $\Omega^\circ$ and all the functionals of the form $f \mapsto \partial^\beta \left[ \mathfrak{R}(f) \right](0) = \omega^\new_\beta(f)$ for all $\beta \in \cM$.
\item Let the list $\Xi \subset (\X(E))^*$ consist of all the functionals $f \mapsto \xi^\circ(f,\mathfrak{R}(f))$ where $\xi^\circ \in \Xi^\circ$. Hence,
\[\sum_{\xi \in \Xi} \lvert \xi(f) \rvert^p = \sum_{\xi^\circ \in \Xi^\circ} \lvert \xi^\circ(f,\mathfrak{R}(f)) \rvert^p.\]
\item Let $T : \X(E) \rightarrow \X$ be defined by the formula
\[ Tf = \theta^\circ \cdot T^\circ(f, \mathfrak{R}(f)) + (1-\theta^\circ) \cdot \mathfrak{R}(f).\]
\end{itemize}
We note that the functionals $\xi \in \Xi$ and the map $T$ have $\Omega$-assisted bounded depth. 

We can list all the functionals in $\Xi$ and $\Omega$, with each functional expressed in short form, using work and storage at most $C N$. (We have already computed the functionals in the lists $\Omega^\circ$ and $\Xi^\circ$, and we have computed the map $f \mapsto \mathfrak{R}(f)$, all expressed in short form.)

We give a query algorithm for $T$. A query consists of a point $\underline{x} \in \R^n$. Then, using property (1) of $\theta^\circ$ we can write
\begin{equation*}
J_{\underline{x}}(Tf) =
\left\{
\begin{aligned}
&J_{\underline{x}} \theta^\circ \odot_{\underline{x}} J_{\underline{x}} T^\circ + J_{\underline{x}}(1 - \theta^\circ) \odot_{\underline{x}} \mathfrak{R}(f) &\mbox{if} \; \underline{x} \in Q^\circ \\
&  \mathfrak{R}(f) &\mbox{if} \; \underline{x} \notin Q^\circ.
\end{aligned}
\right.
\end{equation*}
We test whether $\underline{x} \in Q^\circ$ or $\underline{x} \in \R^n \setminus Q^\circ$. If $\underline{x} \in Q^\circ$, then we compute the map $f \mapsto J_{\underline{x}} (Tf)$ in short form in terms of the assists $\Omega$. This uses the query algorithm for $T^\circ$ and property (4) of $\theta^\circ$. Note that we can computate a short form representation of the $\odot_{\underline{x}}$-product or sum of polynomial-valued maps which are given in short form, using work at most $C$. If $\underline{x} \in \R^n \setminus Q^\circ$, then the map is given by $f \mapsto J_{\underline{x}}(Tf) = \mathfrak{R}(f)$, which is given in short form in terms of the assists $\Omega$. This completes the description of the query algorithm for $T$.  The query work is at most $C \log N$, as promised in Theorem \ref{main_thm_hom}.

\noindent\textbf{\underline{Main conditions:}}

According to the first bullet point in Proposition \ref{mprop}, we have 
\begin{equation}\label{star3}
\left\{
\begin{array}{ll}
& \#(\Xi) \leq C N , \;\; \mbox{and} \\
& \displaystyle \sum_{\omega \in \Omega} \depth(\omega) \leq CN.
\end{array}
\right.
\end{equation}

From property (2) of $\theta^\circ$, and since $T^\circ(f, \mathfrak{R}(f) ) = f$ on $E$, we see that
\begin{equation}
\label{m6}
Tf =f \; \mbox{on} \; E.
\end{equation}

We now estimate $\| Tf \|_{\X}$. A standard argument shows that
\[\| Tf \|_{\X} \leq C \cdot \left[ \| T^\circ(f,\mathfrak{R}(f)) \|_{\X(Q^\circ)} + \| T^\circ(f,\mathfrak{R}(f)) - \mathfrak{R}(f) \|_{L^p(Q^\circ)} \right].\]
(See the proof of Lemma \ref{patch_lem}.)
According to the third bullet point in Proposition \ref{mprop}, we therefore have
\begin{align} \label{m7}
\| Tf \|_{\X}^p & \leq C \cdot \sum_{\xi^\circ \in \Xi^\circ} \lvert \xi^\circ(f,\mathfrak{R}(f)) \rvert^p \\
& = C \cdot \sum_{\xi \in \Xi} \lvert \xi(f) \rvert^p. \notag{}
\end{align}

We now observe that
\begin{align}
\label{m7a}
\sum_{\xi \in \Xi} \lvert \xi(f) \rvert^p = \sum_{\xi^\circ \in \Xi^\circ} \lvert \xi^\circ(f,\mathfrak{R}(f)) \rvert^p & \leq C \inf_{R \in \cP}  \sum_{\xi^\circ \in \Xi^\circ} \lvert \xi^\circ(f,R) \rvert^p  \quad (\mbox{see \eqref{mm1}})\\
\notag{}
& \leq C \inf_{R \in \cP} \|(f,R) \|_{\frac{65}{64}Q^\circ}^p \quad (\mbox{see Proposition \ref{mprop}}).
\end{align}
Moreover, by definition of the trace seminorm,
\[
 \|(f,R)\|_{\frac{65}{64}Q^\circ}  = \inf_{F \in \X} \left\{ \| F \|_{\X(\frac{65}{64}Q^\circ)} + \| F - R \|_{L^p(\frac{65}{64}Q^\circ)} :  F = f \; \mbox{on} \; E \right\}.
\]
Note that $\| F - R \|_{L^p(\frac{65}{64}Q^\circ)} \leq C \| F \|_{\X(\frac{65}{64}Q^\circ)}$ if we choose $R = J_{x_{Q^\circ}} F$ (thanks to the Sobolev inequality). Therefore,
\begin{equation*}
\inf_{R \in \cP} \|(f,R)\|_{\frac{65}{64}Q^\circ}   \leq  C \cdot  \inf_{F \in \X} \bigl\{ \|F \|_{\X} :  F = f \; \mbox{on} \; E \bigr\} = C \cdot \|f\|_{\X(E)}.
\end{equation*}
The previous estimates imply that
\begin{equation}
\label{m8}
\sum_{\xi \in \Xi} \lvert \xi(f) \rvert^p  \leq C \cdot \|f\|_{\X(E)}^p.
\end{equation}

Finally, we have $ \| f\|_{\X(E)} = \inf_{F \in \X} \bigl\{ \|F\|_\X : F = f \; \mbox{on} \; E \bigr\}  \leq \| Tf\|_{\X} $, thanks to \eqref{m6}. This estimate and \eqref{m7} imply that
\begin{equation}
\label{m9}
c \cdot \| f\|_{\X(E)}^p \leq \sum_{\xi \in \Xi} \lvert \xi(f) \rvert^p .
\end{equation}
In view of \eqref{star3}-\eqref{m9}, we have proven Theorem \ref{main_thm_hom}.

\hfill $\qed$

\section{Extension in Inhomogeneous Sobolev Spaces}\label{inhom_sec}

Let $E \subset \R^n$ be finite, and let $N = \#(E)$.

The inhomogeneous Sobolev space $W^{m,p}(\R^n) \subset L^{m,p}(\R^n)$ consists of real-valued functions $F$ on $\R^n$ such that $\partial^\alpha F \in L^p(\R^n)$ for all $\lvert \alpha \rvert \leq m$. This space is equipped with the norm \[\| F \|_{W^{m,p}(\R^n)} = \left( \int_{\R^n}  \sum_{| \alpha| \leq m} \lvert \partial^\alpha F(x) \rvert^p dx \right)^{1/p}.\]

Let $W^{m,p}(E)$ denote the space of functions $f : E \rightarrow \R$, equipped with the trace norm
\[\| f\|_{W^{m,p}(E)} = \inf_{F \in W^{m,p}(\R^n)} \left\{ \| F \|_{W^{m,p}(\R^n)} :  F = f \; \mbox{on} \; E \right\}.\]

We use our extension results for the homogeneous Sobolev space $L^{m,p}(\R^n)$ to obtain analogous results for the inhomogeneous Sobolev space $W^{m,p}(\R^n)$. We will exhibit a query algorithm for a linear extension operator $T : W^{m,p}(E) \rightarrow W^{m,p}(\R^n)$ and we will compute a formula that approximates the $W^{m,p}(E)$ trace norm. We will do so using one-time work at most $C N \log N$ in space at most $C N$. Given $\underline{x} \in \R^n$ and $\lvert \alpha \rvert \leq m-1$, we will explain how to compute $\partial^\alpha Tf (\underline{x})$ using work at most $C \log N$.

\subsection{Case I}
We assume that $N = \#(E) \geq 2$ and that $E \subset \frac{1}{32} Q^\circ$, where $Q^\circ = [0,1)^n$.

\label{pp24}
We apply Proposition \ref{mprop} to define an extension operator $T^\circ : W^{m,p}(E) \oplus \cP \rightarrow W^{m,p}(Q^\circ)$ and lists $\Xi^\circ \subset (W^{m,p}(E) \oplus \cP)^*$ and $\Omega^\circ \subset (W^{m,p}(E))^*$. 

We define a cutoff function $\theta^\circ$ on $\R^n$. As in Section \ref{hom_sec}, we assume that the function $\theta^\circ \in C^m(\R^n)$ satisfies $\theta^\circ = 0$ on $\R^n \setminus Q^\circ$, $\theta^\circ = 1$ on $E$, and $\lvert \partial^\alpha \theta^\circ (x) \rvert \leq C$ for all $x \in \R^n$ and $|\alpha| \leq m$. Furthermore, we assume that we can compute $J_{\underline{x}} \theta^\circ$ for $\underline{x} \in \R^n$ using work at most $C$. We accomplish this by taking $\theta^\circ$ to be an appropriate spline function.

We define a linear map $T : W^{m,p}(E) \rightarrow W^{m,p}(\R^n)$ by
\[T f := \theta^\circ \cdot T^\circ(f,0) \qquad \mbox{for any} \; f \in W^{m,p}(E).\]
Proposition \ref{mprop} states that $T^\circ(f,0) = f$ on $E$. Thus, since $\theta^\circ = 1$ on $E$ we have
\begin{equation}
\label{eeop}
T f = f \text{ on } E.
\end{equation} 
We write $J_{\underline{x}} (T f) = J_{\underline{x}}  \theta^\circ  \odot_{\underline{x}} J_{\underline{x}}  T^\circ(f,0)$. Hence, we compute $J_{\underline{x}} (T f) = 0$ whenever $\underline{x} \in \R^n \setminus Q^\circ$ (since $\theta^\circ \equiv 0$ on $\R^n \setminus Q^\circ$). On the other hand, if $\underline{x} \in Q^\circ$ then we can compute the map $f \mapsto J_{\underline{x}}T^\circ(f,0)$ in short form in terms of the assists $\Omega^\circ$ (see Proposition \ref{mprop}), hence we can compute the map $f \mapsto J_{\underline{x}} (T f)$ in short form by basic algebra (multiplying polynomials). Thus we have given a query algorithm for $T$.

Proposition \ref{mprop} states that
\begin{equation}
\label{z2}
\| T^\circ( f, 0) \|_{L^{m,p}(Q^\circ)}^p  + \| T^\circ(f,0)  \|_{L^p(Q^\circ)}^p \leq C \sum_{\xi \in \Xi^\circ} \lvert \xi(f,0) \rvert^p
\end{equation}
and
\begin{align}
\label{z3}
\sum_{\xi \in \Xi^\circ} \lvert \xi(f,0) \rvert^p & \leq C \cdot \inf_{F} \left\{ \| F \|^p_{L^{m,p}(\frac{65}{64}Q^\circ)} +  \| F \|^p_{L^p(\frac{65}{64}Q^\circ)}: F = \; f \; \mbox{on} \; E \right\} \\
& \leq C \cdot \inf_{F} \left\{ \| F \|^p_{W^{m,p}(\frac{65}{64}Q^\circ)}: F = \; f \; \mbox{on} \; E \right\}.  \notag{}
\end{align}

We write $H = T f = \theta^\circ H^\circ$ with $H^\circ = T^\circ(f,0)$. Recall that $\theta^\circ$ is supported on $Q^\circ$ and that the derivatives of $\theta^\circ$ are bounded by a constant $C$. Hence, applying the Leibniz rule we see that
\begin{align*}
\| H \|_{W^{m,p}(\R^n)}^p &\leq C \cdot \sum_{|\alpha| + |\beta | \leq m} \int_{Q^\circ} \lvert \partial^\alpha H^\circ (x) \rvert^p \cdot \lvert \partial^\beta \theta^\circ(x) \rvert^p dx  \\
& \leq C \cdot \biggl[ \sum_{| \alpha | \leq m} \| \partial^\alpha H^\circ \|^p_{L^p(Q^\circ)} \biggr].
\end{align*}
Thus, Proposition \ref{int_ineq} implies that
\begin{equation*}
\| H \|_{W^{m,p}(\R^n)} \leq C \cdot \biggl[ \| H^\circ \|_{L^{m,p}(Q^\circ)} + \| H^\circ \|_{L^p(Q^\circ)} \biggr].
\end{equation*}
We finally apply \eqref{z2} and insert the definitions of $H$ and $H^\circ$ to see that
\begin{equation}\label{z5}
\| T f \|_{W^{m,p}(\R^n)}^p \leq C \cdot \sum_{\xi \in \Xi^\circ} \lvert \xi(f,0) \rvert^p.
\end{equation}

We set $\Omega = \Omega^\circ$. We define a list $\Xi \subset (\X(E))^*$ consisting of the functionals $f \mapsto \xi(f, 0 )$ for all $\xi \in \Xi^\circ$. The estimates \eqref{z3} and \eqref{z5} imply that 
\[ c \| T f \|_{W^{m,p}(\R^n)}^p \leq \sum_{\xi \in \Xi} \lvert \xi(f) \rvert^p \leq C \inf_{F} \bigl\{ \| F \|^p_{W^{m,p}(\frac{65}{64}Q^\circ)}: F = \; f \; \mbox{on} \; E \bigr\}.
\]
Moreover, note that $\| f \|_{W^{m,p}(E)} \leq \| T f \|_{W^{m,p}(\R^n)}$, since $Tf = f$ on $E$.

All the functionals $\xi \in \Xi$ and the map $T$ have $\Omega$-assisted bounded depth. We have given a query algorithm for $T$, and we have listed the functionals in $\Omega$. We can list the functionals in $\Xi$, expressed in short form in terms of the assists $\Omega$, using work and storage at most $C N$. To see this, note that there are at most $C N$ functionals in $\Xi^\circ$. We  determine  a short form representation of the  functional $f \mapsto \xi(f,0)$ using the short form representation of $(f,P) \mapsto \xi(f,P)$ by just deleting all the coefficients of the variables $(\partial^\alpha P(0))_{\alpha \in \cM}$. This requires work at most $C$ per functional $\xi$.

According to the first bullet point in Proposition \ref{mprop} we have 
\begin{equation*}
\left\{
\begin{array}{ll}
& \#(\Xi) \leq C N , \;\; \mbox{and}  \\
& \displaystyle \sum_{\omega \in \Omega} \depth(\omega) \leq CN.
\end{array}
\right.
\end{equation*}

The preceding argument establishes the case $N = \#(E) \geq 2$ in the result below.

\begin{prop} \label{inhom_prop1}
Assume that we are given a finite subset $E \subset \frac{1}{32} Q^\circ$, with $Q^\circ = [0,1)^n$. Let $N = \#(E)$. 
\begin{itemize}
\item We compute lists $\Omega$ and $\Xi$, consisting of functionals on $W^{m,p}(E) = \{ f : E \rightarrow \R\}$, with the following properties.
\begin{itemize}
\item The sum of $\depth(\omega)$ over all $\omega \in \Omega$ is bounded by $C N$. The number of functionals in $\Xi$ is at most $C N$.
\item Each functional $\xi$ in $\Xi$ has $\Omega$-assisted bounded depth. The functionals in $\Omega$ and $\Xi$ are represented in their short form.
\item For all $f \in W^{m,p}(E)$ we have
\[ c \cdot \| f\|_{W^{m,p}(E)} \leq \left[ \sum_{\xi \in \Xi} \lvert \xi(f) \rvert^p  \right]^{1/p}  \leq C  \cdot \inf \left\{ \| F \|_{W^{m,p}(\frac{65}{64}Q^\circ)}  : F = f \; \mbox{on} \; E \right\}.\]
\end{itemize}
\end{itemize}

Moreover, there exists a linear map $T : W^{m,p}(E) \rightarrow W^{m,p}(\R^n)$ with the following properties.
\begin{itemize}
\item $T$ has $\Omega$-assisted depth at most $C$.
\item $T f = f$ on $E$ and $\| Tf \|_{W^{m,p}(\R^n)}^p \leq C \cdot  \sum_{\xi \in \Xi} \lvert \xi(f) \rvert^p$ for all $f \in W^{m,p}(E)$.
\item We produce a query algorithm that operates as follows.

Given a query point $\underline{x} \in \R^n$, we compute a short form description of the $\Omega$-assisted bounded depth map $ f \mapsto J_{\underline{x}}\left( Tf  \right)$ using work and storage at most $C \log(2+ N)$.
\end{itemize}

The computations above require one-time work at most $C N \log (2+ N) + C$ in space at most $C N + C$.

\end{prop}
\begin{proof}
We have already established the proposition in case $N = \#(E) \geq 2$.

When $E$ is a singleton $\{ x^\circ\}$ (i.e., $N=1$), we define $\Xi = \{ \xi_0 \}$ and $\Omega = \emptyset$, where $\xi_0(f) := f(x^\circ)$. We define
\[ (Tf)(x) = \theta^\circ(x) \cdot f(x^\circ) \quad \mbox{for any} \; f : E \rightarrow \R,\]
where $\theta^\circ \in C^m(\R^n)$ is supported on $Q^\circ$ and $\theta^\circ(x^\circ) = 1$. Moreover, a simple computation shows that $\| Tf \|_{W^{m,p}(\R^n)} \leq C \cdot \lvert \xi_0(f) \rvert$. 

On the other hand,
\[ \lvert \xi_0(f) \rvert = \lvert f(x^\circ) \rvert  \leq C \cdot \| F \|_{W^{m,p}(\frac{65}{64}Q^\circ)}\] 
for any $F \in W^{m,p}(\R^n)$ such that $F(x^\circ) = f(x^\circ)$, thanks to the Sobolev inequality. The above computations require a constant amount of work. This completes the proof of the proposition in the case $N = \#(E) = 1$.

When $E = \emptyset$, we define $\Xi = \Omega = \emptyset$. We define $T : W^{m,p}(E) \rightarrow W^{m,p}(\R^n)$ to be the trivial (zero) map defined on a zero-dimensional space. These objects vacuously satisfy the conclusion of the proposition. This completes the proof of the result.
\end{proof}

\subsection{Case II}
\label{sec_caseii}

Here we strengthen Proposition \ref{inhom_prop1} by removing the hypothesis that $E$ is contained in a unit cube. We employ a standard partition of unity argument.

Assume that $E$ is a finite subset of $\R^n$. Let $N = \#(E)$.

We decompose $\R^n$ into a collection of cubes $\{ Q_{\mathbf{k}} : \mathbf{k} \in \Z^n\}$. For each $\mathbf{k} = (k_1,k_2,\cdots,k_n) \in \Z^n$ we define
\[Q_{\mathbf{k}} =  \left( k_1 \cdot 2^{-10} , k_1 \cdot 2^{-10} + 1 \right] \times \cdots \times \left(k_n \cdot 2^{-10}, k_n \cdot 2^{-10} + 1 \right].\]
Note that the union of all the $Q_{\mathbf{k}}$ is equal to $\R^n$. Moreover, 
\begin{equation}
\label{bdd_overlap}
\text{each point in} \; \R^n \; \mbox{is contained in at most} \; C \; \text{of the cubes} \; 200 Q_{\mathbf{k}}.
\end{equation}

We define $E_{\mathbf{k}} := E \cap Q_{\mathbf{k}}$ for each $\mathbf{k} \in \Z^n$.  According to the above,
\[ \bigcup_{\mathbf{k} \in \Z^n} E_{\mathbf{k}} = E \;\; \text{and} \;\; \sum_{\mathbf{k} \in \Z^n} \# (E_{\mathbf{k}} ) \leq C \cdot N.\]
We define
\[
\cI := \{ \mathbf{k} \in \Z^n : E_{\mathbf{k}} \neq \emptyset\}.
\]
For each $x \in E$ we can easily list all the indices $\mathbf{k} \in \Z^n$ such that $x \in Q_{\mathbf{k}}$. We concatenate these lists and remove duplicate indices by sorting. Thus, we can compute the collection $\cI$. We know that $\# \cI \leq C \cdot N$ by \eqref{bdd_overlap}. Moreover, the computation of $\cI$ requires work at most $C N \log (N+2)$ in space at most $C N$. For each $\mathbf{k}$ arising from some $x \in E$ as above, we include $x$ in a list associated to $\mathbf{k}$. In this way we construct the subsets $E_{\mathbf{k}}$ for each $\mathbf{k} \in \cI$. This computation requires work at most $C N \log(N+2)$.

For each $\mathbf{k} \in \cI$ we do the following. According to Proposition \ref{inhom_prop1}, we can compute lists $\Xi_{\mathbf{k}}$ and $\Omega_{\mathbf{k}}$ of linear functionals on $W^{m,p}(E_{\mathbf{k}})$. We also give a query algorithm for a linear extension operator $T_{\mathbf{k}} : W^{m,p}(E_{\mathbf{k}}) \rightarrow W^{m,p}(\R^n)$ (see below). The following properties hold.

\begin{description}
\item[(a)] Each $\xi \in \Xi_{\mathbf{k}}$ has $\Omega_{\mathbf{k}}$-assisted bounded depth.
\item[(b)] The sum of $\depth(\omega)$ over all $\omega \in \Omega_{\mathbf{k}}$ is bounded by $C \#(E_{\mathbf{k}})$. The number of functionals in $\Xi_{\mathbf{k}}$ is at most $C \cdot \#(E_{\mathbf{k}})$.
\item[(c)] $T_{\mathbf{k}}$ has $\Omega_{\mathbf{k}}$-assisted bounded depth.
\item[(d)] $(T_{\mathbf{k}} f)(x) = f(x)$ for all $x \in E_{\mathbf{k}}$.
\item[(e)] For any $f \in W^{m,p}(E_{\mathbf{k}})$ we have
\begin{equation}\label{zz1} 
\sum_{\xi \in \Xi_{\mathbf{k}}} \lvert \xi(f) \rvert^p \leq C \inf \left\{ \| F \|^p_{W^{m,p}(200Q_{\mathbf{k}} )} : F = f \; \mbox{on} \; E_{\mathbf{k}} \right\}
\end{equation}
and
\begin{equation}
\label{zz2}
\| T_{\mathbf{k}}f \|^p_{W^{m,p}(\R^n)} \leq C \cdot \sum_{\xi \in \Xi_{\mathbf{k}}} \lvert \xi(f) \rvert^p.
\end{equation}
\item[(f)] We can query the extension operator. A query consists of a point $\underline{x} \in \R^n$. We respond to the query $\underline{x}$ with a short form description of the $\Omega_{\mathbf{k}}$-assisted bounded depth map $f \mapsto J_{\underline{x}} (T_{\mathbf{k}} f)$. The query work is at most $C \log (2+ \#(E_{\mathbf{k}}))$.

\end{description}
(Here, we use the fact that $E_{\mathbf{k}} \subset \frac{1}{32} ( 200 Q_{\mathbf{k}} )$. To achieve the preceding results, we apply Proposition \ref{inhom_prop1} to a rescaled and translated copy of $E_{\mathbf{k}}$. We leave details to the reader.)

The above computations require one-time work at most  $C  \#(E_{\mathbf{k}}) \log(\#(E_{\mathbf{k}}) + 1)$ and storage at most $C \# (E_{\mathbf{k}})$ for each $\mathbf{k} \in \cI$. Thus, the total work and space required are at most $C N \log (N+1) + C$ and $C N + C$, respectively.

We will define an extension operator $T :W^{m,p}(E) \rightarrow W^{m,p}(\R^n)$ and lists $\Xi$ and $\Omega$ consisting of linear functionals  on $W^{m,p}(E)$.
\begin{itemize}
\item Let $\Omega \subset (W^{m,p}(E))^*$ be the union of the lists $\Omega_{\mathbf{k}}$ for all $\mathbf{k} \in \cI$. \\
(If $E = \emptyset$ then we define $\Omega = \emptyset$.)
\item Let $\Xi \subset (W^{m,p}(E))^*$ be the union of the lists $\Xi_{\mathbf{k}}$ for all $\mathbf{k} \in \cI$. Hence,
\begin{equation}\label{zz3}
\sum_{\xi \in \Xi} \lvert \xi(f) \rvert^p = \sum_{\mathbf{k} \in \cI} \sum_{\xi \in \Xi_{\mathbf{k}} } \lvert \xi(f) \rvert^p.
\end{equation}
(If $E = \emptyset$ then we define $\Xi = \emptyset$.)
\end{itemize}

\begin{remk}\label{rem_zz}
The sum of $\depth(\omega)$ over $\omega \in \Omega$ is bounded by
\[C \sum_{\mathbf{k} \in \cI} \#(E_{\mathbf{k}}) \leq C N.\]
Also,
\[ \#(\Xi) \leq \sum_{\mathbf{k} \in \cI} \#(\Xi_{\mathbf{k}}) \leq \sum_{\mathbf{k} \in \cI} C \cdot \#(E_{\mathbf{k}}) \leq C N. \]
\end{remk}

We choose a partition of unity $\{ \theta_{\mathbf{k}} \}_{\mathbf{k} \in\Z^n}$ with the following properties.
\begin{itemize}
\item $\theta_{\mathbf{k}} \in C^m(\R^n)$, and $\theta_{\mathbf{k}}$ is supported on $Q_{\mathbf{k}}$.
\item $\lvert \partial^\alpha \theta_{\mathbf{k}}(x) \rvert \leq C$ for all $|\alpha| \leq m$ and $x \in \R^n$.
\item $\sum_{\mathbf{k} \in \Z^n} \theta_{\mathbf{k}} = 1$ on $\R^n$.
\item Given $\underline{x} \in \R^n$ and $\mathbf{k} \in \Z^n$, we can compute $J_{\underline{x}} \theta_{\mathbf{k}}$ using work at most $C \log (2+ N)$.
\end{itemize}
These conditions are easy to arrange. For instance, see the construction of $\{\theta_Q\}_{Q \in \CZ}$ given in Section \ref{hom_sec}. We leave details to the reader.

We define $T : W^{m,p}(E) \rightarrow W^{m,p}(\R^n)$ by the formula
\[(Tf)(x) = \sum_{\mathbf{k} \in \cI} (T_{\mathbf{k}} f)(x)  \cdot \theta_{\mathbf{k}}(x) \qquad (x \in \R^n).\]

Assume that $\underline{x} \in \R^n$ is given. Note that $J_{\underline{x}} \theta_{\mathbf{k}}$ is nonzero only when $Q_{\mathbf{k}}$ contains $\underline{x}$ (since $\theta_{\mathbf{k}}$ is supported on $Q_{\mathbf{k}}$).  We compute a list of all the indices $\mathbf{k} \in \cI$ such that $\underline{x} \in Q_{\mathbf{k}}$ using a binary search; this requires work at most $C \log (2+ N)$. For each such $\mathbf{k}$, we compute a short form description of the linear map $f \mapsto J_{\underline{x}} (T_{\mathbf{k}} f)$. That requires work at most $C \log (2 + N)$. Hence, to compute $J_{\underline{x}} (Tf)$ we can sum the linear maps $J_{\underline{x}} \theta_{\mathbf{k}} \odot_{\underline{x}} J_{\underline{x}} (T_{\mathbf{k}}f)$ over all the relevant indices $\mathbf{k}$. Thus we can compute a short form description of the linear map $f \mapsto J_{\underline{x}} (Tf)$ using work and storage at most $C \log ( 2 + N) $.

Let $x \in E$. 

Let $\mathbf{k} \in \cI$. Recall that $\theta_{\mathbf{k}}(x) = 0$ if $x \notin Q_{\mathbf{k}}$. Also, $T_{\mathbf{k}}f(x) = f(x)$ if $x \in  Q_{\mathbf{k}}$. Thus, we have $\theta_{\mathbf{k}}(x) T_{\mathbf{k}}f(x) = \theta_{\mathbf{k}}(x) f(x) $ unconditionally.

Let $\mathbf{k} \in \Z^n \setminus \cI$. By definition of $\cI$ we know that $x \notin Q_{\mathbf{k}}$, hence $\theta_{\mathbf{k}}(x) = 0$.

Hence,
\begin{align*}
Tf(x) &= \sum_{\mathbf{k} \in \cI} \theta_{\mathbf{k}}(x) T_{\mathbf{k}}f(x) \\
&= \sum_ {\mathbf{k} \in \cI} \theta_{\mathbf{k}}(x) T_{\mathbf{k}}f(x) +  \sum_ {\mathbf{k} \in \Z^n \setminus \cI} \theta_{\mathbf{k}}(x) f(x)  \\
& =  \sum_{\mathbf{k} \in \Z^n} \theta_{\mathbf{k}}(x) f(x) = f(x) \qquad \mbox{for any} \; x \in E.
\end{align*}

\begin{prop}\label{inhom_prop2}
For each $f \in W^{m,p}(E)$ we have
\[\| Tf \|^p_{W^{m,p}(\R^n)} \leq C \cdot \sum_{\xi \in \Xi} \lvert \xi(f) \rvert^p.\]
Furthermore,
\[c \cdot  \| f \|_{W^{m,p}(E)}^p \leq \sum_{\xi \in \Xi} \lvert \xi(f) \rvert^p \leq C \cdot \| f \|_{W^{m,p}(E)}^p.\]
\end{prop}
\begin{proof}

To prove the first estimate, we recall that $Tf = \sum_{\mathbf{k} \in \cI} (T_{\mathbf{k}} f) \cdot \theta_{\mathbf{k}}$. Recall that $\theta_{\mathbf{k}}$ is supported on $Q_{\mathbf{k}}$ and that the derivatives of $\theta_{\mathbf{k}}$ are uniformly bounded. Also, note that each point in $\R^n$ is contained in at most $C$ of the cubes $Q_{\mathbf{k}}$ (see \eqref{bdd_overlap}). Hence, by the Leibniz rule we have
\begin{align*}
\| Tf \|_{W^{m,p}(\R^n)}^p &\leq C \sum_{\mathbf{k} \in \cI} \| T_{\mathbf{k}} f \|^p_{W^{m,p}(Q_{\mathbf{k}})} \\
& \leq C \sum_{\mathbf{k} \in \cI} \sum_{\xi \in \Xi_{\mathbf{k}}} \lvert \xi(f) \rvert^p \qquad\; (\mbox{see \eqref{zz2}}) \\
& = C \sum_{\xi \in \Xi} \lvert \xi(f) \rvert^p \qquad \qquad (\mbox{see \eqref{zz3}}).
\end{align*}
Hence,
\[\sum_{\xi \in \Xi} \lvert \xi(f) \rvert^p \geq c \| Tf\|_{W^{m,p}(\R^n)}^p \geq c\| f \|_{W^{m,p}(E)}^p.\]
In the last inequality above, we use the definition of the seminorm $\| \cdot \|_{W^{m,p}(E)}$ and the fact that $Tf =f $ on $E$.

For the reverse inequality, we use \eqref{zz1} and \eqref{zz3} and deduce that
\begin{align*}
\sum_{\xi \in \Xi} \lvert \xi(f) \rvert^p & = \sum_{\mathbf{k} \in \cI} \sum_{\xi \in \Xi_{\mathbf{k}}} \lvert \xi(f) \rvert^p \\
& \leq C\cdot  \sum_{\mathbf{k} \in \cI} \inf \left\{ \| F \|_{W^{m,p}(200Q_{\mathbf{k}})}^p : F = f \; \mbox{on} \; E \cap Q_{\mathbf{k}} \right\} \\
& \leq C \cdot \inf \left\{ \sum_{\mathbf{k} \in \cI} \| F \|_{W^{m,p}(200Q_{\mathbf{k}})}^p : F = f \; \mbox{on} \; E\right\} \\
& \leq C \cdot \inf \left\{ \| F \|_{W^{m,p}(\R^n)}^p : F = f \; \mbox{on} \; E\right\} \\
& = C \cdot \| f \|_{W^{m,p}(E)}^p.
\end{align*}
Here, we use \eqref{bdd_overlap} to prove the last inequality.

This completes the proof of Proposition \ref{inhom_prop2}.

\end{proof}

The above construction implies our main result for the inhomogeneous Sobolev space. See Proposition \ref{inhom_prop2} and Remark \ref{rem_zz}.
\begin{thm} \label{main_thm_inhom}
Given a finite subset $E \subset \R^n$ with $N = \#(E)$, we perform one-time work at most $C N \log( 2 + N) + C$ in space at most $C N + C$, after which we have achieved the following.
\begin{itemize}
\item We compute lists $\Omega$ and $\Xi$, consisting of functionals on $W^{m,p}(E) = \{ f : E \rightarrow \R\}$, with the following properties.
\begin{itemize}
\item The sum of $\depth(\omega)$ over all $\omega \in \Omega$ is bounded by $C N$. The number of functionals in $\Xi$ is at most $C N$.
\item Each functional $\xi$ in $\Xi$ has $\Omega$-assisted bounded depth. The functionals in $\Omega$ and $\Xi$ are represented in their short form.
\item For all $f \in W^{m,p}(E)$ we have
\[ c \| f\|_{W^{m,p}(E)} \leq \left[ \sum_{\xi \in \Xi} \lvert \xi(f) \rvert^p  \right]^{1/p}  \leq C \| f \|_{W^{m,p}(E)}.\]
\end{itemize}
\end{itemize}

Moreover, there exists a linear map $T : W^{m,p}(E) \rightarrow W^{m,p}(\R^n)$ with the following properties.
\begin{itemize}
\item $T$ has $\Omega$-assisted depth at most $C$.
\item $T f = f$ on $E$ and $\| Tf \|_{W^{m,p}(\R^n)} \leq C \cdot  \| f \|_{W^{m,p}(E)}$ for all $f \in W^{m,p}(E)$.
\item We produce a query algorithm that operates as follows.

Given a query point $\underline{x} \in \R^n$, we respond with a short form description of the $\Omega$-assisted bounded depth map $ f \mapsto J_{\underline{x}}\left( Tf  \right)$ using work and storage at most $C \log(2+ N)$.
\end{itemize}

\end{thm}

At last, note that Theorem \ref{main_thm_hom} and Theorem \ref{main_thm_inhom} imply the main theorem from the introduction (Theorem \ref{introMainTheorem}). This completes the analysis of our algorithms in the infinite-precision model of computation. In the Appendix we present an analogue of our algorithm in a finite-precision model of computation.

\chapter{Appendix: Modifications for Finite-Precision}

\section{The Finite-Precision Model of Computation}
\label{sec_moc2}

Our model of computation in finite-precision is a slight variant of that described in Section 38 of \cite{FK2}. We spell out the details.

For an integer $S \geq 1$, we work with ``machine numbers'' of the form $k \cdot 2^{-S}$, with $k$ an integer and $\lvert k \rvert \leq 2^{+ 2 S}$. Our model of computation consists of an idealized von Neumann computer \cite{V}, able to handle machine numbers. We make the following assumptions:
\begin{itemize}
\item Given two distinct non-negative machine numbers $x$ and $y$, we can compute the most significant digit in which the binary expansions of $x$ and $y$ differ. (That is, for $\displaystyle x = \sum_{j \geq -S} x_j 2^j$ and $\displaystyle y = \sum_{j \geq -S} y_j 2^j$ with each $x_j$ and $y_j$ equal to $0$ or $1$, we compute the largest $j$ for which $x_j$ does not equal $y_j$.) We assume this takes one unit of ``work''. See also the paragraph following \eqref{task}.

This assumption is reasonable, since presumably a machine number is encoded in the computer as the sequence of its binary digits.

\item Two machine numbers $x$ and $y$ satisfying $\lvert x \rvert \leq 2^\ell$ and $\lvert y \rvert \leq 2^{\ell'}$ with $\ell,\ell' \geq 0$ and $\ell + \ell' \leq S$ can be ``multiplied'' to produce a machine number $x \otimes y$ satisfying $\lvert x \otimes y - xy \rvert \leq 2^{- S}$.

We suppose it takes one unit of ``work'' to compute $x \otimes y$.

We assume that $0 \otimes x = x \otimes 0 = 0$ and that $x \otimes 1 = 1 \otimes x = x$.

We assume that if $\lvert x \rvert \leq 2^\ell$ and  $\lvert y \rvert \leq 2^{\ell'}$, for $\ell,\ell'$ integers, then $\lvert x \otimes y \rvert \leq 2^{\ell + \ell'}$.

\item If $x$ is any machine number other than zero, then we suppose we can produce a machine number ``$1/x$'' in one unit of ``work'', such that $\lvert \text{``}1/x\text{''} - 1/x \rvert \leq 2^{- S}$.

We assume that $\text{``}1/x\text{''} = 1$ when $x=1$.

We assume that if $\lvert x \rvert \geq 2^{\ell}$, for an integer $\ell$, then $\lvert \text{``}1/x\text{''} \rvert \leq 2^{-\ell}$.
\item Two machine numbers $x$ and $y$ satisfying $\lvert x \rvert \leq \ell$ and $\lvert y \rvert \leq \ell'$ for integers $\ell$ and $\ell'$ such that $\ell + \ell' \leq 2^{S}$ may be added to produce their exact sum $x + y$, which is again a machine number.

We assume it takes one unit of ``work'' to compute $x + y$.

\item If $x$ is any machine number, then $-x$ is again a machine number.

We assume it takes one unit of ``work'' to compute $-x$.

\item If $x$ and $y$ are machine numbers, then we can decide whether $x < y$, $y < x$, or $x=y$.

We assume this takes one unit of ``work''.
\item If $x$ is a machine number other than zero, then we can compute the greatest integer $\ell$ such that $2^\ell \leq \lvert x \rvert$.

We assume this takes one unit of ``work''

\item If $x$ is a machine number and $\ell$ is an integer with $\lvert \ell \rvert \leq S$, then we can compute the greatest integer $\leq 2^\ell x$. (If this integer lies outside $[-2^{S},+2^{S}]$, then we produce an error message, and abort our computation.)

We assume this takes one unit of ``work''.

\item We assume we can add, subtract, multiply and divide integers of absolute value $\leq 2^{S}$, in one unit of ``work''.

If we compute $x/y$ in integer arithmetic, for integers $x,y$ ($y \neq 0$) of absolute value at most $2^{S}$, then we obtain the greatest integer $\leq$ the real number $x/y$. If our desired answer lies outside $[- 2^{S}, + 2^{S}]$, then we produe an error message and abort our computation.

\item Given integers $x$, $y$ of absolute value $\leq 2^{S}$, we can decide whether $x < y$, $y < x$, or $x = y$.

We assume this takes one unit of ``work''.

\item If $\ell $ is an integer, with $\lvert \ell \rvert \leq S$, then we can compute exactly the machine number $2^\ell$.

We assume this takes one unit of ``work''.

\item If $x$ and $y$ are machine numbers satisfying $2^{-\ell} \leq x \leq 2^{\ell}$ and $\lvert y \rvert \leq \ell'$ for integers $\ell$ and $\ell'$ such that  $\ell \cdot \ell' \leq S$, then we can compute a machine number ``$x^y$'' in one unit of ``work'', such that $\lvert \text{``}x^y\text{''} - x^y \rvert \leq 2^{-S}$.

\item If $x$ is any positive machine number, then we can compute a machine number ``$\log x$'' in one unit of ``work'', such that $\lvert \text{``}\log x\text{''} - \log x \rvert \leq 2^{-S}$. Here, $\log x$ is the base two logarithm.

\item We assume we can read or write a machine number from (to) the RAM with one unit of ``work''. 

\item We assume we can read a machine number from input or write a machine number to output in one unit of ``work''.

\item We assume we can store a single $S$-bit word in memory using one unit of ``storage''.
\item We assume we can store the address of any memory cell in a single $S$-bit word.
\end{itemize}

Under these assumptions, we say that our computer can process ``$S$-bit machine numbers'' (though the actual implementation of those machine numbers seems to require at least $2S + 2$ bits.) We call $\Delta_{\min} = 2^{-S}$ the ``machine precision'' of our computer.

We fix an integer $\overline{S} \geq 1$. We assume that our computer can process $S$-bit machine numbers for $S = K_{\max} \cdot \overline{S}$, where $K_{\max} \in \N$ satisfies
\begin{equation}
\label{largeK}
K_{\max} \geq C, \; \mbox{for a large enough universal constant} \; C.
\end{equation}

We will show that when our algorithms receive their input as $\overline{S}$-bit machine numbers, then the output produced by our algorithm is accurate to at least $\overline{S}$ bits. We will verify that the work and storage required are as promised: at most $CN \log N$ operations for the one-time work,  and at most $C \log N$ operations for the query work, with $CN$ storage, where the constant $C$ depends only on $m$, $n$, and $p$. 

Throughout the remaining sections, $\Delta_{\min} = 2^{-S}$ will denote the precision of our computer, as just described.

\section{Algorithms in Finite-Precision}
\label{fp_not}

We recall that a universal constant is one that depends only on the parameters $m$, $n$, and $p$. We impose the following assumptions in this section.

\noindent\textbf{Main Assumptions}:
\begin{itemize}
\item We set $\Delta_0 := 2^{- \overline{S}}$ for an integer $\overline{S} \geq 1$.
\item We assume our computer can process $S$-bit machine numbers, with $S := K_{\max} \cdot \overline{S}$, where $K_{\max}$ satisfies \eqref{largeK}. Then $\Delta_{\min} = 2^{-S}$ represents the ``machine precision'' of our computer. A ``machine number'' will always denote an $S$-bit machine number.
\item We set $\Delta_g = 2^{-K_1 \overline{S}}$ and $\Delta_\epsilon = 2^{- K_2 \overline{S}}$ for integers $K_1,K_2 \geq 1$.
\item We assume that $\Delta_{\min} \leq \Delta_\epsilon^C$ and $\Delta_\epsilon \leq \Delta_g^C$ for a large enough universal constant $C$.
\end{itemize}

Assume that $w \in \R$ satisfies $ \lv w \rv \leq 2^{S}$. We may not be able to represent $w$ perfectly on a computer, but we can always store an approximation to $w$. We introduce the relevant notation below.
\begin{itemize}
\item We say that $w$ is \underline{specified to precision $\Delta_\epsilon$} if a machine number $w_\fin$ is given with $\lvert w - w_{\fin} \rvert \leq \Delta_\epsilon$. 
\item We say that $w$ is \underline{computed to precision $\Delta_\epsilon$} if there is a finite-precision algorithm that computes a machine number $w_{\fin}$ with $\lvert w - w_{\fin} \rvert \leq \Delta_\epsilon$.
\item We say that $w$ is \underline{specified (computed) with parameters $(\Delta_g,\Delta_\epsilon)$} if $\lvert w \rvert \leq \Delta_g^{-1}$, and if $w$ is specified (computed) to precision $\Delta_\epsilon$.
\end{itemize}

We illustrate this terminology in the next result, which establishes the numerical stability of arithmetic operations.

\begin{lem} 
\label{num_stab}
Suppose that $\Delta_0$, $\Delta_{\min}$, $\Delta_g$, and $\Delta_\epsilon$ are as in the \textbf{Main Assumptions}.

Let  $x,y \in \R$ be specified with parameters $(\Delta_g,\Delta_\epsilon)$. Then the following hold.
\begin{itemize}
\item We can compute $x + y$ with parameters $(c \Delta_g, C\Delta_\epsilon)$.
\item We can compute $x \cdot y$ with parameters $(\Delta_g^2, C \Delta_\epsilon \Delta_g^{-1})$.
\item If $\lv y \rv \geq \Delta_g$, we can compute $x/y$ with parameters $(\Delta_g^2, C \Delta_\epsilon \Delta_g^{-3})$.
\item If $x \in [ \Delta_g , \Delta_g^{-1} ]$, we can compute $\log x$ with parameters $(c \Delta_g,C \Delta_\epsilon \Delta_g^{-1})$.
\item Suppose that $x \in [ \Delta_g , \Delta_g^{-1} ]$ and $ \lv y \rv \leq A$ with $A \geq 1$, and suppose that $K_{\max} \geq 5A \cdot \max \{ K_1, K_2 \}$. Then we can compute $x^y$ with parameters $(\Delta_g^{A}, \Delta_\epsilon \Delta_g^{- C\cdot A})$. '

\end{itemize}
The above computations require work at most $C$. 

Here, the constants $c$ and $C$ are independent of all the parameters.

\end{lem}

\begin{proof}
By hypothesis, we suppose we are given machine numbers $\overline{x}, \overline{y}$ with $\lv x - \overline{x} \rv \leq \Delta_\epsilon$ and $\lv y - \overline{y} \rv \leq \Delta_\epsilon$. Moreover, we have $\lv x \rv \leq \Delta_g^{-1}$ and  $\lv y \rv  \leq  \Delta_g^{-1} $.

Since $\Delta_\epsilon \leq \Delta_g^{-1}$, we learn that $\lv \overline{x} \rv \leq 2 \Delta_g^{-1}$ and $\lv \overline{y} \rv \leq 2 \Delta_g^{-1}$.


\begin{enumerate}
\item Since $\lv \overline{x} + \overline{y} \rv \leq 4 \Delta_g^{-1} \leq \Delta_{\min}^{-1}$ and since $ \Delta_{\min}^{-1} = 2^S$, we can compute the ($S$-bit) machine number $A = \overline{x} + \overline{y}$. This computation requires one unit of work, by assumption on our model of computation. Then
\[
\lv A - (x+y)  \rv \leq \lv   \overline{x} - x \rv + \lv  \overline{y} - y \rv \leq 2 \Delta_\epsilon.
\]
Moreover, $\lv x + y \rv \leq \lv x \rv + \lv y \rv \leq 2 \Delta_g^{-1}$. 

Thus, we can compute the sum $x + y$ with parameters $(\frac{1}{2} \Delta_g, 2 \Delta_\epsilon)$.

\item Since $\lv \overline{x} \cdot \overline{y} \rv \leq 4 \Delta_g^{-2} \leq \Delta_{\min}^{-1}$, we can compute a machine number $P$ with $\lv P - \overline{x} \cdot \overline{y} \rv \leq \Delta_{\min} \leq \Delta_\epsilon$. (Recall that $\Delta_{\min}$ is the ``machine precision''.) We have
\begin{align*}
\lv x \cdot y - \overline{x} \cdot \overline{y} \rv   &\leq    \lv x -  \overline{x} \rv \cdot \lv y \rv + \lv y - \overline{y} \rv \cdot \lv \overline{x} \rv \\
&\leq  \Delta_\epsilon \cdot \Delta_g^{-1} + \Delta_\epsilon \cdot (2 \Delta_g^{-1}) = 3 \Delta_\epsilon \cdot \Delta_g^{-1}.
\end{align*}
Hence, $\lv P - x \cdot y \rv \leq \Delta_\epsilon + 3 \Delta_\epsilon  \Delta_g^{-1} \leq 4 \Delta_\epsilon \Delta_g^{-1}$. Moreover, $\lv x \cdot y \rv \leq \Delta_g^{-1} \cdot \Delta_g^{-1} = \Delta_g^{-2}$. 

Therefore, we can compute the product $x \cdot y$ with parameters $(\Delta_g^2, 4 \Delta_\epsilon \Delta_g^{-1})$. 

Here, we have only used the assumptions $\Delta_{\min} \leq \frac{1}{4} \Delta_g^{2}$ and $\Delta_{\min} \leq \Delta_\epsilon$

\item Suppose that $\lv y \rv \geq \Delta_g$. Since we may assume $\Delta_\epsilon \leq \Delta_g^{10}$, we have
\[
\lv \overline{y} \rv \geq \lv y \rv - \lv y - \overline{y} \rv \geq \Delta_g - \Delta_\epsilon \geq \Delta_g - \Delta_g^{10}.
\] 
Since $\Delta_g \leq 1/2$, we conclude that $\lv \overline{y} \rv \geq \frac{1}{2} \Delta_g$.

Thus, we can compute a machine number $A$ with $\lv A -  (\overline{y})^{-1} \rv \leq \Delta_{\min} \leq \Delta_\epsilon$.

We have
\[
\lv y^{-1} - (\overline{y})^{-1} \rv = \frac{ \lv y - \overline{y} \rv}{ \lv y \rv \cdot \lv \overline{y} \rv } \leq \frac{\Delta_\epsilon}{\Delta_g \cdot \frac{1}{2} \Delta_g} =  2 \Delta_\epsilon \Delta_g^{-2}.
\]
Hence, $\lv A - y^{-1} \rv \leq  \Delta_\epsilon + 2 \Delta_\epsilon \Delta_g^{-2} \leq 3 \Delta_\epsilon \Delta_g^{-2}$. Moreover, $\lv y^{-1} \rv \leq \Delta_g^{-1}$. 

Therefore, we can compute $y^{-1}$ with parameters $(\Delta_g, 4 \Delta_\epsilon \Delta_g^{-2})$.

\item Suppose that $\lv y \rv \geq \Delta_g$.  According to (3), we can compute $y^{-1}$ with parameters $(\Delta_g,\Delta_\epsilon^{\new})$, where $\Delta_\epsilon^{\new} = 4 \Delta_\epsilon \Delta_g^{-2}$. We have $\Delta_\epsilon^{\new} \leq 4 \Delta_g^{8} \leq 1$ (since $\Delta_\epsilon \leq \Delta_g^{10}$) and $\Delta_\epsilon^{\new} \geq \Delta_\epsilon \geq \Delta_{\min}$. Hence, applying (2), we can compute $x \cdot y^{-1}$ with parameters $(\Delta_g^2, 4 \Delta_\epsilon^\new \Delta_g^{-1}) = (\Delta_g^2, 16 \Delta_\epsilon\Delta_g^{-3})$.

\item Suppose that $\Delta_g \leq x \leq \Delta_g^{-1}$. Since $\lv x - \overline{x} \rv \leq \Delta_\epsilon \leq \Delta_g^{10}$, we have $\frac{1}{2} \Delta_g \leq \overline{x} \leq 2 \Delta_g^{-1}$. 

We can compute a machine number $L$ satisfying $\lv L - \log \overline{x} \rv \leq \Delta_{\min} \leq \Delta_\epsilon$.

\[
\lv \log x - \log \overline{x} \rv \leq \frac{1}{\ln 2} \cdot \lv x - \overline{x} \rv \cdot \max \left\{ x^{-1} , (\overline{x})^{-1} \right\} \leq C \Delta_\epsilon \Delta_g^{-1}.
\]
(Recall that $\log x$ denotes the base two logarithm.) Hence, $\lv L - \log x \rv \leq \Delta_\epsilon  + C \Delta_\epsilon \Delta_g^{-1} \leq C' \Delta_\epsilon \Delta_g^{-1}$. Moreover, $\lv \log x \rv \leq \log \Delta_g^{-1} \leq C \Delta_g^{-1}$. 

Therefore, we can compute the number $\log x$ with parameters $(c \Delta_g,C' \Delta_\epsilon \Delta_g^{-1})$.

\item Suppose that $ \Delta_g \leq x \leq \Delta_g^{-1}$ and $ \lv y \rv \leq A$ for some $A \geq 1$.

Since $\lv x - \overline{x} \rv \leq \Delta_\epsilon \leq \Delta_g^{10}$ and $\lv y - \overline{y} \rv \leq \Delta_\epsilon \leq 1$, we conclude that $\frac{1}{2} \Delta_g \leq \overline{x}  \leq 2 \Delta_g^{-1}$ and $\lv \overline{y} \rv \leq \lv y \rv + \lv y - \overline{y} \rv  \leq 2 A$. 

We have $\lv \overline{x}^{\overline{y}} \rv \leq  (2 \Delta_g^{-1})^{2 A} = 2^{2A \cdot (K_2 \overline{S} +1)} \leq \Delta_{\min}^{-1}$, due to the assumption that $\Delta_{\min} = 2^{ - K_{\max} \overline{S}}$ and $K_{\max} \geq 4A K_2$. Thus, we can compute a machine number $B$ with 
\[
\lv B - \overline{x}^{\overline{y}} \rv \leq \Delta_{\min} \leq \Delta_\epsilon.
\]
This requires work at most $C$.

We have
\begin{align*}
\lv x^y - \overline{x}^{\overline{y}} \rv &\leq \lv x^y - \overline{x}^{y} \rv + \lv \overline{x}^y - \overline{x}^{\overline{y}} \rv \\
& = \lv e^{y \ln x} - e^{y \ln \overline{x}} \rv + \lv e^{y \ln \overline{x}} - e^{\overline{y} \ln \overline{x}} \rv \\
&\leq   \lv y \rv \cdot \lv x - \overline{x} \rv   \cdot    \max \left\{ x^{-1}, ( \overline{x})^{-1} \right\} \cdot     \max \left\{  e^{y \ln x} , e^{y \ln \overline{x}} \right\}   \\
& \hspace{2cm} +  \lv y - \overline{y} \rv \cdot \lv \ln \overline{x} \rv \cdot \max \left\{  e^{y \ln \overline{x}}  , e^{\overline{y} \ln \overline{x}} \right\} \\
& \leq  A \cdot \Delta_\epsilon  \cdot 2 \Delta_g^{-1} \cdot \Delta_g^{-C A} \\
& \hspace{2cm} + \Delta_\epsilon \cdot \ln \left(  2 \Delta_g^{-1} \right) \cdot \Delta_g^{-CA} \\
& \leq \Delta_\epsilon \Delta_g^{-C' A}.
\end{align*}
In the above, we use the estimates $ \lv e^w - e^z \rv \leq \lv w - z \rv \cdot \max \{ e^w, e^z \} $ and $\lv \ln x - \ln \overline{x} \rv \leq  \lv x - \overline{x} \rv \cdot \max \left\{ x^{-1}, (\overline{x})^{-1} \right\}$; both $C$ and $C'$ are numerical constants. 

Hence, $\lv B - x^y \rv \leq \Delta_\epsilon + \Delta_\epsilon \Delta_g^{-C' A} \leq \Delta_\epsilon \Delta_g^{- C'' A}$. Moreover, $\lv x^y \rv \leq  \Delta_g^{- A}$.  

Therefore, we can compute $x^y$ with parameters $(\Delta_g^{A},\Delta_\epsilon \Delta_g^{-C'' A})$ for a numerical constant $C''$.

\end{enumerate}

Thanks to (1)-(6), the conclusions of the lemma are verified. This completes the proof.

\end{proof}

We finish the section with a technical lemma concerning the evaluation of the $\ell^p$ norm by a finite-precision algorithm.\begin{lem}
\label{stablelp}
Let $\Delta \in \left[ \Delta_g, 1 \right]$ be a given machine number. Given real numbers $x_j$ ($1 \leq j \leq J$) with parameters $(\Delta_g,\Delta_\epsilon)$, where $J \leq \Delta_g^{-1}$, we define 
\[
A := \left( \sum_{1 \leq j \leq J} \lv x_j \rv^p + \Delta^p \right)^{1/p}.
\] 
Then there is a finite-precision algorithm, requiring work and storage at most $C \cdot J$, which computes a machine number $\widehat{A}$ that satisfies $\frac{1}{2} \cdot A  \leq \widehat{A} \leq 2 \cdot A$.
\end{lem}
\begin{proof}

All constants in the proof denoted by $C,C',$ etc., will depend only on $p$. We write $\Delta_1 \ll \Delta_2$ to indicate that $\Delta_1 \leq \Delta_2^C$ for a sufficiently large universal constant $C$. We set $\Delta_1 = \Delta_g^{C_0}$ for a sufficiently large universal constant $C_0 \in \N$ that will be determined later. Thus, in the recently introduced notation, we have  $\Delta_1 \ll \Delta_g$. By hypothesis, we are given a machine number $x_j^*$ with $\lv x_j^* - x_j \rv \leq \Delta_\epsilon$, and we guarantee that $\lvert x_j \rvert \leq \Delta_g^{-1}$  for each $j$. We define
\begin{equation}
\label{Bdefn}
B := \left(\sum_{\substack{1 \leq j \leq J \\ \lvert x_j^* \rvert \geq \Delta_1}} \lvert x_j^{*} \rvert^p + \Delta^p \right)^{1/p}.
\end{equation}
Note that
\[
\lv A^p - B^p \rv \leq \sum_{1 \leq j \leq J} \lv \lv x_j \rv^p - \lv x_j^* \rv^p \rv + \sum_{\substack{1 \leq j \leq J \\ \lvert x_j^* \rvert < \Delta_1}} \lv x_j \rv^p.
\]
Since $\lv x_j - x_j^* \rv \leq \Delta_\epsilon$ and $\lv x_j \rv, \lv x_j^* \rv \leq \Delta_g^{-C}$, the first sum is bounded by $C J \Delta_\epsilon \cdot \Delta_g^{-C'}$. Since $\lv x_j \rv \leq \Delta_1 + \Delta_\epsilon \leq 2 \Delta_1$ whenever $\lv x_j^*\rv < \Delta_1$, the second sum is bounded by $C J \Delta_1^p$. Thus, we have $\lv A ^p - B^p \rv \leq  C J \Delta_\epsilon \cdot \Delta_g^{-C} + C J \Delta_1^p$. We obtain the bound $\lv A ^p - B^p \rv  \leq  \Delta_g^{-C''} \Delta_1^p$ for a universal constant $C''$, because $J \leq \Delta_g^{-1}$ and because we may assume that $\Delta_\epsilon \leq \Delta_g^{C_0 p} = \Delta_1^p$.  Note that $A^p$ and $B^p$ are at least $\Delta^p$. Thus, by the mean value theorem, we have 
\begin{equation*}
\lv A - B \rv \leq \lv A^p - B^p \rv \cdot \max_{t \in [ \Delta^p,\infty)} \lv \frac{d}{dt}( t^{1/p}) \rv \leq \Delta_g^{-C''} \Delta_1^p \cdot \frac{1}{p} \Delta^{1-p} \leq \Delta_1^p \Delta_g^{-C'''}.
\end{equation*}
Here, in the last estimate we use that $\Delta \geq \Delta_g$. Note that $C', C'', C'''$ above are independent of $C_0$.

All the summands inside the parentheses in \eqref{Bdefn} are at least  $\Delta_1^p$ (recall that $\Delta \geq \Delta_g \geq \Delta_1$). Also, the number of summands is at most $J + 1 \leq C \Delta_g^{-1}  \leq C \Delta_1^{-1}$. Therefore, by the numerical stability of arithmetic (see Lemma \ref{num_stab}) we can compute a machine number $B_\fin$ such that
\[
\lv B - B_\fin \rv \leq \Delta_\epsilon \Delta_1^{-C}.
\]

We  conclude that $\lv A - B_\fin  \rv \leq \Delta_1^p \Delta_g^{-C'''} +   \Delta_\epsilon \Delta_1^{-C}$. 
We recall that $\Delta_1 = \Delta_g^{C_0}$ and $\Delta_\epsilon \ll \Delta_g$. So, if we pick a sufficiently large universal constant $C_0 \in \N$ then we  can guarantee that $\lv A - B_\fin \rv \leq  \frac{1}{2} \Delta_g$. Note that $A  \geq \Delta \geq \Delta_g$. Thus, we conclude that $A$ and $B_\fin$ differ by at most a factor of $2$. We can therefore define $\widehat{A} = B_\fin$ and the conclusion of the lemma follows.
\end{proof}


\section{Short Form}

Let $E = \{z_1,\cdots,z_N\} \subset \frac{1}{32} Q^\circ$, where $Q^\circ = [0,1)^n$. 

We write $\X(E)$ for the space of all real-valued functions $f$ on $E$, equipped with the trace norm induced by $\X = L^{m,p}(\R^n)$.

Recall that $\cP$ denotes the set of all polynomials  on $\R^n$ of degree at most $m-1$, and $\cM$ denotes the set of all multiindices $\alpha = (\alpha_1,\cdots,\alpha_n)$ with $\lv \alpha \rv \leq m -1$.

We let $\Delta_{\min} \leq \Delta_\epsilon \leq \Delta_g \leq \Delta_0$ be defined as in the \textbf{Main Assumptions} in Section \ref{fp_not}. In particular, recall that $\Delta_{\min} = 2^{-S}$ denotes the machine precision of our computer. When we refer to a ``machine number'' we will always mean an $S$-bit machine number.

Any linear functional $\omega : \X(E) \rightarrow \R$ can be expressed in the form
\begin{equation}\label{sform1}
\omega(f) = \sum_{\ell=1}^L \lambda_\ell \cdot f(z_{j_{\ell}}).
\end{equation}
We call \eqref{sform1} a \underline{short form} of $\omega$. We do not promise that the coefficients $\lambda_\ell$ are non-zero. Thus, in contrast to the notation in infinite-precision, a functional can have more than one short form. 
The \emph{depth} of $\omega$, denoted $\depth(\omega)$, is the number $L$.  Note that $\depth(\omega)$ depends on the short form \eqref{sform1} of $\omega$, and not on $\omega$ alone. This abuse of notation should cause no confusion. 

The short form \eqref{sform1} is given with parameters $(\Delta_g, \Delta_\epsilon)$ if  the numbers $\lambda_\ell$ are given with parameters $(\Delta_g,\Delta_\epsilon)$, and if the list $j_1,\cdots,j_L$ is given. Recall that this means we specify machine numbers $\overline{\lambda}_\ell$ with $\lv \lambda_\ell - \overline{\lambda}_\ell \rv \leq \Delta_\epsilon$, and we promise that $\lvert \lambda_\ell \rvert \leq \Delta_g^{-1}$ for each $\ell$. The indices $j_1,\cdots,j_L$ may be represented as pointers to the memory locations in which the corresponding points of $E$ are stored. We assume that each of these pointers is stored using a single unit of memory.

Let $\Omega = \{\omega_1,\cdots,\omega_{M}\}$ be a list of linear functionals on $\X(E)$.

A functional $\xi : \X(E) \rightarrow \R$ has $\Omega$-assisted depth $d$ provided that it can be written in the form
\begin{equation}
\label{sform2}
\xi(f) = \eta(f) + \sum_{\nu=1}^{\nu_{\max}} \mu_\nu \cdot \omega_{k_\nu}(f),
\end{equation}
where $\eta$ is a linear functional and $\depth(\eta) + \nu_{\max} \leq d$. We call \eqref{sform2} a \underline{short form} of $\xi$. Note that perhaps we can write a given $\xi$ in many different ways in short form.

The short form \eqref{sform2} is given with parameters $(\Delta_g, \Delta_\epsilon)$ in terms of assists $\Omega$ if the functional $\eta$ is given in short form with parameters $(\Delta_g,\Delta_\epsilon)$, the numbers  $\mu_\nu$ are given with parameters $(\Delta_g,\Delta_\epsilon)$, and a list of the indices $k_{\nu_1},\cdots,k_{\nu_{\max}}$ is given.

A functional $\xi : \X(E) \oplus \cP \rightarrow \R$ has $\Omega$-assisted depth $d$ provided that it can be written in the form
\begin{equation}
\label{sform3}
\xi(f,P) = \eta(f) + \sum_{\nu=1}^{\nu_{\max}} \mu_\nu \cdot \omega_{k_\nu}(f) + \sum_{\alpha \in \cM} \theta_\alpha \cdot \frac{1}{\alpha!} \partial^\alpha P(0),
\end{equation}
where $\eta$ is a linear functional and $\depth(\eta) + \nu_{\max} + \# (\cM) \leq d$. We call \eqref{sform3} a \underline{short form} of $\xi$.

The short form \eqref{sform3} is given with parameters $(\Delta_g, \Delta_\epsilon)$ in terms of assists $\Omega$ if  the functional $\eta$ is given in short form with parameters $(\Delta_g,\Delta_\epsilon)$, the numbers $\mu_\nu$ and $\theta_\alpha$ are given with parameters $(\Delta_g,\Delta_\epsilon)$, and a list of the indices $k_{\nu_1},\cdots,k_{\nu_{\max}}$ is given.

A linear map $T : \X(E) \oplus \cP \rightarrow \cP$ is given in short form with parameters $(\Delta_g, \Delta_\epsilon)$ in terms of assists $\Omega$, if for each $\beta \in \cM$ we exhibit a formula
\[ \partial^\beta ( T(f,P))(0) = \eta_\beta(f) + \sum_{\nu=1}^{\nu_{\max}} \mu_{\beta \nu} \cdot \omega_{k_\nu}(f) + \sum_{\alpha \in \cM} \theta_{\beta \alpha} \cdot \frac{1}{\alpha!} \partial^\alpha P(0),\]
where the functional $\eta_\beta$ is given in short form with parameters $(\Delta_g,\Delta_\epsilon)$, the numbers $\mu_{\beta \nu}$ and $\theta_{\beta \alpha}$ are given with parameters $(\Delta_g,\Delta_\epsilon)$, and a list of the indices $k_1,\cdots, k_{\nu_{\max}}$ is given.

Similarly, a linear map $T : \X(E) \rightarrow \cP$ is given in short form with parameters $(\Delta_g, \Delta_\epsilon)$ in terms of assists $\Omega$, if for each $\beta \in \cM$ we exhibit a formula
\[
\partial^\beta ( T(f))(0) = \eta_\beta(f) + \sum_{\nu=1}^{\nu_{\max}} \mu_{\beta \nu} \cdot \omega_{k_\nu}(f),
\] 
where the functional $\eta_\beta$ is given in short form with parameters $(\Delta_g,\Delta_\epsilon)$, the numbers $\mu_{\beta \nu}$ are given with parameters $(\Delta_g,\Delta_\epsilon)$, and a list of the indices $k_1,\cdots, k_{\nu_{\max}}$ is given.

We say we have computed a linear map $T : \X(E) \rightarrow \X$ in short form with parameters $(\Delta_g,\Delta_\epsilon)$ in terms of assists $\Omega$ if for each $\overline{S}$-bit machine point $\underline{x} \in Q^\circ $ and each multiindex $\alpha \in \cM$, we can compute a short form of the linear functional
\[
f \mapsto \partial^\alpha Tf(\underline{x})
\]
with parameters $(\Delta_g,\Delta_\epsilon)$ in terms of the assists $\Omega$. If the functional $f \mapsto \partial^\alpha ( Tf)(\underline{x})$ has $\Omega$-assisted depth $d$, for all $\underline{x} \in \R^n$ and $\alpha \in \cM$, then we say that the map $T$ has $\Omega$-assisted depth $d$. We extend this notation to linear maps $T : \X(E) \oplus \cP \rightarrow \X$ in the obvious way. We only answer queries if $\underline{x} \in Q^\circ$ because enormous $\underline{x}$'s might lead to overflow errors.

\section{Main Algorithms in Finite-Precision}

Our main theorem concerns extension operators for homogeneous Sobolev spaces and is stated below. Later, in Section \ref{inhom_finp_sec}, we will present a corresponding result for inhomogeneous Sobolev spaces (see Theorem \ref{main_thm2_fin}).

We write $c$, $C$, $C'$, etc., to denote universal constants, which depend only on $m$, $n$, and $p$.

Let $x = (x_1,\cdots,x_n) \in \R^n$. We call $x$ an $S_0$-bit ``machine point'' if each coordinate $x_k$ is an $S_0$-bit machine number.

\begin{thm} \label{main_thm_hom_fin}

There exists a universal constant $C \geq 1$  such that the following holds.

Let $\overline{S} \geq 1$ be an integer. 

We fix an $\overline{S}$-bit machine number $p > n$.

We also fix a subset $E \subset \frac{1}{32} Q^\circ$ consisting of $\overline{S}$-bit machine points, with $\#(E) = N \geq 2$, where $Q^\circ = [0,1)^n$.

We assume we are given constants $\Delta_{\min} = 2^{- K_{\max} \overline{S}}$, $\Delta_\epsilon^\circ := 2^{ - K_{1} \overline{S} }$, $\Delta_g^\circ := 2^{ - K_{2} \overline{S} }$, and $\Delta_\junk^\circ := 2^{ - K_{3} \overline{S} }$, for integers $K_1,K_2,K_3,K_{\max} \geq 1$ such that $K_{\max} \geq C \cdot K_1 \geq C^2 \cdot K_2 \geq C^3 \cdot K_3  \geq C^4$.

We assume that our computer can perform arithmetic operations on $S$-bit machine numbers with precision $\Delta_{\min} = 2^{-S}$, where  $S = K_{\max} \overline{S}$.

We compute (see below) lists $\Omega$ and $\Xi$, consisting of linear functionals on $\X(E) = \{ f : E \rightarrow \R\}$, with the following properties.
\begin{itemize}
\item The sum of $\depth(\omega)$ over all $\omega \in \Omega$ is bounded by $C N$, and $ \# \left[ \Xi \right] \leq C N$.
\item Each $\xi$ in $\Xi$ has $\Omega$-assisted depth at most $C$.
\item We compute each $\omega \in \Omega$ in short form with parameters $(\Delta_g^\circ, \Delta_\epsilon^\circ)$, and we compute each $\xi \in \Xi$ in short form with parameters $(\Delta_g^\circ, \Delta_\epsilon^\circ)$ in terms of the assists $\Omega$. 
\item We have
\[ C^{-1} \| f\|_{\X(E)} \leq \left[ \sum_{\xi \in \Xi} \lvert \xi(f) \rvert^p  \right]^{1/p}  \leq C \inf \left\{ \| F \|_{\X} + \Delta_{\junk}^\circ \| F\|_{L^p(Q^\circ)} : F \in \X, \; F = f \; \mbox{on} \; E \right\}\]
for every $f \in \X(E)$.
\end{itemize}

Moreover, there exists a linear map $T : \X(E) \rightarrow \X$ with the following properties.
\begin{itemize}
\item $T$ has $\Omega$-assisted depth at most $C$.
\item $T f = f$ on $E$, and 
\[  
\| Tf \|_{\X} \leq C \inf \left\{ \| F \|_\X + \Delta_\junk^\circ \| F \|_{L^p(Q^\circ)} : F  \in \X, \; F = f \; \mbox{on} \; E \right\}
\]
for every $f \in \X(E)$.
\item We produce a query algorithm with the following properties.

Given an $\overline{S}$-bit machine point $\underline{x} \in Q^\circ$, and given $\alpha \in \cM$, we compute a short form of the linear functional $ f \mapsto \partial^\alpha \left( Tf  \right)(\underline{x})$ in terms of the assists $\Omega$. This linear functional is computed with parameters $(\Delta_g^\circ,\Delta_\epsilon^\circ)$. This computation requires work at most $C \log N$.
\end{itemize}
The above computations require one-time work at most $C N \log N$ in space $C N$.

\end{thm}

\section{Bases for the Space of Polynomials}
\label{bases_sec_fin}

We examine the first algorithm in the infinite-precision text, namely the algorithm \textsc{Fit Basis to Convex Body} in Section \ref{sec_compbase}. This is a preparatory algorithm that is used in a later part of the text. We exhibit a finite-precision version of this algorithm by making several additional assumptions given below.

We impose the assumptions in Theorem \ref{main_thm_hom_fin}. In particular, we assume that $\Delta_{\min} \leq \Delta_\epsilon \leq \Delta_g \leq \Delta_0$ are as in the \textbf{Main Assumptions} in Section \ref{fp_not}. In particular, our computer can perform arithmetic operations on $S$-bit machine numbers with precision $\Delta_{\min} = 2^{-S}$, where $S = K_{\max} \cdot \overline{S}$. 

We introduce a few conventions that are used in the rest of the paper. A machine number will mean an $S$-bit machine number, and a machine point will mean an $S$-bit machine point. A bounded interval $I \subset \R$ is called a machine interval if its endpoints are machine numbers. 

We let $p>n$ be an $\overline{S}$-bit machine number.

\begin{equation}
\label{v0} \text{We assume that } x \in \R^n \text{ is an } (S \mbox{-bit}) \;\text{machine point}.
\end{equation}

Recall that $\cP$ is the vector space of polynomials on $\R^n$ of degree at most $m-1$, and $D = \dim \cP$. We identify $\cP$ with $\R^D$, by identifying $P \in \cP$ with $\left(\frac{1}{\alpha!} \partial^\alpha P(x) \right)_{\alpha \in \cM}$, where $\cM$ denotes the set of all multiindices of order at most $m-1$. We define 
\[ \lvert P \rvert_{x} = \left( \sum_{\alpha \in \cM} \lvert \partial^\alpha P(x) \rvert^p \right)^{1/p}.
\]

\comments{
p. 262: Added line
}

We assume we are given $\Lambda \geq 1$. We write $c(\Lambda)$, $C(\Lambda)$, etc. to denote constants depending on $m$, $n$, $p$, and $\Lambda$. We write $c$, $\widetilde{c}$, $C$, etc. to denote constants depending only on $m$, $n$, and $p$.

We are given a quadratic form $q$ on $\cP$. We assume that $q$ is specified as a $D$ x $D$ matrix $(q_{\beta\gamma})_{\beta, \gamma \in \cM}$ acting on the above coordinates:
\begin{equation}
\label{coordspec}
q(P) = \sum_{\beta, \gamma \in \cM} q_{\beta \gamma} \cdot \partial^\beta P(x) \cdot \partial^\gamma P(x).
\end{equation}

We assume that the matrix $(q_{\alpha \beta})_{\alpha, \beta \in \cM}$ satisfies the following conditions:
\begin{align}
&\label{v1} \text{The entries are such that} \; \lvert q_{\alpha \beta} \rvert \leq \Delta_g^{-1} \text{ and} \; q_{\alpha \beta} \; \text{is specified to precision } \Delta_{\epsilon}. \\
\label{v2}
&\text{We have } (q_{\alpha \beta}) \geq \Delta_g \cdot (\delta_{\alpha \beta}).
\end{align}

From  \eqref{v2} we learn that
\begin{equation}
\label{v3}
\lv q(P) \rv \geq c \Delta_g \cdot \lvert P \rvert_{x}^2 \;\; \mbox{for all} \; P \in \cP,
\end{equation}
for a universal constant $c > 0$.

We fix a multiindex set $\cA \subset \cM$. 

The main result of the section is as follows.

\environmentA{Algorithm: Fit Basis to Convex Body (Finite-Precision Version).}

Given $q$, $x$, $\cA$ as above: We compute a partition of $[\Delta_g,\Delta_g^{-1}]$ into machine intervals $I_\ell$, and for each $\ell$ we compute machine numbers $\lambda_\ell,c_\ell$ with $c_\ell \geq 0$, such that the function $\eta_* : \left[ \Delta_g,\Delta_g^{-1} \right] \rightarrow \R$, defined by
\[
\eta_*(\delta) := c_\ell \cdot \delta^{\lambda_\ell}  \;\; \mbox{for} \; \delta \in I_\ell,
\]
has the following properties.
\begin{itemize}
\item Let $\ooline{\sigma}$ satisfy $\{ q \leq \Lambda^{-1} \} \subset \ooline{\sigma} \subset \{ q \leq \Lambda \}$. Then, for any $\delta \in [ \Delta_g,\Delta_g^{-1}]$, 
\begin{itemize}
\item $\ooline{\sigma}$ has an $(\cA,x,\eta^{1/2},\delta)$-basis for any $\eta > C(\Lambda) \cdot \eta_*(\delta)$, 
\item $\ooline{\sigma}$ does not have an $(\cA,x,\eta^{1/2},\delta)$-basis for any  $\eta <  c(\Lambda) \cdot \eta_*(\delta) $.
\end{itemize}
\item Moreover, $c \cdot \eta_*(\delta_1)  \leq \eta_*(\delta_2) \leq C \cdot \eta_*(\delta_1)$ whenever $\frac{1}{10} \delta_1 \leq \delta_2 \leq 10 \delta_1$.
\item Also, $\eta_*(\delta) \geq \Delta_g^{C}$ for any $\delta \in [ \Delta_g,\Delta_g^{-1}]$.
\item The numbers $c_\ell$ belong to the interval $\left[ \Delta_g^C, \Delta_g^{-C} \right]$, and the exponents $\lambda_\ell$ are of the form $\mu + \nu/p$ for integers $\mu,\nu$ with $\lv \mu \rv, \lv \nu \rv \leq C$.
\item The computation of $I_\ell$, $\lambda_\ell$, and $c_\ell$, requires work and storage at most $C$.
\end{itemize}
Here, $c > 0$ and $C \geq 1$ are constants determined by $m$,$n$, and $p$, while $c(\Lambda)$ and $C(\Lambda)$ are constants determined by $m$,$n$,$p$, and $\Lambda$.

\begin{proof}[\underline{Explanation}]

We first define a rational function $\eta_{\min}(\delta)$ as in Section \ref{sec_compbase}, and then explain how to compute an approximation $\eta_*(\delta)$ to $\eta_{\min}(\delta)$. We will need to estimate the numerical stability of the computation with respect to rounding error.

We recall several of the main definitions in Section \ref{sec_compbase}.

We study the quadratic form
\begin{align}
\label{quadfm}
M^\delta ((P_\alpha)_{\alpha \in \cA}) &:= \sum_{\alpha \in \cA} q(\delta^{m - n/p - |\alpha|} P_\alpha) + \sum_{\substack{\alpha \in \cA, \beta \in \cM \\ \beta > \alpha}} (\delta^{|\beta| - |\alpha|} \partial^\beta P_\alpha(x))^2 \\
\notag{}
&= \sum_{\alpha \in \cA} \sum_{\beta, \gamma \in \cM} \delta^{2(m - n/p - |\alpha|)} q_{\beta \gamma } \cdot \partial^\beta P_\alpha(x) \cdot \partial^\gamma P_\alpha(x) + \sum_{\substack{\alpha \in \cA, \beta \in \cM \\ \beta > \alpha}} (\delta^{|\beta| - |\alpha|} \partial^\beta P_\alpha(x) )^2,
\end{align}
restricted to the affine subspace
\begin{equation*}
H := \left\{ \vec{P} = (P_\alpha)_{\alpha \in \cA} : \partial^\beta P_\alpha(x) = \delta_{\beta \alpha} \; \mbox{for} \; \alpha, \beta \in \cA \right\}.
\end{equation*}
We define
\[\eta_{\min}(\delta) := \min_{\vec{P} \in H} M^\delta(\vec{P}), \] 
which we regard as a function of $\delta \in \left[ \Delta_g, \Delta_g^{-1} \right]$.

Recall from \eqref{slowvariance} and \eqref{abc2} that
\begin{equation}
\label{slowvariance_new}
\eta_{\min}(\delta_1) \leq \eta_{\min}(\delta_2) \leq \left(\frac{\delta_2}{\delta_1} \right)^{2m} \eta_{\min}(\delta_1) \;\; \mbox{for} \;\; \delta_1 \leq \delta_2,
\end{equation}
and that
\begin{equation}
\label{abc2_new}
\ooline{\sigma} \mbox{ has a } (\cA,x,\eta^{1/2},\delta)\mbox{-basis if } \eta > C(\Lambda) \cdot \eta_{\min}(\delta)\mbox{, but not if } \eta < c(\Lambda) \cdot \eta_{\min}(\delta).
\end{equation}

Thanks to \eqref{v3}, for $\delta \in \left[ \Delta_g ,  \Delta_g^{-1} \right]$ we have
\begin{equation}
\label{wef1}
\lvert M^\delta(\vec{P}) \rvert \geq c \Delta_{g}^{2m+1} \cdot \sum_{\alpha \in \cA}  \lvert P_\alpha \rvert_{x}^2 \qquad \mbox{for any} \; \vec{P} = (P_\alpha)_{\alpha \in \cA} \in H.
\end{equation}
Furthermore, for all $\vec{P} \in H$ and $\alpha \in \cA$, we have $\partial^\alpha P_\alpha(x) = 1$, hence $\lvert P_\alpha \rvert_x \geq 1$. Thus, since $\eta_{\min}(\delta)$ is the minimum of $M^\delta(\cdot)$ on $H$, we have
\begin{equation}
\label{lb0}
\eta_{\min}(\delta) \geq \widetilde{c} \Delta_g^{2m+1},
\end{equation}
for a universal constant $\widetilde{c}>0$.

We will compute a piecewise rational function $\widetilde{\eta}_{\min}(\delta)$ that approximates $\eta_{\min}(\delta)$. For $w = (w_{\alpha \beta})_{\alpha \in \cA, \beta \in \cM \setminus \cA} \in \R^J$ we set
\begin{equation}
\label{coordspec2}
P^w_\alpha(z) := \frac{1}{\alpha !} \cdot (z - x)^\alpha + \sum_{\beta \in \cM \setminus \cA} \frac{1}{\beta!} \cdot w_{\alpha \beta} \cdot (z-x)^\beta \qquad (\alpha \in \cA).
\end{equation}
This gives a coordinate mapping $w \mapsto \vec{P}^w = (P^w_\alpha)_{\alpha \in \cA} \in H$, and we set
\begin{align}
\label{coordspec3}
\widetilde{M}^\delta(w) &:= M^\delta(\vec{P}^w) \\
\notag{}
&= \langle A^\delta w, w \rangle - 2 \langle b^\delta, w \rangle + m^\delta \quad (w \in \R^J). 
\end{align}
Here, $A^\delta$ is a matrix, $b^\delta$ is a vector, $m^\delta$ is a scalar - all functions of $\delta$ - and $\langle \cdot, \cdot \rangle$ denotes the standard Euclidean inner product on $\R^J$. The entries of $A^\delta$, $b^\delta$, and $m^\delta$ are all sums of monomials of the form $a \cdot \delta^{\mu + \nu/p}$ with $\mu,\nu \in \Z$ and $a \in \R$.

We have $\| w \|^2 =  \sum_{\alpha,\beta} \lvert w_{\alpha \beta} \rvert^2 \leq c \sum_{\alpha} \lvert P^w_\alpha \rvert_x^2$, since $w_{\alpha \beta} = \partial^\beta P_\alpha^w(x)$ for $\alpha \in \cA$, $\beta \in \cM \setminus \cA$. Thus, from \eqref{wef1} we have $\lv \widetilde{M}^\delta(w) \rv \geq c \Delta_g^{2m+1} \cdot \| w \|^2$, hence
\begin{equation}
\label{wef2}
A^\delta \geq c \Delta_g^{2m+1} \cdot (\delta_{ij}),
\end{equation}
for a universal constant $c>0$. In particular, the matrix $A^\delta$ is invertible.

Recall that $A^\delta = (A^\delta_{ij})$, $b^\delta = (b^\delta_i)$, and $m^\delta$ are given in the form
\begin{equation}
\label{we2_new}
\left\{
\begin{aligned}
A_{ij}^\delta &= \sum_{\mu,\nu} c^{ij}_{\mu \nu} \delta^{\mu + \nu/p} \qquad (1 \leq i, j \leq J), \\
b_j^\delta & = \sum_{\mu,\nu} c^{j}_{\mu \nu} \delta^{\mu + \nu/p} \qquad (1 \leq j \leq J), \\
m^\delta & = \sum_{\mu,\nu} c_{\mu \nu} \delta^{\mu + \nu/p}.
\end{aligned}
\right.
\end{equation}
There are at most $C$ pairs $(\mu,\nu) \in \Z \times \Z$ relevant to the above sums, and $\mu,\nu$ are bounded in magnitude by a universal constant $C$. (See \eqref{we2.1} - \eqref{we2.3}.) 

We insert the formula \eqref{coordspec2} for the polynomials $P_\alpha = P^w_\alpha$ ($\alpha \in \cA$) in the second line of the definition \eqref{quadfm} of $M_\delta$ to produce the expression $\widetilde{M}^\delta(w) = \langle A^\delta w, w \rangle - 2 \langle b^\delta, w \rangle + m^\delta$. We compute each of the numbers $c^{ij}_{\mu \nu}$, $c^j_{\mu \nu}$, and $c_{\mu \nu}$ as a linear combination of the entries of $(q_{\alpha \beta})$ (and the constant $1$). Hence, since the $q_{\alpha \beta}$ are given with parameters $(\Delta_g,\Delta_\epsilon)$, the numbers $c^{ij}_{\mu \nu}$, $c^j_{\mu \nu}$, $c_{\mu \nu}$ can be computed with parameters $(\Delta_g^C,\Delta_g^{-C} \Delta_\epsilon)$.

We can compare exponents of the form $\lambda = \mu + \nu/p$ and $\overline{\lambda} = \overline{\mu} + \overline{\nu}/p$ by expressing $\lambda$, $\overline{\lambda}$ as a ratio of integers and cross-multiplying (recall that $p$ is an $\overline{S}$-bit machine number). This comparison requires at most $C$ units of work. By summing the coefficients of the monomials with the same exponent, we may assume that the exponents $\mu + \nu/p$ in \eqref{we2_new} are pairwise distinct. 

We compute a formula for the inverse matrix $(A^{\delta})^{-1}$ by applying Cramer's rule. Hence,
\begin{equation}
\label{wee0}
(A^\delta)^{-1}_{ij} = \frac{[A^\delta]_{ij}}{\det(A^\delta)}  =  \frac{\sum_{k}  a^{ i j}_{k} \cdot \delta^{\lambda_k}} {\sum_{\ell} b_{\ell} \cdot \delta^{\gamma_\ell} } \qquad (1 \leq i ,j \leq J),
\end{equation}
where $[A^{\delta}]_{ij}$ denotes the $(i,j)$-cofactor of the matrix $A^{\delta}$. The number of terms in the sums in the numerator and denominator of \eqref{wee0} is bounded by $C$.

We compute the numbers $a_k^{ij}$ and $b_\ell$ in \eqref{wee0} by evaluating a polynomial function of the coefficients $c^{ij}_{\mu \nu}$ arising in the entries of the matrix $(A^\delta_{ij})$ in \eqref{we2_new}. The numbers $c_{\mu \nu}^{ij}$ are given with parameters $(\Delta_g^C, \Delta_g^{-C} \Delta_\epsilon)$, so we can compute $a_k^{ij}$ and $b_\ell$ with parameters $(\Delta_g^C,\Delta_g^{-C}\Delta_\epsilon)$. 

The exponents $\lambda_k$ and $\gamma_\ell$ in \eqref{wee0} have the form $\mu + \nu/p$, where $\mu,\nu \in \Z$ are bounded in magnitude by a universal constant $C$.

We compute an expression for $\eta_{\min}(\delta) = \min_w \widetilde{M}^\delta(w)$ as follows. Note that the quadratic function $\widetilde{M}^\delta(w)$ in \eqref{coordspec3} achieves its minimum at $w^\delta := (A^{\delta})^{-1} b^\delta$. Thus,
\begin{align*}
\eta_{\min}(\delta) &= \widetilde{M}^\delta(w^\delta) \\
\notag{}
& = - \langle b^\delta, (A^\delta)^{-1} b^\delta \rangle + m^\delta \\
\notag{}
&= \sum_{i,j=1}^J b_i^\delta  \cdot (A^\delta)^{-1}_{ij} \cdot b_j^\delta + m^\delta.
\end{align*}
Inserting the expressions for the entries of $(A^\delta)^{-1}$, $b^\delta$, and the expression for $m^\delta$, we compute (see below) a rational expression
\begin{equation}
\label{rational2}
\eta_{\min}(\delta) = \frac{ \sum_k a_k \cdot \delta^{\lambda_k} }{   \sum_{\ell} b_{\ell} \cdot \delta^{\gamma_\ell}} = \frac{N(\delta)}{D(\delta)}. 
\end{equation}
The number of terms in the sums in the numerator and denominator of \eqref{rational2} is bounded by $C$.

The denominator $D(\delta) = \sum_\ell b_\ell \cdot \delta^{\gamma_\ell}$ in \eqref{rational2} is the same as the common denominator in the expression for $(A^\delta)^{-1}_{ij}$ in \eqref{wee0}, namely $\det (A^\delta)$. From \eqref{wef2} we have $\det(A^\delta) \geq \Delta_g^C$, hence
\begin{equation}
\label{lb1}
\sum_\ell b_{\ell} \cdot \delta^{\gamma_\ell} \geq \Delta_g^C.
\end{equation}

The exponents $\gamma_\ell$ and $\lambda_k$ in \eqref{rational2} have the form $\mu + \nu/p$, where $\mu,\nu \in \Z$ are bounded in magnitude by a universal constant. We can assume that the $\gamma_\ell$ are distinct, as are the $\lambda_k$. (We combine all the monomials in the numerator or denominator that have the same exponent.)

The numbers $a_k$ in the numerator in \eqref{rational2} are defined by evaluating a polynomial function on the coefficients $a_k^{ij}$, $b_\ell$ in $(A^\delta)^{-1}_{ij}$ (see \eqref{wee0}), and the coefficients $c_{\mu \nu}^j$ and $c_{\mu\nu}$ in $b^\delta_j$ and $m^\delta$ (see \eqref{we2_new}). Thus, we can compute $a_k$ with parameters $(\Delta_g^{{C_1}}, \Delta_g^{-C_1} \Delta_\epsilon)$ for a large enough universal constant $C_1$.

As explained before, the numbers $b_\ell$ are given with parameters $(\Delta_g^{C_1},\Delta_g^{-C_1}\Delta_\epsilon)$.

The exponents $\lambda_k$ in \eqref{rational2} are pairwise distinct and have the form $\mu + \nu/p$ for integers $\mu,\nu \in \Z$ with $\lv \mu \rv, \lv \nu \rv \leq C$. The same is true of the exponents $\gamma_\ell$. Hence, 
\begin{equation}
\label{mach_comp3}
\lvert \lambda_k - \lambda_{k'} \rvert \geq c_0, \;\; \mbox{and} \;\; \lvert \gamma_\ell - \gamma_{\ell'} \rvert \geq c_0
\end{equation}
for all $k \neq k'$ and $\ell \neq \ell'$.\footnote{The constant $c_0$ here depends only on $m$,$n$,$p$, but it may depend sensitively on the approximation of $\frac{1}{p}$ by rationals with low denominators.}

We now perform a crucial rounding step.

We introduce a parameter $\Delta = 2^{-S_1}$ of the form $\Delta = \Delta_g^{C_2}$, for  $C_2 \in \N$ that is assumed to be greater than a large enough universal constant. We will later fix $C_2$ to be a universal constant, but not yet. We assume $\Delta_\epsilon \leq \Delta_g^{C_1 + C_2}$, hence
\begin{equation}
\label{delta1} 
\begin{aligned}
&\Delta_g^{-C_1} \Delta_\epsilon \leq \Delta.
\end{aligned}
\end{equation}

The numbers $a_k$ and $b_\ell$ are given with parameters $(\Delta_g^{C_1},\Delta_g^{-C_1} \Delta_\epsilon)$, so we can compute $S_1$-bit machine numbers $\widetilde{a}_k$ and $\widetilde{b}_\ell$ with
\begin{equation}
\label{mach_comp}
\left\{
\begin{aligned}
& \lvert a_k - \widetilde{a}_k  \rvert \leq \Delta, \qquad \lvert b_\ell - \widetilde{b}_\ell\rvert \leq \Delta, \\
& \lvert \widetilde{a}_k \rvert \leq \Delta_g^{-C}, \qquad\quad  \lvert \widetilde{b}_\ell \rvert \leq \Delta_g^{-C}.
\end{aligned}
\right.
\end{equation}

We set
\begin{equation}
\label{roundedeta}
\widetilde{\eta}_{\min}(\delta) = \frac{\sum_k \widetilde{a}_k \delta^{\lambda_k} }{\sum_\ell  \widetilde{b}_\ell \delta^{\gamma_\ell}} = \frac{\til{N}(\delta)}{\til{D}(\delta)}.
\end{equation}

We use \eqref{mach_comp} and the fact that $\lambda_k$ and $\gamma_\ell$ are bounded by $C$ to estimate the difference between $\eta_{\min}(\delta)$ and $\widetilde{\eta}_{\min}(\delta)$. For $\delta \in \left[ \Delta_g, \Delta_g^{-1} \right]$, we have
\begin{equation*}
\lv a_k \delta^{\lambda_k} - \widetilde{a}_k \delta^{\lambda_k} \rv  \leq  \Delta \Delta_g^{-C}.
\end{equation*}
Hence,
\[
\lv N(\delta) - \til{N}(\delta) \rv = \lv  \sum_k a_k \delta^{\lambda_k}  - \sum_k \til{a}_k \delta^{\lambda_k} \rv \leq C  \Delta \Delta_g^{-C} \leq  \Delta \Delta_g^{-C'}.
\]
Moreover,
\[
\lv N(\delta) \rv
= \lv \sum_k a_k \delta^{\lambda_k} \rv
\leq \sum_k \lv a_k \rv \lv \delta^{\lambda_k} \rv \leq C \Delta_g^{-C} \cdot \Delta_g^{-C} \leq \Delta_g^{-C'}.
\]
Similarly,
\[
\lv D(\delta) - \til{D}(\delta) \rv
= \lv 
\sum_\ell b_\ell \delta^{\gamma_\ell}       - 
\sum_\ell \til{b}_\ell \delta^{\gamma_\ell} 
\rv 
\leq  \Delta \Delta_g^{-C'}.
\]
Moreover,
\begin{align*}
 D(\delta) &= \sum_\ell b_\ell \delta^{\gamma_\ell} \ogeq{\eqref{lb1}} \Delta_g^C, \;\; \mbox{and} \; \\
\lv \til{D}(\delta) \rv &\geq \lv D(\delta) \rv - \lv D(\delta) - \til{D}(\delta) \rv
\\
&\geq \Delta_g^C - \Delta \Delta_g^{-C' } \geq \frac{1}{2} \Delta_g^{C},
\end{align*}
since $\Delta \Delta_g^{-C'} = \Delta_g^{C_2 - C'} \leq \frac{1}{2} \Delta_g^C$, for large enough $C_2$.

Using the previous estimates, we have
\begin{align*}
\lv \eta_{\min}(\delta) -  \widetilde{\eta}_{\min}(\delta) \rv  &=  \lv  \frac{N(\delta)}{D(\delta)} - \frac{\til{N}(\delta)}{\til{D}(\delta)} \rv \\
& \leq  \left\lvert \frac{N(\delta) \cdot \left(  D(\delta) - \til{D}(\delta) \right)}{D(\delta) \cdot \til{D}(\delta)}  \right\rvert + \lv \frac{N(\delta) - \til{N}(\delta)}{\til{D}(\delta)} \rv \\
& \leq \frac{\Delta_g^{-C'} \cdot  \Delta \Delta_g^{-C'}}{\Delta_g^{C} \cdot \frac{1}{2} \Delta_g^C} + \frac{\Delta \Delta_g^{-C'} }{\frac{1}{2} \Delta_g^{C}} \leq \Delta \Delta_g^{-C''}.
\end{align*}
Thus,
\begin{align*}
\lv \eta_{\min}(\delta) - \widetilde{\eta}_{\min}(\delta) \rv 
& \leq \Delta \Delta_g^{-C''}  \\
&\leq \widetilde{c} \Delta_g^{2m+2}
\oleq{\eqref{lb0}} \frac{1}{2} \cdot \eta_{\min}(\delta),
\end{align*} 
where we may choose $C_2$ large enough so that $ \Delta \Delta_g^{-C''} = \Delta_g^{C_2 - C''} \leq \widetilde{c} \Delta_g^{2m+2}$ with $\widetilde{c}$ as in \eqref{lb0}. We fix $C_2 \in \N$ to be a universal constant satisfying the previous bounds.

Therefore, thanks to \eqref{lb0} we have
\begin{equation}
\label{uws2a}
\Delta_g^C \leq \frac{1}{2} \cdot \eta_{\min}(\delta) \leq \widetilde{\eta}_{\min}(\delta) \leq 2 \cdot \eta_{\min}(\delta) \qquad (\Delta_g \leq \delta \leq \Delta_g^{-1}).
\end{equation}

We assume that none of the coefficients in the expression \eqref{roundedeta} are equal to zero, for otherwise we could discard the vanishing terms. Since $\til{a}_k$ and $\til{b}_\ell$ are $S_1$-bit machine numbers and $\Delta = 2^{-S_1}$, this means that
\begin{equation}
\label{mach_comp1}
\lvert \widetilde{a}_k \rvert \geq \Delta = \Delta_g^{C_2}, \;\;\; \lvert \widetilde{b}_\ell \rvert \geq \Delta = \Delta_g^{C_2} \quad \text{for all} \; k,\ell.
\end{equation}

We wish to compute a piecewise monomial function $\eta_*(\delta)$ that differs from $\eta_{\min}(\delta)$ by at most a universal constant factor. The first, second, and third bullet points in the statement of the algorithm \textsc{Fit Basis to Convex Body} (finite-precision) are consequences of \eqref{abc2_new}, \eqref{slowvariance_new}, and \eqref{lb0}, respectively. The guarantees in the fourth bullet point follow by examining the construction in the result below.

\environmentA{Procedure: Approximate Rational Function.}

We are given machine numbers $\widetilde{a}_k, \widetilde{b}_\ell$ satisfying 
\[
\lvert \til{a}_k \rvert, \; \lvert \til{b}_\ell \rvert \in \left[ \Delta_g^C, \Delta_g^{-C} \right].
\]
We are given numbers $\lambda_k$ and $\gamma_\ell$ of the form $\mu + \nu/p$, for integers $\mu,\nu$ with $\lv \mu \rv, \lv \nu \rv \leq C$, such that
\[
\lvert \lambda_k - \lambda_{k'} \rvert \geq c_0, \; \lvert \gamma_\ell - \gamma_{\ell'} \rvert \geq c_0 \qquad \mbox{for all} \; k \neq k', \; \ell \neq \ell'.
\]
Let
\[
\widetilde{\eta}_{\min}(\delta) = \frac{\sum_k  \widetilde{a}_k \delta^{\lambda_k} }{\sum_\ell \widetilde{b}_{\ell} \delta^{\gamma_\ell}}.
\]
Assume that the number of summands in the numerator and denominator is bounded by a universal constant $C$. Suppose that there exists a function $\eta_{\min}(\delta)$ satisfying  \eqref{slowvariance_new} and  \eqref{uws2a}.

We compute machine intervals $I_\ell$, machine numbers $d_\ell$, and numbers $\omega_\ell$, such that $\left[ \Delta_g,\Delta_g^{-1} \right]$ is the disjoint union of the $I_\ell$, and such that the function $\eta_* : \left[ \Delta_g, \Delta_g^{-1} \right] \rightarrow \R$, defined by
\[
\eta_*(\delta) :=  d_\ell \cdot \delta^{\omega_\ell} \quad \; \mbox{for} \; \delta \in I_\ell,
\]
satisfies
\begin{equation*}
c \cdot \eta_*(\delta)  \leq \eta_{\min}(\delta) \leq C \cdot \eta_*(\delta) \quad \mbox{for all} \; \delta \in \left[ \Delta_g,\Delta_g^{-1} \right] .
\end{equation*}
Here, $c$ and $C$ are universal constants.

The numbers $\omega_\ell$ are of the form $\mu + \nu/p$ for integers $\mu,\nu$ with $\lv \mu \rv, \lv \nu \rv \leq C$.

The machine numbers $d_\ell$ are contained in the interval $\left[ \Delta_g^C, \Delta_g^{-C} \right]$.

This computation requires work and storage at most $C$.

\begin{proof}[\underline{Explanation}]

We define
\begin{align*}
\cB &:= \bigcup_{k \neq k'} I_{k k'}, \;\; \mbox{where} \\
& \qquad\qquad I_{k k'} := \left\{ \delta \in \left[\Delta_g , \Delta_g^{-1} \right]  : 5^{-1} \cdot \lvert \widetilde{a}_k \delta^{\lambda_k} \rvert  \leq \lvert \widetilde{a}_{k'} \delta^{\lambda_{k'}} \rvert \leq  5 \cdot \lvert \widetilde{a}_k \delta^{\lambda_k} \rvert  \right\},
\end{align*}
and similarly
\begin{align*}
\cC &:= \bigcup_{\ell \neq \ell'} J_{\ell \ell'}, \;\; \mbox{where} \\
& \qquad\qquad J_{\ell \ell'} := \left\{ \delta \in \left[\Delta_g , \Delta_g^{-1} \right]  :  5^{-1} \cdot  \lvert \widetilde{b}_\ell \delta^{\gamma_\ell} \rvert \leq \lvert \widetilde{b}_{\ell'} \delta^{\gamma_{\ell'}} \rvert  \leq  5 \cdot \lvert \widetilde{b}_\ell \delta^{\gamma_\ell} \rvert  \right\}.
\end{align*}

For any interval $I \subset \left[ \Delta_g, \Delta_g^{-1} \right] \setminus (\cB \cup \cC)$, we have
\begin{equation}
\label{pickmon1_fin}
\left\{ \;\;
\begin{aligned}
& \mbox{there exist unique }  k=k(I) \in \{1,\cdots,K\} \; \mbox{and} \; \ell = \ell(I) \in \{1,\cdots,L \} \\
& \; \mbox{such that}\; \lvert \til{a}_k \delta^{\lambda_k} \rvert >  2 \sum_{k' \neq k} \lvert \til{a}_{k'} \delta^{\lambda_{k'}} \rvert \; \mbox{and} \; \lvert \til{b}_{\ell} \delta^{\gamma_{\ell}} \rvert >  2 \sum_{\ell' \neq \ell} \lvert \til{b}_{\ell'} \delta^{\gamma_{\ell'}}\rvert \; \mbox{for all} \; \delta \in I.
\end{aligned}
\right.
\end{equation}
Moreover, $\int_{\cB \cup \cC} \frac{dt}{t} \leq 2A$ for a universal constant $A$. (This follows by the same reasoning used to prove \eqref{pickmon1} and \eqref{logbound}.)

To compute the endpoints of a nonempty interval $I_{k k'} = [h^-_{k k'}, h_{k k'}^+]$ ($k \neq k'$) we solve the equations
\[
\delta^{\lambda_k - \lambda_{k'}} = 5^{\pm 1} \frac{\lv  \til{a}_{k'} \rv}{\lv \til{a}_k \rv}.
\]
The solutions $\delta = \delta_\pm$ are given by
\[
\delta_{\pm} = \left( 5^{\pm 1} \frac{\lv  \til{a}_{k'} \rv}{\lv \til{a}_k \rv} \right)^{(\lambda_k - \lambda_{k'})^{-1}} \qquad (\mbox{for each choice of} \; \pm ).
\]
From the lower/upper bound on $\lv \til{a}_k \rv$ by $\Delta_g^{\pm C}$, we see that we can compute $\lv \til{a}_{k'}/ \til{a}_k \rv$ to precision $\Delta_g^{-C} \Delta_\epsilon$, and $\lv \til{a}_{k'}/ \til{a}_k \rv \in [\Delta^{-C}_g,\Delta^C_g]$. Since $\lv (\lambda_k - \lambda_{k'})^{-1} \rv \leq c_0^{-1} \leq C$, we can compute $\delta_{+}$ and $\delta_-$ to precision $\Delta_g^{-C} \Delta_\epsilon$, due to the numerical stability of exponentiation.

Now, note that
\[
h_{kk'}^- = \min \{ \delta_- , \delta_+, \Delta_g \}, \qquad h_{kk'}^+ = \max\{ \delta_-, \delta_+, \Delta_g^{-1} \}.
\]
Both $h_{kk'}^-$ and $h_{kk'}^+$ can be computed with parameters $(\Delta_g, \Delta_g^{-C} \Delta_\epsilon)$. Thus, we can compute a machine interval $\widetilde{I}_{kk'} \subset [\Delta_g,\Delta_g^{-1}]$ with $I_{kk'} \subset \widetilde{I}_{kk'}$ and
\begin{equation}\label{haus1}
 \dist( I_{kk'}, \widetilde{I}_{kk'} ) \leq \Delta_g^{-C} \Delta_{\epsilon},
\end{equation}
where $\dist(\cdot,\cdot)$ is the Hausdorff distance.  Due to the previous inclusion, we know that the union of the intervals $\til{I}_{kk'}$ contains the set $\cB = \cup_{k \neq k'} I_{k k'}$.

Similarly, we compute a machine interval $\widetilde{J}_{\ell \ell'} \subset [\Delta_g,\Delta_g^{-1}]$ with $J_{\ell \ell'} \subset \widetilde{J}_{\ell \ell'}$ and 
\begin{equation}\label{haus2}
 \dist( J_{\ell \ell'}, \widetilde{J}_{\ell \ell'} ) \leq \Delta_g^{-C} \Delta_{\epsilon}.
\end{equation}
Again, note that the union of the intervals $\til{J}_{\ell \ell'}$ contains the set $\cC$.

We next compute pairwise disjoint machine intervals $I_\nu^\bad \subset \left[\Delta_g,\Delta_g^{-1} \right]$ such that
\begin{align*}
\bigcup_{\nu=1}^{\nu_{\max}} I_\nu^\bad =  \bigcup_{k \neq k'} \widetilde{I}_{k k'} \cup \bigcup_{\ell \neq \ell'} \widetilde{J}_{\ell \ell'}.
\end{align*} 
We form the intervals $I_\nu^{\bad}$ by concatenating the intersecting intervals among $\widetilde{I}_{k k'} $ and $\widetilde{J}_{\ell \ell'}$. Note that the union of the $I_\nu^\bad$ contains the set $\cB \cup \cC$.

Because the intervals below are contained in $\left[ \Delta_g, \Delta_g^{-1} \right]$, for each $\nu$ we have 
\begin{align}
\label{fff1}
 \int_{I^\bad_\nu} \frac{dt}{t} &\leq \sum_{k,k'} \int_{\til{I}_{kk'}} \frac{dt}{t} + \sum_{\ell, \ell'} \int_{\til{J}_{\ell \ell'}} \frac{dt}{t} \\
 \notag{}
 & \oleq{\eqref{haus1}, \eqref{haus2}} \sum_{k,k'} \left[ \int_{I_{kk'}} \frac{dt}{t} + \Delta_g^{-1} \cdot \Delta_g^{-C} \Delta_\epsilon \right] + \sum_{\ell, \ell'} \left[  \int_{J_{\ell \ell'}} \frac{dt}{t} + \Delta_g^{-1} \Delta_g^{-C} \Delta_\epsilon \right] \\
 \notag{}
 &\leq  \int_{\cB \cup \cC} \frac{dt}{t}  +  \Delta_g^{-C'} \Delta_\epsilon \leq 3 A.
 \end{align}
Recall that $A \geq 1$ is a universal constant.

We compute pairwise disjoint machine intervals $I_\mu \subset \left[\Delta_g,\Delta_g^{-1} \right]$ such that 
\[ 
\bigcup_{\mu=1}^{\mu_{\max}} I_\mu = \left[\Delta_g , \Delta_g^{-1} \right]  \setminus \bigcup_{\nu=1}^{\nu_{\max}} I^{\bad}_\nu.
\]
Thus, since the union of the $I_\nu^\bad$ contains $\cB \cup \cC$, we have $I_\mu \subset \left[\Delta_g,\Delta_g^{-1} \right] \setminus ( \cB \cup \cC)$ for each $\mu$. By \eqref{pickmon1_fin}, there exist $k =k(\mu) \in \{1,\cdots, K  \}$ and $\ell = \ell(\mu) \in \{1,\cdots, L \}$ such that
\begin{align*}
&\lvert \widetilde{a}_k \delta^{\lambda_k} \rvert > 2 \sum_{k' \neq k} \lvert \widetilde{a}_{k'} \delta^{\lambda_{k'}} \rvert \qquad  \mbox{and} \;\; \\
&  \lvert \widetilde{b}_{\ell} \delta^{\gamma_{\ell}} \rvert > 2 \sum_{\ell' \neq \ell} \lvert \widetilde{b}_{\ell'} \delta^{\gamma_{\ell'}} \rvert  \qquad \mbox{for all } \delta \in I_\mu.
\end{align*}
We compute $k=k(\mu)$ and $\ell = \ell(\mu)$, for each $\mu$, by searching over all $k$, $\ell$ to determine the maximal value of $\lv \til{a}_k \delta_*^{\lambda_\ell} \rv$ and $\lv\til{b}_\ell \delta_*^{\gamma_\ell} \rv$ for any fixed $\delta_* \in I_\mu$. Then, by definition of $\widetilde{\eta}_{\min}(\delta)$ we see that
\begin{equation}
\label{fff2}
c \cdot \widetilde{\eta}_{\min}(\delta) \leq \frac{ \widetilde{a}_{k} \delta^{\lambda_{k}} }{ \widetilde{b}_{\ell} \delta^{\gamma_{\ell}}  }  \leq C \cdot \widetilde{\eta}_{\min}(\delta) \qquad \mbox{for all}\;  \delta \in I_\mu.
\end{equation}

According to \eqref{uws2a}, we also have $c \cdot \widetilde{\eta}_{\min}(\delta) \geq c \cdot \Delta_g^{C}$. 

We compute machine numbers $d_\mu$ such that $ \lv d_\mu - \til{a}_k/ \til{b}_\ell\rv \leq \Delta_g^{-C} \Delta_\epsilon$, and numbers $\omega_\mu = \lambda_k - \gamma_\ell$, where $k = k(\mu)$ and $\ell = \ell(\mu)$. We claim that
\[c \cdot \eta_{\min}(\delta) \leq  d_\mu \cdot \delta^{\omega_\mu} \leq C \cdot \eta_{\min}(\delta) \qquad \mbox{for all} \; \delta \in I_\mu.\]
Indeed, $d_\mu \cdot \delta^{\omega_\mu}$ differs from $\frac{ \widetilde{a}_{k} \delta^{\lambda_{k}} }{ \widetilde{b}_{\ell} \delta^{\gamma_{\ell}}  }$ by at most an additive error of $\Delta_g^{-C} \Delta_\epsilon$, since $\Delta_g \leq \delta \leq \Delta_g^{-1}$ and $ \lv \omega_\mu \rv \leq C$. This additive error is bounded by $\frac{1}{2} c \cdot \Delta_g^{C} \leq c \cdot \widetilde{\eta}_{\min}(\delta)$ since, by assumption, $\Delta_\epsilon \leq \frac{1}{2} c \Delta_g^{2C}$. Hence, \eqref{fff2} implies the above claim.

We compute a machine number $\delta_\nu$ in each interval $I_\nu^\bad$. We know that $e^{-3A} \leq \frac{\delta}{\delta_\nu} \leq e^{3A}$ for all $\delta \in I_\nu^\bad$, due to \eqref{fff1}. Hence, \eqref{slowvariance_new}  implies that
\begin{equation}
\label{fff3}
c \cdot  \eta_{\min}(\delta) \leq \eta_{\min}(\delta_\nu) \leq  C \cdot \eta_{\min}(\delta) \qquad \mbox{for all} \; \delta \in I_\nu^\bad.
\end{equation}
We then compute a machine number $\Gamma_\nu$ such that $\lvert \Gamma_\nu  - \widetilde{\eta}_{\min}(\delta_\nu) \rvert \leq \Delta_g^{-C} \Delta_\epsilon$. Thus,  from \eqref{uws2a} and \eqref{fff3} we conclude that
\[  
c' \cdot \eta_{\min}(\delta)  \leq \Gamma_\nu \leq C' \cdot \eta_{\min}(\delta)  \qquad \mbox{for all} \; \delta \in I_\nu^\bad.
\]

We define $\eta_* : \left[ \Delta_g, \Delta_g^{-1} \right] \rightarrow \R$ by
\begin{equation}
\label{etastar}
\eta_*(\delta) = 
\left\{
\begin{array}{c}
d_\mu \cdot \delta^{\omega_\mu}  \;\; \mbox{if} \; \delta \in I_\mu \\
 \Gamma_\nu \;\; \mbox{if} \; \delta \in I_\nu^\bad.
\end{array}
\right.
\end{equation}
As shown above, we have $c \cdot \eta_*(\delta) \leq \eta_{\min}(\delta) \leq C \cdot \eta_*(\delta)$ for $\delta \in \left[ \Delta_g , \Delta_g^{-1} \right]$, hence we obtain the main estimate in the conclusion of the procedure \textsc{Approximate Rational Function} (finite-precision). This completes the explanation.
\end{proof}

As mentioned before, by applying the procedure \textsc{Approximate Rational Function} we compute a function $\eta_*(\delta)$ satisfying the conditions of the algorithm \textsc{Fit Basis to Convex Body} (finite-precision). This completes the explanation.

\end{proof}

\section{Compressing Norms in Finite-Precision}
\label{sec_cn_fin}

We assume that $\Delta_{\min} \leq \Delta_\epsilon \leq \Delta_g \leq \Delta_0 \leq 1$ are as in the \textbf{Main Assumptions} in Section \ref{fp_not}. In particular, $\Delta_{\min} = 2^{- S}$ ($S = K_{\max} \overline{S}$) denotes the machine precision of our computer, and $\Delta_0 = 2^{- \overline{S}}$. We assume that $\Delta_\epsilon \leq \Delta_g^C$ for a large enough universal constant $C$.

Let $\mu$ be a linear functional on $\R^D$ given in the form $\mu(v) = v \cdot w$, where $w \in \R^D$ is given as $w=(w_1,\cdots,w_D)$. We define
\[
\| \mu \| := \max_{1 \leq i \leq D} \lv w_i \rv.
\]
We say that $\mu$ is specified with parameters $(\Delta_g,\Delta_\epsilon)$ if $\| \mu \| \leq \Delta_g^{-1}$ and if $w_j$ is specified to precision $\Delta_\epsilon$ for each $i=1,\cdots,D$. This means that machine numbers $w_i^\fin$ are given with $\lv w_i - w_i^\fin \rv \leq \Delta_\epsilon$ for each $1 \leq i \leq D$.

We assume that the following data are given.
\begin{itemize}
\item We fix a machine number $\Delta \in \bigl[ \Delta_g , 1 \bigr]$ of the form $\Delta = 2^{- K \overline{S}}$ for an integer $K \geq 1$.
\item We specify linear functionals $\overline{\mu}_1,\cdots,\overline{\mu}_{\overline{L}}$ on $\R^D$ with parameters $(\Delta_g,\Delta_\epsilon)$. We assume that $\overline{L} \leq \Delta_g^{-1}$, and that $ D \leq \til{C}$ for a universal constant $\til{C}$.

\item We fix an $\overline{S}$-bit machine number $p>1$.
\end{itemize}
We denote $\displaystyle \lvert v \rvert = \left( \sum_{i=1}^D \lvert v_i \rvert^p \right)^{1/p}$ for $v = (v_1,\cdots,v_D) \in \R^D$.

\environmentA{Algorithm: Compress Norms (Finite-Precision Version).}

Fix $1 < p < \infty$, and fix an integer $D \geq 1$ as above. Let $\overline{\mu}_1, \cdots,\overline{\mu}_{\overline{L}}$ be linear functionals on $\R^D$, and let $\Delta \in \bigl[ \Delta_g, 1 \bigr]$ be as above. 

We compute linear functionals $\mu_1^*,\cdots,\mu_D^*$ on $\R^D$ such that
\begin{equation}
\label{desired_est}
c \cdot \sum_{i=1}^D \lvert \mu_i^*(v) \rvert^p \leq \sum_{\ell=1}^{\overline{L}} \lvert \overline{\mu}_\ell(v) \rvert^p + \Delta^p \lvert v \rvert^p \leq C \cdot \sum_{i=1}^D \lvert \mu_i^* (v) \rvert^p   \text{ for all } v \in \R^D.
\end{equation}
The $\mu_i^*$ are represented as $v \mapsto v  \cdot \overset{*}{w}^i$ where $\overset{*}{w}^i = (\overset{*}{w}^i_1,\cdots,\overset{*}{w}^i_D)$ and the $\overset{*}{w}^i_k$ are computed with parameters $(\Delta_g^C,\Delta_g^{-C} \Delta_\epsilon)$. 

This computation requires work and storage at most $C \overline{L}$. 

Here, $c > 0$ and $C \geq 1$ are universal constants.
\begin{proof}[Explanation:]

We proceed by induction on $D$.

First consider the base case $D = 1$. The given functionals on $\R^1$ have the form $ \overline{\mu}_\ell(v) = w_\ell \cdot v $ ($1 \leq \ell \leq \overline{L}$), where the numbers $w_\ell$ are specified with parameters $(\Delta_g,\Delta_\epsilon)$. We define $\gamma := \bigl(\lv w_1 \rv^p + \cdots + \lv w_{\overline{L}} \rv^p + \Delta^p \bigr)^{1/p}$. Using Lemma \ref{stablelp}, we compute a machine number $\widehat{\gamma}$ such that $\frac{1}{2} \gamma \leq \widehat{\gamma} \leq 2 \gamma$. Define the functional $\mu^*_1(v) = \widehat{\gamma} \cdot v$ on $\R^1$. Then the estimate \eqref{desired_est} holds with $c = \frac{1}{2}$ and $C = 2$.

We now treat the induction step. Fix an integer $D \geq 2$. We assume by induction that the algorithm \textsc{Compress Norms} has been established when $D$ is replaced by $D - 1$.

We write $c$, $c'$, $C$, $C'$, etc., to denote constants depending only on $p$ and $D$.

We define the functionals
\begin{equation}
\label{tildemu}
\omega _i(v) := \Delta \cdot v_i \qquad \mbox{for} \; v = (v_1,\cdots,v_D) \in \R^D, \quad \mbox{for each} \; i=1,\cdots,D.
\end{equation}

Let  $\{\mu_1,\cdots,\mu_{L}\}$ denote the collection $\{ \overline{\mu}_1,\cdots \overline{\mu}_{\overline{L}}, \omega_1,\cdots,\omega_D\}$ of linear functionals on $\R^D$. Except for minor modifications, we mimic the computation in the infinite-precision version of \textsc{Compress Norms} (see Section \ref{sec_lf}), using the collection $\{\mu_1,\cdots,\mu_{L}\}$ as input. We include the extra functionals $\omega_i$ in order to ensure that we never encounter division by a small number. This leads to the required numerical stability. We provide details of the computation below.

For each $1 \leq i \leq L$, we write
\begin{align}
\label{fff4}
\mu_i(v_1,\cdots,v_D) &:=  \beta_i^* \cdot v_D + \mu_{i}(v_1,\cdots,v_{D-1},0) \\
\notag{}
& = \epsilon_i \cdot \bigl[ \beta_i v_D - \widetilde{\mu}_i(v_1,\cdots,v_{D-1})  \bigr],
\end{align}
where $\beta_i = \lv \beta^*_i \rv$, $\epsilon_i = \sgn(\beta^*_i)$, and $\widetilde{\mu}_i(v_1,\cdots,v_{D-1}) = - \epsilon_i \cdot \mu_i(v_1,\cdots,v_{D-1},0)$. Here, $\sgn(\cdot)$ denotes the ``signum'' function: $\sgn(\alpha) = 1$ if $\alpha \geq 0$, and $\sgn(\alpha) = -1$ if $\alpha < 0$.

The numbers $\beta^*_i$ in \eqref{fff4} are given with parameters $(\Delta_g,\Delta_\epsilon)$, since the functionals $\overline{\mu}_\ell$ are given with parameters $(\Delta_g,\Delta_\epsilon)$ and the $\omega_i$ are given exactly. Hence, we can compute $\beta_i$ with parameters $(\Delta_g, 10\Delta_\epsilon)$ for each $i$. We cannot compute $\epsilon_i$ or $\widetilde{\mu}_i$ with any accuracy unless $\lv \beta_i^* \rv > \Delta_\epsilon$, but this remark will not cause much difficulty.

We set $\Delta_1 = \Delta_g^{C_0}$, for a universal constant $C_0 \in \N$ that will be determined later. Recall that $\beta_i$ is specified to precision $\Delta_\epsilon$, and that $\Delta_\epsilon \leq \frac{1}{4} \Delta_g^{C_0} = \frac{1}{4} \Delta_1$. Hence, we can compute a subset $I^{\fin} \subset \{1,\cdots,L\}$ such that
\begin{equation}\label{mainp}
\beta_i  \leq 2 \Delta_1 \; \mbox{for} \; i \notin I^{\fin}, \;\; \mbox{and} \;\; \beta_i \geq \Delta_1\; \mbox{for} \; i \in I^{\fin}.
\end{equation}
(Just compare the machine approximation of each $\beta_i$ to $ \frac{3}{2} \Delta_1$.)

We compute $\epsilon_i = \sgn(\beta^*_i)$ exactly if $i \in I^{\fin}$, since then we have $ \lvert \beta^*_i \rvert = \beta_i \geq \Delta_1 \geq 2 \Delta_\epsilon$. (We do not attempt to compute  $\epsilon_i$ for $i \notin I^{\fin}$.) Hence, we can compute the functional $\widetilde{\mu}_i$ with parameters $(\Delta_g,\Delta_\epsilon)$ for each $i \in  I^{\fin}$.

Alternatively, for each $i \notin I^{\fin}$, we define the functional  $\widetilde{\widetilde{\mu}}_i(v_1,\cdots,v_{D-1}) = \mu_i(v_1,\cdots,v_{D-1},0)$, which is given with parameters $(\Delta_g,\Delta_\epsilon)$. We have either $\widetilde{\widetilde{\mu}}_i = - \widetilde{\mu}_i$ or $\widetilde{\widetilde{\mu}}_i = \widetilde{\mu}_i$, though we do not guarantee which case occurs.

We have $\mu_{i_0} = \omega_D$ for some $i_0 \in \{1,\cdots,L\}$. From \eqref{tildemu}, we see that $\beta_{i_0} = \Delta \geq \Delta_g > 2\Delta_1$, since $\beta_{i_0}$ is the magnitude of the coefficient of $v_D$ in $\mu_{i_0} = \omega_D$. Hence, $i_0 \in I^{\fin}$, thanks to \eqref{mainp}. Therefore,
\begin{equation}\label{nondegen1}
\mathbf{B} := \sum_{i\in I^{\fin}} \lvert \beta_i \rvert^p \geq \lv \beta_{i_0} \rv^p  =  \Delta^p
\end{equation}
Each $\beta_i$ in \eqref{nondegen1} is given with parameters $(\Delta_g,\Delta_\epsilon)$. Hence, we can compute $\mathbf{B}$ to precision $L \cdot \Delta_g^{-C} \Delta_\epsilon \leq \Delta_g^{-C'} \Delta_\epsilon$, since $\#(I^{\fin}) \leq L \leq \Delta_g^{-C}$ (note: the error invoked in computing each exponentiation $\lv \beta_i \rv^p$ is bounded by $\Delta_g^{-C} \Delta_\epsilon$). Clearly, also $\mathbf{B} \in \left[ \Delta_g^C, \Delta_g^{-C} \right]$. Hence, for each  $i \in I^{\fin}$, we can compute $\mbox{Prob}(i) := \lv\beta_i\rv^p/\mathbf{B}$ with parameters $(1,\Delta_g^{-C}\Delta_\epsilon)$.

Recall that the coefficients of $\til{\mu}_i : \R^{D-1} \rightarrow \R$ are bounded by $\Delta_g^{-1}$, and $D \leq C$. Therefore,
\begin{equation}
\label{bound1}
\lv \til{\mu}_i(v_1,\cdots,v_{D-1}) \rv^p \leq \Delta_g^{-C} \lv v \rv^p.
\end{equation}

The list $\{ \mu_1,\cdots,\mu_L \}$  consists of the functionals $\overline{\mu}_\ell$ and $\omega_i$ (defined in \eqref{tildemu}). Hence,
\begin{equation}\label{b_exp}
 \sum_{\ell=1}^{\overline{L}} \lv \overline{\mu}_\ell(v_1,\cdots,v_D) \rv^p + \Delta^p \lv v \rv^p = \sum_{i=1}^L \lv \mu_i(v_1,\cdots,v_D) \rv^p
\end{equation}
which differs by at most a factor of $\overline{C}$ from
\begin{align}\label{c_exp}
&\mathbf{B} \cdot \lv v_D - \overline{\mu}(v_1,\cdots,v_{D-1}) \rv^p \\
\notag{}
& + \biggl\{  \mathbf{B} \cdot \sum_{i \in I^{\fin}} \mbox{Prob}(i) \cdot \lv \overline{\mu}(v_1,\cdots,v_{D-1}) - \beta_i^{-1} \til{\mu}_i(v_1,\cdots,v_{D-1}) \rv^p \\
\notag{}
& \qquad + \left[ \sum_{i \notin I^{\fin}} \lv \beta_i v_D - \til{\mu}_i(v_1,\cdots,v_{D-1} ) \rv^p\right]
\biggr\},
\end{align}
where
\begin{align}
\label{mubar_fin}
\overline{\mu}(v_1,\cdots,v_{D-1}) &:= \sum_{i \in I^{\fin}} \mbox{Prob}(i) \cdot \bigl\{ \beta_i^{-1} \widetilde{\mu}_i(v_1,\cdots,v_{D-1}) \bigr\} \\
&= \mathbf{B}^{-1} \cdot \sum_{i \in I^{\fin}} \beta_i^{p-1} \cdot \widetilde{\mu}_i(v_1,\cdots,v_{D-1}). \notag{}
\notag{}
\end{align}
These estimates are proven using the same method as in the derivation following \eqref{mubar}. (In contrast to the prior setting, we no longer guarantee here that $\beta_i = 0$ for $i \notin I^{\fin}$, which explains why the third line in \eqref{c_exp} contains the extra term $\beta_i v_D$.) 

Note that $\cI^\fin \neq \emptyset$, as we saw just before \eqref{nondegen1}. Recall that $\lv \text{Prob}(i) \rv \leq 1$ for each $i$, and that $\beta_i \geq \Delta_1 = \Delta_g^{C_0}$ for each $i \in I^{\fin}$. Therefore, from \eqref{mubar_fin} we see that
\[ 
\| \overline{\mu} \| \leq \#(I^{\fin}) \cdot \Delta_g^{-C} \cdot \max_i \left\{ \| \widetilde{\mu}_i \| \right\} \leq \Delta_g^{-C'}.
\] 
Moreover, we can compute $\overline{\mu}$ in \eqref{mubar_fin} to precision $\Delta_g^{-C} \Delta_\epsilon$. Hence, we can compute $\overline{\mu}$ with parameters $(\Delta_g^{C}, \Delta_g^{-C} \Delta_\epsilon)$.

We next estimate the term inside the brackets in \eqref{c_exp}. Applying the estimate \\
$\lv \lv x + y \rv^p  - \lv x \rv^p \rv \leq p \cdot \lv y \rv \cdot  (\lv x \rv + \lv y \rv)^{p-1}$, we have
\begin{align*}
&\lv
\sum_{i \notin I^{\fin}} \lv \beta_i v_D - \til{\mu}_i(v_1,\cdots,v_{D-1}) \rv^p   - 
\sum_{i \notin I^{\fin}} \lv  \til{\mu}_i(v_1,\cdots,v_{D-1}) \rv^p 
\rv \\
&\leq
p \cdot
\sum_{i \notin I^{\fin}} \lv \beta_i v_D \rv \cdot \bigl\{ \lv \beta_i v_D \rv + \lv \til{\mu}_i(v_1,\cdots,v_{D-1}) \rv \bigr\}^{p-1} \\
& \leq C L \Delta_g^{-C} \Delta_1 \lv v \rv^p \leq \Delta_g^{C_0-C'} \lv v \rv^p.
\end{align*}
The constant $C'$ is independent of $C_0$. Here, we use estimate \eqref{bound1}, that $\lv \beta_i \rv \leq 2 \Delta_1$ for $i \notin I^{\fin}$ (see \eqref{mainp}), and that the number of relevant $i$ is  bounded by $ L \leq \Delta_g^{-C}$. Hence, 
\begin{align*}
\mathfrak{S} - \Delta_g^{C_0-C'} \lv v \rv^p & \leq \biggl[ \text{bracketed expression in \eqref{c_exp}} \biggr] \leq \mathfrak{S}  + \Delta_g^{C_0-C'} \lv v \rv^p, \\
& \mbox{where} \; \mathfrak{S} :=  \sum_{i \notin I^{\fin}} \lv  \widetilde{\mu}_i(v_1,\cdots,v_{D-1}) \rv^p.
\end{align*}
We now fix the constant $C_0$ used to define $\Delta_1 = \Delta_g^{C_0}$. We take $C_0$ much larger than $C'$ above, so that the junk term $\Delta_g^{C_0-C'} \lv v \rv^p$ is bounded by $\frac{1}{10} (\overline{C})^{-1} \Delta_g^p \lv v \rv^p \leq \frac{1}{10} (\overline{C})^{-1} \Delta^p \lv v \rv^p$. Hence, we can replace the expression inside square brackets in \eqref{c_exp} with $\mathfrak{S}$, and we can absorb the junk term $\Delta_g^{C_0-C'} \lv v \rv^p$ into the junk term $\Delta^p \lv v \rv^p$ in \eqref{b_exp}. Consequently, \\
$\displaystyle  \sum_{\ell=1}^{\overline{L}} \lv \overline{\mu}_\ell(v_1,\cdots,v_D) \rv^p + \Delta^p \lv v \rv^p $ differs by at most a factor of $C''$ from
\begin{align*}
\notag{}
& \mathbf{B}\cdot \lv v_D - \overline{\mu}(v_1,\cdots,v_{D-1}) \rv^p \\
\notag{}
& + \biggl\{ \mathbf{B} \cdot \sum_{i \in I^{\fin}} \mbox{Prob}(i) \cdot \lv \overline{\mu}(v_1,\cdots,v_{D-1}) - \beta_i^{-1} \til{\mu}_i(v_1,\cdots,v_{D-1}) \rv^p \\
\notag{}
& \qquad + \left[ \sum_{i \notin I^{\fin}} \lv \til{\mu}_i(v_1,\cdots,v_{D-1}) \rv^p\right]
\biggr\}
\end{align*}
We add $\Delta^p  \lv (v_1,\cdots,v_{D-1}) \rv^p$ to both expressions in the previous sentence. Note that \\
$\Delta^p \lv (v_1,\cdots,v_{D-1}) \rv^p + \Delta^p \lv v \rv^p$ differs by at most a factor of $2$ from $\Delta^p \lv v \rv^p$. Therefore,  \\
$\displaystyle  \sum_{\ell=1}^{\overline{L}} \lv \overline{\mu}_\ell(v_1,\cdots,v_D) \rv^p + \Delta^p \lv v \rv^p$ differs by at most a factor of $C'''$ from
\begin{align}
\label{curly}
& \mathbf{B} \cdot \lv v_D - \overline{\mu}(v_1,\cdots,v_{D-1}) \rv^p \\
\notag{}
& \qquad+ \biggl\{  \sum_{i \in I^{\fin}}  \lv \beta_i \overline{\mu}(v_1,\cdots,v_{D-1}) - \til{\mu}_i(v_1,\cdots,v_{D-1}) \rv^p \\
\notag{}
& \qquad\qquad + \sum_{i \notin I^{\fin}} \lv \til{\til{\mu}}_i(v_1,\cdots,v_{D-1} )\rv^p +\Delta^p \lv (v_1,\cdots,v_{D-1}) \rv^p
\biggr\}
\end{align}
(Recall that $\mbox{Prob}(i) = \lv \beta_i \rv^p/ \mathbf{B}$ and that $\widetilde{\mu}_i = \pm  \widetilde{\widetilde{\mu}}_i$.) We consider the functionals arising inside the curly brackets above, namely
\[
\widehat{\mu}_i(v_1,\cdots,v_{D-1}) := 
\left\{
\begin{array}{ccc}
& \beta_i \overline{\mu}(v_1,\cdots,v_{D-1}) - \til{\mu}_i(v_1,\cdots,v_{D-1}) \;\; &\mbox{if} \; i \in I^{\fin}. \\
& \til{\til{\mu}}_i(v_1,\cdots,v_{D-1}) \; &\mbox{if} \; i \notin I^{\fin}
\end{array}
\right.
\]
Note that $\| \widehat{\mu}_i \| \leq \Delta_g^{-C}$, since the same upper bound holds for $\overline{\mu}$, $\til{\til{\mu}}_i = \pm \til{\mu}_i$, and $\beta_i$. Moreover, each $\widehat{\mu}_i$ can be computed to precision $\Delta_g^{-C} \Delta_\epsilon$. Hence, we can compute $\widehat{\mu}_i$ ($1 \leq i \leq L$) with parameters $(\Delta_g^C, \Delta_g^{-C} \Delta_\epsilon)$.

The functionals $\widehat{\mu}_i$ are given with parameters $(\Delta_g^{C}, \Delta_g^{-C} \Delta_\epsilon)$ and $\Delta_g^C \leq \Delta_g \leq \Delta \leq 1$. Hence, by the induction hypothesis, we can compute functionals $\mu_1^*,\cdots,\mu_{D-1}^* : \R^{D-1} \rightarrow \R$ such that
\[
\sum_{i=1}^{D-1} \lv \mu_i^* (v_1,\cdots,v_{D-1}) \rv^p
\]
differs by at most a factor of $C$ from the expression in curly brackets in \eqref{curly}. The $\mu_1^*,\cdots,\mu_{D-1}^*$ are specified with parameters $(\Delta_g^{C'},\Delta_g^{-C'} \Delta_\epsilon)$.

We define
\[
\mu_D^*(v_1,\cdots,v_D) := \mathbf{B}^{1/p} \cdot \left[ v_D - \overline{\mu}(v_1,\cdots,v_{D-1}) \right].
\]
We can compute $\mu_D^*$ with parameters $(\Delta_g^{C},\Delta_g^{-C} \Delta_\epsilon)$, since the same is true of $\mathbf{B}$ and $\overline{\mu}$, and since $\mathbf{B} \geq \Delta^p \geq \Delta_g^p$ (see \eqref{nondegen1}).

Thus, from \eqref{curly}, we see that 
\begin{align*}
& c \cdot \left[  \sum_{\ell=1}^{\overline{L}} \lv \overline{\mu}_\ell(v_1,\cdots,v_D) \rv^p + \Delta^p \lv v \rv^p\right] \\
& \hspace{2cm}  \leq \lv \mu_D^*(v_1,\cdots,v_D) \rv^p + \sum_{i=1}^{D-1} \lv \mu_i^*(v_1,\cdots,v_{D-1}) \rv^p \\
& \hspace{4cm}   \leq C \cdot \left[ \sum_{\ell=1}^{\overline{L}} \lv \overline{\mu}_\ell(v_1,\cdots,v_D) \rv^p + \Delta^p \lv v \rv^p\right] .
\end{align*}

This completes the explanation of the finite-precision version of \textsc{Compress Norms}.

\end{proof}

\section{Algorithm: Optimize via Matrix}
\label{sec_ovm_fin}

We define $\Delta_{\min} \leq \Delta_\epsilon \leq \Delta_g \leq \Delta_0 \leq 1$ as in the \textbf{Main Assumptions} in Section \ref{fp_not}. In particular, $\Delta_{\min} = 2^{- S}$ ($S = K_{\max} \overline{S}$) denotes the machine precision of our computer, and $\Delta_0 = 2^{- \overline{S}}$. We assume that $\Delta_\epsilon \leq \Delta_g^C$ for a large enough universal constant $C$.

We are given the following data:

\begin{itemize}
\item We fix a machine number $\Delta \in \bigl[ \Delta_g , 1 \bigr]$ of the form $\Delta = 2^{- K \overline{S}}$ for an integer $K \geq 1$.
\item We are given a matrix $A = (a_{\ell j})_{\substack{1 \leq \ell \leq L \\ 1 \leq j \leq J}}$. The numbers $a_{\ell j}$ are specified with parameters $(\Delta_g,\Delta_\epsilon)$. We have $1 \leq L \leq \Delta_g^{-1}$ and $1 \leq J \leq C$ for a universal constant $C$.
\item We fix an $\overline{S}$-bit machine number $p>1$.

\end{itemize}

\environmentA{Algorithm: Optimize via Matrix (Finite-Precision).}

Given $1 < p < \infty$, given $\Delta$, and given a matrix $A = (a_{\ell j})_{\substack{1 \leq \ell \leq L \\ 1 \leq j \leq J}}$ as above, we compute a matrix $B = (b_{j\ell})_{\substack{ 1 \leq j \leq J \\ 1 \leq \ell \leq L}}$.  We guarantee that the following conditions hold.

Let $y_1,\cdots, y_{L}$ be real numbers, and set $x_j^* = \sum_{\ell = 1}^{L} b_{j \ell} y_\ell$ for each $j=1,\cdots,J$.

Then 
\[ 
\sum_{\ell=1}^{L} \lvert y_\ell + \sum_{j=1}^J a_{\ell j} x_j^* \rvert^p  \leq C_1 \cdot \left[ \sum_{\ell=1}^{L} \lvert y_\ell + \sum_{j=1}^J a_{\ell j} x_j \rvert^p + \Delta^p \sum_{j=1}^J \lvert x_j \rvert^p \right]\]
for any real numbers $x_1,\cdots,x_J$.

The numbers $b_{j \ell}$ are computed with parameters $(\Delta_g^{C_1}, \Delta_g^{-C_1} \Delta_\epsilon)$. 

The algorithm requires work and storage at most $C_1 \cdot L$.

Here, $C_1$ is a universal constant.

\begin{proof}[\underline{Explanation}]
We write $c, C, C'$, etc., to denote universal constants. 

We proceed by induction on $J$. We first handle the case $J=1$.

Assume that an $L \times 1$ matrix $(a_{\ell})_{1 \leq \ell \leq L}$ is given, with each number $a_\ell$ specified with parameters $(\Delta_g,\Delta_\epsilon)$.

Let $y_1,\cdots,y_L$ be given real numbers. 

We define $y_{0} = 0$ and $a_{0} = \Delta$.

We compute an index set $\cL \subset \{ 0,\cdots, L\}$ such that $\lv a_\ell \rv \geq \Delta_g^{10}$ for $\ell \in \cL$, and $\lv a_\ell \rv \leq 2 \Delta_g^{10}$ for $\ell \in \{0,\cdots, L \} \setminus \cL$. To do so, we compare the machine approximation of each $\lv a_\ell \rv$ to the machine number $\frac{3}{2} \Delta_g^{10}$.

Note that $a_0 = \Delta \geq \Delta_g > 2 \Delta_g^{10}$, which implies that $0 \in \cL$. In particular, $\cL \neq \emptyset$.

If $\lv a_\ell \rv \leq 2 \Delta_g^{10}$ then the quantities $\lv y_\ell \rv + 2\Delta_g^{10} \cdot \lv x \rv $ and $\lv y_\ell + a_\ell x \rv + 2\Delta_g^{10} \cdot \lv x \rv$ differ by at most a factor of $2$, thanks to the triangle inequality. Thus,
\begin{align*}
\lv y_\ell + a_\ell x \rv^p + \Delta_g^{10p} \cdot \lv x \rv^p &\sim  \lv y_\ell \rv^p + \Delta_g^{10p} \cdot \lv x \rv^p \\
& \mbox{for} \; \ell \in \{0,\cdots,L\} \setminus \cL,
\end{align*}
where $A \sim B$ indicates that $c \cdot A \leq B \leq C \cdot A$ for some universal constants $c > 0$ and $C \geq 1$. Therefore, we have
\begin{align}
\label{bcd1}
& \sum_{\ell=0}^L \lv y_\ell + a_\ell x \rv^p +  \mathcal{E}(x) \\
& \qquad\qquad \sim \sum_{\ell \in \cL} \lv y_\ell + a_\ell x \rv^p + \sum_{\ell \in \{0,\cdots,L\} \setminus \cL} \lv y_\ell \rv^p +  \mathcal{E}(x), \notag{} \\
& \mbox{where} \; \mathcal{E}(x) = \# ( \{ 0,\cdots,L\} \setminus \cL ) \cdot \Delta_g^{10p} \cdot \lv x \rv^p. \notag{}
\end{align}
Since $L \leq \Delta_g^{-1}$, it follows that $\mathcal{E}(x) \leq \Delta_g^p \cdot \lv x \rv^p \leq \Delta^p \cdot \lv x \rv^p = \lv y_{0} + a_{0} x \rv^p$. Therefore, because $\lv y_{0} + a_{0} x \rv^p$ is a summand on both sides of \eqref{bcd1}, we can discard the error term $\mathcal{E}(x)$ and obtain
\begin{equation*}
\sum_{\ell=0}^L \lv y_\ell + a_\ell x \rv^p  \sim \sum_{\ell \in \cL} \lv y_\ell + a_\ell x \rv^p + \sum_{ \ell \in \{0,\cdots,L\} \setminus \cL} \lv y_\ell \rv^p.
\end{equation*}
We write this estimate in the form
\begin{align}
\label{fff7_new}
\sum_{\ell=0}^L \lv y_\ell + a_\ell x \rv^p  \sim \sum_{\ell \in \cL} \lv \overline{y}_\ell + x \rv^p \cdot \lv a_\ell  \rv^p &+ \sum_{ \ell \in \{0,\cdots,L\} \setminus \cL} \lv y_\ell \rv^p,  \\
&\mbox{where we define} \; \overline{y}_\ell := \frac{y_\ell}{a_\ell}. \notag{}
\end{align}

Now, we want to minimize the expression $\mathcal{T}(x) = \sum_{\ell \in \cL} \lv \overline{y}_\ell + x \rv^p \cdot \lv a_\ell \rv^p$ up to a universal constant factor. We define
\[
\left\{
\begin{aligned}
x^* &:= - \sum_{\ell \in \cL} \overline{y}_\ell \cdot  \mbox{Prob}(\ell), \;\; \mbox{where} \\
\mbox{Prob}(\ell) &:= \left( \sum_{\ell' \in \cL} \lv a_{\ell'} \rv^p \right)^{-1} \cdot \lv a_{\ell} \rv^p  \qquad \mbox{for} \;\; \ell \in \cL.\end{aligned}
\right.
\]
Recall that $\cL \neq \emptyset$, hence $\mbox{Prob}(\ell)$ is a well-defined probability measure on $\cL$. By applying \eqref{cn1} we conclude that $\mathcal{T}(x_*) \leq C \cdot \mathcal{T}(x)$ for all $x \in \R$. Therefore, we have
\[
\sum_{\ell=0}^L \lv y_\ell + a_\ell x \rv^p \leq C' \cdot \sum_{\ell=0}^L \lv y_\ell + a_\ell x^* \rv^p.
\]
Because $y_0 = 0$ and $a_0 = \Delta$, this implies that
\[
\sum_{\ell=1}^L \lv y_\ell + a_\ell x \rv^p \leq C' \cdot \left[ \sum_{\ell=1}^L \lv y_\ell + a_\ell x^* \rv^p + \Delta^p \lv x \rv^p \right],
\]
as desired in the case $J=1$ of our algorithm. Note that
\[
\left\{
\begin{aligned}
x^* &= - \sum_{\ell \in \cL}  \overline{y}_\ell \cdot \mbox{Prob}(\ell) = \sum_{\ell \in \cL \setminus \{0\}} y_\ell \cdot b_\ell, \;\;\mbox{where} \\
b_\ell &= - \left( \sum_{\ell' \in \cL} \lv a_{\ell'} \rv^p \right)^{-1} \cdot \lv a_\ell \rv^p \cdot a_\ell^{-1} \;\; \mbox{for} \; \ell \in \cL \setminus \{0\}.
\end{aligned}
\right.
\]
It is safe to discard the $\ell=0$ term in the sum, because $y_0 = \overline{y}_0 = 0$ by definition. Note that $\lv a_\ell \rv$ and $\lv a_{\ell'} \rv$, for $\ell, \ell' \in \cL$, belong to the interval $[ \Delta_g^{10},\Delta_g^{-1}]$. Therefore, we can compute $\lv a_{\ell'} \rv^p$ and $\lv a_\ell \rv^p$ to precision $\Delta_g^{-C} \Delta_\epsilon$; moreover, we can compute the expression $\left( \cdots \right)^{-1}$ - in the formula for $b_\ell$ - with precision $\Delta_g^{-C} \Delta_\epsilon$. Thus, we can compute the coefficients $b_\ell$, for each $\ell \in \cL \setminus \{0\}$, with precision $\Delta_g^{-C} \Delta_\epsilon$. Furthermore, note  that each $\lv b_\ell \rv$ is bounded by $\Delta_g^{-C}$ for a universal constant $C \geq 1$. 

All the remaining coefficients $b_\ell$, for $\ell \in \{1,\cdots,L\} \setminus \cL$, are defined to be $0$. Thus, $x^* = \sum_{\ell=1}^L y_\ell \cdot b_\ell$, and $b_\ell$ can be computed with the desired parameters. Thus, we have established the case $J=1$ of our algorithm.

\underline{For the general case}, we use induction on $J$.

Let $J \geq 2$, and let $1 < p < \infty$ and assume that we are given an $L \times J$ matrix $A = (a_{\ell j})_{\substack{1 \leq \ell \leq L \\ 1 \leq j \leq J}}$. We assume that the numbers $a_{\ell j}$ are specified with parameters $(\Delta_g,\Delta_\epsilon)$.

Let real numbers $y_1,\cdots,y_L$ be given.

We have
\begin{equation}
\label{ww2}
 \sum_{\ell=1}^L \lvert y_\ell + \sum_{j=1}^J a_{\ell j} x_j \rvert^p  = \sum_{\ell=1}^{L} \lvert \widehat{y}_\ell + \sum_{j=1}^{J-1} a_{\ell j} x_j \rvert^p \qquad ((x_1,\cdots,x_J) \in \R^J),
\end{equation}
using new variables
\begin{equation}\label{ww1}
\widehat{y}_\ell = y_\ell + a_{\ell J} \cdot x_J \qquad \mbox{for} \; 1 \leq \ell \leq L.
\end{equation}

By applying the algorithm \textsc{Optimize via Matrix} recursively to $1 < p < \infty$ and the submatrix $(a_{\ell j})_{\substack{1 \leq \ell \leq L \\ 1 \leq j \leq J-1}}$, we compute a matrix $(\widehat{b}_{j \ell})_{\substack{1 \leq j \leq J-1 \\ 1 \leq \ell \leq L}}$ such that the following holds.

\begin{itemize}
\item We compute the numbers $\widehat{b}_{j \ell}$ with parameters $( \Delta_g^C, \Delta_g^{-C}\Delta_\epsilon)$ for a universal constant $C$.
\item Let $\widehat{y}_1,\cdots,\widehat{y}_L$ be given, and set
\begin{equation}\label{ww3}
\widehat{x}_j = \sum_{\ell=1}^L \widehat{b}_{j \ell} \widehat{y}_\ell \qquad \mbox{for} \; 1 \leq j \leq J-1.
\end{equation}
Then, for any real numbers $x_1,\cdots,x_{J-1}$, we have
\begin{equation}\label{ww4}
\sum_{\ell=1}^{L} \lvert \widehat{y}_\ell + \sum_{j=1}^{J-1} a_{\ell j} \widehat{x}_j \rvert^p  \leq C\cdot \left[ \sum_{\ell=1}^{L} \lvert \widehat{y}_\ell + \sum_{j=1}^{J-1} a_{\ell j} x_j \rvert^p  + \Delta^p \sum_{j=1}^{J-1} \lvert x_j \rvert^p \right].
\end{equation}
\end{itemize}

Using \eqref{ww2}-\eqref{ww4}, we draw the following conclusion.

Let real numbers $y_1,\cdots,y_L$ be given, and let $x_1,\cdots,x_J$ be arbitrary. We define $\widehat{y}_1,\cdots, \widehat{y}_{L}$ by \eqref{ww1}, next define $\widehat{x}_1,\cdots,\widehat{x}_{J-1}$ by \eqref{ww3}, and finally set 
\begin{equation}
\label{ww5}
\widehat{x}_J = x_J.
\end{equation}
Then 
\begin{equation*}
\sum_{\ell=1}^{L} \lvert y_\ell + \sum_{j=1}^J a_{\ell j} \widehat{x}_j \rvert^p  \leq C \cdot \left[   \sum_{\ell=1}^{L} \lvert y_\ell + \sum_{j=1}^J a_{\ell j} x_j \rvert^p + \Delta^p \sum_{j=1}^{J-1} \lvert x_j \rvert^p \right],
\end{equation*}
hence
\begin{equation}\label{ww6}
\sum_{\ell=1}^{L} \lvert y_\ell + \sum_{j=1}^J a_{\ell j} \widehat{x}_j \rvert^p + \Delta^p \cdot \lvert \widehat{x}_J \rvert^p  \leq C \cdot \left[   \sum_{\ell=1}^{L} \lvert y_\ell + \sum_{j=1}^J a_{\ell j} x_j \rvert^p + \Delta^p \sum_{j=1}^{J} \lvert x_j \rvert^p \right],
\end{equation}
and moreover
\begin{equation}\label{ww7}
\widehat{x}_j = \sum_{\ell=1}^{L} \widehat{b}_{j \ell} \cdot (y_\ell + a_{\ell J} \widehat{x}_J ) \quad \mbox{for} \; j=1,\cdots,J-1.
\end{equation}
Thus,
\begin{align}\label{ww8}
&\widehat{x}_j = \sum_{\ell=1}^{L} \widehat{b}_{j \ell} y_\ell + g_j \widehat{x}_J, \quad \mbox{where} \\
\label{ww9}
&g_j := \sum_{\ell=1}^{L} \widehat{b}_{j \ell} a_{\ell J} \qquad \mbox{for} \; j=1,\cdots, J-1.
\end{align}
We compute the numbers $g_j$ with parameters $( \Delta_g^{C'}, \Delta_g^{-C'} \Delta_\epsilon)$ using work at most $C L$. This is possible because $L \leq \Delta_g^{-1}$ and because of parameters with which $\widehat{b}_{j \ell}$ and $a_{\ell j}$ are specified. In the above discussion, the numbers $x_1,\cdots,x_J$ are arbitrary, the numbers $\widehat{x}_1,\cdots,\widehat{x}_{J-1}$ are defined from $\widehat{x}_J$ by \eqref{ww7}, and $\widehat{x}_J = x_J$.

Next, note that
\begin{align*}
y_\ell + \sum_{j=1}^{J} a_{\ell j} \widehat{x}_j & = y_\ell + \sum_{j=1}^{J-1} a_{\ell j} \left[ \sum_{\ell' = 1}^L \widehat{b}_{j \ell'} y_{\ell'} +  g_j \widehat{x}_J  \right]  + a_{\ell J} \widehat{x}_J \\
&= \left\{ y_\ell + \sum_{j=1}^{J-1} a_{\ell j} \sum_{\ell'=1}^{L} \widehat{b}_{j \ell'} y_{\ell'} \right\} + \left\{ a_{\ell J} + \sum_{j=1}^{J-1} a_{\ell j} g_j \right\} \widehat{x}_J \\
& =: \hspace{1cm}  y_\ell^{\ouch} \hspace{2.5cm} + \hspace{1cm}  h_\ell \cdot \widehat{x}_J.
\end{align*}
 Here,
\begin{equation}
\label{ww10}
y_\ell^\ouch = y_\ell + \sum_{j=1}^{J-1} a_{\ell j} \sum_{\ell'=1}^{L} \widehat{b}_{j \ell'} y_{\ell'},
\end{equation}
and
\begin{equation}
\label{ww11}
h_\ell = a_{\ell J} + \sum_{j=1}^{J-1} a_{\ell j} g_j \quad \mbox{for} \; \ell=1,\cdots,L.
\end{equation}
Thus,
\begin{equation} \label{ww12} 
\sum_{\ell=1}^{L} \lvert y_\ell + \sum_{j=1}^J a_{\ell j} \widehat{x}_j \rvert^p = \sum_{\ell=1}^{L} \lvert y^\ouch_\ell + h_\ell \widehat{x}_J \rvert^p.
\end{equation}
Here, \eqref{ww12} holds whenever $\widehat{x}_1,\cdots,\widehat{x}_{J-1}$ are determined from $\widehat{x}_J$ via \eqref{ww8}. 

We compute the numbers $h_\ell$ with parameters $(\Delta_g^C, \Delta_g^{-C} \Delta_\epsilon)$, using work at most $C L$.

Note that it is too expensive to compute $y^\ouch_\ell$ for all $\ell$ ($1 \leq \ell \leq L$); that computation would require $\sim L^2 J$ work. However, the $y^\ouch_\ell$ defined above are independent of our choice of $\widehat{x}_J$.

Applying the known case $J=1$ of our algorithm \textsc{Optimize via Matrix}, we compute from the $h_\ell$ a vector of coefficients $\gamma_\ell$ ($1 \leq \ell \leq L$), for which the following holds.
\begin{itemize}
\item We compute the numbers $\gamma_\ell$ with parameters $( \Delta_g^C,\Delta_g^{-C} \Delta_\epsilon)$ for a universal constant $C$.
\item Let
\begin{equation}
\label{ww13}
\check{x}_J = \sum_{\ell=1}^{L} \gamma_\ell y^{\ouch}_\ell.
\end{equation}
Then
\[
\sum_{\ell=1}^{L} \lvert y_\ell^\ouch + h_\ell \check{x}_J \rvert^p    \leq    C \cdot \left[  \sum_{\ell=1}^{L} \lvert y_\ell^\ouch + h_\ell \widehat{x}_J \rvert^p + \Delta^p \cdot \lvert \widehat{x}_J \rvert^p \right]
\]
for any real number $\widehat{x}_J$.
\end{itemize}

We thus learn the following.

Let $\check{x}_1,\cdots,\check{x}_{J-1}$ be defined from $\check{x}_J$ as in \eqref{ww8}, i.e.,
\begin{equation}\label{ww14}
\check{x}_j := \sum_{\ell=1}^L \widehat{b}_{j \ell} y_\ell + g_j \check{x}_J \quad \mbox{for} \; j=1,\cdots,J-1.
\end{equation}
Let $\widehat{x}_J$ be any real number, and let $\widehat{x}_1,\cdots,\widehat{x}_{J-1}$ be determined from $\widehat{x}_J$ by \eqref{ww8}. Then 
\begin{equation}\label{ww15}
 \sum_{\ell=1}^{L} \lvert y_\ell + \sum_{j=1}^J a_{\ell j} \check{x}_j \rvert^p  \leq C \cdot \left[ \sum_{\ell=1}^{L} \lvert y_\ell + \sum_{j=1}^J a_{\ell j} \widehat{x}_j \rvert^p + \Delta^p \cdot \lvert \widehat{x}_J \rvert^p \right]
 \end{equation}
(See \eqref{ww12}.) 

From \eqref{ww6} and \eqref{ww15}, we see that
\begin{equation}\label{ww16}
\sum_{\ell=1}^L \lvert y_\ell + \sum_{j=1}^J a_{\ell j} \check{x}_j \rvert^p 
\leq C \cdot \left[ \sum_{\ell=1}^{L} \lvert y_\ell + \sum_{j=1}^J a_{\ell j} x_j \rvert^p + \Delta^p \sum_{j=1}^{J} \lvert x_j \rvert^p  \right].
\end{equation}
Here, $\check{x}_1,\cdots,\check{x}_J$ are computed from \eqref{ww13},\eqref{ww14}; and $x_1,\cdots,x_J$ are arbitrary.

We produce efficient formulas for the $\check{x}_j$. Putting \eqref{ww10} into \eqref{ww13}, we find that
\begin{align*}
\check{x}_J & = \sum_{\ell=1}^L \gamma_\ell \cdot \left\{ y_\ell + \sum_{j=1}^{J-1} a_{\ell j} \sum_{\ell' = 1}^L \widehat{b}_{j \ell'} y_{\ell'} \right\} \\
& = \sum_{\ell=1}^L \gamma_\ell \cdot y_\ell + \sum_{\ell' = 1}^L \sum_{j=1}^{J-1}  \left[ \sum_{\ell=1}^L   \gamma_\ell a_{\ell j}  \right] \widehat{b}_{j \ell'} y_{\ell'} \\
& = \sum_{\ell=1}^L \left\{ \gamma_\ell + \sum_{j=1}^{J-1}  \left[ \sum_{\ell'=1}^L  \gamma_{\ell'} a_{\ell' j}  \right] \widehat{b}_{j \ell}  \right\} \cdot y_\ell.
\end{align*}
Therefore, setting
\begin{equation}\label{ww17}
\Delta_j = \sum_{\ell=1}^L \gamma_\ell a_{\ell j} \quad \mbox{for} \; j=1,\cdots, J -1
\end{equation}
and
\begin{equation}\label{ww18}
b^{\#\#}_{J\ell} = \gamma_\ell + \sum_{j=1}^{J-1} \Delta_j \widehat{b}_{j \ell} \quad \mbox{for} \; \ell=1,\cdots,L
\end{equation}
we find that
\begin{equation}\label{ww19}
\check{x}_J = \sum_{\ell=1}^L b^{\# \#}_{J \ell} y_\ell.
\end{equation}
Substituting \eqref{ww19} into \eqref{ww14}, we find that
\[\check{x}_j = \sum_{\ell=1}^L \left\{ \widehat{b}_{j \ell} + g_j b_{J \ell}^{\# \#} \right\} y_\ell \quad \mbox{for} \; j=1,\cdots,J-1.\]
Thus, setting
\begin{align}\label{ww20}
b^{\# \#}_{j \ell} = \widehat{b}_{j \ell} + g_j b^{\# \#}_{J \ell} \qquad \mbox{for} \;\; & j=1,\cdots,J-1,\\
& \qquad  \ell = 1,\cdots, L \notag{}
\end{align}
we have
\begin{equation}\label{ww21}
\check{x}_j = \sum_{\ell=1}^L b^{\# \#}_{j \ell} y_\ell \quad \mbox{for} \; j=1,\cdots,J-1.
\end{equation}
Recalling \eqref{ww19}, we see that \eqref{ww21} holds for $j=1,\cdots,J$. Thus, with $\check{x}_1,\cdots,\check{x}_{J}$ defined by \eqref{ww21}, we have
\[
\sum_{\ell=1}^L \lvert y_\ell + \sum_{j=1}^J a_{\ell j} \check{x}_j \rvert^p \leq C \cdot \left[ \sum_{\ell=1}^{L} \lvert y_\ell + \sum_{j=1}^J a_{\ell j} x_j \rvert^p + \Delta^p \sum_{j=1}^{J} \lvert x_j \rvert^p  \right] 
\]
for any real numbers $x_1,\cdots,x_J$. (See \eqref{ww16}.)

So the matrix $B = (b^{\# \#}_{j \ell})_{\substack{1 \leq j \leq J \\ 1 \leq \ell \leq L}}$ is as promised in our algorithm.

We make a few additional remarks on the computation of $(b^{\# \#}_{j \ell})_{\substack{1 \leq j \leq J \\ 1 \leq \ell \leq L}}$.
\begin{itemize}
\item Recall that the numbers $\gamma_\ell$, $a_{\ell j}$, and $\widehat{b}_{j \ell}$ are given with parameters $( \Delta_g^C, \Delta_g^{-C} \Delta_\epsilon)$. Also recall that $L \leq \Delta_g^{-1}$.
\item Thus, the numbers $\Delta_j$ ($1 \leq j \leq J-1$) in \eqref{ww17} can be computed with parameters $( \Delta_g^C,\Delta_g^{-C} \Delta_\epsilon)$.
\item Consequently, the numbers $b^{\# \#}_{J\ell}$ ($1 \leq \ell \leq L$) in \eqref{ww18} can be computed with parameters $( \Delta_g^C,  \Delta_g^{-C} \Delta_\epsilon)$.
\item Recall that the numbers $g_j$ ($1 \leq j \leq J-1$) in \eqref{ww9} can be computed  with parameters $(\Delta_g^C, \Delta_g^{-C} \Delta_\epsilon)$.
\item Therefore, the numbers $b^{\# \#}_{j \ell}$ ($1 \leq j \leq J-1$, $1 \leq \ell \leq L$) in \eqref{ww20} can be computed with parameters $( \Delta_g^C, \Delta_g^{-C} \Delta_\epsilon)$.
\item Thus, the matrix $B = (b^{\# \#}_{j \ell})$ can be computed to the accuracy promised in the algorithm.
\end{itemize}

\end{proof}

\section{Statement of Main Technical Results}\label{sec_mainresults_fin}

We will prove a modified version of the Main Technical Results for $\cA$ (see Chapter \ref{sec_mainresults}), which accounts for the rounding errors that can arise in the computation.

We define a norm $\displaystyle \lvert P \rvert := \left( \sum_{\alpha \in \cM} \lvert \partial^\alpha P(0) \rvert^p \right)^{1/p}$ for $P \in \cP$. Thus, $\lv P \rv$ denotes the $\ell^p$-norm of the vector $(\partial^\alpha P(0))_{\alpha \in \cM}$. We will always use this norm in the course of the proof.

We fix an integer $\overline{S} \geq 1$.

We are given a finite set $E \subset \frac{1}{32} Q^\circ$, with $Q^\circ = [0,1)^n$. We  assume  that $N = \#(E) \geq 2$. We additionally assume that $E$  consists of  $\overline{S}$-bit machine points. Thus, 
\begin{equation}
\label{farapart}
\lvert x - x' \rvert \geq \Delta_0 \;\; \mbox{for distinct} \; x,x' \in E\text{,}
\end{equation}
where $\Delta_0 := 2^{- \overline{S}}$. Hence, 
\begin{equation}
\label{Nbound}
\#(E) = N \leq \Delta_0^{-n}.
\end{equation}

For $\Delta_1 , \Delta_2 \in (0,1]$, we write $\Delta_1 \ll \Delta_2$ to indicate that $\Delta_1 \leq \Delta_2^C$ for a sufficiently large universal constant $C$.

We introduce constants $\Delta_\epsilon^\circ := 2^{ - K_{1} \overline{S} }$, $\Delta_g^\circ := 2^{ - K_{2} \overline{S} }$, and $\Delta_\junk^\circ := 2^{ - K_{3} \overline{S} }$ as in Theorem \ref{main_thm_hom_fin}. Here, $K_1,K_2,K_3$ are positive integers, which are assumed to be sufficiently well-separated in the sense that $K_1 \geq C \cdot K_2 \geq C^2 \cdot K_3$ for a large enough universal constant $C$.

For each $\cA \subset \cM$, we will use parameters  $\Delta_\epsilon(\cA) = \Delta_0^{K_1(\cA)}$, $\Delta_g(\cA) = \Delta_0^{K_2(\cA)}$, and $\Delta_\junk(\cA) = \Delta_0^{K_3(\cA)}$ for integer exponents $K_1(\cA) \geq K_2(\cA) \geq K_3(\cA) \geq 1$. We assume the exponents are chosen so that 
\begin{align}\label{dominate}
 \Delta_\epsilon(\cM) &\ll  \cdots \ll \Delta_\epsilon(\emptyset)  \ll \Delta_\epsilon^\circ \\
 \notag{}
& \ll \Delta_g^\circ \ll \Delta_g(\emptyset) \ll \cdots \ll \Delta_g(\cM)  \\
\notag{}
& \ll \Delta_\junk(\cM) \ll \cdots \ll \Delta_\junk(\emptyset) \ll \Delta_\junk^\circ \\
\notag{}
& \ll \Delta_0.
\end{align}
In particular,
\begin{equation}
\label{constants0}
\left\{
\begin{aligned}
& \Delta_\epsilon(\emptyset) \leq (\Delta_\epsilon^\circ)^C \\
& \Delta_g^\circ \leq (\Delta_g(\emptyset))^C \\
& \Delta_\junk(\emptyset) \leq (\Delta_\junk^\circ)^C
\end{aligned}
\right.
\end{equation}
and
\begin{equation}
\label{constants1}
\left\{
\begin{aligned}
&\Delta_\epsilon(\cA^-) \leq \Delta_{\epsilon}(\cA^+)^{C} \\
&\Delta_g(\cA^+) \leq \Delta_g(\cA^-)^C \\
& \Delta_\junk(\cA^-) \leq \Delta_\junk(\cA^+)^C \\
&\Delta_\epsilon(\emptyset) \leq \Delta_g(\emptyset)^C \\
& \Delta_g(\cM) \leq \Delta_{\junk}(\cM)^C \\
& \Delta_\junk(\emptyset) \leq \Delta_0^C
\end{aligned}
\right.
\end{equation}
for any $\cA^+ > \cA^-$ and for a large enough universal constant $C$. We refer the reader to Section \ref{sec_multi} for the definition of the order relation $>$ on sets of multiindices. The conditions in \eqref{dominate}, \eqref{constants0}, \eqref{constants1} are clearly consistent with one another. We will use these conditions throughout the course of the proof.


We assume throughout the course of the proof that we can perform arithmetic operations on $S$-bit machine numbers to precision $\Delta_{\min} = 2^{-S}$, where $S = K_{\max} \overline{S}$. Here, the parameter $K_{\max} \in \N$ is larger than all the exponents $K_j(\cA)$ (for all $\cA \subset \cM$ and $j=1,2,3$).

The Main Technical Results for $\cA$ are as in Chapter 3, with the following modifications.

\begin{itemize}
\item We define a dyadic decomposition $\CZ(\cA)$ of $Q^\circ$. We guarantee all the properties \textbf{(CZ1)}$\cdots$\textbf{(CZ5)} in Chapter \ref{sec_mainresults}. We further guarantee that
\begin{equation}
\label{lowed} 
\delta_Q \geq \frac{1}{32} \cdot \Delta_0 \; \mbox{for all} \; Q \in \CZ(\cA).
\end{equation}
Hence, each cube in $\CZ(\cA)$ has $\til{S}$-bit machine points as corners, where $\til{S} \leq \overline{S} + 100$. Thus, we can store each cube in $\CZ(\cA)$ on our computer using at most $C$ units of storage (for a universal constant $C$). However, we will not compute all the cubes in $\CZ(\cA)$ for this would require too much work.
\item We let $\CZ_{\main}(\cA)$ consist of all the cubes $Q \in \CZ(\cA)$ such that $\frac{65}{64}Q \cap E \neq \emptyset$. For each $Q \in \CZ_{\main}(\cA)$, we will compute $\Omega(Q,\cA)$, $\Xi(Q,\cA)$, and $T_{(Q,\cA)}$ as in the three bullet points below.
\item The assists $\omega \in \Omega(Q,\cA)$ are to be given in short form with parameters $(\Delta_g(\cA),\Delta_\epsilon(\cA))$.
\item The functionals $\xi \in \Xi(Q,\cA))$ are to be given in short form with parameters $(\Delta_g(\cA),\Delta_\epsilon(\cA))$ in terms of the assists $\Omega(Q,\cA)$.

We define
\[M_{(Q,\cA)}(f,P) := \left( \sum_{\xi \in \Xi(Q,\cA)} \lvert \xi(f,P) \rvert^p \right)^{1/p}.\]

For each $(f,P) \in \X(\frac{65}{64}Q \cap E) \oplus \cP$, we guarantee that
\[c \cdot \| (f,P) \|_{(1+a(\cA))Q} \leq M_{(Q,\cA)}(f,P) \leq C \cdot \left[ \| (f,P) \|_{\frac{65}{64}Q} + \Delta_\junk(\cA) \cdot \lvert P \rvert \right]. \]

\item The operators $T_{(Q,\cA)}$ map $\X(\frac{65}{64}Q \cap E) \oplus \cP$ into $\X$. 
\begin{description}
\item[(E1)] $T_{(Q,\cA)}(f,P) = f$ on $(1+a(\cA))Q \cap E$ for each $(f,P)$.
\item[(E2)] $\| T_{(Q,\cA)}(f,P) \|^p_{\X((1+a(\cA))Q)} + \delta_Q^{-mp} \| T_{(Q,\cA)}(f,P) - P \|^p_{L^p((1+a(\cA))Q)} \leq C \left[ M_{(Q,\cA)}(f,P) \right]^p$ for each $(f,P)$.
\item[(E3)] $T_{(Q,\cA)}$ has $\Omega(Q,\cA)$-assisted depth at most $C$.
\end{description}
\item The only modification to the algorithm $\CZ$-\textsc{Oracle} is as follows: \\
We assume that the query $\underline{x} \in Q^\circ$ is an $S$-bit machine point. We compute a list of all the cubes $Q \in \CZ(\cA^-)$ such that $\underline{x} \in \frac{65}{64} Q$.

(Recall that $S = K_{\max} \overline{S}$ is the maximum bit length of a machine number representable on our computer.)
\item The algorithm \textsc{Compute Main-Cubes} is unchanged. We compute and store all the cubes in $\CZ_{\main}(\cA)$.
\item The only modifications to the algorithm \textsc{Compute Functionals} are as follows: \\
The functionals $\omega \in \Omega(Q, \cA)$ are computed in short form with parameters $(\Delta_g(\cA), \Delta_\epsilon(\cA))$. The functionals $\xi \in \Xi(Q,\cA)$ are computed in short form with parameters $(\Delta_g(\cA) , \Delta_\epsilon(\cA))$ in terms of the assists $\Omega(Q,\cA)$.

\item The only modifications to the algorithm \textsc{Compute Extension Operators} are as follows:  \\
Let $\underline{x} \in Q^\circ$ be an $S$-bit machine point, and let $\alpha \in \cM$. We compute the linear functional $(f,P) \mapsto \partial^\alpha  (T_{(Q,\cA)}(f,P))(\underline{x})$ in short form with parameters $(\Delta_g(\cA) , \Delta_\epsilon(\cA))$ in terms of the assists $\Omega(Q,\cA)$. This requires work at most $C \log N$, as before.

\item All the constants $c_*(\cA), S(\cA),\epsilon_1(\cA),\epsilon_2(\cA),a(\cA),c,C$ depend only on $m,n,p$, and $\cA$. The constant $S(\cA) \geq 1$ is an integer. We further assume that $a(\cA)$ is an integer power of $2$. (This is a new assumption in the finite-precision case.)
\item We perform the above computations using one-time work at most $C N \log N$ and storage at most $CN$.
\end{itemize}

\section{Algorithms for Dyadic Cubes}

\label{adc_fin}

We make the following assumptions.
\begin{itemize}
\item We are given machine numbers $\Delta_\epsilon = 2^{- K_1 \overline{S}}$ and $\Delta_g = 2^{-K_2 \overline{S}}$, for integers $K_1, K_2 \geq 1$. 
\item We assume that our computer can perform arithmetic operations on $S$-bit machine numbers with precision $\Delta_{\min} = 2^{-S}$, where $S = K_{\max} \cdot \overline{S}$. 
\item We assume that $\Delta_{\min} \leq \Delta_\epsilon^C$, $ \Delta_\epsilon  \leq \Delta_g^{C}$, and $\Delta_g \leq 2^{- C \overline{S}}$ for a large enough universal constant $C$.
\end{itemize}

Whenever we refer to a  machine number in this section, we mean an $S$-bit machine number, with $S$ as above.

We call a dyadic cuboid $Q = \prod_{j=1}^n I_j \subset \R^n$ a ``machine cuboid'' if each $I_j$ is an interval of the form $[a_j,b_j)$, where $a_j$ and $b_j$ are machine numbers. Recall that each $I_j$ is contained in $[0,\infty)$, by definition of cuboids (see Section \ref{algs1}).

Let $Q,Q'$ be given machine cuboids. The following task can be performed using one unit of ``work'':
\begin{equation}
\label{task} \text{Compute the smallest machine cuboid } Q \text{ containing both } Q' \; \mbox{and} \; Q''.
\end{equation}
Let us explain why we charge only one unit of work to perform the task \eqref{task}.

We suppose that a non-negative machine number $x$ is represented in the computer by its binary digits $(x_i)_{-S  \leq i \leq S}$, where 
\[
x = \sum_{i = - S}^{ + S} x_i 2^i \;\; \mbox{and each } x_i \in \{0,1\}.
\]
We suppose that the bit pattern $(x_i)_{- S \leq i \leq S}$ fits in a single machine word. Given two distinct non-negative machine numbers $x,y$ with binary digits $(x_i)_{-S \leq i \leq S}$, $(y_i)_{-S \leq i \leq S}$ respectively, we return the largest $i_*$ for which $x_{i_*} \neq y_{i_*}$. Recall that in our model of computation for finite-precision arithmetic, we assume that the computation of $i_*$ from $(x_i)$ and $(y_i)$ takes one unit of ``work''. (See Section \ref{sec_moc2}.) Moreover, there are computers in use today for which the computation of $i_*$ from $(x_i)$ and $(y_i)$ may be accomplished by executing $O(1)$ assembly language instructions.

Note that the smallest smallest dyadic interval containing $x$ and $y$ has length $2^{i_*}$. It follows easily that the task \eqref{task} may be accomplished using at most $C$ operations. That is why we consider it reasonable to charge one unit of ``work'' to carry out \eqref{task}.

Therefore, we can determine whether $Q < Q'$, $Q' < Q$, or $Q = Q'$, using $\cO(1)$ computer operations. We refer here to the order relation $<$ on dyadic cuboids defined in Section \ref{sec_dc}.

We should point out that the task \eqref{task} appears to require more than $\cO(1)$ operations in several standard models of computation (not used here). See the discussion of ``quad trees'' and ``segment trees'' in \cite{BCKO}.

We will obtain modified versions of the algorithms in Section \ref{algs1} for our finite-precision model of computation.

We mention a common modification we will make to the statements of many of the algorithms in Section \ref{algs1}. All the cuboids that are input data to an algorithm are assumed to be machine cuboids. Similarly, all the cuboids that are produced as output data are guaranteed to be machine cuboids. We can clearly store a machine cuboid on our computer using $\mathcal{O}(1)$ units of storage.

\begin{itemize}
\item \textbf{Modification 1:} We introduce a bit of notation relevant to the notion of DTrees and ADTrees. See the discussion in Section \ref{pdef_sec}.

Recall that each node $x$ of a DTree $T$ is marked with a dyadic cuboid $Q_x$. When we speak of a DTree $T$ in this section, it is assumed that $Q_x$ is a machine cuboid for each $x \in T$.

Recall that each node $x$ of an ADTree $T$ is marked with linear functionals $\mu_1^x,\cdots,\mu_D^x$ on $\R^D$. We write $\displaystyle \mu_i^x : (v_1,\cdots,v_D) \mapsto \sum_{j=1}^D \theta_{ij}^x v_j$.

We will assume that $D \leq C$ for a universal constant $C$, in what follows.

We say that $\mu_1^x,\cdots,\mu_D^x$ are specified with parameters $(\Delta_g,\Delta_\epsilon)$ if $\lvert \theta_{ij}^x \rvert \leq \Delta_g^{-1}$ and if each $\theta_{ij}^x$ is specified to precision $\Delta_{\epsilon}$. If that's the case for each node $x$ and if the number of nodes of the ADTree is at most $\Delta_g^{-1}$, then we say that the ADTree $T$ is specified with parameters $(\Delta_g,\Delta_\epsilon)$.

\item All of the algorithms that involve BTrees are combinatorial in nature, hence they remain the  same in our finite-precision model of computation. In particular, \textsc{BTree1} and \textsc{Make Control Tree} (deluxe edition and paperback edition) are unchanged.

\item \textbf{Modification 2:} We make the following changes to the algorithm \textsc{Make Control Tree (Hybrid Version)} (see Section \ref{control_tree}).

Assume that an ADTree $T$ is given with parameters $(\Delta_g,\Delta_\epsilon)$, with each node $x$ in $T$ marked by linear functionals $\mu_1^x,\cdots,\mu_D^x$ on $\R^D$. Also, we are given a machine number $\Delta \in [ \Delta_g,1]$ of the form $\Delta = 2^{-K \overline{S}}$ for an integer $K \geq 1$.

Then we compute the control tree $\CT(T)$, with all its markings except for the trees $\BT(\xi)$ ($\xi \in \CT(T)$). For each node $\xi \in \CT(T)$, we compute functionals $\mu^\xi_1,\cdots,\mu^\xi_D : \R^D \rightarrow \R$ of the form
\[\mu_i^\xi : (v_1,\cdots,v_D) \mapsto \sum_{j=1}^D \theta_{ij}^\xi v_j.\]
The numbers $\theta_{ij}^\xi$ are computed with parameters $(\Delta_g^C,\Delta_g^{-C} \Delta_\epsilon)$. That is, we guarantee that $\lvert \theta_{ij}^\xi \rvert \leq \Delta_g^{-C}$ and each $\theta_{ij}^\xi$ is computed to precision $\Delta_g^{-C} \Delta_\epsilon$. We guarantee that for each $\xi \in \CT(T)$ we have
\begin{equation}
\label{controltree_fin}
c\sum_{i=1}^D \lvert \mu_i^\xi(v) \rvert^p \leq \sum_{x \in \BT(\xi) } \sum_{i=1}^D \lvert \mu_i^x(v) \rvert^p + \Delta^p \lvert v \rvert^p \leq C\sum_{i=1}^D \lvert \mu_i^\xi(v) \rvert^p.
\end{equation}
Recall that we denote $\lvert v \rvert^p = \sum_{j} \lvert v_j \rvert^p$ for $v=(v_1,\cdots,v_D)$.

The work and storage requirements are the same as before.

That completes the list of modifications to the hybrid version of \textsc{Make  Control Tree}.

To obtain this result, we apply the finite-precision version of \textsc{Compress Norms} (see Section \ref{sec_cn_fin}) where before we used its infinite-precision counterpart. The proof of \eqref{controltree_fin} is exactly as before.

\item \textbf{Modification 3:} In the algorithm \textsc{Encapsulate}:  Assume that $T$ is a DTree with $N$ nodes such that each node $x$ in $T$ is marked with a machine cuboid $Q_x$.  We perform $C N ( 1 + \log N)$ one-time work in space $C N$ after which we can answer queries. A query consists of a machine cuboid $Q$. The response to a query is an encapsulation $S$ of $Q$, consisting of at most $C + C \log N $ nodes of $\CT(T)$. The work and storage used to answer a query are at most $C+ C \log N$, where $C$ denotes a constant depending only on the dimension $n$.

For the explanation of the algorithm, just note that one can compare two machine cuboids to determine whether one contains the other, using at most $C$ units of work. Thus, we can proceed as in the infinite-precision version of the algorithm \textsc{Encapsulate} using our finite-precision computer.

\item \textbf{Modification 4:} In the algorithm \textsc{ADProcess} (see Section \ref{encaps}):

We assume our ADTree $T$ is given with parameters $(\Delta_g,\Delta_\epsilon)$. We are given a machine number $\Delta \in [ \Delta_g,1]$ of the form $\Delta = 2^{-K \overline{S}}$ for an integer $K \geq 1$.

A query consists of a machine cuboid $Q$. The response to a query is a list of linear functionals $\mu_1^Q,\cdots,\mu_D^Q$ on $\R^D$ such that
\begin{align}
\label{errbd}
c \sum_{i=1}^D \lvert \mu_i^Q(v) \rvert^p &\leq \sum_{\substack{ x \in T \\ Q_x \subset Q}} \sum_{i=1}^D \lvert \mu_i^x(v) \rvert^p + \Delta^p \log(\Delta_g^{-1}) \lvert v \rvert^p \\
&\leq C  \left[ \sum_{i=1}^D \lvert \mu_i^Q(v) \rvert^p + \Delta^p \log(\Delta_g^{-1}) \lvert v \rvert^p \right] \quad \mbox{for all} \; v \in \R^D. \notag{}
\end{align}

Each $\mu_i^Q$ has the form $\displaystyle \mu_i^Q : (v_1,\cdots,v_D) \rightarrow \sum_{j=1}^D \theta_{ij}^Q v_j$. We compute each $\theta_{ij}^Q$ with parameters $(\Delta_g^C, \Delta_g^{-C} \Delta_\epsilon)$.

The work and storage requirements remain the same as before.

That completes the list of modifications to \textsc{ADProcess}.

We present the modifications needed in the explanation of the algorithm. We use the finite-precision versions of the algorithms \textsc{Make Control Tree (Hybrid Version)} and \textsc{Compress Norms} in place of their infinite-precision counterparts. As part of the one-time work, we compute the control tree $\CT(T)$. Each node $\xi \in \CT(T)$ is marked with linear functionals $\mu_1^\xi,\cdots,\mu_D^\xi$ which satisfy \eqref{controltree_fin}. We compute the functionals $\mu_k^\xi$ with parameters $(\Delta_g^C,\Delta_g^{-C} \Delta_\epsilon)$.

Using the algorithm \textsc{Encapsulate}, we respond to a query as follows.

Given a machine cuboid $Q$, we produce a set $S$ of at most $C + C \log N$ nodes in $\CT(T)$ such that $\left\{ x \in T : Q_x \subset Q \right\}$ is the disjoint union over $\xi \in S$ of $\BT(\xi)$. Therefore, by the finite-precision version of \textsc{Make Control Tree (Hybrid Version)} (see \textbf{Modification 2} above), the expression 
\begin{align*}
\mathfrak{E}_1 & = \sum_{\substack{x \in T \\ Q_x \subset Q}} \sum_{i=1}^D \lvert \mu_i^x(v) \rvert^p + \#(S) \cdot \Delta^p \cdot \lvert v \rvert^p \\
&= \sum_{\xi \in S} \left[ \sum_{x \in \BT(\xi)} \sum_{i=1}^D \lvert \mu_i^x(v) \rvert^p + \Delta^p \cdot \lvert v \rvert^p \right]  \\
\end{align*}
differs by at most a factor of $C$ from 
\[
\mathfrak{E}_2 = \sum_{\xi \in S} \sum_{i=1}^D \lvert \mu_i^\xi(v) \rvert^p.
\]

Applying \textsc{Compress Norms} (finite-precision) (see Section \ref{sec_cn_fin}) to the expression $\mathfrak{E}_2$, we compute linear functionals $\mu_1^Q,\cdots,\mu_D^Q$ such that $\mathfrak{E}_2 + \Delta^{p} \cdot \lvert v \rvert^p$ differs by at most a factor of $C$ from $\displaystyle \sum_{i=1}^D \lvert \mu_i^Q(v) \rvert^p$. Each functional $\mu_i^Q$ is computed with parameters $(\Delta_g^C,\Delta_g^{-C}\Delta_\epsilon)$.

Therefore, $\displaystyle \sum_{i=1}^D \lvert \mu_i^Q(v) \rvert^p$ differs by at most a factor of $C$ from
\[
\mathfrak{E}_1 + \Delta^p \lvert v \rvert^p = \sum_{\substack{x \in T \\ Q_x \subset Q}} \sum_{i=1}^D \lvert \mu_i^x(v) \rvert^p + \left( \#(S) + 1 \right) \cdot \Delta^p \cdot \lvert v \rvert^p
\]
Note that $\#(S) \leq C + C \log ( \#(T)) \leq C \log (\Delta_g^{-1})$, so the junk term $\left(\#(S)+1 \right) \cdot \Delta^p \cdot \lvert v \rvert^p$ is bounded by $C \Delta^p \log (\Delta_g^{-1}) \lvert v \rvert^p$. That concludes the proof of \eqref{errbd}.

The work and storage requirements are as promised.

\item \textbf{Modification 5:} In the algorithms \textsc{Make Forest}, \textsc{Fill in Gaps}, and \textsc{Make DTree}: All dyadic cuboids are assumed to be machine cuboids. The explanations of these algorithms require no modification.

\item \textbf{Modification 6:} In the algorithm \textsc{Compute Norms From Marked Cuboids}:

We suppose our cuboids $Q_1,\cdots,Q_N$ have corners whose coordinates are $\til{S}$-bit machine numbers, with $\til{S} \leq C \overline{S}$ for a universal constant $C$. Hence, $N \leq 2^{C n \overline{S}} \leq \Delta_0^{-C}$, where we set $\Delta_0 = 2^{- \overline{S}}$.

We are given a machine number $\Delta \in [ \Delta_g,1]$ of the form $\Delta = 2^{-K \overline{S}}$ for an integer $K \geq 1$.

Each linear functional $\mu_\ell^{Q_i}$ is given as $\displaystyle \mu_\ell^{Q_i} : (v_1,\cdots,v_D) \mapsto \sum_{j=1}^D \theta_{\ell j}^i v_j$. The numbers $ \theta_{\ell j}^i$ are specified with parameters $(\Delta_g,\Delta_{\epsilon})$. We assume that $\displaystyle \widehat{N} := \sum_{i=1}^N (L_i + 1) \leq \Delta_g^{-1}$.

A query consists of a dyadic cuboid $Q$ whose corners are machine points.

The response to a query $Q$ is a list of linear functionals $\widehat{\mu}_1^Q,\cdots,\widehat{\mu}_D^Q : \R^D \rightarrow \R$ for which we guarantee the estimate
\begin{align*}
c \sum_{j=1}^D \lvert \widehat{\mu}_j^Q(v) \rvert^p & \leq \sum_{Q_i \subset Q} \sum_{j=1}^{L_i} \lvert \mu_j^{Q_i}(v) \rvert^p + \Delta^p \Delta_0^{-C} \log(\Delta_g^{-1}) \cdot \lvert v \rvert^p \\
&\leq C \left[  \sum_{j=1}^D \lvert \widehat{\mu}_j^Q(v) \rvert^p + \Delta^p \Delta_0^{-C} \log(\Delta_g^{-1}) \cdot \lvert v \rvert^p \right] \quad \mbox{for all} \; v \in \R^D.
\end{align*}

Each $\widehat{\mu}_i^Q$ has the form $\displaystyle \mu_i^{Q} : (v_1,\cdots,v_D) \mapsto \sum_{j=1}^D \theta_{i j}^Q v_j$. The numbers $\theta_{i j}^Q$ are computed with parameters $(\Delta_g^C, \Delta_g^{-C}\Delta_{\epsilon})$.

That completes the list of modifications to the algorithm \textsc{Compute Norms From Marked Cuboids}.

The explanation of the algorithm is as follows.

For each $i=1,\cdots,N$ we apply the finite-precision version of \textsc{Compress Norms} (see Section \ref{sec_cn_fin}) to produce  linear functionals $\overline{\mu}_j^{Q_i}$ on $\R^D$ for $1 \leq  j \leq  D$ such that
\begin{equation}
\label{latest1}
c \cdot \sum_{j=1}^D \lv \overline{\mu}_j^{Q_i}(v) \rv^p \leq \sum_{j=1}^{L_i} \lv \mu_j^{Q_i}(v) \rv^p  + \Delta^p \cdot  \lv v \rv^p \leq C \cdot \sum_{j=1}^D \lv \overline{\mu}_j^{Q_i}(v) \rv^p.
\end{equation}
Note that each $L_i \leq \Delta_g^{-1}$ by assumption, so the algorithm may be applied as stated. The functionals $\overline{\mu}_j^{Q_i}$ are given with parameters $(\Delta_g^C,\Delta_g^{-C} \Delta_\epsilon)$.

Using the algorithm \textsc{Make Dtree}, we construct a DTree $T$ with at most $C N$ nodes, such that each $Q_i$ is a node of $T$. We mark each $Q_i$ in $T$ with the list of functionals $\overline{\mu}_1^{Q_i}, \cdots, \overline{\mu}_D^{Q_i}$, and we mark all the other nodes in $T$ with a list of linear functionals that are all zero. When equipped with these markings, $T$ forms an ADTree. Note that $\#(T) \leq C N \leq C \Delta_0^{-C} \leq \Delta_g^{-1}$. Hence, the ADTree $T$ is specified with parameters $(\Delta_g^C,\Delta_g^{-C} \Delta_\epsilon)$ (recall that $\overline{\mu}_j^{Q_i}$ are specified with such parameters).

We apply the algorithm \textsc{ADProcess} to the ADTree $T$. Thus, given a machine cube $Q$, we can compute a list of linear functionals $\widehat{\mu}_1^Q, \cdots, \widehat{\mu}_D^Q$ on $\R^D$ such that
\begin{align*}
c \cdot \sum_{j=1}^D \lv \widehat{\mu}_j^Q(v) \rv^p  &\leq \sum_{\substack{ 1 \leq i \leq N \\ Q_i \subset Q}} \sum_{j=1}^D \lv \overline{\mu}_j^{Q_i}(v) \rv^p + \Delta^p \log(\Delta_g^{-C}) \cdot \lv v \rv^p \\
& \leq C \cdot \left[ \sum_{j=1}^D \lv \widehat{\mu}_j^Q(v) \rv^p + \Delta^p \log(\Delta_g^{-C}) \cdot \lv v \rv^p \right].
\end{align*}
Note that $\Delta \in [ \Delta_g^C,1 ]$, so the algorithm may be applied as stated. Using  \eqref{latest1}, we determine that
\begin{align*}
c \cdot \sum_{j=1}^D \lv \widehat{\mu}_j^Q(v) \rv^p  &\leq \sum_{\substack{ 1 \leq i \leq N \\ Q_i \subset Q}}\sum_{j=1}^{L_i} \lv \mu_j^{Q_i}(v) \rv^p   +  \mathfrak{E}(v)  \\
& \leq C \cdot \left[ \sum_{j=1}^D \lv \widehat{\mu}_j^Q(v) \rv^p + \mathfrak{E}(v) \right],
\end{align*}
where $\displaystyle \mathfrak{E}(v) = \sum_{\substack{ 1 \leq i \leq N \\ Q_i \subset Q}} \Delta^p \cdot \lv v \rv^p + \Delta^p \log(\Delta_g^{-C}) \cdot \lv v \rv^p$. Since $N \leq \Delta_0^{-C}$, we conclude that $\mathfrak{E}(v) \leq \Delta^p \Delta_0^{-C} \log(\Delta_g^{-C}) \cdot \lv v \rv^p \leq \Delta^p \Delta_0^{-C'} \log (\Delta_g^{-1}) \cdot \lv v \rv^p$. Thus, the previous estimate implies the desired condition on the functionals $\widehat{\mu}_1^Q, \cdots, \widehat{\mu}_D^Q$.

This completes the explanation of the finite-precision version of the algorithm \textsc{Compute Norms From Marked Cuboids}.

\item The algorithm \textsc{Placing a Point Inside Target Cuboids} requires minor changes in finite-precision: We assume that  $Q_1,\cdots,Q_N$ are machine cuboids, and that the query $\underline{x}$ is a machine point. The response to a query $\underline{x}$ is either one of the $Q_i$ containing $\underline{x}$, or else a promise that no such $Q_i$ exists. The work to answer a query is at most $C \cdot ( 1 + \log N)$. The explanation of the algorithm requires no modification.
\end{itemize}

The above bullet points conclude the description of the changes needed in Section \ref{algs1}.

We make several additional assumptions in Sections \ref{sec_CK}, \ref{sec_bbd}, \ref{clusters}, \ref{sec_ptkc}. Aside from these assumptions all the relevant algorithms are unchanged in our finite-precision model of computation. The new assumptions are as follows.
\begin{itemize}
\item Each point $x \in E$ is an $\overline{S}$-bit machine point (i.e., the coordinates of $x$ are $\overline{S}$-bit machine numbers). 
\item All numbers are $S$-bit machine numbers and all given points are $S$-bit machine points, where $S = K_{\max} \overline{S}$.
\item Each dyadic cuboid has $\til{S}$-bit machine points as corners, where $\til{S} \leq C \overline{S}$.

\end{itemize}

\section{CZ Decompositions}
\label{sec_cz_fin}

We describe the modifications required in Section \ref{sec_czdecomp}.

We let $\Delta_{\min} = 2^{-S}$ denote the machine precision of our computer, where $S = K_{\max} \overline{S}$.

We call a dyadic cube $Q = \prod_{k=1}^n I_k \subset \R^n$ a ``machine cube'' if each $I_k$ is an interval of the form $[a_k,b_k)$, where $a_k$ and $b_k$ are machine numbers.

\begin{itemize}

\item In Sections \ref{sec_rkr} and \ref{sec_czoracle}, we are given a subset $E \subset Q^\circ = [0,1)^n$. We assume that each point in $E$ is an $\overline{S}$-bit machine point. Hence, $\lvert x - y \rvert \geq \Delta_0 = 2^{-\overline{S}}$ for distinct $x,y \in E$.

\item

In Section \ref{sec_czoracle}, we are given a list  $\vec{\Delta} = (\Delta(x))_{x \in E}$ consisting of positive real numbers. We assume that the numbers $\Delta(x)$ are $\til{S}$-bit machine numbers, where $\til{S} \leq C \overline{S}$ for a universal constant $C \geq 1$. Hence, $\Delta(x) \geq 2^{-\til{S}}  \geq \Delta_0^C$ for all $x \in E$. 
\end{itemize}

We recall the definition of the Calder\'on-Zygmund decomposition $\CZ(\vec{\Delta})$ given in Section \ref{sec_czdecomp}:
\begin{itemize}
\item $\CZ(\vec{\Delta})$ consists of the maximal dyadic cubes $Q \subset Q^\circ$ such that either $\#(E \cap 3Q) \leq 1$ or $\Delta(x) \geq \delta_Q$ for all $x \in E \cap 3Q$.
\end{itemize}

Note that each $Q \in \CZ(\vec{\Delta})$ either satisfies $Q = Q^\circ$ or $\#(9Q \cap E) \geq \#(3 Q^+ \cap E) \geq 2$, and that $\lv x - y \rv \geq \Delta_0$ for any distinct $x,y \in E$. Hence, $\delta_Q \geq c \cdot \Delta_0$ for any $Q \in \CZ(\vec{\Delta})$. Therefore, the cubes in $\CZ(\vec{\Delta})$ have $\til{S}$-bit machine points as corners, where $\til{S} \leq C \overline{S}$.

We modify the \textsc{Plain Vanilla CZ-Oracle} to operate as follows: Given an $S$-bit machine point $ \underline{x} \in Q^\circ$, we respond with the cube $Q \in \CZ(\vec{\Delta})$ that contains $\underline{x}$.

We are given a dyadic decomposition $\CZ_{\old}$ as in Section \ref{sec_czoracle}. We assume that each $Q \in \CZ_{\old}$ is a machine cube. We assume we have available a $\CZ_{\old}$-\textsc{Oracle}: Given an $S$-bit machine point $\underline{x} \in Q^\circ$, we produce the unique cube $Q \in \CZ_{\old}$ that contains $\underline{x}$.

Let us define a dyadic decomposition $\CZ_\new$ as in Section \ref{sec_czoracle}: Let  $\CZ_\new$ consist of the maximal dyadic cubes $Q \subset Q^\circ$ such that either $Q \in \CZ_{\old}$ or $\Delta(x) \geq \delta_Q$ for all $x \in E \cap 3Q$. Thus, the decomposition $\CZ_{\old}$ refines the decomposition $\CZ_{\new}$. Since the cubes in $\CZ_\old$ are  machine cubes, we see that the cubes in $\CZ_\new$ are machine cubes.

We make the following minor changes to the statement of the \textsc{Glorified CZ-Oracle}. A query consists of an $S$-bit machine point $ \underline{x} \in Q^\circ$. The response to a query is a list of the cubes $Q \in \CZ_\new$ such that $\underline{x} \in \frac{65}{64} Q$. The explanation of the finite-precision algorithm is unchanged.

We now turn to Section \ref{czalg_sec} and Section \ref{sec_pou}. 

Let $E \subset \frac{1}{32} Q^\circ$, with $\#(E) = N \geq 2$. We assume that $E$ consists of $\overline{S}$-bit machine points. 

We are given a locally finite collection $\CZ$, consisting of dyadic cubes, that forms a partition of $Q^\circ$ (or $\R^n$).

We do not list all the cubes in $\CZ$. Instead, we are given a $\CZ$-\textsc{Oracle}, which operates as follows: Given an $S$-bit machine point $\underline{x} \in Q^\circ$ (or $\underline{x} \in \R^n$), we list all the cubes $Q \in \CZ$ such that $\underline{x} \in \frac{65}{64} Q$. We guarantee that each such $Q$ is an $\til{S}$-bit machine cube with $\til{S} \leq C \overline{S}$. We charge at most $C \log N$ units of work to answer a query.

Under the previous assumptions, the finite-precision versions of the algorithms \textsc{Find Neighbors} and \textsc{Find Main-Cubes} are unchanged.

We modify the algorithm \textsc{Compute Cutoff Function} as follows: We are given machine numbers $\Delta_\epsilon$ and $\Delta_g$, which are large integer powers of $\Delta_0 = 2^{-\overline{S}}$. Given a machine number $\overline{r} \in (\Delta_0,1/64]$ ($\Delta_0 = 2^{- \overline{S}}$), an $S$-bit machine point $\underline{x} \in Q^\circ$, and a machine cube $Q \in \CZ$, we compute the numbers $\frac{1}{\alpha!} \partial^\alpha \left(J_{\underline{x}} \widetilde{\theta}_Q \right)(0)$ (all $\alpha \in \cM$) with parameters $(\Delta_g,\Delta_\epsilon)$.

\section{Starting the Induction}

We begin the proof of the finite-precision version of the Main Technical Results for $\cA$. (See Section \ref{sec_mainresults_fin}.)

We proceed by induction on the multiindex sets $\cA \subset \cM$. Recall that the the collection of multiindex sets carries a total order relation $<$. (See Section \ref{sec_multi}.)

For the base case of the induction, we must prove the Main Technical Results for $\cA = \cM$.

Recall that $\Delta_g(\cM)$ and $\Delta_\epsilon(\cM)$ are assumed to be integer powers of $\Delta_0 = 2^{-\overline{S}}$ that satisfy $\Delta_\epsilon(\cM) \ll \Delta_g(\cM)$ (see \eqref{dominate}). We denote $\Delta_g = \Delta_g(\cM)$ and $\Delta_\epsilon = \Delta_\epsilon(\cM)$ in the course of this section. Thus, we may assume that $\Delta_\epsilon \leq \Delta_g^C$ for a sufficiently large universal constant $C$.

We define $\CZ(\cM)$ to be the collection of the maximal dyadic cubes $Q \subset Q^\circ$ such that $\#(E \cap 3Q) \leq 1$.

From \eqref{farapart}, we know that 
\begin{equation}
\label{deltamin}
\delta_{\min} = \min \left\{ \lv x - y \rv : x,y \in E, \; x \neq y \right\} \geq  \Delta_0.
\end{equation}
Thus, we can compute a machine number  $\delta$ with $c \cdot \delta_{\min} < \delta < \frac{1}{100} \delta_{\min}$ using the BBD Tree (see Theorem \ref{bbd_thm}). This data structure requires no modifications for finite-precision computation. We construct the \textsc{Plain Vanilla $\CZ$-Oracle} in Section \ref{sec_cz_fin}, where we take $\Delta(x) = \delta$ for all $x \in E$. This yields a $\CZ(\cM)$-\textsc{Oracle} as in the Main Technical Results for $\cM$. See Remark \ref{rem_cw}.

Note that $\#(E \cap 3Q^+) \geq 2$ for each $Q \in \CZ(\cM)$, where $Q^+$ is the dyadic parent of $Q$. We can thus choose distinct points $x,x' \in E \cap 3Q^+$. From \eqref{farapart} we know that $\lv x - x' \rv \geq \Delta_0$. Hence, $6 \delta_Q = \delta_{3Q^+} \geq \Delta_0$. In particular,
\[
\delta_Q \geq \frac{1}{32} \Delta_0 \;\; \mbox{for each} \; Q \in \CZ(\cM).
\]
Thus the decomposition $\CZ(\cM)$ satisfies the additional property \eqref{lowed} required in finite-precision.

Using the algorithm \textsc{Find Main-Cubes} (see Section \ref{sec_cz_fin}), we list all the cubes $Q \in \CZ_{\main}(\cM)$, and we compute a point $x(Q) \in E \cap \frac{65}{64}Q$ associated to each $Q \in \CZ_{\main}(\cM)$.

The definitions of the linear maps $T_{(Q,\cM)}$ and linear functionals $\xi_Q$ in \eqref{base1} and \eqref{base2} are unchanged. We take the assist set $\Omega(Q,\cM)$ to be empty for each $Q \in \CZ_{\main}(\cM)$.

For each $Q \in \CZ_{\main}(\cM)$, the linear functional $\xi_Q$ is computed in short form with parameters $(\Delta_g,\Delta_\epsilon)$. Furthermore, given an $S$-bit machine point $\underline{x} \in  Q^\circ$ and a multiindex $\alpha \in \cM$, we compute the linear functional $(f,P) \mapsto \partial^\alpha(T_{(Q,\cM)}(f,P))(\underline{x})$ with parameters  $(\Delta_g,\Delta_\epsilon)$. This completes the description of the changes to the algorithm \textsc{Find Main-Cubes and Compute Extension Operators (Base Case)}. The explanation of this algorithm is obvious once we examine the formulas for $T_{(Q,\cM)}$ and $\xi_Q$ in \eqref{base1} and \eqref{base2}. 

That concludes the proof of the finite-precision version of the Main Technical Results for $\cA = \cM$. This completes the base case of our induction. Next, we turn to the induction step.

\section{The Induction Step}
\label{sec_ind_step_fin}

Let $\cA \subsetneq \cM$ be a given multiindex set. Let $\cA^- < \cA$ be maximal with respect to the order $<$ on multiindex sets. By induction we assume that the Main Technical Results for $\cA^-$ hold. We list below a few consequences of these results. (See Section \ref{sec_mainresults_fin}.)

We denote $\Delta_g := \Delta_g(\cA^-)$, $\Delta_\epsilon := \Delta_\epsilon(\cA^-)$, and $\Delta_{\junk} := \Delta_\junk(\cA^-)$, which are the parameters arising in the Main Technical Results for $\cA^-$. These parameters are all large integer powers of $\Delta_0 = 2^{- \overline{S}}$. These parameters are not fixed just yet. Hence, in the course of the proof of the Main Technical Results for $\cA$ we may  impose additional assumptions on the relationships between these parameters. From \eqref{constants1}, we may assume estimates of the form
\begin{equation}
\label{constants2}
\Delta_\epsilon \leq \Delta_g^C, \; \Delta_g \leq \Delta_\junk^C, \;\mbox{and} \; \Delta_\junk \leq \Delta_0^C, \;\; \mbox{for a universal constant } C.
\end{equation}
For example, the first estimate in \eqref{constants2} is derived from \eqref{constants1} as follows:
\[
\Delta_\epsilon(\cA^-) \leq \Delta_\epsilon(\emptyset) \leq \Delta_g(\emptyset)^C \leq \Delta_g(\cA^-)^C.
\]
The other conditions can be derived in a similar fashion.

By the induction hypothesis, we have defined a dyadic decomposition $\CZ(\cA^-)$ of the unit cube $Q^\circ$, which satisfies conditions \textbf{(CZ1)}-\textbf{(CZ5)} in the Main Technical Results for $\cA^-$. Furthermore, from the finite-precision version of these results, we have
\[\delta_Q \geq \frac{1}{32} \cdot \Delta_0 \;\; \mbox{for each} \; Q \in \CZ(\cA^-).
\]
Recall that  $\CZ_{\main}(\cA^-)$ is the collection of all the cubes $Q \in \CZ(\cA^-)$ such that $\frac{65}{64}Q \cap E \neq \emptyset$.

According to the Main Technical Results for $\cA^-$, we have computed a list of all the cubes in $\CZ_{\main}(\cA^-)$. Furthermore, we have access to a $\CZ(\cA^-)$-\textsc{Oracle} that operates as follows: Given an $S$-bit machine point  $\underline{x} \in Q^\circ$, we list all the cubes $Q \in \CZ(\cA^-)$ such that $\underline{x} \in \frac{65}{64}Q$, using work at most $C \log N$.

We have computed a list of assists $\Omega(Q,\cA^-)$, and a list of assisted functionals $\Xi(Q,\cA^-)$ for each $Q \in \CZ_{\main}(\cA^-)$. Each $\omega \in \Omega(Q,\cA^-)$ is a linear functional on $\X( \frac{65}{64}Q \cap E)$, and is given in short form with parameters $(\Delta_g,\Delta_\epsilon)$; each $\xi \in \Xi(Q,\cA^-)$ is a linear functional on $\X(\frac{65}{64}Q \cap E) \oplus \cP$, and is given in short form with parameters $(\Delta_g,\Delta_\epsilon)$  in terms of the assists $\Omega(Q,\cA^-)$. We guarantee that
\[c \cdot \|(f,P) \|_{(1+a(\cA^-))Q} \leq M_{(Q,\cA^-)}(f,P) \leq C \cdot \left[ \| (f,P)\|_{\frac{65}{64}Q} + \Delta_\junk \lvert P \rvert \right]\]
where
\[M_{(Q,\cA^-)}(f,P) := \left( \sum_{\xi \in \Xi(Q,\cA^-)}  \lvert \xi(f,P) \rvert^p \right)^{1/p}.\]
We recall that $\displaystyle \lv P \rv = \left( \sum_{\alpha \in \cM} \lv \partial^\alpha P(0) \rv^p \right)^{1/p}$.

We have computed a linear map $T_{(Q,\cA^-)} : \X(\frac{65}{64}Q \cap E) \oplus \cP \rightarrow \X$  for each $Q \in \CZ_{\main}(\cA^-)$ in the following sense: Given an $S$-bit machine point $\underline{x} \in Q^\circ$ and a multiindex $\alpha \in \cM$, we compute the linear functional
\[
(f,P) \mapsto \partial^\alpha( T_{(Q,\cA^-)}(f,P))(\underline{x})
\] 
in short form with parameters $(\Delta_g,\Delta_\epsilon)$ in terms of the assists $\Omega(Q,\cA^-)$. This computation requires work at most $C \log N$.

The previous computations are carried out using one-time work at most $CN \log N$ in space at most $C N$, thanks to the induction hypothesis.

The finite-precision version of the algorithm \textsc{Approximate Old Trace Norm} is as follows. 

\environmentA{Algorithm: Approximate Old Trace Norm (Finite-Precision).}

For each $Q \in \CZ_{\main}(\cA^-)$ we compute linear functionals $\xi_1^Q,\cdots,\xi_D^Q$ on $\cP$, such that
\[
c \cdot \sum_{i=1}^D \lvert \xi_i^Q(P) \rvert^p \leq\sum_{\xi \in \Xi(Q,\cA^-)} \lvert \xi(0,P) \rvert^p + \Delta_g^p \cdot \lvert P \rvert^p \leq C \cdot \sum_{i=1}^D \lvert \xi_i^Q(P) \rvert^p.
\]
The functionals $\xi_i^Q$ have the form $\displaystyle \xi_i^Q : (P) \mapsto \sum_{\alpha \in \cM} \theta_{i \alpha}^Q \cdot \frac{1}{\alpha!} \partial^\alpha P(0)$. The numbers $\theta_{i \alpha}^Q$ are computed with parameters $(\Delta_g^C,\Delta_g^{-C}\Delta_\epsilon)$.

\begin{proof}[\underline{Explanation}]
The explanation proceeds just as in infinite-precision. We simply apply the finite-precision version of the algorithm \textsc{Compress Norms} (with $\Delta = \Delta_g$) instead of the infinite-precision version of this algorithm. See Section \ref{sec_cn_fin}.
\end{proof}

That completes the description of the Main Technical Results for $\cA^-$.

We now begin the proof of the Main Technical Results for $\cA$.

\subsection{The Non-monotonic Case.}

Section \ref{nonmon_sec} requires no change in finite-precision. For a nonmonotonic set $\cA \subsetneq \cM$, we can simply define $\CZ(\cA) = \CZ(\cA^-)$, and 
\begin{align*}
&\Omega(Q,\cA) = \Omega(Q,\cA^-), \; \Xi(Q,\cA) = \Xi(Q,\cA^-), \;\; \mbox{and} \\
& T_{(Q,\cA)} =T_{(Q,\cA^-)} \; \; \text{for each} \; Q \in \CZ_{\main}(\cA) = \CZ_{\main}(\cA^-).
\end{align*}

As before, the Main Technical Results for $\cA$ (finite-precision) follow from the Main Technical Results for $\cA^-$ (finite-precision). 

We can compute everything to the desired precision provided that $\Delta_\epsilon(\cA) \geq \Delta_\epsilon(\cA^-)$, $\Delta_g(\cA) \leq \Delta_g(\cA^-)$, and $\Delta_\junk(\cA) \geq \Delta_\junk(\cA^-)$. These conditions are allowed in view of the assumptions in \eqref{constants1},

This proves the Main Technical Results for a nonmonotonic set $\cA$.

\subsection{The Monotonic Case}
\label{sec_ocz_fin}

Henceforth, we assume that $\cA$ is a monotonic set. 

The statement of Proposition \ref{newcz_prop} remains the same, except for the following changes.

The definition of the dyadic decomposition $\overline{CZ}(\cA^-)$ in Proposition \ref{newcz_prop} is unchanged. The $\overline{CZ}(\cA^-)$-\textsc{Oracle} operates as follows.

\begin{itemize}
\item Given a query consisting of an $S$-bit machine point $\underline{x} \in \R^n$ such that $\lv \underline{x} \rv \leq 2^{\overline{S}}$, we compute \underline{exactly} a list of the cubes $Q \in \overline{CZ}(\cA^-)$ such that $\underline{x} \in \frac{65}{64} Q$. The work required to answer a query is at most $C \log N$.
\end{itemize}

One can check that the explanation of the $\overline{CZ}(\cA^-)$-\textsc{Oracle} in Proposition \ref{newcz_prop} applies equally well in the finite-precision setting (under the additional hypotheses on $\underline{x}$ stated above). Indeed, since $\lvert \underline{x} \rvert \leq 2^{\overline{S}}$, we see that each of the dyadic cubes that is relevant to the explanation of the algorithm is contained in the rectangular box $[-2^{C \overline{S}}, 2^{C \overline{S}}]^n$ and has sidelength in the interval $[2^{-C\overline{S}},2^{C \overline{S}}]$ for a universal constant $C$. Therefore, the relevant dyadic cubes have $C\overline{S} $-bit machine points as corners. We may assume that $C \overline{S} \leq S$, and therefore all the relevant dyadic cubes are $S$-bit machine cubes and can be processed on our finite-precision computer, which allows the previous explanation to apply in the current setting.

\comments{Changed above paragraph
}

This concludes the description of  changes needed in Section \ref{sec_ocz}.

\subsection{Keystone Cubes}

Section \ref{key_sec} requires no change in finite-precision. We define integer constants $S_0,S_1,S_2$ as in \eqref{consts}. The \textsc{Keystone-Oracle} is unchanged. The explanation follows just as before from the \textsc{Main Keystone Cube Algorithm} and the algorithm \textsc{List All Keystone Cubes}.

\section{An Approximation to the Sigma}
\label{sec_appx}

Given a polynomial $P \in \cP$, we define
\[ 
\lvert P \rvert := \left( \sum_{\alpha \in \cM} \lvert \partial^\alpha P(0) \rvert^p \right)^{1/p}.
\]

Recall that the parameters $\Delta_\epsilon(\cA^-)$, $\Delta_g(\cA^-)$, and $\Delta_{\junk}(\cA^-)$, are denoted by $\Delta_\epsilon$, $\Delta_g$, and $\Delta_{\junk}$, respectively.


We denote $\Delta_\new = \Delta_{\junk}(\cA)$, which is the constant arising in the Main Technical Results for $\cA$. We assume that $\Delta_\new$ satisfies
\begin{equation}
\label{constants3}
\left\{
\begin{aligned}
&\Delta_\epsilon \leq \Delta_g \leq \Delta_\junk \leq \Delta_\new^5 \leq \Delta_\new \leq \Delta_0^C, \\
&C \cdot \Delta_\new \leq 1 \qquad \mbox{for a large enough universal constant} \; C.
\end{aligned}
\right.
\end{equation}
The conditions in \eqref{constants3} are easily derived from \eqref{constants1}.

We denote $a = a(\cA^-)$ for the constant in the Main Technical Results for $\cA^-$. Instead of \eqref{mdefn} and \eqref{n_appx} in Section \ref{sec_assign}, we have estimates from the finite-precision version of the Main Technical Results for $\cA^-$. Namely, for each $Q \in \CZ_{\main}(\cA^-)$, the functional
\begin{equation}
\label{mdefn_fin}
M_{(Q,\cA^-)}(f,R) = \left( \sum_{\xi \in \Xi(Q,\cA^-)} \lvert \xi(f,R) \rvert^p \right)^{1/p}
\end{equation}
satisfies
\begin{equation}
\label{n_appx_fin}
c \cdot \| (f,R) \|_{(1+a)Q} \leq M_{(Q,\cA^-)}(f,R) \leq C \cdot \left[ \| (f,R) \|_{\frac{65}{64}Q} + \Delta_\junk \lvert R \rvert \right].
\end{equation}
The estimate \eqref{n_appx_fin} holds for all $Q \in \overline{\CZ}(\cA^-)$, since we assume by definition that $\Xi(Q,\cA^-) = \emptyset$ for all $Q \in \overline{\CZ}(\cA^-) \setminus \CZ_{\main}(\cA^-)$.



We set
\[
\cI(Q^\#) := \bigl\{ Q \in \overline{\CZ}(\cA^-): Q \cap S_0 Q^\# \neq \emptyset \bigr\}.
\]
(Recall that  $S_0 = S(\cA^-)$.)

Lemma \ref{touch_lem} is unchanged in finite-precision. Similarly, the conditions \eqref{touch1_a}, \eqref{touch2_a} continue to hold.

The finite-precision version of the algorithm \textsc{Make New Assists and Assign Keystone Jets} is as follows.

\environmentA{Algorithm: Make New Assists and Assign Keystone Jets (Finite-Precision)}

For each keystone cube $Q^\#$, we compute a list of new assists $\Omega^\new(Q^\#)$ and we compute an $\Omega^{\new}(Q^\#)$-assisted bounded depth linear map $R_{Q^\#}^\# : \X(S_1 Q^\# \cap E) \oplus \cP \rightarrow \cP$. 

Each of the new assists $\omega \in \Omega^\new(Q^\#)$ is given in the form
\[
\omega : f \mapsto \sum_{x \in S} \lambda_x \cdot f(x).
\] 
Here, the set $S \subset E$ may depend on $\omega$. The coefficients $\lambda_x$ are computed with parameters $(\Delta_g^{C}, \Delta_g^{-C} \Delta_\epsilon)$. The sum of $\depth(\omega) = \#(S)$ over all the new assists $\omega \in \Omega^{\new}(Q^\#)$ and over all keystone cubes $Q^\#$, is bounded by $C N$.

Similarly, the maps $(f,P) \mapsto R^\# = R^\#_{Q^\#}(f,P)$ are such that
\begin{align*}
&(f,P) \mapsto \partial^\alpha R^\#(0) \quad (\text{for any } \alpha \in \cA) \;\; \text{has the form} \\
& \sum_{x \in S} \lambda_x \cdot f(x) + \sum_{\omega \in \Omega'} \mu_\omega \cdot \omega(f) + \sum_{\beta \in \cM} \theta_\beta \frac{1}{\beta!} \partial^\beta P(0).
\end{align*}
Here, the subsets $S \subset E$ and $\Omega' \subset \Omega$ may depend on $Q^\#$, and the coefficients $\lambda_x,  \mu_\omega, \theta_\beta$ are computed with parameters $( \Delta_g^{C}, \Delta_g^{-C} \Delta_\epsilon)$. The number  $\#(S) + \#(\Omega') + \# (\cM)$ is bounded by a universal constant $C$.

The polynomial $R^\#$ satisfies the following properties.
\begin{itemize}
\item $\partial^\alpha ( R^\# - P ) \equiv 0$ for all $\alpha \in \cA$. (Recall that $\cA$ is monotonic; see Remark \ref{mon_rem}.)
\item Let $R \in \cP$ with $\partial^\alpha (R - P) \equiv 0$ for all $\alpha \in \cA$. Then
\begin{equation} \label{crab}
\sum_{Q \in \cI(Q^\#)} \sum_{\xi \in \Xi(Q,\cA^-)} \lvert \xi(f,R^\#) \rvert^p  \leq C \cdot \left[ \sum_{Q \in \cI(Q^\#)}  \sum_{\xi \in \Xi(Q,\cA^-)} \lvert \xi(f,R) \rvert^p + \Delta_g^p \cdot \lvert R \rvert^p \right].
\end{equation}
\end{itemize}

\begin{proof}[\underline{Explanation}] 

As before, we define coordinates on $V_P$, which is the space of all polynomials $R \in \cP$ such that $\partial^\alpha \left[ R - P \right] \equiv 0$ for all $\alpha \in \cA$. The coordinate map $w \mapsto R_w$ is given by
\[
R_w(x) = \sum_{\alpha \in \cA} \frac{1}{\alpha!} \partial^\alpha P(0) \cdot x^\alpha + \sum_{j=1}^k \frac{1}{\alpha_j!} w_j \cdot x^{\alpha_j} \qquad \mbox{for} \; w=(w_1,\cdots,w_k) \in \R^k,
\]
where $\cM \setminus \cA = \{\alpha_1,\cdots,\alpha_k\}$. Note that
\begin{equation}
\label{pnorm1}
\lv R_w \rv^p = \sum_{j=1}^k \lvert w_j \rvert^p + \sum_{\beta \in \cA} \lvert \partial^\beta P(0) \rvert^p \geq \sum_{j=1}^k \lvert w_j \rvert^p.
\end{equation}

We compute a list $\mathfrak{L}$ of all the $Q \in \CZ_{\main}(\cA^-)$ such that $Q \cap S_0 Q^\# \neq \emptyset$. We produce this list by the same method used in infinite-precision (recall that $S_0 = S(\cA^-) \in \N$ is a universal constant, as stated in the Main Technical Results for $\cA^-$.)

From the Main Technical Results for $\cA^-$ (finite-precision), we can compute a list of the functionals $\xi_\ell : (f,R_w) \mapsto \xi(f,R_w)$, with $1 \leq \ell \leq L$, where $\xi$ is an arbitrary element of the list $\Xi(Q,\cA^-)$ for some $Q \in \mathfrak{L}$. Each $\xi_\ell$ is expressed in short form with parameters $(\Delta_g,\Delta_\epsilon)$ in terms of the assists $\Omega(Q,\cA^-)$. This means that we have \eqref{xi_short}, where
\begin{itemize}
\item The numbers $\check{\mu}_{\ell j}$ in \eqref{xi_short} are specified with parameters $(\Delta_g,\Delta_\epsilon)$;
\item The functionals $\lambda_\ell$, $\check{\lambda}_\ell$, and the coefficients $\mu_{\ell a}$ in \eqref{xi_short} are specified with parameters $(\Delta_g,\Delta_\epsilon)$. 
\item We have $L \leq C N \leq \Delta_g^{-C}$ for a large enough universal constant $C$. (Recall that $N \leq \Delta_0^{-n} \leq \Delta_g^{-1}$; see \eqref{Nbound} and \eqref{constants2}.)
\end{itemize}

We process the functionals $w \mapsto \xi_\ell(f,R_w)$, with $f$ and $(\partial^\alpha P(0))_{\alpha \in \cA}$ held fixed, using the algorithm \textsc{Optimize via Matrix} (finite-precision), where we set $\Delta = \Delta_g$ (see Section \ref{sec_ovm_fin}). Thus, we can compute (see below) numbers $b_{j \ell}$, such that, for 
\begin{equation}
\label{newfunc}
\left\{
\begin{aligned}
&\omega_j^\new(f) =   \sum_{\ell=1}^L b_{j\ell} \left[ \lambda_\ell(f) + \sum_{a=1}^{I_\ell} \mu_{\ell a} \omega_{\ell a}(f) \right] \\
&\lambda_j^\new((\partial^\alpha P(0))_{\alpha \in \cA}) = \sum_{\ell=1}^L b_{j\ell} \cdot \check{\lambda}_\ell((\partial^\alpha P(0))_{\alpha \in \cA}),
\end{aligned}
\right.
\end{equation}
we have
\begin{equation}
\label{eq2000}
\sum_{Q \in \cI(Q^\#)} \sum_{\xi \in \Xi(Q,\cA^-)} \lvert \xi(f,R_{w^*}) \rvert^p \leq  C \cdot \left[ \sum_{Q \in \cI(Q^\#)} \sum_{\xi \in \Xi(Q,\cA^-)} \lvert \xi(f,R_{w}) \rvert^p + \Delta_g^p\sum_{j=1}^k \lvert w_j \rvert^p \right],
\end{equation}
for all $w=(w_1,\cdots,w_k) \in \R^k$, where 
\begin{align}
\label{wstar}
&w^* = (w^*_1,\cdots,w^*_k), \;\; \mbox{with} \\
&w^*_j = \omega_j^\new(f) + \lambda_j^\new((\partial^\alpha P(0))_{\alpha \in \cA}) \quad \text{for all} \; 1 \leq j \leq k.
\notag{}
\end{align}
This is a consequence of the finite-precision version of the algorithm \textsc{Optimize via Matrix}. We compute the numbers $b_{j \ell}$ with parameters $(\Delta_g^C, \Delta_g^{-C}\Delta_\epsilon)$.

Recall that $w \mapsto R_w$ parametrizes the space $V_P$ of all polynomials $R \in \cP$ with $\partial^\alpha \left[ R - P \right] \equiv 0$. Thus, using  \eqref{pnorm1} and \eqref{eq2000} we see that
\[
\sum_{Q \in \cI(Q^\#)} \sum_{\xi \in \Xi(Q,\cA^-)} \lvert \xi(f,R_{w^*}) \rvert^p \leq  C \cdot \left[ \sum_{Q \in \cI(Q^\#)} \sum_{\xi \in \Xi(Q,\cA^-)} \lvert \xi(f,R) \rvert^p + \Delta_g^p \cdot \lvert R \rvert^p \right],
\]
for any polynomial $R \in V_P$. We can thus set $R^\#_{Q^\#} = R_{w^*}$ and we obtain the estimate \eqref{crab}.

We now produce a numerically accurate formula for the new assists $\omega_j^{\new}$ ($j=1,\cdots,k$) and for the functionals $\lambda_j^\new((\partial^\alpha P(0))_{\alpha \in \cA})$. We examine the relevant definitions.

In the expression for $\lambda_j^\new$ in \eqref{newfunc}, the numbers $b_{j \ell}$ are given with parameters $(\Delta_g^C, \Delta_g^{-C}\Delta_\epsilon)$, the functionals $\check{\lambda}_\ell$ are given with parameters $(\Delta_g,\Delta_\epsilon)$, and $L \leq \Delta_g^{-C}$. Thus, we can compute the functionals $\lambda_j^\new$ with parameters $(\Delta_g^C, \Delta_g^{-C}\Delta_\epsilon)$.

We will review our computation of a short form representation of each $\omega_j^\new(f)$ in \eqref{newfunc}, following the infinite-precision text. We need to document roundoff errors  at each stage of the computation, and estimate the size of the relevant numbers.

We write $\omega_j^\new = \omega_j^{\new,1} + \omega_j^{\new,2}$, with $\omega_j^{\new,1}$  and $\omega_j^{\new,2}$ defined via \eqref{first_func} and \eqref{last_fnc}, respectively.

We first review the computation of $\omega_j^{\new,1}$ in \eqref{first_func}.

\begin{itemize}
\item The numbers $c_\ell(x)$  ($x \in S_\ell \subset E$) in \eqref{cl} are given with parameters $(\Delta_g,\Delta_\epsilon)$, since each $\lambda_\ell$ is given with the same parameters by assumption. The weights $d_j(x)$ in \eqref{new_form1} are computed by evaluating the sum $d_j(x) = \sum_{\ell} b_{j \ell} \cdot c_\ell(x)$. Each term $b_{j \ell} \cdot c_\ell(x)$ is computed to precision $\Delta_g^{-C} \Delta_\epsilon$. Hence, each weight $d_j(x)$ is computed to precision $L \Delta_g^{-C} \Delta_\epsilon \leq C N \Delta_g^{-C} \Delta_\epsilon \leq \Delta_g^{-C'} \Delta_\epsilon$; we compute each $d_j(x)$ by sorting, just as before. Moreover, each $d_j(x)$ satisfies $\lv d_j(x) \rv \leq L \Delta_g^{-C} \leq \Delta_g^{-C'}$. Therefore, we can compute each $d_j(x)$ with parameters $(\Delta_g^{C},\Delta_g^{-C} \Delta_\epsilon)$. The bounds on the work and storage required by this computation are the same as before.
\end{itemize}

We can thus express the functionals $\omega_j^{\new,1}$ in the form \eqref{new_form1}, where the coefficients $d_j(x)$ are given with parameters $(\Delta_g^C,\Delta_g^{-C}\Delta_\epsilon)$. Thus, by definition, we can compute the functional $\omega_j^{\new,1}$ with parameters $(\Delta_g^C,\Delta_g^{-C}\Delta_\epsilon)$.

We now review the computation of $\omega_j^{\new,2}$ in \eqref{last_fnc}.

We express $\omega_j^{\new,2}$ in the form \eqref{new_form2} by sorting. The coefficients $q_{j\omega}$ in  \eqref{new_form2} are computed by evaluating $q_{j \omega} = \sum_{(\ell,a)} b_{j \ell} \cdot \mu_{\ell a}$, where the sum is over certain pairs $(\ell,a)$; see \eqref{q_defn}. The numbers $b_{j\ell}$ are given with parameters $(\Delta_g^{C}, \Delta_g^{-C} \Delta_\epsilon)$, the numbers $\mu_{\ell a}$ are given with parameters $(\Delta_g,\Delta_\epsilon)$, and the number of terms in the sum is bounded by $C N \leq \Delta_g^{-C}$. Hence, the numbers $q_{\ell \omega}$ can be computed with parameters $(\Delta_g^{C}, \Delta_g^{-C} \Delta_\epsilon)$.

We finally express $\omega_j^{\new,2}$ in the form \eqref{new_form3}. As before, the coefficients $k_j(x)$, which we compute by sorting, are given with parameters $(\Delta_g^C,\Delta_g^{-C}\Delta_\epsilon)$.

We have shown how to compute the new assists $\omega^\new_j = \omega^{\new,1}_j + \omega^{\new,j}_2$ in short form with parameters $(\Delta_g^{C},\Delta_g^{-C}\Delta_\epsilon)$. We have seen that computation follows by the same method as in infinite-precision, and by making careful note of the roundoff errors and the size of numbers involved in the computation, we verified that the computation could be carried out on our finite-precision computer.

The functionals $(f,P) \mapsto \partial^\alpha \bigl[ R^\#_{Q^\#}(f,P) \bigr](0)$ can clearly be computed in short form in terms of the assists $\omega_j^{\new}$ ($j=1,\cdots,k$) with parameters $(\Delta_g^C,\Delta_g^{-C}\Delta_\epsilon)$, as desired. (See \eqref{derivativesform}.)

This completes the explanation of the algorithm \textsc{Make New Assists and Assign Keystone Jets} (finite-precision).

\end{proof}

For each $(f,R) \in \X(S_1 Q^\# \cap E) \oplus \cP$ we set
\begin{align} \label{sharp_norm_fin} \left[ M^\#_{Q^\#}(f,R) \right]^p &:= \sum_{\substack{Q \in \cI(Q^\#)}} \sum_{\xi \in \Xi(Q,\cA^-)} \lvert \xi(f,R) \rvert^p  = \sum_{\substack{ Q \in \cI(Q^\#)}}  \left[ M_{(Q,\cA^-)}(f,R)\right]^p.
\end{align}

Let $P \in \cP$. From the previous algorithm, we see that the polynomial $R^\# = R_{Q^\#}^\#(f,P)$ satisfies
\begin{equation}
\label{crab2}
\left\{
\begin{aligned}
&\partial^\alpha( R^\# - P) \equiv 0 \;\; \mbox{for all} \; \alpha \in \cA, \\
&M^\#_{Q^\#}(f,R^\#)  \leq  C \cdot \left[ M^\#_{Q^\#}(f,R) + \Delta_{\junk} \lv R \rv \right],
\end{aligned}
\right.
\end{equation}
for any polynomial $R \in \cP$ such that $\partial^\alpha \left[ R - P \right] \equiv 0$ for all $\alpha \in \cA$. \\
(Recall that $\Delta_g \leq \Delta_{\junk}$.)

We replace Lemma \ref{key_lem1} with the following result.
\begin{lem}\label{key_lem1_fin}
Let $Q^\#$ be a keystone cube. Then
\begin{align*} 
c \cdot \|(f, R)\|_{S_0 Q^\#} &\leq M^\#_{Q^\#}(f,R)  \\
& \leq C \cdot \left[  \|(f,R)\|_{S_1 Q^\#} + \Delta_\junk \lvert R \rvert \right].
\end{align*}
for all $(f,R) \in \X(S_1 Q^\# \cap E) \oplus \cP$. Here, $c>0$ and $C \geq 1$ are universal constants.
\end{lem}
\begin{proof}
The first inequality $c \|(f, R)\|_{S_0 Q^\#} \leq M^\#_{Q^\#}(f,R)$ is proven using the argument in Lemma \ref{key_lem1}. Where before we referred to \eqref{n_appx}, we now refer to \eqref{n_appx_fin}.

We prove the second inequality as follows. From \eqref{n_appx_fin}, we have
\begin{equation*}
\left[ M_{Q^\#}^\#(f,R) \right]^p = \sum_{Q \in \cI(Q^\#) } \left[ M_{(Q,\cA^-)}(f,R) \right]^p \leq C \cdot \sum_{Q \in \cI(Q^\#)} \left[ \|(f,R)\|_{\frac{65}{64}Q}^p + \Delta_\junk^p \lvert R \rvert^p \right]
\end{equation*}
Since $\frac{65}{64}Q \subset S_1 Q^\#$ and $\delta_{Q^\#} \leq \delta_Q \leq C \delta_{Q^\#}$ for each $Q \in \cI(Q^\#)$, and since $\# \bigl[ \cI(Q^\#)  \bigr] \leq C$, by Lemma \ref{lem_normmon}, we can estimate the above line by
\[ C \cdot \left[ \|(f,R)\|^p_{S_1 Q^\#} + \Delta^p_\junk \lvert R \rvert^p \right] \leq C' \cdot \left[ \|(f,R)\|_{S_1 Q^\#} + \Delta_\junk \lvert R \rvert \right]^p.\]
This completes the proof of Lemma \ref{key_lem1_fin}.
\end{proof}

Proposition \ref{key_prop1} requires a few changes in finite-precision. Here is the modified statement:

\begin{prop}\label{key_prop1_fin}
Let $\hQ$ be a dyadic subcube of $Q^\circ$, such that $3\hQ$ is tagged with $(\cA,\epsilon)$.
Assume also that $Q^\# \in \CZ(\cA^-)$ is a keystone cube, and that $S_1 Q^\# \subseteq \frac{65}{64}\hQ$.

If $H \in \X$, $H = f$ on $E \cap S_1Q^\#$, and $\partial^\alpha H(x_{Q^\#}) = \partial^\alpha P(x_{Q^\#})$ for all $\alpha \in \cA$, then
\begin{equation} \label{imp_fin} 
\delta_{Q^\#}^{-m} \| H - R_{Q^\#}^\#  \|_{L^p(S_1 Q^\#)} \leq C \cdot \left[ \| H \|_{\X(S_1 Q^\#)} + \Delta_\junk \lvert J_{x_{Q^\#}} H \rvert \right].
\end{equation}
Here, $C \geq 1$ is a universal constant; and $R^\#_{Q^\#} = R^\#_{Q^\#}(f,P)$. \\
(See the algorithm \textsc{Make New Assists and Assign Keystone Jets}.)
\end{prop}

\begin{proof}

The proof of \eqref{smallsigma}, assuming \eqref{insigma0} and \eqref{zeroder}, is unchanged in the present setting. We need only examine and fix the last paragraph in the proof, right after the paragraph containing \eqref{e732}.

We modify the text that begins after the sentence ``We now prove the main assertion in Proposition \ref{key_prop1}''. The revised  text is as follows.
Suppose that $H \in \X$ satisfies $H = f$ on $E \cap S_1 Q^\#$ and $\partial^\alpha H(x_{Q^\#}) = \partial^\alpha P(x_{Q^\#})$ for $\alpha \in \cA$. Then $\partial^\alpha (J_{x_{Q^\#}} H - P ) \equiv 0$ for all $\alpha \in \cA$. (Recall, $\cA$ is monotonic.) We apply the estimate in \eqref{crab2}, and then we apply Lemma \ref{key_lem1_fin}. Therefore,
\begin{align*}
M_{Q^\#}^\#(f,R^\#_{Q^\#}) &\leq C \left[ M_{Q^\#}^\#(f,J_{x_{Q^\#}}H) + \Delta_\junk \lvert J_{x_{Q^\#}} H \rvert \right] \\ 
& \leq C' \left[ \| (f,J_{x_{Q^\#}} H) \|_{S_1 Q^\#} + \Delta_\junk \lvert J_{x_{Q^\#}} H \rvert \right] \\
& \leq C'' \left[  \| H \|_{\X(S_1 Q^\#)} + \Delta_\junk \lvert J_{x_{Q^\#}} H \rvert \right].
\end{align*}
Thus,
\[M_{Q^\#}^\#(0, R_{Q^\#}^\# - J_{x_{Q^\#}} H) \leq C \left[ \| H \|_{\X(S_1 Q^\#)} + \Delta_\junk \lvert J_{x_{Q^\#}} H \rvert \right].\]
Lemma \ref{key_lem1_fin} implies that $\| (0, R_{Q^\#}^\# - J_{x_{Q^\#}} H) \|_{S_0 Q^\#} \leq C \left[ \| H \|_{\X(S_1 Q^\#)} + \Delta_\junk \lvert J_{x_{Q^\#}} H \rvert \right]$, hence
\[
R_{Q^\#} - J_{x_{Q^\#}} H    \in      C    \left[ \| H \|_{\X(S_1 Q^\#)} + \Delta_\junk \lvert J_{x_{Q^\#}}H \rvert  \right] \cdot \sigma(S_0 Q^\#). \]
Since $\partial^\alpha (R_{Q^\#}^\# - J_{x_{Q^\#}} H) \equiv 0$  for $\alpha \in \cA$, \eqref{smallsigma} implies that
\begin{equation}
\label{eqqq1}
\left\lvert \partial^\beta \left[ J_{x_{Q^\#}} H - R^\#_{Q^\#} \right] (x_{Q^\#}) \right\rvert \leq C \cdot \delta_{Q^\#}^{m-n/p - |\beta|} \cdot \left[ \| H \|_{\X(S_1 Q^\#)} + \Delta_\junk \lvert J_{x_{Q^\#}}H \rvert \right]
\end{equation}
for all $\beta \in \cM$.

Hence, by the Sobolev inequality we have
\[\delta_{Q^\#}^{-m} \| H - R^\#_{Q^\#} \|_{L^p(S_1 Q^\#)} \leq C \left[ \| H \|_{\X(S_1 Q^\#)} + \Delta_\junk \lvert J_{x_{Q^\#}}H \rvert \right] \]
This is the desired estimate. (See \eqref{imp_fin}.) This completes the proof of Proposition \ref{key_prop1_fin}.
\end{proof}
We derive a few consequences of Proposition \ref{key_prop1_fin}. We make all the assumptions in the hypothesis of Proposition \ref{key_prop1_fin}. First note that Lemma \ref{pnorm} and \eqref{eqqq1} give
\begin{align}
\label{imp2_fin}
\lvert J_{x_{Q^\#} } H - R_{Q^\#}^\# \rvert   \leq C \cdot \sum_{\beta \in \cM} \lvert  \partial^\beta (J_{x_{Q^\#}} H - R_{Q^\#}^\#)(0)  \rvert  &\leq C' \sum_{\beta \in \cM} \lvert  \partial^\beta (J_{x_{Q^\#}}H - R_{Q^\#}^\#)(x_{Q^\#}) \rvert \notag{} \\
&\leq C'' \left[ \| H \|_{\X(S_1 Q)} + \Delta_\junk \lvert J_{x_{Q^\#}} H \rvert  \right]. \notag{}
\end{align}
(Recall that $\delta_{Q^\#} \leq 1$.) Hence,
\begin{align*}
&\lvert J_{x_{Q^\#} } H \rvert \leq \lvert R_{Q^\#}^\# \rvert + C \cdot \left[ \| H \|_{\X(S_1 Q^\#)} + \Delta_\junk \lvert J_{x_{Q^\#}} H \rvert \right] \implies \\
& \lvert J_{x_{Q^\#} } H \rvert \leq 2 \cdot \lvert R_{Q^\#}^\# \rvert + 2C \cdot \| H \|_{\X(S_1 Q^\#)},
\end{align*}
where we have used that
\[C \cdot \Delta_\junk \leq 1/2.\]
Hence, using the above estimate in \eqref{imp_fin}, we see that
\begin{align*}
\delta_{Q^\#}^{-m} \| H - R_{Q^\#}^\#  \|_{L^p(S_1 Q^\#)} & \leq  C' \left[ \| H \|_{\X(S_1 Q^\#)} + \Delta_\junk \lvert R_{Q^\#}^\# \rvert + \Delta_\junk \| H \|_{\X(S_1 Q^\#)}\right] \\
& \leq C'' \left[ \| H \|_{\X(S_1 Q^\#)} + \Delta_\junk \lvert R^\#_{Q^\#} \rvert \right]. 
\end{align*}
In summary, we have proven the following result.

\begin{lem}
\label{goodversion_lem}
Under the assumptions of Proposition \ref{key_prop1_fin}, we have
\begin{equation}
\label{goodversion}
\delta_{Q^\#}^{-m} \| H - R_{Q^\#}^\#  \|_{L^p(S_1 Q^\#)} \leq C'' \left[ \| H \|_{\X(S_1 Q^\#)} + \Delta_\junk \lvert R^\#_{Q^\#} \rvert \right]. 
\end{equation}
\end{lem}

We will prove one more lemma before returning to the main line of our argument.

\begin{lem}\label{new_fin_lem}
Under the assumptions of Proposition \ref{key_prop1_fin}, we have
\[ 
\lvert R^\# - P \rvert \leq C \cdot \bigl[  \| (f,P) \|_{S_1 Q^\#} +  \Delta_\junk \lvert P \rvert \bigr] 
\]
where $R^\# = R^\#_{Q^\#}(f,P)$ for a keystone cube $Q^\#$.
\end{lem}
\begin{proof}
Let $Q^\#$ be a keystone cube, and denote $R^\# = R^\#_{Q^\#}(f,P)$.

Suppose that $\widetilde{R} \in \cP$ satisfies $\partial^\alpha \widetilde{R} \equiv 0$ for all $\alpha \in \cA$. Recall that \eqref{insigma0} and \eqref{zeroder} imply \eqref{smallsigma}. Put another way, this means that
\[
\lvert \partial^\beta \widetilde{R} (x_{Q^\#}) \rvert \leq C \cdot (\delta_{Q^\#})^{m - n/p - |\beta|}  \cdot \| (0,\widetilde{R}) \|_{S_0 Q^\#} \;\; \mbox{for all} \; \beta \in \cM.
\]
Thus, noting that $\delta_{Q^\#} \leq 1$, we apply Lemma \ref{key_lem1_fin} to deduce that
\begin{align*}
\lvert \widetilde{R} \rvert & \leq C \cdot \sum_{\beta \in \cM} \lvert \partial^\beta \widetilde{R} (x_{Q^\#}) \rvert  \\
& \leq C' \cdot M_{Q^\#}^\#(0,\widetilde{R}). 
\end{align*}
Taking $R = P$ in \eqref{crab2}, we see that the polynomial $R^\# = R^\#_{Q^\#}(f,P)$ satisfies
\begin{equation*}
\left\{
\begin{aligned}
& \partial^\alpha (R^\# - P) \equiv 0 \qquad \mbox{for all} \; \alpha \in \cA \\
& M_{Q^\#}^\#(f,R^\#) \leq  C \cdot \left[ M_{Q^\#}^\#(f,P) + \Delta_\junk \cdot \lvert P \rvert\right].
\end{aligned}
\right.
\end{equation*}
Thus, the $\cA$-derivatives of $\widetilde{R} = R^\# - P$ vanish, so we may apply the previous estimate to give
\begin{align*}
\lvert R^\# - P \rvert &\leq C' \cdot M_{{Q^\#}}^\#(0 , R^\# - P) \\
& \leq C'' \cdot \left[ M_{Q^\#}^\#(f,R^\#) + M_{Q^\#}^\#(f,P) \right] \\
& \leq C''' \cdot \left[ M_{Q^\#}^\#(f,P) + \Delta_\junk \lvert P \rvert \right] \\
& \leq \overline{C} \cdot \left[ \| (f,P) \|_{S_1 {Q^\#}} + \Delta_\junk \lvert P \rvert \right].
\end{align*}
Here, in the last estimate, we use Lemma \ref{key_lem1_fin}.

This completes the proof of Lemma \ref{new_fin_lem}.

\end{proof}

In Section \ref{mk_cubes}, all of the marked cubes are assumed to be $\til{S}$-bit machine cubes, with $\til{S} \leq C \overline{S}$. All the functionals are to be given in short form with parameters $(\Delta_g^C,\Delta_g^{-C} \Delta_\epsilon)$. This concludes the description of the changes to Section \ref{mk_cubes}.

Let $\hQ \subset Q^\circ$ be a dyadic cube. Recall that we say that $\hQ$ is a testing cube if $\hQ$ can be written as the disjoint union of cubes from $\CZ(\cA^-)$ (see Section \ref{sec_testcube}). 

\begin{remk}
\label{test_fp}
Recall that each cube $Q$ in $\CZ(\cA^-)$ satisfies $\delta_Q \geq c \cdot \Delta_0$ with $\Delta_0 = 2^{-\overline{S}}$ for a universal constant $c > 0$, by the Main Technical Results for $\cA^-$. Hence, every testing cube $\hQ$ has $\til{S}$-bit machine points as corners, where $\til{S} \leq C \overline{S}$ for a universal constant $C$.
\end{remk}

Recall, in \eqref{tdefn}, we introduce a parameter $t_G$, which is an integer power of $2$. Furthermore, we assume that 
\begin{equation}
\label{tprop}
t_G = 2^{-\til{S}} \;\mbox{for an integer} \; \til{S} \; \mbox{with} \; 1 \leq \til{S} \leq C \overline{S}.
\end{equation}

In finite-precision, we make one slight change to Lemma \ref{lem_cover}. 
We assume that the constant $a_{\new}$ is picked to satisfy\begin{equation}
\label{aprop} 
a_{\new} = 2^{-\til{S}} \;\mbox{for an integer} \; \til{S} \; \mbox{with} \; 1 \leq \til{S} \leq C \overline{S}.
\end{equation}
To see that this is possible, we examine the proof of Lemma \ref{lem_cover}. We observe that it suffices to choose $a_{\new} = \frac{a \cdot t_G}{512}$, where $a = a(\cA^-)$ arises in the Main Technical Results for $\cA^-$. There it is stated that $a$ is a universal constant and an integer power of $2$. Because $t_G$ satisfies \eqref{tprop}, we conclude that the constant $a_{\new}$ written above satisfies \eqref{aprop}.

This concludes the description of the changes required in Section \ref{sec_testcube}.

\subsection{Testing Functionals}
\label{sec_newfunc}

We continue on to  Section \ref{sec_as}.

Recall that $\Delta_\new$ is a machine number that satisfies \eqref{constants3}. We will make use of the condition
\begin{equation}\label{constants4}
\Delta_\new \leq c(\epsilon,t_G),
\end{equation}
where $c(\epsilon,t_G) \in (0,1)$ is a small constant depending on $m,n,p,t_G,$ and $\epsilon$. We later choose $\epsilon$ and $t_G$ to depend only on $m$, $n$, and $p$ - hence, \eqref{constants4} is consistent with \eqref{constants1}, \eqref{constants2}, and \eqref{constants3}.

We assume we are given a testing cube $\hQ \subset Q^\circ$.

Just as in \eqref{jet1}, for $Q \in \CZ(\cA^-)$ with $Q \subset (1+100t_G)\hQ$, we define the map
\begin{equation}
\label{jet1_fin}
(f,P) \mapsto  R^{\hQ}_Q := \left\{
     \begin{array}{lr}
       P & : \delta_Q \geq t_G \delta_{\hQ}\\
       R^\#_{\mathcal{K}(Q)}(f,P) & : \delta_Q < t_G \delta_{\hQ}
     \end{array}  
   \right.
\end{equation}
for any $(f,P) \in \X(\frac{65}{64} \hQ \cap E) \oplus \cP$. Recall that $S_1 \mathcal{K}(Q) \subset \frac{65}{64}\hQ$, hence this map is well-defined.

\begin{itemize}
\item Recalling the precision with which we compute the maps $R^\#_{Q^\#}$, we see that $(f,P) \mapsto \partial^\alpha \left[ R^\hQ_Q(f,P) \right](0)$ has the form 
\begin{equation}\label{keymaps1}
\sum_{x \in E} \lambda_x f(x) + \sum_\omega \mu_\omega \omega(f) + \sum_{\beta \in \cM} \theta_\beta \cdot \frac{1}{\beta !} \partial^\beta P(0),
\end{equation}
where the possibly nonzero coefficients $\lambda_x, \mu_\omega, \theta_\beta$ are computed with parameters $(\Delta_g^C, \Delta_g^{-C}\Delta_\epsilon)$. The number of possibly nonzero coefficients is at most a universal constant.
\item Recalling the precision with which we compute each $\xi$ in $\Xi(Q,\cA^-)$, we see that $(f,P) \mapsto \xi(f,R^\hQ_Q)$ has the form
\begin{equation}
\label{keymaps2}
\sum_{x \in E} \widetilde{\lambda}_x f(x) + \sum_\omega \widetilde{\mu}_\omega \omega(f) + \sum_{\beta \in \cM} \widetilde{\theta}_\beta \cdot \frac{1}{\beta !} \partial^\beta P(0),
\end{equation}
where the possibly nonzero coefficients $\widetilde{\lambda}_x, \widetilde{\mu}_\omega, \widetilde{\theta}_\beta$ are computed with parameters $(\Delta_g^C, \Delta_g^{-C}\Delta_\epsilon)$. The number of possibly nonzero coefficients is at most a universal constant.
\end{itemize}

We will need to modify the definition of $M_\hQ(f,P)$ in \eqref{i}-\eqref{iv}.

We define $\left[ M_{\hQ}(f,P) \right]^p$ to be the sum of the terms \textbf{(I)}-\textbf{(IV)} (see \eqref{i}-\eqref{iv}) and the sum of the terms

\begin{align}
\label{newterms}
&\text{\bf (V)} = \Delta_\new^{2p} \sum_{x \in \frac{65}{64}\hQ \cap E}  \lvert f(x) - P(x) \rvert^p, \qquad \mbox{and} \\
\notag{}
&\text{\bf (VI)} = \Delta_\new^{p} \lvert P \rvert^p = \Delta_\new^{p} \sum_{\beta \in \cM}  \lvert \partial^\beta P(0) \rvert^p.
\end{align}
Recall that $\Delta_\new$ is a machine number satisfying \eqref{constants3}.

Each of the linear functionals arising in \textbf{(I)}-\textbf{(IV)} (see \eqref{i}-\eqref{iv}) and in \textbf{(V)}-\textbf{(VI)}, can be computed with precision $(\Delta_g^C,\Delta_\epsilon \Delta_g^{-C})$. In Section \ref{supp_data_fin}, we will analyze the work required to compute all these functionals.

As in \eqref{s1} from the infinite-precision text, we define
\[\ooline{\sigma}(\hQ) := \left\{ P \in \cP : M_{\hQ}(0,P) \leq 1 \right\}.\]
We replace the algorithm \textsc{Approximate New Trace Norm} from the infinite-precision text with the finite-precision version below.

\environmentA{Approximate New Trace Norm (Finite-Precision).}

We are given a machine number $t_G > 0$ as in \eqref{tprop}.

We perform one-time work at most $C(t_G) N \log N$ in space $C(t_G)N$, after which we can answer queries as follows.

A query consists of a testing cube $\hQ$. The response to the query $\hQ$ is a list $\mu_1^{\hQ},\ldots,\mu_D^{\hQ}$ of linear functionals on $\cP$ such that 
\begin{equation}
\label{wc12}
c(t_G) \cdot \sum_{i=1}^D \lvert \mu_i^{\hQ}(P) \rvert^p \leq  \left[ M_{\hQ}(0,P) \right]^p \leq  C(t_G) \cdot \left[ \sum_{i=1}^D \lvert \mu_i^{\hQ}(P) \rvert^p  \right].
\end{equation}
The functionals $\mu_1^{\hQ},\cdots,\mu_D^{\hQ}$ have the form 
\[P \mapsto \sum_{\beta \in \cM} \text{coeff}_\beta \frac{1}{\beta!} \partial^\beta P(0)\]
where the coefficients $\text{coeff}_\beta$ are computed with parameters $(\Delta_g^{C}, \Delta_g^{-C} \Delta_\epsilon)$.

We define a quadratic form on $\cP$ by 
\begin{equation}
\label{quadform_fin}
q_{\hQ}(P) := \sum_{i=1}^D \lvert \mu_i^{\hQ}(P) \rvert^2.
\end{equation}
This quadratic form satisfies
\begin{equation}
\label{qform_bd1_fin}
c(t_G) \cdot \left[ M_{\hQ}(0,P) \right]^2 \leq  q_{\hQ}(P) \leq C(t_G) \cdot \left[ M_{\hQ}(0,P) \right]^2.
\end{equation}
In particular,
\begin{equation}
\label{qform_bd2_fin}
\{ q_\hQ \leq c(t_G) \} \subset \ooline{\sigma}(\hQ) \subset \{ q_\hQ \leq C(t_G)\}.
\end{equation}

The quadratic form $q_{\hQ}$ is given in the form
\[q_{\hQ}(P) = \sum_{\alpha,\beta \in \cM} q_{\alpha \beta} \cdot \frac{1}{\alpha!}  \partial^\alpha P(0) \cdot \frac{1}{\beta!} \partial^\beta P(0).\]
The $q_{\alpha \beta}$ form a symmetric matrix. For each $\alpha,\beta \in \cM$ we compute the number $q_{\alpha \beta}$ with parameters $(\Delta_g^C, \Delta_g^{-C} \Delta_\epsilon)$. 

The work required to answer a query is at most $C(t_G) \log N$.

Here, $c(t_g) > 0$ and $C(t_g) \geq 1$ are constants depending on $t_G$, $m$, $n$, and $p$.

\begin{proof}[\underline{Explanation}]  

We use the finite-precision versions of the algorithms \textsc{Approximate Old Trace Norm}, \textsc{Compute Norms From Marked Cuboids}, and \textsc{Compress Norms}. The explanation otherwise remains the same. We describe a few minor changes that are needed.

Given a keystone cube $Q^\#$, the polynomial map $P \mapsto R^\#_{Q^\#}(0,P)$ in \eqref{polymap} is given with parameters $(\Delta_g^C,\Delta_g^{-C}\Delta_\epsilon)$; indeed, the functionals $\overline{\lambda}_{(Q^\#,\beta)}$ are given with such parameters, as stated in the finite-precision version of the algorithm \textsc{Make New Assists and Assign Keystone Jets}.

We will apply a marking procedure, just as in infinite-precision. We provide details below.

\begin{itemize}
\item For each $Q \in \CZ_{\main}(\cA^-)$ and $1 \leq i \leq D$, we mark the cube $Q$ with the functional $P \mapsto \xi_{(Q,i)}(P)$, defined by
\[
\xi_{(Q,i)}(P) = \xi_i^Q(R^\#_{\mathcal{K}(Q)}(0,P)).
\]
Note that each functional $ R \mapsto \xi_i^Q(R)$ is given with parameters $(\Delta_g^C,\Delta_g^{-C}\Delta_\epsilon)$; see the statement of the algorithm \textsc{Approximate Old Trace Norm} (finite-precision). Also, the polynomial map $P \mapsto R^\#_{\mathcal{K}(Q)}(0,P)$ is given with parameters $(\Delta_g^C,\Delta_g^{-C}\Delta_\epsilon)$. We can stably compute the composition of a linear functional on $\R^D$ with a linear map on $\R^D$, hence we can compute each $\xi_{(Q,i)}$ with parameters $(\Delta_g^C,\Delta_g^{-C}\Delta_\epsilon)$, for a possibly larger constant $C$.

\item For each $(Q',Q'') \in \BD(\cA^-)$ and $\beta \in \cM$, we mark  the cube $Q'$ with the functional $P \mapsto \xi_{(Q',Q'',\beta)}(P)$, defined by
\[
\xi_{(Q',Q'',\beta)}(P) =  \delta_{Q'}^{n/p - m+|\beta|} \cdot \partial^\beta \left\{ R^\#_{\mathcal{K}(Q')}(0,P) - R^\#_{\mathcal{K}(Q'')}(0,P) \right\}(x_{Q'}).
\] This functional is given with parameters $(\Delta_g^C,\Delta_g^{-C}\Delta_\epsilon)$ because the polynomial maps $P \mapsto R^\#_{\mathcal{K}(Q')}(0,P)$ and $P \mapsto R^\#_{\mathcal{K}(Q'')}(0,P)$ are given with the same parameters, because $\delta_{Q'}$ is a machine number with $c \cdot \Delta_0 \leq \delta_{Q'} \leq 1$ (since $Q' \in \CZ(\cA^-)$, this is a consequence of the Main Technical Results for $\cA^-$), and because $x_{Q'}$ is a $\til{S}$-bit machine point with $\til{S} \leq C \overline{S}$.
\end{itemize}

Each of the functionals $\xi$ that is associated to a marked cube has the form 
\[
P \mapsto \xi(P) = \sum_{\beta \in \cM} \mbox{coeff}_\beta \cdot \partial^\beta P(0),
\]
where the coefficients $\text{coeff}_\beta$ are specified with parameters $(\Delta_g^C,\Delta_g^{-C}\Delta_\epsilon)$. We perform one-time work for the algorithm \textsc{Compute Norms From Marked Cubes} (finite-precision) on the marked cubes described in the above bullet points. (See \textbf{Modification 6} in Section \ref{adc_fin}.) All of the marked cubes belong to $\CZ(\cA^-)$. By assumption, all cubes in $\CZ(\cA^-)$ have sidelength in the interval $\left[ c \cdot \Delta_0, 1 \right]$, where $\Delta_0 = 2^{- \overline{S}}$. Hence, the marked cubes have $\til{S}$-bit machine points as corners, where $\til{S} \leq C \overline{S}$. Thus, the finite-precision version of the algorithm applies.

Next, we explain the query work for the algorithm \textsc{Approximate New Trace Norm}.

We are given a testing cube $\hQ$. As in infinite-precision, we partition $(1+t_G)\hQ$ into dyadic cubes $Q_1,\cdots,Q_L \subset \R^n$ such that $\delta_{Q_\ell} = (t_G/4)\delta_\hQ$.  Hence, $\delta_{Q_\ell} \geq 2^{- C \overline{S}}$ for a universal constant $C$, for each $\ell=1,\cdots,L$, thanks to \eqref{tprop} and Remark \ref{test_fp}. Consequently, each $Q_\ell$ is a $\til{S}$-bit machine cube with $\til{S} \leq C \overline{S}$. Moreover, note that $L \leq C(t_G)$.

Thus, we can apply the query algorithm in \textsc{Compute Norms From Marked Cuboids} (finite-precision)  with $\Delta = \Delta_g$ for each cube $Q_\ell$ ($1 \leq \ell \leq L$). We refer the reader to \textbf{Modification 6} in Section \ref{adc_fin} for the statement of the relevant algorithm. Thus, with work at most $C(t_G) \log N$ we can compute linear functionals $\mu_1^{Q_\ell},\cdots,\mu_{D}^{Q_\ell}$ such that 
\begin{align}
\label{fin_1}
c \sum_{k=1}^D \lvert \mu_k^{Q_\ell}(P) \rvert^p &\leq \sum_{\substack{Q \in \CZ(\cA^-) \\ \text{linear functional} \; \xi}} \bigl\{ \lvert \xi(P) \rvert^p : Q \subset Q_\ell, \; Q \; \mbox{marked with} \; \xi \bigr\} + \Delta_g^{p/2} \lvert P \rvert^p \\
& \notag{} \leq C \left[  \sum_{k=1}^D \lvert \mu_k^{Q_\ell}(P) \rvert^p + \Delta_g^{p/2} \lvert P \rvert^p \right],
\end{align}
where each functional $\mu_k^{Q_\ell}$  in \eqref{fin_1} is given with parameters $(\Delta_g^C,\Delta_g^{-C} \Delta_\epsilon)$. Here, we use the estimate $\Delta_g^p  \Delta_0^{-C} \log(\Delta_g^{-1}) \leq \Delta_g^{p/2}$, which follows from the assumption $\Delta_g \ll \Delta_0$. Recall that $\lv P \rv$ is defined to be the $\ell^p$ norm of the vector $(\partial^\alpha P(0))_{\alpha \in \cM}$.

We sum \eqref{fin_1} from $\ell = 1,\cdots,L$. The sum of the junk terms is equal to $L \Delta_g^{p/2} \lvert P \rvert^p \leq C(t_G) \Delta_g^{p/2} \lvert P \rvert^p$. Hence, as in infinite-precision, in analogy with \eqref{e763}, we have
\begin{align}
\label{e763_fin}
c(t_G) \sum_{\ell=1}^L \sum_{k=1}^D \lvert \mu_k^{Q_\ell}(P) \rvert^p &\leq \left[ \mathfrak{S}_1 + \mathfrak{S}_2 \right] + \Delta_g^{p/2} \lvert P \rvert^p \\
\notag{}
& \leq C(t_G) \left[ \sum_{\ell=1}^L \sum_{k=1}^D \lvert \mu_k^{Q_\ell}(P) \rvert^p + \Delta_g^{p/2} \lvert P \rvert^p \right].
\end{align}
For the definition of the terms $\mathfrak{S}_1$ and $\mathfrak{S}_2$, see \eqref{e763}.

We compute all the functionals in $\mathbf{(F_1)}$-$\mathbf{(F_6)}$ (see the text following \eqref{e763}), as in infinite-precision, by looping over relevant cubes and listing the relevant functionals. We have described in the above how to compute  the functionals in $\mathbf{(F_1)}$, $\mathbf{(F_2)}$, $\mathbf{(F_3)}$, and $\mathbf{(F_4)}$; each such functional is given with parameters $(\Delta_g^{C}, \Delta_g^{-C} \Delta_\epsilon)$. Additionally, note that the maps $P \mapsto R_Q^\hQ(0,P)$ are given with parameters $(\Delta_g^{C},\Delta_g^{-C}\Delta_\epsilon)$; see \eqref{jet1_fin}. Hence, we can compute each functional in $\mathbf{(F_5)}$ and $\mathbf{(F_6)}$ with parameters $(\Delta_g^{C}, \Delta_g^{-C} \Delta_\epsilon)$.

Therefore, each functional listed in $\mathbf{(F_1)}$-$\mathbf{(F_6)}$ has the form 
\[P \mapsto \sum_{\beta \in \cM} d_\beta \cdot \frac{1}{\beta!} \partial^\beta P(0),\]
where the numbers $d_\beta$ are computed with parameters $(\Delta_g^{C}, \Delta_g^{-C} \Delta_\epsilon)$. 

In addition, we compute the functionals
\begin{itemize}
\item[$\mathbf{(F_7)}$] $\boxed{ \lambda_\beta(P) := \Delta_\new \cdot \partial^\beta P(0)}$ \qquad for $\beta \in \cM$.
\end{itemize}
Each functional in $\mathbf{(F_7)}$ is given with parameters $(\Delta_g^{C}, \Delta_g^{-C} \Delta_\epsilon)$.

We define $\left[ X(P) \right]^p$ to be the sum of the $p$-th powers of the functionals  in $\mathbf{(F_1)}$-$\mathbf{(F_7)}$. 

Note that the sum of the $p$th powers of the functionals in $\mathbf{(F_7)}$ is equal to the term \textbf{(VI)} (see \eqref{newterms}). 

Hence, we have 
\[ 
\left[ X(P) \right]^p = \left[ X_\old(P) \right]^p + \text{\textbf{(VI)}},
\] 
where $\left[ X_\old(P) \right]^p$ is the sum of the $p$-th powers of the functionals in $\mathbf{(F_1)}$-$\mathbf{(F_6)}$. 

We have
\begin{align}
\label{Xeq_fin}
c(t_G) \cdot \left[ X_\old(P) \right]^p &\leq \left[ \text{sum of terms \textbf{(I)}-\textbf{(IV)} with} \; f \equiv 0 \right] + \Delta_g^{p/2} \lv P \rv^p \\
\notag{}
&\leq  C(t_G) \cdot \left[ \left[ X_\old(P) \right]^p + \Delta_g^{p/2} \lv P \rv^p \right].
\end{align}
To obtain the above estimate, we reason as in the paragraph containing \eqref{Xeq} and the paragraph following \eqref{Xeq}, making sure to use the estimate \eqref{e763_fin}  in place of \eqref{e763}. (Recall that $X_\old$ corresponds to $X$ in our earlier notation.)

We examine the term \textbf{(V)} (see \eqref{newterms}), which arises in the definition of $M_{\hQ}(0,P)$. When $f \equiv 0$, we have
\[
\text{\textbf{(V)}} = \Delta_\new^{2p} \sum_{x \in \frac{65}{64}\hQ \cap E} \lv P(x) \rv^p \leq C \Delta_\new^{2p} N \sum_{\beta \in \cM} \lv \partial^\beta P(0) \rv^p \leq \Delta_\new^p \sum_{\beta \in \cM} \lv \partial^\beta P(0) \rv^p.
\]
Here, we use the estimates $N \leq \Delta_0^{-n} \leq \Delta_\new^{-p/10}$ (see \eqref{Nbound} and \eqref{constants3}) and $C \Delta_\new^{19p/10} \leq \Delta_\new^p$ (see \eqref{constants3}). Hence, when $f=0$, the term \textbf{(V)} is bounded by \textbf{(VI)}.  Thus, up to constant factors, $\bigl[ M_{\hQ}(0,P) \bigr]^p$ is equivalent to the sum of the terms \textbf{(I)}-\textbf{(IV)} and \textbf{(VI)} (with $f \equiv 0$). 

Therefore, by adding the term $\text{\textbf{(VI)}} + \Delta_g^{p/2} \lvert P \rvert^p$ to the chain of inequalities \eqref{Xeq_fin}, we learn that
\begin{align*}
c(t_G) \cdot \left\{ \bigl[ X_{\old}(P) \bigr]^p + \text{\textbf{(VI)}} + \Delta_g^{p/2} \lvert P \rvert^p \right\} & \leq \left[ M_{\hQ}(0,P) \right]^p + \Delta_g^{p/2} \lv P \rv^p \\
&\leq C(t_G) \cdot \left\{ \left[ X_{\old}(P) \right]^p + \text{\textbf{(VI)}}  + \Delta_g^{p/2} \lv P \rv^p \right\}.
\end{align*}
Note that the middle term above is comparable to $\left[ M_\hQ(0,P)\right]^p$, since 
\[
\left[ M_\hQ(0,P)\right]^p \geq \text{\textbf{(VI)}} = \Delta_\new^{p} \lvert P \rvert^p \geq \Delta_g^{p/2} \cdot \lv P \rv^p.
\]
Here, we use that $\Delta_g \leq \Delta^2_\new$ (see \eqref{constants3}). Since $\left[ X(P) \right]^p = \left[ X_\old(P) \right]^p + \text{\textbf{(VI)}}$, we conclude that
\[
c(t_G) \cdot \left\{ \bigl[ X(P) \bigr]^p + \Delta_g^{p/2}  \lv P \rv^p\right\} \leq \left[ M_\hQ(0,P)\right]^p  \leq C(t_G) \cdot  \left\{ \bigl[ X(P) \bigr]^p + \Delta_g^{p/2} \lv P \rv^p \right\}.
\]

Recall that $\bigl[ X(P) \bigr]^p$ is the sum of the $p$th powers of the functionals in $\mathbf{(F_1)}$-$\mathbf{(F_7)}$. Processing the functionals in $\mathbf{(F_1)}$-$\mathbf{(F_7)}$ using \textsc{Compress Norms} (finite-precision), we compute functionals $\mu_1^\hQ,\cdots,\mu_D^\hQ$ on $\cP$ such that
\[c \cdot \sum_{i=1}^D \lvert \mu_i^\hQ(P) \rvert^p \leq  \bigl[ X(P) \bigr]^p + \Delta_g^{p/2} \lvert P \rvert^p \leq C \cdot \sum_{i=1}^D \lvert \mu_i^\hQ(P) \rvert^p.\]
The functionals $\mu_i^\hQ$ are given with parameters $(\Delta_g^{C},\Delta_g^{-C} \Delta_\epsilon)$. 

The previous two estimates establish \eqref{wc12}. 

Moreover, properties \eqref{qform_bd1_fin} and \eqref{qform_bd2_fin} are immediate from the definition of $q_{\hQ}$ in \eqref{quadform_fin} and the equivalence of the $\ell_p$ and $\ell_2$ norms on a finite-dimensional space. Each of the $\mu_i^\hQ$ is given with parameters $(\Delta_g^{C},\Delta_g^{-C} \Delta_\epsilon)$, hence the coefficients $q_{\alpha \beta}$ of the quadratic form $q_{\hQ}$ can be computed with parameters $(\Delta_g^{C},\Delta_g^{-C} \Delta_\epsilon)$ (for a possibly larger constant $C$).

This concludes the explanation of the query algorithm. It is easy to check that the query work at most $C(t_G) \log N$.

\end{proof}

\subsection{Supporting Data}
\label{supp_data_fin}

We assume we are given a testing cube $\hQ \subset Q^\circ$.

We explain the main modifications to Section \ref{sec_computingstuff} needed here.

\begin{itemize}
\item \textbf{Modification 1:} As part of the supporting data for $\hQ$, we include a list of all the points $x \in \frac{65}{64}\hQ \cap E$, in addition to all the other data described in Section \ref{sec_computingstuff}. We call this the \emph{modified supporting data} for $\hQ$.
\end{itemize} 

The list $\Omega(\hQ)$ of the new assist functionals is defined as in \eqref{test_assists}.

\begin{itemize}
\item \textbf{Modification 2:} The algorithm \textsc{Compute New Assists} operates as follows. Given a testing cube $\hQ$, and given the supporting data for $\hQ$, we compute a list of all the functionals in $\Omega(\hQ)$. We compute a short form of each $\omega \in \Omega(\hQ)$ with parameters $(\Delta_g^C, \Delta_g^{-C} \Delta_\epsilon )$.

\item \textbf{Modification 3:} We make only minor changes to the algorithm \textsc{Compute Supporting Map}. The linear maps $R_Q^\hQ$ are to be computed in short form in terms of the assists $\Omega(\hQ)$ with precision $(\Delta_g^C,\Delta_g^{-C}\Delta_\epsilon)$. The explanation of the algorithm is unchanged.

\item \textbf{Modification 4:} We replace the algorithm \textsc{Compute New Assisted Functionals} with the finite-precision version described below. 
\end{itemize}

\environmentA{Algorithm: Compute New Assisted Functionals (Finite-Precision).}

Given a testing cube $\hQ$ and its modified supporting data, we define
\begin{equation}
\label{work2_fin}
\Work_2^\fin(\hQ) := \Work_2(\hQ) + C(t_G) \cdot \#\left(\frac{65}{64}\hQ \cap E \right)
\end{equation}
and
\begin{equation}
\label{space2_fin}
\Space_2^\fin(\hQ) := \Space_2(\hQ) + C(t_G) \cdot \#\left(\frac{65}{64}\hQ \cap E\right),
\end{equation}
where $\Work_2(\hQ)$ and $\Space_2(\hQ)$ are defined in \eqref{work2} and \eqref{space2}, respectively.

We compute a list $\Xi(\hQ)$ of functionals on $\X(E) \oplus \cP$, such that
\[
\left[ M_\hQ(f,P) \right]^p = \sum_{\xi \in \Xi(\hQ)} \lvert \xi(f,P) \rvert^p.
\]
Each functional $\xi$ in $\Xi(\hQ)$ is given in short form in terms of assists $\Omega(\hQ)$ with parameters $(\Delta_g^C,\Delta_g^{-C} \Delta_\epsilon)$.

This computation requires work at most $\Work_2^\fin(\hQ)$ in space $\Space_2^\fin(\hQ)$.

\begin{proof}[\underline{Explanation}]
We include in the list $\Xi(\hQ)$ all the same functionals as before, namely, the functionals in \eqref{i1}-\eqref{iv1}, as well as a few additional functionals  described below. Each ``assisted functional'' in \eqref{i1}-\eqref{iv1} is given in short form with parameters $(\Delta_g^C,\Delta_g^{-C} \Delta_\epsilon)$ in terms of the assists $\Omega(\hQ)$. Indeed, all the functionals $\xi$ and maps $R_Q^\hQ$, $R_{Q'}^\hQ$, $R_{Q''}^\hQ$, $R_{Q_{\spec}}^\hQ$, which are relevant to \eqref{i1}-\eqref{iv1}, are given in short form with parameters $(\Delta_g^C,\Delta_g^{-C}\Delta_\epsilon)$. See \eqref{keymaps1} and \eqref{keymaps2}.

Hence, we can compute all the ``assisted functionals'' in \eqref{i1}-\eqref{iv1} with parameters $(\Delta_g^C,\Delta_g^{-C}\Delta_\epsilon)$, as claimed. 

In addition, we include in the list $\Xi(\hQ)$ the additional functionals
\[\lambda_x(f,P) := \Delta_\new^2 \cdot (f(x) - P(x)) \quad \mbox{for each} \; x \in \frac{65}{64}\hQ \cap E\]
and
\[\lambda_\beta(f,P) := \Delta_\new \cdot \partial^\beta P(0) \quad \mbox{for each} \; \beta \in \cM.\]
There are the new functionals needed in finite-precision. That completes the definition of $\widehat{\Xi}(\hQ)$. Note that $\Delta_\new \leq 1 \leq \Delta_g^{-1}$. Hence, each functional $\lambda_x$ and $\lambda_\beta$ can be expressed in short form (without assists) using coefficients that are bounded in magnitude by $\Delta_g^{-C}$. Moreover, each coefficient can be computed to precision $\Delta_\epsilon$, which is within the precision  available to our computer. (Recall: The computer works with precision $\Delta_{\min} \ll \Delta_\epsilon$.) Hence, the functionals $\lambda_x$ and $\lambda_\beta$ can be computed in short form with parameters $(\Delta_g^C,\Delta_\epsilon)$ (without assists). 

The sum of $\lvert \xi(f,P) \rvert^p$ over all $\xi$ in $\Xi(\hQ)$ is equal to $\left[M_\hQ(f,P)\right]^p$, by definition.

The additional term $\#(\frac{65}{64}\hQ \cap E)$ in \eqref{work2_fin} and \eqref{space2_fin} accounts for the additional work and space, respectively, needed to compute the functionals $\lambda_x$, $\lambda_\beta$.

This completes the explanation of the algorithm.
\end{proof}

Given a testing cube $\hQ$, the covering cubes $\cov(\hQ) \subset \CZ(\cA^-)$ are defined as in \eqref{scrI}, namely
\[
\cov(\hQ) := \bigl\{ Q \in \CZ(\cA^-) : Q \subset (1+t_G)\hQ\bigr\}.
\]
We introduce a family of cutoff functions $\theta_Q^\hQ$ (for $Q \in \cov(\hQ)$) that satisfy \eqref{pou701}-\eqref{pou703}. 

The finite-precision version of the algorithm \textsc{Compute POU} is as follows.

\environmentA{Compute POU (finite-precision).}

After one-time work at most $C N \log N$ in space $CN$, we can answer queries as follows.

A query consists of a testing cube $\hQ$ and an $S$-bit machine point $\underline{x} \in Q^\circ$. 

Notice that $\hQ$ has $\til{S}$-bit machine points as corners, with $\til{S} \leq C \overline{S}$ for a universal constant $C$, so we can safely process $\hQ$ on our finite-precision computer. (See Remark \ref{test_fp}.)

We respond to the query with a list of all the cubes $Q_1,\cdots,Q_L \in \cov(\hQ)$ (with $Q_1,\cdots,Q_L$ all distinct) such that $\underline{x} \in \frac{65}{64} Q_\ell$. Futhermore, we compute the numbers $\frac{1}{\alpha !} \partial^\alpha J_{\underline{x}} \theta_{Q_\ell}^{\hQ}(0)$ (for all $\ell=1,\cdots,L$ and $\alpha \in \cM$) with parameters $(\Delta_g^C,\Delta_g^{-C} \Delta_\epsilon)$.

To answer a query requires work and storage at most $C \log N$.
\begin{proof}[\underline{Explanation}]
The explanation is just as in the infinite-precision case. We refer the reader to the explanation given just after equation \eqref{pou703} for more details.  
\end{proof}

The definitions and conditions in \eqref{test_eo_aux}-\eqref{extopdefn} are unchanged; as before, the Main Technical Results for $\cA^-$ (finite-precision) yield the algorithm \textsc{Compute New Extension Operator}, with the following modifications:

\begin{itemize}
\item \textbf{Modification 5:} We assume that $\underline{x} \in Q^\circ$ is an $S$-bit machine point. We compute the functionals  $(f,P) \mapsto \partial^\beta ( J_{\underline{x}} T_{\hQ}(f,P))(0)$ which have the form
\[\sum_{\ell = 1}^L \gamma_{\ell} \cdot \omega_\ell(f) + \sum_{j=1}^J \lambda_j  \cdot f(x_j) + \sum_{\gamma \in \cM} \theta_\gamma \cdot \frac{1}{\gamma!} \partial^\gamma P(0), \]
where each $\omega_\ell$ is in $\Omega(\hQ)$ and each $x_j$ is in $E \cap \frac{65}{64}\hQ$; each real number $\gamma_\ell$, $\lambda_j$, and $\theta_\gamma$, is computed with parameters $(\Delta_g^C,\Delta_g^{-C} \Delta_\epsilon)$; and $L + J + \#(\cM) \leq C$, for a universal constant $C$.
\end{itemize}

\section{Inequalities for Testing Functionals}

Let $a_\new = a_\new(t_G)$ be the constant from Lemma \ref{lem_cover}. We recall  that $a_{\new} = 2^{-\til{S}}$ for an integer $\til{S}$ with $1 \leq \til{S} \leq C \overline{S}$. (See \eqref{aprop}.)

First, in Proposition \ref{prop_bddextop}, we state and prove some properties of the extension operator $(f,P) \mapsto T_\hQ(f,P)$ defined in \eqref{extopdefn}. The assertion and proof of Proposition \ref{prop_bddextop} are unchanged in the current setting.

Next, we prove a few estimates to show that the testing functional $M_{\hQ}$ defined before well-approximates the trace seminorm near the testing cube $\hQ$. Such estimates were stated before in Proposition \ref{prop_normappx} (the conditional/unconditional inequalities). In the present setting, the statement and proof of the corresponding estimates will need to be modified. The next result contains the relevant estimates.

\begin{prop}\label{prop_normappx_fin}
Let $\hQ$ be a testing cube, and let $(f,P) \in \X(\frac{65}{64}\hQ \cap E) \oplus \cP$. Then the following estimates hold.
\begin{description}
\item[Unconditional inequality] $\| (f,P) \|_{(1+a_\new)\hQ} \leq C(t_G) \cdot M_{\hQ}(f,P)$.
\item[Conditional inequality] If $3\hQ$ is tagged with $(\cA,\epsilon)$, then 
\[ M_{\hQ}(f,P) \leq C(t_G) \cdot (1/\epsilon) \cdot \left[ \|(f,P) \|_{\frac{65}{64}\hQ}  + \Delta_\new \cdot \lvert P \rvert \right].\]
\end{description}
\end{prop}

The unconditional inequality is a direct consequence of Proposition \ref{prop_bddextop} as in the infinite-precision proof of Proposition \ref{prop_normappx}. We now prove the conditional inequality.

The finite-precision version of Lemma \ref{lem_simple} reads as follows.
\begin{lem}\label{lem_simple_fin}
Suppose that the testing cube $\hQ$ is $\eta$-simple for some $\eta \geq t_G$. \\
Then 
\[M_{\hQ}(f,P) \leq C(t_G) \cdot \left[ \| (f,P) \|_{\frac{65}{64}\hQ} + \Delta_\new \cdot \lvert P \rvert \right],\]
where $C(t_G)$ depends only on $m$, $n$, $p$, and $t_G$.
\end{lem}
\begin{proof}
The proof requires minor modifications. If $\hQ$ is $\eta$-simple ($\eta \geq t_G$), then the terms \textbf{(II)}, \textbf{(III)}, \textbf{(IV)} vanish, as explained in the outset of the proof of Lemma \ref{lem_simple}. This leaves us with the original term \textbf{(I)}, and the new terms \textbf{(V)}, and \textbf{(VI)}. Recall, the definition of \textbf{(I)} is given in Section \ref{sec_as}. and the definitions of \textbf{(V),(VI)} are given in Section \ref{sec_newfunc}.

We consider an arbitrary summand $\left[ M_{(Q,\cA^-)}(f,P) \right]^p$ in the term \textbf{(I)} (see \eqref{i}). From \eqref{n_appx_fin}, we have
\[
\left[ M_{(Q,\cA^-)}(f,P) \right]^p \leq C \cdot \left[ \| (f,P) \|_{\frac{65}{64}Q}^p + \Delta_\junk^p \lvert P \rvert^p \right],
\]
hence
\[ 
\left[ M_{(Q,\cA^-)}(f,P) \right]^p \leq  C(t_G) \cdot \left[ \| (f,P) \|^p_{\frac{65}{64}\hQ} + \Delta_\junk^p  \lvert P \rvert^p \right],
\]
where we have used that $\frac{65}{64} Q \subset \frac{65}{64}\hQ$ and that $\hQ$ is $\eta$-simple with $\eta \geq t_G$; for more details, see the proof of Lemma \ref{lem_simple}. Since the number of cubes $Q$ relevant to term \textbf{(I)} is bounded by $C(t_G)$, we conclude that 
\[
\text{ \bf{(I)}} \leq C(t_G) \cdot \left[ \| (f,P) \|_{\frac{65}{64}\hQ}^p + \Delta_\junk^p  \lvert P \rvert^p \right].
\]

Note that term \textbf{(VI)} is equal to $\Delta_\new^p \lvert P \rvert^p$.

It remains to estimate term \textbf{(V)}. Recall that $\#(\frac{65}{64}\hQ \cap E) \leq N \leq \Delta_0^{-n}$ (see \eqref{Nbound}). Hence,
\begin{align}
\label{wc2a}
\text{\bf (V)} =  \Delta_\new^{2p} \sum_{x \in \frac{65}{64}\hQ \cap E} \lvert f(x) - P(x) \rvert^p & \leq \Delta_\new^{2p} \Delta_0^{-C} \| (f,P) \|_{\frac{65}{64}\hQ}^p \\
\label{wc2b}
&\leq \| (f,P) \|_{\frac{65}{64}\hQ}^p.
\end{align}
We explain below the previous two estimates.

We deduce the estimate in \eqref{wc2a} by picking a function $\widetilde{F}$ that satisfies $\| \widetilde{F} \|_{\X(\frac{65}{64}\hQ)}^p + \delta_\hQ^{-mp} \| \widetilde{F} - P \|_{L^p(\frac{65}{64} \hQ)} \leq 2 \cdot \| (f,P) \|_{\frac{65}{64}\hQ}^p$ and $\widetilde{F} = f$ on $E \cap \frac{65}{64}\hQ$. Then, for each $x \in E \cap \frac{65}{64}\hQ$, we apply the estimate \eqref{s_ineq0} from Lemma \ref{si2} to the function $F = \widetilde{F} - P$. Hence, we have
\[
\lvert f(x) - P(x) \rvert = \lv \widetilde{F}(x) - P(x) \rv \leq C \cdot \delta_{\hQ}^{-n/p} \|\widetilde{F} - P \|_{L^p(\frac{65}{64} \hQ)} + \delta_{\hQ}^{m-n/p} \| \widetilde{F} \|_{\X(\frac{65}{64} \hQ)}.
\]
Recall that $\delta_{\hQ} \geq c \Delta_0$, since $\hQ$ is a testing cube. Hence, we have $\lvert f(x) - P(x) \rvert \leq C \Delta_0^{-C} \lvert (f,P) \rvert_{\frac{65}{64}\hQ}$. Summing over $x \in E \cap \frac{65}{64}\hQ$ and using the fact that $\#(E \cap \frac{65}{64}\hQ) \leq \Delta_0^{-C}$, we obtain the stated estimate. 

We deduce the estimate in \eqref{wc2b} from  the estimate $\Delta_\new \leq \Delta_0^{C/(2p)}$; see \eqref{constants3}.

Therefore,
\[
\left[M_\hQ(f,P)\right]^p = \text{\bf (I)} + \text{\bf (V)} + \text{\bf (VI)} \leq C(t_G) \cdot \left[ \| (f,P) \|_{\frac{65}{64}\hQ}^p + \Delta_\new^p \lvert P \rvert^p + \Delta_\junk^p  \lvert P \rvert^p \right].
\]
Since $\Delta_\junk \leq \Delta_\new$ (see \eqref{constants3}), the above estimate implies the conclusion of Lemma \ref{lem_simple_fin}. 
\end{proof}

Lemma \ref{lem_simple_fin} implies the conditional inequality in the $\eta$-simple case.

So we may assume that $\hQ$ is not $\eta$-simple as in \eqref{notsimple}.

Both Proposition \ref{prop_A} and Proposition \ref{sob_prop} hold in the present setting, without change. The proofs are as before.

We now prove the conditional inequality. We describe how the estimates from before will need to be changed in the present setting.

On the right-hand side of \eqref{e772} we add the terms $\text{\bf(V)}$ and $\text{\bf(VI)}$. Note that
\begin{align*}
\text{\textbf{(V)}} + \text{\textbf{(VI)}} &=\Delta_\new^{2p} \sum_{x \in \frac{65}{64}\hQ \cap E} \lvert f(x) - P(x) \rvert^p + \Delta_\new^{p}  \lvert P \rvert^p \\
& \leq \| (f,P) \|_{\frac{65}{64}\hQ}^p + \Delta_\new^{p} \lvert P \rvert^p \quad (\text{thanks to \eqref{wc2a}}) \\
& \leq \left[ \text{RHS of the conditional inequality} \right]^p.
\end{align*}

Therefore, the extra terms in \eqref{e772} don't hurt.

Instead of \eqref{e772a}, our inductive assumption now states that
\begin{align*}
M_{(Q,\cA^-)}(f, R_Q^\hQ) & \leq C \cdot \left[ \| (f , R_Q^\hQ) \|_{\frac{65}{64}Q} + \Delta_\junk \lvert R_Q^\hQ \rvert \right] \\
&\leq C \cdot  \left[\|H \|_{\X(\frac{65}{64}Q)} + \delta_Q^{-m} \|H - R_Q^\hQ \|_{L^p(\frac{65}{64}Q)} +  \Delta_\junk \lvert R_Q^\hQ \rvert  \right].
\end{align*}
Hence, in place of \eqref{e773}, we now have
\begin{align}
\label{eq1_fin}
\left[ M_{\hQ}(f,P) \right]^p \leq C(t_G) \cdot \biggl[ &\| H \|_{\X(\frac{65}{64}\hQ)}^p + \delta_\hQ^{-mp} \| H - P \|_{L^p(\frac{65}{64}\hQ)}^p \\
& + \sum_{Q \subset (1+100t_G)\hQ} \left[ \| H \|_{\X(\frac{65}{64}Q)}^p + \delta_Q^{-mp} \| H - R_Q^\hQ \|_{L^p(\frac{65}{64}Q)}^p \right] \notag{} \\
& + \Delta_\junk^p \sum_{Q \subset (1+100t_G)\hQ}  \lvert R_Q^\hQ \rvert^p \;\; (=: \mathfrak{S})  \notag{}\\
& + \bigl[ \text{RHS of conditional inequality} \bigr]^p \biggr]. \notag{}
\end{align}
The third and fourth lines contain new terms not present in the original estimate.

We will now estimate the extra term $\mathfrak{S}$ in the third line of \eqref{eq1_fin}

We write $\mathfrak{S} = \mathfrak{S}_1 + \mathfrak{S}_2$, with
\begin{align*}
\mathfrak{S}_1 &=  \Delta_\junk^p \sum_{\substack{Q \subset (1+100t_G)\hQ \\\delta_Q \geq t_G \delta_\hQ}} \lvert R_Q^\hQ \rvert^p, \\
\mathfrak{S}_2 &=   \Delta_\junk^p \sum_{\substack{Q \subset (1+100t_G)\hQ \\ \delta_Q < t_G \delta_\hQ}} \lvert R_Q^\hQ \rvert^p
\end{align*}

We estimate the term $\mathfrak{S}_1$. The number of cubes in $\CZ(\cA^-)$ is at most $\Delta_0^{-C}$. Moreover, by definition, $R_Q^\hQ  = P$ for each $Q \in \CZ(\cA^-)$ relevant to the $\mathfrak{S}_1$ (see \eqref{jet1_fin}).  Hence,
\begin{align}\label{temp11}
\mathfrak{S}_1 & \leq \Delta_\junk^p \cdot \Delta_0^{-C} \cdot  \lvert P \rvert^p \\
& \leq \Delta_\new^p \cdot \lvert P \rvert^p \leq \bigl[ \text{RHS of conditional inequality} \bigr]^p. \notag{}
\end{align}
(Here, we use that $\Delta_\junk \leq \Delta_\new^2 \leq \Delta_\new \cdot \Delta_0^{C/p}$; see \eqref{constants3}.)

We will now estimate the term $\mathfrak{S}_2$. Let $Q \in \CZ(\cA^-)$ satisfy $Q \subset (1+100t_G) \hQ$ and $\delta_Q < t_G \delta_\hQ$. Note that the keystone cube $Q^\# = \mathcal{K}(Q)$ associated to $Q$ satisfies $S_1 Q^\# \subset \frac{65}{64}\hQ$. (See \eqref{edp2} in Proposition \ref{sob_prop}.) Furthermore, by definition \eqref{jet1} we have $R_Q^\hQ = R^\#_{Q^\#}$. Moreover, since the number of cubes in $\CZ(\cA^-)$ is at most $\Delta_0^{-C}$, we have
\begin{equation}
\label{eq2_fin}
\mathfrak{S}_2 \leq  \max_{Q^\# \; \text{keystone}} \left\{   \Delta_\junk^p \Delta_0^{-C} \lvert R^\#_{Q^\#} \rvert^p :  S_1 Q^\# \subset \frac{65}{64}\hQ \right\}.
\end{equation}

Let $Q^\#$ be a keystone cube with $S_1 Q^\# \subset \frac{65}{64}\hQ$. Lemma \ref{new_fin_lem} states that
\[\lvert R_{Q^\#}^\# - P \rvert \leq C \cdot \left[  \| (f,P) \|_{S_1 Q^\#} + \Delta_\junk \lvert P \rvert \right].\]
Hence,
\begin{equation}
\label{eq3_fin}
\lvert R_{Q^\#}^\# \rvert^p \leq C \cdot \left[ \| (f,P) \|_{S_1 Q^\#}^p +   \lvert P \rvert^p \right].
\end{equation}
We define
\[ \widetilde{F} = P + \sum_{x \in S_1 Q^\# \cap E} \theta_x \cdot (f(x)-P(x))
\]
where $\theta_x(y)$ ($x \in E$) are cutoff functions satisfying (a) $\theta_x \equiv 1$ on a neighborhood of $x$, (b)  $\theta_x$ is supported on a ball $B(x,c \Delta_0)$ for a small universal constant $c$, and (c) $\| \partial^\alpha \theta_x \|_{L^\infty} \leq \Delta_0^{-C}$ for all $|\alpha | \leq m$. Indeed note that we may take $\theta_x(y) = \theta(y-x)$ for a fixed cutoff function $\theta$ supported on a small ball about the origin. From (a) and (b) we deduce that $\theta_x(z) \equiv 0$ for any $z \in E \setminus \{x\}$, since $\lvert x - y \rvert \geq \Delta_0$ for distinct points $x,y \in E$. Thus,  we have $\widetilde{F}(x) = f(x)$ for each $x \in E$. Moreover,
\[
\| \widetilde{F} \|^p_{\X(S_1 Q^\#)} \leq \Delta_0^{-C} \sum_{x \in S_1 Q^\# \cap E} \lvert f(x)-P(x) \rvert^p
\]
and
\[
\| \widetilde{F} - P \|_{L^p(S_1 Q^\#)}^p \leq \Delta_0^{-C} \sum_{x \in S_1 Q^\# \cap E} \lvert f(x)-P(x) \rvert^p.
\]
Hence, by definition of the trace seminorm, $\|(f,P)\|_{S_1 Q^\#}^p \leq \Delta_0^{-C} \sum_{x} \lvert f(x)-P(x) \rvert^p$. Thus, estimate \eqref{eq3_fin} implies that
\begin{equation}
\label{eq4_fin}
\lvert R_{Q^\#}^\# \rvert^p \leq C   \left[ \Delta_0^{-C'} \sum_{x \in S_1 Q^\# \cap E} \lvert f(x) - P(x) \rvert^p +  \lvert P \rvert^p \right].
\end{equation}

Therefore, returning to \eqref{eq2_fin}, we have
\begin{equation*}
\mathfrak{S}_2  \leq  \Delta_\junk^p \Delta_0^{-C''}  \left[ \sum_{x \in \frac{65}{64} \hQ \cap E} \lvert f(x) - P(x) \rvert^p +  \lvert P \rvert^p  \right].
\end{equation*}

Since $\Delta_\junk^p \Delta_0^{-C''} \leq \Delta_\new^{4p} \cdot \left[ \Delta_\new^p \Delta_0^{-C''} \right] \leq \Delta_\new^{4p}$ (see \eqref{constants3}), we conclude that
\begin{align*}
\mathfrak{S}_2 \leq \Delta_\new^{4p} \sum_{x \in \frac{65}{64} \hQ \cap E} \lvert f(x) - P(x) \rvert^p +  \Delta_\new^{4p} \lvert P \rvert^p &\leq \Delta_\new^{2p} \left[ M_\hQ(f,P) \right]^p\\
&\leq \frac{1}{2 C(t_G)} \left[ M_\hQ(f,P) \right]^p,
\end{align*} 
where $C(t_G)$ is the constant in \eqref{eq1_fin}. To obtain the previous estimates, we make sure to choose $\Delta_\new^{2p} \leq \frac{1}{2 C(t_G)}$ (see \eqref{constants4}).

This completes our estimation of the term $\mathfrak{S}_2$.

In the estimate \eqref{eq1_fin}, we consider the term $C(t_G) \cdot \mathfrak{S} = C(t_G) \cdot \mathfrak{S}_1 +  C(t_G) \cdot \mathfrak{S}_2$ on the right-hand side, and note that $C(t_G) \cdot \mathfrak{S}_2$ is irrelevant since it is bounded by a half of the left-hand side; moreover, the term $ C(t_G) \cdot \mathfrak{S}_1$ is bounded from above by $C(t_G) \cdot \left[ \text{RHS of conditional inequality} \right]^p$, thanks to \eqref{temp11}. Therefore, in place of \eqref{eq1_fin}, we have the simpler estimate:
\begin{align} 
\label{e773_fin}
\bigl[ M_{\hQ}(f,P)\bigr]^p \leq C(t_G) \cdot & \biggl[ \|H \|_{\X(\frac{65}{64}\hQ)}^p  +  \delta_\hQ^{-mp}  \| H - P \|_{L^p(\frac{65}{64}\hQ)}^p  \\
\notag{}
& + \sum_{Q \subset (1+100t_G)\hQ } \bigl[ \| H\|_{\X(\frac{65}{64}Q)}^p + \delta_Q^{-mp} \| H - R^\hQ_Q \|_{L^p(\frac{65}{64}Q)}^p \bigr] \\
\notag{}
& +  \left[ \text{RHS of conditional inequality} \right]^p \biggr].
\end{align}
The difference between the estimates \eqref{e773_fin} and  \eqref{e773} is that the right-hand side of \eqref{e773_fin}  contains an extra term: $\left[ \text{RHS of conditional inequality} \right]^p$.

The estimates preceding \eqref{e774} in \textbf{Stage II} are unchanged. Using these estimates in \eqref{e773_fin}, we obtain
\begin{align}
\label{e774_fin}
\left[M_{\hQ}(f,P)\right]^p & \leq C(t_G) \cdot \biggl(  \|H \|_{\X(\frac{65}{64}\hQ)}^p + \delta_\hQ^{-mp}  \| H - P \|_{L^p(\frac{65}{64}\hQ)}^p  + \sum_{\substack{ Q \subset (1+100t_G) \hQ \\ \delta_Q < t_G \delta_\hQ}} \delta_Q^{-mp} \| H - R^\#_{\mathcal{K}(Q)} \|_{L^p(\frac{65}{64}Q)}^p \\
\notag{}
& + \bigl[ \text{RHS of conditional inequality} \bigr]^p \biggr) \notag{} \\
\notag{}
& \leq C(t_G) \cdot \biggl(  \|H \|_{\X(\frac{65}{64}\hQ)}^p +  \delta_\hQ^{-mp} \| H - P \|_{L^p(\frac{65}{64}\hQ)}^p  + \sum_{\substack{ Q^\# \; \tiny{\mbox{keystone}} \\ \; S_1 Q^\# \subset \frac{65}{64}\hQ}} (\delta_{Q^\#})^{-mp} \| H - R^\#_{Q^\#} \|^p_{L^{p}(S_1Q^\#)} \\
\notag{}
& + \bigl[ \text{RHS of conditional inequality} \bigr]^p \biggr).
\end{align}
The difference between the estimates \eqref{e774_fin} and \eqref{e774} is the extra  term \\
$\bigl[ \text{RHS of conditional inequality} \bigr]^p$ that appears in \eqref{e774_fin}.

We now examine the estimates in \textbf{Stage III}.

In place of the infinite-precision inequality \eqref{wc3}, which reads
\[
\delta_{Q^\#}^{-mp} \| H - R_{Q^\#}^\# \|_{L^p(S_1 Q^\#)}^p \leq C \| H \|_{\X(S_1 Q^\#)}^p,
\]
we now apply \eqref{goodversion} which reads
\[
\delta_{Q^\#}^{-mp} \| H - R_{Q^\#}^\# \|_{L^p(S_1 Q^\#)}^p \leq C \left[ \| H \|_{\X(S_1 Q^\#)}^p  + \Delta_ \junk^p \lvert R_{Q^\#}^\# \rvert^p \right].
\] 
From \eqref{eq4_fin} we thus have
\[
\delta_{Q^\#}^{-mp} \| H - R_{Q^\#}^\# \|_{L^p(S_1 Q^\#)}^p  \leq C \left[ \| H \|_{\X(S_1 Q^\#)}^p  +  \Delta_ \junk^p \cdot  \left[ \Delta_0^{-C'} \sum_{x \in S_1 Q^\# \cap E} \lvert f(x) - P(x) \rvert^p +  \lvert P \rvert^p \right] \right].
\]
There are at most $C N \leq \Delta_0^{-C}$ keystone cubes in $\CZ(\cA^-)$. Hence, since the collection $\{ S_1 Q^\# : Q^\# \text{ keystone} \}$ has bounded overlap, we have
\[
\sum_{\substack{Q^\# \text{ keystone} \\ S_1 Q^\# \subset \frac{65}{64}\hQ}}  \delta_{Q^\#}^{-mp} \| H - R_{Q^\#}^\# \|_{L^p(S_1 Q^\#)}^p \leq C \| H \|_{\X(\frac{65}{64}\hQ)}^p  + \left\{ C \Delta_ \junk^p \Delta_0^{-C}  \left[ \sum_{x \in \frac{65}{64}\hQ \cap E} \lvert f(x) - P(x) \rvert^p + \lvert P \rvert^p \right] \right\}.
\]
We have $C \Delta_\junk^p \Delta_0^{-C} \leq C \Delta_\new^{3p} \leq \frac{1}{2 C(t_G)}\Delta_\new^{2p}$, due to the assumptions \eqref{constants3} and \eqref{constants4}. Thus, the term inside the curly braces in the above estimate is bounded by
\[ 
\frac{1}{2 C(t_G)} \left[ \text{ \bf (V)} + \text{ \bf (VI) } \right] \leq \frac{1}{2 C(t_G)} \left[ M_\hQ(f,P) \right]^p.
\]
We put the previous estimates into \eqref{e774_fin} to obtain
\begin{align*}
\left[ M_{\hQ}(f,P) \right]^p & \leq C(t_G) \cdot \left[C \|H \|_{\X(\frac{65}{64}\hQ)}^p + \delta_\hQ^{-mp}  \| H - P \|_{L^p(\frac{65}{64}\hQ)}^p  + \frac{1}{2 C(t_G)} \left[ M_{\hQ}(f,P) \right]^p \right] \\
& \qquad + \bigl[ \text{RHS of conditional inequality} \bigr]^p.
\end{align*}
From the third bullet point in Proposition \ref{prop_A}, we deduce that
\[ 
\left[ M_{\hQ}(f,P) \right]^p \leq C(t_G) \cdot \Lambda^{(2D+1)p} \| (f,P) \|_{\frac{65}{64}\hQ}^p + \bigl[ \text{RHS of conditional inequality} \bigr]^p.
\]
Since $\Lambda^{2D+1} \leq 1/ \epsilon$, this estimate implies the conditional inequality in Proposition \ref{prop_normappx_fin}. This completes the proof of Proposition \ref{prop_normappx_fin}.

We fix $t_G > 0$ to be a universal constant, small enough so that the preceding results hold. We define the universal constant $a(\cA) = a_\new$, with $a_\new$ defined as in Lemma \ref{lem_cover}.

For a moment, we fix $\epsilon = \epsilon_0$ in Proposition \ref{prop_normappx_fin} for a small universal constant $\epsilon_0$. This implies the following result.

\begin{prop} \label{inc_prop_new}

There exist universal constants $\epsilon_0 > 0$ and $C \geq 1$ such that the following estimates hold.

\begin{description}
\item[Unconditional Inequality] $\| (f,P) \|_{(1+a(\cA))\hQ} \leq C \cdot M_{\hQ}(f,P)$.
\item[Conditional Inequality] If $3\hQ$ is tagged with $(\cA,\epsilon_0)$, then 
\[ M_{\hQ}(f,P) \leq C \left[ \|(f,P) \|_{\frac{65}{64}\hQ} + \Delta_\new \cdot \lvert P \rvert \right].\]
\end{description}

\end{prop}

We no longer fix $\epsilon = \epsilon_0$. Once again, we assume that $\epsilon$ is a small parameter, less than a small enough universal constant.

We assume that 
\begin{equation}\label{constants5}
\Delta_\new \leq c(\epsilon),
\end{equation}
for a small enough constant $c(\epsilon)$, depending only on $\epsilon$, $m$, $n$, and $p$.

We will need new proofs of Propositions \ref{tool1},\ref{tag_prop1},\ref{tag_prop2}. We recall the statements of these results and give the new proofs.

\noindent\textbf{Proposition \ref{tool1}}. \emph{ Let $\hQ$ be a testing cube. If 
$$\left[ \#\left( \frac{65}{64} \hQ \cap E \right) \leq 1 \;\; \mbox{or} \;\; \ooline{\sigma}(\hQ) \;\; \mbox{has an} \;\; (\cA',x_{\hQ},\epsilon,\delta_{\hQ})\mbox{-basis for some} \; \cA' \leq \cA \right]$$
then $(1+ a(\cA)) \hQ$ is tagged with $(\cA,\epsilon^\kappa)$. Otherwise, no cube containing $3 \hQ$ is tagged with $(\cA,\epsilon^{1/\kappa})$. Here, $\kappa>0$ is a universal constant.
}

\begin{proof}
If $\#(\frac{65}{64}\hQ \cap E) \leq 1$, then $(1+a(\cA))\hQ$ is tagged with $(\cA,\epsilon)$.

Suppose $\ooline{\sigma}(\hQ)$ has an $(\cA',x_{\hQ},\epsilon,\delta_{\hQ})$-basis with $\cA' \leq \cA$. Call this basis $(P_\alpha)_{\alpha \in \cA'}$. Then
\begin{itemize}
\item $P_\alpha \in \epsilon \cdot \delta_\hQ^{|\alpha| + n/p - m} \cdot \ooline{\sigma}(\hQ)$ \quad for all $\alpha \in \cA'$.
\item $\partial^\beta P_\alpha(x_\hQ) = \delta_{\alpha \beta} $ \qquad \quad for all $\alpha, \beta \in \cA'$.
\item $|\partial^\beta P_\alpha(x_\hQ)| \leq \epsilon \cdot \delta_\hQ^{|\alpha| - |\beta|}$ \;\; for all $\alpha \in \cA'$, $\beta \in \cM$, $\beta > \alpha$.
\end{itemize}
Since $\ooline{\sigma}(\hQ) = \{ P : M_{\hQ}(0,P) \leq 1 \}$, we have $M_{\hQ}(0,P_\alpha) \leq \epsilon \delta_\hQ^{|\alpha| + n/p - m}$ for $\alpha \in \cA'$. So, the Unconditional Inequality gives
\[\|(0,P_\alpha)\|_{(1+a(\cA))\hQ} \leq  C' \epsilon \delta_\hQ^{|\alpha| + n/p - m} \quad \text{for all} \; \alpha \in \cA'.\]
Thus,
\[P_\alpha \in C' \epsilon \delta_\hQ^{|\alpha| + n/p - m} \sigma((1+a(\cA))\hQ) \quad \text{for all} \; \alpha \in \cA'.\]
With the second and third bullet points above, this shows that $(P_\alpha)_{\alpha \in \cA'}$ is an $(\cA',x_{\hQ},C' \epsilon,\delta_\hQ)$-basis for $\sigma((1+a(\cA))\hQ)$. Hence, $(P_\alpha)_{\alpha \in \cA'}$ is an $(\cA',x_{\hQ},\epsilon^\kappa,\delta_{ (1+a(\cA)) \hQ})$-basis for $\sigma((1+a(\cA))\hQ)$, for a small enough universal constant $\kappa$. Since $\cA' \leq \cA$, it follows that $(1+a(\cA)) \hQ$ is tagged with $(\cA,\epsilon^\kappa)$, as claimed. That  proves the first half of Proposition \ref{tool1}.

On the other hand, suppose $Q \supset 3 \hQ$ and suppose $Q$ is tagged with $(\cA, \epsilon^{1/\kappa'})$, for a small enough universal constant $\kappa' > 0$, to be chosen below (not any previous $\kappa'$). Then $3 \hQ$ is tagged with $(\cA, \epsilon^{\kappa/\kappa'})$ for some universal constant $\kappa>0$, thanks to Lemma \ref{pre_lem5}. Hence, as long as $\epsilon$ is small enough so that $\epsilon^{\kappa/\kappa'} \leq \epsilon_0$, the Conditional Inequality applies:
\begin{equation*} M_{\hQ}(0,P) \leq C \left[ \| (0,P) \|_{\frac{65}{64}\hQ} +  \Delta_\new \lvert P \rvert \right] \qquad \mbox{for any} \; P \in \cP
\end{equation*}
Also, by Lemma \ref{pre_lem5}, $\frac{65}{64}\hQ$ is tagged with $(\cA,\epsilon^{\kappa/\kappa'})$. So either $\#(\frac{65}{64}Q \cap E) \leq 1$ (in which case we have finished the proof of Proposition \ref{tool1}) or else $\sigma(\frac{65}{64}Q)$ has an $(\cA',x_\hQ,\epsilon^{\kappa/\kappa'},\delta_\hQ)$-basis for some $\cA' \leq \cA$. 

In the latter case, Lemma \ref{pre_lem2} gives an $(\cA'', x_\hQ,\epsilon^{\overline{\kappa}/\kappa' },\delta_\hQ,\Lambda)$-basis, with 
\begin{align*}
\cA'' \leq \cA' \leq \cA, \; \widetilde{\kappa} \leq \overline{\kappa} \leq \widetilde{\widetilde{\kappa}}, \; \mbox{with} \; \widetilde{\kappa}, \widetilde{\widetilde{\kappa}} > 0  \; \mbox{universal constants} \; \mbox{independent of} \; \kappa', \\
\mbox{ and } \epsilon^{\overline{\kappa}/\kappa'} \Lambda^{100D} \leq \epsilon^{\overline{\kappa}/2\kappa'}.
\end{align*} 
Call this basis $(P_\alpha)_{\alpha \in \cA''}$. Then
\begin{itemize}
\item $P_\alpha \in \epsilon^{\overline{\kappa}/\kappa'} \cdot \delta_\hQ^{|\alpha| + n/p - m} \cdot \sigma(\frac{65}{64}Q)$ \quad for all $\alpha \in \cA''$.
\item $\partial^\beta P_\alpha(x_\hQ) = \delta_{\alpha \beta} $ \qquad \quad for all $\alpha, \beta \in \cA''$.
\item $|\partial^\beta P_\alpha(x_\hQ)| \leq \epsilon^{\overline{\kappa}/\kappa'} \cdot \delta_\hQ^{|\alpha| - |\beta|}$ \;\; for all $\alpha \in \cA''$, $\beta \in \cM$, $\beta > \alpha$.
\item $|\partial^\beta P_\alpha(x_\hQ)| \leq \Lambda \cdot \delta_\hQ^{|\alpha| - |\beta|}$ \;\; for all $\alpha \in \cA''$, $\beta \in \cM$.
\end{itemize}
We deduce a few conclusions from the above bullet points. The first bullet point implies that $\|(0,P_\alpha)\|_{\frac{65}{64}\hQ} \leq \epsilon^{\overline{\kappa}/\kappa'} \delta_\hQ^{|\alpha|+n/p-m}$; the last bullet point implies that $|\partial^\beta P_\alpha(x_\hQ)| \leq \Lambda \cdot \Delta_0^{-C}$ for all $\alpha \in \cA''$, $\beta \in \cM$ (since $\Delta_0 \leq \delta_\hQ \leq 1$), hence
\begin{align*}
\lvert P_\alpha \rvert  = \left(\sum_{\beta \in \cM} \lvert \partial^\beta P_\alpha(0) \rvert^p \right)^{1/p} &\leq C \cdot \left(\sum_{\beta \in \cM} \lvert \partial^\beta P_\alpha(x_{\hQ}) \rvert^p \right)^{1/p} \qquad (\mbox{since} \; \lvert x_{\hQ} \rvert \leq C) \\
& \leq C' \Lambda \cdot \Delta_0^{-C}.
\end{align*}
Hence, the (known) Conditional Inequality implies the estimate
\begin{align}
\label{bbb2}
M_{\hQ}(0,P_\alpha) &\leq  C \cdot \left[ \| (0,P_\alpha) \|_{\frac{65}{64}\hQ} +  \Delta_\new \cdot \lvert P_\alpha \rvert \right] \\
& \leq C \epsilon^{\overline{\kappa}/\kappa'} \delta_\hQ^{|\alpha| + n/p - m} + C \Lambda \Delta_0^{-C} \Delta_\new \leq C' \epsilon^{\overline{\kappa}/\kappa'} \delta_\hQ^{|\alpha| + n/p - m} \quad \mbox{for} \; \alpha \in \cA'',
\notag{}
\end{align}
where we use that $\delta_\hQ \leq 1$ and $\lv \alpha \rv + n/p - m < 0$, and
\[ 
\Lambda \cdot  \Delta_0^{-C} \Delta_\new \oleq{\eqref{constants3}} \Lambda \cdot \Delta_\new^{1/2} \oleq{\eqref{constants5}} \Lambda \cdot  \epsilon^{2\cdot \widetilde{\widetilde{\kappa}}/\kappa'} \leq \epsilon^{\overline{\kappa}/\kappa'}.
\]
Here, in the last inequality, we use that $\Lambda \leq \Lambda^{100D} \leq \epsilon^{-\overline{\kappa}/2\kappa'}$, where $\overline{\kappa} \leq \widetilde{\widetilde{\kappa}}$.

Now, the estimate \eqref{bbb2} implies that
\[P_\alpha \in C' \epsilon^{\overline{\kappa}/\kappa'} \delta_\hQ^{|\alpha| + n/p - m} \cdot \ooline{\sigma}(\hQ) \quad \mbox{for all} \; \alpha \in \cA''.\]
This estimate, together with the second and third bullet points above, shows that 
\[
(P_\alpha)_{\alpha \in \cA''} \; \mbox{is an} \; (\cA'',x_{\hQ},C \epsilon^{\overline{\kappa}/\kappa'},\delta_\hQ)\mbox{-basis for} \; \ooline{\sigma}(\hQ).
\]
We ensure that $C \epsilon^{\overline{\kappa}/\kappa'} \leq \epsilon$ by choosing $\kappa'$ to be a small enough universal constant. Hence, $\ooline{\sigma}(\hQ)$ has an $(\cA'',x_\hQ,\epsilon,\delta_\hQ)$-basis with $\cA'' \leq \cA' \leq \cA$. This completes the proof of Proposition \ref{tool1}.

\end{proof}

\noindent\textbf{Proposition \ref{tag_prop1}}. \emph{
Suppose $\hQ_1 \subset \hQ_2$ are testing cubes with $\#(3\hQ_2 \cap E) \geq 2$, and $(1+a(\cA))\hQ_1 \cap E = 3 \hQ_2 \cap E$.
Suppose $\ooline{\sigma}(\hQ_1)$ has an $(\cA',x_{\hQ_1},\epsilon,\delta_{\hQ_2})$-basis. Then $3 \hQ_2$ is tagged with $(\cA',\epsilon^\kappa)$ for a universal constant $\kappa$.
}

\begin{proof}

By Lemma \ref{pre_lem4}, $\ooline{\sigma}(\hQ_1)$ has an $(\cA'',x_{\hQ_2},\epsilon^\kappa,\delta_{\hQ_2})$-basis, with $\cA'' \leq \cA'$, for some universal constant $\kappa$. Call that basis $(P_\alpha)_{\alpha \in \cA''}$. Then for each $\alpha \in \cA''$, we have
\begin{itemize}
\item $P_\alpha \in \epsilon^{{\kappa}} \cdot \delta_{\hQ_2}^{|\alpha| + n/p - m} \cdot \ooline{\sigma}(\hQ_1)$.
\item $\partial^\beta P_\alpha(x_{\hQ_2}) = \delta_{\alpha \beta} $ \qquad \quad for all $\beta \in \cA''$.
\item $|\partial^\beta P_\alpha(x_{\hQ_2})| \leq \epsilon^{{\kappa}} \cdot \delta_{\hQ_2}^{|\alpha| - |\beta|}$ \;\; for $\beta \in \cM$, $\beta > \alpha$.
\end{itemize}
The first condition here gives $M_{\hQ_1}(0,P_\alpha) \leq \epsilon^{{\kappa}} \delta_{\hQ_2}^{| \alpha | + n/p - m}$. So, by the Unconditional Inequality,
\[\| (0,P_\alpha) \|_{(1+a(\cA))\hQ_1} \leq C \epsilon^{{\kappa}} \delta_{\hQ_2}^{|\alpha| + n/p - m}.\]
Hence, $P_\alpha \in C \epsilon^{{\kappa}} \delta_{\hQ_2}^{|\alpha| + n/p - m} \sigma((1+a(\cA))\hQ_1)$. By Lemma \ref{pre_lem0}, $\sigma(3 \hQ_2)$ is comparable to $\sigma((1+a(\cA))\hQ_1) + \cB(x_{\hQ_2},\delta_{\hQ_2})$, so
\[\sigma((1+a(\cA)) \hQ_1) \subset C \sigma(3 \hQ_2).\]
Thus, $P_\alpha \in C \epsilon^{{\kappa}} \delta_{\hQ_2}^{|\alpha| + n/p - m} \sigma(3 \hQ_2)$ for all $\alpha \in \cA''$. With the second and third bullet points above, this shows that $\sigma(3\hQ_2)$ has an $(\cA'',x_{\hQ_2},C \epsilon^{{\kappa}} \delta_{\hQ_2})$-basis, with $ \cA'' \leq \cA'$. Therefore, $3 \hQ_2$ is tagged with $(\cA',\epsilon^{\kappa/2})$, if $\epsilon$ is less than a small enough universal constant. This completes the proof of Proposition \ref{tag_prop1}.
\end{proof}

Corollary \ref{cor_tag1} is a direct consequence  of Proposition \ref{tag_prop1}, just as before.

\noindent \textbf{Proposition \ref{tag_prop2}}. \emph{
Suppose that $\hQ_1 \subset \hQ_2$ are testing cubes,  $\#(3\hQ_2 \cap E) \geq 2$, and $(1+a(\cA))\hQ_1 \cap E = 3 \hQ_2 \cap E$. Suppose  $3 \hQ_2$ is tagged with $(\cA,\epsilon)$. Then $\ooline{\sigma}(\hQ_1)$ has an $(\cA',x_{\hQ_1},\epsilon^{\kappa'},\delta_{\hQ_2})$-basis for some $\cA' \leq \cA$ and for some universal constant $\kappa'$.
}

\begin{proof}
Since $3 \hQ_1 \subset 3 \hQ_2$ and $3 \hQ_2$ is tagged with $(\cA,\epsilon)$, Lemma \ref{pre_lem5} shows that $3 \hQ_1$ is tagged with $(\cA,\epsilon^\kappa)$ for a universal constant $\kappa$. Hence, the Conditional Inequality holds for $\hQ_1$. Hence,
\begin{equation}\label{smiley}
M_{\hQ_1}(0,P) \leq C \left[ \| (0,P)\|_{\frac{65}{64}\hQ_1} +  \Delta_\new \lvert P \rvert \right] \qquad \mbox{for} \; P \in \cP.
\end{equation}
Now, since $\frac{65}{64}\hQ_1 \cap E = 3 \hQ_2 \cap E$ and $\frac{65}{64}\hQ_1 \subset 3 \hQ_2$, we know from Lemma \ref{pre_lem0} that
\begin{equation*}
\sigma(3 \hQ_2) \subset C \cdot \left[ \sigma\left(\frac{65}{64}\hQ_1\right) + \cB(x_{\hQ_2},\delta_{3 \hQ_2})  \right] .
\end{equation*}
(We have $\cB(x_{\hQ_1}, \delta_{3 \hQ_2}) \subset C \cB(x_{\hQ_2},\delta_{3 \hQ_2})$, because $\lvert x_{\hQ_1} - x_{\hQ_2} \rvert \leq \delta_{\hQ_2}$. Hence, the above inclusion follows from Lemma \ref{pre_lem0}.)

Recall that $3 \hQ_2$ is tagged with $(\cA,\epsilon)$ and $\#(3 \hQ_2 \cap E) \geq 2$. Hence, $\sigma(3 \hQ_2)$ has an $(\cA',x_{\hQ_2},\epsilon,\delta_{3 \hQ_2})$-basis, for some $\cA' \leq \cA$. By Lemma \ref{pre_lem2}, there exist a multiindex set $\cA'' \leq \cA' \leq \cA$ and numbers $\Lambda \geq 1$, $\kappa_1 \leq \overline{\kappa}  \leq \kappa_2$, such that
\[ \sigma(3 \hQ_2) \; \mbox{has an} \; (\cA'',x_{\hQ_2},\epsilon^{\overline{\kappa}}, \delta_{3 \hQ_2},\Lambda)\mbox{-basis, where} \; \epsilon^{\overline{\kappa}} \Lambda^{100D} \leq \epsilon^{\overline{\kappa}/2},
\]
for some universal constants $\kappa_1,\kappa_2 \in (0,1]$. Therefore,
\[\sigma\left( \frac{65}{64}\hQ_1 \right) + \cB(x_{\hQ_2},\delta_{3 \hQ_2}) \; \mbox{has an} \; (\cA'',x_{\hQ_2},C \epsilon^{\overline{\kappa}}, \delta_{3 \hQ_2},\Lambda)\mbox{-basis}.\]
From Lemma \ref{pre_lem1}, we see that $\sigma \left(\frac{65}{64}\hQ_1 \right)$ has an $(\cA'',x_{\hQ_2},C' \epsilon^{\overline{\kappa}} \Lambda, \delta_{3 \hQ_2}, C \Lambda)$-basis. Here, $C' \epsilon^{\overline{\kappa}} \Lambda \leq C' \epsilon^{\overline{\kappa}/2} \leq \epsilon^{\overline{\kappa}/4}$, for sufficiently small $\epsilon$. Let $(P_\alpha)_{\alpha \in \cA}$ be that basis. Thus, for each $\alpha \in \cA''$,
\begin{itemize}
\item $P_\alpha \in \epsilon^{\overline{\kappa}/4} \cdot \delta_{\hQ_2}^{|\alpha| + n/p - m} \cdot \sigma(\frac{65}{64}\hQ_1)$.
\item $\partial^\beta P_\alpha(x_{\hQ_2}) = \delta_{\alpha \beta} $ \qquad \quad for all $\beta \in \cA''$.
\item $|\partial^\beta P_\alpha(x_{\hQ_2})| \leq \epsilon^{ \overline{\kappa}/4} \cdot \delta_{\hQ_2}^{|\alpha| - |\beta|}$ \;\; for all $\beta \in \cM$, $\beta > \alpha$.
\item $|\partial^\beta P_\alpha(x_{\hQ_2})| \leq C \Lambda \cdot \delta_{\hQ_2}^{|\alpha| - |\beta|}$ \;\; for all $\beta \in \cM$.
\end{itemize}
The first and fourth bullet points imply that $\| (0,P_\alpha) \|_{\frac{65}{64}\hQ_1} \leq C \epsilon^{\overline{\kappa}/4} \delta_{\hQ_2}^{|\alpha| + n/p - m}$ and $\lvert P_\alpha \rvert \leq C \Lambda \Delta_0^{-C}$, hence \eqref{smiley} gives 
\begin{equation}
\label{bbb3} 
M_{\hQ_1}(0,P_\alpha) \leq C' \epsilon^{ \overline{\kappa}/4} \delta_{\hQ_2}^{|\alpha| + n/p - m} + C \Lambda \Delta_0^{-C} \Delta_\new \leq C' \epsilon^{ \overline{\kappa}/4} \delta_{\hQ_2}^{|\alpha| + n/p - m},
\end{equation}
which implies that 
\begin{equation}
\label{bbb4}
P_\alpha \in C' \epsilon^{ \overline{\kappa}/4} \delta_{\hQ_2}^{|\alpha| + n/p - m} \cdot \ooline{\sigma}(\hQ_1) \qquad \mbox{for} \; \alpha \in \cA''.
\end{equation}
Here, to prove \eqref{bbb3}, we use that $\delta_{\hQ_2} \leq 1$ and $\lv \alpha \rv + n/p - m < 0$, and
\[ 
\Lambda \cdot  \Delta_0^{-C} \Delta_\new \oleq{\eqref{constants3}} \Lambda \cdot \Delta_\new^{1/2} \oleq{\eqref{constants5}} \Lambda \cdot  \epsilon^{2 \kappa_2} \leq \epsilon^{\overline{\kappa}}.
\]
Here, in the last inequality, we use that $\Lambda \leq \Lambda^{100D} \leq \epsilon^{-\overline{\kappa}/2}$, where $\overline{\kappa} \leq \kappa_2$.

With the second and third bullet points, \eqref{bbb4} shows that $(P_\alpha)_{\alpha \in \cA''}$ forms an \\
$(\cA'',x_{\hQ_2},C \epsilon^{\overline{\kappa}/4},\delta_{\hQ_2}$)-basis for $\ooline{\sigma}(\hQ_1)$. Hence, by Lemma \ref{pre_lem4}, it follows that $\ooline{\sigma}(\hQ_1)$ has an $(\cA''',x_{\hQ_1},\epsilon^{\kappa'},\delta_{\hQ_2})$-basis for some $\cA''' \leq \cA'' \leq \cA' \leq \cA$ and for a small enough universal constant $\kappa'$. This completes the proof of Proposition \ref{tag_prop2}.
\end{proof}

The proofs of Propositions \ref{tag_prop3} and \ref{tool2} are unchanged.

The statement of the algorithm \textsc{Optimize Basis} requires modification.

\environmentA{Algorithm: Optimize Basis (Finite-Precision)}

We perform one time work at most $C N \log N$ in space $C N$, after which we can answer queries as follows.

A query consists of a testing cube $\hQ$ and a set $\cA \subset \cM$

We respond to the query $(\hQ,\cA)$ by producing the following.
\begin{itemize}
\item A collection of machine intervals $I_\ell$ ($1 \leq \ell \leq \ell_{\max}$). The intervals $I_\ell$ are pairwise disjoint, the union of the $I_\ell$ is $\left[ \Delta_g, \Delta_g^{-1} \right]$, and $\ell_{\max} \leq C$.
\item A list of non-negative machine numbers $a_\ell$ ($\ell=1,\cdots,\ell_{\max}$). The numbers $a_\ell$ are bounded in magnitude by $\Delta_g^{-C}$.
\item A list of numbers $\lambda_\ell$. Each $\lambda_\ell$ has the form $\mu_\ell + \nu_\ell/p$, with $\mu_\ell,\nu_\ell \in \Z$ and $\lv \mu_\ell \rv, \lv \nu_\ell \rv \leq C$.
\item Let $\eta^{(\hQ,\cA)}(\delta) := a_\ell \delta^{\lambda_\ell}$ for $\delta \in I_\ell$. Then we have:
\begin{description}
\item[(A1)] For each $\delta \in [ \Delta_g , \Delta_g^{-1}]$ there exists $\cA' \leq \cA$ such that $\ooline{\sigma}(\hQ)$ has an $(\cA',x_{\hQ},\eta^{1/2},\delta)$-basis for all $\eta > C \cdot \eta^{(\hQ,\cA)}(\delta)$.

\item[(A2)] For each $\delta \in [ \Delta_g , \Delta_g^{-1}]$ and any $\cA' \leq \cA$, $\ooline{\sigma}(\hQ)$ \underline{does not} have an $(\cA',x_{\hQ}, \eta^{1/2},\delta)$-basis for any $ \eta < c \cdot \eta^{(\hQ,\cA)}(\delta)$.

\item[(A3)] Moreover, $c \cdot \eta^{(\hQ,\cA)}(\delta_1) \leq \eta^{(\hQ,\cA)}(\delta_2) \leq C \cdot \eta^{(\hQ,\cA)}(\delta_1)$ whenever $\frac{1}{10} \delta_1 \leq \delta_2 \leq 10 \delta_1$ and $\delta_1,\delta_2 \in [ \Delta_g,\Delta_g^{-1} ]$.

\item[(A4)] Also, $\eta^{(\hQ,\cA)}(\delta) \geq \Delta_g^C$, for all $\delta \in [ \Delta_g,\Delta_g^{-1}]$.
\end{description}
\item To answer a query requires work at most $C \log N$.
\end{itemize}

\begin{proof}[\underline{Explanation}] Recall that we can perform arithmetic operations to within precision $\Delta_\epsilon$.

We denote $\mathbb{Z}\left[\frac{1}{p} \right] = \left\{ \lambda = k + \ell \cdot \frac{1}{p} : k,\ell \in \Z \right\}$. If $\lambda = k + \ell \cdot \frac{1}{p} \in \mathbb{Z}\left[\frac{1}{p}\right]$, with $k$ and $\ell$ bounded by a universal constant, then we say that $\lambda$ is a machine element of $\mathbb{Z} \left[ \frac{1}{p} \right]$. Such a $\lambda$ can be stored on our computer  using at most $C$ units of storage.


Recall that we defined $\displaystyle \lvert P \rvert_x = \left( \sum_{\alpha \in \cM} \lvert \partial^\alpha P(x) \rvert^p  \right)^{1/p}$ for $P \in \cP$ and $x \in \R^n$, and $\lvert P \rvert = \lvert P \rvert_0$. The vectors $(\partial^\alpha P(x))_{\alpha \in \cM}$ and $(\partial^\alpha P(0) )_{\alpha \in \cM}$ are related by multiplication against an invertible matrix $A(x) = (A_{\alpha \beta}(x) )_{\alpha, \beta \in \cM}$. This is a consequence of  Taylor's formula. Note that the operator norm of the matrix $A(x)$ is bounded by a universal constant if $\lvert x \rvert \leq 1$. Thus,
\begin{equation}
\label{simnorm}
C^{-1} \lvert P \rvert _x \leq \lvert P \rvert \leq C \lvert P \rvert_x \qquad \mbox{for} \; P \in \cP \; \mbox{and} \;  \lvert x \rvert \leq 1.
\end{equation}

Using \textsc{Approximate New Trace Norm} (see Section \ref{sec_newfunc}), we compute a quadratic form $q_\hQ$ on $\cP$ such that $\{ q_\hQ \leq c \} \subset \ooline{\sigma}(\hQ) \subset \{ q_\hQ \leq C \}$, where $c > 0$ and $C \geq 1$ are universal constants.\footnote{Recall that by now we have fixed $t_G$ to be a universal constant; hence, the constants $c(t_G)$ and $C(t_G)$ in \textsc{Approximate New Trace Norm} are now universal constants $c$, $C$.} The quadratic form $q_{\hQ}$ is given in the form
\[
q_{\hQ}(P) = \sum_{\alpha,\beta \in \cM} \widetilde{q}_{\alpha \beta} \cdot \frac{1}{\alpha!} \partial^\alpha P(0) \cdot \frac{1}{\beta!} \partial^\beta P(0),
\]
where we compute the numbers $\widetilde{q}_{\alpha \beta}$ with parameters $(\Delta_g^C,\Delta_g^{-C} \Delta_\epsilon)$. Using a linear change of basis, we write
\[
q_{\hQ}(P) = \sum_{\alpha , \beta \in \cM} q_{\alpha \beta} \cdot \frac{1}{\alpha!} \partial^\alpha P(x_{\hQ}) \cdot \frac{1}{\beta!} \cdot \partial^\beta P(x_{\hQ}).
\]
Each $q_{\alpha \beta}$ is a  linear combination of all the numbers $\widetilde{q}_{\alpha \beta}$. Thus,  we can compute each $q_{\alpha \beta}$ with parameters $(\Delta_g^C,\Delta_g^{-C} \Delta_\epsilon)$.

From the conditions in the algorithm \textsc{Approximate New Trace Norm}, we know that $q_{\hQ}(P) \geq c \cdot (M_\hQ(0,P))^2$. Furthermore, the term $\text{\textbf{(VI)}} = \Delta_\new^p \cdot \lvert P \rvert^p$ is a summand in  $\left[ M_{\hQ}(0,P) \right]^p$, hence $M_{\hQ}(0,P) \geq  \Delta_\new \cdot \lvert P \rvert$ (see \eqref{newterms}).  Hence, using \eqref{simnorm}, we see that
\[q_\hQ(P) \geq c' \Delta_\new^2 \lvert P \rvert^2 \geq c'' \Delta_{\new}^2 \lvert P \rvert_{x_{\hQ}}^2.\]
(Note that  $\lvert x_{\hQ} \rvert \leq 1$, since $\hQ \subset Q^\circ = [0,1)^n$.)

Therefore, the matrix $(q_{\alpha \beta})$ satisfies 
\[ 
(q_{\alpha \beta}) \geq c \Delta_\new^2 \cdot (\delta_{\alpha \beta}) \geq \Delta_g \cdot (\delta_{\alpha \beta}).
\]
This means that we can apply the finite-precision version of \textsc{Fit Basis to Convex Body} to the matrix $(q_{\alpha \beta})$ and the convex body $\ooline{\sigma}(\hQ)$ (see Section \ref{bases_sec_fin}). We can therefore compute a piecewise monomial function $\eta^{(\hQ,\cA')}_*(\delta)$ for each $\cA' \leq \cA$. We guarantee that
\begin{itemize}
\item For any $\delta \in [ \Delta_g,\Delta_g^{-1}]$,
\begin{itemize}
\item \textbf{(P1)} $\ooline{\sigma}(\hQ)$ has an $(\cA',x_{\hQ},\eta^{1/2},\delta)$-basis, for any $\eta > C \cdot \eta_*^{(\hQ,\cA')}(\delta)$,
\item \textbf{(P2)} $\ooline{\sigma}(\hQ)$ does not have an $(\cA',x_{\hQ},\eta^{1/2},\delta)$-basis, for any $\eta < c \cdot \eta_*^{(\hQ,\cA')}(\delta)$.

\end{itemize}
\item \textbf{(P3)} Moreover, $c \cdot \eta_*^{(\hQ,\cA')}(\delta_1) \leq \eta_*^{(\hQ,\cA')}(\delta_2) \leq C \cdot \eta_*^{(\hQ,\cA')}(\delta_1)$, whenever $\frac{1}{10} \delta_1 \leq \delta_2 \leq 10 \delta_1$, for $\delta_1,\delta_2  \in [ \Delta_g,\Delta_g^{-1}]$.
\item \textbf{(P4)} Also, $\eta_*^{(\hQ,\cA')}(\delta) \geq \Delta_g^C$ for any $\delta \in [ \Delta_g,\Delta_g^{-1}]$.
\item The function $\eta_*^{(\hQ,\cA')} : \left[ \Delta_g, \Delta_g^{-1} \right] \rightarrow \R$ is given in the form
\[
\eta^{(\hQ,\cA')}_*(\delta) = a_{\ell, \cA'} \cdot \delta^{\lambda_{\ell,\cA'}} \quad \mbox{for} \; \delta \in I_{\ell, \cA'}.
\]
To represent $\eta_*^{(\hQ,\cA')}$ we store the following data: pairwise disjoint machine intervals $I_{\ell,\cA'}$ ($1 \leq \ell \leq \ell_{\max}(\cA')$) that  form a partition  of $\left[ \Delta_g, \Delta_g^{-1} \right]$; machine numbers $a_{\ell,\cA'} \in [\Delta_g^C, \Delta_g^{-C}]$; and exponents $\lambda_{\ell,\cA'}$ that are machine elements of $\mathbb{Z}\left[ \frac{1}{p} \right]$. We guarantee that $\ell_{\max}(\cA') \leq C$ for each $\cA' \leq \cA$.
\end{itemize}
By computing all the nonempty  intersections of the intervals $I_{\ell,\cA'}$, we write each
$\eta^{(\hQ,\cA')}_*(\delta)$ in the form
\[
\eta^{(\hQ,\cA')}_*(\delta) = c_{\ell, \cA'} \cdot \delta^{\gamma_{\ell,\cA'}} \quad \mbox{for} \; \delta \in I_{\ell} \;\; (\ell=1,2,\cdots,\ell_{\max}).
\]
Here, we compute the following: machine intervals $I_\ell$ ($1 \leq \ell \leq \ell_{\max}$) that partition $\left[ \Delta_g, \Delta_g^{-1} \right]$; machine numbers  $c_{\ell,\cA'} \in [\Delta_g^C, \Delta_g^{-C}]$; and exponents $\gamma_{\ell,\cA'}$ that are machine elements in $\mathbb{Z}\left[ \frac{1}{p} \right]$. Moreover,  $\ell_{\max} \leq C$.

We define 
\begin{equation}
\label{defnofeta}
\eta(\delta) := \min_{\cA' \leq \cA} \eta_*^{(\hQ,\cA')}(\delta).
\end{equation}
We will compute a piecewise-monomial approximation to the function $\eta(\delta)$ using the following procedure.

\environmentA{Procedure: Process Monomials}. Assume that we are given the following: A machine interval $I \subset [ \Delta_g, \Delta_g^{-1} ]$, machine numbers $a_1,a_2 \in [\Delta_g^C,\Delta_g^{-C}]$, and machine elements  $\gamma_1, \gamma_2 \in \mathbb{Z} \left[ \frac{1}{p} \right]$. We define monomial functions $m_1(\delta) = a_1 \delta^{\gamma_1}$ and $m_2(\delta) = a_2 \delta^{\gamma_2}$. Then we produce one of three outcomes
\begin{enumerate}
\item We guarantee that $m_1(\delta) \leq m_2(\delta) + \Delta_\epsilon^{1/2}$ for all $\delta \in I$.
\item We guarantee that $m_2(\delta) \leq m_1(\delta) + \Delta_\epsilon^{1/2}$ for all $\delta \in I$.
\item We compute a  machine number $\delta_* \in I$, and distinct indices $j,k \in \{1,2\}$, such that
\[
\left\{
\begin{aligned}
m_j(\delta) & \leq m_k(\delta) + \Delta_\epsilon^{1/2} \;\; \mbox{for} \;\; \delta \in I \cap (0,\delta_*] \\
m_k(\delta) & \leq m_j (\delta) + \Delta_\epsilon^{1/2} \;\; \mbox{for} \;\; \delta \in I \cap [\delta_*,\infty).
\end{aligned}
\right.
\]
\end{enumerate}
This computation requires work and storage  at most $C$.

\begin{proof}[Explanation]
If $\gamma_1 = \gamma_2$ then outcome (1) occurs if $a_1 \leq a_2$, and outcome (2) occurs if $a_1 > a_2 $. Thus, we can respond in the case when $\gamma_1 = \gamma_2$

Assume instead that $\gamma_1 \neq \gamma_2$. Note that in this case we have $c_0 \leq \lvert \gamma_1 - \gamma_2 \rvert \leq C_0$ for universal constants $c_0$ and $C_0$, since $\gamma_1$ and $\gamma_2$ are machine elements in $\mathbb{Z}\left[ \frac{1}{p} \right]$.\footnote{The constant $c_0$ here depends only on $m$,$n$,$p$, but it may depend sensitively on the approximation of $\frac{1}{p}$ by rationals with low denominators.} 

We define a monomial function $m(\delta) := \frac{m_1(\delta)}{m_2(\delta)} = a \cdot \delta^\gamma$, where $a = \frac{a_1}{a_2}$ and $\gamma = \gamma_1 - \gamma_2$. The unique solution to the equation $m(\delta) = 1$ is given by
\begin{equation}
\label{deltasol}
\delta_{\text{sol}} := a^{\frac{1}{\gamma}}.
\end{equation}
Since monomial functions are monotonic, we have either
\begin{enumerate}[(a)]
\item $m(\delta) < 1$ for $\Delta_g \leq \delta < \delta_{\text{sol}}$, and $m(\delta) > 1$ for $\delta_{\text{sol}} < \delta \leq \Delta_g^{-1}$; or
\item $m(\delta) > 1$ for $\Delta_g \leq \delta <  \delta_{\text{sol}}$, and $m(\delta) < 1$ for $\delta_{\text{sol}} < \delta \leq \Delta_g^{-1}$.
\end{enumerate}
We know that (a) holds if $\gamma > 0$, and (b) holds if $\gamma < 0$. We  can determine which case occurs because the rational number $\gamma$ is given to exact precision. 

Since $a \in [\Delta_g^C,\Delta_g^{-C}]$ and $c_0 \leq \lv \gamma \rv \leq C_0$, due to the numerical stability of exponentiation we can compute a machine number $\delta_*$ such that
\begin{equation}
\label{appx123}
\delta_* \in [\Delta_g^C,\Delta_g^{-C} ]  \; \mbox{and} \; \lv \delta_* - \delta_{\text{sol}}\rv \leq \Delta_g^{-C} \Delta_\epsilon \;\; (\mbox{see \eqref{deltasol}}).
\end{equation}
From \eqref{appx123} and the Lipschitz continuity of $m(\delta) $, we have $\lv m(\delta) - 1 \rv \leq \Delta_g^{-C} \Delta_\epsilon$ for all $\delta$ in the interval between $\delta_*$ and $\delta_{\text{sol}}$. 

Therefore, in case (a) we have $m(\delta) = \frac{m_1(\delta)}{m_2(\delta)} \leq 1 + \Delta_g^{-C} \Delta_\epsilon$ for all $\Delta_g \leq \delta \leq \delta_*$, and $m(\delta) = \frac{m_1(\delta)}{m_2(\delta)} \geq 1 - \Delta_g^{-C} \Delta_\epsilon$ for all $\delta_* \leq \delta \leq \Delta_g^{-1}$. Note that both $m_1(\delta)$ and $m_2(\delta)$ are in the range $[\Delta_g^C,\Delta_g^{-C}]$ if $\delta \in [\Delta_g,\Delta_g^{-1}]$. Thus, in case (a) we determine that 
\[
m_1(\delta) \leq m_2(\delta) \cdot (1 + \Delta_g^{-C} \Delta_\epsilon) \leq m_2(\delta) + \Delta_g^{-2C} \Delta_\epsilon \leq m_2(\delta) + \Delta_\epsilon^{1/2} \;\;\; \mbox{for} \; \Delta_g \leq \delta \leq \delta_*,
\]
and similarly,  $m_1(\delta) \geq m_2(\delta) - \Delta_\epsilon^{1/2}$ for $\delta_* \leq \delta \leq \Delta_g^{-1}$. Thus, we can respond as follows:
\begin{itemize}
\item If $\delta_*$ is to the left of the interval $I$ then outcome (2) occurs.
\item If $\delta_*$ is to the right of the interval $I$ then outcome (1) occurs.
\item If $\delta_*$ belongs to the interval $I$ then outcome (3) occurs with $j=1$ and $k  = 2$.
\end{itemize}

Similarly, in case (b) we determine that $m_1(\delta) \geq m_2(\delta) -  \Delta_\epsilon^{1/2}$  for all $\Delta_g \leq \delta \leq \delta_*$, and similarly,  $m_1(\delta) \leq m_2(\delta) +  \Delta_\epsilon^{1/2}$ for all $\delta_* \leq \delta \leq \Delta_g^{-1}$. Thus, we can respond as follows:
\begin{itemize}
\item If $\delta_*$ is to the left of the interval $I$ then outcome (1) occurs.
\item If $\delta_*$ is to the right of the interval $I$ then outcome (2) occurs.
\item If $\delta_*$ belongs to the interval $I$ then outcome (3) occurs with $j=2$ and $k  = 1$.
\end{itemize}

That completes the explanation of the procedure \textsc{Process Monomials}. Clearly, the work and storage are at most $C$ for a universal constant $C$.

\end{proof}

We return to the setting before the above procedure. 

Fix $\ell \in \{1,\cdots,\ell_{\max}\}$. Applying the procedure \textsc{Process Monomials}, for each pair $(\cA',\cA'')$ such that $\cA' \leq \cA$ and $\cA'' \leq \cA$,  we produce one of three outcomes.

In outcome (1), we guarantee  that $\eta_*^{(\hQ,\cA')} \leq \eta_*^{(\hQ,\cA'')} + \Delta_\epsilon^{1/2}$, uniformly on the interval $I_\ell$. 

In outcome (2), we guarantee  that $\eta_*^{(\hQ,\cA'')} \leq \eta_*^{(\hQ,\cA')} + \Delta_\epsilon^{1/2}$, uniformly on the interval $I_\ell$.

In outcome (3), we divide the interval $I_\ell$ at the point $\delta_{\ell,\cA',\cA''} =  \delta_* \in I_\ell$ to obtain \emph{split subintervals} $I_\ell^{- } = I_\ell \cap (0,\delta_*]$ and $I_\ell^+ = I_\ell \cap ( \delta_*, \infty)$. (A subinterval may contain only a single point or be empty.) We guarantee that $\eta_*^{(\hQ,\cA')} \leq \eta_*^{(\hQ,\cA'')} + \Delta_\epsilon^{1/2}$ on one of the split subintervals, and $\eta_*^{(\hQ,\cA'')} \leq \eta_*^{(\hQ,\cA')} + \Delta_\epsilon^{1/2}$ on the other.  We determine which inequality is satisfied on each subinterval.

For each pair $(\cA',\cA'')$ such that outcome (3) occurs, we have computed a machine number $\delta_{\ell,\cA',\cA''}$ in $I_\ell$. We sort these numbers and remove duplicates to obtain a list
\[
\delta_1 < \delta_2 < \cdots < \delta_{K_\ell}.
\]
Note that $K_\ell \leq \# \{ (\cA',\cA'') : \cA' \leq \cA, \; \cA'' \leq \cA\} \leq C$ for a universal constant $C$. We define $\delta_0$ and $\delta_{K_\ell+1}$ to be the left and right endpoints of $I_\ell$,  respectively. We let $I_\ell^k := [\delta_k,\delta_{k+1}]$ for each $0 \leq k \leq K_\ell$. We thus obtain a (possibly trivial) partition of $I_\ell$ into subintervals $I_{\ell}^0,\cdots,I_\ell^{K_\ell}$.

For each interval $I_\ell^k$ and each pair $(\cA', \cA'')$ such that $\cA' \leq \cA$ and $\cA'' \leq \cA$, we guarantee either that $\eta_*^{(\hQ,\cA')} \leq \eta_*^{(\hQ,\cA'')} + \Delta_\epsilon^{1/2}$ on $I_\ell^k$ ($\cA'$ beats $\cA''$ on $I_\ell^k$), or  that $\eta_*^{(\hQ,\cA'')} \leq \eta_*^{(\hQ,\cA')} + \Delta_\epsilon^{1/2}$ on $I_\ell^k$ ($\cA''$ beats $\cA'$ on $I_\ell^k$). To make such a guarantee, we look at the previous outcomes. If outcome (1) occurs, then $\cA'$ beats $\cA''$ on $I_\ell^k$. If outcome (2) occurs, then  $\cA''$ beats $\cA'$ on $I_\ell^k$. If outcome (3) occurs, then we determine which of the split subintervals of $I_\ell$ contains $I_\ell^k$. Once we have done that, we can make a correct guarantee by using the guarantee made in outcome (3) for the split subinterval.

For each of the intervals $I_\ell^k$ we perform the following computation. We initialize $\cS = \{ \cA' \subset \cM : \cA' \leq \cA \}$. We initialize $\overline{\cA}$ to be any member of $\cS$. Then we run the following loop.
\begin{itemize}
\item \textsc{While}:  $\cS \neq \{ \overline{\cA} \}$
\item
\begin{itemize}
\item Select an arbitrary $\cA' \in \cS \setminus \{ \overline{\cA}\}$. 
\item If we guarantee that $\cA'$ beats $\overline{\cA}$ on $I_\ell^k$, then discard $\overline{\cA}$ from $\cS$ and set $\overline{\cA} = \cA'$.
\item If we guarantee that $\overline{\cA}$ beats $\cA'$ on $I_\ell^k$, then discard $\cA'$ from $\cS$. Do not modify $\overline{\cA}$.
\item (Note that we make at most one guarantee.)
\end{itemize}
\end{itemize} 
Let $\cA_1$ denote the sole member remaining in $\cS$ once the loop is complete. For any $\cA'  \subset \cM$ with $\cA' \leq \cA$, there is a sequence of ``competitors'' $\cA_2,\cdots,\cA_J$ with $\cA_J = \cA'$, such that $\cA_j$ beats $\cA_{j+1}$ on $I_\ell^k$ for $j=1,\cdots, J - 1$. This is clear because $\cA'$ is selected in the loop at some iteration, and as long as $\cA' \neq \cA_1$ we can be certain that $\cA'$ is beaten by some competitor, who in turn is beaten by another competitor, and so on until the loop terminates with the final competitor $\cA_1$. Clearly, the number of competitors $J$ is bounded by a universal constant $C$. By combining the estimates coming from each competition, we learn that $\eta^{(\hQ,\cA_1)}_* \leq \eta^{(\hQ,\cA_J)}_* + J \Delta_\epsilon^{1/2}$. Therefore, $\eta^{(\hQ,\cA_1)}_* \leq \eta^{(\hQ,\cA')}_* + C \Delta_\epsilon^{1/2}$.

Thus, for each $0 \leq k \leq K_\ell$, we can compute a multiindex set $\cA_\ell^k \leq \cA$ such that
\begin{equation}
\label{nearmax}
\eta_*^{(\hQ,\cA_{\ell}^k)}(\delta) \leq \eta_*^{(\hQ,\cA')}(\delta) + C \Delta_\epsilon^{1/2} \qquad \mbox{for all} \; \delta \in I_\ell^k, \;\; \mbox{for all} \; \cA' \leq \cA.
\end{equation}

We repeat the previous construction for each $\ell\in \{1,\cdots,\ell_{\max}\}$.

We therefore obtain machine intervals $I_\ell^k$  ($0 \leq k  \leq K_{\ell}$, $1 \leq \ell \leq \ell_{\max}$), which form a partition of $[\Delta_g,\Delta_g^{-1}]$, and multiindex sets $\cA_\ell^k$ as in \eqref{nearmax}.

We define a function $\widetilde{\eta} : \left[ \Delta_g,\Delta_g^{-1} \right] \rightarrow \R$ by
\[
\widetilde{\eta}(\delta) :=  \eta_*^{(\hQ,\cA_\ell^k)}(\delta) = c_{\ell,\cA_\ell^k} \cdot \delta^{\lambda_{\ell,\cA_\ell^k}} \qquad \mbox{if} \; \delta \in I_\ell^k.
\]

Since $\ell_{\max} \leq C$, the previous construction can be executed using work and storage at most a universal constant $C'$.

We will make use of the properties \textbf{(P1)}-\textbf{(P4)} of the functions $\eta_*^{(\hQ,\cA')}$ that were stated earlier in this section.

Recall that $\eta(\delta)$ is the minimum of $\eta_*^{(\hQ,\cA')}(\delta)$ over all $\cA' \leq \cA$ (see \eqref{defnofeta}). Since $\cA_\ell^k \leq \cA$ for all $(k,\ell)$, we have $\widetilde{\eta}(\delta) \geq  \eta(\delta)$. Moreover, taking the minimum with respect to $\cA'$ in  \eqref{nearmax}, we conclude that $\widetilde{\eta}(\delta) \leq \eta(\delta) + C\Delta_\epsilon^{1/2}$. Thanks to \textbf{(P4)}, we have $\eta(\delta) \geq \Delta_g^C \geq C \Delta_\epsilon^{1/2}$. Thus, we learn that
\begin{equation}\label{comp_fin}
\eta(\delta) \leq \widetilde{\eta}(\delta) \leq 2 \cdot  \eta(\delta).
\end{equation}

We next prove that the the function $\eta^{(\hQ,\cA)}(\delta) = \widetilde{\eta}(\delta)$ satisfies \textbf{(A1)}-\textbf{(A4)}.

\noindent\underline{\textbf{Proof of (A1).}}

Let $\delta \in [\Delta_g,\Delta_g^{-1}]$. Also, let $\eta > C \cdot \widetilde{\eta}(\delta)$, with $C$ as in \textbf{(P1)}. Then, thanks to \eqref{comp_fin}, we have
\[
\eta > C \cdot \eta(\delta) = C \cdot \min_{\cA' \leq \cA} \eta_*^{(\hQ,\cA')}(\delta).
\]
Hence, $\eta > C \cdot  \eta_*^{(\hQ,\cA')}(\delta)$ for some $\cA' \leq \cA$. According to \textbf{(P1)}, we learn that $\ooline{\sigma}(\hQ)$ has an $(\cA',x_\hQ,\eta^{1/2},\delta)$-basis. This completes the proof of \textbf{(A1)}.

\noindent\underline{\textbf{Proof of (A2).}}

Let $\delta \in [\Delta_g,\Delta_g^{-1}]$. Also, let $\eta < \frac{c}{2} \cdot \widetilde{\eta}(\delta)$, with $c > 0$ as in \textbf{(P2)}. Then, thanks to \eqref{comp_fin}, we have
\[
\eta  \leq c \cdot \eta(\delta) = c \cdot \min_{\cA' \leq \cA} \eta_*^{(\hQ,\cA')}(\delta).
\]
Hence, $\eta < c \cdot  \eta_*^{(\hQ,\cA')}(\delta)$ for all $\cA' \leq \cA$. According to \textbf{(P2)}, we learn that $\ooline{\sigma}(\hQ)$ does not have an $(\cA',x_\hQ,\eta^{1/2},\delta)$-basis for any $\cA' \leq \cA$. This completes the proof of \textbf{(A2)}.

\noindent\underline{\textbf{Proof of (A3).}}

Let $\delta_1,\delta_2 \in [\Delta_g,\Delta_g^{-1}]$, with $\frac{1}{10} \delta_1 \leq \delta_2 \leq 10 \delta_1$.

 According to \textbf{(P3)}, for each $\cA' \leq \cA$ we have
\[
 c \cdot \eta_*^{(\hQ,\cA')}(\delta_1) \leq \eta_*^{(\hQ,\cA')}(\delta_2) \leq  C \cdot \eta_*^{(\hQ,\cA')}(\delta_1).
\]
Taking the minimum with respect to $\cA' \leq \cA$, we learn that $c \cdot \eta(\delta_1) \leq \eta(\delta_2) \leq C \cdot \eta(\delta_1)$. According to \eqref{comp_fin}, we therefore have $\frac{1}{4} c \cdot \widetilde{\eta}(\delta_1) \leq \widetilde{\eta}(\delta_2) \leq 4C \cdot \widetilde{\eta}(\delta_1)$. This completes the proof of \textbf{(A3)}.

\noindent\underline{\textbf{Proof of (A4).}}

We have
\[
\widetilde{\eta}(\delta) \ogeq{\eqref{comp_fin}} \eta(\delta) = \min_{\cA' \leq \cA} \eta_*^{(\hQ,\cA')}(\delta) \ogeq{\textbf{(P4)}} \Delta_g^{C}.
\]
This completes the proof of \textbf{(A4)}.

Thus, properties \textbf{(A1)}, \textbf{(A2)}, \textbf{(A3)}, \textbf{(A4)} are satisfied for the function  $\eta^{(\hQ,\cA)}(\delta) = \widetilde{\eta}(\delta)$.

This concludes the explanation of the algorithm.

\end{proof}

\section{Computing Lengthscales}

Each point $x \in E$ is assumed to be an $\overline{S}$-bit machine point. Recall that $\Delta_0 = 2^{- \overline{S}}$. Hence, 
\begin{equation}
\label{Econd}
\lv x' - x'' \rv \geq \Delta_0  \;\; \text{for distinct} \; x',x'' \in E.
\end{equation}

Recall that $\CZ(\cA^-)$ consists of disjoint dyadic cubes that form a partition of $Q^\circ = [0,1)^n$. According to the Main Technical Results for $\cA^-$, we have $\delta_Q \geq c \cdot \Delta_0$ for each $Q$ in $\CZ(\cA^-)$, for a universal constant $c$. Therefore, each $Q$ in $\CZ(\cA^-)$ is an $\til{S}$-bit machine cube, where $\til{S} \leq C \overline{S}$ for a universal constant $C$.

Recall that a testing cube is a dyadic cube $\hQ \subset Q^\circ$ that can be written as a disjoint union of  cubes in $\CZ(\cA^-)$. We then have $\delta_{\hQ} \geq c \cdot \Delta_0$ for a universal constant $c > 0$ (see Remark \ref{test_fp}).

We set $\lambda := 1/40$.

\environmentA{Algorithm: Compute Interesting Cubes (Finite-Precision)}

We compute a tree $T$ consisting of testing cubes. The nodes in $T$ consist of all the cubes $Q \in \CZ(\cA^-)$ that contain points of $E$, all the testing cubes $\hQ$ for which $\diam(3 \hQ \cap E) \geq \lambda \cdot \delta_{\hQ}$, and the unit cube $Q^\circ$. 

Here, $T$ is a tree with respect to inclusion. We mark each internal node $Q$ in $T$ with pointers to its children, and we mark each node $Q$ in $T$ (except for the root) with a pointer to its parent.

The number of nodes in $T$ is at most $C N$, and $T$ can be  computed with work at most $C N \log N$ in space $C N$.

\begin{proof}[\underline{Explanation:}]
We use the explanation in Section \ref{sec_cl}. We need to check that the computation is valid in our finite-precision model of computation.

We compute representative pairs from the well-separated pairs decomposition of $E$ using the algorithm \textsc{Make WSPD} (see Section \ref{sec_CK}). The representative pairs $(x_\nu',x_\nu'') \in E \times E \setminus \{ (x,x) : x \in E \}$ ($1 \leq \nu \leq \nu_{\max}$) satisfy $\lv x_\nu' - x_{\nu}'' \rv \geq \Delta_0$, thanks to \eqref{Econd}. 

Next, we loop over all $\nu$ and list all the dyadic cubes $\widetilde{Q}$ with $x_\nu' , x_\nu'' \in 5 \widetilde{Q}$ and $\lv x_\nu' - x_\nu'' \rv \geq \frac{\lambda}{2} \delta_{\widetilde{Q}}$. We call this list $Q_1,\cdots,Q_K$. Since $5Q_k$ contains some representative pair $(x_\nu',x_\nu'')$, we have $\delta_{Q_k} \geq \frac{1}{5} \lv x_\nu' - x_\nu'' \rv \geq \frac{1}{5} \Delta_0$ for each $k=1,\cdots,K$.

Note that the ``BBD Tree algorithm'' in Theorem \ref{bbd_thm} is unchanged in finite-precision, so we can compute $\diam(3Q_k \cap E)$ for each $k=1,\cdots,K$. We remove any cubes from our list that satisfy $\diam(3 Q_k \cap E) < \lambda \delta_{Q_k}$, which occurs if and only if  $\delta_{Q_k} > 40 \cdot \diam(3 Q_k \cap E)$. We also compute the cube in $\CZ(\cA^-)$ that contains the center of each $Q_k$, using the $\CZ(\cA^-)$-\textsc{Oracle}. If $Q_k$ is strictly contained in this cube, then we remove $Q_k$ from our list. Denote the surviving cubes by $\widetilde{Q}_1,\cdots,\widetilde{Q}_{\widetilde{K}}$.

We list all the cubes $Q \in \CZ(\cA^-)$ that contain points of $E$ (take all the cubes $Q$ in $\CZ_{\main}(\cA^-)$ that satisfy $E \cap Q \neq \emptyset$), the cubes $\widetilde{Q}_1,\cdots,\widetilde{Q}_{\widetilde{K}}$, and the unit cube $Q^\circ$. We sort this list to remove duplicates, and organize it in a tree $T$ using the algorithm \textsc{Make Forest} (see Section \ref{maketree_sec}).

That completes the explanation of the algorithm.
\end{proof}

\environmentA{Algorithm: Compute Critical Testing Cubes (Finite-Precision).}

Given $\epsilon > 0$, which is less than a small enough universal constant, we produce a collection $\widehat{\mathcal{Q}}_\epsilon$ of testing cubes with the following properties.
\begin{enumerate}[(a)]
\item Each point $x \in E$ belongs to some cube $\hQ_x \in \widehat{\mathcal{Q}}_\epsilon$.
\item The cardinality of $\widehat{\mathcal{Q}}_\epsilon$ is at most $C \cdot N$.
\item If $\widehat{Q} \in \widehat{\mathcal{Q}}_\epsilon$ strictly contains a cube in $\CZ(\cA^-)$, then $(1+a(\cA))\widehat{Q}$ is tagged with $(\cA, \epsilon^\kappa)$.
\item If $\widehat{Q} \in \widehat{\mathcal{Q}}_\epsilon$ and $\delta_{\widehat{Q}} \leq c^*$, then no cube containing $S \widehat{Q}$ is tagged with $(\cA, \epsilon^{1/\kappa})$.
\item Each cube $Q$ in $\mathcal{Q}_\epsilon$ satisfies $\delta_Q \geq c \cdot \Delta_0$.
\end{enumerate}

The algorithm requires work at most $C N \log N$ in space $CN$.

Here, $c^* > 0$ and $S \geq 1$ are integer powers of $2$, which depend only on $m$, $n$, $p$; also, $\kappa \in (0,1)$ and $C \geq 1$ are universal constants.

\begin{proof}[\underline{Explanation}]

The main change to the explanation is that we use the finite-precision version of \textsc{Optimize Basis} instead of the infinite-precision version. We also need to show that the roundoff errors that can arise have little effect.

Note that condition (e) will hold for each $Q$ in $\mathcal{Q}_\epsilon$, since we promise that $\mathcal{Q}_\epsilon$ contains only testing cubes. (See Remark \ref{test_fp}.)

We let $\Lambda \geq 1$ be a sufficiently large integer power of two, as before. We will later choose $\Lambda$ to be bounded by a universal constant, but not yet. We assume that $\Lambda$ is a machine number.

We construct a tree $T$ of interesting cubes with the algorithm \textsc{Compute Interesting Cubes}.

We next explain the construction of the collection $\widehat{\cQ}_\epsilon$.

We proceed with Steps 0-6. The construction is almost identical to that in infinite-precision. We refer the reader to the earlier text. We will only record the necessary changes

We assume we have carried out the one-time work of the BBD Tree in Section \ref{sec_bbd}. Thus, given an $\til{S}$-bit machine cube $Q$, with $\til{S} \leq C \overline{S}$, we can compute $\# \left( \frac{65}{64}Q \cap E \right)$ using work at most $C \cdot \log N$.

Therefore, we can compute $\# \left( \frac{65}{64}Q \cap E \right)$ for each $Q$ in $T$.

For each cube $Q_1$ in $T$ we perform Steps 0-3.

Step 0 is unchanged: We find the parent $Q_2$ of $Q_1$ in $T$.

In Step 1 in the infinite-precision text, it says ``We determine whether or not there exists a number $\delta \in [ \Lambda^{10} \delta_{Q_1}, \Lambda^{-10} \delta_{Q_2}]$ with the property that $\epsilon^{1/\kappa_5} \leq \eta^{(\upQ_1,\cA)}(\delta) \leq \epsilon^{\kappa_5}$. If such a $\delta$ exists, we can easily find one.'' We can no longer make such an accurate determination because of inevitable roundoff errors. We will need to make the modifications listed below.

\begin{itemize}
\item \underline{Step 1 (Modified) :} In finite-precision, we compute a piecewise-monomial representation for the function $\eta^{(Q^\up_1,\cA)}(\delta)$ using  the finite-precision version of \textsc{Optimize Basis}, where $Q_1^\up$ is the dyadic cube with $Q_1 \subset Q_1^\up$ and $\delta_{Q_1^\up} = \Lambda \cdot \delta_Q$. We produce one of two outcomes. Either we guarantee that there does not exist a $\delta \in [ \Lambda^{10} \delta_{Q_1}, \Lambda^{-10} \delta_{Q_2}]$ such that 
\begin{equation}
\label{step1neg_fin}
\epsilon^{1/\kappa_5} \leq \eta^{(\upQ_1,\cA)} (\delta) \leq \epsilon^{\kappa_5},
\end{equation}
or else we compute a machine number $\delta \in [ \Lambda^{10} \delta_{Q_1}, \Lambda^{-10} \delta_{Q_2}]$ satisfying 
\begin{equation}
\label{step1_fin}
\frac{1}{2} \epsilon^{1 / \kappa_5} \leq \eta^{(\upQ_1,\cA)}(\delta) \leq 2 \epsilon^{\kappa_5}.
\end{equation}
The number $\delta$ is computed \underline{exactly}.

The factors of $2$ in the above estimate arise because of roundoff errors in the computation of $\delta$. Indeed, we can bound any roundoff error by $\Delta_\epsilon \Delta_g^{-C}$, which is at most $100^{-1} \cdot \epsilon^{1/\kappa_5}$, since $\Delta_\epsilon \Delta_g^{-C} \leq \Delta_\epsilon^{1/2}$ (see \eqref{constants2}) and $\Delta_\epsilon^{1/2} \leq \Delta_\new \leq 100^{-1} \cdot \epsilon^{1/\kappa_5}$ (see \eqref{constants3} and \eqref{constants5}).

As before, in the second alternative we can find a dyadic cube $Q$ with $Q_1 \subset Q \subset Q_2$, $\Lambda^{10} \delta_{Q_1} \leq \delta_Q \leq \Lambda^{-10} \delta_{Q_2}$, and such that 
\begin{equation} \label{step1_fin}\left[ \epsilon^{1/\kappa_6} \leq \eta^{(\upQ_1,\cA)}( \delta_Q)  \right]\;\; \mbox{and} \;\; \left[ \eta^{(\upQ_1,\cA)}(\delta_Q) \leq \epsilon^{\kappa_6} \right].
\end{equation}
(Compare to \eqref{step1}.) 

Here, by choosing $\kappa_6$ sufficiently small, we can make the extra factors of $2$ disappear.

In the second alternative, we add $Q$ to the collection $\hQ_\epsilon$. That completes the computation in Step 1.

\end{itemize}
Note that $[ \delta_{Q_1}, \delta_{Q_2}] \subset [\Delta_g,\Delta_g^{-1}]$, since each cube in $T$ has sidelength in $[ c \cdot \Delta_0, 1]$, and since $\Delta_g \leq c \cdot \Delta_0$. This comment justifies the previous computation, since the function $\eta^{(Q_1^\up,\cA)}(\delta)$ is defined only for $\delta \in [\Delta_g,\Delta_g^{-1}]$.

Similarly, in Steps 2 - 6, we make the following changes.

\begin{itemize}
\item \underline{Step 2 (Modified) :} We examine each dyadic cube $Q$ with $Q_1 \subset Q \subset Q_2$, $\delta_Q \leq \Lambda^{-10}$, and [$\delta_Q \leq \Lambda^{10} \delta_{Q_1}$ or $\delta_Q \geq \Lambda^{-10} \delta_{Q_2}$]. We compute a piecewise-monomial function $\eta^{(Q^\up,\cA)}(\delta)$ using the finite-precision version of \textsc{Optimize Basis}. We produce one of two outcomes. Either we guarantee that 
\begin{equation} \label{step2neg_fin}
\left[ \epsilon^{1/\kappa_5} >  \eta^{(Q^\up,\cA)}(\delta_{Q^\up}) \right] \; \mbox{or} \; \left[ \#\left( \frac{65}{64} Q \cap E\right) \geq 2 \;\; \mbox{and} \;\; \eta^{(Q,\cA)}(\delta_Q) > \epsilon^{\kappa_5} \right],
\end{equation}
where $Q^\up$ is the unique dyadic cube with $Q \subset Q^\up$ and $\delta_{Q^\up} = \Lambda \delta_Q$, or else we guarantee that
\begin{equation}
\label{step2_fin}
\left[ \frac{1}{2} \epsilon^{1/\kappa_5} \leq  \eta^{(Q^\up,\cA)}(\delta_{Q^\up}) \right] \; \mbox{and} \; \left[ \#\left( \frac{65}{64} Q \cap E\right) \leq 1 \;\; \mbox{or} \;\; \eta^{(Q,\cA)}(\delta_Q) \leq 2 \epsilon^{\kappa_5} \right].
\end{equation}
The extra factors of $2$ allow for roundoff errors in the computation of $\eta^{(Q^\up,\cA)}(\delta_{Q^\up})$.

We add $Q$ to the collection $\hQ_\epsilon$ in the second alternative.

\item \underline{Step 3 (Modified) :} We examine each dyadic cube $Q$ with $Q_1 \subset Q \subset Q_2$ and $\delta_Q \geq \Lambda^{-10}$. We apply the finite-precision version of \textsc{Optimize Basis} to compute a function $\eta^{(Q,\cA)}(\delta)$. We produce one of two outcomes. Either we guarantee that 
\begin{equation} \label{step3neg_fin}
\left[ \#\left( \frac{65}{64} Q \cap E\right) \geq 2 \;\; \mbox{and} \;\; \eta^{(Q,\cA)}(\delta_Q) > \epsilon^{\kappa_5} \right], 
\end{equation}
or else we guarantee that
\begin{equation} \label{step3_fin}
\left[ \#\left( \frac{65}{64} Q \cap E\right) \leq 1 \;\; \mbox{or} \;\; \eta^{(Q,\cA)}(\delta_Q) \leq 2 \epsilon^{\kappa_5} \right].
\end{equation}
The extra factors of $2$ allow for roundoff errors in the computation of $\eta^{(Q,\cA)}(\delta_{Q})$.

We add $Q$ to the collection $\hQ_\epsilon$ in the second alternative.

\item \underline{Step 4 (Modified) :} We apply the finite-precision version of \textsc{Optimize Basis} to compute a function $\eta^{(Q^\circ,\cA)}(\delta)$. We produce one of two outcomes. Either we guarantee that 
\begin{equation} \label{step4neg_fin}
\left[ \eta^{(Q^\circ,\cA)}(\delta_{Q^\circ}) > \epsilon^{\kappa_5} \right],
\end{equation}
or else we guarantee that
\begin{equation} \label{step4_fin}
\left[ \eta^{(Q^\circ,\cA)}(\delta_{Q^\circ}) \leq 2 \epsilon^{\kappa_5} \right].
\end{equation}
We add $Q^\circ$ to the collection $\hQ_\epsilon$ in the second alternative.

\item \underline{Step 5 (Modified) :} We examine all dyadic cubes $Q \subset Q^\circ$ such that $\delta_Q \geq \Lambda^{-10}$. We add $Q$ to the collection $\widehat{\mathcal{Q}}_\epsilon$ if and only if $Q \in \CZ(\cA^-)$.

\item \underline{Step 6 (Modified) :} We examine all cubes $Q \in \CZ(\cA^-)$ such that $\delta_Q \leq \Lambda^{-10}$ and $Q \cap E \neq \emptyset$. We apply the finite-precision version of \textsc{Optimize Basis} to compute a function $\eta^{(\upQ,\cA)}(\delta)$, where $Q^\up$ is the dyadic cube with $Q \subset Q^\up$ and $\delta_{Q^\up} = \Lambda \delta_Q$. We produce one of two outcomes. Either we guarantee that 
\begin{equation} \label{step6neg_fin}
\left[ \epsilon^{1/\kappa_5} > \eta^{(\upQ,\cA)}(\delta_{\upQ}) \right], 
\end{equation}
or else we guarantee that
\begin{equation} \label{step6_fin}
\left[ \frac{1}{2} \epsilon^{1/\kappa_5} \leq \eta^{(\upQ,\cA)}(\delta_{\upQ}) \right].
\end{equation}
We add $Q$ to the collection $\hQ_\epsilon$ in the second alternative.
\end{itemize}

As before, we see that $\# ( \widehat{\cQ}_\epsilon) \leq C(\Lambda) \cdot N$, hence property (b) holds.

Recall that Propositions \ref{tool1},\ref{tag_prop1}, and \ref{tag_prop2} are unchanged in the finite-precision case - only their proofs required modification. Hence, the analysis that the above algorithm works proceeds as before. In place of the conditions \eqref{step1}, \eqref{step2}, \eqref{step3}, \eqref{step4}, and \eqref{step6} we use the conditions presented in the above bullet points. 

The proof of properties (c) and (d) requires minor changes to reflect the loss of factors of $2$. By choosing smaller values for $\kappa_1,\cdots,\kappa_{20}$, we arrange that the extra factors of $2$ can be absorbed into relevant estimates in the proof. Thus, we can prove properties (c) and (d) for each cube in $\cQ_\epsilon$ using the same argument as before.

The proof of property (a) requires minor changes to reflect the loss of factors of $2$.

We fix a point $x \in E$. 

As before, we consider the increasing chain of cubes $Q_0 \subset Q_1 \subset \cdots \subset Q_{\nu_{\max}}$ in $T$, such that $Q_{\ell+1}$ is a parent of $Q_\ell$ in $T$, $x \in Q_0$, and $Q_0 \in \CZ(\cA^-)$.

As before, we consider the First Extreme Case, the Second Extreme Case, and the Main Case.

To prove (a), we will show that there exists a cube $Q \in \widehat{\cQ}_\epsilon$ such that $x \in Q$.

In the \underline{First Extreme Case}: We assume that $3 Q^\circ$ is tagged with $(\cA,\epsilon)$ and deduce that $\eta^{(Q^\circ,\cA)}(\delta_{Q^\circ}) \leq \epsilon^{\kappa_5}$. Hence, according to the above construction in Step 4, we included $Q^\circ$ in $\widehat{\cQ}_\epsilon$.

In the \underline{Second Extreme Case}: We assume that $3 Q_0$ is not tagged with $(\cA,\epsilon)$ and we deduce that $\eta^{(\upQ_0,\cA)}(\delta_{\upQ_0}) \geq \epsilon^{1/\kappa_5}$. Hence, according to the construction in Step 6, we included $Q_0$ in $\widehat{\cQ}_\epsilon$. 

In the \underline{GI subcase} of the \underline{Main Case}: From the assumptions in the GI subcase we prove \eqref{eq01} and \eqref{eq02} (see the analysis in infinite-precision). This means that \eqref{step2neg_fin} does not hold for the cube $Q$. Hence, according to the construction in Step 2, we included $Q$ in $\widehat{\cQ}_\epsilon$. 

In the \underline{GUI subcase} of the \underline{Main Case}: From the assumptions in the GUI subcase we prove \eqref{eq4} and \eqref{eq5}. Hence, \eqref{step1neg_fin} holds with $\delta = \delta_Q$. Thus, we pass to the second alternative in our construction in Step 1 (for the cube $Q_\nu \in T$). Hence, we decided to include in $\widehat{\cQ}_\epsilon$ a cube $Q'$ with $Q_\nu \subset Q' \subset Q_{\nu+1}$.

In the \underline{NM subcase} of the \underline{Main Case}: From the assumptions in the NM subcase we prove \eqref{eq6}. Hence, \eqref{step3neg_fin} fails to hold for the cube $Q$. Hence, in the construction in Step 3, we included $Q$ in $\widehat{\cQ}_\epsilon$. 

Thus, as in infinite-precision, we see that there exists $Q' \in \widehat{\cQ}_\epsilon$ with $Q_0 \subset Q' \subset Q_{\nu_{\max}}$, and hence $x \in Q'$. This completes the proof of (a).

We choose $\Lambda \geq 1$ to a be a large enough universal constant so that the above holds. That concludes the explanation of the algorithm. 
\end{proof}

According to our construction, each $Q$ in $\mathcal{Q}_\epsilon$ satisfies $\delta_Q \geq c \cdot \Delta_0$. Furthermore, by hypothesis, each $x \in E$ is an $\overline{S}$-bit machine point. 

Thus, we can apply the algorithm \textsc{Placing a Point Inside Target Cuboids} to compute a cube $Q_x \in \mathcal{Q}_\epsilon$ containing each $x \in E$. This requires work at most $C N \log N$ in space $CN$. Thus, the algorithm \textsc{Compute Lengthscales} is unchanged in finite-precision (see Section \ref{sec_cl3}). 

Proposition \ref{lengthscales_prop} still holds in the finite-precision setting.

\section{Passing from Lengthscales to CZ Decompositions}

We explain how to define a decomposition $\CZ(\cA)$ of $Q^\circ$ into machine cubes, and how to define a $\CZ(\cA)$-\textsc{Oracle}.

For each $x \in E$, we compute the machine numbers
\[
\Delta_\cA(x) := \delta_{Q_x}.
\]

We say that a testing cube $Q \subset Q^\circ$ is $OK(\cA)$ if either $Q \in \CZ(\cA^-)$ or $\Delta_{\cA}(x) \geq K \delta_Q$ for all $x \in E \cap 3Q$, where $K := \frac{2^{30}}{a(\cA)}$ (here, the constant $10^{9}$ in Section \ref{sec_cz} is replaced by $2^{30}$).

We define a Calder\'on-Zygmund decomposition $\CZ(\cA)$ of the unit cube $Q^\circ$ to consist of the maximal dyadic subcubes $Q \subset Q^\circ$ that are $OK(\cA)$.

Clearly, $\CZ(\cA^-)$ refines the decomposition $\CZ(\cA)$, namely, each cube in $\CZ(\cA)$ is a disjoint union of the cubes in $\CZ(\cA^-)$. Since $\delta_Q \geq \frac{1}{32} \cdot \Delta_0$ for each $Q \in \CZ(\cA^-)$ (by the finite-precision version of the Main Technical Results for $\cA^-$), we have
\begin{equation}
\label{hh111}
\delta_Q \geq \frac{1}{32} \cdot \Delta_0 \quad \mbox{for each} \;  Q \in \CZ(\cA).
\end{equation}
This implies an additional property of $\CZ(\cA)$ that is required in the finite-precision version of Main Technical Results for $\cA$.

We construct a $\CZ(\cA)$-\textsc{Oracle} using the \textsc{Glorified CZ-Oracle} in Section \ref{sec_czdecomp}, where we take $\Delta(x) := \Delta_\cA(x)/K = \Delta_\cA(x) \cdot a(\cA) \cdot 2^{-30}$. Note that $a(\cA) = a_{\new} = 2^{- \til{S}}$, where $\til{S} \leq C \overline{S}$ for a universal constant $C$ (see \eqref{aprop}). Note also that $\Delta_\cA(x) = \delta_{Q_x}$ is an $\til{S}$-bit machine number (recall that $Q_x$ is a testing cube, and use Remark \ref{test_fp}). Thus, $\Delta(x)$ is an $\til{S}$-bit machine number for each $x \in E$, where $\til{S} \leq C' \overline{S}$ for a universal constant $C'$. Thus, the extra hypotheses required for the finite-precision version of the \textsc{Glorified CZ-Oracle}  are valid (see Section \ref{sec_cz_fin}).

The remaining properties \textbf{(CZ1)}-\textbf{(CZ5)} of the decomposition $\CZ(\cA)$  are proven in Section \ref{sec_cz}. See Propositions \ref{prop561}, \ref{prop_gg}, and \ref{mainprops}.

We have thus proven all the properties of the decomposition $\CZ(\cA)$ stated in the Main Technical Results for $\cA$.

\section{Completing the Induction}

In executing the algorithm \textsc{Produce All Supporting Data} in finite-precision, we need to produce extra stuff, since we added stuff to the definition of \emph{modified supporting data} (see \textbf{Modification 1} in Section \ref{supp_data_fin}). For each $Q \in \CZ_{\main}(\cA)$, we need to list all the points $x \in  E \cap \frac{65}{64}Q$. However, it's easy to do that. The procedure is as follows: We loop over all points $x \in E$. For each $x$, we use the $\CZ(\cA)$-\textsc{Oracle} to find all the $Q \in \CZ_{\main}(\cA)$ such that $x \in \frac{65}{64}Q$, and we then add $x$ to a list associated to each relevant $Q$. Any given $x$ is associated to at most $C$ cubes $Q$, and we can find each cube in the list $\CZ_{\main}(\cA)$ by binary search that requires work at most $C \log N$.  Therefore, this procedure requires work at most $C N \log N$ in space $ C N$. Thus, the work and storage used by the finite-precision version of the algorithm \textsc{Produce All Supporting Data} are bounded as required.

In place of \eqref{q1} and \eqref{q2}, we have to prove the estimates.
\[\| (f,P)\|_{(1+a(\cA))\hQ} \leq C M_{\hQ}(f,P)\]
and
\[M_{\hQ}(f,P) \leq C \| (f,P) \|_{\frac{65}{64}\hQ} + C \Delta_\new \lvert P \rvert. \]
We prove these estimates using the finite-precision Unconditional and Conditional Inequalities, just as in the infinite-precision case.

We separately treat the \emph{simple} and \emph{non-simple} cubes $\hQ \in \CZ(\cA)$ in Sections \ref{notsimple_sec} and \ref{simple_sec}. We make a few small changes to the analysis. which are documented below.

\begin{itemize}
\item In Section \ref{notsimple_sec}: We defined lists $\Xi(\hQ,\cA)$ and $\Omega(\hQ,\cA)$ of linear functionals, and a linear map $T_{(\hQ,\cA)}$ for each of the \emph{non-simple cubes} $\hQ \in \CZ(\cA)$. The definitions are unchanged. See the versions of the algorithms \textsc{Compute New Assists}, \textsc{Compute New Assisted Functionals}, and \textsc{Compute New Extension Operator} in Section \ref{supp_data_fin}. The linear functionals and linear maps here are all computed with parameters $(\Delta_g^C, \Delta_g^{-C} \Delta_\epsilon)$.
\item In Section \ref{notsimple_sec}: We need to control an extra sum when evaluating the upper bound on the work and storage. Namely, we have to control the sum
\[\sum_{\hQ \in \CZ_{\main}(\cA)} \left\{ \# \left( \frac{65}{64}\hQ \cap E \right) \right\}.\]
This extra term arises from the work of applying the finite-precision version of \textsc{Compute New Assisted Functionals} (see Section \ref{supp_data_fin}). This sum is bounded by $CN$, thanks to the bounded overlap of the cubes $\frac{65}{64}\hQ$, for $\hQ \in \CZ(\cA)$. Hence, the work and storage needed to compute all the functionals defined in Section \ref{notsimple_sec} are bounded as required.

\item In Section \ref{simple_sec}: We defined lists $\Xi(\hQ,\cA)$ and $\Omega(\hQ,\cA)$, and a linear map $T_{(\hQ,\cA)}$ for each of the \emph{simple cubes} $\hQ \in \CZ(\cA)$.  The definitions are unchanged. See the relevant text. The linear functionals and linear maps here are all computed with parameters $(\Delta_g^C, \Delta_g^{-C} \Delta_\epsilon)$.
\item In Section \ref{simple_sec}: The finite-precision version of \eqref{simple_2} (from the Main Technical Results for $\cA^-$) states that
\begin{equation}
\label{simple_2_fin}
C^{-1}  \| (f,R) \|_{(1+a)Q} \leq M_{(Q,\cA^-)}(f,P) \rvert \leq C  \left[ \| (f,P) \|_{\frac{65}{64}Q} + \Delta_\junk \lvert P \rvert \right].
\end{equation}

\item In Section \ref{simple_sec}: The statement and proof of Proposition \ref{simple_extprops} are unchanged.

\item In Section \ref{simple_sec}: The finite-precision version of Lemma \ref{simple_lem1} states that
\[
C^{-1} \|(f,P) \|_{(1+a_\new)\hQ} \leq M_{(\hQ,\cA)}(f,P) \leq C \left[ \| (f,P) \|_{\frac{65}{64}\hQ} + \Delta_\new \lvert P \rvert \right].
\]
We prove this estimate as follows. From \eqref{simple_2_fin} we have
\[
\left[ M_{(\hQ,\cA)}(f,P) \right]^p \leq C \sum_{\substack{ Q \in \CZ_{\main}(\cA^-) \\  Q \subset (1+t_G)\hQ} }  \left[ \| (f,P) \|_{\frac{65}{64} Q}^p + \Delta_\junk^p \lvert P \rvert^p \right].\]
The number of $Q$ arising in the above sum is bounded by $C$. Hence, $M_{(\hQ,\cA)}(f,P) \leq C \left[ \| (f,P) \|_{\frac{65}{64}\hQ} + \Delta_\new \lvert P \rvert \right]$, just as in the proof of Lemma \ref{simple_lem1} in the infinite-precision setting. Here, we use that $\Delta_{\junk} \leq \Delta_\new$; see \eqref{constants3}.

\item In the \textbf{Closing Remarks}: We fix $\epsilon$ to be a small enough universal constant. The parameters $\Delta_g = \Delta_g(\cA^-)$, $\Delta_\epsilon = \Delta_\epsilon(\cA^-)$, and $\Delta_\new$ are assumed to satisfy \eqref{constants3}, \eqref{constants4} and \eqref{constants5}.\footnote{Recall that we have fixed $t_G$ and $\epsilon$ to be universal constants. Hence, the conditions \eqref{constants4} and \eqref{constants5} state that $\Delta_\new$ is less than a small enough universal constant. These are among the conditions \eqref{constants0} and \eqref{constants1} imposed before. } We also impose the assumptions $\Delta_{\junk}(\cA) \geq \Delta_{\new}$, $\Delta_{g}(\cA) \leq \Delta_g^{C}$, and $\Delta_\epsilon(\cA) \geq \Delta_g^{-C} \Delta_\epsilon$, for a large enough universal constant $C$. Thus, we obtain the Main Technical Results for $\cA$  from the above bullet points. 
\item If $\cA = \emptyset$ (the maximal multiindex set) then the induction is complete. We do not fix a choice of the parameters $\Delta_\epsilon(\cA)$, $\Delta_g(\cA)$, $\Delta_\junk(\cA)$ (for $\cA \subset \cM$) just yet. These parameters are determined later in the proofs of our Main Theorems.
\end{itemize}

\section{Main Theorems}
\label{mainthm_fin_sec}
\subsection{Homogeneous Sobolev spaces}
\label{hom_fp}

In this section we prove Theorem \ref{main_thm_hom_fin} using the Main Technical Results for $\cA = \emptyset$.

We assume we are given parameters $\Delta_{\min} = 2^{- K_{\max} \overline{S}}$, $\Delta_\epsilon^\circ := 2^{ - K_{1} \overline{S} }$, $\Delta_g^\circ := 2^{ - K_{2} \overline{S} }$, and $\Delta_\junk^\circ := 2^{ - K_{3} \overline{S} }$, for integers $K_1,K_2,K_3,K_{\max} \geq 1$ as in Theorem \ref{main_thm_hom_fin}.

The proof is identical to the argument in Section \ref{hom_sec}, except for minor changes, which we describe below.

We start from the sentence ``By translating and rescaling, we may assume $\cdots$,'' which follows the statement of Theorem \ref{main_thm_hom}.

We let the parameters $\Delta_\epsilon = \Delta_\epsilon(\emptyset)$, $ \Delta_g = \Delta_g(\emptyset)$, and $\Delta_\junk =  \Delta_\junk(\emptyset)$ be as in the Main Technical Results for $\cA = \emptyset$.

According to the Main Technical Results for $\cA = \emptyset$, we are given the following objects.

There is a dyadic decomposition $\CZ$ of the unit cube $Q^\circ$. The \textsc{$\CZ$-Oracle} operates as before, except that the query point $\underline{x} \in Q^\circ$ is required to be an $S$-bit machine point. We can list all the cubes $Q \in \CZ$ such that $\underline{x} \in \frac{65}{64}Q$, using work at most $C \log N$.

For each $Q \in \CZ$ with $\frac{65}{64}Q \cap E \neq \emptyset$, we are given a collection $\Omega(Q) \subset \left[ \X(E \cap \frac{65}{64}Q) \right]^*$ of assist functionals, a collection $\Xi(Q) \subset \left[ \X(E \cap \frac{65}{64}Q) \oplus \cP \right]^*$ of assisted functionals, and a linear map  $T_Q :  \X(E \cap \frac{65}{64}Q) \oplus \cP \rightarrow \X$. 

We recall some of the main properties of these objects in the bullet points below.

\begin{itemize}
\item \textbf{Modification 1:} For each $Q \in \CZ$ with $\frac{65}{64}Q \cap E \neq \emptyset$, the linear functionals $\omega \in \Omega(Q)$ are given with parameters $(\Delta_g,\Delta_\epsilon)$; also, the linear functionals $\xi \in \Xi(Q)$ are given in short form with parameters $(\Delta_g,\Delta_\epsilon)$ in terms of the assists $\Omega(Q)$.

Given $Q \in \CZ$ with $\frac{65}{64}Q \cap E \neq \emptyset$, given an $S$-bit machine point $\underline{x} \in Q^\circ$, and given $\alpha \in \cM$, we compute the linear functional $(f,P) \mapsto \partial^\alpha( T_Q(f,P))(\underline{x})$ in short form with parameters $(\Delta_g,\Delta_\epsilon)$ in terms of the assists $\Omega(Q)$.

\item \textbf{Modification 2:} We replace \eqref{mm0} with the corresponding estimate from the finite-precision version of the Main Technical Results for $\cA = \emptyset$, namely:
\begin{equation}
\label{mm0_fin}
\sum_{\xi \in \Xi(Q)} \lvert \xi(f,P) \rvert^p \leq C \left[ \|(f,P)\|_{\frac{65}{64}Q}^p + \Delta_\junk^p \lvert P \rvert^p \right].
\end{equation}
\item The linear maps $T_Q$ satisfy \eqref{m1} and \eqref{m2} just as before.
\item From the conditions in the Main Technical Results we learn that $\delta_Q > c_*$ for every $Q \in \CZ$, for the universal constant $c_* = c_*(\emptyset)$. Using the \textsc{$\CZ$-Oracle}, we can list all the cubes in $\CZ$ using work at most $C \log N$. The algorithm is as before.

\item \textbf{Modification 3:} As before, we let $a$ denote the universal constant $a(\emptyset)$. According to the finite-precision version of the Main Technical Results, we know that $a$ is an integer power of $2$. Thus, $a$ is a machine number. We define a family of cutoff functions $\widetilde{\theta}_Q$ (for $Q \in \CZ$) as before. We refer the reader to the text in Section \ref{hom_sec} for a statement of the relevant properties of $\widetilde{\theta}_Q$. The finite-precision version of the algorithm \textsc{Compute Auxiliary Functions} requires slight modification to allow for roundoff errors. Given $Q \in \CZ$ and given an $S$-bit machine point $\underline{x} \in Q^\circ$, we compute the numbers $\partial^\alpha ( \widetilde{\theta}_Q)(\underline{x})$ for all $\alpha \in \cM$. We guarantee that the numbers $\partial^\alpha ( \widetilde{\theta}_Q)(\underline{x})$ have magnitude  at most $\Delta_g^{-C}$ and are computed to precision $\Delta_g^{-C} \Delta_\epsilon$ for a universal constant $C$. For the explanation, we define a spline function $\widetilde{\theta}$ (depending on $a$) with $\widetilde{\theta} \geq 1/2$ on $Q^\circ = [0,1)^n$, $\widetilde{\theta} \equiv 0$ outside $(1+a)Q^\circ$, $0 \leq \widetilde{\theta}_Q \leq 1$ on $\R^n$, and $\lv \partial^\beta \widetilde{\theta}_Q(x) \rv \leq C$ (for $\beta \in \cM$, $x \in \R^n$). We also assume that the derivatives of $\widetilde{\theta}$ at a general $S$-bit machine point in $\R^n$ can be computed to precision $\Delta_\epsilon$. This is possible because the machine precision of our computer is $\Delta_{\min} \ll \Delta_\epsilon$. We define $\widetilde{\theta}_Q$ to be an appropriately shifted and rescaled version of $\widetilde{\theta}$ that is supported on the cube $(1+a)Q$. Since $\delta_Q \geq c^*$ for all $Q \in \CZ$, we learn that $\lv \partial^\beta \widetilde{\theta}_Q(x) \rv \leq C' \leq \Delta_g^{-C'}$ for a large enough universal constant $C'$. We can compute $\partial^\alpha \widetilde{\theta}_Q(\underline{x})$ (for $\alpha \in \cM$) with precision $\Delta_g^{-C} \Delta_\epsilon$ by rescaling the $\alpha$-derivative of $\widetilde{\theta}$ at a suitable machine point in $\R^n$ (determined by  $\underline{x}$).

\item \textbf{Modification 4:} We modify \textsc{Compute POU2} to take into account roundoff errors. Given $Q \in \CZ$ and given an $S$-bit machine point $\underline{x} \in Q^\circ$, we compute the numbers $\partial^\alpha ( \theta_Q)(\underline{x})$ for each $\alpha \in \cM$. The numbers $\partial^\alpha ( \theta_Q)(\underline{x})$ are bounded in magnitude by $\Delta_g^{-C}$ and are computed to precision $\Delta_g^{-C} \Delta_\epsilon$ for a universal constant $C$. Here, $\theta_Q$ is defined in terms of $\widetilde{\theta}_Q$ as in Section \ref{hom_sec}. The explanation is obvious. We choose the function $\eta(t)$ to be a spline function whose derivatives can be computed to precision $\Delta_\epsilon$, and then we compute the derivatives of $\theta_Q$ using the Leibniz rule. Of course, we still have the properties (1)-(4) of the partition of unity $(\theta_Q)$.
\item The definitions of $\Xi^\circ$, $\Omega^\circ$, and $T^\circ$, are unchanged. We define $\Xi^\circ$ to be the union of the lists $\Xi(Q)$, and we define $\Omega^\circ$ to be the union of the lists $\Omega(Q)$. As before, we define $T^\circ(f,P)$ as in \eqref{star2}, namely:
\[
T^\circ(f,P) = \sum_{\substack{ Q \in \CZ \\ \frac{65}{64}Q \cap E \neq \emptyset }} \theta_Q \cdot T_Q(f,P) + \sum_{\substack{ Q \in \CZ \\ \frac{65}{64}Q \cap E = \emptyset }} \theta_Q \cdot P. 
\]

\item \textbf{Modification 5:} The second bullet point in Proposition \ref{mprop} is changed to account for roundoff errors. Given an $S$-bit machine point $\underline{x}  \in Q^\circ$ and given $\alpha \in \cM$, we compute the linear functional $(f,P) \mapsto \partial^\alpha ( T^\circ(f,P))(\underline{x})$ in short form with parameters $(\Delta_g^C,\Delta_g^{-C} \Delta_\epsilon)$ in terms of the assists $\Omega^\circ$. The explanation is an obvious consequence of the Leibniz rule, since the functionals $(f,P) \mapsto \partial^\beta(  T_Q(f,P))(\underline{x})$ can be computed with parameters $(\Delta_g^C,\Delta_g^{-C}\Delta_\epsilon)$, and the numbers $\partial^\beta ( \theta_Q)(\underline{x})$ can be computed with parameters $(\Delta_g^C,\Delta_g^{-C}\Delta_\epsilon)$.
\item \textbf{Modification 6:}  The fourth bullet point of Proposition \ref{mprop} is changed to instead consist of the estimate
\begin{equation}
\label{bullet4}
 \sum_{\xi \in \Xi^\circ} \lvert \xi(f,P) \rvert^p \leq C \cdot \left[ \| (f,P) \|_{\frac{65}{64}Q^\circ}^p + \Delta_\junk^p \lvert P \rvert^p \right].
 \end{equation}
Next, we explain how to modify the proof of Proposition \ref{mprop}.
\item \textbf{Modification 7:} We replace \eqref{snuff22} with
\[ \sum_{\substack{ Q \in \CZ \\ \frac{65}{64}Q \cap E \neq \emptyset}}  \sum_{\xi \in \Xi(Q)} \lvert \xi(f,P) \rvert^p \leq C \cdot \sum_{\substack{ Q \in \CZ \\ \frac{65}{64}Q \cap E \neq \emptyset}} \left[ \| (f,P) \|_{\frac{65}{64}Q}^p + \Delta_\junk^p \lvert P \rvert^p \right],\]
which follows from \eqref{mm0_fin}.

Now, the cardinality of $\CZ$ is at most a universal constant and $\| (f,P) \|_{\frac{65}{64}Q} \leq C \| (f,P) \|_{\frac{65}{64}Q^\circ}$, just as before. Hence, we have
\[\sum_{\substack{ Q \in \CZ \\ \frac{65}{64}Q \cap E \neq \emptyset}}  \sum_{\xi \in \Xi(Q)} \lvert \xi(f,P) \rvert^p \leq C \cdot \left[ \| (f,P) \|^p_{\frac{65}{64}Q^\circ} + \Delta_\junk^p \lvert P \rvert^p \right].\]
But this is just the estimate in the modified fourth bullet point of Proposition \ref{mprop} (see \textbf{Modification 6}). The proof of Proposition \ref{mprop} is otherwise unchanged. This completes the proof of the modified version of Proposition \ref{mprop}.

\item \textbf{Modification 8:}  Now we introduce a linear map $\mathfrak{R} : \X(E) \mapsto \cP$ using the finite-precision version of \textsc{Optimize via Matrix} with $\Delta = \Delta_\junk$. We compute the map $\mathfrak{R}$ in short form with parameters $(\Delta_g,\Delta_\epsilon)$ in the following sense: For each $\alpha \in \cM$, we compute the linear functional $f \mapsto \partial^\alpha ( \mathfrak{R}(f))(0)$ in short form with parameters $(\Delta_g^C,\Delta_g^{-C}\Delta_\epsilon)$ (without assists). We guarantee that
\begin{equation} \label{mm1_fin}
\sum_{\xi \in \Xi^\circ} \lvert \xi(f, \mathfrak{R}(f)) \rvert^p \leq C \cdot \left[ \sum_{\xi \in \Xi^\circ} \lvert \xi(f,R) \rvert^p + \Delta_\junk^p \lvert R \rvert^p \right] \quad \mbox{for any} \; R \in \cP.
\end{equation}
(This estimate is the finite-precision analogue of \eqref{mm1}.) 

We can answer slightly more general queries: Given an $S$-bit machine point $\underline{x} \in Q^\circ$, and given $\alpha \in \cM$, we compute the linear functional $f \mapsto \partial^\alpha (\mathfrak{R}(f))(\underline{x})$. This follows because of Taylor's formula, which allows us to express the functional $f \mapsto \partial^\alpha( \mathfrak{R}(f))(\underline{x})$ as a weighted combination
\[
\sum_{\lv \beta \rv \leq m - 1 - \lv \alpha \rv} \frac{1}{\beta!} \cdot \left( \underline{x} \right)^{\beta} \partial^{\alpha + \beta} (\mathfrak{R}(f))(0)
\]
of the linear functionals $f \mapsto \partial^\gamma (\mathfrak{R}(f))(0)$ ($\gamma \in \cM$). The coefficients in this combination can be computed to precision $(\Delta_g,\Delta_\epsilon)$, and so the claim follows.

\item \textbf{Modification 9:} The list $\Xi$ consists of all the functionals $\xi: f \mapsto \xi^\circ(f,\mathfrak{R}(f))$ with $\xi^\circ \in \Xi^\circ$.  We compute each $\xi \in \Xi$ with parameters $(\Delta_g^C,\Delta_g^{-C}\Delta_\epsilon)$ by composing a linear functional $\xi^\circ \in \Xi^\circ$ with the linear map $f \mapsto \mathfrak{R}(f)$.

\item \textbf{Modification 10:} The cutoff function $\theta^\circ$ is defined as before. The same properties (1)-(4) hold. For the construction, we choose $\theta^\circ$ to be an appropriate spline function. The computation of $\theta^\circ$ is modified to take into account roundoff errors. Given an $S$-bit machine point $\underline{x} \in Q^\circ$, we compute the numbers $ \partial^\alpha (\theta^\circ)(\underline{x})$ (all $\alpha \in \cM$) to within precision $\Delta_g^{-C} \Delta_\epsilon$; these numbers have magnitude at most $\Delta_g^{-C}$. This computation requires work at most $C$.

\item \textbf{Modification 11:} Just as before, we define $T: \X(E) \rightarrow \X$ by the formula $Tf = \theta^\circ \cdot T^\circ(f,\mathfrak{R}(f)) + (1 - \theta^\circ) \cdot \mathfrak{R}(f)$. We need to modify the query algorithm for $T$ to take into account roundoff error. Given an $S$-bit machine point $\underline{x} \in Q^\circ$, and given $\alpha \in \cM$, we compute the linear functional $f \mapsto \partial^\alpha(T(f))(\underline{x})$ in short form with parameters $(\Delta_g^C,\Delta_g^{-C} \Delta_\epsilon)$ in terms of the assists $\Omega$. The explanation is an obvious consequence of the Leibniz rule, since the linear maps $\mathfrak{R}$, $T^\circ$, and the cutoff function $\theta^\circ$ have been computed with parameters $(\Delta_g^C,\Delta_g^{-C} \Delta_\epsilon)$, as described in the previous bullet points.

\item For the same reason as before, we have
\[
\| Tf\|_\X^p \leq C \cdot \sum_{\xi \in \Xi} \lv \xi(f) \rv^p.
\] 
(See the proof of \eqref{m7}).
\item \textbf{Modification 12:} The estimate \eqref{m7a} no longer holds. Instead, we have
\begin{align*}
\sum_{\xi \in \Xi} \lv \xi(f) \rv^p &= \sum_{\xi^\circ \in \Xi^\circ} \lvert \xi^\circ(f,\mathfrak{R}(f)) \rvert^p \\
& \leq C \inf_{R \in \cP}  \left\{ \sum_{\xi^\circ \in \Xi^\circ} \lvert \xi^\circ(f,R) \rvert^p + \Delta_\junk^p \lvert R \rvert^p  \right\} \\
& \leq C \inf_{R \in \cP} \left\{ \|(f,R) \|_{\frac{65}{64}Q^\circ}^p + \Delta_\junk^p \lvert R \rvert^p \right\}.
\end{align*}
(See \eqref{mm1_fin}.)

For an arbitrary $F \in \X$ with $F = f$ on $E$, set $R = J_x F$, and estimate $\| F - R \|_{L^p(\frac{65}{64}Q^\circ)} \leq C \| F \|_\X$ using the Sobolev inequality. Also, by the Sobolev inequality, $\lv J_x F \rv \leq \| F \|_\X + \| F \|_{L^p(Q^\circ)}$. So the last infimum above is dominated by $C \cdot \left[ \| F \|_\X^p  + \Delta_\junk^p \| F \|_{L^p(Q^\circ)}^p \right]$. Hence,
\begin{equation*}
\sum_{\xi \in \Xi} \lv \xi(f) \rv^p \leq C \cdot \inf \left\{ \| F \|^p_\X  + \Delta_\junk^p \| F \|_{L^p(Q^\circ)}^p  : F \in \X, \; F = f \; \mbox{on} \; E \right\}.
\end{equation*}
\item Just as before, we prove that
\[
c \cdot \| f\|_{\X(E)}^p \leq \sum_{\xi \in \Xi} \lv \xi(f) \rv^p.
\]
(See \eqref{m9}.)
\end{itemize}

Recall that we have set $\Delta_g = \Delta_g(\emptyset)$, $\Delta_\epsilon = \Delta_\epsilon(\emptyset)$, and $\Delta_{\junk} = \Delta_\junk(\emptyset)$ in the above bullet points.

From \eqref{constants0}, we may impose the assumption $\Delta_\junk \leq \Delta_\junk^\circ$. Thus, from the last three bullet points we learn that
\[
c \cdot \| f \|_{\X(E)} \leq \left( \sum_{\xi \in \Xi} \lv \xi(f) \rv^p \right)^{1/p} \leq C \cdot \inf \left\{ \| F \|_\X + \Delta_\junk^\circ \cdot \| F \|_{L^p(Q^\circ)} : F \in \X, \; F = f \; \mbox{on} \; E \right\}
\]
and
\[
\| Tf \|_\X \leq C \cdot \inf \left\{ \| F \|_\X + \Delta_\junk^\circ \cdot \| F \|_{L^p(Q^\circ)} : F \in \X, \; F = f \; \mbox{on} \; E \right\},
\]
as desired (see Theorem \ref{main_thm_hom_fin}).

All of the functionals $f \mapsto \omega(f)$, $f \mapsto \xi (f)$, and $f \mapsto \partial^\alpha (Tf)(\underline{x})$ in the above bullet points, which arise in the statement of Theorem \ref{main_thm_hom_fin}, are specified with parameters $(\Delta_g^{C_0}, \Delta_g^{-C_0} \Delta_\epsilon)$ for a universal constant $C_0$. According to \eqref{constants0} and \eqref{constants1}, we may assume that $\Delta_g^\circ \leq (\Delta_g)^{C_0}$ and $\Delta_g^{-C_0} \Delta_\epsilon \leq \Delta_\epsilon^{1/2} \leq \Delta_\epsilon^\circ$. Thus, we can compute all of the functionals relevant to Theorem \ref{main_thm_hom_fin} with parameters $(\Delta_g^\circ,\Delta_\epsilon^\circ)$.

This completes the proof of Theorem \ref{main_thm_hom_fin}.

\subsection{Inhomogeneous Sobolev Spaces}
\label{inhom_finp_sec}

Once we pass from Homogeneous $L^{m,p}(\R^n)$ to Inhomogeneous $W^{m,p}(\R^n)$, the error terms $\Delta_\junk^\circ \| F \|_{L^p(Q^\circ)}$ in Theorem \ref{main_thm_hom_fin} will become irrelevant.

Our main result for inhomogeneous Sobolev spaces is Theorem \ref{main_thm2_fin} written below.

We follow the argument in Section \ref{inhom_sec}, with the following changes.

\begin{itemize}
\item \textbf{Modification 1:} We let $T^\circ, \Xi^\circ,\Omega^\circ$ be defined as in the previous section. We will use the finite-precision version of Proposition \ref{mprop}, which guarantees the following:
\begin{itemize} 
\item We list the functionals in $\Omega^\circ$. Each $\omega^\circ \in \Omega^\circ$ is specified in short form with parameters $(\Delta_g^C,\Delta_g^{-C} \Delta_\epsilon)$.
\item We list the functionals in $\Xi^\circ$. Each $\xi^\circ \in \Xi^\circ$ is specified in short form in terms of the assists $\Omega^\circ$ with parameters $(\Delta_g^C,\Delta_g^{-C} \Delta_\epsilon)$. The functionals in $\Xi^\circ$ satisfy the modified estimate \eqref{bullet4}.

\item Given an $S$-bit machine point $\underline{x} \in Q^\circ$ and given $\alpha \in \cM$, we compute the linear functional $(f,P) \mapsto \partial^\alpha ( T^\circ(f,P))(\underline{x})$ in short form with parameters $(\Delta_g^C,\Delta_g^{-C} \Delta_\epsilon)$ in terms of the assists $\Omega^\circ$, using work at most $C \log N$.
\end{itemize}

\item \textbf{Modification 2:} We introduce a cutoff function $\theta^\circ$. Let $\underline{x} \in Q^\circ$ be a given point with $S$-bit machine numbers as coordinates. We compute the numbers $ \partial^\alpha ( \theta^\circ)(\underline{x})$ (all $\alpha \in \cM$) up to an additive error of absolute value at most $\Delta_g^{-C} \Delta_\epsilon$; these numbers have absolute value at most $\Delta_g^{-C}$. This requires work at most $C$.

\item \textbf{Modification 3:} Instead of \eqref{z3}, we use the modified estimate from the fourth bullet point in the finite-precision version of Proposition \ref{mprop}. (See \eqref{bullet4}.) This gives
\begin{equation}
\label{z3_fin}
\sum_{\xi \in \Xi^\circ} \lvert \xi(f,0) \rvert^p \leq C \inf \biggl\{ \| F \|^p_{L^{m,p}(\frac{65}{64}Q^\circ)}  +  \| F  \|^p_{L^p(\frac{65}{64}Q^\circ)} : \; F \in \X, \;\; F = \; f \; \mbox{on} \; E \biggr\}. \notag{}
\end{equation}
The junk term in \eqref{bullet4} disappears because we set $P=0$. We needn't modify \eqref{z5}.

\item \textbf{Modification 4:} The rest of the content of the section is unchanged. In particular, the collections $\Xi$ and $\Omega$ consisting of linear functionals on $\X(E)$, and the linear map $T : \X(E) \mapsto \X$ are defined as before. The functionals in $\Omega$ are computed in short form with parameters $(\Delta_g^C, \Delta_g^{-C} \Delta_\epsilon)$, and the functionals in $\Xi$ are computed in short form with parameters $(\Delta_g^C,\Delta_g^{-C} \Delta_\epsilon)$ in terms of the assists $\Omega$.  Given an $S$-bit machine point $\underline{x} \in Q^\circ$ and given $\alpha \in \cM$, we can compute the linear functional $f \mapsto \partial^\alpha ( T(f))(\underline{x})$ in short form with parameters $(\Delta_g^C,\Delta_g^{-C} \Delta_\epsilon)$ in terms of the assists $\Omega$, using work at most $C \log N$.
\end{itemize}

All the functionals in the above bullet points are computed with parameters $(\Delta_g^{C_0}, \Delta_g^{-C_0} \Delta_\epsilon)$ for a universal constant $C_0$. According to \eqref{constants0} and \eqref{constants1}, we may assume that $\Delta_g^\circ \leq (\Delta_g)^{C_0}$ and $\Delta_g^{-C_0} \Delta_\epsilon \leq \Delta_\epsilon^{1/2} \leq \Delta_\epsilon^\circ$, for parameters $\Delta_g^\circ$ and $\Delta_\epsilon^\circ$ as in the statement of Theorem \ref{main_thm_hom_fin}. Thus, we can compute our functionals with parameters $(\Delta_g^\circ,\Delta_\epsilon^\circ)$ for suitable $\Delta_g^\circ$ and $\Delta_\epsilon^\circ$ (see below).

We have proven the following theorem, which is our main extension theorem for inhomogeneous Sobolev spaces in a finite-precision model of computation.

\begin{thm} \label{main_thm2_fin} There exists $C = C(m,n,p) \geq 1$ such that the following holds.

Let $\overline{S} \geq 1$ be an integer.

Assume $E \subset \frac{1}{32} Q^\circ$ satisfies $\#(E) = N \geq 2$, where $Q^\circ = [0,1)^n$. Assume that the points of $E$ have $\overline{S}$-bit machine numbers as coordinates. 

Assume that $K_1,K_2,K_{\max} \in \N$ satisfy $K_{\max} \geq C \cdot K_1 \geq  C^2 \cdot K_2 \geq C^3$.

Let $\Delta^\circ_{\min} = 2^{- K_{\max} \overline{S}}$, $\Delta^\circ_\epsilon := 2^{ - K_1 \overline{S}}$, and $\Delta^\circ_g = 2^{- K_2 \overline{S}}$.

We assume that our computer can perform arithmetic operations on $S$-bit machine numbers with precision $\Delta^\circ_{\min}$, where $S = K_{\max} \cdot \overline{S}$.

Then we can compute lists $\Omega$ and $\Xi$, consisting of linear functionals on $W^{m,p}(E) = \{ f : E \rightarrow \R\}$, with the following properties.

\begin{itemize}
\item The sum of $\depth(\omega)$ over all $\omega \in \Omega$ is bounded by $C N$. The number of functionals in $\Xi$ is at most $C N$.
\item Each functional $\xi$ in $\Xi$ has $\Omega$-assisted bounded depth.
\item The functionals $\omega \in \Omega$ and $\xi \in \Xi$ are computed in short form with parameters $(\Delta^\circ_g, \Delta^\circ_\epsilon)$.
\item For all $f \in W^{m,p}(E)$ we have
\[ C^{-1} \| f\|_{W^{m,p}(E)} \leq \left[ \sum_{\xi \in \Xi} \lvert \xi(f) \rvert^p  \right]^{1/p}  \leq C \| f \|_{W^{m,p}(E)}.\]
\end{itemize}

Moreover, there exists a linear map $T : W^{m,p}(E) \rightarrow W^{m,p}(\R^n)$ with the following properties.
\begin{itemize}
\item $T$ has $\Omega$-assisted bounded depth.
\item $T f = f$ on $E$ and 
\[  
\| Tf \|_{W^{m,p}(\R^n)} \leq C \cdot \| f \|_{W^{m,p}(E)}
\]
for all $f \in \X(E)$.
\item We produce a query algorithm that operates as follows.

Given an $S$-bit machine point $\underline{x} \in Q^\circ$ and given $\alpha \in \cM$, we compute a short form description of the $\Omega$-assisted bounded depth linear functional $ W^{m,p}(E) \ni f \mapsto \partial^\alpha ( Tf )(\underline{x})$. We compute this functional in short form with parameters $(\Delta^\circ_g,\Delta^\circ_\epsilon)$. This requires work at most $C \log N$.

\end{itemize}

The  computations above require one-time work at most $C N \log N$ in space $C N$.

\end{thm}


\end{document}